\pdfoutput=1

\documentclass[spanningrule,twoside]{cambridge7A}

\title{Introduction to Homotopy Type Theory}
\author{Egbert Rijke}
\date{\today}

\usepackage{amsmath,amsthm}
\usepackage[T1]{fontenc}
\usepackage[utf8]{inputenc}
\usepackage[english]{babel}
\usepackage{fontawesome}
\usepackage{newtxtext}
\usepackage[vvarbb]{newtxmath} 
\usepackage[final]{microtype}
\usepackage{booktabs}

\usepackage{atbegshi,picture} 
\usepackage{totcount} 
\usepackage{comment}
\usepackage{xspace}
\usepackage{ifthen}
\usepackage{mathtools}
\usepackage{tikz,tikz-cd}
\usepackage{csquotes}
\usepackage{chngcntr}
\usepackage{caption}

\usepackage{enumitem}
\newlength{\listindent}
\setlist{leftmargin=\listindent}
\setitemize[1]{leftmargin=\listindent}
\setenumerate[1]{leftmargin=\listindent}

\usepackage{bussproofs}

\usepackage{makeidx}
\usepackage[hyphens]{url}
\usepackage[style=numeric,maxnames=10,backend=biber]{biblatex}
\definecolor{darkgreen}{rgb}{0,0.4,0.4}
\usepackage[toc]{appendix}

\usepackage[pdfpagelabels,colorlinks,citecolor=orange,linkcolor=darkgreen,unicode,bookmarksopen=true,bookmarksopenlevel=2]{hyperref}
\usepackage{aliascnt}

\counterwithout{section}{chapter}
\renewcommand{\thechapter}{\Roman{chapter}}

\newtheorem{thm}{Theorem}[subsection]
\def\defthm#1#2{%
  \newaliascnt{#1}{thm}
  \newtheorem{#1}[#1]{#2}
  \aliascntresetthe{#1}
}

\defthm{cor}{Corollary}
\defthm{lem}{Lemma}
\defthm{prp}{Proposition}
\theoremstyle{definition}
\defthm{defn}{Definition}
\defthm{quasidefn}{Quasi-definition}
\defthm{rmk}{Remark}
\defthm{eg}{Example}
\defthm{axiom}{Axiom}
\defthm{postulate}{Postulate}

\usepackage[capitalize,english]{cleveref}

\crefname{thm}{Theorem}{Theorems}
\crefname{cor}{Corollary}{Corollaries}
\crefname{lem}{Lemma}{Lemmas}
\crefname{prp}{Proposition}{Propositions}
\theoremstyle{definition}
\crefname{defn}{Definition}{Definitions}
\crefname{rmk}{Remark}{Remarks}
\crefname{eg}{Example}{Examples}
\crefname{axiom}{Axiom}{Axioms}
\crefname{postulate}{Postulate}{Postulates}
\crefname{quasidefn}{Quasi-definition}{Quasi-definitions}

\input{hott}

\newtotcounter{exercisecounter}
\newcommand{\exitem}{\item\stepcounter{exercisecounter}}

\newlist{exenum}{enumerate}{1}
\setlist[exenum]{noitemsep,leftmargin=\listindent,label=\thesection.\arabic*}

\crefname{exenumi}{Exercise}{Exercises}
  
\newlist{subexenum}{enumerate}{1}
\setlist[subexenum]{noitemsep,leftmargin=\listindent,label=(\alph*),ref=\theexenumi~(\alph*)}
\crefname{subexenumi}{Exercise}{Exercises}
  
\renewenvironment{exercises}
{%
\subsection*{Exercises}%
\addcontentsline{toc}{subsection}{Exercises}%
\sectionmark{Exercises}%
\begin{exenum}}
{%
\end{exenum}}

\AtBeginShipout{\AtBeginShipoutUpperLeft{%
    \put(\dimexpr1.5cm\relax,-\paperheight+1.5cm){%
      \begin{minipage}[b]{18cm}%
        \textnormal{\small{%
          This material will be published by Cambridge University Press as
          \emph{Introduction to Homotopy Type Theory} by Egbert Rijke.
          This pre-publication version is free to view and download for
          personal use only. Not for re-distribution, re-sale or use in
          derivative works. \textcopyright{} Egbert Rijke 2020\\[\baselineskip]
          The best way to support this book while it is still being prepared for publication, is to let the author know if you find inaccuracies of any kind at \faEnvelopeO~\href{mailto:egbert.rijke@fmf.uni-lj.si}{egbert.rijke@fmf.uni-lj.si}, or at the homotopy type theory chat at \faCommentsO~\href{https://hott.zulipchat.com}{https://hott.zulipchat.com}.}}
      \end{minipage}%
    }}%
}
      
\makeindex

\begin{document}

\hypersetup{pageanchor=false}

\setcounter{tocdepth}{2}

\maketitle

\frontmatter

\hypersetup{pageanchor=true}

\tableofcontents

\chapter*{Preface}

This book started as a set of lecture notes for an Introduction to Homotopy Type Theory course that I taught at Carnegie Mellon University in the spring of 2018. The goal of that course was to give students from a wide variety of backgrounds, including mathematics, computer science, and philosophy majors, a solid foundational understanding of what univalent mathematics is about, and that is also the purpose of this book.

This book would not exist without the consistent and generous support of three people: Steve Awodey, Dan Grayson, and Andrej Bauer. Back in 2017, Steve proposed that I would teach an introductory course to homotopy type theory at Carnegie Mellon University, where I was a PhD student at the time. When the course was finished he took the course notes to Cambridge University Press to propose a book project. Throughout the entire process of writing this book I have relied on Steve's advice, and I owe him a big univalent thanks. After my graduation from CMU, I held a postdoc position at the University of Illinois at Urbana-Champaign, where Dan Grayson was my mentor. The formalization project of the book started to take off during this year, and the book has benefited from many enjoyable discussions with Dan about writing, about formalization, and univalent mathematics, and his support extended much beyond the book project. After one year in Illinois I moved to Ljubljana, where Andrej generously offered me a stable, multiple year postdoc position. What once was a set of course notes now transformed into a textbook, and some of the clearest explanations on how to think in type theory and how to define concepts correctly have their origin in the many conversations with Andrej. I cannot thank you enough.

Furthermore, this book project has benefited from insightful discussions and many people using an early version of this book. I would like to thank
William Barnett,
Katja Ber\v{c}i\v{c},
Marc Bezem,
Ulrik Buchholtz,
Ali Caglayan,
Dan Christensen,
Thierry Coquand,
Peter Dybjer,
Jacob Ender,
Mart\'in Escard\'o,
Sam van Gool,
Kerem G\"une\c{s},
Bob Harper,
Matej Jazbec,
Urban Jezernik,
Tom de Jong,
Ivan Kobe,
Anders Mortberg,
Clive Newstead, 
Charles Rezk,
Emily Riehl,
Mike Shulman,
Elif Uskuplu,
Chetan Vuppulury, and
Blaž Zupančič
for sharing their thoughts about the book, their support, finding typos, and generally for their many helpful comments on the early drafts of this book.
\\[\baselineskip]
\begin{flushright}
  Egbert Rijke \\*
  December 19th 2022, Ljubljana
\end{flushright}

\vfill
The author gratefully acknowledges the support by the Air Force Office of Scientific Research through MURI grant FA9550-15-1-0053, grants FA9550-17-1-0326 and FA9550-21-1-0024, and support by the Slovenian Research Agency research programme P1-0294.


\chapter*{Introduction}

This is an introductory textbook to univalent mathematics and homotopy type theory, a mathematical foundation that takes advantage of the structural nature of mathematical definitions and constructions. It is common in mathematical practice to consider equivalent objects to be the same, for example, to identify isomorphic groups. In set theory it is not possible to make this common practice formal. For example, there are as many distinct trivial groups in set theory as there are distinct singleton sets. Type theory, on the other hand, takes a more structural approach to the foundations of mathematics that accommodates the univalence axiom. This, however, requires us to rethink what it means for two objects to be equal.

\section*{The origins of homotopy type theory}

Homotopy type theory emerged about 10 years ago, following the discovery of the homotopy interpretation of Martin-L\"of's dependent type theory by Awodey and Warren \cite{AwodeyWarren} and independently by Voevodsky \cite{Voevodsky06}, and Voevodsky's discovery of the univalence axiom \cite{Voevodsky10}. Martin-L\"of's dependent type theory \cite{MartinLof84} is a foundational language for mathematics which is used in many of today's computer proof assistants.

In dependent type theory there are primitive objects called types and primitive objects called elements. Martin-L\"of's dependent type theory contains type formation rules for many operations that we are familiar with from traditional mathematics, such as products, sums, and inductive types such as the type of natural numbers. It is called \emph{dependent} type theory because both types and elements may be parametrized by elements of other types.

One of the distinguishing features of Martin-L\"of's dependent type theory is the \emph{identity type}. The identity type
\begin{equation*}
  \idtypevar{A}(a,b)
\end{equation*}
is an example of a dependent type because it is parametrized by two elements $a,b:A$. The elements of the identity type are called identifications, and the type theoretical way to assert that $a$ and $b$ are equal elements of type $A$ is to assert that there is an element in the identity type $\idtypevar{A}(a,b)$. In other words, to prove in type theory that two elements $a$ and $b$ in a type $A$ are equal, one has to define an identification $p:\idtypevar{A}(a,b)$. It is therefore common to write $a= b$ for the identity type $\idtypevar{A}(a,b)$, or if we want to be explicit about the ambient type we can write $a=_Ab$.

The rules for the identity type postulate that it is inductively generated by one single element
\begin{equation*}
  \refl{a}:a=_A a
\end{equation*}
parametrized by $a:A$. This raises an important question. Since the identity type is a type, it could be possible that there are many identifications $p:a=_A b$ between any two elements $a,b:A$. On the other hand, the identity type is generated inductively by only one element, i.e., by reflexivity. Is it possible to prove, using the rules of the identity type, that there is indeed at most one identification between any two elements of a type?

The pioneers of type theory have already been aware that this seemed impossible, but it was not until Hofmann and Streicher constructed the groupoid model of Martin-L\"of's dependent type theory \cite{hs:gpd-typethy} that the question was settled. In their model, types are groupoids and the type of identifications between two objects in a groupoid is the set of isomorphisms between them. Since there can be multiple isomorphisms between two objects in a groupoid, there can be multiple identifications between two objects. Furthermore, they showed that under this interpretation, the identity type indeed satisfies the rule that it is inductively generated by reflexivity. In other words, they soundly refuted the idea that identity types have at most one element. This is quite unlike ordinary mathematics, where two elements of a set are either equal or they aren't. At the end of their paper they even wondered whether there could be a similar model of type theory using higher groupoids, but the theory of higher groupoids had still been underdeveloped at that point in time. Nevertheless, the stage was set for the homotopy interpretation of type theory to emerge.

In the homotopy interpretation of type theory, we think of types as spaces. Their elements are points in those spaces, and for any two points in a space there is a \emph{space of paths} from one to the other. Analogously, for any two elements in a type there is a \emph{type of identifications} from one to the other. This way of thinking about types turned out to be very fruitful, and it opened the door to rethinking the foundations of mathematics with a prominent role for homotopy theory. It is important, however, to step back and ask the question:
\\[\baselineskip]
\emph{How is it possible that type theoretic foundations for mathematics can be so different from the usual set theoretic foundations for mathematics?}
\\[\baselineskip]
At first glance, types seem to be objects that contain stuff just like sets. There is a type $\N$ of natural numbers, a type $\Z$ of integers, standard finite types $\Fin{k}$, function types $A\to B$ and product types $A\times B$, and all of them are not very different from their set theoretic counterparts. The type of natural numbers contains the natural numbers; the type $A\to B$ contains functions from $A$ to $B$; the type $A\times B$ contains pairs $(a,b)$ consisting of elements $a:A$ and $b:B$, and so on. A big difference between type theory and set theory, however, is that in type theory types and elements are separate entities, whereas in classical set theory there is a global elementhood relation: everything in set theory is a set, and for any two sets $x$ and $y$ we can ask the question whether the proposition $x\in y$ holds. In type theory, on the other hand, there are things called types and separately there are things called elements. Furthermore, every element in type theory has a designated type. For example, $\btrue$ and $\bfalse$ are specified to be elements of type $\bool$; the numbers $\zeroN$, $1_\N$, $2_\N$, and so on, are specified to be elements of type $\N$; the successor function $\succN$ is an element of type $\N\to\N$, and the identification $\refl{a}$ is an element of type $a=_A a$. In other words, types in Martin L\"of's dependent type theory don't share their elements. Whereas in set theory a set $x$ is uniquely determined by its relation $y\in x$ with respect to all other sets $y$, types in type theory are constructed out of a small set of type forming operations, and each type forming operation comes with set of structural rules that postulate how to construct elements of that type and how elements of that type can be used. This simple change of setup has deep implications for how the foundational system works, and ultimately it opens the door to new ways of thinking about the foundations of mathematics, including the homotopy interpretation of type theory.

First of all it turns out to be extremely useful for computers if we keep track of the types of elements. Most of the widely used computer proof assistants such as Agda, Coq, or Lean, are based on a type theory. Given that every element comes equipped with a designated type, the computer can verify whether a function has been applied to elements of the correct type and outputs elements of the specified type. Such type checking algorithms are at the heart of every proof assistant, and they can be used to verify the correctness of mathematical constructions as well as proofs.

Furthermore, the identity type only compares elements in the same type. The question whether $\btrue=1_\N$ simply doesn't make sense in type theory, because $\btrue$ is a boolean and $1_\N$ is a natural number. However, if the identity type can only compare two elements in the same type, how can we hope to prove that two types are the same? This is possible with universes, which are types of which the elements themselves encode types. Given two types $A$ and $B$ in the same universe $\UU$, we can ask whether there is an identification $A=_\UU B$ in the universe $\UU$. The univalence axiom gives a characterization of this identity type. It asserts that an identification of types is equivalently described as an equivalence of types: There is an equivalence
\begin{equation*}
  (A=_\UU B) \simeq (A\simeq B)
\end{equation*}
for any two types $A$ and $B$ in the same universe $\UU$. Roughly speaking, there are as many identifications between $A$ and $B$ in $\UU$ as there are equivalences between them. For example, there are $k!$ equivalences from the standard finite type $\Fin{k}$ to itself, so the univalence axiom implies that there must be $k!$ identifications from $\Fin{k}$ to itself. Extending this example, the type $\Set_\UU$ of all sets in a universe $\UU$ should be thought of as the groupoid of all sets, because the identifications in $\Set_\UU$ correspond to equivalences between sets. This directly violates the principles of Zermelo-Fraenkel set theory, because in set theory two sets are equal if and only if they contain exactly the same elements, whereas in type theory the identifications between two types are equivalently described by equivalences between them. For example, there are many distinct singleton sets in Zermelo-Fraenkel set theory, but they are all the same in univalent mathematics.

That raises the question: Is the univalence axiom consistent? The answer is a resounding yes. Voevodsky proposed a model of Martin-L\"of's dependent type theory in which he interpreted types as Kan simplicial sets. Kan simplicial sets are the higher groupoids that Hofmann and Streicher alluded to at the end of their paper about the groupoid model of type theory. The simplicial model of type theory with the univalence axiom was later published in \cite{KapulkinLeFanuLumsdaine}.

At this point it became clear that Martin-L\"of's dependent type theory together with the univalence axiom could serve as a new foundational system for mathematics, which has homotopy theory built into its core. The famous HoTT book \cite{hottbook} was the first textbook exploring this exciting new subject. It was written during the special year 2012-2013 at the Institute for Advanced Study in Princeton as a collaborative effort by over 50 participants. The HoTT book opened up many new avenues of research, including general mathematics from a univalent point of view, (higher) group theory \cite{symmetrybook}, synthetic homotopy theory \cite{BruneriePhD}, and modal homotopy type theory \cite{Corfield,RSS}.

It has now been 10 years since the HoTT book was published. Since the publication of the HoTT book, some important open problems have been solved. It was conjectured that homotopy type theory should be modeled by all higher toposes. Higher toposes are $\infty$-categories in which the objects resemble homotopy types of spaces. The simplest $\infty$-topos is the $\infty$-category of simplicial sets, i.e., Voevodsky's model of univalence. The question whether all $\infty$-toposes model Martin-L\"of's dependent type theory with univalent universes was settled affirmatively by Michael Shulman in \cite{Shulman19}. In other words, all theorems proven in Homotopy Type Theory are valid in all $\infty$-toposes.

Another problem was whether it is possible to find a constructive model of univalence. In the simplicial model of type theory, Voevodsky used the axiom of choice to construct univalent universes. However, type theory is traditionally considered a foundation for constructive mathematics, so it was natural to ask whether it was possible to justify univalence constructively. This question was solved when Bezem, Coquand, and Huber found a model of dependent type theory and univalence in cubical sets \cite{BezemCoquandHuber,BCH19}. The cubical extension of the Agda proof assistant is based on this model. By the constructive interpretation of univalence in the cubical model it becomes possible to compute with the univalence axiom. Axel Ljungstr\"om has recently used cubical Agda to compute and formally verify that Brunerie's number \cite{BruneriePhD}, which is a number $n$ such that $\pi_4(\sphere{3})\cong\Z/n$, is $2$.


\section*{About this book}

Type theory can be confusing for people who are new to the subject, since mathematical training traditionally focuses on sets, and the differences between set theory and type theory may appear to be rather subtle to the untrained eye. The book therefore starts with a chapter that focuses on Martin-L\"of's dependent type theory, without going into homotopy theory. We first introduce the system of type dependency, gradually introduce all the type formers with their rules, and show how to get some basic mathematics off the ground in type theory.

In the second chapter we build the univalent foundation of mathematics. The central concepts of univalent mathematics are the notion of equivalence, contractibility, the hierarchy of truncation levels which includes propositions and sets, and eventually the univalence axiom. It should be noted that the univalence axiom can technically be introduced as soon as equivalences are defined, but this tends to be confusing rather than enlightening. For a good understanding of the univalence axiom, the student should have a good working knowledge of type theory, and in order to use univalence effectively they should be familiar with some of the subtleties in introducing mathematical concepts in univalent mathematics. I have three particular examples in mind: the definitions of the image of a map, surjectivity of a map, and finiteness of a type all require some type theoretical finesse. We cover those topics before we cover the univalence axiom, which will then also serve as a source of illustrative applications of the univalence axiom.

In the final chapter of the book we define the circle. The circle was the first example of a higher inductive type, and Shulman's proof using the univalence axiom of the fact that its fundamental group is $\Z$ is a pure gem \cite{LicataShulman}. It led to the realization that the methods of algebraic topology equally apply to univalent type theory, and perhaps it is because of this proof that our subject is called \emph{homotopy type theory}.

Each chapter is divided into sections that are roughly the length of one lecture, and at the end of each section there is a set of exercises. There is a total of \total{exercisecounter} exercises in this book. These exercises are an essential part of the material, and they will be referred to throughout the text. We encourage the reader to read through them, and make sure that they understand what the exercises are asking. When you see an exercise referred to in the text, we hope that you will feel encouraged to try it, or feel rewarded if you have already put in the hard work.

The more ambitious student may even try to formalize the solutions of some of the exercises in a computer proof assistant. Proof assistants provide an excellent way to become familiar with type theory, because they give instant feedback on your work. This book, including the solutions to most of its exercises, has also been formalized in the agda-unimath library \cite{Agda-UniMath}. For practice with formalization, especially the exercises in the first chapter on Martin-L\"of's dependent type theory are all very suitable.


\mainmatter

\setcounter{tocdepth}{2}

\chapter{Martin-L\"of's Dependent Type Theory}
\label{chap:type-theory}%

Dependent type theory is a formal system to organize all mathematical objects, structure, and knowledge. Dependent type theory is about types, or more generally dependent types, and their elements. There are many ways to think about type theory, types, and its elements. Types can be interpreted as sets, i.e., there is an interpretation of type theory into Zermelo-Fraenkel set theory, but there are some important differences between type theory and set theory, and the interpretation of types as sets has significant limitations. One of the differences is that in type theory, every element comes equipped with its type. We will write $a:A$ for the judgment that $a$ is an element of type $A$. This leads us to a second important difference between type theory and set theory. Set theory is axiomatized in the formal system of first order logic, whereas type theory is its own formal system. Types and their elements are constructed by following the rules of this formal system, and the only way to construct an element is to construct it as an element of a previously constructed type. The expression $a:A$ is therefore not considered to be a proposition, i.e., something which one can assert about an arbitrary element and an arbitrary type, but it is considered to be a judgment, i.e., an assessment that is part of the construction of the element $a:A$.

In type theory there is a much stronger focus on equality of elements than there is in set theory. It is said that a type is not fully understood until (i) one understands how to construct an element of the type and (ii) one understands precisely how to show that two elements of the type are equal. Equality in type theory is governed by the identity type. Unlike in classical set theory, where equality is a decidable proposition of first order logic, the \emph{type} $x=y$ of identifications of two elements $x,y:A$ is itself a type, and therefore it could possess intricate further structure. 

Dependent type theory is built up in several stages. At the first stage we give structural rules, which express the general theory of type dependency. There is no ambient deductive system of first order logic in type theory. Type theory is its own deductive system, and the structural rules are at the heart of this system. The basic operations that are governed by the structural rules are substitution and weakening operations. After the general system of dependent type theory has been set up, we introduce the ways in which we can form types. The most fundamental class of types are dependent function types, or $\Pi$-types. They are used for practically everything. Next, we introduce the type of natural numbers, where we use type-dependency to formulate a type-theoretic version of the induction principle. By the type-theoretic nature of this induction principle, it can be used in two ways: it can be used to construct the many familiar operations on $\N$, such as addition and multiplication, and it can also be used to prove properties about those operations.

The next idea is that we can consider induction principles for many other types as well. This leads to the idea of more general inductive types. In \cref{sec:inductive} we introduce the unit type, the empty type, the booleans, coproducts, dependent pair types, and cartesian products. All of these are examples of inductive types, and their induction principles can be used to construct the basic operations on them, as well as to prove properties about those operations.

Then we come to the most characteristic ingredient of Martin L\"of's dependent type theory: the identity type. The identity type $x=_Ay$ is an example of a \emph{dependent} type, because it is indexed by $x,y:A$, and it is inductively generated by the reflexivity element $\refl{x}:x=_Ax$. The catch is, however, that the identity type $x=_Ay$ is just another type, and it could potentially have many different elements.

The last class of types that we introduce are universes. Universes are type families that are closed under the operations of type theory: $\Pi$-types, $\Sigma$-types, identity types, and so on. Universes play a fundamental role in the theory. One important reason for introducing universes is that they can be used to define type families over inductive types via their induction principles. For example, this allows us to define the ordering relations $\leq$ and $<$ on the natural numbers. We will also use the universes to show the Peano axioms asserting that $\succN$ is injective, and that $\zeroN$ is not a successor.

In the final two sections of this chapter, we start developing mathematics in type theory. In \cref{sec:modular-arithmetic} we study the Curry-Howard interpretation, and use it to develop modular arithmetic in type theory. In \cref{sec:decidability} we study the concept of decidability, and use it to obtain basic theorems in elementary number theory, such as the well-ordering theorem, the construction of the greatest common divisor, and the infinitude of primes. Both of these sections can be viewed as tutorials in type theory, designed to give you some practical experience with type theory before diving into the intricacies of the univalent foundations of mathematics.

\section{Dependent type theory}%
\label{sec:dtt}%
\index{dependent type theory|(}%

Dependent type theory is a system of inference rules that can be combined to make \emph{derivations}. In these derivations, the goal is often to construct an element of a certain type. Such an element can be a function if the type of the constructed element is a function type; a proof of a property if the type of the constructed element is a proposition; but it can also be an identification if the type of the constructed element is an identity type, and so on. In some respect, a type is just a collection of mathematical objects and constructing elements of a type is the everyday mathematical task or challenge. The system of inference rules that we call type theory offers a principled way of engaging in mathematical activity.

\subsection{Judgments and contexts in type theory}%
\index{judgment|(}%
\index{context|(}%

A mathematical argument or construction consists of a sequence of deductive steps, each one using finitely many premises in order to get to the next stage in the proof or construction. Such steps can be represented by \define{inference rules}\index{inference rule|see {rule}}, which are written in the form
\begin{prooftree}
  \AxiomC{$\mathcal{H}_1$\quad $\mathcal{H}_2$ \quad \dots \quad $\mathcal{H}_n$}
  \UnaryInfC{$\mathcal{C}$.}
\end{prooftree}
Inference rules contain above the horizontal line\index{horizontal line|see {inference rule}} a finite list $\mathcal{H}_1$, $\mathcal{H}_2$, \dots, $\mathcal{H}_n$ of \emph{judgments} for the \define{premises}\index{inference rule!premises|textbf}\index{premise of an inference rule|textbf}, and below the horizontal line a single judgment $\mathcal{C}$ for the \define{conclusion}\index{inference rule!conclusion|textbf}\index{conclusion of an inference rule|textbf}. The system of dependent type theory is described by a set of such inference rules.

A straightforward example of an inference rule that we will encounter in \cref{sec:pi} when we introduce function types\index{function type}, is the inference rule
\begin{prooftree}
  \AxiomC{$\Gamma\vdash a:A$}
  \AxiomC{$\Gamma\vdash f:A\to B$}
  \BinaryInfC{$\Gamma\vdash f(a):B$.}
\end{prooftree}
This rule asserts that in any context $\Gamma$ we may use an element $a:A$ and a function $f:A\to B$ to obtain an element $f(a):B$. Each of the expressions
\begin{align*}
  \Gamma & \vdash a :A \\*
  \Gamma & \vdash f : A \to B \\*
  \Gamma & \vdash f(a):B
\end{align*}
are examples of judgments.

\begin{defn}\label{defn:judgments}
  There are four kinds of \define{judgments} in Martin-L\"of's dependent type theory:
  \begin{enumerate}
  \item \emph{$A$ is a (well-formed) \define{type} in context $\Gamma$.}
    \index{well-formed type}\index{type}
    We express this judgment as\index{Gamma turnstile A type@{$\Gamma\vdash A~\type$}}\index{judgment!Gamma turnstile A type@{$\Gamma\vdash A~\type$}}
    \begin{equation*}
      \Gamma\vdash A~\type.
    \end{equation*}
  \item \emph{$A$ and $B$ are \define{judgmentally equal types} in context $\Gamma$.}
    \index{judgmental equality!of types} We express this judgment as\index{Gamma turnstile A is B type@{$\Gamma\vdash A\jdeq B~\type$}}\index{judgment!Gamma turnstile A is B type@{$\Gamma\vdash A\jdeq B~\type$}}
    \begin{equation*}
      \Gamma\vdash A \jdeq B~\type.
    \end{equation*}
  \item \emph{$a$ is an \define{element}\index{element|textbf} of type $A$ in context $\Gamma$.} We express this judgment as\index{Gamma turnstile a in A@{$\Gamma\vdash a:A$}}\index{judgment!Gamma turnstile a in A@{$\Gamma\vdash a:A$}}
    \begin{equation*}
      \Gamma \vdash a:A.
    \end{equation*}
  \item \emph{$a$ and $b$ are \define{judgmentally equal elements} of type $A$ in context $\Gamma$.}\index{judgmental equality!of elements} We express this judgment as\index{Gamma turnstile a is b in A@{$\Gamma\vdash a\jdeq b:A$}}\index{judgment!Gamma turnstile a is b in A@{$\Gamma\vdash a\jdeq b:A$}}
    \begin{equation*}
      \Gamma\vdash a\jdeq b:A.
    \end{equation*}
  \end{enumerate}
\end{defn}

We see that any judgment is of the form $\Gamma\vdash\mathcal{J}$, consisting of a \emph{context} $\Gamma$ and a \emph{judgment thesis} $\mathcal{J}$ asserting either that $A$ is a type, that $A$ and $B$ are equal types, that $a$ is an element of type $A$, or that $a$ and $b$ are equal elements of type $A$. The role of a context is to declare what \define{hypothetical elements}\index{hypothetical elements|textbf} are assumed, along with their types. Hypothetical elements are commonly called \define{variables}\index{variable|textbf}.

\begin{defn}\label{defn:context}
  A \define{context}\index{context|textbf} is a finite list of \define{variable declarations}\index{variable declaration|textbf}
\begin{equation}\label{eq:context}
x_1:A_1,~x_2:A_2(x_1),~\ldots,~x_n:A_n(x_1,\ldots,x_{n-1})
\end{equation}
satisfying the condition that for each $1\leq k\leq n$ we can derive the judgment
\begin{equation*}
  x_1:A_1,~\ldots,~x_{k-1}:A_{k-1}(x_1,\ldots,x_{k-2})\vdash A_k(x_1,\ldots,x_{k-1})~\type,
\end{equation*}
using the inference rules of type theory. We may use variable names other than $x_1,\ldots,x_n$, as long as no variable is declared more than once.
\end{defn}

The condition in \cref{defn:context} that each of the hypothetical elements is assigned a type, is checked recursively. In other words, to check that a list of variable declarations as in \cref{eq:context} is a context, one starts on the left and works their way to the right, verifying that each hypothetical elements $x_k$ is assigned a type. 

Note that there is a context of length $0$, the \define{empty context}\index{context!empty context|textbf}\index{empty context|textbf}, which declares no variables. This context satisfies the requirement in \cref{defn:context} vacuously. A list of variable declarations $x_1:A_1$ of length one is a context if and only if $A_1$ is a type in the empty context. We will soon encounter the type $\N$ of natural numbers\index{natural numbers}, which is an example of a type in the empty context.

The next case is that a list of variable declarations $x_1:A_1,~x_2:A_2(x_1)$ of length two is a context if and only if $A_1$ is a type in the empty context, and $A_2(x_1)$ is a type in context $x_1:A_1$. This process repeats itself for longer contexts.
\index{judgment|)}
\index{context|)}

\subsection{Type families}
\index{type family|(}
It is a feature of \emph{dependent} type theory that all judgments are context dependent, and indeed that even the types of the variables in a context may depend on any previously declared variables. For example, if $n$ is a natural number and we know from the context that $n$ is prime, then we don't have enough information yet to decide whether or not $n$ is odd. However, if we also know from the context that $n+2$ is prime, then we can derive that $n$ must be odd. Context dependency is everywhere -- not only in mathematics, but also in language and in everyday life -- and it gives rise to the notion of \emph{type families} and their \emph{sections}.

\begin{defn}
  Consider a type $A$ in context $\Gamma$. A \define{family}\index{family of types|see{type family}}\index{type family|textbf}\index{family of types|textbf} of types over $A$ in context $\Gamma$ is a type $B(x)$ in context $\Gamma,x:A$. In other words, in the situation where
\begin{equation*}
  \Gamma,~x:A\vdash B(x)~\type,
\end{equation*}
we say that $B$ is a family of types over $A$ in context $\Gamma$. Alternatively, we say that $B(x)$ is a type \define{indexed}\index{indexed type|textbf}\index{type!indexed type|textbf} by $x:A$, in context $\Gamma$.
\end{defn}

We think of a type family $B$ over $A$ in context $\Gamma$ as a type $B(x)$ varying along $x:A$. A basic example of a type family occurs when we introduce \emph{identity types}\index{identity type} in \cref{sec:identity}. They are introduced as follows:
\begin{prooftree}
  \AxiomC{$\Gamma\vdash a:A$}
  \UnaryInfC{$\Gamma,~x:A\vdash a=x~\type$.}
\end{prooftree}
This rule asserts that given an element $a:A$ in context $\Gamma$, we may form the type $a=x$ in context $\Gamma,~x:A$. The type $a=x$ in context $\Gamma,~x:A$ is an example of a type family over $A$ in context $\Gamma$.

\begin{defn}
Consider a type family $B$ over $A$ in context $\Gamma$. A \define{section}\index{section of a type family} of the family $B$ over $A$ in context $\Gamma$ is an element of type $B(x)$ in context $\Gamma,x:A$, i.e., in the judgment
\begin{equation*}
  \Gamma,~x:A\vdash b(x):B(x)
\end{equation*}
we say that $b$ is a section of the family $B$ over $A$ in context $\Gamma$. Alternatively, we say that $b(x)$ is an element of type $B(x)$ \define{indexed}\index{indexed element|textbf}\index{element!indexed element|textbf} by $x:A$ in context $\Gamma$.
\end{defn}

Note that in the above situations $A$, $B$, and $b$ also depend on the variables declared in the context $\Gamma$, even though we have not explicitly mentioned them. It is indeed common practice to not mention every variable in the context $\Gamma$ in such situations.
\index{type family|)}

\subsection{Inference rules}\label{sec:rules}

We are now ready to present the system of inference rules that underlies dependent type theory. These rules are known as the \define{structural rules} of type theory, since they establish the basic mathematical framework for type dependency. There are six sets of inference rules:
\begin{enumerate}
\item Rules about the formation of contexts, types, and their elements
\item Rules postulating that judgmental equality is an equivalence relation.
\item Variable conversion rules.
\item Substitution rules.
\item Weakening rules.
\item The generic element.
\end{enumerate}

\subsubsection*{Rules about the formation of contexts, types, and their elements}
In the definition of well-formed contexts, types, and elements we specified that for a type $B(x)$ to be well-formed in context $\Gamma,x:A$, it must be the case that $A$ is a well-formed type in context $\Gamma$. The following rules follow from the presuppositions about contexts, types, and their elements, and may be used freely in derivations:

\begin{center}
  \begin{minipage}{.35\textwidth}
    \begin{prooftree}
      \AxiomC{$\Gamma,x:A\vdash B(x)~\type$}
      \UnaryInfC{$\Gamma\vdash A~\type$}
    \end{prooftree}
    
    \begin{prooftree}
      \AxiomC{$\Gamma\vdash a:A$}
      \UnaryInfC{$\Gamma\vdash A~\type$}
    \end{prooftree}
  \end{minipage}
  \begin{minipage}{.25\textwidth}
    \begin{prooftree}
      \AxiomC{$\Gamma\vdash A\jdeq B~\type$}
      \UnaryInfC{$\Gamma\vdash A~\type$}
    \end{prooftree}
    
    \begin{prooftree}
      \AxiomC{$\Gamma\vdash a\jdeq b:A$}
      \UnaryInfC{$\Gamma\vdash a:A$}
    \end{prooftree}
  \end{minipage}
  \begin{minipage}{.25\textwidth}
    \begin{prooftree}
      \AxiomC{$\Gamma\vdash A\jdeq B~\type$}
      \UnaryInfC{$\Gamma\vdash B~\type$}
    \end{prooftree}
    
    \begin{prooftree}
      \AxiomC{$\Gamma\vdash a\jdeq b:A$}
      \UnaryInfC{$\Gamma\vdash b:A$}
    \end{prooftree}
  \end{minipage}
\end{center}

\subsubsection*{Judgmental equality is an equivalence relation}

\index{rules!for type dependency!judgmental equality is an equivalence relation|(}
The rules postulating that judgmental equality on types and on elements is an equivalence relation simply postulate that these relations are reflexive, symmetric, and transitive\index{judgmental equality!is an equivalence relation}:
\begin{center}
  \begin{small}
    \begin{minipage}{.22\textwidth}
      \begin{center}
        \begin{prooftree}
          \AxiomC{$\Gamma\vdash A~\textrm{type}$}
          \UnaryInfC{$\Gamma\vdash A\jdeq A~\textrm{type}$}
        \end{prooftree}
      \end{center}
    \end{minipage}
    \begin{minipage}{.28\textwidth}
      \begin{center}
        \begin{prooftree}
          \AxiomC{$\Gamma\vdash A\jdeq B~\textrm{type}$}
          \UnaryInfC{$\Gamma\vdash B\jdeq A~\textrm{type}$}
        \end{prooftree}
      \end{center}
    \end{minipage}
    \begin{minipage}{.48\textwidth}
      \begin{prooftree}
        \AxiomC{$\Gamma\vdash A\jdeq B~\textrm{type}$}
        \AxiomC{$\Gamma\vdash B\jdeq C~\textrm{type}$}
        \BinaryInfC{$\Gamma\vdash A\jdeq C~\textrm{type}$}
      \end{prooftree}
    \end{minipage}
    \\
    \bigskip
    \begin{minipage}{.22\textwidth}
      \begin{prooftree}
        \AxiomC{$\Gamma\vdash a:A$}
        \UnaryInfC{$\Gamma\vdash a\jdeq a : A$}
      \end{prooftree}
    \end{minipage}
    \begin{minipage}{.28\textwidth}
      \begin{prooftree}
        \AxiomC{$\Gamma\vdash a\jdeq b:A$}
        \UnaryInfC{$\Gamma\vdash b\jdeq a: A$}
      \end{prooftree}
    \end{minipage}
    \begin{minipage}{.48\textwidth}
      \begin{prooftree}
        \AxiomC{$\Gamma\vdash a\jdeq b : A$}
        \AxiomC{$\Gamma\vdash b\jdeq c: A$}
        \BinaryInfC{$\Gamma\vdash a\jdeq c: A$.}
      \end{prooftree}
    \end{minipage}
  \end{small}
\end{center}
\index{rules!for type dependency!judgmental equality is an equivalence relation|)}

\subsubsection*{Variable conversion rules}
\index{rules!for type dependency!variable conversion|(}
The \define{variable conversion rules}\index{judgmental equality!conversion rules}\index{variable conversion rules|textbf}\index{conversion rule!variable|textbf}\index{rules!for type dependency!variable conversion|textbf} are rules postulating that we can convert the type of a variable to a judgmentally equal type. The first variable conversion rule states that
\begin{prooftree}
\AxiomC{$\Gamma\vdash A\jdeq A'~\textrm{type}$}
\AxiomC{$\Gamma,~x:A,~\Delta\vdash B(x)~\type$}
\BinaryInfC{$\Gamma,~x:A',~\Delta\vdash B(x)~\type$.}
\end{prooftree}
In this conversion rule, the context $\Gamma,~x:A,~\Delta$ is just any extension of the context $\Gamma,~x:A$, i.e., it is a context of the form
\begin{equation*}
  x_1:A_1,~\ldots,~x_{n-1}:A_{n-1},~x:A,~x_{n+1}:A_{n+1},~\ldots,~x_{n+m}:A_{n+m}.
\end{equation*}

Similarly, there are variable conversion rules for judgmental equality of types, for elements, and for judgmental equality of elements. To avoid having to state essentially the same rule four times, we state all four variable conversion rules at once using a \emph{generic judgment thesis} $\mathcal{J}$, which can be any of the four kinds described in \cref{defn:judgments}:
\begin{prooftree}
\AxiomC{$\Gamma\vdash A\jdeq A'~\textrm{type}$}
\AxiomC{$\Gamma,~x:A,~\Delta\vdash \mathcal{J}$}
\BinaryInfC{$\Gamma,~x:A',~\Delta\vdash \mathcal{J}$.}
\end{prooftree}
An analogous \emph{element conversion rule}, stated in \cref{ex:term_conversion}, converting the type of an element to a judgmentally equal type, is derivable using the rules from the rules presented in this section.
\index{rules!for type dependency!variable conversion|)}

\subsubsection*{Substitution}
\index{substitution|(}\index{rules!for type dependency!rules for substitution|(}

Consider an element $f(x):B(x)$ indexed by $x:A$ in context $\Gamma$, and suppose we also have an element $a:A$. Then we can simultaneously substitute $a$ for all occurrences of $x$ in $f(x)$ to obtain a new element $f[a/x]$, which has type $B[a/x]$. A precise definition of substitution requires us to get too deep into the theory of the syntax of type theory, but a mathematician is of course no stranger to substitution. For example, substituting $0$ for $x$ in the polynomial
\begin{equation*}
  1+x+x^2+x^3
\end{equation*}
results in the number $1+0+0^2+0^3$, which can be computed to the value $1$.

Type theoretic substitution is similar. Type theoretic substitution is in fact a bit more general than what we have described above. Suppose we have a type
\begin{equation*}
  \Gamma,~x:A,~y_{1}:B_{1},~\ldots,~y_{n}:B_{n}\vdash C~\textrm{type}
\end{equation*}
and an element $a:A$ in context $\Gamma$. Then we can simultaneously substitute $a$ for all occurrences of $x$ in the types $B_1,\ldots,B_n$ and $C$, to obtain
\begin{equation*}
  \Gamma,~y_{1}:B_{1}[a/x],~\ldots,~y_{n}:B_{n}[a/x]\vdash C[a/x]~\mathrm{type}.
\end{equation*}
Note that the variables $y_{1},~\ldots,y_{n}$ are assigned new types after performing the substitution of $a$ for $x$. Similarly, we can substitute $a$ for $x$ in an element $c:C$ to obtain the element $c[a/x]:C[a/x]$, and we can substitute $a$ for $x$ in a judgmental equality thesis, either of types or elements, by simply substituting on both sides of the equation. The \define{substitution rule} are therefore stated using a generic judgment $\mathcal{J}$:
\begin{prooftree}
\AxiomC{$\Gamma\vdash a:A$}
\AxiomC{$\Gamma,~x:A,~\Delta\vdash \mathcal{J}$}
\RightLabel{$S$.}
\BinaryInfC{$\Gamma,~\Delta[a/x]\vdash \mathcal{J}[a/x]$}
\end{prooftree}
Furthermore, we add two more `congruence rules' for substitution, postulating that substitution by judgmentally equal elements results in judgmentally equal types and elements:
\begin{prooftree}
\AxiomC{$\Gamma\vdash a\jdeq a':A$}
\AxiomC{$\Gamma,~x:A,~\Delta\vdash B~\type$}
\BinaryInfC{$\Gamma,~\Delta[a/x]\vdash B[a/x]\jdeq B[a'/x]~\type$}
\end{prooftree}
\begin{prooftree}
\AxiomC{$\Gamma\vdash a\jdeq a':A$}
\AxiomC{$\Gamma,~x:A,~\Delta\vdash b:B$}
\BinaryInfC{$\Gamma,~\Delta[a/x]\vdash b[a/x]\jdeq b[a'/x]:B[a/x]$.}
\end{prooftree}
To see that these rules make sense, we observe that both $B[a/x]$ and $B[a'/x]$ are types in context $\Delta[a/x]$, provided that $a\jdeq a'$. This is immediate by recursion on the length of $\Delta$.

\begin{defn}
  When $B$ is a family of types over $A$ in context $\Gamma$, and if we have $a:A$, then we also say that $B[a/x]$ is the \define{fiber}\index{type family!fiber of a type family|textbf}\index{fiber of a type family|textbf} of $B$ at $a$. We will usually write $B(a)$ for the fiber of $B$ at $a$.

  When $b$ is a section of the family $B$ over $A$ in context $\Gamma$, we call the element $b[a/x]$ the \define{value} of $b$ at $a$. Again, we will usually write $b(a)$ for the value of $b$ at $a$.
\end{defn}
\index{substitution|)}\index{rules!for type dependency!rules for substitution|)}

\subsubsection*{Weakening}
\index{weakening|(}\index{rules!for type dependency!rules for weakening|(}
If we are given a type $A$ in context $\Gamma$, then any judgment made in a longer context $\Gamma,~\Delta$ can also be made in the context $\Gamma,~x:A,~\Delta$, for a fresh variable $x$. The \define{weakening rule}\index{weakening} asserts that weakening by a type $A$ in context preserves well-formedness and judgmental equality of types and elements.
\begin{prooftree}
\AxiomC{$\Gamma\vdash A~\textrm{type}$}
\AxiomC{$\Gamma,~\Delta\vdash \mathcal{J}$}
\RightLabel{$W$.}
\BinaryInfC{$\Gamma,~x:A,~\Delta \vdash \mathcal{J}$}
\end{prooftree}
This process of expanding the context by a fresh variable of type $A$ is called \define{weakening} (by $A$).

In the simplest situation where weakening applies, we have two types $A$ and $B$ in context $\Gamma$. Then we can weaken $B$ by $A$ as follows
\begin{prooftree}
  \AxiomC{$\Gamma\vdash A~\textrm{type}$}
  \AxiomC{$\Gamma\vdash B~\textrm{type}$}
  \RightLabel{$W$}
  \BinaryInfC{$\Gamma,~x:A\vdash B~\type$}
\end{prooftree}
in order to form the type $B$ in context $\Gamma,~x:A$. The type $B$ in context $\Gamma,~x:A$ is called the \define{constant family}\index{type family!constant family|textbf}\index{constant family|textbf} $B$, or the \define{trivial family}\index{type family!trivial family|textbf}\index{trivial family|textbf} $B$.
\index{weakening|)}\index{rules!for type dependency!rules for weakening|)}

\subsubsection*{The generic elements}
If we are given a type $A$ in context $\Gamma$, then we can weaken $A$ by itself to obtain that $A$ is a type in context $\Gamma,~x:A$. The rule for the \define{generic element}\index{generic element|textbf}\index{rules!for type dependency!generic element|textbf} now asserts that any hypothetical element $x:A$ in the context $\Gamma,~x:A$ is also an element of type $A$ in context $\Gamma,~x:A$.
\begin{prooftree}
\AxiomC{$\Gamma\vdash A~\textrm{type}$}
\RightLabel{$\delta$.}
\UnaryInfC{$\Gamma,~x:A\vdash x:A$}
\end{prooftree}
This rule is also known as the \define{variable rule}\index{variable rule|textbf}\index{rules!for type dependency!variable rule|textbf}. One of the reasons for including the generic element is to make sure that the variables declared in a context---i.e., the hypothetical elements---are indeed \emph{elements}. It also provides the \emph{identity function}\index{identity function} on the type $A$ in context $\Gamma$.

\subsection{Derivations}\label{sec:derivations}

\index{derivation|(}
A \define{derivation}\index{derivation|textbf} in type theory is a finite tree in which each node is a valid rule of inference. At the root of the tree we find the conclusion, and in the leaves of the tree we find the hypotheses. We give two examples of derivations: a derivation showing that any variable can be changed to a fresh one, and a derivation showing that any two variables that do not mutually depend on one another can be swapped in order.

Given a derivation with hypotheses $\mathcal{H}_1,\ldots,\mathcal{H}_n$ and conclusion $\mathcal{C}$, we can form a new inference rule
\begin{prooftree}
  \AxiomC{$\mathcal{H}_1$}
  \AxiomC{$\cdots$}
  \AxiomC{$\mathcal{H}_n$}
  \TrinaryInfC{$\mathcal{C}$.}
\end{prooftree}
Such a rule is called \define{derivable}, because we have a derivation for it. In order to keep proof trees reasonably short and manageable, we use the convention that any derived rules can be used in future derivations.

\subsubsection*{Changing variables}

\index{change of variables}
Variables can always be changed to fresh variables. We show that this is the case by showing that the inference rule\index{rules!for type dependency!change of variables}
\begin{prooftree}
  \AxiomC{$\Gamma,~x:A,~\Delta\vdash \mathcal{J}$}
  \RightLabel{$x'/x$}
  \UnaryInfC{$\Gamma,~x':A,~\Delta[x'/x]\vdash \mathcal{J}[x'/x]$}
\end{prooftree}
is derivable, where $x'$ is a variable that does not occur in the context $\Gamma,~x:A,~\Delta$. 

Indeed, we have the following derivation using substitution, weakening, and the generic element:
\begin{prooftree}
  \AxiomC{$\Gamma\vdash A~\type$}
  \RightLabel{$\delta$}
  \UnaryInfC{$\Gamma,~x':A\vdash x':A$}
  \AxiomC{$\Gamma\vdash A~\type$}
  \AxiomC{$\Gamma,~x:A,~\Delta\vdash \mathcal{J}$}
  \RightLabel{$W$}
  \BinaryInfC{$\Gamma,~x':A,~x:A,~\Delta\vdash \mathcal{J}$}
  \RightLabel{$S$.}
  \BinaryInfC{$\Gamma,~x':A,~\Delta[x'/x]\vdash \mathcal{J}[x'/x]$}
\end{prooftree}
In this derivation it is the application of the weakening rule where we have to check that $x'$ does not occur in the context $\Gamma,~x:A,~\Delta$.

\subsubsection*{Interchanging variables}

The \define{interchange rule}\index{rules!for type dependency!interchange}\index{interchange rule} states that if we have two types $A$ and $B$ in context $\Gamma$, and we make a judgment in context $\Gamma,~x:A,~y:B,~\Delta$, then we can make that same judgment in context $\Gamma,~y:B,~x:A,~\Delta$ where the order of $x:A$ and $y:B$ is swapped. More formally, the interchange rule is the following inference rule
\begin{prooftree}
\AxiomC{$\Gamma\vdash B~\textrm{type}$}
\AxiomC{$\Gamma,~x:A,~y:B,~\Delta\vdash \mathcal{J}$}
\BinaryInfC{$\Gamma,~y:B,~x:A,~\Delta\vdash \mathcal{J}$.}
\end{prooftree}
Just as the rule for changing variables, we claim that the interchange rule is a derivable rule.

The idea of the derivation for the interchange rule is as follows: If we have a judgment
\begin{equation*}
  \Gamma,~x:A,~y:B,~\Delta\vdash\mathcal{J},
\end{equation*}
then we can change the variable $y$ to a fresh variable $y'$ and weaken the judgment to obtain the judgment
\begin{equation*}
  \Gamma,~y:B,~x:A,~y':B,~\Delta[y'/y]\vdash\mathcal{J}[y'/y].
\end{equation*}
Now we can substitute $y$ for $y'$ to obtain the desired judgment $\Gamma,~y:B,~x:A,~\Delta\vdash\mathcal{J}$. The formal derivation is as follows:
\begin{small}
  \begin{prooftree}
    \AxiomC{$\Gamma\vdash B~\textrm{type}$}
    \UnaryInfC{$\Gamma,~y:B\vdash y:B$}
    \UnaryInfC{$\Gamma,~y:B,~x:A\vdash y:B$}
    \AxiomC{$\Gamma\vdash B~\textrm{type}$}
    \AxiomC{$\Gamma,~x:A,~y:B,~\Delta\vdash \mathcal{J}$}
    \UnaryInfC{$\Gamma,~x:A,~y':B,~\Delta[y'/y]\vdash \mathcal{J}[y'/y]$}
    \BinaryInfC{$\Gamma,~y:B,~x:A,~y':B,~\Delta[y'/y]\vdash \mathcal{J}[y'/y]$}
    \BinaryInfC{$\Gamma,~y:B,~x:A,~\Delta\vdash \mathcal{J}$}
  \end{prooftree}%
\end{small}
\index{derivation|)}

\begin{exercises}
  \exitem \label{ex:term_conversion}
  \begin{subexenum}
  \item Give a derivation for the following \define{element conversion rule}\index{element conversion rule|textbf}\index{rules!for type dependency!element conversion|textbf}\index{conversion rule!element|textbf}:
    \begin{prooftree}
      \AxiomC{$\Gamma\vdash A\jdeq A'~\textrm{type}$}
      \AxiomC{$\Gamma\vdash a:A$}
      \BinaryInfC{$\Gamma\vdash a:A'$.}
    \end{prooftree}
  \item Give a derivation for the following \define{congruence rule} for element conversion:
    \begin{prooftree}
      \AxiomC{$\Gamma\vdash A\jdeq A'~\textrm{type}$}
      \AxiomC{$\Gamma\vdash a\jdeq b:A$}
      \BinaryInfC{$\Gamma\vdash a\jdeq b:A'$.}
    \end{prooftree}
  \end{subexenum}
\end{exercises}
\index{dependent type theory|)}


\section{Dependent function types}
\label{sec:pi}

\index{Pi-type@{$\Pi$-type}|see {dependent function type}}
\index{dependent function type|(}
A fundamental concept of dependent type theory is that of a dependent function. A dependent function is a function of which the type of the output may depend on the input. For example, when we concatenate a vector of length $m$ with a vector of length $n$, we obtain a vector of length $m+n$. Dependent functions are a generalization of ordinary functions, because an ordinary function $f:A\to B$ is a function of which the output $f(x)$ has type $B$ regardless of the value of $x$.

\subsection{The rules for dependent function types}
Consider a section $b$ of a family $B$ over $A$ in context $\Gamma$, i.e., consider
\begin{equation*}
  \Gamma,x:A\vdash b(x):B(x).
\end{equation*}
From one point of view, such a section $b$ is an operation or assignment $x\mapsto b(x)$, or a program\index{program}, that takes as input $x:A$ and produces a term $b(x):B(x)$. From a more mathematical point of view we see $b$ as a choice of an element of each $B(x)$. In other words, we may see $b$ as a function that takes $x:A$ to $b(x):B(x)$. Note that the type $B(x)$ of the output may depend on $x:A$. The assignment $x\mapsto b(x)$ is in this sense a \emph{dependent} function. The type of all such dependent functions is called the \define{dependent function type}, and we will write
\begin{equation*}
  \prd{x:A}B(x)
\end{equation*}
for the type of dependent functions. There are four principal rules for $\Pi$-types:
\begin{enumerate}
\item The \emph{formation rule}, which tells us how we may form dependent function types.
\item The \emph{introduction rule}, which tells us how to introduce new terms of dependent function types.
\item The \emph{elimination rule}, which tells us how to use arbitrary terms of dependent function types.
\item The \emph{computation rules}, which tell us how the introduction and elimination rules interact. These computation rules guarantee that every term of a dependent function type is indeed a dependent function taking the values by which it is defined.
\end{enumerate}
In the cases of the formation rule, the introduction rule, and the elimination rule, we also need rules that assert that all the constructions respect judgmental equality. Those rules are called \define{congruence rules}, and they are part of the specification of dependent function types.

\subsubsection{The $\Pi$-formation rule}
\index{dependent function type!formation rule|textbf}
The \define{$\Pi$-formation rule} tells us how $\Pi$-types are constructed. The idea of $\Pi$-types is that $\prd{x:A}B(x)$ is a type of \define{dependent functions}\index{dependent function type|textbf}, for any type family $B$ of types over $A$, so the $\Pi$-formation rule is as follows:\index{rules!for dependent function types!formation|textbf}\index{P (x:A) B(x)@{$\prd{x:A}B(x)$}|see{dependent function type}}\index{P (x:A) B(x)@{$\prd{x:A}B(x)$}|textbf}
\begin{prooftree}
\AxiomC{$\Gamma,x:A\vdash B(x)~\type$}
\RightLabel{$\Pi$.}
\UnaryInfC{$\Gamma\vdash \prd{x:A}B(x)~\type$}
\end{prooftree}
This rule simply states that in order to form the type $\prd{x:A}B(x)$ in context $\Gamma$, we must have a type family $B$ over $A$ in context $\Gamma$.

We also require that the operation of forming dependent function types respects judgmental equality. This is postulated in the \define{congruence rule} for $\Pi$-types:
\index{rules!for dependent function types!congruence|textbf}
\index{dependent function type!congruence rule|textbf}
\begin{prooftree}
\AxiomC{$\Gamma\vdash A\jdeq A'~\type$}
\AxiomC{$\Gamma,x:A\vdash B(x)\jdeq B'(x)~\textrm{type}$}
\RightLabel{$\Pi$-eq.}
\BinaryInfC{$\Gamma\vdash \prd{x:A}B(x)\jdeq\prd{x:A'}B'(x)~\type$}
\end{prooftree}

\subsubsection{The $\Pi$-introduction rule}
The introduction rule for dependent functions tells us how we may construct dependent functions of type $\prd{x:A}B(x)$. The idea is that a dependent function $f:\prd{x:A}B(x)$ is an operation that takes an $x:A$ to $f(x):B(x)$. Hence the introduction rule of dependent functions postulates that, in order to construct a dependent function one has to construct a term $b(x):B(x)$ indexed by $x:A$ in context $\Gamma$, i.e.:
\begin{prooftree}
  \AxiomC{$\Gamma,x:A \vdash b(x) : B(x)$}
  \RightLabel{$\lambda$.}
  \UnaryInfC{$\Gamma\vdash \lam{x}b(x) : \prd{x:A}B(x)$}
\end{prooftree}
This introduction rule%
\index{dependent function type!introduction rule|see {$\lambda$-abstraction}}
for dependent functions is also called the \define{$\lambda$-abstraction rule}%
\index{lambda-abstraction@{$\lambda$-abstraction}|textbf}%
\index{rules!for dependent function types!lambda-abstraction@{$\lambda$-abstraction}|textbf}%
\index{dependent function type!lambda-abstraction@{$\lambda$-abstraction}|textbf}, and we also say that the $\lambda$-abstraction $\lam{x}b(x)$ \define{binds} the variable $x$ in $b$. Just like ordinary mathematicians, we will sometimes write $x\mapsto b(x)$ for a function $\lam{x}b(x)$. The map $n\mapsto n^2$ is an example.

We will also require that $\lambda$-abstraction respects judgmental equality. Therefore we postulate the \define{congruence rule} for $\lambda$-abstraction,
\index{rules!for dependent function types!lambda-congruence@{$\lambda$-congruence}}
\index{lambda-congruence@{$\lambda$-congruence}}
\index{dependent function type!lambda-congruence@{$\lambda$-congruence}}
which asserts that\label{page:lambda-eq}
\begin{prooftree}
  \AxiomC{$\Gamma,x:A \vdash b(x)\jdeq b'(x) : B(x)$}
  \RightLabel{$\lambda$-eq.}
  \UnaryInfC{$\Gamma\vdash \lam{x}b(x)\jdeq \lam{x}b'(x) : \prd{x:A}B(x)$}
\end{prooftree}

\subsubsection{The $\Pi$-elimination rule}

\index{dependent function type!elimination rule|see {evaluation}}
The elimination rule for dependent function types provides us with a way to \emph{use} dependent functions. The way to use a dependent function is to evaluate it at an argument of the domain type. The $\Pi$-elimination rule is therefore also called the \define{evaluation rule}\index{evaluation|textbf}\index{rules!for dependent function types!evaluation|textbf}\index{dependent function type!evaluation|textbf}:
\begin{prooftree}
\AxiomC{$\Gamma\vdash f:\prd{x:A}B(x)$}
\RightLabel{$ev$.}
\UnaryInfC{$\Gamma,x:A\vdash f(x) : B(x)$}
\end{prooftree}
This rule asserts that given a dependent function $f:\prd{x:A}B(x)$ in context $\Gamma$ we obtain a term $f(x)$ of type $B(x)$ indexed by $x:A$ in context $\Gamma$. Again we require that evaluation respects judgmental equality:
\begin{prooftree}
  \AxiomC{$\Gamma\vdash f\jdeq f':\prd{x:A}B(x)$}
  \RightLabel{$ev$-eq.}
  \UnaryInfC{$\Gamma,x:A\vdash f(x)\jdeq f'(x):B(x)$}
\end{prooftree}

\subsubsection{The $\Pi$-computation rules}

\index{dependent function type!computation rules|see {$\beta$- and $\eta$-rules}}
We now postulate rules that specify the behavior of functions. First, we have a rule that asserts that a function of the form $\lam{x}b(x)$ behaves as expected: when we evaluate it at $x:A$, then we obtain the value $b(x):B(x)$. This rule is called the \define{$\beta$-rule}\index{b-rule for P-types@{$\beta$-rule for $\Pi$-types}|textbf}\index{rules!for dependent function types!b-rule@{$\beta$-rule}|textbf}\index{dependent function type!b-rule@{$\beta$-rule}|textbf}
\begin{prooftree}
\AxiomC{$\Gamma,x:A \vdash b(x) : B(x)$}
\RightLabel{$\beta$.}
\UnaryInfC{$\Gamma,x:A \vdash (\lambda y.b(y))(x)\jdeq b(x) : B(x)$}
\end{prooftree}
Second, we postulate a rule that asserts that all elements of a $\Pi$-type are (dependent) functions. This rule is known as the \define{$\eta$-rule}\index{eta-rule@{$\eta$-rule} for Pi-types@{for $\Pi$-types}|textbf}\index{rules!for dependent function types!eta-rule@{$\eta$-rule}|textbf}\index{dependent function type!eta-rule@{$\eta$-rule}|textbf}
\begin{prooftree}
\AxiomC{$\Gamma\vdash f:\prd{x:A}B(x)$}
\RightLabel{$\eta$.}
\UnaryInfC{$\Gamma \vdash \lam{x}f(x) \jdeq f : \prd{x:A}B(x)$}
\end{prooftree}
In other words, the computation rules ($\beta$ and $\eta$) for dependent function types postulate that $\lambda$-abstraction rule and the evaluation rule are mutual inverses. This completes the specification of dependent function types.

\subsection{Ordinary function types}

An important special case of $\Pi$-types arises when both $A$ and $B$ are types in context $\Gamma$. In this case, we can first weaken $B$ by $A$ and then apply the $\Pi$-formation rule to obtain the type $A\to B$ of \emph{ordinary} functions from $A$ to $B$, as in the following derivation:
\begin{prooftree}
\AxiomC{$\Gamma\vdash A~\textrm{type}$}
\AxiomC{$\Gamma\vdash B~\textrm{type}$}
\RightLabel{$W$}
\BinaryInfC{$\Gamma,x:A\vdash B~\textrm{type}$}
\RightLabel{$\Pi$}
\UnaryInfC{$\Gamma\vdash \prd{x:A}B~\textrm{type}$.}
\end{prooftree}
A term $f:\prd{x:A}B$ is a function that takes an argument $x:A$ and returns $f(x):B$. In other words, terms of type $\prd{x:A}B$ are indeed ordinary functions from $A$ to $B$. Therefore, we define the type $A\to B$\index{A arrow B@{$A\to B$}|see {function type}}\index{A arrow B@{$A\to B$}|textbf} of \define{(ordinary) functions}\index{function type|textbf} from $A$ to $B$ by
\begin{equation*}
  A\to B\defeq\prd{x:A}B.
\end{equation*}
If $f:A\to B$ is a function, then the type $A$ is also called the \define{domain} of $f$, and the type $B$ is also called the \define{codomain} of $f$.

Sometimes we will also write $B^A$\index{B^A@{$B^A$}|see {function type}} for the type $A\to B$.  Formally, we make such definitions by adding one more line to the above derivation:
\begin{prooftree}
\AxiomC{$\Gamma\vdash A~\textrm{type}$}
\AxiomC{$\Gamma\vdash B~\textrm{type}$}
\RightLabel{$W$}
\BinaryInfC{$\Gamma,x:A\vdash B~\textrm{type}$}
\RightLabel{$\Pi$}
\UnaryInfC{$\Gamma\vdash \prd{x:A}B~\textrm{type}$}
\UnaryInfC{$\Gamma\vdash A\to B \defeq \prd{x:A}B~\textrm{type}$.}
\end{prooftree}

\begin{rmk}
  More generally, we can make definitions at the end of a derivation if the conclusion is a certain type in context, or if the conclusion is a certain term of a type in context. Suppose, for instance, that we have a derivation
  \begin{prooftree}
    \AxiomC{$\mathcal{D}$}
    \UnaryInfC{$\Gamma\vdash a:A$,}
  \end{prooftree}
  in which the derivation $\mathcal{D}$ makes use of the premises $\mathcal{H}_1$, \ldots,$\mathcal{H}_n$. If we wish to make a definition $\newdef\defeq a$, then we can extend the derivation tree with
  \begin{prooftree}
    \AxiomC{$\mathcal{D}$}
    \UnaryInfC{$\Gamma\vdash a:A$}
    \UnaryInfC{$\Gamma\vdash\newdef\defeq a:A$.}
  \end{prooftree}
  The effect of such a definition is that we have extended our type theory with a new constant $\newdef$, for which the following inference rules are valid
  \begin{center}
    \begin{minipage}{.45\textwidth}
      \begin{prooftree}
        \AxiomC{$\mathcal{H}_1$\quad $\mathcal{H}_2$ \quad \dots \quad $\mathcal{H}_n$}
        \UnaryInfC{$\Gamma\vdash\newdef:A$}
      \end{prooftree}
    \end{minipage}
    \begin{minipage}{.45\textwidth}
      \begin{prooftree}
        \AxiomC{$\mathcal{H}_1$\quad $\mathcal{H}_2$ \quad \dots \quad $\mathcal{H}_n$}
        \UnaryInfC{$\Gamma\vdash\newdef\jdeq a:A$.}
      \end{prooftree}
    \end{minipage}
  \end{center}
  In our example of the definition of the ordinary function type $A\to B$, we therefore have by definition the following valid inference rules
  \begin{center}
    \begin{minipage}{.45\textwidth}
      \begin{prooftree}
        \AxiomC{$\Gamma\vdash A~\textrm{type}$}
        \AxiomC{$\Gamma\vdash B~\textrm{type}$}
        \BinaryInfC{$\Gamma\vdash A\to B~\textrm{type}$}
      \end{prooftree}
    \end{minipage}
    \begin{minipage}{.45\textwidth}
      \begin{prooftree}
        \AxiomC{$\Gamma\vdash A~\textrm{type}$}
        \AxiomC{$\Gamma\vdash B~\textrm{type}$}
        \BinaryInfC{$\Gamma\vdash A\to B\jdeq \prd{x:A}B~\textrm{type}$.}
      \end{prooftree}
    \end{minipage}
  \end{center}
  There are of course many such definitions throughout the development of dependent type theory, the univalent foundations of mathematics, and synthetic homotopy theory. They are all included in the index at the end of this book.
\end{rmk}

\begin{rmk}
  By the term conversion rules of \cref{ex:term_conversion} we can now use the rules for $\lambda$-abstraction, evaluation, and so on, to obtain corresponding rules for the ordinary function type $A\to B$. We give a brief summary of these rules, omitting the congruence rules.\index{rules!for function types}
  \begin{prooftree}
    \AxiomC{$\Gamma\vdash A~\textrm{type}$}
    \AxiomC{$\Gamma\vdash B~\textrm{type}$}
    \RightLabel{$\to$}
    \BinaryInfC{$\Gamma\vdash A\to B~\textrm{type}$}
  \end{prooftree}%
  \begin{center}
    \begin{minipage}{.55\textwidth}
      \begin{prooftree}
        \AxiomC{$\Gamma\vdash B~\textrm{type}$}
        \AxiomC{$\Gamma,x:A\vdash b(x):B$}
        \RightLabel{$\lambda$}
        \BinaryInfC{$\Gamma\vdash \lam{x}b(x):A\to B$}
      \end{prooftree}%
    \end{minipage}
    \begin{minipage}{.35\textwidth}
      \begin{prooftree}
        \AxiomC{$\Gamma\vdash f:A\to B$}
        \RightLabel{$ev$}
        \UnaryInfC{$\Gamma,x:A\vdash f(x):B$}
      \end{prooftree}%
    \end{minipage}
  \end{center}
  \begin{center}
    \begin{minipage}{.55\textwidth}
      \begin{prooftree}
        \AxiomC{$\Gamma\vdash B~\textrm{type}$}
        \AxiomC{$\Gamma,x:A\vdash b(x):B$}
        \RightLabel{$\beta$}
        \BinaryInfC{$\Gamma,x:A\vdash(\lam{y}b(y))(x)\jdeq b(x):B$}
      \end{prooftree}%
    \end{minipage}
    \begin{minipage}{.40\textwidth}
      \begin{prooftree}
        \AxiomC{$\Gamma\vdash f:A\to B$}
        \RightLabel{$\eta$}
        \UnaryInfC{$\Gamma\vdash\lam{x} f(x)\jdeq f:A\to B$}
      \end{prooftree}
    \end{minipage}
  \end{center}
\end{rmk}

Now we can use these rules to construct some familiar functions, such as the identity function $\idfunc:A\to A$ on an arbitrary type $A$, and the composition $g\circ f:A\to C$ of any two functions $f:A\to B$ and $g:B\to C$. 

\begin{defn}
For any type $A$ in context $\Gamma$, we define the \define{identity function}\index{identity function|textbf}\index{function!identity function|textbf} $\idfunc[A]:A\to A$\index{id A@{$\idfunc[A]$}|see {identity function}}\index{id A@{$\idfunc[A]$}|textbf} using the generic term:
\begin{prooftree}
\AxiomC{$\Gamma\vdash A~\textrm{type}$}
\UnaryInfC{$\Gamma,x:A\vdash x:A$}
\UnaryInfC{$\Gamma\vdash \lam{x}x:A\to A$}
\UnaryInfC{$\Gamma\vdash \idfunc[A]\defeq\lam{x}x:A\to A$.}
\end{prooftree}
\end{defn}

The identity function therefore satisfies the following inference rules:
  \begin{center}
    \begin{minipage}{.45\textwidth}
      \begin{prooftree}
        \AxiomC{$\Gamma\vdash A~\textrm{type}$}
        \UnaryInfC{$\Gamma\vdash \idfunc[A]:A\to A$}
      \end{prooftree}
    \end{minipage}
    \begin{minipage}{.45\textwidth}
      \begin{prooftree}
        \AxiomC{$\Gamma\vdash A~\textrm{type}$}
        \UnaryInfC{$\Gamma\vdash \idfunc[A]\jdeq\lam{x}x:A\to A$.}
      \end{prooftree}
    \end{minipage}
  \end{center}

Next, we define the composition of functions. We will introduce the composition operation itself as a function $\comp$ that takes two arguments: the first argument is a function $g:B\to C$, and the second argument is a function $f:A\to B$. The output is a function $\comp(g,f):A\to C$, for which we often write $g\circ f$.

\begin{rmk}
  Since composition is a function that takes multiple arguments, we need to know how to represent such functions. Types of functions with multiple arguments can be formed by iterating the $\Pi$-formation rule or the $\to$-formation rule. For example, a function
  \begin{equation*}
    f:A\to (B\to C)
  \end{equation*}
  takes two arguments: first it takes an argument $x:A$, and the output $f(x)$ has type $B\to C$. This is again a function type, so $f(x)$ is a function that takes an argument $y:B$, and its output $f(x)(y)$ has type $C$. We will usually write $f(x,y)$ for $f(x)(y)$.

  Similarly, when $C(x,y)$ is a family of types indexed by $x:A$ and $y:B(x)$, then we can form the dependent function type $\prd{x:A}\prd{y:B(x)}C(x,y)$. In the special case where $C(x,y)$ is a family of types indexed by two elements $x,y:A$ of the same type, then we often write
  \begin{equation*}
    \prd{x,y:A}C(x,y)
  \end{equation*}
  for the type $\prd{x:A}\prd{y:A}C(x,y)$.

  With the idea of iterating function types, we see that type of the composition operation $\comp$ is
  \begin{equation*}
    (B\to C)\to ((A\to B)\to (A\to C)).
  \end{equation*}
  It is the type of functions, taking a function $g:B\to C$, to the type of functions $(A\to B)\to (A\to C)$. Thus, $\comp(g)$ is again a function, mapping a function $f:A\to B$ to a function of type $A\to C$.
\end{rmk}

\begin{defn}
For any three types $A$, $B$, and $C$ in context $\Gamma$, there is a \define{composition}\index{function!composition|textbf}\index{composition!of functions|textbf} operation
\begin{equation*}
\comp:(B\to C)\to ((A\to B)\to (A\to C)).
\end{equation*}
We will usually write $g\circ f$\index{g o f@{$g\circ f$}|textbf} for $\comp(g,f)$\index{comp(g,f)@{$\comp(g,f)$}|textbf}\index{comp(g,f)@{$\comp(g,f)$}|see {composition, of functions}}.
\end{defn}

\begin{constr}
  The idea of the definition is to define $\comp(g,f)$ to be the function $\lam{x}g(f(x))$. The function $\comp$ is therefore defined as
  \begin{equation*}
    \comp\defeq \lam{g}\lam{f}\lam{x}g(f(x)).
  \end{equation*}
  The derivation we use to construct $\comp$ is as follows:
  \begin{small}
  \begin{prooftree}
    \AxiomC{$\Gamma\vdash A~\type$}
    \AxiomC{$\Gamma\vdash B~\type$}
    \RightLabel{(a)}
    \BinaryInfC{$\Gamma,f:B^A,x:A\vdash f(x):B$}
    \UnaryInfC{$\Gamma,g:C^B,f:B^A,x:A\vdash f(x):B$}
    \AxiomC{$\Gamma\vdash B~\type$}
    \AxiomC{$\Gamma\vdash C~\type$}
    \RightLabel{(b)}
    \BinaryInfC{$\Gamma,g:C^B,y:B\vdash g(y):C$}
    \UnaryInfC{$\Gamma,g:C^B,f:B^A,y:B\vdash g(y):C$}
    \UnaryInfC{$\Gamma,g:C^B,f:B^A,x:A,y:B\vdash g(y):C$}
    \BinaryInfC{$\Gamma,g:C^B,f:B^A,x:A\vdash g(f(x)) : C$}
    \UnaryInfC{$\Gamma,g:C^B,f:B^A\vdash \lam{x}g(f(x)):C^A$}
    \UnaryInfC{$\Gamma,g:B\to C\vdash \lam{f}\lam{x}g(f(x)):B^A\to C^A$}
    \UnaryInfC{$\Gamma\vdash\lam{g}\lam{f}\lam{x}g(f(x)):C^B\to (B^A\to C^A)$}
    \UnaryInfC{$\Gamma\vdash\comp\defeq \lam{g}\lam{f}\lam{x}g(f(x)):C^B\to (B^A\to C^A)$.}
  \end{prooftree}
  \end{small}
  Note, however, that we haven't derived the rules (a) and (b) yet. These rules assert that the \emph{generic functions} of $A\to B$ and $B\to C$ can also be evaluated. The formal derivation of this fact is as follows:
  \begin{prooftree}
    \AxiomC{$\Gamma\vdash A~\type$}
    \AxiomC{$\Gamma\vdash B~\type$}
    \BinaryInfC{$\Gamma\vdash A \to B~\type$}
    \UnaryInfC{$\Gamma,f:A\to B\vdash f:A\to B$}
    \UnaryInfC{$\Gamma,f:A\to B,x:A\vdash f(x):B$.}
  \end{prooftree}
  This completes the construction of $\comp$.
\end{constr}

In the remainder of this section we will see how to use the given rules for function types to derive the laws of a category\index{category laws!for functions} for functions. These are the laws that assert that function composition is associative and that the identity function satisfies the unit laws.

\begin{lem}
Composition of functions is associative\index{associativity!of function composition}, i.e., we can derive
\begin{prooftree}
\AxiomC{$\Gamma\vdash f:A\to B$}
\AxiomC{$\Gamma\vdash g:B\to C$}
\AxiomC{$\Gamma\vdash h:C\to D$}
\TrinaryInfC{$\Gamma \vdash (h\circ g)\circ f\jdeq h\circ(g\circ f):A\to D$.}
\end{prooftree}
\end{lem}

\begin{proof}
  The main idea of the proof is that both $((h\circ g)\circ f)(x)$ and $(h\circ (g\circ f))(x)$ evaluate to $h(g(f(x))$, and therefore $(h\circ g)\circ f$ and $h\circ(g\circ f)$ must be judgmentally equal. This idea is made formal in the following derivation:
  \begin{prooftree}
    \AxiomC{$\Gamma\vdash f:A\to B$}
    \UnaryInfC{$\Gamma,x:A\vdash f(x):B$}
    \AxiomC{$\Gamma\vdash g:B\to C$}
    \UnaryInfC{$\Gamma,y:B\vdash g(y):C$}
    \UnaryInfC{$\Gamma,x:A,y:B\vdash g(y):C$}
    \BinaryInfC{$\Gamma,x:A\vdash g(f(x)):C$}
    \AxiomC{$\Gamma\vdash h:C\to D$}
    \UnaryInfC{$\Gamma,z:C\vdash h(z):D$}
    \UnaryInfC{$\Gamma,x:A,z:C\vdash h(z):D$}
    \BinaryInfC{$\Gamma,x:A\vdash h(g(f(x))):D$}
    \UnaryInfC{$\Gamma,x:A\vdash h(g(f(x)))\jdeq h(g(f(x))):D$}
    \UnaryInfC{$\Gamma,x:A\vdash (h\circ g)(f(x))\jdeq h((g\circ f)(x)):D$}
    \UnaryInfC{$\Gamma,x:A\vdash ((h\circ g)\circ f)(x)\jdeq (h\circ (g \circ f))(x):D$}
    \UnaryInfC{$\Gamma\vdash (h\circ g)\circ f\jdeq h\circ(g\circ f):A\to D$.}
  \end{prooftree}
\end{proof}

\begin{lem}\label{lem:fun_unit}
Composition of functions satisfies the left and right unit laws\index{left unit law|see {unit laws}}\index{right unit law|see {unit laws}}\index{unit laws!for function composition}, i.e., we can derive
\begin{prooftree}
\AxiomC{$\Gamma\vdash f:A\to B$}
\UnaryInfC{$\Gamma\vdash \idfunc[B]\circ f\jdeq f:A\to B$}
\end{prooftree}
and
\begin{prooftree}
\AxiomC{$\Gamma\vdash f:A\to B$}
\UnaryInfC{$\Gamma\vdash f\circ\idfunc[A]\jdeq f:A\to B$.}
\end{prooftree}
\end{lem}

\begin{proof}
  Note that it suffices to derive that $\idfunc(f(x))\jdeq f(x)$ in context $\Gamma,x:A$, because once we derived this equality we can finish the derivation with
  \begin{prooftree}
    \AxiomC{$\vdots$}
    \UnaryInfC{$\Gamma,x:A\vdash\idfunc(f(x))\jdeq f(x):B$}
    \UnaryInfC{$\Gamma\vdash\lam{x}\idfunc(f(x))\jdeq\lam{x}f(x):A\to B$}
    \AxiomC{$\Gamma\vdash f:A\to B$}
    \UnaryInfC{$\Gamma\vdash\lam{x}f(x)\jdeq f:A\to B$}
    \BinaryInfC{$\Gamma\vdash\idfunc\circ f\jdeq f:A\to B$.}  
  \end{prooftree}
  The derivation of the equality $\idfunc(f(x))\jdeq f(x)$ in context $\Gamma,x:A$ is as follows:
  \begin{prooftree}
    \AxiomC{$\Gamma\vdash f:A\to B$}
    \UnaryInfC{$\Gamma,x:A\vdash f(x):B$}
    \AxiomC{$\Gamma\vdash A~\type$}
    \AxiomC{$\Gamma\vdash B~\type$}
    \UnaryInfC{$\Gamma,y:B\vdash\idfunc(y)\jdeq y:B$}
    \BinaryInfC{$\Gamma,x:A,y:B\vdash\idfunc(y)\jdeq y:B$}
    \BinaryInfC{$\Gamma,x:A\vdash\idfunc(f(x))\jdeq f(x):B$.}
  \end{prooftree}
  We leave the right unit law as \cref{ex:fun_right_unit}.
\end{proof}

\begin{exercises}
  \exitem \label{ex:eta_ext}The $\eta$-rule is often seen as a judgmental extensionality principle. Use the $\eta$-rule to show that if $f$ and $g$ take equal values, then they must be equal, i.e., give a derivation for the rule
  \begin{prooftree}
    \def\fCenter{\Gamma}
    \Axiom$\fCenter\vdash f:\prd{x:A}B(x)$
    \noLine
    \UnaryInf$\fCenter\vdash g:\prd{x:A}B(x)$
    \noLine
    \UnaryInf$\fCenter,x:A\vdash f(x)\jdeq g(x):B(x)$
    \UnaryInf$\fCenter\vdash f\jdeq g:\prd{x:A}B(x).$
  \end{prooftree}
  \exitem \label{ex:fun_right_unit}Give a derivation for the right unit law of \cref{lem:fun_unit}.\index{unit laws!for function composition}
  \exitem 
  \begin{subexenum}
  \item Construct the \define{constant map}\index{constant map|textbf}\index{function!constant map|textbf}\index{const x@{$\const_x$}|textbf}\index{function!const@{$\const$}|textbf}
    \begin{prooftree}
      \AxiomC{$\Gamma\vdash A~\textrm{type}$}
      \UnaryInfC{$\Gamma,y:B\vdash \const_y:A\to B$.}
    \end{prooftree}
  \item Show that
    \begin{prooftree}
      \AxiomC{$\Gamma\vdash f:A\to B$}
      \UnaryInfC{$\Gamma,z:C\vdash \const_z\circ f\jdeq\const_z : A\to C$.}
    \end{prooftree}
  \item Show that
    \begin{prooftree}
      \AxiomC{$\Gamma\vdash A~\textrm{type}$}
      \AxiomC{$\Gamma\vdash g:B\to C$}
      \BinaryInfC{$\Gamma,y:B\vdash g\circ\const_y\jdeq \const_{g(y)}:A\to C$.}
    \end{prooftree}
  \end{subexenum}
  \exitem \label{ex:swap}
  \begin{subexenum}
  \item Define the \define{swap function}\index{function!swap|textbf}\index{swap function|textbf}
    \begin{prooftree}
      \AxiomC{$\Gamma\vdash A~\mathrm{type}$}
      \AxiomC{$\Gamma\vdash B~\mathrm{type}$}
      \AxiomC{$\Gamma,x:A,y:B\vdash C(x,y)~\mathrm{type}$}
      \TrinaryInfC{$\Gamma\vdash \sigma:\Big(\prd{x:A}\prd{y:B}C(x,y)\Big)\to\Big(\prd{y:B}\prd{x:A}C(x,y)\Big)$}
    \end{prooftree}
    that swaps the order of the arguments.
  \item Show that
  \end{subexenum}
  \vspace{-.5\baselineskip}
  \begin{small}
    \begin{prooftree}
      \AxiomC{$\Gamma\vdash A~\mathrm{type}$}
      \AxiomC{$\Gamma\vdash B~\mathrm{type}$}
      \AxiomC{$\Gamma,x:A,y:B\vdash C(x,y)~\mathrm{type}$}
      \TrinaryInfC{$\Gamma\vdash \sigma\circ\sigma\jdeq\idfunc:\Big(\prd{x:A}\prd{y:B}C(x,y)\Big)\to \Big(\prd{x:A}\prd{y:B}C(x,y)\Big).$}
    \end{prooftree}
  \end{small}
\end{exercises}
\index{dependent function type|)}


\section{The natural numbers}
\label{sec:nat}

\index{inductive type|(}
\index{natural numbers|(}

The set of natural numbers is the most important object in mathematics. We quote Bishop\index{Bishop on the positive integers}, from his Constructivist Manifesto, the first chapter in Foundations of Constructive Analysis \cite{Bishop1967}, where he gives a colorful illustration of its importance to mathematics.

\begin{quote}
  ``The primary concern of mathematics is number, and this means the positive integers. We feel about number the way Kant felt about space. The positive integers and their arithmetic are presupposed by the very nature of our intelligence and, we are tempted to believe, by the very nature of intelligence in general. The development of the theory of the positive integers from the primitive concept of the unit, the concept of adjoining a unit, and the process of mathematical induction carries complete conviction. In the words of Kronecker, the positive integers were created by God. Kronecker would have expressed it even better if he had said that the positive integers were created by God for the benefit of man (and other finite beings). Mathematics belongs to man, not to God. We are not interested in properties of the positive integers that have no descriptive meaning for finite man. When a man proves a positive integer to exist, he should show how to find it. If God has mathematics of his own that needs to be done, let him do it himself.''
\end{quote}

A bit later in the same chapter, he continues:

\begin{quote}
  ``Building on the positive integers, weaving a web of ever more sets and ever more functions, we get the basic structures of mathematics: the rational number system, the real number system, the euclidean spaces, the complex number system, the algebraic number fields, Hilbert space, the classical groups, and so forth. Within the framework of these structures, most mathematics is done. Everything attaches itself to number, and every mathematical statement ultimately expresses the fact that if we perform certain computations within the set of positive integers, we shall get certain results.''
\end{quote}

\subsection{The formal specification of the type of natural numbers}
The type $\N$\index{N@{$\N$}|see {natural numbers}} of \define{natural numbers} is the archetypal example of an inductive type\index{inductive type!natural numbers}. The rules we postulate for the type of natural numbers come in four sets, just as the rules for $\Pi$-types:
\begin{enumerate}
\item The formation rule, which asserts that the type $\N$ can be formed.
\item The introduction rules, which provide the zero element $\zeroN$ and the successor function $\succN$.
\item The elimination rule. This rule is the type theoretic version of the induction principle for $\N$.
\item The computation rules, which assert that any application of the elimination rule behaves as expected on the constructors $\zeroN$ and $\succN$ of $\N$.
\end{enumerate}

\subsubsection{The formation rule of $\N$}

\index{natural numbers!rules for N@{rules for $\N$}!formation}
\index{rules!for N@{for $\N$}!formation}
The type $\N$ is formed by the \define{$\N$-formation} rule
\begin{prooftree}
  \AxiomC{}
  \RightLabel{$\N$-form.}
  \UnaryInfC{$\vdash \N~\type$}
\end{prooftree}
In other words, $\N$ is postulated to be a type in the empty context.

\subsubsection{The introduction rules of $\N$}
Unlike the set of positive integers in Bishop's remarks, Peano's first axiom postulates that $0$ is a natural number. The introduction rules for $\N$ equip it with the \define{zero element} and the \define{successor function}.
\index{natural numbers!rules for N@{rules for $\N$}!introduction rules}
\index{rules!for N@{for $\N$}!introduction rules}
\index{natural numbers!operations on N@{operations on $\N$}!0 N@{$\zeroN$}}
\index{0 N@{$\zeroN$}}
\index{successor function!on N@{on $\N$}}
\index{natural numbers!operations on N@{operations on $\N$}!succ N@{$\succN$}}
\index{succ N@{$\succN$}}

\bigskip
\begin{minipage}{.45\textwidth}
  \begin{prooftree}
    \AxiomC{}
    \UnaryInfC{$\vdash \zeroN:\N$}
  \end{prooftree}
\end{minipage}
\begin{minipage}{.45\textwidth}
  \begin{prooftree}
    \AxiomC{}
    \UnaryInfC{$\vdash \succN:\N\to\N$}
  \end{prooftree}
\end{minipage}

\bigskip
\begin{rmk}
  Every element in type theory always comes equipped with its type. Therefore it is possible in type theory that all elements have a \emph{unique} type. In general, it is therefore good practice to make sure that every element is given a unique name, and in formalized mathematics in computer proof assistants this is even required. For example, the element $\zeroN$ has type $\N$, and it is not also a type of $\Z$. This is why we annotate the terms $\zeroN$ and $\succN$ with their type in the subscript. The type $\Z$ of the integers will be introduced in the next section, which will come equipped with a zero element $\zeroZ$ and a successor function $\succZ$.
\end{rmk}

\subsubsection{The induction principle of $\N$}

\index{natural numbers!rules for N@{rules for $\N$}!elimination|see {induction}}
\index{natural numbers!rules for N@{rules for $\N$}!induction principle|(}
\index{induction principle!of N@{of $\N$}|(}
The classical induction principle of the natural numbers tells us what we have to do in order to show that $\forall_{(n\in\N)}P(n)$ holds, for a predicate $P$ over $\N$. Recall that a predicate $P$ on a set $X$ is just a proposition $P(x)$ about an arbitrary $x\in X$. For example, the assertion that `$n$ is divisible by five' is a predicate on the natural numbers.

In dependent type theory we may think of a type family $P$ over $\N$ as a predicate over $\N$. The type theoretical induction principle of $\N$ is therefore formulated using a type family $P$ over $\N$:\index{ind N@{$\indN$}|textbf}\index{rules!for N@{for $\N$}!induction principle|textbf}\index{natural numbers!indN@{$\indN$}|textbf}
\begin{prooftree}
  \def\fCenter{\Gamma}
  \Axiom$\fCenter, n:\N\vdash P(n)~\type$
  \noLine
  \UnaryInf$\fCenter\ \vdash p_0:P(\zeroN)$
  \noLine
  \UnaryInf$\fCenter\ \vdash p_S:\prd{n:\N}P(n)\to P(\succN(n))$
  \RightLabel{$\N$-ind.}
  \UnaryInf$\fCenter\ \vdash \indN(p_0,p_S):\prd{n:\N} P(n)$
\end{prooftree}
In other words, the type theoretical induction principle of $\N$ tells us what we need to do in order to construct a dependent function $\prd{n:\N}P(n)$. Just as in the classical induction principle, there are two things to be constructed given a type family $P$ over $\N$: in the \define{base case}\index{base case} we need to construct an element $p_0:P(\zeroN)$, and for the \define{inductive step}\index{inductive step} we need to construct a function of type $P(n)\to P(\succN(n))$ for all $n:\N$.

\begin{rmk}
  We might alternatively present the induction principle of $\N$ as the following inference rule
  \begin{prooftree}
    \AxiomC{$\Gamma,n:\N\vdash P(n)~\type$}
    \UnaryInfC{$\Gamma\vdash \indN : P(\zeroN)\to \Big(\Big(\prd{n:\N}P(n)\to P(\succN(n))\Big)\to \prd{n:\N}P(n)\Big)$.}
  \end{prooftree}
  In other words, for any type family $P$ over $\N$ there is a \emph{function} $\indN$ that takes two arguments, one for the base case and one for the inductive step, and returns a section of $P$. We claim that this rule is \emph{interderivable} with the rule $\N$-ind above.
  
  To see that indeed we get such a function from the rule $\N$-ind, we use generic elements. First, we let $\Gamma'$ be the context
  \begin{equation*}
    \Gamma,~p_0:P(\zeroN),~p_S:\prd{n:\N}P(n)\to P(\succN(n)).
  \end{equation*}
  By weakening we obtain that
  \begin{align*}
    & \Gamma',~n:\N\vdash P(n)~\type \\
    & \Gamma'\vdash p_0 : P(\zeroN) \\
    & \Gamma'\vdash p_S : \prd{n:\N}P(n)\to P(\succN(n)).
  \end{align*}
  Therefore, the induction principle of $\N$ provides us with a dependent function
  \begin{equation*}
    \Gamma' \vdash \indN(p_0,p_S) : \prd{n:\N}P(n).
  \end{equation*}
  Now we proceed by $\lambda$-abstraction twice to obtain a function
  \begin{equation*}
    \indN : P(\zeroN)\to \Big(\Big(\prd{n:\N}P(n)\to P(\succN(n))\Big) \to \prd{n:\N}P(n)\Big)
  \end{equation*}
  in the original context $\Gamma$. This shows that we can define the function $\indN$ from the rule $\N$-ind. Conversely, we can derive the rule $\N$-ind from the rule that presents $\indN$ as a function. We conclude that the ``official'' rule $\N$-ind and the rule that presents $\indN$ as a function are indeed interderivable.
\end{rmk}
\index{natural numbers!rules for N@{rules for $\N$}!induction|)}
\index{induction principle!of N@{of $\N$}|)}

\subsubsection{The computation rules of $\N$}

\index{computation rules!for N@{for $\N$}|(}
\index{natural numbers!rules for N@{rules for $\N$}!computation rules|(}
The computation rules for $\N$ postulate that the dependent function
\begin{equation*}
  \indN(p_0,p_S):\prd{n:\N}P(n)
\end{equation*}
behaves as expected when it is applied to $\zeroN$ or a successor. There is one computation rule for each step in the induction principle, covering the base case and the inductive step.

The computation rule for the base case is\index{rules!for N@{for $\N$}!computation rules|(}\index{computation rules!for N@{for $\N$}|textbf}
\begin{prooftree}
    \def\fCenter{\Gamma}
  \Axiom$\fCenter, n:\N\vdash P(n)~\type$
  \noLine
  \UnaryInf$\fCenter\ \vdash p_0:P(\zeroN)$
  \noLine
  \UnaryInf$\fCenter\ \vdash p_S:\prd{n:\N}P(n)\to P(\succN(n))$
  \UnaryInf$\fCenter\ \vdash \indN(p_0,p_S,\zeroN)\jdeq p_0 : P(\zeroN).$
\end{prooftree}
\begin{samepage}
  The computation rule for the inductive step has the same premises as the computation rule for the base case:
  \begin{prooftree}
    \AxiomC{$\cdots$}
    \UnaryInfC{$\Gamma, n:\N \vdash  \indN(p_0,p_S,\succN(n))\jdeq p_S(n,\indN(p_0,p_S,n)) : P(\succN(n))$.}
  \end{prooftree}
\end{samepage}
This completes the formal specification of the type $\N$ of natural numbers.
\index{rules!for N@{for $\N$}!computation rules|)}
\index{computation rules!for N@{for $\N$}|)}
\index{natural numbers!rules for N@{rules for $\N$}!computation rules|)}

\subsection{Addition on the natural numbers}

\index{addition on N@{addition on $\N$}|(}
\index{natural numbers!operations on N@{operations on $\N$}!addition|(}
The type theoretic induction principle of $\N$ can be used to do all the usual constructions of operations on $\N$, and to derive all the familiar properties about natural numbers. Many of those properties, however, require a few more ingredients of Martin-L\"of's dependent type theory. For example, the traditional inductive proof that the triangular numbers can be calculated by
\begin{equation*}
  1+\cdots+n = \frac{n(n+1)}{2}
\end{equation*}
is analogous in type theory, but it requires the identity type to state this equation. We will introduce the identity type in \cref{sec:identity}. Until we have fully specified all the ways of forming types in Martin-L\"of's dependent type theory, we are a bit limited in what we can do with the natural numbers, but at the present stage we can define some of the familiar operations on $\N$. We give in this section the type theoretical construction the \define{addition operation} by induction on $\N$, along with the complete derivation tree.

\begin{defn}\label{defn:addN}
  We define a function\index{add N@{$\addN$}|textbf}\index{natural numbers!operations on N@{operations on $\N$}!add N@{$\addN$}|textbf}\index{addition on N@{addition on $\N$}|textbf}\index{natural numbers!addition|textbf}
  \begin{equation*}
    \addN:\N\to (\N\to\N)
  \end{equation*}
  satisfying the specification
  \begin{align*}
    \addN(m,\zeroN) & \jdeq m \\
    \addN(m,\succN(n)) & \jdeq \succN(\addN(m,n)).
  \end{align*}
  Usually we will write $m+n$ for $\addN(m,n)$.
\end{defn}

\begin{proof}[Construction.]
  We will construct the binary operation $\addN:\N\to(\N\to\N)$ by induction on the second variable. In other words, we will construct an element
  \begin{equation*}
    m:\N \vdash \addN(m):\N\to\N.
  \end{equation*}
  The context $\Gamma$ we work in is therefore $m:\N$. The induction principle of $\N$ is used with the family of types $P(n)\defeq \N$ indexed by $n:\N$ in context $m:\N$. Therefore we need to construct
  \begin{align*}
    m:\N & \vdash \addzeroN(m) : \N\\
    m:\N & \vdash \addsuccN(m) : \N\to(\N\to\N),\\
    \intertext{in order to obtain}
    m:\N & \vdash \addN(m)\defeq\indN(\addzeroN(m),\addsuccN(m)):\N\to\N.
  \end{align*}
  The element $\addzeroN(m):\N$ in context $m:\N$ is of course defined to be $m:\N$, i.e., by the generic element, because adding zero should just be the identity function.
  To see how the function $\addsuccN(m):\N\to(\N\to\N)$ should be defined, we look at the specification of $\addN(m)$ when it is applied to a successor:
  \begin{equation*}
    \addN(m,\succN(n))\jdeq \succN(\addN(m,n)).
  \end{equation*}
  This shows us that we should define
  \begin{equation*}
    \addsuccN(m,n,x)\jdeq \succN(x),
  \end{equation*}
  because with this definition we will have
  \begin{align*}
    \addN(m,\succN(n)) & \jdeq \indN(\addzeroN(m),\addsuccN(m),\succN(n)) \\
                       & \jdeq \addsuccN(m,n,\addN(m,n)) \\
                       & \jdeq \succN(\addN(m,n)).
  \end{align*}
  The formal derivation for the construction of $\addsuccN$ is as follows:
  \begin{prooftree}
    \AxiomC{}
    \UnaryInfC{$\vdash\N~\type$}
    \AxiomC{}
    \UnaryInfC{$\vdash\N~\type$}
    \AxiomC{}
    \UnaryInfC{$\vdash \succN:\N\to\N$}
    \BinaryInfC{$n:\N\vdash \succN:\N\to\N$}
    \BinaryInfC{$m:\N,n:\N \vdash \succN:\N\to\N$}
    \UnaryInfC{$m:\N \vdash \lam{n}\succN:\N\to (\N \to \N)$}
    \UnaryInfC{$m:\N \vdash \addsuccN(m) \defeq \lam{n}\succN:\N\to (\N \to \N)$.}
  \end{prooftree}
  We combine this derivation with the induction principle of $\N$ to complete the construction of addition:
  \begin{prooftree}
    \AxiomC{$\vdots$}
    \UnaryInfC{$m:\N\vdash \addzeroN(m) \defeq m:\N$}
    \AxiomC{$\vdots$}
    \UnaryInfC{$m:\N\vdash \addsuccN(m):\N\to (\N \to \N)$}
    \BinaryInfC{$m:\N\vdash\indN(\addzeroN(m),\addsuccN(m)):\N\to \N$}
    \UnaryInfC{$m:\N\vdash\addN(m)\defeq\indN(\addzeroN(m),\addsuccN(m)):\N\to \N$.}
  \end{prooftree}
  The asserted judgmental equalities then hold by the computation rules for $\N$.
\end{proof}

\begin{rmk}
  By the computation rules for $\N$ it follows that
  \begin{equation*}
    m+\zeroN\jdeq m,\qquad\text{and}\qquad m+\succN(n)\jdeq\succN(m+n).
  \end{equation*}
  A simple consequence of this definition is that $\succN(n)\jdeq n+1$, as one would expect. However, the rules that we provided so far are not sufficient to also conclude that $\zeroN+n\jdeq n$ and $\succN(m) + n\jdeq \succN(m+n)$. In fact, dependent type theory with its inductive types does not provide any means to prove such judgmental equalities.

  Nevertheless, once we have introduced the \emph{identity type} in \cref{sec:identity} we will be able to \emph{identify} $\zeroN+n$ with $n$, and $\succN(m)+n$ with $\succN(m+n)$. See \cref{prp:unit-laws-add-N,prp:successor-laws-add-N}. 
\end{rmk}
\index{addition on N@{addition on $\N$}|)}
\index{natural numbers!operations on N@{operations on $\N$}!addition|)}

\subsection{Pattern matching}

Note that in definition \cref{defn:addN} we stated that $\addN$ is a function of type $\N\to (\N\to\N)$ satisfying the specification
\begin{align*}
    \addN(m,\zeroN) & \jdeq m \\
  \addN(m,\succN(n)) & \jdeq \succN(\addN(m,n)).
\end{align*}
Such a specification is enough to characterize the function $\addN(m)$ entirely, because it postulates the behaviour of $\addN(m)$ at the constructors of $\N$.
It is therefore convenient to present the definition of $\addN$ recursively in the following way:
\begin{align*}
  \addN(m,\zeroN) & \defeq m \\
  \addN(m,\succN(n)) & \defeq \succN(\addN(m,n)).
\end{align*}

More generally, if we want to define a dependent function $f:\prd{n:\N}P(n)$ by induction on $n$, using
\begin{align*}
  p_0 & : P(\zeroN) \\
  p_S & : \prd{n:\N} P(n)\to P(\succN(n)),
\end{align*}
we can present that definition by writing
\begin{align*}
  f(\zeroN) & \defeq p_0 \\
  f(\succN(n)) & \defeq p_S(n,f(n)). 
\end{align*}
When the definition of $f$ is presented in this way, we say that $f$ is defined by \define{pattern matching}\index{pattern matching} on the variable $n$. To see that $f$ is fully specified when it is defined by pattern matching, we have to recover the dependent function
\begin{equation*}
  p_S:\prd{n:\N}P(n)\to P(\succN(n))
\end{equation*}
from the expression $p_S(n,f(n))$ that was used in the definition of $f$. This can of course be done by replacing all occurrences of the term $f(n)$ in the expression $p_S(n,f(n))$ with a fresh variable $x:P(n)$. In other words, when a subexpression of $p_S(n,f(n))$ \emph{matches} $f(n)$, we replace that subexpression by $x$. This is where the name \emph{pattern matching} comes from. Many computer proof assistants have the pattern matching mechanism built in, because it is a concise way of presenting a recursive definition. Another advantage of presenting definitions by pattern matching is that the judgmental equalities by which the object is defined are immediately displayed. Those judgmental equalities are all that is known about the defined object, and often proving things about it amounts to finding a way to apply those judgmental equalities.

Pattern matching can also be used in more complicated situations, such as defining a function by pattern matching on multiple variables, or by iterated pattern matching. For example, an alternative definition of addition on $\N$ could be given by pattern matching on both variables:
\begin{align*}
  \addpN(\zeroN,\zeroN) & \defeq \zeroN \\*
  \addpN(\zeroN,\succN(n)) & \defeq \succN(n) \\*
  \addpN(\succN(m),\zeroN) & \defeq \succN(m) \\*
  \addpN(\succN(m),\succN(n)) & \defeq \succN(\succN(\addpN(m,n)).
\end{align*}

An example of a definition by iterated pattern matching is the \define{Fibonacci function} $F:\N\to\N$. This function is defined by
\begin{align*}
  F(\zeroN) & \defeq \zeroN \\
  F(\oneN) & \defeq \oneN \\
  F(\succN(\succN(n))) & \defeq F(\succN(n))+F(n).
\end{align*}
However, since $F(\succN(\succN(n)))$ is defined using both $F(\succN(n))$ and $F(n)$, it is not immediately clear how to present $F$ by the usual induction principle of $\N$. It is a nice puzzle, which we leave as \cref{ex:fibonacci}, to find a definition of the Fibonacci sequence with the usual induction principle of $\N$. 

\begin{exercises}
  \exitem
  \begin{subexenum}
  \item Define the \define{multiplication} operation
    \index{multiplication!on N@{on $\N$}|textbf}
    \index{natural numbers!operations on N@{operations on $\N$}!mul N@{$\mulN$}|textbf}
    \index{mul N@{$\mulN$}|textbf}
    \begin{equation*}
      \mulN :\N\to(\N\to\N).
    \end{equation*}
  \item Define the \define{exponentiation function} $n,m\mapsto m^n$ of type $\N\to (\N\to \N)$.
  \index{exponentiation function on N@{exponentiation function on $\N$}|textbf}
  \index{natural numbers!operations on N@{operations on $\N$}!exponentiation|textbf}
  \end{subexenum}
  \exitem Define the binary \define{min} and \define{max} functions
  \index{minimum function|textbf}
  \index{maximum function|textbf}
  \index{natural numbers!operations on N@{operations on $\N$}!minN@{$\minN$}|textbf}
  \index{natural numbers!operations on N@{operations on $\N$}!maxN@{$\maxN$}|textbf}
  \begin{equation*}
    \minN,\maxN:\N\to(\N\to\N).
  \end{equation*}
  \exitem
  \begin{subexenum}
  \item Define the \define{triangular numbers}
    \begin{equation*}
      1+\cdots+n.
    \end{equation*}
    \index{triangle number|textbf}
    \index{natural numbers!operations on N@{operations on $\N$}!triangle number|textbf}
  \item Define the \define{factorial} operation $n\mapsto n!$.
  \index{factorial operation|textbf}
  \index{natural numbers!operations on N@{operations on $\N$}!n factorial@{$n"!$}|textbf}
  \end{subexenum}
  \exitem Define the \define{binomial coefficient} $\binom{n}{k}$\index{(n k)@{$\binom{n}{k}$}|see {binomial coefficient}}\index{(n k)@{$\binom{n}{k}$}|textbf}\index{binomial coefficient|textbf}\index{natural numbers@operations on N@{operations on $\N$}!binomial coefficient|textbf} for any $n,k:\N$, making sure that $\binom{n}{k}\jdeq 0$ when $n<k$.
  \index{binomial coefficient|textbf}
  \index{natural numbers!operations on N@{operations on $\N$}!binomial coefficient|textbf}
  \exitem \label{ex:fibonacci}Use the induction principle of $\N$ to define the \define{Fibonacci sequence} as a function $F:\N\to\N$ that satisfies the equations\index{Fibonacci sequence|textbf}\index{natural numbers!operations on N@{operations on $\N$}!Fibonacci sequence|textbf}
  \begin{samepage}
    \begin{align*}
      F(\zeroN) & \jdeq \zeroN \\
      F(\oneN) & \jdeq \oneN \\
      F(\succN(\succN(n))) & \jdeq F(\succN(n))+F(n).
    \end{align*}
  \end{samepage}
  \exitem Define division by two rounded down as a function $\N\to\N$ in two ways: first by pattern matching, and then directly by the induction principle of $\N$.
\end{exercises}
\index{natural numbers|)}


\section{More inductive types}\label{sec:inductive}

In the previous section we introduced the type of natural numbers. Many other types can also be introduced as inductive types. In this section we will see by example how that works. We will introduce the unit type, the empty type, coproducts, dependent pair types, and cartesian products as inductive types, and in the next section the identity type will be introduced as an inductive family of types.

From this section on, we will also start using a more informal style. The inductive types will be specified by a description of their constructors and induction principles in terms of operations on dependent function types, which is more tightly connected with how we will use them, but we will not display the formal rules. It is a good exercise for the reader to formally specify at least some of the inductive types of this section by stating their formal rules.

\subsection{The idea of general inductive types}

Just like the type of natural numbers, other inductive types are also specified by their \emph{constructors}, an \emph{induction principle}, and their \emph{computation rules}: 
\begin{enumerate}
\item The constructors tell what structure the inductive type comes equipped with. There may be any finite number of constructors, even no constructors at all, in the specification of an inductive type. 
\item The induction principle specifies the data that should be provided in order to construct a section of an arbitrary type family over the inductive type. The idea of the induction principle is always the same: in order to define a dependent function $f:\prd{x:A}B(x)$, one has to specify the behaviour of $f$ at the constructors of $A$.
\item The computation rules assert that the inductively defined section agrees on the constructors with the data that was used to define the section. Thus, there is a computation rule for every constructor.
\end{enumerate}
Since any inductively defined function is entirely determined by its behavior on the constructors, we can again present such inductive definitions by pattern matching. Therefore, we will also specify for each inductive type how to give definitions by pattern matching.

\subsection{The unit type}
\index{unit type|(}
\index{inductive type!unit type|(}
A straightforward example of an inductive type is the \emph{unit type}, which has just one constructor. 
Its induction principle is analogous to just the base case of induction on the natural numbers.

\begin{defn}
We define the \define{unit type}\index{1 @{$\unit$}|see {unit type}}\index{1 @{$\unit$}|textbf}\index{unit type|textbf} to be a type $\unit$ equipped with a term\index{unit type!star@{$\ttt$}|textbf}
\begin{equation*}
\ttt:\unit,
\end{equation*}
satisfying the induction principle\index{induction principle!of the unit type|textbf}\index{unit type!induction principle|textbf} that for any family of types $P(x)$ indexed by $x:\unit$, there is a function\index{ind 1@{$\indunit$}|textbf}\index{unit type!indunit@{$\indunit$}|textbf}
\begin{equation*}
\indunit : P(\ttt)\to\prd{x:\unit}P(x)
\end{equation*}
for which the computation rule\index{computation rules!for the unit type|textbf}\index{unit type!computation rules|textbf}
\begin{equation*}
\indunit(p,\ttt) \jdeq p
\end{equation*}
holds. Alternatively, a definition of a dependent function $f:\prd{x:\unit}P(x)$ by induction using $p:P(\ttt)$ can be presented by pattern matching as
\begin{equation*}
  f(\ttt)\defeq p.
\end{equation*}
\end{defn}

A special case of the induction principle arises when $P$ does not actually depend on $\unit$. If we are given a type $A$, then we can first weaken it to obtain the constant family over $\unit$, with value $A$. Then the induction principle of the unit type provides a function
\begin{equation*}
  \indunit : A \to (\unit\to A).
\end{equation*}
In other words, by the induction principle for the unit type we obtain for every $x:A$ a function $\pt_x\defeq\indunit(x):\unit\to A$.\index{pt x@{$\pt_x$}|textbf}
\index{unit type|)}
\index{inductive type!unit type|)}

\subsection{The empty type}
\index{empty type|(}
\index{inductive type!empty type|(}
The empty type is a degenerate example of an inductive type. It does \emph{not} come equipped with any constructors, and therefore there are also no computation rules. The induction principle merely asserts that any type family has a section. In other words: if we assume the empty type has a term, then we can prove anything.

\begin{defn}
We define the \define{empty type}\index{0 @{$\emptyt$}|see {empty type}}\index{0 @{$\emptyt$}|textbf}\index{empty type|textbf} to be a type $\emptyt$ satisfying the induction principle\index{induction principle!of the empty type|textbf}\index{empty type!induction principle|textbf} that for any family of types $P(x)$ indexed by $x:\emptyt$, there is a term\index{ind 0@{$\indempty$}|textbf}\index{empty type!indempty@{$\indempty$}|textbf}
\begin{equation*}
\indempty : \prd{x:\emptyt}P(x).
\end{equation*}
\end{defn}

It is again a special case of the induction principle that we have a function
\begin{equation*}
  \exfalso\defeq\indempty:\emptyt\to A
\end{equation*}
for any type $A$. Indeed, to obtain this function one first weakens $A$ to obtain the constant family over $\emptyt$ with value $A$, and then the induction principle gives the desired function. The function $\exfalso$ can be used to draw any conclusion after deriving a contradiction, because \emph{ex falso quodlibet}.

We can also use the empty type to define the negation operation on types.

\begin{defn}
  For any type $A$ we define \define{negation}\index{negation!of types|textbf}\index{n A@{$\neg A$}|see {negation}|textbf} of $A$ by
  \begin{align*}
    \neg A & \defeq A\to\emptyt.
  \intertext{We also say that a type $A$ \define{is empty} if it comes equipped with an element of type $\neg A$. Therefore, we also define}
    \isempty(A) & \defeq A\to\emptyt.
  \end{align*}
\end{defn}

\begin{rmk}
  Since $\neg A$ is the type of functions from $A$ to $\emptyt$, a proof of $\neg A$ is given by assuming that $A$ holds, and then constructing an element of the empty type. In other words, we prove $\neg A$ by assuming $A$ and deriving a contradiction. This proof technique is called \define{proof of negation}\index{proof of negation}.

  Proofs of negation should not be confused with proofs by contradiction\index{proof by contradiction}. Even though a proof of negation involves deriving a contradiction, in logic a \define{proof by contradiction} of a proposition $P$ is an argument where we conclude that $P$ holds after showing that $\neg P$ implies a contradiction. In other words, a proof by contradiction uses the logical step $\neg\neg P \Rightarrow P$, which is also called \define{double negation elimination}.

  In type theory, however, note that the type $\neg\neg A$ is the type of functions
  \begin{equation*}
    (A\to\emptyt)\to\emptyt.
  \end{equation*}
  This type is quite different from the type $A$ itself, and with the given rules of type theory it is not possible to construct a function $\neg\neg A \to A$ unless more is known about the type $A$. In other words, before one can prove by contradiction that there is an element in $A$, one has to construct a function $\neg\neg A\to A$, and it depends on the specific type $A$ whether this is possible at all. In \cref{ex:dne-is-decidable} we will see a situation where we can indeed construct a function $\neg\neg A\to A$. In practice, however, we will rarely use double negation elimination.
\end{rmk}

In the following proposition we illustrate how to work with the type theoretic definition of negation.

\begin{prp}\label{prp:contravariant-neg}
  For any two types $P$ and $Q$, there is a function
  \begin{equation*}
    (P\to Q)\to (\neg Q \to \neg P).
  \end{equation*}
\end{prp}

\begin{proof}
  The desired function is defined by $\lambda$-abstraction, so we begin by assuming that we have a function $f:P\to Q$. Then we have to construct a function $\neg Q\to\neg P$, which is again constructed by $\lambda$-abstraction. We assume that we have $\tilde{q}:\neg Q$. By our definition of $\neg Q$ we see that $\tilde{q}$ is a function $Q\to\emptyt$. Now we have to construct a term of type $\neg P$, which is the type of functions $P\to\emptyt$. We apply $\lambda$-abstraction once more, so we assume $p:P$. Now we have
  \begin{align*}
    f & : P \to Q \\*
    \tilde{q} & : Q\to \emptyt \\*
    p & : P,
  \end{align*}
  and our goal is to construct a term of the empty type.

  Since we have $f:P\to Q$ and $p:P$, we obtain $f(p):Q$. Moreover, we have $\tilde{q}:Q\to\emptyt$, so we obtain $\tilde{q}(f(p)):\emptyt$. This completes the proof. The function we have constructed is
  \begin{equation*}
    \lam{f}\lam{\tilde{q}}\lam{p}\tilde{q}(f(p)):(P\to Q)\to(\neg Q\to\neg P).\qedhere
  \end{equation*}
\end{proof}

We leave it to the reader to construct the corresponding natural deduction tree, that formally constructs a function
\begin{equation*}
  (P\to Q)\to(\neg Q\to \neg P).
\end{equation*}
\index{empty type|)}
\index{inductive type!empty type|)}

\subsection{Coproducts}\label{sec:coprod}
\index{coproduct|(}
\index{inductive type!coproduct|(}
\begin{defn}
Let $A$ and $B$ be types. We define the \define{coproduct}\index{disjoint sum|see {coproduct}}\index{coproduct|textbf} $A+B$\index{A + B@{$A+B$}|see {coproduct}} to be a type that comes equipped with\index{inl@{$\inl$}|textbf}\index{coproduct!inl@{$\inl$}|textbf}\index{inr@{$\inr$}|textbf}\index{coproduct!inr@{$\inr$}|textbf}
\begin{align*}
\inl & : A \to A+B \\*
\inr & : B \to A+B,
\end{align*}
satisfying the induction principle\index{induction principle!of coproducts|textbf}\index{coproduct!induction principle|textbf} that for any family of types $P(x)$ indexed by $x:A+B$, there is a term\index{ind +@{$\indcoprod$}|textbf}\index{coproduct!ind+@{$\indcoprod$}|textbf}
\begin{equation*}
\indcoprod : \Big(\prd{x:A}P(\inl(x))\Big)\to\Big(\Big(\prd{y:B}P(\inr(y))\Big)\to\prd{z:A+B}P(z)\Big)
\end{equation*}
for which the computation rules\index{computation rules!for coproducts|textbf}\index{coproduct!computation rules|textbf}
\begin{align*}
\indcoprod(f,g,\inl(x)) & \jdeq f(x) \\*
\indcoprod(f,g,\inr(y)) & \jdeq g(y)
\end{align*}
hold. Alternatively, a definition of a dependent function $h:\prd{x:A+B}P(x)$ by induction using $f:\prd{x:A}P(\inl(x))$ and $g:\prd{y:B}P(\inr(y))$ can be presented by pattern matching as
\begin{align*}
  h(\inl(x)) & \defeq f(x) \\*
  h(\inr(y)) & \defeq g(y).
\end{align*}
Sometimes we write $[f,g]$ for the function $\indcoprod(f,g)$. The coproduct of two types is sometimes also called the \define{disjoint sum}.
\end{defn}

By the induction principle of coproducts we obtain a function
\begin{equation*}
  \indcoprod:(A\to X) \to \big((B\to X) \to (A+B\to X)\big)
\end{equation*}
for any type $X$. Note that this special case of the induction principle of coproducts is very similar to the elimination rule of disjunction in first order logic: if $P$, $P'$, and $Q$ are propositions, then we have
\begin{equation*}
  (P\to Q)\to \big((P'\to Q)\to (P\lor P'\to Q)\big).
\end{equation*}
Indeed, we can think of \emph{propositions as types} and of terms as their constructive proofs. Under this interpretation of type theory the coproduct is indeed the disjunction.

\begin{rmk}\label{rmk:functor-coprod}
  A simple application of the induction principle for coproducts gives us a map\index{coproduct!functorial action|textbf}\index{functorial action!of coproducts|textbf}\index{f + g@{$f+g$}|see {functorial action, of coproducts}}\index{f + g@{$f+g$}|textbf}
  \begin{equation*}
    f+g:A+B\to A'+B'
  \end{equation*}
  for every $f:A\to A'$ and $g:B\to B'$. Indeed, the map $f+g$ is defined by
  \begin{align*}
    (f+g)(\inl(x)) & \defeq \inl(f(x)) \\*
    (f+g)(\inr(y)) & \defeq \inr(g(y)).
  \end{align*}
\end{rmk}

\begin{prp}
  Consider two types $A$ and $B$, and suppose that $B$ is empty. Then there is a function
  \begin{equation*}
    (A+B)\to A.
  \end{equation*}
\end{prp}

\begin{rmk}
  In other words, there is a function
  \begin{equation*}
    \isempty(B) \to ((A+B)\to A),
  \end{equation*}
  for any two types $A$ and $B$. Similarly, there is a function
  \begin{equation*}
    \isempty(A)\to ((A+B)\to B),
  \end{equation*}
  for any two types $A$ and $B$.
\end{rmk}

\begin{proof}
  We will construct the function $(A+B)\to A$ with the induction principle of the coproduct $A+B$. Therefore, we must construct two functions:
  \begin{align*}
    f & : A\to A \\*
    g & : B\to A.
  \end{align*}
  The function $f$ is simply defined to be the identity function $\idfunc:A\to A$. Recall that we have assumed that $B$ is empty, so we have a function $\tilde{b}:B\to\emptyt$. Furthermore, we always have the function $\exfalso:\emptyt\to A$. Therefore, we can define $g\defeq \exfalso\circ \tilde{b}$ to complete the proof.
\end{proof}
\index{coproduct|)}
\index{inductive type!coproduct|)}

\subsection{The type of integers}
\index{integers|(}
The set of integers is usually defined as a quotient of the set $\N\times\N$, by the equivalence relation
\begin{equation*}
  ((n,m)\sim (n',m')) := (n+m' = n'+m).
\end{equation*}
We haven't introduced the identity type yet, in order to consider the type of identifications $n+m'=n'+m$, but more importantly there are no quotient types in Martin-L\"of's dependent type theory. We will only discuss quotient types in \cref{sec:set-quotients} after we have assumed the univalence axiom and propositional truncations, because we will use the univalence axiom and propositional truncations to define them and derive their basic properties. Nevertheless, the type of integers is also definable in dependent type theory without set quotients, but we have to settle for a more pedestrian version of the integers that is defined using coproducts.

\begin{defn}
  We define the \define{integers}\index{Z@{$\Z$}|see {integers}} to be the type $\Z\defeq\N+(\unit+\N)$. The type of integers comes equipped with inclusion functions of the positive and negative integers\index{integers!in-pos@{$\inpos$}}\index{integers!in-neg@{$\inneg$}}
  \begin{alignat*}{2}
    \inpos & \defeq \inr\circ\inr\quad & & : \N\to \Z \\*
    \inneg & \defeq \inl\quad & & : \N \to \Z
  \end{alignat*}
  and with the constants\index{integers!-1 Z@{$-1_\Z$}}\index{integers!0 Z@{$0_\Z$}}\index{integers!1 Z@{$1_\Z$}}\index{-1 Z@{$-1_\Z$}}\index{0 Z@{$0_\Z$}}\index{1 Z@{$1_{\Z}$}}
  \begin{align*}
    -1_\Z & \defeq \inneg(0)\\*
    0_\Z & \defeq \inr(\inl(\ttt))\\*
    1_\Z & \defeq \inpos(0).
  \end{align*}
\end{defn}

The definition of the integers as the coproduct $\N+(\unit+\N)$ can be pictured as follows:
\begin{equation*}
  \begin{tikzcd}[column sep=0]
    \phantom{\unit+\N} & \unit \arrow[dr] & & \N \arrow[dl] \\
    \N \arrow[dr] & \phantom{\unit+\N} & \unit+\N \arrow[dl] & \phantom{\unit+\N} \\
    & \Z
  \end{tikzcd}
\end{equation*}

\begin{rmk}\label{lem:Z_ind}
  The type of integers is built entirely out of inductive types. Therefore it is possible to derive an induction principle especially tailored for the type $\Z$, which can be used to define the basic operations on $\Z$, such as the successor map, addition and multiplication. This induction principle asserts that for any type family $P$ over $\Z$,  we can define a dependent function $f:\prd{k:\Z}P(k)$ recursively by
  \begin{align*}
    f(-1_\Z) & \defeq p_{-1} \\*
    f(\inneg(\succN(n))) & \defeq p_{-S}(n,f(\inneg(n))) \\*
    f(0_\Z) & \defeq p_{0} \\*
    f(1_\Z) & \defeq p_{1} \\*
    f(\inpos(\succN(n))) & \defeq p_S(n,f(\inpos(n))),
  \end{align*}
  where the types of $p_{-1}$, $p_{-S}$, $p_0$, $p_1$, and $p_S$ are 
  \begin{align*}
    p_{-1} & :P(-1_\Z) \\*
    p_{-S} & : \prd{n:\N}P(\inneg(n))\to P(\inneg(\succN(n)))\\*
    p_{0} & : P(0_\Z) \\*
    p_{1} & : P(1_\Z) \\*
    p_{S} & : \prd{n:\N}P(\inpos(n))\to P(\inpos(\succN(n))).
  \end{align*}
\end{rmk}

\begin{defn}
We define the \define{successor function}\index{successor function!on Z@{on $\Z$}|textbf} on the integers $\succZ:\Z\to\Z$\index{succ Z@{$\succZ$}|textbf}\index{integers!succ Z@{$\succZ$}|textbf} using the induction principle of \cref{lem:Z_ind}, taking
\begin{align*}
\succZ(-1_\Z) & \defeq \zeroZ \\*
\succZ(\inneg(\succN(n))) & \defeq \inneg(n) \\*
\succZ(0_\Z) & \defeq \oneZ \\*
\succZ(1_\Z) & \defeq \inpos(1_\N) \\*
\succZ(\inpos(\succN(n))) & \defeq \inpos(\succN(\succN(n))).
\end{align*}
\end{defn}
\index{integers|)}

\subsection{Dependent pair types}

\index{dependent pair type|(}
\index{inductive type!dependent pair type|(}

Given a type family $B$ over $A$, we may consider pairs $(a,b)$ of terms, where $a:A$ and $b:B(a)$. Note that the type of $b$ depends on the first term in the pair. Therefore we call such a pair a \define{dependent pair}\index{dependent pair|textbf}. The type of such dependent pairs is the inductive type that is generated by the dependent pairs.

\begin{defn}
  Consider a type family $B$ over $A$.
  The \define{dependent pair type} (or \define{$\Sigma$-type}) \index{dependent pair type|textbf}\index{S-type@{$\Sigma$-type}|see {dependent pair type}}\index{S-type@{$\Sigma$-type}|textbf}is defined to be the inductive type $\sm{x:A}B(x)$ equipped with a \define{pairing function}\index{pairing function|textbf}\index{pair@{$\pair$}|textbf}\index{dependent pair type!pair@{$\pair$}|textbf}
  \begin{equation*}
    \pair :\prd{x:A} \Big(B(x)\to \sm{y:A}B(y)\Big).
  \end{equation*}
  The induction principle\index{induction principle!of Sigma types@{of $\Sigma$-types}|textbf}\index{dependent pair type!induction principle|textbf} for $\sm{x:A}B(x)$ asserts that for any family of types $P(p)$ indexed by $p:\sm{x:A}B(x)$, there is a function\index{dependent pair type!indSigma@{$\indSigma$}|textbf}\index{ind Sigma@{$\indSigma$}|textbf}
  \begin{equation*}
    \indSigma:\Big(\prd{x:A}\prd{y:B(x)}P(\pair(x,y))\Big)\to\Big(\prd{z:\sm{x:A}B(x)}P(z)\Big).
  \end{equation*}
  satisfying the computation rule\index{computation rules!for S-types@{for $\Sigma$-types}|textbf}\index{dependent pair type!computation rule|textbf}
  \begin{equation*}
    \indSigma(g,\pair(x,y))\jdeq g(x,y).
  \end{equation*}
  Alternatively, a definition of a dependent function $f:\prd{z:\sm{x:A}B(x)}P(z)$ by induction using a function $g:\prd{x:A}\prd{y:B(x)}P((x,y))$ can be presented by pattern matching as
  \begin{equation*}
    f(\pair(x,y))\defeq g(x,y).
  \end{equation*}
  We will usually write $(x,y)$ for $\pair(x,y)$\index{(x,y)@{$(x,y)$}|see {dependent pair}}\index{(x,y)@{$(x,y)$}|textbf}.
\end{defn}

The induction principle of $\Sigma$-types can be used to define the projection functions.

\begin{defn}
  Consider a type $A$ and a type family $B$ over $A$.
  \begin{enumerate}
  \item The \define{first projection map}\index{first projection map|textbf}\index{projection map!first projection|textbf}\index{dependent pair type!pr 1@{$\proj 1$}|textbf}\index{pr 1@{$\proj 1$}|textbf}\index{pr 1@{$\proj 1$}|see{first projection map}}
    \begin{equation*}
      \proj 1:\Big(\sm{x:A}B(x)\Big)\to A
    \end{equation*}
    is defined by induction as
    \begin{equation*}
      \proj 1(x,y) \defeq x.
    \end{equation*}
  \item The \define{second projection map}\index{second projection map|textbf}\index{projection map!second projection|textbf}\index{dependent pair type!pr 2@{$\proj 2$}|textbf}\index{pr 2@{$\proj 2$}|textbf}\index{pr 2@{$\proj 2$}|see{second projection map}} is a dependent function
    \begin{equation*}
      \proj 2 : \prd{p:\sm{x:A}B(x)} B(\proj 1(p)),
    \end{equation*}
    defined by induction as
    \begin{equation*}
      \proj 2(x,y) \defeq y.
    \end{equation*}
  \end{enumerate}
\end{defn}
\index{dependent pair type|)}
\index{inductive type!dependent pair type|)}

\begin{rmk}
  If we want to construct a function
  \begin{equation*}
    f:\prd{z:\sm{x:A}B(x)}P(z)
  \end{equation*}
  by $\Sigma$-induction, then we get to assume a pair $(x,y)$ consisting of $x:A$ and $y:B(x)$ and our goal will be to construct an element of type $P(x,y)$. The induction principle of $\Sigma$-types is therefore converse to the \define{currying operation}, a familiar concept from the theory of programming languages, which is given by the function
  \begin{equation*}
    \evpair : \Big(\prd{z:\sm{x:A}B(x)}P(z)\Big)\to \Big(\prd{x:A}\prd{y:B(x)}P(x,y)\Big)
  \end{equation*}
  given by $f\mapsto\lam{x}\lam{y}f(x,y)$. The induction principle $\indSigma$ is therefore also known as the \define{uncurrying operation}. 
\end{rmk}

\index{cartesian product type|(}
\index{inductive type!cartesian product|(}
A common special case of the $\Sigma$-type occurs when the $B$ is a constant family over $A$, i.e., when $B$ is just a type weakened by $A$.
In this case, the type $\sm{x:A}B$ is the type of \emph{ordinary} pairs $(x,y)$ where $x:A$ and $y:B$, where the type of $y$ does not depend on $x$. The type of ordinary pairs $(x,y)$ consisting of $x:A$ and $y:B$ is of course the \emph{product} of $A$ and $B$, so we see that product types arise as a special case of $\Sigma$-types in a similar way to how function types were defined as a special case of $\Pi$-types.

\begin{defn}
  Consider two types $A$ and $B$. Then we define the \define{(cartesian) product}\index{cartesian product type|textbf}\index{product of types|textbf}\index{A x B@{$A\times B$}|see {cartesian product}}\index{A x B@{$A\times B$}|textbf} $A\times B$ of $A$ and $B$ by
  \begin{equation*}
    A\times B \defeq \sm{x:A}B.
  \end{equation*}
\end{defn}

\begin{rmk}
  Since $A\times B$ is defined as a $\Sigma$-type, it follows that cartesian products also satisfy the induction principle of $\Sigma$-types. In this special case, the induction principle\index{induction principle!of cartesian products|textbf}\index{cartesian product type!induction principle|textbf} for $A\times B$ asserts that for any type family $P$ over $A\times B$ there is a function\index{ind times@{$\ind{\times}$}|textbf}\index{cartesian product type!indtimes@{$\ind{\times}$}|textbf}
\begin{equation*}
\ind{\times} : \Big(\prd{x:A}\prd{y:B}P(x,y)\Big)\to\Big(\prd{z:A\times B} P(z)\Big)
\end{equation*}
that satisfies the computation rule\index{computation rules!for cartesian products}\index{cartesian product type!computation rule|textbf}
\begin{align*}
\ind{\times}(g,(x,y)) & \jdeq g(x,y).
\end{align*}
\end{rmk}

The projection maps are defined similarly to the projection maps of $\Sigma$-types. When one thinks of types as propositions\index{propositions as types!conjunction}, then $A\times B$ is interpreted as the conjunction of $A$ and $B$.
\index{cartesian product type|)}
\index{inductive type!cartesian product|)}

\begin{exercises}
  \exitem
  \begin{subexenum}
  \item \label{ex:int_pred}\index{integers|(}\index{predecessor function|textbf}\index{integers!pred Z@{$\predZ$}|textbf}\index{pred Z@{$\predZ$}|textbf}Define the predecessor function $\predZ:\Z\to \Z$.
  \item \label{ex:int_group_ops}Define the group operations\index{add Z@{$\addZ$}|textbf}\index{integers!add Z@{$\addZ$}|textbf}\index{neg Z@{$\negZ$}}\index{integers!neg Z@{$\negZ$}}
    \begin{align*}
      \addZ & : \Z \to (\Z \to \Z) \\*
      \negZ & : \Z \to \Z.
    \end{align*}
  \item \label{ex:mulZ}Define the multiplication operation $\mulZ : \Z \to (\Z \to \Z)$.\index{mul Z@{$\mulZ$}}\index{integers!mul Z@{$\mulZ$}}
  \end{subexenum}
  \exitem \label{ex:bool}The type of \define{booleans}\index{booleans|textbf}\index{bool@{$\bool$}|see {booleans}}\index{bool@{$\bool$}|textbf} is defined to be an inductive type $\bool$ that comes equipped with\index{booleans!true@{$\btrue$}|textbf}\index{booleans!false@{$\bfalse$}|textbf}\index{false@{$\bfalse$}|textbf}\index{true@{$\btrue$}|textbf}
  \begin{equation*}
    \bfalse : \bool\qquad\text{and}\qquad\btrue : \bool.
  \end{equation*}
  The induction principle\index{induction principle!of the booleans|textbf}\index{booleans!induction principle|textbf} of the booleans asserts that for any family of types $P(x)$ indexed by $x:\bool$, there is a term\index{ind_bool@{$\indbool$}|textbf}
  \begin{equation*}
    \indbool : P(\bfalse)\to \Big(P(\btrue)\to \prd{x:\bool}P(x)\Big)
  \end{equation*}
  for which the computation rules\index{computation rules!for the booleans|textbf}\index{booleans!computation rules|textbf}
  \begin{align*}
    \indbool(p_0,p_1,\bfalse) & \jdeq p_0 \\*
    \indbool(p_0,p_1,\btrue) & \jdeq p_1
  \end{align*}
  hold.
  \begin{subexenum}
  \item Construct the \define{boolean negation} function $\negbool:\bool\to\bool$\index{neg-bool@{$\negbool$}|textbf}\index{booleans!neg-bool@{$\negbool$}|textbf}.
  \item Construct the \define{boolean conjunction} operation $\blank\land\blank : \bool\to(\bool\to\bool)$.\index{boolean conjunction|textbf}\index{booleans!conjunction|textbf}
  \item Construct the \define{boolean disjunction} operation $\blank\lor\blank : \bool\to(\bool\to\bool)$.\index{boolean disjunction|textbf}\index{booleans!disjunction|textbf}
  \end{subexenum}
  \exitem Let $P$ and $Q$ be types. We will write $P\leftrightarrow Q$ for the type of \define{bi-implications}\index{bi-implication|textbf} ${(P\to Q)}\times {(Q\to P)}$. Use the fact that $\neg P$\index{negation} is defined as the type $P\to\emptyt$ of functions from $P$ to the empty type to give type theoretic proofs of the constructive tautologies in this exercise.\label{ex:dne-dec}
  \begin{subexenum}
  \item \label{ex:no-fixed-points-neg}Show that
    \begin{enumerate}
    \item $\neg(P\times \neg P)$
    \item $\neg(P\leftrightarrow \neg P)$.
    \end{enumerate}
  \item \label{ex:dn-monad}Construct the following maps in the structure of the \define{double negation monad}:
    \begin{enumerate}
    \item $P\to\neg\neg P$
    \item $(P\to Q)\to(\neg\neg P\to\neg\neg Q)$
    \item $(P\to \neg\neg Q)\to (\neg\neg P \to\neg\neg Q)$.
    \end{enumerate}
  \item Prove that the following double negations of classical laws hold:
    \begin{enumerate}
    \item $\neg\neg(\neg\neg P \to P)$
    \item $\neg\neg(((P\to Q)\to P)\to P)$
    \item $\neg\neg((P\to Q)+(Q\to P))$
    \item $\neg\neg(P+\neg P)$.
    \end{enumerate}
  \item \label{ex:dne-is-decidable}Show that
    \begin{enumerate}
    \item $(P+\neg P)\to(\neg\neg P\to P)$
    \item $\neg\neg(Q\to P)\leftrightarrow ((P+\neg P)\to (Q\to P))$.    
    \end{enumerate}
  \item Prove the following tautologies, showing that $\neg P$, $P\to\neg\neg Q$, and $\neg\neg P\times\neg\neg Q$ are \define{double negation stable}:
    \begin{enumerate}
    \item $\neg\neg\neg P \to \neg P$
    \item $\neg\neg(P \to \neg\neg Q)\to (P\to\neg\neg Q)$
    \item $\neg\neg((\neg\neg P)\times(\neg\neg Q))\to (\neg\neg P)\times(\neg\neg Q)$.
    \end{enumerate}
  \item Show that
    \begin{enumerate}
    \item $\neg\neg(P\times Q)\leftrightarrow (\neg\neg P)\times(\neg\neg Q)$
    \item $\neg\neg(P+Q)\leftrightarrow \neg (\neg P \times \neg Q)$
    \item $\neg\neg(P\to Q)\leftrightarrow (\neg\neg P\to\neg\neg Q)$.
    \end{enumerate}
  \end{subexenum}
\exitem \label{ex:lists}For any type $A$ we can define the type $\lst(A)$\index{list A@{$\lst(A)$}|see {lists in $A$}}\index{list A@{$lst(A)$}|textbf} of \define{lists}\index{lists in A@{lists in $A$}|textbf}\index{inductive type!lists of elements of A@{lists of elements of $A$}|textbf} of elements of $A$ as the inductive type with constructors\index{lists in A@{lists in $A$}!nil@{$\nil$}|textbf}\index{nil@{$\nil$}|textbf}\index{cons(a,l)@{$\cons(a,l)$}|textbf}\index{lists in A@{lists in $A$}!cons@{$\cons$}|textbf}
  \begin{align*}
    \nil & : \lst(A) \\*
    \cons & : A \to (\lst(A) \to \lst(A)).
  \end{align*}
  \begin{subexenum}
  \item Write down the induction principle and the computation rules for $\lst(A)$.\index{induction principle!of list A@{of $\lst(A)$}|textbf}\index{lists in A@{lists in $A$}!induction principle|textbf}
  \item Let $A$ and $B$ be types, suppose that $b:B$, and consider a binary operation $\mu:A\to (B \to B)$. Define a function\index{fold-list@{$\foldlist$}|textbf}\index{lists in A@{lists in $A$}!fold-list@{$\foldlist$}|textbf}
    \begin{equation*}
      \foldlist(\mu) : \lst(A)\to B
    \end{equation*}
    that iterates the operation $\mu$, starting with $\foldlist(\mu,\nil)\defeq b$.
  \item Define the operation
    \begin{equation*}
      \maplist : (A\to B) \to (\lst(A)\to\lst(B))
    \end{equation*}
    for any two types $A$ and $B$.
  \item Define a function $\lengthlist:\lst(A)\to\N$.\index{length-list@{$\lengthlist$}|textbf}\index{lists in A@{lists in $A$}!length-list@{$\lengthlist$}|textbf}
  \item Define the functions\index{sum-list@{$\sumlist$}|textbf}\index{lists in A@{lists in $A$}!sum-list@{$\sumlist$}|textbf}
    \begin{align*}
      \sumlist & : \lst(\N) \to \N \\
      \productlist & : \lst(\N)\to\N,
    \end{align*}
    where $\sumlist$ adds all the elements in a list of natural numbers, and $\productlist$ takes their product.
  \item Define a function\index{concat-list@{$\concatlist$}|textbf}\index{lists in A@{lists in $A$}!concat-list@{$\concatlist$}|textbf}\index{concatenation!of lists|textbf}
    \begin{equation*}
      \concatlist : \lst(A) \to (\lst(A) \to \lst(A))
    \end{equation*}
    that concatenates any two lists of elements in $A$.
  \item Define a function\index{flatten-list@{$\flattenlist$}|textbf}\index{lists in A@{lists in $A$}!flatten-list@{$\flattenlist$}|textbf}
    \begin{equation*}
      \flattenlist : \lst(\lst(A)) \to \lst(A)
    \end{equation*}
    that concatenates all the lists in a lists of lists in $A$.
  \item Define a function $\reverselist : \lst(A) \to \lst(A)$ that reverses the order of the elements in any list.\index{reverse-list@{$\reverselist$}|textbf}\index{lists in A@{lists in $A$}!reverse-list@{$\reverselist$}|textbf}
  \end{subexenum}
\end{exercises}


\section{Identity types}\label{sec:identity}

\index{identity type|(}
\index{inductive type!identity type|(}
From the perspective of types as proof-relevant propositions, how should we think of \emph{equality} in type theory? Given a type $A$, and two elements $x,y:A$, the equality $\id{x}{y}$ should again be a type. Indeed, we want to \emph{use} type theory to prove equalities. \emph{Dependent} type theory provides us with a convenient setting for this: the identity type $\id{x}{y}$ is dependent on $x,y:A$. 

Then, if $\id{x}{y}$ is to be a type, how should we think of the elements of $\id{x}{y}$. An element $p:\id{x}{y}$ witnesses that $x$ and $y$ are equal elements of type $A$. In other words $p:\id{x}{y}$ is an \emph{identification} of $x$ and $y$. In a proof-relevant world, there might be many elements of type $\id{x}{y}$. I.e., there might be many identifications of $x$ and $y$. And, since $\id{x}{y}$ is itself a type, we can form the type $\id{p}{q}$ for any two identifications $p,q:\id{x}{y}$. That is, since $\id{x}{y}$ is a type, we may also use the type theory to prove things \emph{about} identifications (for instance, that two given such identifications can themselves be identified), and we may use the type theory to perform constructions with them. As we will see in this section, we can give every type a groupoidal structure.

Clearly, the equality type should not just be any type dependent on $x,y:A$. Then how do we form the equality type, and what ways are there to use identifications in constructions in type theory? The answer to both these questions is that we will form the identity type as an \emph{inductive} type, generated by just a reflexivity identification providing an identification of $x$ to itself. The induction principle then provides us with a way of performing constructions with identifications, such as concatenating them, inverting them, and so on. Thus, the identity type is equipped with a reflexivity element, and further possesses the structure that are generated by its induction principle and by the type theory. This inductive construction of the identity type is elegant, beautifully simple, but far from trivial!

The situation where two elements can be identified in possibly more than one way is analogous to the situation in \emph{homotopy theory}, where two points of a space can be connected by possibly more than one \emph{path}. Indeed, for any two points $x,y$ in a space, there is a \emph{space of paths} from $x$ to $y$. Moreover, between any two paths from $x$ to $y$ there is a space of \emph{homotopies} between them, and so on. From \cref{chap:uf} on we will take full advantage of this idea in order to develop the univalent foundations of mathematics.

\subsection{The inductive definition of identity types}

\begin{defn}
  Consider a type $A$ and let $a:A$. Then we define the \define{identity type}\index{identity type|textbf} of $A$ at $a$ as an inductive family of types $a =_A x$\index{a = x@{$a = x$}|see {identity type}} indexed by $x:A$, of which the constructor is\index{refl@{$\refl{}$}|textbf}\index{identity type!refl@{$\refl{}$}|textbf}
  \begin{equation*}
    \refl{a}:a=_Aa.
  \end{equation*}
  The induction principle of the identity type\index{identity type!induction principle|textbf}\index{induction principle!of the identity type|textbf} postulates that for any family of types $P(x,p)$ indexed by $x:A$ and $p:a=_A x$, there is a function\index{path-ind@{$\pathind$}|textbf}\index{identity type!path-ind@{$\pathind$}|textbf}
  \begin{equation*}
    \pathind_a:P(a,\refl{a}) \to \prd{x:A}\prd{p:a=_A x} P(x,p)
  \end{equation*}
  that satisfies $\pathind_a(u,a,\refl{a})\jdeq u$, give $u:P(a,\refl{a})$.

  An element of type $a=_A x$ is also called an \define{identification}\index{identification|textbf}\index{identity type!identification|textbf} of $a$ with $x$, and sometimes it is called a \define{path}\index{path|textbf}\index{identity type!path|textbf} from $a$ to $x$.
The induction principle for identity types is sometimes called \define{identification elimination}\index{identification elimination|textbf}\index{induction principle!identification elimination|textbf}\index{identity type!identification elimination|textbf} or \define{path induction}\index{path induction|textbf}\index{identity type!path induction|textbf}\index{induction principle!path induction|textbf}. We also write $\idtypevar{A}$\index{Id A@{$\idtypevar{A}$}|see {identity type}}\index{Id A@{$\idtypevar{A}$}|textbf} for the identity type on $A$, and often we write $a=x$ for the type of identifications of $a$ with $x$, omitting reference to the ambient type $A$.
\end{defn}

\begin{rmk}
  We see that the identity type is not just an inductive type, like the inductive types $\N$, $\emptyt$, and $\unit$ for example, but it is an inductive \emph{family} of types. Even though we have a type $a=_A x$ for any $x:A$, the constructor only provides an element $\refl{a}:a=_A a$, identifying $a$ with itself. The induction principle then asserts that in order to prove something about all identifications of $a$ with some $x:A$, it suffices to prove this assertion about $\refl{a}$ only. We will see in the next sections that this induction principle is strong enough to derive many familiar facts about equality, namely that it is a symmetric and transitive relation, and that all functions preserve equality.
\end{rmk}

\begin{rmk}
  \index{rules!identity type|(}\index{identity type!rules|(}
  Since the identity types require getting used to, we provide the formal rules
  for identity types. The identity type is formed by the formation rule:
  \begin{prooftree}
    \AxiomC{$\Gamma\vdash a:A$}
    \UnaryInfC{$\Gamma,x:A\vdash a=_A x~\type$}
  \end{prooftree}
  The constructor of the identity type is then given by the introduction rule:
  \begin{prooftree}
    \AxiomC{$\Gamma\vdash a:A$}
    \UnaryInfC{$\Gamma\vdash \refl{a}:a=_A a$}
  \end{prooftree}
  The induction principle is now given by the elimination rule:
  \begin{prooftree}
    \AxiomC{$\Gamma\vdash a:A$}
    \AxiomC{$\Gamma,x:A,p:a=_A x\vdash P(x,p)~\type$}
    \BinaryInfC{$\Gamma\vdash \pathind_a:P(a,\refl{a})\to\prd{x:A}\prd{p:a=_A x}P(x,p)$}
  \end{prooftree}
  And finally the computation rule is:
  \begin{prooftree}
    \AxiomC{$\Gamma\vdash a:A$}
    \AxiomC{$\Gamma,x:A,p:a=_A x\vdash P(x,p)~\type$}
    \BinaryInfC{$\Gamma,u:P(a,\refl{a}) \vdash \pathind_a(u,a,\refl{a})\jdeq u : P(a,\refl{a})$}
  \end{prooftree}
  \index{rules!identity type|)}\index{identity type!rules|)}
\end{rmk}

\begin{rmk}
  One might wonder whether it is also possible to form the identity type at a \emph{variable} of type $A$, rather than at an element. This is certainly possible: since we can form the identity type in \emph{any} context, we can form the identity type at a variable $x:A$ as follows:
  \begin{prooftree}
    \AxiomC{$\Gamma,x:A\vdash x:A$}
    \UnaryInfC{$\Gamma,x:A,y:A\vdash x=_A y~\type$}
  \end{prooftree}
  In this way we obtain the `binary' identity type. Its constructor is then also indexed by $x:A$. We have the following introduction rule
  \begin{prooftree}
    \AxiomC{$\Gamma,x:A\vdash x:A$}
    \UnaryInfC{$\Gamma,x:A\vdash \refl{x}:x=_A x$}
  \end{prooftree}
  and similarly we have elimination and computation rules.
\end{rmk}

\subsection{The groupoidal structure of types}\label{sec:groupoid}
\index{groupoid laws!of identifications|(}
We show that identifications can be \emph{concatenated} and \emph{inverted}, which corresponds to the transitivity and symmetry of the identity type.

\begin{defn}\label{defn:id_concat}
Let $A$ be a type. We define the \define{concatenation}\index{concatenation!of identifications|textbf}\index{concat@{$\concat$}|textbf}\index{identity type!concatenation|textbf} operation
\begin{equation*}
\concat : \prd{x,y,z:A} (\id{x}{y})\to ((\id{y}{z})\to (\id{x}{z})).
\end{equation*}
We will write $\ct{p}{q}$ for $\concat(p,q)$.
\end{defn}

\begin{constr}
  We first construct a function
  \begin{equation*}
    f(x):\prd{y:A}(x=y)\to\prd{z:A}(y=z)\to(x=z)
  \end{equation*}
  for any $x:A$. By the induction principle for identity types, it suffices to construct
  \begin{equation*}
    f(x,x,\refl{x}):\prd{z:A} (x=z)\to(x=z).
  \end{equation*}
  Here we have the function $\lam{z}\idfunc[(x=z)]$. The function $f(x)$ we obtain via identity elimination is explicitly thus defined as
  \begin{equation*}
    f(x)\defeq\pathind_x(\lam{z}\idfunc):\prd{y:A} (x=y)\to \prd{z:A} (y=z)\to (x=z).
  \end{equation*}
  To finish the construction of $\concat$, we use \cref{ex:swap} to swap the order of the third and fourth variable of $f$, i.e., we define
  \begin{equation*}
    \concat_{x,y,z}(p,q):=f(x,y,p,z,q).\qedhere
  \end{equation*}
\end{constr}

\begin{defn}\label{defn:id_inv}
Let $A$ be a type. We define the \define{inverse operation}\index{inverse operation!for identifications|textbf}\index{inv@{$\invfunc$}|textbf}\index{identity type!inverse operation|textbf}
\begin{equation*}
\invfunc:\prd{x,y:A} (x=y)\to (y=x).
\end{equation*}
Most of the time we will write $p^{-1}$ for $\invfunc(p)$.
\end{defn}

\begin{constr}
By the induction principle for identity types, it suffices to construct
\begin{equation*}
\invfunc(\refl{x}): x=x,
\end{equation*}
for any $x:A$. Here we take $\invfunc(\refl{x})\defeq \refl{x}$.
\end{constr}

The next question is whether the concatenation and inverting operations on identifications behave as expected. More concretely: is concatenation of identifications associative, does it satisfy the unit laws, and is the inverse of an identification indeed a two-sided inverse?

For example, in the case of associativity we are asking to compare the identifications
\begin{equation*}
  \ct{(\ct{p}{q})}{r}\qquad\text{and}\qquad\ct{p}{(\ct{q}{r})}
\end{equation*}
for any $p:x=y$, $q:y=z$, and $r:z=w$ in a type $A$. The computation rules of the identity type are not strong enough to conclude that $\ct{(\ct{p}{q})}{r}$ and $\ct{p}{(\ct{q}{r})}$ are judgmentally equal. However, both $\ct{(\ct{p}{q})}{r}$ and $\ct{p}{(\ct{q}{r})}$ are elements of the same type: they are identifications of type $x=w$. Since the identity type is a type like any other, we can ask whether there is an \emph{identification}
\begin{equation*}
\ct{(\ct{p}{q})}{r}=\ct{p}{(\ct{q}{r})}.
\end{equation*}
This is a very useful idea: while it is often impossible to show that two elements of the same type are judgmentally equal, it may be the case that those two elements can be \emph{identified}. Indeed, we identify two elements by constructing an element of the identity type, and we can use all the type theory at our disposal in order to construct such an element. In this way we can show, for example, that addition on the natural numbers or on the integers is associative and satisfies the unit laws. And indeed, here we will show that concatenation of identifications is associative and satisfies the unit laws.

\begin{defn}\label{defn:id_assoc}
  Let $A$ be a type and consider three consecutive identifications
  \begin{equation*}
    \begin{tikzcd}
      x \arrow[r,equals,"p"] & y \arrow[r,equals,"q"] & z \arrow[r,equals,"r"] & w
    \end{tikzcd}
  \end{equation*}
  in $A$. We define the \define{associator}\index{associativity!of concatenation of identifications}
  \begin{equation*}
    \assoc(p,q,r) : \ct{(\ct{p}{q})}{r}=\ct{p}{(\ct{q}{r})}.
  \end{equation*}
\end{defn}

\begin{constr}
By the induction principle for identity types it suffices to show that
\begin{equation*}
\prd{z:A}\prd{q:x=z}\prd{w:A}\prd{r:z=w} \ct{(\ct{\refl{x}}{q})}{r}= \ct{\refl{x}}{(\ct{q}{r})}.
\end{equation*}
Let $q:x=z$ and $r:z=w$. Note that by the computation rule of identity types we have a judgmental equality $\ct{\refl{x}}{q}\jdeq q$. Therefore we conclude that
\begin{equation*}
  \ct{(\ct{\refl{x}}{q})}{r}\jdeq \ct{q}{r}.
\end{equation*}
Similarly we have a judgmental equality $\ct{\refl{x}}{(\ct{q}{r})}\jdeq \ct{q}{r}$. Thus we see that the left-hand side and the right-hand side in
\begin{equation*}
  \ct{(\ct{\refl{x}}{q})}{r}=\ct{\refl{x}}{(\ct{q}{r})}
\end{equation*}
are judgmentally equal, so we can simply define $\assoc(\refl{x},q,r)\defeq\refl{\ct{q}{r}}$.
\end{constr}

\begin{defn}\label{defn:id_unit}
Let $A$ be a type. We define the left and right \define{unit law operations}\index{unit laws!for concatenation of identifications}, which assigns to each $p:x=y$ the identifications\index{left unit@{$\leftunit$}|textbf}\index{right unit@{$\rightunit$}|textbf}
\begin{align*}
\leftunit(p) & : \ct{\refl{x}}{p}=p \\
\rightunit(p) & : \ct{p}{\refl{y}}=p,
\end{align*}
respectively.
\end{defn}

\begin{constr}
By identification elimination it suffices to construct
\begin{align*}
\leftunit(\refl{x}) & : \ct{\refl{x}}{\refl{x}} = \refl{x} \\
\rightunit(\refl{x}) & : \ct{\refl{x}}{\refl{x}} = \refl{x}.
\end{align*}
In both cases we take $\refl{\refl{x}}$.
\end{constr}

\begin{defn}\label{defn:id_invlaw}
Let $A$ be a type. We define left and right \define{inverse law operations}\index{inverse law operations!for identifications}\index{left inv@{$\leftinv$}|textbf}\index{right inv@{$\rightinv$}|textbf}
\begin{align*}
\leftinv(p) & : \ct{p^{-1}}{p} = \refl{y} \\
\rightinv(p) & : \ct{p}{p^{-1}} = \refl{x}.
\end{align*}
\end{defn}

\begin{constr}
By identification elimination it suffices to construct
\begin{align*}
\leftinv(\refl{x}) & : \ct{\refl{x}^{-1}}{\refl{x}} = \refl{x} \\
\rightinv(\refl{x}) & : \ct{\refl{x}}{\refl{x}^{-1}} = \refl{x}.
\end{align*}
Using the computation rules we see that
\begin{equation*}
\ct{\refl{x}^{-1}}{\refl{x}}\jdeq \ct{\refl{x}}{\refl{x}}\jdeq\refl{x},
\end{equation*}
so we define $\leftinv(\refl{x})\defeq \refl{\refl{x}}$. Similarly it follows from the computation rules that
\begin{equation*}
\ct{\refl{x}}{\refl{x}^{-1}} \jdeq \refl{x}^{-1}\jdeq \refl{x}
\end{equation*}
so we again define $\rightinv(\refl{x})\defeq\refl{\refl{x}}$. 
\end{constr}

\begin{rmk}
  We have seen that the associator, the unit laws, and the inverse laws, are all proven by constructing an identification of identifications. And indeed, there is nothing that would stop us from considering identifications of those identifications of identifications. We can go up as far as we like in the \emph{tower of identity types}\index{tower of identity types}\index{identity type!tower of identity types}, which is obtained by iteratively taking identity types.

  The iterated identity types give types in homotopy type theory a very intricate structure. One important way of studying this structure is via the homotopy groups of types, a subject that we will gradually be working towards.
\end{rmk}
\index{groupoid laws!of identifications|)}

\subsection{The action on identifications of functions}

\index{action on paths|(}
\index{identity type!action on paths|(}
Using the induction principle of the identity type we can show that every function preserves identifications.
In other words, every function sends identified elements to identified elements.
Note that this is a form of continuity for functions in type theory: if there is an identification that identifies two points $x$ and $y$ of a type $A$, then there also is an identification that identifies the values $f(x)$ and $f(y)$ in the codomain of $f$. 

\begin{defn}\label{defn:ap}
Let $f:A\to B$ be a map. We define the \define{action on paths}\index{function!action on paths|textbf}\index{identity type!action on paths|textbf}\index{action on paths|textbf} of $f$ as an operation\index{ap f@{$\apfunc{f}$}|see {action on paths}}\index{ap f@{$\apfunc{f}$}|textbf}
\begin{equation*}
\apfunc{f} : \prd{x,y:A} (\id{x}{y})\to(\id{f(x)}{f(y)}).
\end{equation*}
Moreover, there are operations\index{ap-id@{$\apid$}|textbf}\index{action on paths!ap-id@{$\apid$}|textbf}\index{ap-comp@{$\apcomp$}|textbf}\index{action on paths!ap-comp@{$\apcomp$}|textbf}
\begin{align*}
\apid_A & : \prd{x,y:A}\prd{p:\id{x}{y}} \id{p}{\ap{\idfunc[A]}{p}} \\
\apcomp(f,g) & : \prd{x,y:A}\prd{p:\id{x}{y}} \id{\ap{g}{\ap{f}{p}}}{\ap{g\circ f}{p}}.
\end{align*}
\end{defn}

\begin{constr}
First we define $\apfunc{f}$ by the induction principle of identity types, taking
\begin{equation*}
\apfunc{f}(\refl{x})\defeq \refl{f(x)}.
\end{equation*}
Next, we construct $\apid_A$ by the induction principle of identity types, taking
\begin{equation*}
\apid_A(\refl{x}) \defeq \refl{\refl{x}}.
\end{equation*}
Finally, we construct $\apcomp(f,g)$ by the induction principle of identity types, taking
\begin{equation*}
\apcomp(f,g,\refl{x}) \defeq \refl{\refl{g(f(x))}}.\qedhere
\end{equation*}
\end{constr}

\begin{defn}\label{defn:ap-preserve}
Let $f:A\to B$ be a map. Then there are identifications\index{ap-refl@{$\aprefl$}|textbf}\index{ap-inv@{$\apinv$}|textbf}\index{ap-concat@{$\apconcat$}|textbf}\index{action on paths!ap-refl@{$\aprefl$}|textbf}\index{action on paths!ap-inv@{$\apinv$}|textbf}\index{action on paths!ap-concat@{$\apconcat$}|textbf}
\begin{align*}
\aprefl(f,x) & : \id{\ap{f}{\refl{x}}}{\refl{f(x)}} \\
\apinv(f,p) & : \id{\ap{f}{p^{-1}}}{\ap{f}{p}^{-1}} \\
\apconcat(f,p,q) & : \id{\ap{f}{\ct{p}{q}}}{\ct{\ap{f}{p}}{\ap{f}{q}}}
\end{align*}
for every $p:\id{x}{y}$ and $q:\id{x}{y}$.
\end{defn}

\begin{constr}
To construct $\aprefl(f,x)$ we simply observe that ${\ap{f}{\refl{x}}}\jdeq {\refl{f(x)}}$, so we take
\begin{equation*}
\aprefl(f,x)\defeq\refl{\refl{f(x)}}.
\end{equation*}
We construct $\apinv(f,p)$ by identification elimination on $p$, taking
\begin{equation*}
\apinv(f,\refl{x}) \defeq \refl{\ap{f}{\refl{x}}}.
\end{equation*}
Finally we construct $\apconcat(f,p,q)$ by identification elimination on $p$, taking
\begin{equation*}
\apconcat(f,\refl{x},q)  \defeq \refl{\ap{f}{q}}.\qedhere
\end{equation*}
\end{constr}
\index{action on paths|)}
\index{identity type!action on paths|)}

\subsection{Transport}

\index{transport|(}
Dependent types also come with an action on identifications: the \emph{transport} functions.
Given an identification $p:\id{x}{y}$ in the base type $A$, we can transport any element $b:B(x)$ to the fiber $B(y)$.

\begin{defn}
Let $A$ be a type, and let $B$ be a type family over $A$.
We will construct a \define{transport}\index{transport|textbf}\index{type family!transport|textbf}\index{identity type!transport|textbf} operation\index{tr B@{$\tr_B$}|textbf}
\begin{equation*}
\tr_B:\prd{x,y:A} (\id{x}{y})\to (B(x)\to B(y)).
\end{equation*}
\end{defn}

\begin{constr}
We construct $\tr_B(p)$ by induction on $p:x=_A y$, taking
\begin{equation*}
\tr_B(\refl{x}) \defeq \idfunc[B(x)].\qedhere
\end{equation*}
\end{constr}

Thus we see that type theory cannot distinguish between identified elements $x$ and $y$, because for any type family $B$ over $A$ one obtains an element of $B(y)$ from the elements of $B(x)$.

As an application of the transport function we construct the \emph{dependent} action on paths\index{dependent action on paths|textbf} of a dependent function $f:\prd{x:A}B(x)$. Note that for such a dependent function $f$, and an identification $p:\id[A]{x}{y}$, it does not make sense to directly compare $f(x)$ and $f(y)$, since the type of $f(x)$ is $B(x)$ whereas the type of $f(y)$ is $B(y)$, which might not be exactly the same type. However, we can first \emph{transport} $f(x)$ along $p$, so that we obtain the element $\tr_B(p,f(x))$ which is of type $B(y)$. Now we can ask whether it is the case that $\tr_B(p,f(x))=f(y)$. The dependent action on paths of $f$ establishes this identification.

\begin{defn}\label{defn:apd}
Given a dependent function $f:\prd{a:A}B(a)$ and an identification $p:\id{x}{y}$ in $A$, we construct an identification\index{apd f@{$\apdfunc{f}$}|textbf}
\begin{equation*}
\apd{f}{p} : \id{\tr_B(p,f(x))}{f(y)}.
\end{equation*}
\end{defn}

\begin{constr}
The identification $\apd{f}{p}$ is constructed by the induction principle for identity types. Thus, it suffices to construct an identification
\begin{equation*}
\apd{f}{\refl{x}}:\id{\tr_B(\refl{x},f(x))}{f(x)}.
\end{equation*}
Since transporting along $\refl{x}$ is the identity function on $B(x)$, we simply take $\apd{f}{\refl{x}}\defeq\refl{f(x)}$. 
\end{constr}
\index{transport|)}

\subsection{The uniqueness of \texorpdfstring{$\refl{}$}{refl}}\label{sec:refl-unique}%

The identity type is an inductive \emph{family} of types. This has some subtle, but important implications. For instance, while the type $a=x$ indexed by $x:A$ is inductively generated by $\refl{a}$, the type $a=a$ is \emph{not} inductively generated by $\refl{a}$. Hence we cannot use the induction principle of identity types to show that $p=\refl{a}$ for any $p:a=a$. The obstacle, which prevents us from applying the induction principle of identity types in this case, is that the endpoint of $p:a=a$ is not free.

Nevertheless, the identity type $a=x$ is generated by a single element $\refl{a}:a=a$, so it is natural to wonder in what sense the reflexivity identification is unique. An identification with an element $a$ is specified by first giving the endpoint $x$ with which we seek to identify $a$, and then giving the identification $p:a=x$. It is therefore only the pair $(a,\refl{a})$ which is unique in the type of all pairs
\begin{equation*}
  (x,p):\sm{x:A}a=x.
\end{equation*}
We prove this fact in the following proposition.

\begin{prp}\label{prp:contraction-total-space-id}
  Consider an element $a:A$. Then there is an identification
  \begin{equation*}
    (a,\refl{a})=y
  \end{equation*}
  in the type $\sm{x:A}a=x$, for any $y:\sm{x:A}a=x$.
\end{prp}

\begin{proof}
  By $\Sigma$-induction it suffices to show that there is an identification
  \begin{equation*}
    (a,\refl{a})=(x,p)
  \end{equation*}
  for any $x:A$ and $p:a=x$. We proceed by the induction principle of identity types.
  Therefore it suffices to show that
  \begin{equation*}
    (a,\refl{a})=(a,\refl{a}).
  \end{equation*}
  We obtain such an identification by reflexivity.
\end{proof}

\cref{prp:contraction-total-space-id} shows that there is, up to identification, only one element in $\Sigma$-type of the identity type. Such types are called contractible, and they are the subject of \cref{sec:contractible}.

\subsection{The laws of addition on \texorpdfstring{$\N$}{ℕ}}\label{subsec:addN}

Now that we have introduced the identity type, we can start proving equations. We will prove here that there are identifications\index{unit laws!for addition on N@{for addition on $\N$}}\index{successor laws!for addition on N@{for addition on $\N$}}\index{associativity!of addition on N@{of addition on $\N$}}\index{commutativity!of addition on N@{of addition on $\N$}}\index{natural numbers!unit laws for addition}\index{natural numbers!successor laws for addition}\index{natural numbers!associativity of addition}\index{natural numbers!commutativity of addition}
\begin{align*}
  0+n & = n & m+0 & = m \\
  \succN(m)+n & = \succN(m+n) & m+\succN(n) & = \succN(m+n) \\
  (m+n)+k & = m+(n+k) & m+n & = n+m.
\end{align*}
The unit laws, associativity, and commutativity of addition are of course familiar. The successor laws will be useful to prove commutativity. In \cref{ex:semi-ring-laws-N} you will be asked to prove the laws of multiplication on $\N$. There will again be \emph{successor laws} as part of this exercise, because they are useful intermediate steps in the more complicated laws.

Recall that addition on the natural numbers is defined in such a way that
\begin{align*}
  m+0 & \jdeq m & m+\succN(n) & \jdeq \succN(m+n).
\end{align*}
These two judgmental equalities are all we currently know about the function $m,n\mapsto m+n$ on $\N$. Consequently, we will have to find ways to apply these two judgmental equalities in our proofs of the laws of addition. Of course, the judgmental equalities coincide with two of the six laws. For the remaining four laws, we will have to proceed by induction on $\N$.

\begin{prp}\label{prp:unit-laws-add-N}
  For any natural number $n$, there are identifications
  \begin{align*}
    \leftunitlawaddN(n) & : 0+n=n \\
    \rightunitlawaddN(n) & : n+0=n.
  \end{align*}
\end{prp}

\begin{proof}
  We can define
  \begin{equation*}
    \rightunitlawaddN(n)\defeq\refl{n},
  \end{equation*}
  because the computation rule for addition gives us that $n+0\jdeq n$.

  It remains to define the left unit law. We proceed by induction on $n$. In the base case we have to show that $0+0=0$, which holds by reflexivity. For the inductive step, assume that we have an identification $p:0+n=n$. Our goal is to show that $0+\succN(n)=\succN(n)$. However, it suffices to construct an identification
  \begin{equation*}
    \succN(0+n)=\succN(n),
  \end{equation*}
  because by the computation rule for addition we have that $0+\succN(n)\jdeq\succN(0+n)$. Now we use the action on paths of $\succN:\N\to\N$ to obtain
  \begin{equation*}
    \ap{\succN}{p}:\succN(0+n)=\succN(n).
  \end{equation*}
  The left unit law is therefore defined by
  \begin{equation*}
    \leftunitlawaddN(n)\defeq\indN(\refl{0},\lam{p}\ap{\succN}{p}).\qedhere
  \end{equation*}
\end{proof}

\begin{prp}\label{prp:successor-laws-add-N}
  For any natural numbers $m$ and $n$, there are identifications
  \begin{align*}
    \leftsuccessorlawaddN(m,n) & : \succN(m)+n=\succN(m+n) \\
    \rightsuccessorlawaddN(m,n) & : m+\succN(n)=\succN(m+n).
  \end{align*}
\end{prp}

\begin{proof}
  We can define
  \begin{equation*}
    \rightsuccessorlawaddN(m,n)\defeq\refl{\succN(m+n)}
  \end{equation*}
  because we have a judgmental equality $m+\succN(n)\jdeq\succN(m+n)$ by the computation rules for $\addN$.

  The left successor law is constructed by induction on $n$. In the base case we have to construct an identification $\succN(m)+0=\succN(m+0)$, which is obtained by reflexivity. For the inductive step, assume that we have an identification $p:\succN(m)+n=\succN(m+n)$. Our goal is to show that
  \begin{equation*}
    \succN(m)+\succN(n)=\succN(m+\succN(n)). 
  \end{equation*}
  Note that we have the judgmental equalities
  \begin{align*}
    \succN(m)+\succN(n) & \jdeq\succN(\succN(m)+n) \\
    \succN(m+\succN(n)) & \jdeq\succN(\succN(m+n))
  \end{align*}
  Therefore it suffices to construct an identification
  \begin{equation*}
    \succN(\succN(m)+n)=\succN(\succN(m+n)).
  \end{equation*}
  Such an identification is given by $\ap{\succN}{p}$.
\end{proof}

\begin{prp}
  Addition on the natural numbers is associative, i.e., for any three natural numbers $m$, $n$, and $k$, there is an identification
  \begin{equation*}
    \associativeaddN(m,n,k):(m+n)+k=m+(n+k).
  \end{equation*}
\end{prp}

\begin{proof}
  We construct $\associativeaddN(m,n,k)$ by induction on $k$. In the base case we have the judgmental equalities
  \begin{equation*}
    (m+n)+0\jdeq m+n\jdeq m+(n+0).
  \end{equation*}
  Therefore we define $\associativeaddN(m,n,0)\defeq\refl{m+n}$.

  For the inductive step, let $p:(m+n)+k=m+(n+k)$. Our goal is to show that
  \begin{equation*}
    (m+n)+\succN(k)=m+(n+\succN(k)).
  \end{equation*}
  Note that we have the judgmental equalities
  \begin{align*}
    (m+n)+\succN(k) & \jdeq \succN((m+n)+k) \\
    m+(n+\succN(k)) & \jdeq m+(\succN(n+k)) \\
                    & \jdeq \succN(m+(n+k))
  \end{align*}
  Therefore it suffices to construct an identification
  \begin{equation*}
    \succN((m+n)+k)=\succN(m+(n+k)),
  \end{equation*}
  which we have by $\ap{\succN}{p}$.
\end{proof}

\begin{prp}
  Addition on the natural numbers is commutative, i.e., for any two natural numbers $m$ and $n$ there is an identification
  \begin{equation*}
    \commutativeaddN(m,n) : m+n=n+m.
  \end{equation*}
\end{prp}

\begin{proof}
  We construct $\commutativeaddN(m,n)$ by induction on $m$. In the base case we have to show that $0+n=n+0$, which holds by the unit laws for $n$, proven in \cref{prp:unit-laws-add-N}.

  For the inductive step, let $p:m+n=n+m$. Our goal is to construct an identification $\succN(m)+n=n+\succN(m)$. Now it is clear why we first proved the successor laws: we compute
  \begin{align*}
    \succN(m)+n & = \succN(m+n) \\
                & = \succN(n+m) \\
                & \jdeq n+\succN(m).
  \end{align*}
  The first identification is obtained by \cref{prp:successor-laws-add-N}, and the second identification is the identification $\ap{\succN}{p}$.
\end{proof}

\begin{exercises}
  \exitem \label{ex:inv_assoc}Show that the operation inverting identifications distributes over the concatenation operation, i.e., construct an identification
  \index{distributivity!of inv over concat@{of $\invfunc$ over $\concat$}}
  \index{identity type!distributive-inv-concat@{$\distributiveinvconcat$}|textbf}
  \begin{align*}
    \distributiveinvconcat(p,q):\id{(\ct{p}{q})^{-1}}{\ct{q^{-1}}{p^{-1}}}.
  \end{align*}
  for any $p:\id{x}{y}$ and $q:\id{y}{z}$.
  \exitem \label{ex:inv_con}For any $p:x=y$, $q:y=z$, and $r:x=z$, construct maps
  \index{identity type!inv-con@{$\invcon$}|textbf}
  \index{inv-con@{$\invcon$}|textbf}
  \index{identity type!con-inv@{$\coninv$}|textbf}
  \index{con-inv@{$\coninv$}|textbf}
  \begin{align*}
    \invcon(p,q,r) & : (\ct{p}{q}=r)\to (q=\ct{p^{-1}}{r}) \\
    \coninv(p,q,r) & : (\ct{p}{q}=r)\to (p=\ct{r}{q^{-1}}).
  \end{align*}
  \exitem Let $B$ be a type family over $A$, and consider an identification $p:\id{a}{x}$ in $A$. Construct for any $b:B(a)$ an identification\index{lift@{$\lift$}|textbf}\index{identity type!lift@{$\lift$}|textbf}
  \begin{equation*}
    \lift_B(p,b) : \id{(a,b)}{(x,\tr_B(p,b))}.
  \end{equation*}
  In other words, an identification $p:x=y$ in the \emph{base type} $A$ \emph{lifts} to an identification in $\sm{x:A}B(x)$ for every element in $B(x)$, analogous to the path lifting property for fibrations in homotopy theory.
  \exitem Consider four consecutive identifications
  \begin{equation*}
    \begin{tikzcd}
      a \arrow[r,equals,"p"] & b \arrow[r,equals,"q"] & c \arrow[r,equals,"r"] & d \arrow[r,equals,"s"] & e
    \end{tikzcd}
  \end{equation*}
  in a type $A$. In this exercise we will show that the \define{Mac Lane pentagon}\index{Mac Lane pentagon|textbf}\index{identity type!Mac Lane pentagon|textbf} for identifications commutes.
  \begin{subexenum}
  \item Construct the five identifications $\alpha_1,\ldots,\alpha_5$ in the pentagon
    \begin{equation*}
      \begin{tikzcd}[column sep=-1.5em]
        &[-2em] \ct{(\ct{(\ct{p}{q})}{r})}{s} \arrow[rr,equals,"\alpha_4"] \arrow[dl,equals,swap,"\alpha_1"] & & \ct{(\ct{p}{q})}{(\ct{r}{s})} \arrow[dr,equals,"\alpha_5"] &[-2em] \\
        \ct{(\ct{p}{(\ct{q}{r})})}{s} \arrow[drr,equals,swap,"\alpha_2"] & & & & \ct{p}{(\ct{q}{(\ct{r}{s})})}, \\
        & & \ct{p}{(\ct{(\ct{q}{r})}{s})} \arrow[urr,equals,swap,"\alpha_3"]
      \end{tikzcd}
    \end{equation*}
    where $\alpha_1$, $\alpha_2$, and $\alpha_3$ run counter-clockwise, and $\alpha_4$ and $\alpha_5$ run clockwise.
  \item Show that
    \begin{equation*}
      \ct{(\ct{\alpha_1}{\alpha_2})}{\alpha_3} = \ct{\alpha_4}{\alpha_5}.
    \end{equation*}
  \end{subexenum}
  \exitem \label{ex:semi-ring-laws-N}In this exercise we show that the operations of addition and multiplication on the natural numbers satisfy the laws of a commutative \define{semi-ring}.%
  \index{semi-ring laws!for N@{for $\N$}}%
  \index{natural numbers!semi-ring laws}%
  \index{associativity!of multiplication on N@{of multiplication on $\N$}}%
  \index{unit laws!for multiplication on N@{for multiplication on $\N$}}%
  \index{commutativity!of multiplication on N@{of multiplication on $\N$}}%
  \index{distributivity!of mulN over addN@{of $\mulN$ over $\addN$}}%
  \index{natural numbers!associativity of multiplicatoin@{associativity of multiplication}}
  \index{natural numbers!unit laws for multiplication}
  \index{natural numbers!zero laws for multiplication}
  \index{natural numbers!commutativity of multiplication}
  \index{natural numbers!distributivity of multiplication over addition}
  \begin{subexenum}
  \item Show that multiplication satisfies the following laws:
    \begin{align*}
      m\cdot 0 & = 0 & m\cdot 1 & = m & m\cdot \succN(n) & = m+m\cdot n \\
      0\cdot m & = 0 & 1\cdot m & = m & \succN(m)\cdot n & = m\cdot n+n.
    \end{align*}
  \item Show that multiplication on $\N$ is commutative:
    \begin{equation*}
      m\cdot n=n\cdot m.
    \end{equation*}
  \item \label{ex:distributive-mul-addN}Show that multiplication on $\N$ distributes over addition from the left and from the right, i.e., show that we have identifications
    \begin{align*}
      m\cdot (n+k) & = m\cdot n + m\cdot k \\
      (m+n)\cdot k & = m\cdot k + n\cdot k.
    \end{align*}
  \item Show that multiplication on $\N$ is associative:
    \begin{align*}
      (m\cdot n)\cdot k & = m\cdot (n\cdot k).
    \end{align*}
  \end{subexenum}
  \exitem \label{ex:is-equiv-succ-Z}Show that
  \begin{equation*}
    \succZ(\predZ(k))=k \qquad\text{and}\qquad \predZ(\succZ(k))=k
  \end{equation*}
  for any $k:\Z$, where $\predZ$ is the predecessor function on the integers, defined in \cref{ex:int_pred}.
  \exitem \label{ex:int_group_laws}\index{integers!group laws} In this exercise we will show that the laws for abelian groups hold for addition on the integers, using the group operations on $\Z$ defined in \cref{ex:int_group_ops}.
  \begin{subexenum}
  \item Show that addition satisfies the left and right unit laws, i.e., show that\index{unit laws!for addition on Z@{for addition on $\Z$}}\index{integers!unit laws for addition}
    \begin{align*}
      0+x & = x \\
      x+0 & = x.
    \end{align*}
  \item Show that the following successor and predecessor laws hold for addition on $\Z$.
    \begin{align*}
      \predZ(x)+y & = \predZ(x+y) & \succZ(x)+y & = \succZ(x+y) \\
      x+\predZ(y) & = \predZ(x+y) & x+\succZ(y) & = \succZ(x+y).
    \end{align*}
  \item Use part (b) to show that addition on the integers is associative and commutative, show that\index{associativity!of addition on Z@{of addition on $\Z$}}\index{commutativity!of addition on Z@{of addition on $\Z$}}\index{integers!associativity of addition}\index{integers!commutativity of addition}
    \begin{align*}
      (x+y)+z & = x + (y+z) \\
      x+y & = y+x.
    \end{align*}
  \item Show that addition satisfies the left and right inverse laws:\index{inverse laws!for addition on Z@{for addition on $\Z$}}\index{integers!inverse laws for addition}
    \begin{align*}
      (-x)+x & =0 \\
      x+(-x) &=0.
    \end{align*}
  \end{subexenum}
  \exitem \label{ex:ring-Z}In this exercise we will show that $\Z$ satisfies the axioms of a \define{ring}\index{ring!integers}\index{integers!is a ring}, using the multiplication operation defined in \cref{ex:mulZ}.
  \begin{subexenum}
  \item Show that multiplication on $\Z$ satisfies the following laws for $0$ and $1$\index{zero laws!for mulZ@{for $\mulZ$}}\index{unit laws!for multiplication on Z@{for multiplication on $\Z$}}\index{mul Z@{$\mulZ$}!unit laws}\index{mul Z@{$\mulZ$}!zero laws}\index{integers!zero laws for multiplication}\index{unit laws for multiplication}:
    \begin{align*}
      0\cdot x & = 0 & 1\cdot x & = x \\
      x\cdot 0 & = 0 & x\cdot 1 & = x.
    \end{align*}
  \item Show that multiplication on $\Z$ satisfies the predecessor and successor laws\index{mul Z@{$\mulZ$}!predecessor laws}\index{mul Z@{$\mulZ$}!successor laws}\index{integers!successor laws for addition}\index{integers!predecessor laws for addition}:
    \begin{align*}
      \predZ(x)\cdot y & = x\cdot y-y & \succZ(x)\cdot y & = x\cdot y + y \\
      x\cdot \predZ(y) & = x\cdot y-x & y\cdot \succZ(y) & = x\cdot y + x.
    \end{align*}
  \item Show that multiplication on $\Z$ distributes over addition, both from the left and from the right\index{mul Z@{$\mulZ$}!distributive over addZ@{distributive over $\addZ$}}\index{distributivity!of mulZ over addZ@{of $\mulZ$ over $\addZ$}}\index{integers!distributivity of multiplication over addition}:
    \begin{align*}
      x\cdot(y+z) & = x\cdot y+ x\cdot z \\
      (x+y)\cdot z & = x\cdot z + y\cdot z.
    \end{align*}
  \item Show that multiplication on $Z$ is associative and commutative\index{associativity!of multiplication on Z@{of multiplication on $\Z$}}\index{mul Z@{$\mulZ$}!associativity}\index{commutativity!of multiplication on Z@{of multiplication on $\Z$}}\index{mul Z@{$\mulZ$}!commutativity}\index{integers!associativity of multiplication}\index{integers!commutativity of multiplication}:
    \begin{align*}
      (x\cdot y)\cdot z & = x\cdot (y\cdot z) \\
      x\cdot y & = y\cdot x.
    \end{align*}
  \end{subexenum}
\end{exercises}

\index{identity type|)}
\index{inductive type!identity type|)}
\index{inductive type|)}


\section{Universes}\label{sec:universes}

\index{universe|(}
\index{universal family|(}
To complete our specification of dependent type theory, we introduce type theoretic \emph{universes}. Universes can be thought of as types that consist of types. In reality, however, a universe consists of a type $\UU$ equipped with a type family $\Ty$ over $\UU$. For any $X:\UU$ we think of $X$ as an \emph{encoding}\index{encoding of a type in a universe} of the type $\Ty(X)$. The type family $\Ty$ is called a \emph{universal type family}.

There are several reasons to equip type theory with universes. One important reason is that it enables us to define new type families over inductive types, using their induction principle. We use this way of defining type families to define many familiar relations over $\N$, such as the ordering relations $\leq$ and $<$. We also introduce a relation $\EqN$ called the \emph{observational equality} on $\N$. This equivalence relation can be used to show that $\zeroN\neq\oneN$.

The idea of introducing an observational equality relation for a particular type is that it should help us thinking about the identity type. The identity type has been introduced in a very generic and uniform way. In specific cases, however, we have a clear idea of what the equality relation \emph{should be}. In the case of the natural numbers, for instance, we will use the observational equality $\EqN$ to characterize the identity type of $\N$. Characterizing identity types is one of the main themes in homotopy type theory.

A second reason to introduce universes is that it allows us to define many types of types equipped with structure. One of the most important examples is the type of groups\index{group}, which is the type of types equipped with the group operations satisfying the group laws, and for which the underlying type is a set\index{set}. We won't discuss the condition for a type to be a set until \cref{chap:hierarchy}, so the definition of groups in type theory will be given much later.

\subsection{Specification of type theoretic universes}

A universe consists of a type $\UU$ of which the elements can be thought of as `codes' for types. A universe also comes equipped with a type family $\Ty$ indexed by $\UU$. Given an element $X:\UU$, we think of the type $\Ty(X)$ as the type of elements of $X$. The family $\Ty$ is called the \define{universal type family}.

One of the distinguishing features of universes is that they are closed under all the type constructors. Given a universe $\UU$ with universal type family $\Ty$, how do we express that $\UU$ is closed under $\Sigma$-types, for example? Recall that a $\Sigma$-type is formed using a type $A$ and a type family $B$ over $A$. Thus, if $A$ is a type in $\UU$ and $B$ is a family of types in $\UU$, we would like to express that the $\Sigma$-type is also a type in $\UU$. However, we cannot just assert that $\sm{x:A}B(x)$ is an element of the universe, because type theory carefully distinguishes between types and elements.

We express that $\UU$ is closed under $\Sigma$-types using a new operation $\check{\Sigma}$, which takes two arguments. The first argument is an element $X:\UU$, and the second argument is a family of types in $\UU$ indexed by the elements of $X$, i.e., a map $\Ty(X)\to\UU$. Thus we say that $\UU$ is closed under $\Sigma$-types by asserting that $\UU$ comes equipped with an operation
\begin{equation*}
  \check{\Sigma} : \prd{X:\UU} (\Ty(X)\to\UU)\to\UU
\end{equation*}
Furthermore, we ask that the element $\check{\Sigma}(X,Y):\UU$ satisfies the judgmental equality
\begin{equation*}
  \Ty(\check{\Sigma}(X,Y))\jdeq\sm{x:\Ty(X)}\Ty(Y(x)).
\end{equation*}
This judgmental equality asserts that the element $\check{\Sigma}(X,Y)$ of the universe $\UU$ \emph{represents} the $\Sigma$-type $\sm{x:\Ty(X)}\Ty(Y(x))$.

We will similarly assume that universes are closed under $\Pi$-types and the other ways of forming types. However, there is an important restriction: it would be inconsistent to assume that the universe is contained in itself. One way of thinking about this is that universes are types of \emph{small} types, and it cannot be the case that the universe is small with respect to itself. In \cref{subsec:russell} we will use a variant of Russell's paradox to derive a contradiction when $\UU$ is assumed to be (equivalent to a type) in $\UU$. Instead of assuming that the universe contains itself, we will assume that there are plenty of universes: enough universes so that any type family can be obtained by substituting into the universal type family of some universe.

\begin{defn}\label{defn:universe}
  A \define{universe}\index{universe|textbf} in type theory is a type $\UU$\index{U@{$\UU$}|see {universe}}\index{V@{$\VV$}|see {universe}} in the empty context, equipped with a type family $\Ty$\index{T@{$\Ty$}|see {universal family}} over $\UU$ called a \define{universal family}\index{type family!universal family|textbf}\index{universal family|textbf}, that is closed under the type forming operations in the sense that it comes equipped with the following structure:
  \begin{enumerate}
  \item $\UU$ is closed under $\Pi$, in the sense that it comes equipped with a function
    \begin{equation*}
      \check{\Pi} :\prd{X:\UU}(\Ty(X)\to\UU)\to\UU
    \end{equation*}
    for which the judgmental equality
    \begin{equation*}
      \Ty\big(\check{\Pi}(X,Y)\big)\jdeq \prd{x:\Ty(X)}\Ty(Y(x)).
    \end{equation*}
    holds, for every $X:\UU$ and $Y:\Ty(X)\to\UU$.
  \item $\UU$ is closed under $\Sigma$ in the sense that it comes equipped with a function
    \begin{equation*}
      \check{\Sigma} :\prd{X:\UU}(\Ty(X)\to\UU)\to\UU
    \end{equation*}
    for which the judgmental equality
    \begin{equation*}
      \Ty\big(\check{\Sigma}(X,Y)\big) \jdeq \sm{x:\Ty(X)}\Ty(Y(x))
    \end{equation*}
    holds, for every $X:\UU$ and $Y:\Ty(X)\to\UU$.
  \item $\UU$ is closed under identity types, in the sense that it comes equipped with a function
    \begin{equation*}
      \check{\mathrm{I}} : \prd{X:\UU}\Ty(X)\to(\Ty(X)\to\UU)
    \end{equation*}
    for which the judgmental equality
    \begin{equation*}
      \Ty\big(\check{\mathrm{I}}(X,x,y)\big)\jdeq (\id{x}{y})
    \end{equation*}
    holds, for every $X:\UU$ and $x,y:\Ty(X)$.
  \item $\UU$ is closed under coproducts, in the sense that it comes equipped with a function
    \begin{equation*}
      \mathbin{\check{+}}:\UU\to (\UU\to\UU)
    \end{equation*}
    that satisfies $\Ty\big(X\mathbin{\check{+}}Y\big)\jdeq \Ty(X)+\Ty(Y)$.
  \item $\UU$ contains elements $\check{\emptyt},\check{\unit},\check{\N}:\UU$
    that satisfy the judgmental equalities
    \begin{align*}
      \Ty(\check{\emptyt}) & \jdeq \emptyt \\
      \Ty(\check{\unit}) & \jdeq \unit \\
      \Ty(\check{\N}) & \jdeq \N.
    \end{align*}
  \end{enumerate}
  Consider a universe $\UU$ and a type $A$ in context $\Gamma$. We say that $A$ is a type in $\UU$, or that $\UU$ \define{contains} $A$, if $\UU$ comes equipped with an element $\check{A}:\UU$ in context $\Gamma$, for which the judgment
  \begin{equation*}
    \Gamma\vdash\Ty\big(\check{A}\big)\jdeq A~\type
  \end{equation*}
  holds. If $A$ is a type in $\UU$, we usually write simply $A$ for $\check{A}$ and also $A$ for $\Ty(\check{A})$.
\end{defn}

\begin{rmk}
  Since ordinary function types are defined as a special case of dependent function types, we don't have to assume separately that universes are closed under ordinary function types. Similarly, it follows from the assumption that universes are closed under dependent pair types that universes are closed under cartesian product types.
\end{rmk}

\subsection{Assuming enough universes}

\index{enough universes|(}
\index{universe!enough universes|(}
  Most of the time we will get by with assuming one universe $\UU$, and indeed we recommend on a first reading of this text to simply assume that there is one universe $\UU$. However, sometimes we might want to consider the universe $\UU$ itself to be a type in some universe. In such situations we cannot get by with a single universe, because the assumption that $\UU$ is a element of itself would lead to inconsistencies like the Russell's paradox.

  Russell's paradox is the famous argument that there cannot be a set of all sets. If there were such a set $S$, then we could consider Russell's subset
  \begin{equation*}
    R:=\{x\in S\mid x\notin x\}.
  \end{equation*}
  Russell then observed that $R\in R$ if and only if $R\notin R$, so we reach a contradiction. A variant of this argument reaches a similar contradiction when we assume that $\UU$ is a universe that contains an element $\check{\UU}:\UU$ such that $\mathcal{T}\big(\check{\UU}\big)\jdeq \UU$. In order to avoid such paradoxes, Russell and Whitehead formulated the \emph{ramified theory of types} in their book \emph{Principia Mathematica}. The ramified theory of types is a precursor of Martin L\"of's type theory that we are studying in this book.  

  Even though the universe is not an element of itself, it is still convenient if every type, including any universe, is in \emph{some} universe. Therefore we will assume that there are sufficiently many universes:
  
  \begin{postulate}\label{enough-universes}
    We assume that there are \define{enough universes}\index{universe!enough universes|textbf}\index{enough universes|textbf}, i.e., that for every finite list of types in context
    \begin{equation*}
      \Gamma_1\vdash A_1~\type\qquad\cdots\qquad\Gamma_n\vdash A_n~\type,
    \end{equation*}
    there is a universe $\UU$ that contains each $A_i$ in the sense that $\UU$ comes equipped with
    \begin{align*}
      \Gamma_i\vdash \check{A}_i:\UU
    \end{align*}
    for which the judgment
    \begin{equation*}
      \Gamma_i\vdash \Ty\big(\check{A}_i\big)\jdeq A_i~\type
    \end{equation*}
    holds.
  \end{postulate}

With this assumption it will rarely be necessary to work with more than one universe at the same time. Using the assumption that for any finite list of types in context there is a universe that contains those types, we obtain many specific universes.

\begin{defn}
  The \define{base universe}\index{base universe|textbf}\index{universe!base universe|textbf} $\UU_0$\index{U 0@{$\UU_0$}|see {base universe}} is the universe that we obtain using \cref{enough-universes} with the empty list of types in context.
\end{defn}

In other words, the base universe is a universe that is closed under all the ways of forming types, but it isn't specified to contain any further types.

\begin{defn}\label{defn:successor-universe}
  The \define{successor universe}\index{successor universe|textbf}\index{universe!successor universe|textbf} of a universe $\UU$ is the universe $\UU^+$\index{U +@{$\UU^+$}|see {successor universe}} obtained using \cref{enough-universes} with the finite list
  \begin{align*}
    & \vdash \UU~\type \\
    X:\UU & \vdash \mathcal{T}(X)~\type.
  \end{align*}
\end{defn}

\begin{rmk}\label{rmk:successor-universe}
  The successor universe $\UU^+$ of $\UU$ therefore contains the type $\UU$ as well as every type in $\UU$, in the following sense
  \begin{align*}
    & \vdash \check{\UU}:\UU^+ & & \vdash \mathcal{T}^+(\check{\UU})\jdeq\UU~\type \\
    X:\UU & \vdash \check{\mathcal{T}}(X) :\UU^+ & X:\UU & \vdash \mathcal{T}^+(\check{\mathcal{T}}(X))\jdeq \mathcal{T}(X)~\type.
  \end{align*}
  In particular, we obtain a function $i:\UU\to\UU^+$ that includes the types in $\UU$ into $\UU^+$, given by
  \begin{equation*}
    i\defeq \lam{X}\check{\mathcal{T}}(X).
  \end{equation*}

  Using successor universes we can create an infinite tower
  \begin{equation*}
    \UU,\ \UU^+,\ \UU^{++},\ \ldots
  \end{equation*}
  of universes, starting at any universe $\UU$, in which each universe is contained in the next. However, such towers of universes need not be exhaustive in the sense that it might not be the case that every type is contained in a universe in this tower.
\end{rmk}

\begin{defn}\label{defn:join-universe}
  The \define{join}\index{join of universes|textbf}\index{universe!join of universes|textbf} of two universes $\UU$ and $\VV$ is the universe $\UU\sqcup\VV$\index{U u V@{$U\sqcup V$}|see {join of universes}} that we obtain using \cref{enough-universes} with the two types
  \begin{align*}
    X:\UU & \vdash \mathcal{T}_{\UU}(X)~\type \\
    Y:\VV & \vdash \mathcal{T}_{\VV}(Y) ~\type.
  \end{align*}
\end{defn}

\begin{rmk}\label{rmk:join-universe}
  Since the join $\UU\sqcup\VV$ contains all the types in $\UU$ and $\VV$, there are maps
  \begin{align*}
    i : \UU\to\UU\sqcup\VV \\
    j : \VV\to\UU\sqcup\VV
  \end{align*}
  Note that we don't postulate any relations between the universes. In general it will therefore be the case that the universes $(\UU\sqcup\VV)\sqcup\mathcal{W}$ and $\UU\sqcup(\VV\sqcup\mathcal{W})$ will be unrelated.
\end{rmk}
\index{enough universes|)}
\index{universe!enough universes|)}

\subsection{Observational equality of the natural numbers}
\index{natural numbers!observational equality|(}
\index{observational equality!on N@{on $\N$}|(}

Using universes, we can define many relations on the natural numbers. We give here the example of \emph{observational equality} of $\N$. The idea of observational equality is that, if we want to prove that $m$ and $n$ are observationally equal, we may do so by looking at $m$ and $n$:
\begin{enumerate}
\item If both $m$ and $n$ are $\zeroN$, then they are observationally equal.
\item If one of them is $\zeroN$ and the other is a successor, then they are not observationally equal.
\item If both $m$ and $n$ are successors, say $m\jdeq\succN(m')$ and $n\jdeq \succN(n')$, then $m$ and $n$ are observationally equal if and only if their predecessors $m'$ and $n'$ are observationally equal.
\end{enumerate}
Thus, observational equality is an inductively defined relation, which gives us an algorithm for checking equality on $\N$. Indeed, it can be used to show that equality of natural numbers is \emph{decidable}, i.e., there is a program that decides for any two natural numbers $m$ and $n$ whether they are equal or not.

\begin{defn}\label{defn:obs_nat}
We define the \define{observational equality}\index{observational equality!on N@{on $\N$}|textbf}\index{natural numbers!observational equality|textbf} of $\N$ as binary relation $\EqN:\N\to(\N\to\UU_0)$\index{Eq N@{$\EqN$}|textbf}\index{natural numbers!Eq N@{$\EqN$}|textbf} satisfying
\begin{align*}
\EqN(\zeroN,\zeroN) & \jdeq \unit & \EqN(\succN(n),\zeroN) & \jdeq \emptyt \\
\EqN(\zeroN,\succN(n)) & \jdeq \emptyt & \EqN(\succN(n),\succN(m)) & \jdeq \EqN(n,m).
\end{align*}
\end{defn}

\begin{constr}
We define $\EqN$ by double induction on $\N$. By the first application of induction it suffices to provide
\begin{align*}
E_0 & : \N\to\UU_0 \\
E_S & : \N\to ((\N\to\UU_0)\to(\N\to\UU_0))
\end{align*}
We define $E_0$ by induction, taking $E_{00}\defeq \unit$ and $E_{0S}(n,X,m)\defeq \emptyt$. The resulting family $E_0$ satisfies
\begin{align*}
E_0(\zeroN) & \jdeq \unit \\
E_0(\succN(n)) & \jdeq \emptyt.
\end{align*} 
We define $E_S$ by induction, taking $E_{S0}\defeq \emptyt$ and $E_{SS}(n,X,m)\defeq X(m)$. The resulting family $E_S$ satisfies
\begin{align*}
E_S(n,X,\zeroN) & \jdeq \emptyt \\
E_S(n,X,\succN(m)) & \jdeq X(m) 
\end{align*}
Therefore we have by the computation rule for the first induction that the judgmental equality
\begin{align*}
\EqN(\zeroN,m) & \jdeq E_0(m) \\
\EqN(\succN(n),m) & \jdeq E_S(n,\EqN(n),m)
\end{align*}
holds, from which the judgmental equalities in the statement of the definition follow.
\end{constr}

The observational equality of the natural numbers is important because it can be used to prove equalities and negations of equalities. \cref{prp:Eq-eq-N} enables us to do so.

\begin{lem}
  Observational equality of $\N$ is a reflexive relation, i.e., we have
  \begin{equation*}
    \reflEqN : \prd{n:\N}\EqN(n,n).
  \end{equation*}
\end{lem}

\begin{proof}
  The function $\reflEqN$ is defined by induction on $n$, taking
  \begin{align*}
    \reflEqN(\zeroN) & \defeq \ttt \\
    \reflEqN(\succN(n)) & \defeq \reflEqN(n).\qedhere
  \end{align*}
\end{proof}

\begin{prp}\label{prp:Eq-eq-N}
  For any two natural numbers $m$ and $n$, we have
  \begin{equation*}
    (m=n)\leftrightarrow \EqN(m,n).
  \end{equation*}
\end{prp}

\begin{proof}
  The function $(m=n)\to\EqN(m,n)$ is defined by the induction principle of identity types, using the reflexivity of $\EqN$.

  The converse $\EqN(m,n)\to (m=n)$ is defined by induction on $m$ and $n$. If both $m$ and $n$ are zero, we have $\refl{\zeroN}:\zeroN=\zeroN$. If one of $m$ and $n$ is zero and the other is a successor, then $\EqN(m,n)$ is empty and we have a function $\emptyt\to (m=n)$ by the induction principle of the empty type. In the inductive step, suppose we have a function $f:\EqN(m,n)\to (m=n)$. Then we can define a function
  \begin{equation*}
    \EqN(\succN(m),\succN(n))\to (\succN(m)=\succN(n))
  \end{equation*}
  as the composite
  \begin{equation*}
    \begin{tikzcd}
      \EqN(\succN(m),\succN(n)) \arrow[r,dashed] \arrow[d,swap,"\idfunc"] & (\succN(m)=\succN(n)). \\
      \EqN(m,n) \arrow[r,swap,"f"] & (m=n) \arrow[u,swap,"\apfunc{\succN}"]
    \end{tikzcd}
  \end{equation*}
  Note that the map on the left is the identity function, because we have the judgmental equality $\EqN(\succN(m),\succN(n))\jdeq\EqN(m,n)$ by definition of $\EqN$.
\end{proof}
\index{natural numbers!observational equality|)}
\index{observational equality!on N@{on $\N$}|)}

\subsection{Peano's seventh and eighth axioms}\label{sec:peano-axioms}
Using the observational equality of $\N$, we can prove Peano's seventh and eighth axioms. In his \emph{Arithmetices Principia} \cite{Peano}, the natural numbers are based at $1$, but today it is customary to have the natural numbers based at $0$. Adapting for this, the seventh and eighth axioms assert that
\begin{enumerate}
\item[(P7)] For any two natural numbers $m$ and $n$, we have
  \begin{equation*}
    (m=n)\leftrightarrow (\succN(m)=\succN(n)).
  \end{equation*}
\item[(P8)] For any natural number $n$, we have $\zeroN\neq\succN(n)$.
\end{enumerate}

\begin{thm}\label{thm:is-injective-succ-N}
  For any two natural numbers $m$ and $n$, we have\index{is injective!succN@{$\succN$}}\index{succN@{$\succN$}!is injective}
  \begin{equation*}
    (m=n)\leftrightarrow (\succN(m)=\succN(n)).
  \end{equation*}
\end{thm}

\begin{proof}
  The forward implication is given by the action on paths of the successor function
  \begin{equation*}
    \apfunc{\succN}:(m=n)\to(\succN(m)=\succN(n)).
  \end{equation*}
  The direction of interest is the converse, which asserts that the successor function is injective.

  Here we use \cref{prp:Eq-eq-N}, which asserts that $(m=n)\leftrightarrow \EqN(m,n)$ for all $m,n:\N$. Furthermore, we have $\EqN(\succN(m),\succN(n))\jdeq \EqN(m,n)$. Therefore, we obtain
  \begin{equation*}
    \begin{tikzcd}
      (\succN(m)=\succN(n)) \arrow[r,dashed] \arrow[d] & (m=n) \\
      \EqN(\succN(m),\succN(n)) \arrow[r,swap,"\idfunc"] & \EqN(m,n), \arrow[u]
    \end{tikzcd}
  \end{equation*}
  and we define the function $(\succN(m)=\succN(n))\to(m=n)$ as the composite of the maps going down, then right, and then up.
\end{proof}

\begin{thm}\label{prp:zero-one}
  For any natural number $n$, we have $\zeroN\neq\succN(n)$.\index{natural numbers!zero is not a successor}
\end{thm}

\begin{proof}
  By \cref{prp:Eq-eq-N} it follows that there is a family of maps
  \begin{equation*}
    (\zeroN=n)\to \EqN(\zeroN,n).
  \end{equation*}
  indexed by $n:\N$. Since $\EqN(\zeroN,\succN(n))\jdeq\emptyt$ it follows that
  \begin{equation*}
    (\zeroN=\succN(n))\to \emptyt,
  \end{equation*}
  which is precisely the claim.
\end{proof}

\begin{exercises}
  \exitem
  \begin{subexenum}
  \item Show that
    \begin{align*}
      (m=n) & \leftrightarrow (m+k=n+k) \\*
      (m=n) & \leftrightarrow (m\cdot(k+1)=n\cdot(k+1))
    \end{align*}
    for all $m,n,k:\N$. In other words, adding $k$ and multiplying by $k+1$ are injective functions.
  \item \label{ex:is-zero-summand-is-zero-sum-N}Show that
    \begin{align*}
      (m+n=0) & \leftrightarrow (m=0)\times (n=0)\\*
      (mn=0) & \leftrightarrow (m=0)+(n=0)\\*
      (mn=1) & \leftrightarrow (m=1)\times (n=1)
    \end{align*}
    for all $m,n:\N$.
  \item Show that
    \begin{align*}
      m & \neq m+(n+1) \\*
      m+1 & \neq (m+1)(n+2)
    \end{align*}
    for all $m,n:\N$.
  \end{subexenum}
  \exitem \label{ex:obs_bool}\index{observational equality!on bool@{on $\bool$}|textbf}\index{booleans!observational equality|textbf}
  \begin{subexenum}
  \item Define observational equality $\EqBool$\index{Eq-bool@{$\EqBool$}|textbf} by induction on the booleans.
  \item Show that
    \begin{equation*}
      (x=y)\leftrightarrow \EqBool(x,y)
    \end{equation*}
    for any $x,y:\bool$.
  \item \label{ex:zero-neq-one-bool}Show that $b\neq\negbool(b)$ for any $b:\bool$. Conclude that $\bfalse\neq\btrue$. \index{false neq true@{$\bfalse\neq\btrue$}}\index{booleans!false neq true@{$\bfalse\neq\btrue$}}
  \end{subexenum}
  \exitem \label{ex:order_N}The ordering relation\index{relation!order}\index{order relation!leq on N@{$\leq$ on $\N$}}\index{leq@{$\leq$}!on N@{on $\N$}} $\leq$ on $\N$ is defined recursively by
    \begin{align*}
      (\zeroN\leq\zeroN) & \defeq \unit & (\zeroN\leq n+1) & \defeq \unit \\
      (m+1\leq\zeroN) & \defeq \emptyt & (m+1\leq n+1) & \defeq (m\leq n).
    \end{align*}
  \begin{subexenum}
  \item Show that $\leq$ satisfies the axioms of a \emph{poset}, i.e., show that $\leq$ is
    \begin{enumerate}
    \item reflexive,\index{reflexive!leq on N@{$\leq$ on $\N$}}
    \item antisymmetric, and
    \item transitive.
    \end{enumerate}
  \item Show that
    \begin{equation*}
      (m\leq n)+(n\leq m)
    \end{equation*}
    for any $m,n:\N$.
  \item Show that
    \begin{align*}
      (m \leq n) & \leftrightarrow (m+k \leq n+k)
    \end{align*}
    holds for any $m,n,k:\N$.
  \item Show that
    \begin{align*}
      (m \leq n) & \leftrightarrow (m\cdot(k+1) \leq n\cdot(k+1))
    \end{align*}
    holds for any $m,n,k:\N$.
  \item Show that $k\leq \minN(m,n)$ holds if and only if both $k\leq m$ and $k\leq n$ hold, and show that $\maxN(m,n)\leq k$ holds if and only if both $m\leq k$ and $n\leq k$ hold.
  \end{subexenum}
  \exitem \label{ex:strict-order-N}The strict ordering relation\index{relation!strict order}\index{order relation!le on N@{$<$ on $\N$}}\index{le@{$<$}!on N@{on $\N$}} $<$ on $\N$ is defined recursively by
    \begin{align*}
      (\zeroN<\zeroN) & \defeq \emptyt & (\zeroN< n+1) & \defeq \unit \\
      (m+1<\zeroN) & \defeq \emptyt & (m+1< n+1) & \defeq (m< n).
    \end{align*}
  \begin{subexenum}
  \item Show that the strict ordering relation is
    \begin{enumerate}
    \item antireflexive,
    \item antisymmetric, and
    \item transitive.
    \end{enumerate}
  \item Show that $n<n+1$ and
    \begin{equation*}
      (m<n)\to (m<n+1)
    \end{equation*}
    for any $m,n:\N$.
  \item \label{ex:contradiction-le}Show that
    \begin{align*}
      (m<n) & \leftrightarrow (m+1\leq n) \\
      (m<n) & \leftrightarrow (n \nleq m)
    \end{align*}
    for any $m,n :\N$.
  \end{subexenum}
  \exitem \label{ex:distN}The \define{distance function}\index{distance function on N@{distance function on $\N$}|textbf}\index{dist N x y@{$\distN(x,y)$}|textbf}\index{natural numbers!distance function}
  \begin{equation*}
    \distN : \N \to (\N \to \N)
  \end{equation*}
  is defined recursively by
  \begin{align*}
    \distN(0,0) & \defeq 0 & \distN(0,n+1) & \defeq n+1 \\
    \distN(m+1,0) & \defeq m+1 & \distN(m+1,n+1) & \defeq \distN(m,n).
  \end{align*}
  In other words, the distance between two natural numbers is the \emph{symmetric difference} between them.
  \begin{subexenum}
  \item \label{ex:is-metric-distN}Show that $\distN$ satisfies the axioms of a metric:
    \begin{enumerate}
    \item $(m=n)\leftrightarrow (\distN(m,n)=0)$,
    \item $\distN(m,n) = \distN(n,m)$,
    \item $\distN(m,n) \leq\distN(m,k)+\distN(k,n)$.
    \end{enumerate}
  \item \label{ex:distN-triangle-equality}Show that $\distN(m,n)=\distN(m,k)+\distN(k,n)$ if and only if either both $m\leq k$ and $k\leq n$ hold or both $n\leq k$ and $k\leq m$ hold.
  \item \label{ex:translation-invariant-distN}Show that $\distN$ is translation invariant and linear:
    \begin{align*}
      \distN(a+m,a+n) & = \distN(m,n),\\
      \distN(k\cdot m,k\cdot n) & =k\cdot\distN(m,n).
    \end{align*}
  \item Show that $x+\distN(x,y)=y$ for any $x\leq y$. 
  \end{subexenum}
  \exitem Construct the \define{absolute value function}\index{absolute value function on Z@{absolute value function on $\Z$}|textbf}\index{integers!absolute value function|textbf}
  \begin{equation*}
    |\blank|:\Z\to\N
  \end{equation*}
  and show that it satisfies the following three properties:
  \begin{enumerate}
  \item $(x=0)\leftrightarrow (|x|=0)$,
  \item $|x+y|\leq |x|+|y|$,
  \item $|xy|=|x||y|$.
  \end{enumerate}
\end{exercises}
\index{universe|)}
\index{universal family|)}


\section{Modular arithmetic via the Curry-Howard interpretation}\label{sec:modular-arithmetic}

We have now fully described Martin-L\"of's dependent type theory. It is now up to us to start developing some mathematics in it, and Martin-L\"of's dependent type theory is great for elementary mathematics, such as basic number theory, some algebra, and combinatorics. The fundamental idea that is used to develop basic mathematics in type theory is the Curry-Howard interpretation. This is a translation of logic into type theory, which we will use to express concepts of mathematics such as divisibility, the congruence relations, and so on.

We will also introduce the family $\Fin{}$ of the standard finite types, indexed by $\N$, and show how each $\Fin{k+1}$ can be equipped with the group structure of integers modulo $k+1$. Our goal here is to demonstrate how to do those things in type theory, so we will aim for a high degree of accuracy.

\subsection{The Curry-Howard interpretation}\label{sec:Curry-Howard}

The \emph{Curry-Howard interpretation} is an interpretation of logic into type theory. Recall that in type theory there is no separation between the logical framework and the general theory of collections of mathematical objects the way there is in the more traditional setup with Zermelo-Fraenkel set theory, which is postulated by axioms in first order logic. These two aspects of the foundations of mathematics are unified in type theory. The idea of the Curry-Howard interpretation is therefore to express propositions as types, and to think of the elements of those types as their proofs. We illustrate this idea with an example.

\begin{eg}
  A natural number $d$ is said to divide a natural number $n$ if there exists a natural number $k$ such that $d\cdot k=n$. To represent the divisibility predicate in type theory, we need to define a \emph{type}
  \begin{equation*}
    d\mid n,
  \end{equation*}
  of which the elements are witnesses that $d$ divides $n$. In other words, $d\mid n$ should be the type that consists of natural numbers $k$ equipped with an identification $d\cdot k=n$. In general, the type of $x:A$ equipped with $y:B(x)$ is represented as the type $\sm{x:A}B(x)$. The interpretation of the existential quantification ($\exists$) into type theory via the Curry-Howard interpretation is therefore using $\Sigma$-types.
\end{eg}

\begin{defn}
  Consider two natural numbers $d$ and $n$. We say that $d$ \define{divides}\index{divisibility on N@{divisibility on $\N$}|textbf}\index{natural numbers!divisibility|textbf} $n$ if there is a element of type\index{d {"|" n}@{$d\mid n$}|textbf}\index{d {"|" n}@{$d\mid n$}|see{divisibility on $\N$}}
  \begin{equation*}
    d\mid n\defeq \sm{k:\N}d\cdot k=n.
  \end{equation*}
\end{defn}

\begin{rmk}
  This type-theoretical definition of the divisibility relation using $\Sigma$-types has two important consequences:
  \begin{enumerate}
  \item The principal way to show that $d\mid n$ holds is to construct a pair $(k,p)$ consisting of a natural number $k$ and an identification $p:d\cdot k=n$.
  \item The principal way to use a hypothesis $H:d\mid n$ in a proof is to proceed by $\Sigma$-induction on the variable $H$. We then get to assume a natural number $k$ and an identification $p:d\cdot k=n$, in order to proceed with the proof.
  \end{enumerate}
\end{rmk}

\begin{eg}\label{rmk:elementary-facts-div}
  Just as existential quantification ($\exists$) is translated via the Curry-Howard interpretation to $\Sigma$-types, the translation of the universal quantification ($\forall$) in type theory via the Curry-Howard interpretation is to $\Pi$-types. For example, the assertion that every natural number is divisible by $1$ is expressed in type theory as
  \begin{equation*}
    \prd{x:\N} 1\mid x.
  \end{equation*}
  In other words, in order to show that every number $x:\N$ is divisible by $1$ we need to construct a dependent function
  \begin{equation*}
    \lam{x}p(x):\prd{x:\N}1\mid x.
  \end{equation*}
  We do this by constructing an element
  \begin{equation*}
    p(x):\sm{k:\N}1\cdot k=x
  \end{equation*}
  indexed by $x:\N$. Such an element $p(x)$ is constructed as the pair $(x,q(x))$, where the identification $q(x):1\cdot x=x$ is obtained from the left unit law of multiplication on $\N$, which was constructed in \cref{ex:semi-ring-laws-N}.

  Similarly, the type theoretic proof that every natural number $k$ divides $0$, i.e., that $k\mid 0$, is the pair $(0,p)$ consisting of the natural number $0$ and the identification $p:k\cdot 0=0$ obtained from the right annihilation law of multiplication on $\N$. This identification was also constructed in \cref{ex:semi-ring-laws-N}.
\end{eg}

In the following proposition we will see examples of how a hypothesis of type $d\mid x$ can be used.

\begin{prp}\label{prp:div-3-for-2}
  Consider three natural numbers $d$, $x$ and $y$. If $d$ divides any two of the three numbers $x$, $y$, and $x+y$, then it also divides the third.
\end{prp}

\begin{proof}
  We will only show that if $d$ divides $x$ and $y$, then it divides $x+y$. The remaining two claims, that if $d$ divides $y$ and $x+y$ then it divides $x$, and that if $d$ divides $x$ and $x+y$ then it divides $y$, are left as \cref{ex:div-3-for-2}.

  Suppose that $d$ divides both $x$ and $y$. By assumption we have elements
  \begin{equation*}
    H:\sm{k:\N}d\cdot k=x,\qquad\text{and}\qquad K:\sm{k:\N}d\cdot k=y.
  \end{equation*}
  Since the types of the variables $H$ and $K$ are $\Sigma$-types, we proceed by $\Sigma$-induction on $H$ and $K$. Therefore we get to assume a natural number $k:\N$ equipped with an identification $p:d\cdot k=x$, and a natural number $l:\N$ equipped with an identification $q:d\cdot l=y$. Our goal is now to construct an identification
  \begin{equation*}
    d\cdot (k+l)=x+y.
  \end{equation*}
We construct such an identification as a concatenation $\ct{\alpha}{(\ct{\beta}{\gamma})}$, where the types of the identifications $\alpha$, $\beta$, and $\gamma$ are as follows:
  \begin{equation*}
    \begin{tikzcd}
      d\cdot(k+l) \arrow[r,equals,"\alpha"] & d\cdot k+d\cdot l \arrow[r,equals,"\beta"] & x+d\cdot l \arrow[r,equals,"\gamma"] & x+y.
    \end{tikzcd}
  \end{equation*}
  The identification $\alpha$ is obtained from the fact that multiplication on $\N$ distributes over addition, which was shown in \cref{ex:distributive-mul-addN}. The identifications $\beta$ and $\gamma$ are constructed using the action on paths of a function:
  \begin{equation*}
    \beta\defeq\ap{(\lam{t}t+d\cdot l)}{p},\qquad\text{and}\qquad \gamma\defeq \ap{(\lam{t}x+t)}{q}
  \end{equation*}
  To conclude the proof that $d\mid x+y$, note that we have constructed the pair
  \begin{equation*}
    (k+l,\ct{\alpha}{(\ct{\beta}{\gamma})}):\sm{k:\N}d\cdot k=x+y.\qedhere
  \end{equation*}
\end{proof}

The full Curry-Howard interpretation of logic into type theory also involves interpretations of disjunction, conjunction, implication, and equality.

The introduction and elimination rules for disjunction are, for instance,
\begin{equation*}
  \AxiomC{$P$}
  \UnaryInfC{$P\lor Q$}
  \DisplayProof
  \qquad
  \AxiomC{$Q$}
  \UnaryInfC{$P\lor Q$}
  \DisplayProof
  \qquad
  \text{and}
  \qquad
  \AxiomC{$P\Rightarrow R$}
  \AxiomC{$Q\Rightarrow R$}
  \BinaryInfC{$P\lor Q\Rightarrow R$}
  \DisplayProof
\end{equation*}
The two introduction rules assert that $P\lor Q$ holds provided that $P$ holds, and that $P\lor Q$ holds provided that $Q$ holds. These rules are analogous to the introduction rules for coproduct, which assert that there are functions $\inl : A\to A+B$ and $\inr : B \to A+B$. Furthermore, the non-dependent elimination principle for coproducts gives a function
\begin{equation*}
  (A\to C) \to ((B \to C) \to (A+B \to C))
\end{equation*}
for any type $C$, which is again analogous to the elimination rule of disjunction. The Curry-Howard interpretation of disjunction into type theory is therefore as coproducts.

To interpret conjunction into type theory we observe that the introduction rule and elimination rules for conjunction are
\begin{equation*}
  \AxiomC{$P$}
  \AxiomC{$Q$}
  \BinaryInfC{$P\land Q$}
  \DisplayProof
  \qquad
  \text{and}
  \qquad
  \AxiomC{$P\land Q$}
  \UnaryInfC{$P$}
  \DisplayProof
  \qquad
  \AxiomC{$P\land Q$}
  \UnaryInfC{$Q$}
  \DisplayProof
\end{equation*}
Product types possess such structure, where we have the pairing operation $\pair:A\to (B\to A\times B)$ and the projections $\proj 1:A\times B\to A$ and $\proj 2 : A\times B\to B$ give interpretations of the introduction and elimination rules for conjunction. The Curry-Howard interpretation of conjunction into type theory is therefore by products. We summarize the full Curry-Howard interpretation in \cref{table:Curry-Howard}.

\begin{table}[t]
  \begin{tabular}{ll}
    \toprule
    \multicolumn{2}{c}{The Curry-Howard interpretation} \\
    \midrule
    Propositions & Types \\
    Proofs & Elements \\
    Predicates & Type families \\
    $\top$ & $\unit$ \\
    $\bot$ & $\emptyt$ \\
    $P\lor Q$ & $A+B$ \\
    $P\land Q$ & $A\times B$ \\
    $P\Rightarrow Q$ & $A\to B$ \\
    $\neg P$ & $A\to \emptyt$ \\
    $\exists_{x}P(x)$ & $\sm{x:A}B(x)$ \\
    $\forall_{x}P(x)$ & $\prd{x:A}B(x)$ \\
    $x=y$ & $x=y$ \\
    \bottomrule
  \end{tabular}
  \caption{\label{table:Curry-Howard}The Curry-Howard interpretation of logic into type theory.}
\end{table}

\begin{rmk}
  We should note, however, that despite the similarities between logic and type theory that are highlighted in the Curry-Howard interpretation, there are also some differences. One important difference is that types may contain many elements, whereas in logic, propositions are usually considered to be \emph{proof irrelevant}. This means that to establish the truth of a proposition it only matters \emph{whether} it can be proven, not in how many different ways it can be proven. To address this dissimilarity between general types and logic, we will introduce in \cref{chap:uf} a more refined way of interpreting logic into type theory. In \cref{chap:hierarchy} we will define the type $\isprop(A)$, which expresses the property that the type $A$ is a proposition. Furthermore, we will introduce the \emph{propositional truncation} operation in \cref{sec:propositional-truncation}, which we will use to interpret logic into type theory in such a way that all logical assertions are interpreted as types that satisfy the condition of being a proposition.
\end{rmk}

\subsection{The congruence relations on \texorpdfstring{$\N$}{ℕ}}

Relations in the Curry-Howard interpretation of logic into type theory are also type valued. More specifically, a binary relation on a type $A$ is a family of types $R(x,y)$ indexed by $x,y:A$. Such relations are sometimes called \emph{typal}.

\begin{defn}
  Consider a type $A$. A \define{(typal) binary relation} on $A$ is defined to be a family of types $R(x,y)$ indexed by $x,y:A$. Given a binary relation $R$ on $A$, we say that $R$ is \define{reflexive} if it comes equipped with
  \begin{align*}
    \rho & : \prd{x:A}R(x,x), \\
    \intertext{we say that $R$ is \define{symmetric} if it comes equipped with}
    \sigma & : \prd{x,y:A} R(x,y)\to R(y,x), \\
    \intertext{and we say that $R$ is \define{transitive} if it comes equipped with}
    \tau & : \prd{x,y,z:A} R(x,y)\to (R(y,z)\to R(x,z)).
  \end{align*}
  A \define{(typal) equivalence relation} on $A$ is a reflexive, symmetric, and transitive binary typal relation on $A$.
\end{defn}

To define the congruence relation modulo $k$ in type theory using the Curry-Howard interpretation, we will define for any three natural numbers $x$, $y$, and $k$, a \emph{type}
\begin{equation*}
  x\equiv y\mod k
\end{equation*}
consisting of the proofs that $x$ is congruent to $y$ modulo $k$. We will define this type by directly interpreting Gauss' definition of the congruence relations in his \emph{Disquisitiones Arithmeticae} \cite{Gauss}: two numbers $x$ and $y$ are congruent modulo $k$ if $k$ divides the symmetric difference $\distN(x,y)$ between $x$ and $y$. Recall that $\distN(x,y)$ was defined in \cref{ex:distN} recursively by
  \begin{align*}
    \distN(0,0) & \defeq 0 & \distN(0,y+1) & \defeq y+1 \\
    \distN(x+1,0) & \defeq x+1 & \distN(x+1,y+1) & \defeq \distN(x,y).
  \end{align*}

\begin{defn}
  Consider three natural numbers $k,x,y:\N$. We say that $x$ is \define{congruent to $y$ modulo $k$}\index{congruence relations on N@{congruence relations on $\N$}|textbf}\index{natural numbers!congruence relations|textbf} if it comes equipped with an element of type
  \begin{equation*}
    x\equiv y \mod k \defeq k\mid\distN(x,y).
  \end{equation*}
\end{defn}

\begin{eg}
  For example, $k\equiv 0\mod k$. To see this, we have to show that $k\mid\distN(k,0)$. Since $\distN(k,0)=k$ it suffices to show that $k\mid k$. That is, we have to construct a natural number $l$ equipped with an identification $p:kl=k$. Of course, we choose $l\defeq 1$, and the equation $k1=k$ holds by the right unit law for multiplication on $\N$, which was shown in \cref{ex:semi-ring-laws-N}.
\end{eg}

\begin{prp}\label{prp:congruence-eqrel}
  For each $k:\N$, the congruence relation modulo $k$ is an equivalence relation.
\end{prp}

\begin{proof}
  Reflexivity follows from the fact that $\distN(x,x)=0$, and any number divides $0$. Symmetry follows from the fact that $\distN(x,y)=\distN(y,x)$ for any two natural numbers $x$ and $y$.

  The non-trivial part of the claim is therefore transitivity. Here we use the fact that for any three natural numbers $x$, $y$, and $z$, at least one of the equalities
  \begin{align*}
    \distN(x,y)+\distN(y,z) & =\distN(x,z) \\
    \distN(y,z)+\distN(x,z) & =\distN(x,y) \\
    \distN(x,z)+\distN(x,y) & =\distN(y,z)
  \end{align*}
  holds. A formal proof of this fact is given by case analysis on the six possible ways in which $x$, $y$, and $z$ can be ordered:
  \begin{align*}
    x\leq y & \text{ and }y\leq z, & x\leq z & \text{ and }z\leq y, \\
    y\leq z & \text{ and }z\leq x, & y\leq x & \text{ and }x\leq z, \\
    z\leq x & \text{ and }x\leq y, & z\leq y & \text{ and }y\leq x.
  \end{align*}
  Therefore it follows by \cref{ex:distN-triangle-equality} and \cref{prp:div-3-for-2} that ${k\mid\distN(x,z)}$ if ${k\mid\distN(x,y)}$ and ${k\mid\distN(y,z)}$.
\end{proof}

\subsection{The standard finite types}\label{sec:Fin}

The standard finite sets are classically defined as the sets $\{x\in\N\mid x<k\}$. This leads to the question of how to interpret a subset $\{x\in A\mid P(x)\}$ in type theory.

Since type theory is set up in such a way that elements come equipped with their types, subsets aren't formed the same way as in set theory, where the comprehension axiom is used to form the set $\{x\in A\mid P(x)\}$ for any predicate $P$ over $A$. The Curry-Howard interpretation dictates that predicates are interpreted as dependent types. Therefore, a set of elements $x\in A$ such that $P(x)$ holds is interpreted in type theory as the type of terms $x:A$ equipped with an element (a proof) $p:P(x)$. In other words, we interpret a subset $\{x\in A\mid P(x)\}$ as the type $\sm{x:A}P(x)$.

\begin{rmk}
  The alert reader may now have observed that the interpretation of a subset $\{x\in A\mid P(x)\}$ in type theory is the same as the interpretation of the proposition $\exists_{(x\in A)}P(x)$, while indeed the subset $\{x\in A\mid P(x)\}$ has a substantially different role in mathematics than the proposition $\exists_{(x\in A)}P(x)$. This points at a slight problem of the Curry-Howard interpretation of the existential quantifier. While the Curry-Howard interpretation of the existential quantifier is nevertheless useful and important, we will reinterpret the existential quantifier in type theory in \cref{sec:logic}.
\end{rmk}

Since subsets are interpreted as $\Sigma$-types, the `classical' definition of the standard finite types is
\begin{equation*}
  \classicalFin_k:=\sm{x:\N}x<k.
\end{equation*}
This is a perfectly fine definition of the standard finite types. However, the usual definition of the standard finite types in Martin-L\"of's dependent type theory is a more direct, recursive definition, which takes full advantage of the inductive constructions of dependent type theory. 

\begin{defn}\label{defn:fin}
  We define the type family $\Fin{}$ of the \define{standard finite types}\index{Fin k@{$\Fin{k}$}|see {standard finite type}}\index{Fin k@{$\Fin{k}$}|textbf}\index{standard finite type}\index{type family!of standard finite types} over $\N$ recursively by
  \begin{align*}
    \Fin{0} & \defeq \emptyt \\*
    \Fin{k+1} & \defeq \Fin{k}+\unit.
  \end{align*}
  We will write $i$ for the inclusion $\inl:\Fin{k}\to\Fin{k+1}$ and we will write $\ttt$ for the point $\inr(\ttt)$.
\end{defn}

In \cref{ex:classical-Fin} you will be asked to show that the types $\classicalFin_k$ and $\Fin{k}$ are isomorphic.

\begin{rmk}
The type family $\Fin{}$ over $\N$ can be given its own induction principle, which is, at least for the time being, the principal way to make constructions on $\Fin{k}$ for arbitrary $k:\N$ and to prove properties about those constructions. The induction principle of the standard finite types tells us that the family of standard finite types is inductively generated by
\begin{align*}
  i & : \Fin{k}\to\Fin{k+1} \\*
  \ttt & : \Fin{k+1}. 
\end{align*}
In other words, we can define a dependent function $f:\prd{k:\N}\prd{x:\Fin{k}}P_k(x)$ by defining
\begin{align*}
  g_k & : \prd{x:\Fin{k}}P_k(x)\to P_{k+1}(i(x)) \\*
  p_k & : P_{k+1}(\ttt)
\end{align*}
for each $k:\N$. The function $f$ defined in this way then satisfies the judgmental equalities
\begin{align*}
  f_{k+1}(i(x)) & \jdeq g_k(x,f_k(x)) \\*
  f_{k+1}(\ttt) & \jdeq p_k.
\end{align*}
These judgmental equalities completely determine the function $f$, and therefore we may also present such inductive definitions by pattern matching:
  \begin{align*}
    f_{k+1}(i(x)) & \defeq g_k(x,f_k(x)) \\*
    f_{k+1}(\ttt) & \defeq p_k.
  \end{align*}
\end{rmk}

We will often use definitions by pattern matching for two reasons: (i) such definitions are concise, and (ii) they display the judgmental equalities that hold for the defined object. Those judgmental equalities are the only thing we know about that object, and proving a claim about it often amounts to finding a way to apply these judgmental equalities.

To illustrate this way of working with the standard finite types, we define the inclusion functions $\Fin{k}\to\N$, and show that these are injective. In order to show that $\natFin_k$ is injective, we will also show that $\natFin_k$ is bounded.

\begin{defn}\label{defn:natFin}
  We define the inclusion $\natFin_k : \Fin{k}\to\N$ inductively by
  \begin{align*}
    \natFin_{k+1}(i(x)) & \defeq \natFin_{k}(x) \\
    \natFin_{k+1}(\ttt) & \defeq k.
  \end{align*}
\end{defn}

\begin{lem}\label{lem:is-bounded-natFin}
  The function $\natFin:\Fin{k}\to\N$ is bounded, in the sense that $\natFin(x)< k$ for each $x:\Fin{k}$.
\end{lem}

\begin{proof}
  The proof is by induction. In the base case there is nothing to show. In the inductive step, we have the inequalities $\natFin_{k+1}(i(x))\jdeq\natFin_{k}(x)<k<k+1$, where the first inequality holds by the inductive hypothesis, and we also have
  \begin{equation*}
    \natFin_{k+1}(\ttt)\jdeq k<k+1.\qedhere
  \end{equation*}
\end{proof}

\begin{prp}\label{prp:is-injective-natFin}
  The inclusion function $\natFin_k : \Fin{k}\to \N$ is injective, for each $k:\N$.
\end{prp}

\begin{proof}
  We define a function $\alpha_k(x,y):(\natFin_k(x)=\natFin_k(y))\to (x=y)$ recursively by
  \begin{align*}
    \alpha_{k+1}(i(x),i(y),p) & \defeq \ap{i}{\alpha_k(x,y,p)} & \alpha_{k+1}(i(x),\ttt,p) & \defeq \exfalso(f(p)) \\
    \alpha_{k+1}(\ttt,i(y),p) & \defeq \exfalso(g(p)) & \alpha_{k+1}(\ttt,\ttt,p) & \defeq \refl{},
  \end{align*}
  where $f:(\natFin_{k+1}(i(x))=\natFin_{k+1}(\ttt))\to\emptyt$ and $g:(\natFin_{k+1}(\ttt)=\natFin_{k+1}(i(y)))\to\emptyt$ are obtained from the fact that $\natFin_{k+1}(i(z))\jdeq\natFin_k(z)<k$ for any $z:\Fin{k}$, and the fact that $\natFin_{k+1}(\ttt)\jdeq k$.
\end{proof}

\subsection{The natural numbers modulo \texorpdfstring{$k+1$}{k+1}}\label{subsec:finite-types-quotient-maps}

Given an equivalence relation $\sim$ on a set $A$ in classical mathematics, the quotient $A/{\sim}$ comes equipped with a quotient map $q:A\to A/{\sim}$ that satisfies two important properties: (1) The map $q$ satisfies the condition
\begin{equation*}
  q(x)=q(y)\leftrightarrow x\sim y,
\end{equation*}
and (2) the map $q$ is surjective. The first condition is called the \define{effectiveness} of the quotient map.

In classical mathematics, a map $f:A\to B$ is said to be surjective if for every $b\in B$ there exists an element $a\in A$ such that $f(a)=b$. Following the Curry-Howard interpretation, a map $f:A\to B$ is therefore surjective if it comes equipped with a dependent function
\begin{equation*}
  \prd{b:B}\sm{a:A}f(a)=b.
\end{equation*}
However, there is a subtle issue with this interpretation of surjectivity. It is somewhat stronger than the classical notion of surjectivity, because a dependent function $\prd{b:B}\sm{a:A}f(a)=b$ provides for every element $b:B$ an \emph{explicit} element $a:A$ equipped with an explicit identification $p:f(a)=b$, whereas in the classical notion of surjectivity such an element $a\in A$ is merely asserted to exist. To emphasize that the Curry-Howard interpretation of surjectivity is stronger than intended we make the following definition, and we will properly introduce surjective maps in \cref{subsec:surjective}.

\begin{defn}
  Consider a function $f:A\to B$. We say that $f$ is \define{split surjective} if it comes equipped with an element of type
  \begin{equation*}
    \issplitsurjective(f):=\prd{b:B}\sm{a:A}f(a)=b.
  \end{equation*}
\end{defn}

Martin-L\"of's dependent type theory doesn't have a general way of forming quotients of types. However, in the specific case of the congruence relations on $\N$ we can define the type of natural numbers modulo $k+1$ as the standard finite type $\Fin{k+1}$. We will show that $\Fin{k+1}$ comes equipped with a map
\begin{equation*}
  [\blank]_{k+1}:\N\to \Fin{k+1}
\end{equation*}
for each $k:\N$, and we will show in \cref{thm:effective-mod-k,thm:issec-nat-Fin} that this map satisfies conditions (1) and (2) in the split surjective sense.

To prepare for the definition of the quotient map $[\blank]_{k+1}$, we will first define a zero element of $\Fin{k+1}$ and successor function on each $\Fin{k}$. We will also define an auxiliary function $\skipzeroFin_k:\Fin{k}\to\Fin{k+1}$, which is used in the definition of the successor function. The map $[\blank]_{k+1}$ is then defined by iterating the successor function. 

\begin{defn} ~\nopagebreak
  \begin{enumerate}
  \item We define the \define{zero element} $\zeroFin_k:\Fin{k+1}$ recursively by
    \begin{align*}
      \zeroFin_0 & \defeq\ttt \\*
      \zeroFin_{k+1} & \defeq i(\zeroFin_k).
                       \intertext{Since there is a mismatch between the index of $\zeroFin_k$ and the index of its type, we will often simply write $\zeroFin$ or $0$ for the zero element of $\Fin{k+1}$.
    \item We define the function $\skipzeroFin_k:\Fin{k}\to\Fin{k+1}$ recursively by}
      \skipzeroFin_{k+1}(i(x)) & \defeq i(\skipzeroFin_k(x)) \\*
      \skipzeroFin_{k+1}(\ttt) & \defeq \ttt.
    \intertext{\item We define the \define{successor function} $\succFin_k:\Fin{k}\to\Fin{k}$ recursively by}
      \succFin_{k+1}(i(x)) & \defeq \skipzeroFin_k(x) \\*                       
      \succFin_{k+1}(\ttt)    & \defeq \zeroFin_k.
    \end{align*}
  \end{enumerate}
\end{defn}

\begin{defn}
  For any $k:\N$, we define the map $[\blank]_{k+1}:\N\to\Fin{k+1}$ recursively on $x$ by
  \begin{align*}
    [0]_{k+1} & \defeq 0 \\*
    [x+1]_{k+1} & \defeq \succFin_{k+1}[x]_{k+1}.
  \end{align*}
\end{defn}

Our next intermediate goal is to show that $x\equiv \natFin[x]_{k+1}\mod k+1$ for any natural number $x$. This fact is a consequence of the following simple lemma, that will help us compute with the maps $\natFin : \Fin{k}\to\N$.

\begin{lem}\label{lem:nat-Fin}
  We make three claims:
  \begin{enumerate}
  \item For any $k:\N$ there is an identification
    \begin{align*}
      \natFin(\zeroFin_k) & = 0
  \intertext{\item For any $k:\N$ and any $x:\Fin{k}$, we have}
      \natFin(\skipzeroFin_k(x)) & = \natFin(x)+1.
  \intertext{\item For any $k:\N$ and any $x:\Fin{k}$, we have}
      \natFin(\succFin_k(x)) & \equiv \natFin(x)+1 \mod k.
    \end{align*}
  \end{enumerate}
\end{lem}

\begin{proof}
  For the first claim, we define an identification $\alpha_k:\natFin(\zeroFin_k)=0$ recursively by
  \begin{align*}
    \alpha_0 & \defeq \refl{} \\
    \alpha_{k+1} & \defeq \alpha_k.
  \intertext{For the second claim, we define an identification $\beta_k(x):\natFin(\skipzeroFin_k(x))=\natFin(x)+1$ recursively by}
    \beta_{k+1}(i(x)) & \defeq \beta_k(x) \\
    \beta_{k+1}(\ttt) & \defeq \refl{}.
  \end{align*}
  For the third claim, we again define an element $\gamma_k(x):\natFin(\succFin_k(x)) \equiv \natFin(x)+1\mod{k}$ recursively. To obtain
  \begin{equation*}
    \gamma_{k+1}(i(x)) : \natFin(\succFin_{k+1}(i(x))) \equiv\natFin(i(x))+1\mod{k+1},
  \end{equation*}
  we calculate
  \begin{align*}
    \natFin(\succFin_{k+1}(i(x))) & \jdeq \natFin(\skipzeroFin(x)) & & \text{by definition of }\succFin\\
                                  & = \natFin(x)+1 & & \text{by claim (ii).}
  \end{align*}
  Since the congruence relation modulo $k+1$ is reflexive, we obtain $\gamma_{k+1}(i(x))$ from the identification of the above calculation. To obtain
  \begin{equation*}
    \gamma_{k+1}(\ttt) : \natFin(\succFin_{k+1}(\ttt)) \equiv \natFin(\ttt)+1\mod{k+1},
  \end{equation*}
  we calculate
  \begin{align*}
    \natFin(\succFin_{k+1}(\ttt)) & \jdeq \natFin(0) & & \text{by definition of }\succFin \\
                                  & = 0 & & \text{by claim (i)} \\
                                  & \equiv k+1 & & \text{by \cref{rmk:elementary-facts-div}} \\
                                  & \jdeq \natFin(\ttt)+1 & & \text{by definition of }\natFin.\qedhere
  \end{align*}
\end{proof}

\begin{prp}\label{prp:cong-nat-mod-succ}
  For any $x:\N$ we have
  \begin{equation*}
    \natFin[x]_{k+1}\equiv x \mod k+1.
  \end{equation*}
\end{prp}

\begin{proof}
  The proof by induction on $x$. The fact that
  \begin{equation*}
    \natFin[0]_{k+1}\equiv 0 \mod {k+1}
  \end{equation*}
  is immediate from the fact that $\natFin[0]_{k+1}\jdeq\natFin(0)=0$, which was shown in \cref{lem:nat-Fin}. In the inductive step, we have to show that
  \begin{equation*}
    \natFin[x+1]_{k+1}\equiv x+1\mod k+1.
  \end{equation*}
  This follows from the following computation
  \begin{align*}
    \natFin[x+1]_{k+1} & \jdeq \natFin(\succFin_{k+1}[x]_{k+1}) & & \text{by definition of }[\blank]_{k+1} \\
                       & \equiv \natFin[x]_{k+1}+1 & & \text{by \cref{lem:nat-Fin}} \\
                       & \equiv x+1 & & \text{by the inductive hypothesis.}\qedhere
  \end{align*}
\end{proof}

We need one more fact before we can prove \cref{thm:effective-mod-k,thm:issec-nat-Fin}.

\begin{prp}\label{cor:eq-congN}
  For any natural number $x<d$ we have
  \begin{equation*}
  d\mid x\leftrightarrow x=0.  
  \end{equation*}
  Consequently, for any two natural numbers $x$ and $y$ such that $\distN(x,y)<k$, we have
  \begin{equation*}
    x\equiv y\mod k\leftrightarrow x=y.
  \end{equation*}
\end{prp}

\begin{proof}
  Note that the implication $x=0\to d\mid x$ is trivial, so it suffices to prove the forward implication
  \begin{equation*}
    d\mid x \to x=0.
  \end{equation*}
  This implication clearly holds if $x\jdeq 0$. Therefore we only have to show that $d\mid x+1$ implies $x+1=0$, if we assume that $x+1<d$. In other words, we will derive a contradiction from the hypotheses that $x+1<d$ and $d\mid x+1$. To reach a contradiction we use \cref{ex:contradiction-le}, by which it suffices to show that $d\leq x+1$.
  
  We proceed by $\Sigma$-induction on the (unnamed) variable of type $d\mid x+1$, so we get to assume a natural number $k$ equipped with an identification $p:dk=x+1$. In the case where $k\jdeq 0$ we reach an immediate contradiction via \cref{prp:zero-one}, because we obtain that $0=d\cdot 0=x+1$. In the case where $k\jdeq\succN(k')$ it follows that
  \begin{equation*}
    d\leq dk'+ d\jdeq dk = x+1.\qedhere
  \end{equation*}
\end{proof}

\begin{thm}\label{thm:effective-mod-k}
  Consider a natural number $k$. Then we have
  \begin{equation*}
    [x]_{k+1}=[y]_{k+1} \leftrightarrow x\equiv y\mod k+1,
  \end{equation*}
  for any $x,y:\N$.
\end{thm}

\begin{proof}
  First note that, since $\natFin$ is injective by \cref{prp:is-injective-natFin}, we have
  \begin{align*}
    [x]_{k+1}=[y]_{k+1} & \leftrightarrow \natFin[x]_{k+1}=\natFin[y]_{k+1}.
  \end{align*}
  Since the inequalities $\natFin[x]_{k+1}<k+1$ and $\natFin[y]_{k+1}<k+1$ hold by \cref{lem:is-bounded-natFin}, it follows by \cref{cor:eq-congN} that
  \begin{equation*}
    \natFin[x]_{k+1}=\natFin[y]_{k+1}\leftrightarrow \natFin[x]_{k+1}\equiv\natFin[y]_{k+1}\mod k+1.   
  \end{equation*}
  The latter condition is by \cref{prp:cong-nat-mod-succ} equivalent to the condition that $x\equiv y\mod k+1$.
\end{proof}

\begin{thm}\label{thm:issec-nat-Fin}
  For any $x:\Fin{k+1}$ there is an identification
  \begin{equation*}
    [\natFin(x)]_{k+1}=x.
  \end{equation*}
  In other words, the map $[\blank]_{k+1}:\N\to \Fin{k+1}$ is split surjective.
\end{thm}

\begin{proof}
  Since $\natFin:\Fin{k+1}\to\N$ is injective by \cref{prp:is-injective-natFin}, it suffices to show that
  \begin{equation*}
    \natFin[\natFin(x)]_{k+1}=\natFin(x).
  \end{equation*}
  Now observe that $\natFin[\natFin(x)]_{k+1}<k+1$ and $\natFin(x)<k+1$. By \cref{cor:eq-congN} it therefore suffices to show that
  \begin{equation*}
    \natFin[\natFin(x)]_{k+1}\equiv\natFin(x)\mod{k+1}.
  \end{equation*}
  This fact is an instance of \cref{prp:cong-nat-mod-succ}.
\end{proof}

\subsection{The cyclic groups}
We can now define the cyclic groups $\Z/k$ for each $k:\N$. Note that $\Z/k$ must come equipped with the structure of a quotient $\Z/{\equiv}$ of $\Z$ by the congruence relation modulo $k$. In the case where $k\jdeq 0$, we have that $x\equiv y\mod{0}$ if and only if $x=y$. This motivates the following definition:

\begin{defn}\label{defn:Zk}
  We define the type $\Z/k$ for each $k:\N$ by
  \begin{equation*}
    \Z/0\defeq \Z\qquad\text{and}\qquad \Z/{(k+1)}\defeq\Fin{k+1}.
  \end{equation*}
\end{defn}

Recall from \cref{ex:int_group_laws} that $\Z/0$ already comes equipped with the structure of a group, but the group structure on $\Z/{(k+1)}$ remains to be defined.

\begin{defn}
  We define the \define{addition} operation on $\Z/{(k+1)}$ by
  \begin{equation*}
    x+y\defeq[\natFin(x)+\natFin(y)]_{k+1},
  \end{equation*}
  and we define the \define{additive inverse} operation on $\Z/{(k+1)}$ by
  \begin{equation*}
    -x\defeq[\distN(\natFin(x),k+1)]_{k+1}.
  \end{equation*}
\end{defn}

\begin{rmk}
  The following congruences modulo $k+1$ follow immediately from \cref{prp:cong-nat-mod-succ}:
  \begin{align*}
    \natFin(0) & \equiv 0 \\
    \natFin(x+y) & \equiv \natFin(x)+\natFin(y) \\
    \natFin(-x) & \equiv \distN(\natFin(x),k+1).
  \end{align*}
\end{rmk}

Before we show that addition on $\Z/{k}$ satisfies the group laws, we have to show that addition on $\N$ preserves the congruence relation.

\begin{prp}
  Consider $x,y,x',y':\N$. If any two of the following three properties hold, then so does the third:
  \begin{enumerate}
  \item $x\equiv x'\mod k$,
  \item $y\equiv y'\mod k$,
  \item $x+y\equiv x'+y'\mod k$.
  \end{enumerate}
\end{prp}

\begin{proof}
  Recall that the distance function $\distN$ is translation invariant by \cref{ex:translation-invariant-distN}. Therefore it follows that
  \begin{equation}\label{eq:translation-invariant-congN}
    a\equiv b\mod k \leftrightarrow a+c\equiv b+c\mod k.\tag{\textasteriskcentered}
  \end{equation}
  We will use this observation to prove the claim.
  
  First, suppose that $x\equiv x'$ and $y\equiv y'$ modulo $k$. Then it follows by \cref{eq:translation-invariant-congN} that
  \begin{equation*}
    x+y\equiv x'+y\equiv x'+y'.
  \end{equation*}
  This shows that (i) and (ii) together imply (iii).

  Next, suppose that $x\equiv x'$ and $x+y\equiv x'+y'$ modulo $k$. Then it follows that
  \begin{equation*}
    x+y\equiv x'+y'\equiv x+y'.
  \end{equation*}
  Applying \cref{eq:translation-invariant-congN} once more in the reverse direction, we obtain that $y\equiv y'$ modulo $k$. This shows that (i) and (iii) together imply (ii).

  The remaining claim, that (ii) and (iii) together imply (i), follows by commutativity of addition from the fact that (i) and (iii) together imply (ii).
\end{proof}

\begin{thm}
  The addition operation on $\Z/{k}$ satisfies the laws of an abelian group:
  \begin{align*}
    0+x & = x & x+0 & = x \\
    (-x)+x & = 0 & x+(-x) & = 0 \\
    (x+y)+z & = x+(y+z) & x+y & = y+x. 
  \end{align*}
\end{thm}

\begin{proof}
  The fact that the addition operation on $\Z/0$ satisfies the laws of an abelian group was stated as \cref{ex:int_group_laws}. Therefore we will only show that addition on $\Z/{(k+1)}$ satisfies the laws of an abelian group.

  We first note that by commutativity of addition on $\N$, it follows immediately that addition on $\Z/{(k+1)}$ is commutative.

  To prove associativity, note that by \cref{thm:effective-mod-k} it suffices to show that
  \begin{equation*}
    \natFin(x+y)+\natFin(z)\equiv\natFin(x)+\natFin(y+z)\mod k+1.
  \end{equation*}
  Since addition on $\Z/{(k+1)}$ maps preserves the congruence relation, and since we have the congruences
  \begin{align*}
    \natFin(x+y) & \equiv \natFin(x)+\natFin(y) \mod k+1 \\
    \natFin(y+z) & \equiv \natFin(y)+\natFin(z) \mod k+1,
  \end{align*}
  it suffices to show that
  \begin{equation*}
    (\natFin(x)+\natFin(y))+\natFin(z) \equiv \natFin(x)+(\natFin(y)+\natFin(z)) \mod k+1.
  \end{equation*}
  This follows immediately by associativity of addition on $\N$.

  To show that addition on $\Z/{(k+1)}$ satisfies the right unit law, we first observe that it suffices to show that
  \begin{equation*}
    [\natFin(x)+\natFin(0)]_{k+1}=[\natFin(x)]_{k+1}
  \end{equation*}
  because there is an identification $[\natFin(x)]_{k+1}=x$ by \cref{thm:issec-nat-Fin}. By \cref{thm:effective-mod-k} it now suffices tho show that
  \begin{equation*}
    \natFin(x)+\natFin(0)\equiv\natFin(x)\mod k+1. 
  \end{equation*}
  This follows immediately from the fact that $\natFin(0)=0$. The left unit law now follows from the right unit law by commutativity. We leave the inverse laws as an exercise.
\end{proof}

\begin{exercises}
  \exitem \label{ex:div-3-for-2}Complete the proof of \cref{prp:div-3-for-2}.
  \exitem \label{ex:is-poset-div}Show that the divisibility relation satisfies the axioms of a poset, i.e., that it is reflexive, antisymmetric, and transitive.
  \exitem \label{ex:div-factorial}Construct a dependent function
  \begin{equation*}
    \prd{x:\N}(x\neq 0)\to ((x\leq n)\to (x\mid n!))
  \end{equation*}
  for every $n:\N$.
  \exitem Define $1\defeq[1]_{k+1}:\Fin{k+1}$. Show that
  \begin{equation*}
    \succFin_{k+1}(x)=x+1
  \end{equation*}
  for any $x:\Fin{k+1}$.
  \exitem \label{ex:Eq-Fin}The observational equality on $\Fin{k}$ is a binary relation
  \begin{equation*}
    \EqFin_{k}:\Fin{k}\to(\Fin{k}\to\UU_0)
  \end{equation*}
  defined recursively by
  \begin{align*}
    \EqFin_{k+1}(i(x),i(y)) & \defeq \EqFin_k(x,y) & \EqFin_{k+1}(i(x),\ttt) & \defeq \emptyt \\*
    \EqFin_{k+1}(\ttt,i(y)) & \defeq \emptyt & \EqFin_{k+1}(\ttt,\ttt) & \defeq \unit.
  \end{align*}
  \begin{subexenum}
  \item \label{ex:eq-iff-Eq-Fin}Show that
  \begin{equation*}
    (x=y)\leftrightarrow \EqFin_k(x,y)
  \end{equation*}
  for any two elements $x,y:\Fin{k}$.
  \item \label{ex:is-injective-i-Fin}Show that the function $i:\Fin{k}\to\Fin{k+1}$ is injective, for each $k:\N$.
  \item \label{ex:neq-zero-succ-Fin}
  Show that
  \begin{equation*}
    \succFin_{k+1}(i(x))\neq 0
  \end{equation*}
  for any $x:\Fin{k}$.
  \item Show that function $\succFin_k:\Fin{k}\to\Fin{k}$ is injective, for each $k:\N$.
  \end{subexenum}
  \exitem \label{ex:has-inverse-succ-Fin}The predecessor function $\predFin_k:\Fin{k}\to\Fin{k}$ is defined in three steps, just as in the definition of the successor function on $\Fin{k}$.
  \begin{enumerate}
  \item We define the element $\negtwoFin_k:\Fin{k+1}$ by
    \begin{align*}
      \negtwoFin_0 & \defeq\ttt \\*
      \negtwoFin_{k+1} & \defeq i(\ttt).
    \intertext{\item We define the function $\skipnegtwoFin_k:\Fin{k}\to\Fin{k+1}$ recursively by}
      \skipnegtwoFin_{k+1}(i(x)) & \defeq i(i(x)) \\*
      \skipnegtwoFin_{k+1}(\ttt) & \defeq \ttt.
    \intertext{\item Finally, we define the \define{predecessor function} $\predFin_k:\Fin{k}\to\Fin{k}$ recursively by}
      \predFin_{k+1}(i(x)) & \defeq \skipnegtwoFin_k(\predFin_k(x)) \\*                       
      \predFin_{k+1}(\ttt)    & \defeq \negtwoFin_k.
    \end{align*}
  \end{enumerate}
  Show that $\predFin_k$ is an inverse to $\succFin_k$, i.e., construct identifications
  \begin{equation*}
    \succFin_k(\predFin_k(x))=x,\qquad\text{and}\qquad\predFin_k(\succFin_k(x))=x
  \end{equation*}
  for each $x:\Fin{k}$.
  \exitem \label{ex:classical-Fin}Recall that
  \begin{equation*}
    \classicalFin_k:=\sm{x:\N}x<k.
  \end{equation*}
  \begin{subexenum}
  \item Show that
    \begin{equation*}
      (x=y)\leftrightarrow (\proj 1(x)=\proj 1(y))
    \end{equation*}
    for each $x,y:\classicalFin_k$.
  \item By \cref{lem:is-bounded-natFin} it follows that the map $\natFin :\Fin{k}\to\N$ induces a map $\natFin:\Fin{k}\to\classicalFin_k$. Construct a map
    \begin{equation*}
      \alpha_k:\classicalFin_k \to \Fin{k}  
    \end{equation*}
    for each $k:\N$, and show that
  \begin{equation*}
    \alpha_k(\natFin(x)) = x \qquad\text{and}\qquad \natFin(\alpha_k(y)) = y
  \end{equation*}
  for each $x:\Fin{k}$ and each $y:\classicalFin_k$. 
  \end{subexenum}
  \exitem \label{ex:ring-Fin}The multiplication operation $x,y\mapsto xy$ on $\Z/{(k+1)}$ is defined by
  \begin{equation*}
    xy \defeq [\natFin(x)\natFin(y)]_{k+1}.
  \end{equation*}
  \begin{subexenum}
  \item Show that $\natFin(xy)\equiv\natFin(x)\natFin(y)\mod{k+1}$ for each $x,y:\Z/{(k+1)}$.
  \item \label{ex:congruence-mulN}Show that
    \begin{equation*}
      xy\equiv x'y'\mod k
    \end{equation*}
    for any $x,y,x',y':\N$ such that $x\equiv x'$ and $y\equiv y' \mod k$.
  \item Show that multiplication on $\Z/{(k+1)}$ satisfies the laws of a commutative ring:
    \begin{align*}
      (xy)z & = x(yz) & xy & = yx \\
      1x & = x & x1 & = x \\
      x(y+z) & = xy+xz & (x+y)z & = xz+yz.
    \end{align*}
  \end{subexenum}
  \exitem \label{ex:euclidean-division}(Euclidean division) Consider two natural numbers $a$ and $b$.
  \begin{subexenum}
  \item Construct two natural numbers $q$ and $r$ such that $(b\neq 0) \to (r<b)$, along with an identification
    \begin{equation*}
      a=qb+r.
    \end{equation*}
  \item Show that for any four natural numbers $q,q'$ and $r,r'$ such that the implications $(b\neq 0) \to (r<b)$ and $(b\neq 0)\to (r'<b)$ hold, and for which there are identifications
    \begin{equation*}
      a=qb+r\qquad\text{and}\qquad a=q'b+r',
    \end{equation*}
    we have $q=q'$ and $r=r'$.
  \end{subexenum}
  \exitem The type $\N_k$ of \define{$k$-ary natural numbers} is an inductive type with the following constructors:
  \begin{align*}
    \constantbasedN{k} & : \Fin{k}\to\basedN{k} \\
    \unaryopbasedN{k} & : \Fin{k}\to (\basedN{k}\to\basedN{k}).
  \end{align*}
  A $k$-ary natural number can be converted back into an ordinary natural number via the function $\convertbasedN{k}:\basedN{k}\to\N$, which is defined recursively by
  \begin{align*}
    \convertbasedN{k}(\constantbasedN{k}(x)) & \defeq \natFin(x) \\
    \convertbasedN{k}(\unaryopbasedN{k}(x,n)) & \defeq k(\convertbasedN{k}(n)+1)+\natFin(x).
  \end{align*}
  \begin{subexenum}
  \item Show that the type $\basedN{0}$ is empty.
  \item Show that the function $\convertbasedN{k}:\basedN{k}\to\N$ is injective.
  \item Show that the function $\convertbasedN{k+1}:\basedN{k+1}\to\N$ has an inverse, i.e. construct a function
    \begin{equation*}
      g_{k} : \N\to\basedN{k+1}
    \end{equation*}
    equipped with identifications
    \begin{align*}
      \convertbasedN{k+1}(g_k(n)) & = n \\
      g_{k}(\convertbasedN{k+1}(x)) & = x
    \end{align*}
    for each $n:\N$ and each $x:\basedN{k+1}$.
  \end{subexenum}
\end{exercises}

\section{Decidability in elementary number theory}\label{sec:decidability}

Martin-L\"of's dependent type theory is a foundation for constructive mathematics, but in constructive mathematics there is no way to show that $P\lor\neg P$ holds for an arbitrary proposition $P$. Likewise, in type theory there is no way to construct an element of type $A+\neg A$ for an arbitrary type $A$. Consequently, if we want to reason by case analysis over whether $A$ is empty or nonempty, we first have to \emph{show} that $A+\neg A$ holds.

A type $A$ that comes equipped with an element of type $A+\neg A$ is said to be \emph{decidable}. Even though we cannot show that all types are decidable, many types are indeed decidable. Examples include the empty type and any type that comes equipped with a point, such as the type of natural numbers.

Decidability is an important concept with many applications in number theory and finite mathematics, and in this section we will explore the applications of decidability to elementary number theory. For example, the natural numbers satisfy a well-ordering principle with respect to decidable type families over the natural numbers; decidability can be used to construct the greatest common divisor of any two natural numbers; and it can also be used to show that there are infinitely many prime numbers.

\subsection{Decidability and decidable equality}

\begin{defn}
  A type $A$ is said to be \define{decidable}\index{decidable type|textbf} if it comes equipped with an element of type\index{is-decidable@{$\isdecidable(A)$}}
  \begin{equation*}
    \isdecidable(A)\defeq A+\neg A.
  \end{equation*}
  A family $P$ over a type $A$ is said to be \define{decidable}\index{decidable family of types|textbf} if $P(x)$ is decidable for every $x:A$.
\end{defn}

\begin{eg}
  The principal way to show that a type $A$ is decidable is to either construct an element $a:A$, or to construct a function $A\to\emptyt$. For example, the types $\unit$ and $\emptyt$ are decidable. Indeed, we have\index{decidable type!unit type}\index{unit type!is decidable}\index{decidable type!empty type}\index{empty type!is decidable}\index{decidable type!type with an element}
  \begin{align*}
    \inl(\ttt) & :\isdecidable(\unit) \\
    \inr(\idfunc) & : \isdecidable(\emptyt).
  \end{align*}
  Furthermore, any type $A$ equipped with an element $a:A$ is decidable because we have $\inl(a):\isdecidable(A)$ for such $A$.
\end{eg}

\begin{eg}\label{eg:decidability-closure}
  The principal way to use a hypothesis that $A$ is decidable is to proceed by the induction principle of coproducts, i.e., to proceed by case analysis.
  
  For example, if $A$ and $B$ are decidable types, then the types $A+B$, $A\times B$, and $A\to B$ are also decidable. This is straightforward to prove directly by pattern-matching on the variables of type $\isdecidable(A)$ and $\isdecidable(B)$. When we go through these proofs, the familiar truth table emerges:
  \begin{center}
    \begin{tabular}{lllll}
      \toprule
      \multicolumn{5}{c}{$\isdecidable$} \\ \cmidrule{1-5}
      $A$ & $B$ & $A+B$ & $A\times B$ & $A\to B$ \\
      \midrule
      $\inl(a)$ & $\inl(b)$ & $\inl(\inl(a))$ & $\inl(a,b)$ & $\inl(\lam{x}b)$ \\
      $\inl(a)$ & $\inr(g)$ & $\inl(\inl(a))$ & $\inr(g\circ \proj 2)$ & $\inr(\lam{h}g(h(a)))$ \\
      $\inr(f)$ & $\inl(b)$ & $\inl(\inr(b))$ & $\inr(f\circ\proj 1)$ & $\inl(\exfalso\circ f)$ \\
      $\inr(f)$ & $\inr(g)$ & $\inr[f,g]$ & $\inr(f\circ\proj 1)$ & $\inl(\exfalso\circ f)$ \\
      \bottomrule
    \end{tabular}
  \end{center}
  Since $A\to B$ is decidable whenever both $A$ and $B$ are decidable, it also follows that the negation $\neg A$ of any decidable type $A$ is decidable.
\end{eg}

\begin{eg}\label{eg:is-decidable-EqN}
  Since the empty type and the unit type are both decidable types, it also follows that the types $\EqN(m,n)$, $m\leq n$ and $m<n$ are decidable for each $m,n:\N$. The proofs in each of the three cases is by induction on $m$ and $n$.

  For instance, to show that $\EqN(m,n)$ is decidable for each $m,n:\N$, we simply note that the types
  \begin{align*}
    \EqN(\zeroN,\zeroN) & \jdeq \unit \\
    \EqN(\zeroN,\succN(n)) & \jdeq \emptyt \\
    \EqN(\succN(m),\zeroN) & \jdeq \emptyt 
  \end{align*}
  are all decidable, and that the type $\EqN(\succN(m),\succN(n))\jdeq \EqN(m,n)$ is decidable by the inductive hypothesis.
\end{eg}

The fact that $\N$ has decidable observational equality also implies that equality itself is decidable on $\N$. This leads to the general concept of decidable equality, which is important in many results about decidability.

\begin{defn}
  We say that a type $A$ has \define{decidable equality}\index{decidable equality|textbf} if the identity type $x=y$ is decidable for every $x,y:A$. We will write\index{has-decidable-equality(A)@{$\hasdecidableequality(A)$}|textbf}
  \begin{equation*}
    \hasdecidableequality(A)\defeq \prd{x,y:A}\isdecidable(x=y).
  \end{equation*}
\end{defn}

Before we show that $\N$ has decidable equality, let us show that if $A\leftrightarrow B$ and $A$ is decidable, then $B$ must be decidable.

\begin{lem}\label{lem:is-decidable-iff}
  Consider two types $A$ and $B$, and suppose that $A\leftrightarrow B$. Then $A$ is decidable if and only if $B$ is decidable.
\end{lem}

\begin{proof}
  Since we have functions $f:A\to B$ and $g:B\to A$ by assumption, we obtain by \cref{prp:contravariant-neg} the functions
  \begin{align*}
    \tilde{f} & : \neg B\to\neg A \\
    \tilde{g} & : \neg A \to \neg B.
  \end{align*}
  By \cref{rmk:functor-coprod} we have therefore the functions
  \begin{align*}
    f+\tilde{g} & : (A+\neg A) \to (B+\neg B) \\
    g+\tilde{f} & : (B+\neg B) \to (A+\neg A).\qedhere
  \end{align*}
\end{proof}

\begin{prp}\label{prp:has-decidable-equality-N}
  Equality on the natural numbers is decidable.\index{natural numbers!decidable equality}\index{N@{$\N$}!has decidable equality}\index{has decidable equality!natural numbers}\index{decidable equality!of N@{of $\N$}}
\end{prp}

\begin{proof}
  Recall from \cref{prp:Eq-eq-N} that we have
  \begin{equation*}
    (m=n)\leftrightarrow \EqN(m,n).
  \end{equation*}
  The claim therefore follows by \cref{lem:is-decidable-iff}, since we have observed in \cref{eg:is-decidable-EqN} that $\EqN(m,n)$ is decidable for every $m,n:\N$.
\end{proof}

It is certainly not provable with the given rules of type theory that every type has decidable equality. In fact, we will show in \cref{thm:hedberg} that if a type has decidable equality, then it is a \emph{set}. However, it is also not provable that every set has decidable equality unless one assumes the \emph{law of excluded middle}. We will discuss this principle in \cref{sec:logic}. For now, it is important to remember that in order to use decidability, we must first \emph{prove that it holds}, and many familiar types do indeed have decidable equality.

\begin{prp}\label{prp:has-decidable-equality-Fin}
  The standard finite type $\Fin{k}$ has decidable equality for each $k:\N$.\index{Fin k@{$\Fin{k}$}!has decidable equality}\index{has decidable equality!Fin k@{$\Fin{k}$}}\index{decidable equality!of Fin k@{of $\Fin{k}$}}
\end{prp}

\begin{proof}
  Recall from \cref{ex:Eq-Fin} that we constructed an observational equality relation $\EqFin_k$ on $\Fin{k}$ for each $k:\N$, which satisfies
  \begin{equation*}
    (x=y)\leftrightarrow \EqFin_k(x,y).
  \end{equation*}
  The type $\EqFin_k(x,y)$ is decidable, since it is recursively defined using the decidable types $\emptyt$ and $\unit$.
\end{proof}

We can use the fact that the finite types $\Fin{k}$ have decidable equality to show that the divisibility relation on $\N$ is decidable. 

\begin{thm}\label{thm:is-decidable-div-N}
  For any $d,x:\N$, the type $d\mid x$ is decidable.\index{divisibility on N@{divisibility on $\N$}!is decidable}
\end{thm}

\begin{proof}
  Note that $0\mid x$ is decidable because $0\mid x$ if and only if $x=0$, which is decidable by \cref{prp:has-decidable-equality-N}. Therefore it suffices to show that $d+1\mid x$ is decidable.

  By \cref{thm:effective-mod-k} it follows that $d+1\mid x$ holds if and only if we have an identification $[x]_{d+1}=0$ in $\Fin{d+1}$. Therefore the claim follows from the fact that $\Fin{d+1}$ has decidable equality.
\end{proof}

\subsection{Constructions by case analysis}
\index{case analysis|(}
\index{with-abstraction|(}

A common way to construct functions and to prove properties about them is by case analysis. For example, a famous function of Collatz is specified by case analysis on whether $n$ is even or odd:
  \begin{equation*}
  \collatz(n) =
  \begin{cases}
    n/2 & \text{if $n$ is even}\\
    3n+1 & \text{if $n$ is odd.}
  \end{cases}
\end{equation*}
The Collatz function is of course uniquely determined by this specification, but it is important to note that there is a bit of work to be done in order to define the Collatz function according to the rules of dependent type theory. First we note that, since the Collatz function is specified by case analysis on whether $n$ is even or odd, we will have to use a dependent function witnessing the fact that every number is either even or odd. In other words, we will make use of the dependent function
\begin{equation*}
  d:\prd{n:\N}\isdecidable(2\mid n),
\end{equation*}
which we have by \cref{thm:is-decidable-div-N}. The type $\isdecidable(2\mid n)$ is the coproduct $(2\mid n)+(2\nmid n)$, so the idea is to proceed by case analysis on whether $d(n)$ is of the form $\inl(x)$ or $\inr(x)$, i.e., by the induction principle of coproducts. However, $d(n)$ is not a free variable of type $\isdecidable(2\mid n)$. Before we can proceed by induction, we must therefore first \emph{generalize} the element $d(n)$ to a free variable $y:\isdecidable(2\mid n)$. In other words, we will first define a function
\begin{equation*}
  h:\prd{n:\N} (\isdecidable(2\mid n)\to\N)
\end{equation*}
by the induction principle of coproducts, and then we obtain the Collatz function by substituting $d(n)$ for $y$ in $h(n,y)$. Putting these ideas together, we obtain the following type theoretical definition of the Collatz function.

\begin{defn}
  Write $d:\prd{n:\N}\isdecidable(2\mid n)$ for the function deciding $2\mid n$, given in \cref{thm:is-decidable-div-N}.
  \begin{enumerate}
  \item We define a function $h:\prd{n:\N}(\isdecidable(2\mid n)\to \N)$ by
    \begin{align*}
      h(n,\inl(m,p)) & \defeq m \\
      h(n,\inr(f)) & \defeq 3n+1.
    \end{align*}
  \item We define the \define{collatz function}\index{Collatz function|textbf} $\collatz:\N\to \N$ by\index{collatz@{$\collatz$}|textbf}
    \begin{equation*}
      \collatz(n)\defeq h(n,d(n)).
    \end{equation*}
  \end{enumerate}
\end{defn}

\begin{rmk}
  The general ideas behind the formal construction of the Collatz function lead to the type theoretic concept of \emph{with-abstraction}. With-abstraction is a type-theoretically precise generalization of case analysis.
  
  In full generality, if our goal is to define a dependent function $f:\prd{x:A}C(x)$, and we already have a function $g:\prd{x:A}B(x)$, then it suffices to define a dependent function
  \begin{equation*}
    h:\prd{x:A}B(x)\to C(x).
  \end{equation*}
  Indeed, given $g$ and $h$ as above, we can define $f\defeq \lam{x}h(x,g(x))$. In other words, to define $f(x)$ using $g(x):B(x)$, we generalize $g(x)$ to an arbitrary element $y:B(x)$ and proceed to define an element $h(x,y):C(x)$.

  With-abstraction is a concise way to present such a definition. In a definition by with-abstraction, we may write
  \begin{equation*}
    f(x)\with [g(x)/y]\defeq h(x,y),
  \end{equation*}
  to define a function $f:\prd{x:A}C(x)$ that satisfies the judgmental equality $f\jdeq \lam{x}h(x,g(x))$. In other words, $f(x)$ is defined to be $h(x,y)$ with $g(x)$ for $y$.
  
  The definition of the Collatz function can therefore be given by with-abstraction as
  \begin{equation*}
    \collatz(n)\with [d(n)/y]\defeq h(x,y).
  \end{equation*}
  However, recall that the function $h$ was defined by pattern matching on $y$. We can combine with-abstraction and pattern matching to obtain a \emph{direct} definition of the Collatz function that doesn't explicitly mention the function $h$ anymore. This gives us the following concise way to define the Collatz function:
  \begin{align*}
    \collatz(n)\with [d(n)/\inl(m,p)] & \defeq m \\
    \collatz(n)\with [d(n)/\inr(f)] & \defeq 3n+1.
  \end{align*}
  Notice that in addition to the information in the specification of the Collatz function, the definition by with-abstraction also tells us which decision procedure was used to decide whether $n$ is even or not. The combination of with-abstraction and pattern matching, which allows us to skip the explicit definition of the function $h$, is what makes with-abstraction so useful.
  \end{rmk}
  
  Using with-abstraction we can find a slight improvement of the decidability results of $A\to B$ and $A\times B$ in \cref{eg:decidability-closure}, and we will use these improved claims in the construction of the greatest common divisor.

\begin{prp}\label{prp:is-decidable-function-type}
  Consider a decidable type $A$, and let $B$ be a type equipped with a function
  \begin{equation*}
    A\to\isdecidable(B).
  \end{equation*}
  Then the types $A\times B$ and $A\to B$ are also decidable.
\end{prp}

\begin{proof}
  We only prove the claim about the decidability of $A\to B$, since the claim about the decidability of $A\times B$ is proven similarly. Since $A$ is assumed to be decidable, we proceed by case analysis on $A+\neg A$. In the case where we have $f:\neg A$, we have the functions
  \begin{equation*}
    \begin{tikzcd}[column sep=large]
      A \arrow[r,"f"] & \emptyt \arrow[r,"\exfalso"] & B.
    \end{tikzcd}
  \end{equation*}
  Therefore we obtain the element $\inl(\exfalso\circ f):\isdec(A\to B)$. In the case where we have an element $a:A$, we have to construct a function
  \begin{equation*}
    d:(A\to\isdecidable(B))\to\isdecidable(A\to B)
  \end{equation*}
  Given $H:A\to\isdecidable(B)$, we can use with-abstraction to proceed by case analysis on $H(a):B+\neg B$. The function $d$ is therefore defined as
  \begin{align*}
    d(H)\with [H(a)/\inl(b)] & \defeq \inl(\lam{x}b) \\
    d(H)\with [H(a)/\inr(g)] & \defeq \inr(\lam{h}g(h(a))).\qedhere
  \end{align*}
\end{proof}

For a general family of decidable types $P$ over $\N$, we cannot prove that the type
\begin{equation*}
  \prd{x:\N}P(x)
\end{equation*}
is decidable. However, if we know in advance that $P(x)$ holds for any $m\leq x$, then we can decide $\prd{x:\N}P(x)$ by checking the decidability of each $P(x)$ until $m$. 

\begin{prp}\label{prp:is-decidable-pi-type}
  Consider a decidable type family $P$ over $\N$ equipped with a natural number $m$ such that the type
  \begin{equation*}
    \prd{x:\N}(m\leq x)\to P(x)
  \end{equation*}
  is decidable. Then the type $\prd{x:\N}P(x)$ is decidable. 
\end{prp}

\begin{proof}
  Our proof is by induction on $m$, but we will first make sure that the inductive hypothesis will be strong enough by quantifying over all decidable type families over $\N$. Of course, we cannot do this directly. However, by the assumption that there are enough universes (\cref{enough-universes}), there is a universe $\UU$ that contains $P$. We fix this universe, and we will prove by induction on $m$ that for every decidable type family $Q:\N\to\UU$ for which the type
  \begin{equation*}
    \prd{x:\N}(m\leq x)\to Q(x),
  \end{equation*}
  is decidable, the type $\prd{x:\N}Q(x)$ is again decidable.

  In the base case, it follows by assumption that the type $\prd{x:A}Q(x)$ is decidable. For the inductive step, let $Q:\N\to\UU$ be a decidable type family for which the type
  \begin{equation*}
    \prd{x:\N}(m+1\leq x)\to Q(x)
  \end{equation*}
  is decidable. Since $Q$ is assumed to be decidable, we can proceed by case analysis on $Q(0)+\neg Q(0)$. In the case of $\neg Q(0)$, it follows that $\neg \prd{x:\N}Q(x)$. In the case where we have $q:Q(0)$, consider the type family $Q':\N\to\UU$ given by
  \begin{equation*}
    Q'(x)\defeq Q(x+1).
  \end{equation*}
  Then $Q'$ is decidable since $Q$ is decidable, and moreover it follows that the type $\prd{x:\N} ({m\leq x})\to Q'(x)$ is decidable. The inductive hypothesis implies therefore that the type $\prd{x:\N}Q'(x)$ is decidable. In the case where $\neg\prd{x:\N}Q'(x)$, it follows that $\neg\prd{x:\N}Q(x)$, and in the case where we have a function $g:\prd{x:\N}Q'(x)$, we can construct a function $f:\prd{x:\N}Q(x)$ by
  \begin{align*}
    f(0) & \defeq q \\
    f(x+1) & \defeq g(x)\qedhere.
  \end{align*}
\end{proof}

\begin{cor}\label{cor:is-decidable-bounded-pi}
  Consider two decidable families $P$ and $Q$ over $\N$, and suppose that $P$ comes equipped with an upper bound $m$. Then the type
  \begin{equation*}
    \prd{n:\N}P(n)\to Q(n)
  \end{equation*}
  is decidable.
\end{cor}

\begin{proof}
  Since $m$ is assumed to be an upper bound for $P$, it follows $P(n)\to Q(n)$ for any $m\leq n$. With this observation we apply \cref{prp:is-decidable-pi-type}.
\end{proof}
\index{case analysis|)}
\index{with-abstraction|)}

\subsection{The well-ordering principle of \texorpdfstring{$\N$}{ℕ}}
\index{well-ordering principle of N@{well-ordering principle of $\N$}|(}
\index{natural numbers!well-ordering principle|(}

The well-ordering principle of the natural numbers in classical mathematics asserts that any nonempty subset of $\N$ has a least element. To formulate the well-ordering principle in type theory, we will use type families over $\N$ instead of subsets of $\N$. Moreover, the classical well-ordering principle tacitly assumes that subsets are decidable. The type theoretic well-ordering principle of $\N$ is therefore formulated using \emph{decidable} families over $\N$.

\begin{defn}
  Let $P$ be a family over $\N$, not necessarily decidable.
  \begin{enumerate}
  \item We say that a natural number $n$ is a \define{lower bound}\index{lower bound|textbf} for $P$ if it comes equipped with an element of type\index{is-lower-bound P (n)@{$\islowerbound_P(n)$}|textbf}
    \begin{equation*}
      \islowerbound_P(n)\defeq \prd{x:\N}P(x)\to (n\leq x).
    \end{equation*}
  \item We say that a natural number $n$ is an \define{upper bound}\index{upper bound|textbf} for $P$ if it comes equipped with an element of type\index{is-upper-bound P (n)@{$\isupperbound_P(n)$}|textbf}
    \begin{equation*}
      \isupperbound_P(n)\defeq \prd{x:\N}P(x)\to (x\leq n).
    \end{equation*}
  \end{enumerate}
\end{defn}

  A minimal element of $P$ is therefore a natural number $n$ for which $P(n)$ holds, and which is also a lower bound for $P$. The well-ordering principle of $\N$ asserts that such an element exists for any decidable family $P$, as soon as $P(n)$ holds for some $n$.

  \begin{thm}[Well-ordering principle of $\N$]\label{thm:well-ordering-principle-N}\index{well-ordering principle of N@{well-ordering principle of $\N$}|textbf}\index{natural numbers!well-ordering principle|textbf}
    Let $P$ be a decidable family over $\N$, where $d$ witnesses that $P$ is decidable. Then there is a function\index{well-ordering-principle P d@{$\wellorderingprinciple(P,d)$}|textbf}
  \begin{equation*}
    \wellorderingprinciple(P,d):\Big(\sm{n:\N}P(n)\Big)\to\Big(\sm{m:\N}P(m)\times\islowerbound_P(m)\Big).
  \end{equation*}
\end{thm}

\begin{proof}
  By the assumption that there are enough universes (\cref{enough-universes}), there is a universe $\UU$ that contains $P$. Instead of proving the claim for the given type family $P$, we will show by induction on $n:\N$ that there is a function
  \begin{equation*}\label{eq:well-ordering}
    Q(n)\to \Big(\sm{m:\N}Q(m)\times\islowerbound_Q(m)\Big)\tag{\textasteriskcentered}
  \end{equation*}
  for every decidable family $Q:\N\to\UU$. Note that we are now also quantifying over the decidable families $Q:\N\to\UU$. This slightly strengthens the inductive hypothesis, which we will be able to exploit.

  The base case is trivial, since $\zeroN$ is a lower bound of every type family over $\N$. For the inductive step, assume that \cref{eq:well-ordering} holds for every decidable type family $Q:\N\to \UU$. Furthermore, let $Q:\N\to\UU$ be a decidable type family equipped with an element $q:Q(\succN(n))$. Our goal is to construct an element of type
  \begin{equation*}
    \sm{m:\N}Q(m)\times\islowerbound_Q(m).
  \end{equation*}
  Since $Q(\zeroN)$ is assumed to be decidable, it suffices to construct a function
  \begin{equation*}
    (Q(\zeroN)+\neg Q(\zeroN))\to \sm{m:\N}Q(m)\times\islowerbound_Q(m).
  \end{equation*}
  Therefore we can proceed by case analysis on $Q(\zeroN)+\neg Q(\zeroN)$. In the case where we have an element of type $Q(\zeroN)$, it follows immediately that $\zeroN$ must be minimal. In the case where $\neg Q(\zeroN)$, we consider the decidable subset $Q'$ of $\N$ given by
  \begin{equation*}
    Q'(n)\defeq Q(\succN(n)).
  \end{equation*}
  Since we have $q:Q'(n)$, we obtain a minimal element in $Q'$ by the inductive hypothesis. Of course, by the assumption that $Q(\zeroN)$ doesn't hold, the minimal element of $Q'$ is also the minimal element of $Q$.
\end{proof}
\index{well-ordering principle of N@{well-ordering principle of $\N$}|)}
\index{natural numbers!well-ordering principle|)}

\subsection{The greatest common divisor}
\index{natural numbers!greatest common divisor|(}
\index{greatest common divisor|(}

The greatest common divisor of two natural numbers $a$ and $b$ is a natural number $\gcd(a,b)$ that satisfies the property that
\begin{equation*}
  x\mid a\ \text{and}\ x\mid b\qquad\text{if and only if}\qquad x\mid\gcd(a,b)
\end{equation*}
for any $x:\N$. In other words, any number $x:\N$ that divides both $a$ and $b$ also divides the greatest common divisor. Moreover, since $\gcd(a,b)$ divides itself, it follows from the reverse implication that $\gcd(a,b)$ divides both $a$ and $b$.

This property can also be seen as the \emph{specification} of what it means to be a greatest common divisor of $a$ and $b$. In formal developments of mathematics, when you're about to construct an object that satisfies a certain specification, it can be useful to start out with that specification. For example, there is more than one way to define the greatest common divisor. We will define it here in \cref{defn:gcd} using the well-ordering principle, but an alternative definition using Euclid's algorithm is of course just as good, since both definitions satisfy the specification that uniquely characterizes it. Hence we make the following specification of the greatest common divisor.

\begin{defn}\label{defn:is-gcd}
  Consider three natural numbers $a$, $b$, and $d$. We say that $d$ is a \define{greatest common divisor}\index{is a greatest common divisor|textbf}\index{greatest common divisor|textbf} of $a$ and $b$ if it comes equipped with an element of type\index{is-gcd a b d@{$\isgcd_{a,b}(d)$}|textbf}
  \begin{equation*}
    \isgcd_{a,b}(d) \defeq \prd{x:\N} (x\mid a)\times (x\mid b)\leftrightarrow (x\mid d).
  \end{equation*}
\end{defn}

The property of being a greatest common divisor uniquely characterizes the greatest common divisor, in the following sense.

\begin{prp}
  Suppose $d$ and $d'$ are both a greatest common divisor of $a$ and $b$. Then $d=d'$.
\end{prp}

\begin{proof}
  If both $d$ and $d'$ are a greatest common divisor of $a$ and $b$, then both $d$ and $d'$ divide both $a$ and $b$, and hence it follows that $d\mid d'$ and $d'\mid d$. Since the divisibility relation was shown to be a partial order in \cref{ex:is-poset-div}, it follows by antisymmetry that $d=d'$.
\end{proof}

Note that for any two natural numbers $a$ and $b$, the type
\begin{equation*}\label{eq:multiples-of-gcd}
  \sm{n:\N}\prd{x:\N} (x\mid a)\times (x\mid b)\to (x\mid n)\tag{\textasteriskcentered}
\end{equation*}
consists of all the multiples of the common divisors of $a$ and $b$, including $0$. On the other hand, the type
\begin{equation*}\label{eq:common-divisors}
  \sm{n:\N}\prd{x:\N} (x\mid n)\to (x\mid a)\times (x\mid b)\tag{\textasteriskcentered\textasteriskcentered}
\end{equation*}
consists of all the common divisors of $a$ and $b$ except in the case where $a=0$ and $b=0$. In this case, the type in \cref{eq:common-divisors} consists of all natural numbers.

\cref{eq:multiples-of-gcd,eq:common-divisors} provide us with two ways to define the greatest common divisor. We can either define the greatest common divisor of $a$ and $b$ as the greatest natural number in the type in \cref{eq:common-divisors} or we can define it as the least \emph{nonzero} natural number of the type in \cref{eq:multiples-of-gcd}, provided that we make an exception in the case where both $a=0$ and $b=0$. Since we already have established the well-ordering principle of $\N$, we will opt for the second approach. In \cref{ex:maximal-element} you will be asked to show that any \emph{bounded} decidable family over $\N$ has a maximum as soon as it contains some natural number. 

In order to correctly define the greatest common divisor using well-ordering principle of $\N$, we need a slight modification of the type family in \cref{eq:multiples-of-gcd}. We define this family as follows:

\begin{defn}\label{defn:fam-gcd}
  Given $a,b:\N$, we define the type family $\ismultipleofgcd(a,b)$ over $\N$ by
  \begin{align*}
    \ismultipleofgcd(a,b,n) & \defeq (a+b\neq 0) \to (n\neq 0)\times \Big(\prd{x:\N} (x\mid a)\times (x\mid b) \to (x\mid n)\Big).
  \end{align*}
\end{defn}

In other words, if $a+b=0$ then the type $\sm{n:\N}M(a,b,n)$ consist of all the natural numbers. On the other hand, if $a+b\neq 0$ it consists of the nonzero natural numbers $n$ with the property that any common divisor of $a$ and $b$ also divides $n$. These are exactly the nonzero multiples of the greatest common divisor of $a$ and $b$.

Since we intend to apply the well-ordering principle, we must show that the family $\ismultipleofgcd(a,b)$ is decidable. This is a step that one can skip in classical mathematics, because all the subsets of $\N$ are decidable there. However, in our current setting we have no choice but to prove it.

\begin{prp}\label{prp:is-decidable-is-multiple-of-gcd}
  The type family $\ismultipleofgcd(a,b)$ is decidable for each $a,b:\N$.
\end{prp}

\begin{proof}
  The type $a+b\neq 0$ is decidable because it is the negation of the type $a+b=0$, which is decidable by \cref{prp:has-decidable-equality-N}. Therefore it suffices to show that the type
  \begin{equation*}
    (n\neq 0)\times \prd{x:\N} (x\mid a)\times (x\mid b)\to (x\mid n)
  \end{equation*}
  is decidable, and by \cref{prp:is-decidable-function-type} we also get to assume that $a+b\neq 0$. The type $n\neq 0$ is again decidable by \cref{prp:has-decidable-equality-N}, so it suffices to show that the type
  \begin{equation*}
    \prd{x:\N}(x\mid a)\times (x\mid b)\to (x\mid n)
  \end{equation*}
  is decidable. The types $(x\mid a)\times(x\mid b)$ and $(x\mid n)$ are decidable by \cref{thm:is-decidable-div-N}, so by \cref{cor:is-decidable-bounded-pi} it suffices to check that the family of types $(x\mid a)\times (x\mid b)$ indexed by $x:\N$ has an upper bound. If $x$ is a common divisor of $a$ and $b$, then it follows that $x$ divides $a+b$. Furthermore, since we have assumed that $a+b\neq 0$, it follows that $x\leq a+b$. This provides the upper bound.
  \end{proof}

We are almost in position to apply the well-ordering principle of $\N$ to define the greatest common divisor. It just remains to show that there is some $n:\N$ for which $M(a,b,n)$ holds. We prove this in the following lemma.

\begin{lem}\label{lem:exists-multiple-of-gcd}
  There is an element of type $\ismultipleofgcd(a,b,a+b)$. 
\end{lem}

\begin{proof}
  To construct an element of type $\ismultipleofgcd(a,b,a+b)$, assume that $a+b\neq 0$. Then we have tautologically that $a+b\neq 0$, and any common divisor of $a$ and $b$ is also a divisor of $a+b$.
\end{proof}

\begin{defn}\label{defn:gcd}
  We define the \define{greatest common divisor}\index{greatest common divisor|textbf} $\gcd:\N\to (\N\to\N)$\index{gcd@{$\gcd$}|textbf} by the well-ordering principle of $\N$ (\cref{thm:well-ordering-principle-N}) as the least natural number $n$ for which $M(a,b,n)$ holds, using the fact that $M(a,b)$ is a decidable type family (\cref{prp:is-decidable-is-multiple-of-gcd}) and that $M(a,b,a+b)$ always holds (\cref{lem:exists-multiple-of-gcd}).
\end{defn}

\begin{lem}\label{lem:is-zero-gcd}
  For any two natural numbers $a$ and $b$, we have $\gcd(a,b)=0$ if and only if $a+b=0$.
\end{lem}

\begin{proof}
  To prove the forward direction, assume that $\gcd(a,b)=0$. By definition of $\gcd(a,b)$ we have that $\ismultipleofgcd(a,b,\gcd(a,b))$ holds. More explicitly, the implication
  \begin{equation*}
    (a+b\neq 0)\to (\gcd(a,b)\neq 0)\times \prd{x:\N}(x\mid a)\times(x\mid b)\to (x\mid\gcd(a,b))
  \end{equation*}
  holds. However, we have assumed that $\gcd(a,b)=0$, so it follows from the above implication that $\neg(a+b\neq 0)$. In other words, we have $\neg\neg(a+b=0)$. The fact that equality on $\N$ is decidable implies via \cref{ex:dne-is-decidable} that $\neg\neg(a+b=0)\to (a+b=0)$, so we conclude that $a+b=0$.

  For the converse direction, recall that the inequality $\gcd(a,b)\leq a+b$ holds by minimality, since $\ismultipleofgcd(a,b,a+b)$ holds by \cref{lem:exists-multiple-of-gcd}. If $a+b=0$, it therefore follows that $\gcd(a,b)\leq 0$, which implies that $\gcd(a,b)=0$.
\end{proof}

\begin{thm}
  For any two natural numbers $a$ and $b$, the number $\gcd(a,b)$ is a greatest common divisor of $a$ and $b$ in the sense of \cref{defn:is-gcd}.
\end{thm}

\begin{proof}
  We give the proof by case analysis on whether $a+b=0$.
  If we assume that $a+b=0$, then it follows that both $a=0$ and $b=0$, and by \cref{lem:is-zero-gcd} it also follows that $\gcd(a,b)=0$. Since any number divides $0$, the claim follows immediately.

  In the case where $a+b\neq 0$, it follows from \cref{lem:is-zero-gcd} that also $\gcd(a,b)\neq 0$. From the fact that $\ismultipleofgcd(a,b,\gcd(a,b))$ we therefore immediately obtain that
  \begin{equation*}
    \prd{x:\N} (x\mid a)\times (x\mid b)\to (x\mid \gcd(a,b)).
  \end{equation*}
  Therefore it remains to show that if $x$ divides $\gcd(a,b)$, then $x$ divides both $a$ and $b$. By transitivity of the divisibility relation it suffices to show that $\gcd(a,b)$ divides both $a$ and $b$. We will show only that $\gcd(a,b)$ divides $a$, the proof that $\gcd(a,b)$ divides $b$ is similar.

  Since $\gcd(a,b)$ is nonzero, it follows by Euclidean division (\cref{ex:euclidean-division}) that there are numbers $q$ and $r<\gcd(a,b)$ such that
  \begin{equation*}
    a = q\cdot\gcd(a,b)+r.
  \end{equation*}
  From this equation and \cref{prp:div-3-for-2} it follows that any number $x$ which divides both $a$ and $b$ also divides $r$, because we have already noted that any such $x$ divides $\gcd(a,b)$. This observation implies that $r=0$, because we have $r<\gcd(a,b)$ by construction and $\gcd(a,b)$ is minimal. Therefore we conclude that $\gcd(a,b)$ divides $a$.
\end{proof}
\index{natural numbers!greatest common divisor|)}
\index{greatest common divisor|)}

\subsection{The infinitude of primes}
\index{prime number|(}
\index{natural number!prime number|(}

When the natural numbers are ordered by the divisibility relation, the number $1$ is at the bottom. Directly above $1$ are the prime numbers. Above the prime numbers are the multiples of two primes, then the multiples of three primes, and so on. At the top of this ordering we find $0$. For any natural number $n$, the numbers strictly below $n$ are the proper divisors of $n$. A prime number is therefore a number of which has exactly one proper divisor.

\begin{defn}
  ~
  \begin{enumerate}
  \item Consider two natural numbers $d$ and $n$. Then $d$ is said to be a \define{proper divisor}\index{proper divisor|textbf} of $n$ if it comes equipped with an element of type\index{is-proper-divisor n d@{$\isproperdivisor(n,d)$}|textbf}
    \begin{equation*}
      \isproperdivisor(n,d)\defeq (d\neq n)\times (d\mid n).
    \end{equation*}
  \item A natural number $n$ is said to be \define{prime}\index{prime number|textbf}\index{natural numbers!prime number|textbf} if it comes equipped with an element of type\index{is-prime(n)@{$\isprime(n)$}|textbf}
  \begin{equation*}
    \isprime(n)\defeq \prd{x:\N}\isproperdivisor(n,x)\leftrightarrow (x=1).
  \end{equation*}
  \end{enumerate}
\end{defn}

\begin{prp}
  For any $n:\N$, the type $\isprime(n)$ is decidable.\index{is-prime(n)@{$\isprime(n)$}!is decidable}
\end{prp}

\begin{proof}
  We will first show that $\isprime(n)\leftrightarrow\isprime'(n)$, where
  \begin{equation*}
    \isprime'(n)\defeq (n\neq 1)\times \prd{x:\N}\isproperdivisor(n,x)\to (x=1).
  \end{equation*}
  For the forward direction, simply note that $1$ is not a proper divisor of itself, and therefore $1$ is not a prime. For the converse direction, suppose that $n\neq 1$ and that any proper divisor of $n$ is $1$. Then it follows that $1$ is a proper divisor of $n$, which implies that $n$ is prime.

  Now we proceed by showing that the type $\isprime'(n)$ is decidable for every $n:\N$. The proof is by case analysis on whether $n=0$ or $n\neq 0$. In the case where $n=0$, note that any nonzero number is a proper divisor of $0$, and therefore $\isprime'(0)$ doesn't hold. In particular, $\isprime'(0)$ is decidable.
  
  Now suppose that $n\neq 0$. In order to show that the type $\isprime'(n)$ is decidable, note that the type $n\neq 1$ is decidable since it is the negation of the decidable type $n=1$. Therefore it suffices to show that the type
  \begin{equation*}
    \prd{x:\N}\isproperdivisor(n,x)\to (x=1)
  \end{equation*}
  is decidable. Since the types $(x\neq n)\times (x\mid n)$ and $x=1$ are decidable, it follows from \cref{cor:is-decidable-bounded-pi} that it suffices to check that
  \begin{equation*}
    ((x\neq n)\times (x\mid n))\to (x\leq n)
  \end{equation*}
  for any $x:\N$. This follows from the implication $(x\mid n)\to (x\leq n)$, which holds because we have assumed that $n\neq 0$.
\end{proof}

The proof that there are infinitely many primes proceeds by constructing a prime number larger than $n$, for any $n:\N$. The number $n!+1$ is relatively prime with any number $x\leq n$. Therefore there is a least number $n<m$ that is relatively prime with any number $x\leq n$, and it follows that this number $m$ must be prime.

\begin{defn}
  For any two natural numbers $n$ and $m$, we define the type
  \begin{equation*}
    R(n,m)\defeq (n<m)\times \prd{x:\N}(x\leq n)\to ((x\mid m)\to (x=1)). 
  \end{equation*}
\end{defn}

\begin{lem}
  The type $R(n,m)$ is decidable for each $n,m:\N$.
\end{lem}

\begin{proof}
  The type $n<m$ and, and for each $x:\N$ both types $x\leq n$ and $(x\mid m)\to (x=1)$ are decidable, so it follows via \cref{cor:is-decidable-bounded-pi} that the product
  \begin{equation*}
    \prd{x:\N}(x\leq n)\to ((x\mid m)\to (x=1))
  \end{equation*}
  is decidable.
\end{proof}

\begin{lem}\label{lem:succ-factorial-has-one-bounded-divisor}
  There is an element of type $R(n,{n!}+1)$ for each $n:\N$.
\end{lem}

\begin{proof}
  The fact that $n<{n!}+1$ follows from the fact that $n\leq n!$, which is shown by induction. We leave this to the reader, and focus on the second aspect of the claim: that every $x\leq n$ that divides ${n!}+1$ must be equal to $1$.

  To see this, note that any divisor of ${n!}+1$ is automatically nonzero, and recall that any nonzero $x\leq n$ divides $n!$ by \cref{ex:div-factorial}. Therefore it follows that any $x\leq n$ that divides ${n!}+1$ also divides $n!$, and consequently it divides $1$ as well. Now we are done, because if $x$ divides $1$ then $x=1$.
\end{proof}

We finally show that there are infinitely many primes.

\begin{thm}
  For each $n:\N$, there is a prime number $p:\N$ such that $n< p$.\index{infinitude of primes}\index{prime number!infinitude of primes}\index{natural numbers!infinitude of primes}
\end{thm}

\begin{proof}
  It suffices to show that for each \emph{nonzero} $n:\N$, there is a prime number $p:\N$ such that $n\leq p$. Let $n$ be a nonzero natural number.

  Since the type $R(n,m)$ is decidable for each $m:\N$, and since $R(n,{n!}+1)$ holds by \cref{lem:succ-factorial-has-one-bounded-divisor}, it follows by the well-ordering principle of $\N$ (\cref{thm:well-ordering-principle-N}) that there is a minimal $m:\N$ such that $R(n,m)$ holds. In order to prove the theorem, we will show that this number $m$ is prime, i.e., that there is an element of type
  \begin{equation*}
    (m\neq 1)\times \prd{x:\N} \isproperdivisor(m,x)\to (x=1).
  \end{equation*}

  First, we note that $m\neq 1$ because $n<m$ holds by construction, and $n$ is assumed to be nonzero. Therefore it suffices to show that $1$ is the only proper divisor of $m$. Let $x$ be a proper divisor of $m$. Since $R(n,m)$ holds by construction, we will prove that $x=1$ by showing that $x\leq n$ holds.

  Since $m$ is nonzero, it follows from the assumption that $x\mid m$ that $x<m$. By minimality of $m$, it therefore follows that $\neg R(n,x)$ holds. However, any divisor of $x$ is also a divisor of $m$ by transitivity of the divisibility relation. Therefore it follows that any $y\leq n$ that divides $x$ must be $1$. In other words:
  \begin{equation*}
    \prd{y:\N}(y\leq n)\to ((y\mid x)\to (y=1))
  \end{equation*}
  holds. Since $\neg R(n,x)$ holds, we conclude now that $n\nless x$. To finish the proof, it follows that $x\leq n$.
\end{proof}
\index{prime number|)}
\index{natural number!prime number|)}

\subsection{Boolean reflection}\label{sec:boolean-reflection}
\index{boolean reflection|(}

We have shown that the type $\isprime(n)$ is decidable for every $n$. In other words, there is an element $d(n):\isdecidable(\isprime(n))$ for every $n$. In principle, we can therefore check whether any \emph{specific} natural number $n$ is prime by inspecting the element $d(n)$: if it is of the form $\inl(x)$ for some $x:\isprime(n)$, then $n$ is prime; if it is of the form $\inr(f)$ for some $f:\neg\isprime(n)$, then $n$ is not prime. In other words, we evaluate the element $d(n)$ using the computation rules of type theory, and then we see whether $n$ is prime or not.

Computers can perform such evaluations, but it is often unfeasible to carry out such evaluations by hand. Moreover, even for computers the task of evaluating a proof term like $\isdecidableisprime(n)$ may quickly get out of hand. With the formalization of the material in this book, the proof assistant Agda returns a proof term of 430 lines of code when we simply ask it to evaluate the term $\isdecidableisprime(7)$, and it returned a proof term of 69373 lines of code when we asked it to evaluate the term $\isdecidableisprime(37)$. There is a much better way to do this: \emph{boolean reflection}.

\begin{defn}
  For any type $A$ we define the map\index{booleanization@{$\booleanization$}|textbf}
  \begin{equation*}
    \booleanization:\isdecidable(A)\to\bool
  \end{equation*}
  by
  \begin{align*}
    \booleanization(\inl(a)) & \defeq \btrue \\*
    \booleanization(\inr(f)) & \defeq \bfalse.
  \end{align*}
\end{defn}

\begin{thm}[Boolean reflection principle]
  For any type $A$ and any decision $d:\isdecidable(A)$, there is a map\index{boolean reflection|textbf}\index{reflect@{$\booleanreflection$}|textbf}
  \begin{equation*}
    \booleanreflection:(\booleanization(d)=\btrue)\to A
  \end{equation*}
  such that $\booleanreflection(\inl(a))\jdeq a$.
\end{thm}

\begin{proof}
  First, recall that by \cref{ex:obs_bool} there is a map $\gamma:(\bfalse=\btrue)\to \emptyt$. We use this to construct $\booleanreflection$ by pattern matching as follows:
  \begin{align*}
    \booleanreflection(\inl(a),p) & \defeq a \\*
    \booleanreflection(\inr(f),p) & \defeq \exfalso(\gamma(p)).\qedhere
  \end{align*}
\end{proof}

\begin{rmk}
  Since the number 37 is a prime, it follows that the booleanization of the term
  \begin{equation*}
    d(37):\isdecidable(\isprime(37))
  \end{equation*}
  has the value $\booleanization(d(37))\jdeq\btrue$. By boolean reflection it therefore follows that
  \begin{equation}\label{eq:is-prime-37}
    \isprimethirtyseven\defeq \booleanreflection(d(37),\refl{}):\isprime(37).\tag{\textasteriskcentered}
  \end{equation}
  The term in $\isprimethirtyseven$ does not, however, contain any explicit information as to why the number 37 is prime. The reason that it type checks is simply that $d(37)$ is judgmentally equal to some term of the form $\inl(t):\isdecidable(\isprime(37))$ and therefore it follows that $\refl{}$ is an identification of type
  \begin{equation*}
    \booleanization(d(37))=\btrue.
  \end{equation*}
  To see that $\isprimethirtyseven$ is indeed an element of type $\isprime(37)$ therefore requires us to evaluate the term $d(37)$. This is not doable by hand. Computer proof assistants, however, are capable of performing this task. In a proof assistant, we may therefore use boolean reflection to offload the task of evaluating the decision algorithm of a decidable type to the computer. This technique has been essential in the formalization of the Feit-Thompson theorem in Coq \cite{Gonthier}. The book \emph{Mathematical Components} \cite{mathematical-components} contains more information about using boolean reflection effectively in formalized mathematics.

  Do not, however, "solve" your homework problems with boolean reflection. If your teaching assistant cannot evaluate your solution, they will conclude that you haven't demonstrated your clear understanding of the problem.
\end{rmk}
\index{boolean reflection|)}

\begin{exercises}
  \exitem
  \begin{subexenum}
  \item State Goldbach's conjecture\index{Goldbach's conjecture} in type theory.
  \item State the twin prime conjecture\index{twin prime conjecture} in type theory.
  \item State the Collatz conjecture\index{Collatz conjecture} in type theory.
  \end{subexenum}
  \noindent If you have a solution to any of these open problems, you should certainly formalize it before you submit it to the Annals of Mathematics.
  \exitem Show that
  \begin{equation*}
    \isdecidable(\isdecidable(P))\to\isdecidable(P)
  \end{equation*}
  for any type $P$.
  \exitem For any family $P$ of decidable types indexed by $\Fin{k}$, construct a function
  \begin{equation*}
    \neg\Big(\prd{x:\Fin{k}}P(x)\Big)\to\sm{x:\Fin{k}}\neg P(x).
  \end{equation*}
  \exitem
  \begin{subexenum}
  \item Define the \define{prime function} $\primefunction:\N\to\N$\index{prime@{$\primefunction$}|see {prime function}}\index{prime function|textbf}\index{natural numbers!prime function} for which $\primefunction(n)$ is the $n$-th prime.

  \item Define the \define{prime-counting function}\index{prime counting function|textbf} $\pi:\N\to\N$\index{p@{$\pi$}|see {prime counting function}}, which counts for each $n:\N$ the number of primes $p\leq n$.
  \end{subexenum}
  \exitem For any natural number $n$, show that
  \begin{equation*}
    \isprime(n)\leftrightarrow (2\leq n)\times \prd{x:\N} (x\mid n)\to (x=1)+(x=n).
  \end{equation*}
  \exitem Consider two types $A$ and $B$. Show that the following are equivalent:
  \begin{enumerate}
  \item There are functions
    \begin{align*}
      & B \to \hasdecidableequality(A) \\
      & A \to \hasdecidableequality(B).
    \end{align*}
  \item The product $A\times B$ has decidable equality.\index{decidable equality!of cartesian products}
  \end{enumerate}
  Conclude that if both $A$ and $B$ have decidable equality, then so does $A\times B$.
  \exitem Consider two types $A$ and $B$, and consider the observational equality $\Eqcoprod$ on the coproduct $A+B$ defined by
  \begin{align*}
    \Eqcoprod(\inl(x),\inl(x')) & \defeq x= x' & \Eqcoprod(\inl(x),\inr(y')) & \defeq \emptyt \\
    \Eqcoprod(\inr(y),\inl(x')) & \defeq \emptyt & \Eqcoprod(\inr(y),\inr(y')) & \defeq y = y'.
  \end{align*}
  \begin{subexenum}
  \item Show that $(x=y)\leftrightarrow\Eqcoprod(x,y)$ for every $x,y:A+B$.
  \item Show that the following are equivalent:
    \begin{enumerate}
    \item Both $A$ and $B$ have decidable equality.\index{decidable equality!of coproducts}
    \item The coproduct $A+B$ has decidable equality.
    \end{enumerate}
    Conclude that $\Z$ has decidable equality.\index{decidable equality!of Z@{of $\Z$}}\index{Z@{$\Z$}!has decidable equality}\index{integers!decidable equality}\index{has decidable equality!integers}
  \end{subexenum}
  \exitem \label{ex:has-decidable-equality-Sigma}Consider a family $B$ over $A$, and consider the following three conditions:
  \begin{enumerate}
  \item The type $A$ has decidable equality.
  \item The type $B(x)$ has decidable equality for each $x:A$.
  \item The type $\sm{x:A}B(x)$ has decidable equality.\index{decidable equality!of S-types@{of $\Sigma$-types}}
  \end{enumerate}
  Show that if (i) holds, then (ii) and (iii) are equivalent, and show that if $B$ has a section $b:\prd{x:A}B(x)$, then (ii) and (iii) together imply (i).
  \exitem Consider a family $B$ of types over $\Fin{k}$, for some $k:\N$.
  \begin{subexenum}
  \item Show that if each $B(x)$ is decidable, then $\prd{x:\Fin{k}}B(x)$ is again decidable.
  \item Show that if each $B(x)$ has decidable equality, then $\prd{x:\Fin{k}}B(x)$ also has decidable equality.
  \end{subexenum}
  \exitem \label{ex:maximal-element}Consider a decidable type family $P$ over $\N$ equipped with an upper bound $m$.
  \begin{subexenum}
  \item Show that the type $\sm{n:\N}P(n)$ is decidable.
  \item Construct a function
  \begin{equation*}
    \Big(\sm{n:\N}P(n)\Big)\to\Big(\sm{n:\N}P(n)\times\isupperbound_P(n)\Big).
  \end{equation*}
  \item Use the function of part (b) to give a second construction of the greatest common divisor, and verify that it satisfies the specification of \cref{defn:is-gcd}.
  \end{subexenum}
  \exitem \label{ex:bezouts-identity-N}
  \begin{subexenum}
  \item For any three natural numbers $x$, $y$, and $z$, show that the type
    \begin{equation*}
      \sm{k:\N}\sm{l:\N}\distN(kx,ly)=z
    \end{equation*}
    is decidable.
  \item (B\'ezout's identity)\index{Bezout's identity@{B\'ezout's identity}}\index{natural numbers!Bezout's identity@{B\'ezout's identity}} For any two natural numbers $x$ and $y$, construct two natural numbers $k$ and $l$ equipped with an identification
  \begin{equation*}
    \distN(kx,ly)=\gcd(x,y).
  \end{equation*}
  \end{subexenum}
  \exitem
  \begin{subexenum}
  \item Show that every natural number $n\geq 2$ has a prime factor.
  \item Define a function
    \begin{equation*}
      \primefactors : \Big(\sm{n:\N}2\leq n\Big)\to \lst(\N)
    \end{equation*}
    such that $\primefactors(n)$ is an increasing list of primes, and $n$ is the product of the primes in the list $\primefactors(n)$.
  \item Show that any increasing list $l$ of primes of which the product is $n$ is equal to the list $\primefactors(n)$.
  \end{subexenum}
  \exitem Show that there are infinitely many primes $p\equiv 3\mod 4$.
  \exitem Show that for each prime $p$, the ring $\Z/p$ of integers modulo $p$ is a field, i.e., construct a multiplicative inverse
  \begin{equation*}
    (\blank)^{-1} : \prd{x:\Z/p}\to (x\neq 0) \to \Z/p
  \end{equation*}
  equipped with identifications
  \begin{align*}
    x^{-1}x & = 1 & xx^{-1} & = 1.
  \end{align*}
  \exitem Let $F:\N\to\N$ be the Fibonacci sequence. Construct the \define{cofibonacci sequence}\index{cofibonacci sequence|textbf}\index{natural numbers!cofibonacci sequence}, i.e., the function $G:\N\to\N$ such that\index{Fibonacci sequence!has left adjoint}
  \begin{equation*}
    (G_m\mid n) \leftrightarrow (m\mid F_n)
  \end{equation*}
  for all $m,n:\N$. Hint: for $m>0$, $G_m$ is the least $x>0$ such that $m\mid F_x$. 
\end{exercises}



\cleardoublepage

\chapter{The Univalent Foundations of Mathematics}\label{chap:uf}

The univalent foundations program is an approach to mathematics in which mathematics is formalized in dependent type theory, using the homotopy interpretation and the univalence axiom. The homotopy interpretation of type theory fully embraces the idea that between any two elements of a type there is a \emph{type} of identifications, much like between any two points in a topological space there is a \emph{space} of paths between them. This idea was first explored by Awodey and Warren in \cite{AwodeyWarren} and independently by Voevodsky in \cite{Voevodsky06}. With the homotopy interpretation of type theory\index{homotopy interpretation}, outlined in the table below, we think of types as spaces, type families as fibrations, and identifications as paths.
\begin{table}
\begin{tabular}{ll}
\toprule
\emph{Type theory} &  \emph{Homotopy theory} \\
\midrule
Types  & Spaces \\
Dependent types & Fibrations \\
Elements & Points \\
Dependent pair type & Total space \\
Identity type & Path fibration\\
\bottomrule
\end{tabular}
\end{table}

Voevodsky's univalence axiom characterizes the identity type of the universes in type theory, asserting that for any two types $A$ and $B$ in a universe $\UU$, we have an equivalence
\begin{equation*}
  (A=_{\UU}B)\simeq (A\simeq B).
\end{equation*}
In other words, identifications of types are equivalent to equivalences of types. A consequence of the univalence axiom is that many kinds of isomorphic objects in mathematics, such as isomorphic groups or isomorphic rings, can be identified.

The concept of equivalences generalizes the concept of set-isomorphisms to type theory in a way that is suitable for the homotopy interpretation of type theory. Equivalent types are the same for all practical purposes, just as isomorphic objects are practically the same in everyday mathematics. By the univalence axiom, isomorphic objects get identified.

However, the informal practice of identifying isomorphic objects is technically inconsistent with the set theoretic foundations of mathematics. The extensionality axiom of Zermelo-Fraenkel set theory implies, for instance, that there are many different singleton sets $\{x\}$. All those singleton sets are isomorphic, so the univalence axiom identifies them, which would be inconsistent within Zermelo-Fraenkel set theory. The assumption of the univalence axiom therefore marks our definitive departure from the set-theoretic foundations of mathematics.

Since the univalence axiom characterizes the identity type of the universe, it is important to understand the general task of characterizing the identity type of any given type. It is a crucial observation, which we already made when we discussed the uniqueness of $\refl{}$ in \cref{sec:refl-unique}, that for any $a:A$, the type
\begin{equation*}
  \sm{x:A}a=x  
\end{equation*}
is contractible. Contractible types are types that are singletons up to homotopy, i.e., they are types $A$ that come equipped with a point $a:A$ such that $a=x$ for every $x:A$. We have seen in \cref{prp:contraction-total-space-id} that the total space of all paths starting at $a$ is such a type, so it is an example of a contractible type. The fundamental theorem of identity types asserts that a type family $B$ over $A$ with $b:B(a)$ has a contractible total space
\begin{equation*}
  \sm{x:A}B(x)
\end{equation*}
if and only if $(a=x)\simeq B(x)$ for all $x:A$. The fundamental theorem of identity types can be used to characterize the identity types of virtually any type that we will encounter. Since types are only fully understood if we also have a clear understanding of their identity types, it is one of the core tasks of a homotopy type theorist to characterize identity types, and the fundamental theorem (\cref{thm:id_fundamental}) is the main tool.

Not all types have very complicated identity types. For example, some types have the property that all their identity types are contractible. For example, the types $\emptyt$ and $\unit$ satisfy this condition. Any two terms of such a type can therefore be identified, so in this sense they are \emph{proof irrelevant}. The only thing that matters about such types is whether or not they are inhabited by a term. Analogously, this is also the case for propositions in the propositional calculus or first order logic. Therefore we call such types propositions, and we see that propositions are present in type theory as certain types.

Next, there are the types of which the identity types are propositions. In other words, the identity types of such types have the property of proof irrelevance. We are familiar with this situation from set theory, because equality in set theory is a proposition. Therefore we call such types sets. The types $\N$, $\Z$, and $\Fin{k}$ are all examples of sets.

It is now clear that there is a hierarchy arising: at the bottom of the hierarchy we have the contractible types; then we have the propositions, of which the identity types are contractible; after the propositions we have the sets, of which the identity types are propositions. Defining sets to be of truncation level $0$, we define a type to be of truncation level $k+1$ if its identity types are of truncation level $k$. Types of truncation level $k$ for $k\geq 1$ are also called $k$-types or $k$-groupoids.

This hierarchy of truncation levels is due to Voevodsky, who recognized that, when you are formalizing mathematics in type theory, it is important to specify the truncation level in which you are working. Most mathematics, for example, takes place at truncation level $0$, the level of sets. Groups, rings, posets, and so on are all set-level objects. Categories, on the other hand, are objects of truncation level $1$, the level of the $1$-groupoids. This is because two objects in a category are considered equal if they are isomorphic, and between any two objects in a category there is a set of isomorphisms.

The fundamental theorem of identity types and the basic facts about truncation levels are proved without assuming any axioms. In other words, they are theorems of Martin-L\"of's dependent type theory, as introduced in \cref{chap:type-theory}. In particular, the rules of dependent type theory are sufficient to characterize the identity types of $\Sigma$-types and of the type of natural numbers, and also to prove the disjointness of coproducts. However, there are still two important characterizations of identity types missing: those of $\Pi$-types and those of universes. For those two cases we need axioms:
\begin{enumerate}
\item For any two dependent functions $f,g:\prd{x:A}B(x)$, the canonical map
  \begin{equation*}
    (f=g)\to (f\htpy g)
  \end{equation*}
  that maps $\refl{f}$ to the constant homotopy, is an equivalence.
\item For any two types $A$ and $B$ in a universe $\UU$, the canonical map
  \begin{equation*}
    (A=B)\to (A\simeq B)
  \end{equation*}
  that maps $\refl{A}$ to the identity equivalence, is an equivalence.
\end{enumerate}
The function extensionality axiom (i) characterizes the identity types of $\Pi$-types, and the univalence axiom (ii) characterizes the identity types of universes.

With the addition of the function extensionality axiom and the univalence axiom, we have almost fully specified the univalent foundations of mathematics. The one ingredient missing is that of quotients. In order to obtain quotients, we will postulate two more axioms:
\begin{enumerate}
  \addtocounter{enumi}{2}
\item We will assume that every type $A$ has a propositional truncation.
\item We will assume the type theoretic \emph{replacement axiom}. 
\end{enumerate}

Propositional truncations are the simplest kind of quotients, identifying all the elements in a type $A$. In other words, the propositional truncation of a type $A$ is a proposition $\brck{A}$ that is true if and only if $A$ is inhabited. Using propositional truncations we can construct the homotopy image of a map. A quotient of a type $A$ by an equivalence relation $R$ can then be constructed as the type of all equivalence classes of $R$, i.e., as the image of the map $R:A\to (A\to \prop_\UU)$. With this construction of the quotient, we immediately obtain a surjective map $q:A\to A/R$, and by the univalence axiom it follows that the quotient is \emph{effective}, i.e., that for any $x,y:A$ we have
\begin{equation*}
  (q(x)=q(y))\simeq R(x,y).
\end{equation*}
However, this construction does not guarantee that the quotient $A/R$ is small with respect to the universe $\UU$, because it is constructed as a subtype of the type $A\to\prop_\UU$. This is why we assume the replacement axiom, which will imply that the quotient $A/R$ is \emph{essentially} small. Essentially small types are types that are equivalent to a small type, and the replacement axiom asserts that if $f:A\to B$ is a map from an essentially small type $A$ into a type $B$ of which the identity types are essentially small, then the image of $f$ is also essentially small. The role of the replacement axiom in type theory is similar to the role of the replacement axiom in Zermelo-Fraenkel set theory: to ensure that quotients are small.

We have two goals in this chapter. The first goal is to fully describe the univalent foundations of mathematics and its most important concepts. Our second goal is to show how to we can start doing ordinary mathematics from a univalent point of view. We therefore show how to derive the strong induction principle for the natural numbers using function extensionality; we give a new interpretation of logic in univalent mathematics using our definition of propositions and the propositional truncations; we show how Cantor's diagonal argument works in univalent mathematics; and we introduce finite types, binomial types, set quotients, the univalent type of all groups. We end this chapter with a variant of Russell's paradox, showing that for any univalent universe $\UU$ there can be no type $U:\UU$ that is equivalent to $\UU$. We hope that after seeing these familiar examples, you will be able to do your own mathematics from a univalent point of view.

\section{Equivalences}\label{sec:equivalences}
\index{equivalence|(} 

In this section we will define equivalences of types. However, we have to be a bit careful in how we define the condition for a map to be an equivalence. It turns out to be important that being an equivalence is a \emph{property} of maps, and not a \emph{structure} on maps. In other words, we want to define the type
\begin{equation*}
  \isequiv(f)
\end{equation*}
in such a way that we will be able to prove that the type $\isequiv(f)$ is a \emph{proposition}. Propositions will be defined in \cref{chap:hierarchy}, and in \cref{chap:funext} we will be able to prove that $\isequiv(f)$ is indeed a proposition.

It turns out that if we naively define a function $f$ to be an equivalence if it has an inverse, then we won't be able to show that $\isequiv(f)$ is a property. We will therefore say that $f$ is an equivalence if it has a separate left and right inverse. This may look odd, but when we define equivalences in this way we will be able to show that $\isequiv(f)$ is a property.

\subsection{Homotopies}
\index{homotopy|(}

In type theory we are very limited in constructing identifications of functions. The following example illustrates a case where type theory provides no rules to construct an identification between two maps, even though they are pointwise equal.

\begin{rmk}\label{rmk:negnegbool}
  Consider the negation function $\negbool : \bool\to\bool$\index{neg bool@{$\negbool$}} on the booleans, which was defined in \cref{ex:bool}. Type theory does not provide any means to show that
  \begin{equation*}
    \negbool\circ\negbool=\idfunc.
  \end{equation*}
  The best we can do is to construct an identification\index{neg neg bool@{$\negnegbool$}|textbf}
  \begin{equation*}
    \negnegbool(b) : \negbool(\negbool(b))=b
  \end{equation*}
  for any $b:\bool$. Indeed, $\negnegbool$ is defined using the induction principle of $\bool$, by
  \begin{align*}
    \negnegbool(\btrue) & \defeq \refl{\btrue} \\
    \negnegbool(\bfalse) & \defeq \refl{\bfalse}.
  \end{align*}
  Therefore we see that, while we cannot identify $\negbool\circ\negbool$ with $\idfunc$, we can define a \emph{pointwise identification} between the values of $\negbool\circ\negbool$ and $\idfunc$.
\end{rmk}

The observations in \cref{rmk:negnegbool} are an instance of a general phenomenon in type theory: it is often much easier to construct a \emph{pointwise identification} between the values of two maps, than it is to construct an identification between those two maps. In fact, the prevalent notion of sameness of maps is the notion of pointwise identification. Since they are so important, we will give them a name and call them \emph{homotopies}\index{pointwise identification|see {homotopy}}\index{pointwise identification|textbf}.

\begin{defn}
Let $f,g:\prd{x:A}B(x)$ be two dependent functions. The type of \define{homotopies}\index{homotopy|textbf} from $f$ to $g$ is defined as the type of pointwise identifications, i.e., we define\index{f htpy g@{$f\htpy g$}|see {homotopy}}\index{f htpy g@{$f\htpy g$}|textbf}
\begin{equation*}
f\htpy g \defeq \prd{x:A} f(x)=g(x).
\end{equation*}
\end{defn}

\begin{eg}
  By \cref{rmk:negnegbool} we have a homotopy
  \begin{equation*}
    \negnegbool : \negbool\circ\negbool\htpy\idfunc.
  \end{equation*}
\end{eg}

\begin{rmk}\label{rmk:commuting-diagrams}
  We will use homotopies, for example, to express the commutativity of diagrams. For example, we say that a triangle\index{homotopy!commutative diagram|textbf}
  \begin{equation*}
    \begin{tikzcd}[column sep=tiny]
      A \arrow[rr,"h"] \arrow[dr,swap,"f"] & & B \arrow[dl,"g"] \\
      & X
    \end{tikzcd}
  \end{equation*}
  \define{commutes}\index{commutative diagram|textbf} if it comes equipped with a homotopy $H:f\htpy g\circ h$. Similarly, we say that a square
  \begin{equation*}
    \begin{tikzcd}
      A \arrow[r,"g"] \arrow[d,swap,"f"] & A' \arrow[d,"{f'}"] \\
      B \arrow[r,swap,"h"] & B'
    \end{tikzcd}
  \end{equation*}
  commutes if it comes equipped with a homotopy $h \circ f\htpy f'\circ g$.
\end{rmk}

Note that the type of homotopies $f\htpy g$ is defined for dependent functions, and moreover the type of homotopies is itself a dependent function type. The definition of homotopies is therefore set up in such a way that we may also consider homotopies \emph{between}\index{homotopy!iterated homotopy}\index{iterated homotopies} homotopies, and even further homotopies between those higher homotopies. More concretely, if $H,K:f\htpy g$ are two homotopies, then the type of homotopies $H\htpy K$ between them is just the type
\begin{equation*}
\prd{x:A} H(x)=K(x).
\end{equation*}

\index{groupoid laws!of homotopies|(}
\index{homotopy!groupoid laws|(}
Since homotopies are pointwise identifications, we can use the groupoidal structure of identity types to also define the groupoidal structure of homotopies. In this case, however, we state the groupoid laws as \emph{homotopies} and \emph{homotopies between homotopies} rather than as identifications.

\begin{defn}\label{defn:htpy_groupoid}\index{groupoid laws!of homotopies}
  For any type family $B$ over $A$ we define the operations on homotopies
  \index{homotopy!refl-htpy@{$\reflhtpy$}|textbf}
  \index{refl-htpy@{$\reflhtpy$}|textbf}
  \index{homotopy!inv-htpy@{$\invhtpy$}|textbf}
  \index{inv-htpy@{$\invhtpy$}|textbf}
  \index{homotopy!concathtpy@{$\concathtpy$}|textbf}
  \index{concat-htpy@{$\concathtpy$}|textbf}
  \begin{align*}
    \reflhtpy & : \prd{f:\prd{x:A}B(x)}f\htpy f \\
    \invhtpy & : \prd{f,g:\prd{x:A}B(x)} (f\htpy g)\to(g\htpy f) \\
    \concathtpy & : \prd{f,g,h:\prd{x:A}B(x)} (f\htpy g)\to ((g\htpy h)\to (f\htpy h))
  \end{align*}
  pointwise by
  \begin{align*}
    \reflhtpy(f) & \defeq \lam{x} \refl{f(x)} \\
    \invhtpy(H) & \defeq \lam{x} H(x)^{-1} \\
    \concathtpy(H,K) & \defeq \lam{x}\ct{H(x)}{K(x)}.
  \end{align*}
  We will often write $H^{-1}$ for $\invhtpy(H)$, and $\ct{H}{K}$ for $\concathtpy(H,K)$.
\end{defn}

\begin{prp}
  Homotopies satisfy the groupoid laws:
  \begin{enumerate}
  \item Concatenation of homotopies is associative up to homotopy, i.e., there is a homotopy
    \begin{equation*}
      \assochtpy(H,K,L) : \ct{(\ct{H}{K})}{L}\htpy\ct{H}{(\ct{K}{L})}
    \end{equation*}
    for any homotopies $H:f\htpy g$, $K:g\htpy h$ and $L:h\htpy i$.
  \item Homotopies satisfy the left and right unit laws up to homotopy, i.e., there are homotopies
    \begin{align*}
    \leftunithtpy(H) & : \ct{\reflhtpy_f}{H}\htpy H \\
    \rightunithtpy(H) & : \ct{H}{\reflhtpy_g}\htpy H 
    \end{align*}
    for any homotopy $H$.
  \item Homotopies satisfy the left and right inverse laws up to homotopy, i.e., there are homotopies
    \begin{align*}
      \leftinvhtpy(H) & : \ct{H^{-1}}{H} \htpy \reflhtpy_g \\
      \rightinvhtpy(H) & : \ct{H}{H^{-1}} \htpy \reflhtpy_f
    \end{align*}
    for any homotopy $H$.
  \end{enumerate}
\end{prp}

\begin{proof}
  The homotopy $\assochtpy(H,K,L)$ is defined pointwise by
  \begin{equation*}
    \assochtpy(H,K,L,x) \defeq \assoc(H(x),K(x),L(x)).
  \end{equation*}
  The other homotopies are similarly defined pointwise.
\end{proof}
\index{groupoid laws!of homotopies|)}
\index{homotopy!groupoid laws|)}

Apart from the groupoid operations and their laws, we will occasionally need \emph{whiskering} operations. Whiskering operations are operations that allow us to compose homotopies with functions. There are two situations where we want this:
\begin{equation*}
  \begin{tikzcd}
    A \arrow[r,bend left,""{name=A,below}] \arrow[r,bend right,""{name=B,above}] \arrow[from=A,to=B,draw=none,"\Downarrow" description] & B \arrow[r] & C & A \arrow[r] & B \arrow[r,bend left,""{name=C,below}] \arrow[r,bend right,""{name=D,above}] \arrow[from=C,to=D,draw=none,"\Downarrow" description] & C.
  \end{tikzcd}
\end{equation*}

\begin{defn}
We define the following \define{whiskering}\index{homotopy!whiskering operations|textbf}\index{whiskering operations!of homotopies|textbf} operations on homotopies:
\begin{enumerate}
\item Suppose $H:f\htpy g$ for two functions $f,g:A\to B$, and let $h:B\to C$. We define\index{h . H@{$h\cdot H$}|see {homotopy, whiskering operations}}\index{h . H@{$h\cdot H$}|textbf}
\begin{equation*}
h\cdot H\defeq \lam{x}\ap{h}{H(x)}:h\circ f\htpy h\circ g.
\end{equation*}
\item Suppose $f:A\to B$ and $H:g\htpy h$ for two functions $g,h:B\to C$. We define\index{H . f@{$H\cdot f$}|see {homotopy, whiskering operations}}\index{H . f@{$H\cdot f$}|textbf}
\begin{equation*}
H\cdot f\defeq\lam{x}H(f(x)):g\circ f\htpy h\circ f.
\end{equation*}
\end{enumerate}
\end{defn}
\index{homotopy|)}

\subsection{Bi-invertible maps}

We use homotopies to define sections and retractions of a map $f$, and to define what it means for a map $f$ to be an equivalence.

\begin{defn}
  Let $f:A\to B$ be a function.
  \begin{enumerate}
  \item The type of \define{sections} of $f$\index{section!of a map|textbf}\index{function!section of a map|textbf} is defined to be the type\index{sec(f)@{$\sections(f)$}|textbf}
    \begin{equation*}
      \sections(f) \defeq \sm{g:B\to A} f\circ g\htpy \idfunc[B].
    \end{equation*}
    In other words, a \define{section} of $f$ is a map $g:B\to A$ equipped with a homotopy $f\circ g\htpy \idfunc$. 
  \item The type of \define{retractions} of $f$\index{retraction|textbf}\index{function!has a retraction|textbf} is defined to be the type\index{retr(f)@{$\retractions(f)$}|textbf}
    \begin{equation*}
      \retractions(f) \defeq \sm{h:B\to A} h\circ f\htpy \idfunc[A].
    \end{equation*}
    If a map $f:A \to B$ has a retraction, we also say that $A$ is a \define{retract}\index{retract!of a type|textbf} of $B$.
  \item We say that a function $f:A\to B$ is an \define{equivalence}\index{equivalence|textbf}\index{is an equivalence|textbf}\index{function!is an equivalence|textbf} if it has both a section and a retraction, i.e., if it comes equipped with an element of type\index{is-equiv(f)@{$\isequiv(f)$}|textbf}
    \begin{equation*}
      \isequiv(f)\defeq\sections(f)\times\retractions(f).
    \end{equation*}
    We will write $\eqv{A}{B}$\index{A simeq B@{$\eqv{A}{B}$}|see {equivalence}} for the type $\sm{f:A\to B}\isequiv(f)$ of all equivalences from $A$ to $B$.
    For any equivalence $e:A\simeq B$ we define $e^{-1}$ to be the section of $e$.\index{equivalence!inverse|textbf}\index{inverse!of an equivalence|textbf}
  \end{enumerate}
\end{defn}

\begin{rmk}
An equivalence, as we defined it here, can be thought of as a \emph{bi-invertible map}\index{bi-invertible map|see {equivalence}}, since it comes equipped with a separate left and right inverse. Explicitly, if $f$ is an equivalence, then there are
\begin{align*}
g & : B\to A & h & : B\to A \\
G & : f\circ g \htpy \idfunc[B] & H & : h\circ f \htpy \idfunc[A].
\end{align*}
\end{rmk}

\begin{eg}\label{thm:id_equiv}
  For any type $A$, the identity function $\idfunc:A\to A$ is an equivalence, since it is its own section and its own retraction\index{identity function!is an equivalence}\index{is an equivalence!identity function}
\end{eg}

\begin{eg}\label{ex:neg_equiv}
  Since we have seen in \cref{rmk:negnegbool} that the negation function $\negbool:\bool\to\bool$ on the booleans is its own inverse, it follows that $\negbool$ is an equivalence.\index{neg bool@{$\negbool$}!is an equivalence}\index{is an equivalence!neg bool@{$\negbool$}}
\end{eg}

\begin{eg}\label{eg:is-equiv-succ-Z}
  The successor and predecessor functions on $\Z$ are equivalences by \cref{ex:is-equiv-succ-Z}\index{successor function!on Z@{on $\Z$}!is an equivalence}\index{succ Z@{$\succZ$}!is an equivalence}\index{is an equivalence!succ Z@{$\succZ$}}. Furthermore, the function
  \begin{equation*}
    x\mapsto x+k
  \end{equation*}
  is an equivalence from $\Z$ to $\Z$, for each $k:\Z$. This follows from the group laws on $\Z$, proven in \cref{ex:int_group_laws}. Indeed, the inverse of $x\mapsto x+k$ is the map $x\mapsto x+(-k)$. Finally, it also follows from the group laws on $\Z$ that the map $x\mapsto -x$ is an equivalence.

  The same holds for the finite types: the maps $\succFin_{k}$, $\predFin_{k}$, $\addFin_{k}(x)$ and $\negFin_{k}$ are all equivalences on $\Fin{k}$.
\end{eg}

\begin{rmk}\label{rmk:has-inverse}
  More generally, if $f$ \define{has an inverse}\index{has an inverse|textbf}\index{function!has an inverse|textbf} in the sense that we have a function $g:B\to A$ equipped with homotopies $f\circ g\htpy\idfunc[B]$ and $g\circ f\htpy\idfunc[A]$, then $f$ is an equivalence. We write\index{has-inverse(f)@{$\hasinverse(f)$}}
  \begin{equation*}
    \hasinverse(f)\defeq\sm{g:B\to A} (f\circ g\htpy \idfunc[B])\times (g\circ f\htpy\idfunc[A]).
  \end{equation*}
  However, we did \emph{not} define equivalences to be functions that have inverses. The reason is that we would like that being an equivalence is a \emph{property}, not a non-trivial structure on the map $f$. This fact requires the function extensionality axiom, but we can already say that if a map $f$ is an equivalence, then it has up to homotopy only one section and only one retraction (see \cref{ex:isprop_isequiv}).

  The type $\hasinverse(f)$ on the other hand, turns out to be homotopically complicated. In \cref{ex:is_invertible_id_S1} we will see that the identity function $\idfunc[\sphere{1}]:\sphere{1}\to\sphere{1}$ on the circle is an example of a map for which
  \begin{equation*}
    \hasinverse(\idfunc[\sphere{1}])\simeq \Z.
  \end{equation*}
\end{rmk}

Even though $\isequiv(f)$ and $\hasinverse(f)$ can be wildly different types, there are maps back and forth between the two. We have already observed in \cref{rmk:has-inverse} that there is a map
\begin{equation*}
  \hasinverse(f)\to\isequiv(f).
\end{equation*}
The following proposition gives the converse implication.

\begin{prp}\label{lem:inv_equiv}
  Any map $f:A\to B$ which is an equivalence, can be given the structure of an invertible map\index{equivalence!has an inverse} i.e., there is a map
  \begin{equation*}
    \isequiv(f)\to\hasinverse(f).
  \end{equation*}
\end{prp}

\begin{proof}
First we construct for any equivalence $f$ with right inverse $g$ and left inverse $h$ a homotopy $K:g\htpy h$. For any $y:B$, we have 
\begin{equation*}
\begin{tikzcd}[column sep=huge]
g(y) \arrow[r,equals,"H(g(y))^{-1}"] & hfg(y) \arrow[r,equals,"\ap{h}{G(y)}"] & h(y).
\end{tikzcd}
\end{equation*} 
In other words, the homotopy $K:g\htpy h$ is defined to be $\ct{(H\cdot g)^{-1}}{(h\cdot G)}$.
Using the homotopy $K$ we are able to show that $g$ is also a left inverse of $f$. For $x:A$ we have the identification
\begin{equation*}
\begin{tikzcd}[column sep=large]
gf(x) \arrow[r,equals,"K(f(x))"] & hf(x) \arrow[r,equals,"H(x)"] & x.
\end{tikzcd}\qedhere
\end{equation*}
\end{proof}

\begin{cor}
The inverse of an equivalence is again an equivalence.\index{inverse!of an equivalence!is an equivalence}\index{is an equivalence!inverse of an equivalence}
\end{cor}

\begin{proof}
Let $f:A\to B$ be an equivalence. By \cref{lem:inv_equiv} it follows that the section of $f$ is also a retraction. Therefore it follows that the section is itself an invertible map, with inverse $f$. Hence it is an equivalence.
\end{proof}

\begin{eg}\label{eg:laws-products-coproducts}
  Types, just as sets in classical mathematics, satisfy the usual laws of coproducts and products, such as unit laws, commutativity, and associativity. These laws are formulated as equivalences:\index{unit laws!for coproducts}\index{associativity!of coproducts}\index{zero laws!for cartesian products}\index{unit laws!for cartesian products}\index{commutativity!of coproducts}\index{commutativity!of cartesian products}\index{associativity!of cartesian products}\index{distributivity!of cartesian product over coproduct}\index{coproduct!unit laws}\index{coproduct!associativity}\index{coproduct!commutativity}\index{cartesian product type!zero laws}\index{cartesian product type!unit laws}\index{cartesian product type!commutativity}\index{cartesian product type!associativity}\index{cartesian product type!distributivity over coproducts}
  \begin{align*}
    \emptyt+B & \simeq B & A+\emptyt & \simeq A \\
    A+B & \simeq B+A & (A+B)+C & \simeq A+(B+C) \\
    \emptyt\times B & \simeq \emptyt & A\times\emptyt & \simeq \emptyt\\
    \unit\times B & \simeq B & A\times\unit & \simeq A \\
    A\times B & \simeq B\times A & (A \times B) \times C & \simeq A \times (B \times C) \\
    A\times (B+C) & \simeq (A\times B)+(A\times C) & (A+B)\times C & \simeq (A\times C)+(B\times C).
  \end{align*}
  All of these equivalences are constructed in a similar way: the maps back and forth as well as the required homotopies are constructed using induction, or, more efficiently, using pattern matching. For example, to show that cartesian products distribute from the left over coproducts, we construct maps
  \begin{align*}
    \alpha & : A\times(B+C)\to (A\times B)+(A\times C) \\
    \beta & : (A\times B)+(A\times C)\to A\times(B+C)
  \end{align*}
  as follows:
  \begin{align*}
    \alpha(x,\inl(y)) & \defeq \inl(x,y) & \beta(\inl(x,y)) & \defeq (x,\inl(y)) \\
    \alpha(x,\inr(z)) & \defeq \inr(x,z) & \beta(\inr(x,z)) & \defeq (x,\inr(z)).
  \end{align*}
  The homotopies $G:\alpha\circ\beta\htpy\idfunc$ and $H:\beta\circ\alpha\htpy \idfunc$ are then defined by
  \begin{align*}
    G(\inl(x,y)) & \defeq \refl{} & H(x,\inl(y)) & \defeq \refl{} \\
    G(\inr(x,z)) & \defeq \refl{} & H(x,\inr(z)) & \defeq \refl{}.
  \end{align*}
  We encourage the reader to write out the definitions of at least a few of these equivalences.
\end{eg}

\begin{eg}\label{eg:laws-Sigma-types}
  The laws for cartesian products can be generalized to arbitrary $\Sigma$-types. The absorption laws and unit laws, for instance, are as follows:
  \index{absorption laws!of dependent pair types}\index{dependent pair type!absorption laws}\index{unit laws!for dependent pair types}\index{dependent pair type!unit laws}
  \begin{align*}
    \sm{x:\emptyt}B(x) & \simeq \emptyt & \sm{x:A}\emptyt & \simeq \emptyt \\
    \sm{x:\unit}B(x) & \simeq B(\ttt) & \sm{x:A}\unit & \simeq A.
  \end{align*}
  Note that the right absorption law and the right unit law are exactly the same as the right absorption and unit laws for cartesian products. The left absorption and unit laws are, however, formulated with a type family $B$ over $\emptyt$ and over $\unit$, and therefore they are slightly more general.
  
  Commutativity cannot be generalized to $\Sigma$-types. Associativity, on the other hand, can be expressed in two ways:\index{associativity!of dependent pair types}\index{dependent pair type!associativity}
  \begin{align*}
    \sm{w:\sm{x:A}B(x)}C(w) & \simeq\sm{x:A}\sm{y:B}C(\pair(x,y)) \\
    \sm{w:\sm{x:A}B(x)}C(\proj 1(w),\proj 2(w)) & \simeq \sm{x:A}\sm{y:B(x)}C(x,y). 
  \end{align*}
  In the first of these equivalences associativity is stated using a type family $C$ over $\sm{x:A}B(x)$ while in the second it is stated using a family of types $C(x,y)$ indexed by $x:A$ and $y:B(x)$.
  
  Finally, we note that $\Sigma$ also distributes over coproducts. In other words, there are the following two equivalences:\index{dependent pair type!distributivity over coproducts}\index{distributivity!of S-types over coproducts@{of $\Sigma$-types over coproducts}}
  \begin{align*}
    \sm{x:A}B(x)+C(x) & \simeq \Big(\sm{x:A}B(x)\Big)+\Big(\sm{x:A}C(x)\Big) \\
    \sm{w:A+B}C(w) & \simeq \Big(\sm{x:A}C(\inl(x))\Big)+\Big(\sm{y:B}C(\inr(y))\Big).
  \end{align*}
\end{eg}

\begin{rmk}
    We haven't stated any laws involving function types or dependent function types, because it requires the function extensionality principle to prove them.
\end{rmk}

\subsection{Characterizing the identity types of \texorpdfstring{$\Sigma$-}{dependent pair }types}

\index{identity type!of a Sigma-type@{of a $\Sigma$-type}|(}
\index{dependent pair type!identity type|(}
\index{characterization of identity type!of S-types@{of $\Sigma$-types}|(}
In this section we characterize the identity type of a $\Sigma$-type as a $\Sigma$-type of identity types. Characterizing identity types is a task that a homotopy type theorist routinely performs, so we will follow the general outline of how such a characterization goes:
\begin{enumerate}
\item First we define a binary relation $R:A\to A\to \UU$ on the type $A$ that we are interested in. This binary relation is intended to be equivalent to its identity type.
\item Then we will show that this binary relation is reflexive, by constructing a dependent function of type
  \begin{equation*}
    \prd{x:A}R(x,x)
  \end{equation*}
\item Using the reflexivity we will show that there is a canonical map
  \begin{equation*}
    (x=y)\to R(x,y)
  \end{equation*}
  for every $x,y:A$. This map is just constructed by path induction, using the reflexivity of $R$.
\item Finally, it has to be shown that the map
  \begin{equation*}
    (x=y)\to R(x,y)
  \end{equation*}
  is an equivalence for each $x,y:A$. 
\end{enumerate}
The last step is usually the most difficult, and we will refine our methods for this step in \cref{chap:fundamental}, where we establish the fundamental theorem of identity types.

In this section we consider a type family $B$ over $A$. Given two pairs
\begin{equation*}
  (x,y),(x',y'):\sm{x:A}B(x),
\end{equation*}
if we have a path $\alpha:x=x'$ then we can compare $y:B(x)$ to $y':B(x')$ by first transporting $y$ along $\alpha$, i.e., we consider the identity type
\begin{equation*}
  \tr_B(\alpha,y)=y'.
\end{equation*}
Thus it makes sense to think of $(x,y)$ to be identical to $(x',y')$ if there is an identification $\alpha:x=x'$ and an identification $\beta:\tr_B(\alpha,y)=y'$. In the following definition we turn this idea into a binary relation on the $\Sigma$-type.

\begin{defn}
  We will define a relation\index{Eq Sigma@{$\EqSigma$}|textbf}\index{dependent pair type!Eq Sigma@{$\EqSigma$}|textbf}\index{dependent pair type!observational equality|textbf}\index{observational equality!on S-types@{on $\Sigma$-types}|textbf}
  \begin{equation*}
    \EqSigma : \Big(\sm{x:A}B(x)\Big)\to\Big(\sm{x:A}B(x)\Big)\to\UU
  \end{equation*}
  by defining
  \begin{equation*}
    \EqSigma(s,t)\defeq\sm{\alpha:\proj 1(s)=\proj 1(t)}\tr_B(\alpha,\proj 2(s))=\proj 2 (t).
  \end{equation*}
\end{defn}

\begin{lem}
  The relation $\EqSigma$ is reflexive, i.e., we can construct
  \begin{equation*}
    \reflexiveEqSigma:\prd{s:\sm{x:A}B(x)}\EqSigma(s,s).
  \end{equation*}
\end{lem}

\begin{constr}
  The element $\reflexiveEqSigma$ is constructed by $\Sigma$-induction on $s:\sm{x:A}B(x)$. Thus, it suffices to construct a dependent function of type
  \begin{equation*}
    \prd{x:A}\prd{y:B(x)}\sm{\alpha:x=x}\tr_B(\alpha,y)=y.
  \end{equation*}
  Here we take $\lam{x}\lam{y}(\refl{x},\refl{y})$.
\end{constr}

\begin{defn}
  Consider a type family $B$ over $A$. Then for any $s,t:\sm{x:A}B(x)$ we define a map\index{pair-eq@{$\paireq$}|textbf}
  \begin{equation*}
    \paireq: (s=t)\to \EqSigma(s,t)
  \end{equation*}
  by path induction, taking $\paireq(\refl{s})\defeq\reflexiveEqSigma(s)$.
\end{defn}

\begin{thm}\label{thm:eq_sigma}
  Let $B$ be a type family over $A$. Then the map\index{pair-eq@{$\paireq$}!is an equivalence}\index{is an equivalence!pair-eq@{$\paireq$}}
  \begin{equation*}
    \paireq: (s=t)\to \EqSigma(s,t)
  \end{equation*}
  is an equivalence for every $s,t:\sm{x:A}B(x)$.\index{is an equivalence!pair-eq@{$\paireq$}}
\end{thm}

\begin{proof}
The maps in the converse direction\index{eq-pair@{$\eqpair$}|textbf}
\begin{equation*}
\eqpair : \EqSigma(s,t)\to(\id{s}{t})
\end{equation*}
are defined by repeated $\Sigma$-induction. By $\Sigma$-induction on $s$ and $t$  we see that it suffices to define a map
\begin{equation*}
\eqpair : \Big(\sm{p:x=x'}\id{\tr_B(p,y)}{y'}\Big)\to(\id{(x,y)}{(x',y')}).
\end{equation*}
A map of this type is again defined by $\Sigma$-induction. Thus it suffices to define a dependent function of type
\begin{equation*}
\prd{p:x=x'} (\id{\tr_B(p,y)}{y'}) \to (\id{(x,y)}{(x',y')}).
\end{equation*}
Such a dependent function is defined by double path induction by sending $\pairr{\refl{x},\refl{y}}$ to $\refl{\pairr{x,y}}$. This completes the definition of the function $\eqpair$.

Next, we must show that $\eqpair$ is a section of $\paireq$. In other words, we must construct an identification
\begin{equation*}
\paireq(\eqpair(\alpha,\beta))=\pairr{\alpha,\beta}
\end{equation*}
for each $\pairr{\alpha,\beta}:\sm{\alpha:x=x'}\id{\tr_B(\alpha,y)}{y'}$. We proceed by path induction on $\alpha$, followed by path induction on $\beta$. Then our goal becomes to construct an identification of type
\begin{equation*}
\paireq(\eqpair\pairr{\refl{x},\refl{y}})=\pairr{\refl{x},\refl{y}}
\end{equation*}
By the definition of $\eqpair$ we have $\eqpair\pairr{\refl{x},\refl{y}}\jdeq \refl{\pairr{x,y}}$, and by the definition of $\paireq$ we have $\paireq(\refl{\pairr{x,y}})\jdeq\pairr{\refl{x},\refl{y}}$. Thus we may take $\refl{\pairr{\refl{x},\refl{y}}}$ to complete the construction of the homotopy $\paireq\circ\eqpair\htpy\idfunc$.

To complete the proof, we must show that $\eqpair$ is a retraction of $\paireq$. In other words, we must construct an identification
\begin{equation*}
\eqpair(\paireq(p))=p
\end{equation*}
for each $p:s=t$. We proceed by path induction on $p:s=t$, so it suffices to construct an identification 
\begin{equation*}
\eqpair\pairr{\refl{\proj 1(s)},\refl{\proj 2(s)}}=\refl{s}.
\end{equation*}
Now we proceed by $\Sigma$-induction on $s:\sm{x:A}B(x)$, so it suffices to construct an identification
\begin{equation*}
\eqpair\pairr{\refl{x},\refl{y}}=\refl{(x,y)}.
\end{equation*}
Since $\eqpair\pairr{\refl{x},\refl{y}}$ computes to $\refl{(x,y)}$, we may simply take $\refl{\refl{(x,y)}}$.
\end{proof}
\index{identity type!of a Sigma-type@{of a $\Sigma$-type}|)}
\index{dependent pair type!identity type|)}
\index{characterization of identity type!of S-types@{of $\Sigma$-types}|)}

\begin{exercises}
  \exitem \label{ex:equiv_grpd_ops}Show that the functions\index{inv@{$\invfunc$}!is an equivalence}\index{is an equivalence!inv@{$\invfunc$}}\index{concat@{$\concat$}!is a family of equivalences}\index{is an equivalence!concat(p)@{$\concat(p)$}}\index{concat'@{$\concat'$}!is a family of equivalences}\index{is an equivalence!concat'@{$\concat'(q)$}}\index{tr B@{$\tr_B$}!is a family of equivalences}\index{is an equivalence!tr B(p)@{$\tr_B(p)$}}
  \begin{align*}
    \invfunc & :(\id{x}{y})\to(\id{y}{x}) \\
    \concat(p) & : (\id{y}{z})\to(\id{x}{z}) \\
    \concat'(q) & : (\id{x}{y}) \to (\id{x}{z}) \\
    \tr_B(p) & :B(x)\to B(y)
  \end{align*}
  are equivalences, where $\concat'(q,p)\defeq \ct{p}{q}$\index{concat'@{$\concat'$}|textbf}. Give their inverses explicitly.
  \exitem
  \begin{subexenum}
  \item Use \cref{ex:zero-neq-one-bool} to show that the constant function $\const_b:\bool\to\bool$ is not an equivalence, for any $b:\bool$.\index{booleans!const b is not an equivalence@{$\const_b$ is not an equivalence}}
  \item Show that $\bool\not\simeq \unit$.
  \item Show that $\N\not\simeq \Fin{k}$ for any $k:\N$. 
  \end{subexenum}
  \exitem
  \begin{subexenum}
  \item \label{ex:htpy_equiv}\index{equivalence!closed under homotopies} Consider two functions $f,g:A\to B$ and a homotopy $H:f\htpy g$. Then
    \begin{equation*}
      \isequiv(f)\leftrightarrow\isequiv(g).
    \end{equation*}
  \item Show that for any two homotopic equivalences $e,e':\eqv{A}{B}$, their inverses are also homotopic.
  \end{subexenum}
  \exitem \label{ex:3_for_2}
  Consider a commuting triangle
  \begin{equation*}
    \begin{tikzcd}[column sep=tiny]
      A \arrow[rr,"h"] \arrow[dr,swap,"f"] & & B \arrow[dl,"g"] \\
      & X.
    \end{tikzcd}
  \end{equation*}
  with $H:f\htpy g\circ h$.
  \begin{subexenum}
  \item Suppose that the map $h$ has a section $s:B \to A$. Show that the triangle
    \begin{equation*}
      \begin{tikzcd}[column sep=tiny]
        B \arrow[rr,"s"] \arrow[dr,swap,"g"] & & A \arrow[dl,"f"] \\
        & X.
      \end{tikzcd}
    \end{equation*}
    commutes, and that $f$ has a section if and only if $g$ has a section.
  \item Suppose that the map $g$ has a retraction $r:X\to B$. Show that the triangle
    \begin{equation*}
      \begin{tikzcd}[column sep=tiny]
        A \arrow[rr,"f"] \arrow[dr,swap,"h"] & & X \arrow[dl,"r"] \\
        & B.
      \end{tikzcd}
    \end{equation*}
    commutes, and that $f$ has a retraction if and only if $h$ has a retraction.
  \item (The \define{3-for-2 property} for equivalences.)\index{equivalence!3-for-2 property}\index{3-for-2 property!of equivalences}\index{composition!of equivalences|textbf}\index{equivalence!composition|textbf} Show that if any two of the functions
    \begin{equation*}
      f,\qquad g,\qquad h
    \end{equation*}
    are equivalences, then so is the third. Conclude that any section and any retraction of an equivalence is again an equivalence.
  \end{subexenum}
  \exitem \label{ex:sigma_swap}
  \begin{subexenum}
  \item Let $A$ and $B$ be types, and let $C$ be a family over $x:A,y:B$. Construct an equivalence
    \begin{equation*}
      \eqv{\Big(\sm{x:A}\sm{y:B}C(x,y)\Big)}{\Big(\sm{y:B}\sm{x:A}C(x,y)\Big)}.
    \end{equation*}
  \item Let $A$ be a type, and let $B$ and $C$ be type families over $A$. Construct an equivalence
    \begin{equation*}
      \Big(\sm{u:\sm{x:A}B(x)}C(\proj 1(u))\Big) \simeq \Big(\sm{v:\sm{x:A}C(x)}B(\proj 1(v))\Big).
    \end{equation*}
  \end{subexenum}
  \exitem \label{ex:coproduct_functor}Recall from \cref{rmk:functor-coprod} that coproducts have a \define{functorial action}\index{functorial action!of coproducts}\index{coproduct!functorial action}, i.e., that for every $f:A\to A'$ and every $g:B\to B'$ we have a map
  \begin{equation*}
    f+g:(A+B)\to (A'+B').
  \end{equation*}
  \begin{subexenum}
  \item Show that $\idfunc[A]+\idfunc[B]\htpy \idfunc[A+B]$.
  \item Show that for any two pairs of composable functions
    \begin{equation*}
      \begin{tikzcd}
        A \arrow[r,"f"] & {A'} \arrow[r,"{f'}"] & {A''}
      \end{tikzcd}
      \qquad\text{and}\qquad
      \begin{tikzcd}
        B \arrow[r,"g"] & {B'} \arrow[r,"{g'}"] & {B''}
      \end{tikzcd}
    \end{equation*}
    there is a homotopy $(f'\circ f)+(g'\circ g) \htpy (f'+g')\circ (f+g)$.
  \item Show that if $H:f\htpy f'$ and $K:g\htpy g'$, then there is a homotopy
    \begin{equation*}
      H+K:(f+g)\htpy (f'+g').
    \end{equation*}
  \item \label{ex:coproduct_functor_equivalence}Show that if both $f$ and $g$ are equivalences, then so is $f+g$. (The converse of this statement also holds, see \cref{ex:is-equiv-is-equiv-functor-coprod}.)
  \end{subexenum}
  \exitem
  \begin{subexenum}
  \item Construct for any two maps $f:A \to A'$ and $g:B\to B'$, a map
    \begin{equation*}
      f\times g : A\times B \to A'\times B'
    \end{equation*}
  \item Show that $\idfunc[A]\times\idfunc[B]\htpy\idfunc[A\times B]$.
  \item Show that for any two pairs of composable functions
    \begin{equation*}
      \begin{tikzcd}
        A \arrow[r,"f"] & {A'} \arrow[r,"{f'}"] & {A''}
      \end{tikzcd}
      \qquad\text{and}\qquad
      \begin{tikzcd}
        B \arrow[r,"g"] & {B'} \arrow[r,"{g'}"] & {B''}
      \end{tikzcd}
    \end{equation*}
    there is a homotopy $(f'\circ f)\times(g'\circ g) \htpy (f'\times g')\circ (f\times g)$.
  \item Show that if $H:f\htpy f'$ and $K:g\htpy g'$, then there is a homotopy
    \begin{equation*}
      H\times K:(f\times g)\htpy (f'\times g').
    \end{equation*}
  \item Show that for any two maps $f:A\to A'$ and $g:B\to B'$, the following are equivalent:
    \begin{enumerate}
    \item The map $f\times g$ is an equivalence.
    \item There are functions
      \begin{align*}
        \alpha & : B \to \isequiv(f) \\
        \beta & : A \to \isequiv(g).
      \end{align*}
    \end{enumerate}
  \end{subexenum}
  \exitem\label{ex:laws-Fin} Construct equivalences
  \begin{align*}
    \Fin{k+l} & \simeq \Fin{k}+\Fin{l} \\
    \Fin{kl} & \simeq \Fin{k}\times\Fin{l}.
  \end{align*}
  \exitem A map $f:X\to X$ is said to be \define{finitely cyclic}\index{finitely cyclic type|textbf} if it comes equipped with an element of type
  \begin{equation*}
    \isfinitelycyclic(f)\defeq\prd{x,y:X}\sm{k:\N}f^k(x)=y.
  \end{equation*}
  \begin{subexenum}
  \item Show that any finitely cyclic map is an equivalence.
  \item Show that $\succFin:\Fin{k}\to\Fin{k}$ is finitely cyclic for any $k:\N$.
  \end{subexenum}
\end{exercises}
\index{equivalence|)}



\section{Contractible types and contractible maps}
\sectionmark{Contractible types and maps}
\label{sec:contractible}

\index{contractible type|(} 
A contractible type is a type which has, up to identification, only one element. In other words, a contractible type is a type that comes equipped with a point, and an identification of this point with any point.

We may think of contractible types as singletons up to homotopy, and indeed we show that the unit type is an example of a contractible type. Moreover, we show that contractible types satisfy an induction principle that is very similar to the induction principle of the unit type.

Another case of an inductive type with a single constructor is the type of identifications $p:a=x$ with a fixed starting point $a:A$. To specify such an identification, we have to give its end point $x:A$ as well as the identification $p:a=x$, and the path induction principle asserts that in order to show something about all such identifications, it suffices to show that thing in the case where the end point is $a$, and the path is $\refl{a}$. This suggests that the total space
\begin{equation*}
  \sm{x:A}a=x
\end{equation*}
of all paths with starting point $a:A$ is contractible. This important fact will be shown in \cref{thm:total_path}, and it is the basis for the fundamental theorem of identity types (\cref{chap:fundamental}).

Next, we introduce the \emph{fiber} of a map $f:A\to B$. The fiber of $f$ at $b:B$ consists of the type of elements $a:A$ equipped with an identification $p:f(a)=b$. In other words, the fiber of $f$ at $b$ is the preimage of $f$ at $b$. In \cref{thm:equiv_contr,thm:contr_equiv} we show that a map is an equivalence if and only if its fibers are contractible. The condition that the fibers of a map are contractible is analogous to the set theoretic notion of bijective map, or $1$-to-$1$-correspondence.

\subsection{Contractible types}

\begin{defn}
  We say that a type $A$ is \define{contractible}\index{contractible type|textbf} if it comes equipped with an element of type\index{is-cont(A)r@{$\iscontr(A)$}|see {contractible type}}
  \begin{equation*}
    \iscontr(A) \defeq \sm{c:A}\prd{x:A}c=x.
  \end{equation*}
  Given a pair $(c,C):\iscontr(A)$, we call $c:A$ the \define{center of contraction}\index{center of contraction|textbf}\index{contractible type!center of contraction|textbf} of $A$, and we call $C:\prd{x:A}c=x$ the \define{contraction}\index{contraction|textbf}\index{contractible type!contraction|textbf} of $A$.
\end{defn}

\begin{rmk}
Suppose $A$ is a contractible type with center of contraction $c$ and contraction $C$. Then the type of $C$ is (judgmentally) equal to the type
\begin{equation*}
\const_c\htpy\idfunc[A].
\end{equation*}
In other words, the contraction $C$ is a \emph{homotopy} from the constant function to the identity function.
\end{rmk}

\begin{eg}
  The unit type is easily seen to be contractible.\index{unit type!is contractible}\index{is contractible!unit type} For the center of contraction we take $\ttt:\unit$. Then we define a contraction $\prd{x:\unit}\ttt=x$ by the induction principle of $\unit$. Applying the induction principle, it suffices to construct an identification of type $\ttt = \ttt$, for which we just take $\refl{\ttt}$.
\end{eg}

\begin{thm}\label{thm:total_path}
For any $a:A$, the type
\begin{equation*}
\sm{x:A}a=x
\end{equation*}
is contractible.\index{identity type!total space is contractible}\index{is contractible!total space of identity type}
\end{thm}

\begin{proof}
  For the center of contraction we take
  \begin{equation*}
    (a,\refl{a}):\sm{x:A}a=x.
  \end{equation*}
  The contraction is constructed in \cref{prp:contraction-total-space-id}.
\end{proof}

\subsection{Singleton induction}
\index{singleton induction|(}
\index{induction principle!singleton induction|(}

Contractible types are singletons up to homotopy. Indeed, every element of a contractible type can be identified with the center of contraction. Therefore we can prove an induction principle for contractible types that is similar to the induction principle of the unit type.

\begin{defn}\label{defn:singleton-induction}
  Suppose $A$ comes equipped with an element $a:A$. Then we say that $A$ satisfies \define{singleton induction}\index{singleton induction|textbf}\index{induction principle!singleton induction|textbf} if for every type family $B$ over $A$, the map\index{ev-pt@{$\evpt$}|textbf}
  \begin{equation*}
    \evpt:\Big(\prd{x:A}B(x)\Big)\to B(a)
  \end{equation*}
  defined by $\evpt(f)\defeq f(a)$ has a section. In other words, if $A$ satisfies singleton induction we have a function and a homotopy\index{ind-sing@{$\singind$}|textbf}\index{comp-sing@{$\singcomp$}|textbf}
  \begin{align*}
    \singind_{a} & : B(a)\to \prd{x:A}B(x) \\
    \singcomp_{a} & : \evpt\circ \singind_{a} \htpy \idfunc
  \end{align*}
  for any type family $B$ over $A$.
\end{defn}

\begin{eg}
  Note that the singleton induction principle is almost the same as the induction principle for the unit type, the difference being that the `computation rule' in the singleton induction for $A$ is stated using an \emph{identification} rather than as a judgmental equality. The unit type\index{unit type!singleton induction} $\unit$ comes equipped with a function
  \begin{equation*}
    \indunit:B(\ttt)\to \prd{x:\unit}B(x)
  \end{equation*}
  for every type family $B$ over $\unit$, satisfying the judgmental equality $\indunit(b,\ttt)\jdeq b$ for every $b:B(\ttt)$ by the computation rule. Therefore, we obtain the homotopy
  \begin{equation*}
    \lam{b}\refl{b}:\evpt\circ\indunit \htpy\idfunc,
  \end{equation*}
  and we conclude that the unit type satisfies singleton induction. 
\end{eg}

\begin{thm}\label{thm:contractible}
Let $A$ be a type. The following are equivalent:\index{is contractible!satisfies singleton induction}\index{singleton induction!is contractible}\index{contractible type!satisfies singleton induction}
\begin{enumerate}
\item The type $A$ is contractible.
\item The type $A$ comes equipped with an element $a:A$, and satisfies singleton induction.
\end{enumerate}
\end{thm}

\begin{proof}
Suppose $A$ is contractible with center of contraction $a$ and contraction $C$. 
First we observe that, without loss of generality, we may assume that $C$ comes equipped with an identification $p:C(a)=\refl{a}$.
To see this, note that we can always define a new contraction $C'$ by
\begin{equation*}
C'(x)\defeq\ct{C(a)^{-1}}{C(x)},
\end{equation*}
which satisfies the requirement by the left inverse law, constructed in \cref{defn:id_invlaw}.

To show that $A$ satisfies singleton induction let $B$ be a type family over $A$, and suppose we have $b:B(a)$. Our goal is to define
\begin{equation*}
  \indsing_a(b):\prd{x:A}B(x).
\end{equation*}
Let $x:A$. Since we have an identification $C(x):a=x$, and an element $b$ in $B(a)$, we may transport $b$ along the path $C(x)$ to obtain
\begin{equation*}
  \indsing_a(b,x)\defeq \tr_B(C(x),b):B(x).
\end{equation*}
Therefore, the function $\indsing_a(b)$ is defined to be the dependent function $\lam{x}\tr_B(C(x),b)$. Now we have to show that $\indsing_a(b,a)=b$. Then we have the identifications
\begin{equation*}
\begin{tikzcd}
\tr_B(C(a),b) \arrow[r,equals,"\ap{\lam{\omega}\tr_B(\omega,b)}{p}"] &[4em] \tr_B(\refl{a},b) \arrow[r,equals,"\refl{b}"] & b.
\end{tikzcd}
\end{equation*}
This shows that the computation rule is satisfied, which completes the proof that $A$ satisfies singleton induction.

For the converse, suppose that $a:A$ and that $A$ satisfies singleton induction. Our goal is to show that $A$ is contractible. For the center of contraction we take the element $a:A$. By singleton induction applied to $B(x)\defeq a=x$ we have the map 
\begin{equation*}
\indsing_{a} : a=a \to \prd{x:A}a=x.
\end{equation*}
Therefore $\indsing_{a}(\refl{a})$ is a contraction.
\end{proof}
\index{singleton induction|)}
\index{induction principle!singleton induction|)}

\subsection{Contractible maps}

\index{contractible map|(}
\begin{defn}
  Let $f:A\to B$ be a function, and let $b:B$. The \define{fiber}\index{fiber|textbf}\index{homotopy fiber|see {fiber}} of $f$ at $b$ is defined to be the type\index{fib f b@{$\fib{f}{b}$}|textbf}
  \begin{equation*}
    \fib{f}{b}\defeq\sm{a:A}f(a)=b.
  \end{equation*}
\end{defn}

In other words, the fiber of $f$ at $b$ is the type of $a:A$ that get mapped by $f$ to $b$.
One may think of the fiber as a type theoretic version of the preimage\index{pre-image|see {fiber}} of a point.

\index{fiber!characterization of identity type|(}
\index{characterization of identity type!of the fiber of a map|(}
\index{identity type!of a fiber|(}
It will be useful to have a characterization of the identity type of a fiber. In order to identify any $(x,p)$ and $(x',p')$ in $\fib{f}{y}$, we may first construct an identification $\alpha:x=x'$. Then we obtain a triangle
\begin{equation*}
  \begin{tikzcd}[column sep=tiny]
    f(x) \arrow[dr,equals,swap,"p"] \arrow[rr,equals,"{\ap{f}{\alpha}}"] & & f(x') \arrow[dl,equals,"{p'}"] \\
    \phantom{f(x')} & y,
  \end{tikzcd}
\end{equation*}
so we may consider the type of identifications $\beta:p=\ct{\ap{f}{\alpha}}{p'}$. We will show that the type of all identifications $(x,p)=(x',p')$ is equivalent to the type of such pairs $(\alpha,\beta)$. 

\begin{defn}
  Let $f:A \to B$ be a map, and let $(x,p),(x',p'):\fib{f}{y}$ for some $y:B$.
  Then we define\index{Eq-fib@{$\Eqfib$}|textbf}\index{fiber!Eq-fib@{$\Eqfib$}|textbf}\index{observational equality!fiber}
  \begin{equation*}
    \Eqfib_f((x,p),(x',p'))\defeq \sm{\alpha:x=x'}p=\ct{\ap{f}{\alpha}}{p'}
  \end{equation*}
  The relation $\Eqfib_f:\fib{f}{y}\to\fib{f}{y}\to\UU$ is a reflexive relation, since we have
  \begin{equation*}
    \lam{(x,p)}(\refl{x},\refl{p}):\prd{(x,p):\fib{f}{y}}\Eqfib_f((x,p),(x,p)).
  \end{equation*}
\end{defn}

\begin{prp}
  Consider a map $f:A\to B$ and let $y:B$. The canonical map
  \begin{equation*}
    ((x,p)=(x',p'))\to\Eqfib_f((x,p),(x',p'))
  \end{equation*}
  induced by the reflexivity of $\Eqfib_f$ is an equivalence for any $(x,p),(x',p'):\fib{f}{y}$.
\end{prp}

\begin{proof}
  The converse map
  \begin{equation*}
    \Eqfib_f((x,p),(x',p'))\to ((x,p)=(x',p'))
  \end{equation*}
  is easily defined by $\Sigma$-induction, and then path induction twice. The homotopies witnessing that this converse map is indeed a right inverse as well as a left inverse are similarly constructed by induction.
\end{proof}
\index{fiber!characterization of identity type|)}
\index{characterization of identity type!of the fiber of a map|)}
\index{identity type!of a fiber|)}

Now we define at the notion of contractible map.

\begin{defn}
We say that a function $f:A\to B$ is \define{contractible}\index{contractible map|textbf} if it comes equipped with an element of type\index{is-contr(f)@{$\iscontr(f)$}|see {contractible map}}\index{is a contractible map|textbf}
\begin{equation*}
\iscontr(f)\defeq\prd{b:B}\iscontr(\fib{f}{b}).
\end{equation*}
\end{defn}

\begin{thm}\label{thm:equiv_contr}
Any contractible map is an equivalence.\index{contractible map!is an equivalence}\index{is an equivalence!contractible map}
\end{thm}

\begin{proof}
Let $f:A\to B$ be a contractible map. Using the center of contraction of each $\fib{f}{y}$, we obtain the dependent function
\begin{align*}
\lam{y}\pairr{g(y),G(y)}:\prd{y:B}\fib{f}{y}.
\end{align*}
Thus, we get map $g:B\to A$, and a homotopy $G:\prd{y:B} f(g(y))=y$. In other words, we get a section of $f$.

It remains to construct a retraction of $f$. Taking $g$ as our retraction, we have to show that $\prd{x:A} g(f(x))=x$. Note that we get an identification $p:f(g(f(x)))=f(x)$ since $g$ is a section of $f$. Therefore, it follows that $(g(f(x)),p):\fib{f}{f(x)}$. Moreover, since $\fib{f}{f(x)}$ is contractible we get an identification $q:\pairr{g(f(x)),p}=\pairr{x,\refl{f(x)}}$. The base path $\ap{\proj 1}{q}$ of this identification is an identification of type $g(f(x))=x$, as desired.
\end{proof}

\subsection{Equivalences are contractible maps}\label{sec:is-contr-map-is-equiv}

In \cref{thm:contr_equiv} we will show the converse to \cref{thm:equiv_contr}, i.e., we will show that any equivalence is a contractible map. We will do this in two steps.

First we introduce a new notion of \emph{coherently invertible map}, for which we can easily show that such maps have contractible fibers. Then we show that any equivalence is a coherently invertible map.

  Recall that an invertible map is a map $f:A\to B$ equipped with $g:B\to A$ and homotopies
  \begin{equation*}
    G : f\circ g \htpy \idfunc\qquad\text{and}\qquad H:g\circ f\htpy \idfunc.
  \end{equation*}
  Then we observe that both $G \cdot f$ and $f \cdot H$ are homotopies of the same type
  \begin{equation*}
    f\circ g\circ f \htpy f.
  \end{equation*}
  A coherently invertible map is an invertible map for which there is a further homotopy $G \cdot f\htpy f\cdot H$.

  \begin{defn}
    Consider a map $f:A\to B$. We say that $f$ is \define{coherently invertible}\index{coherently invertible map|textbf} if it comes equipped with
    \begin{align*}
      g & : B \to A \\
      G & : f \circ g \htpy \idfunc \\
      H & : g \circ f \htpy \idfunc \\
      K & : G \cdot f \htpy f \cdot H.
    \end{align*}
    We will write $\iscohinvertible(f)$\index{is-coh-invertible(f)@{$\iscohinvertible(f)$}|textbf} for the type of quadruples $(g,G,H,K)$.
  \end{defn}

  Although we will encounter the notion of coherently invertible map on some further occasions, the following proposition is our main motivation for considering it.

  \begin{prp}\label{lem:contr-inv}
    Any coherently invertible map has contractible fibers.\index{coherently invertible map!is a contractible map}
  \end{prp}

  \begin{proof}
    Consider a map $f:A\to B$ equipped with
    \begin{align*}
      g & : B \to A \\
      G & : f \circ g \htpy \idfunc \\
      H & : g \circ f \htpy \idfunc \\
      K & : G \cdot f \htpy f \cdot H,
    \end{align*}
    and let $y:B$. Our goal is to show that $\fib{f}{y}$ is contractible. For the center of contraction we take $(g(y),G(y))$. In order to construct a contraction, it suffices to construct a dependent function of type
    \begin{equation*}
      \prd{x:A}\prd{p:f(x)=y}\Eqfib_f((g(y),G(y)),(x,p)).
    \end{equation*}
    By path induction on $p:f(x)=y$ it suffices to construct a dependent function of type
    \begin{equation*}
      \prd{x:A}\Eqfib_f((g(f(x)),G(f(x))),(x,\refl{f(x)})).
    \end{equation*}
    By definition of $\Eqfib_f$, we have to construct for each $x:A$ an identification $\alpha:g(f(x))=x$ equipped with a further identification
    \begin{equation*}
      G(f(x))=\ct{\ap{f}{\alpha}}{\refl{f(x)}}.
    \end{equation*}
    Such a dependent function is constructed as $\lam{x}(H(x),K'(x))$, where the homotopy $H:g\circ f\htpy \idfunc$ is given by assumption, and the homotopy
    \begin{align*}
      K' & : \prd{x:A}G(f(x))=\ct{\ap{f}{H(x)}}{\refl{f(x)}}
    \end{align*}
    is defined as
    \begin{equation*}
      K'\defeq \ct{K}{\rightunithtpy(f\cdot H)^{-1}}.\qedhere
    \end{equation*}
  \end{proof}

  Our next goal is to show that for any map $f:A\to B$ equipped with
  \begin{equation*}
    g:B\to A,\qquad G:f\circ g \htpy \idfunc,\qquad\text{and}\qquad H:g\circ f\htpy \idfunc,
  \end{equation*}
  we can improve the homotopy $G$ to a new homotopy $G':f\circ g\htpy \idfunc$ for which there is a further homotopy
  \begin{equation*}
    f\cdot H\htpy G'\cdot f.
  \end{equation*}
  Note that this situation is analogous to the situation in the proof of \cref{thm:contractible}, where we improved the contraction $C$ so that it satisfied $C(c)=\refl{}$. The extra coherence $f\cdot H\htpy G'\cdot f$ is then used in the proof that the fibers of an equivalence are contractible.

\begin{defn}\label{defn:htpy_nat}\index{homotopy!naturality|textbf}
Let $f,g:A\to B$ be functions, and consider $H:f\htpy g$ and $p:x=y$ in $A$. We define the identification\index{nat-nat@{$\nathtpy$}|textbf}\index{homotopy!nat-htpy@{$\nathtpy$}|textbf}
\begin{equation*}
\nathtpy(H,p) \defeq  \ct{\ap{f}{p}}{H(y)}=\ct{H(x)}{\ap{g}{p}}
\end{equation*}
witnessing that the square
\begin{equation*}
\begin{tikzcd}
f(x) \arrow[r,equals,"H(x)"] \arrow[d,equals,swap,"\ap{f}{p}"] & g(x) \arrow[d,equals,"\ap{g}{p}"] \\
f(y) \arrow[r,equals,swap,"H(y)"] & g(y)
\end{tikzcd}
\end{equation*}
commutes. This square is also called the \define{naturality square}\index{naturality square of homotopies|textbf} of the homotopy $H$ at $p$.
\end{defn}

\begin{constr}
  By path induction on $p$ it suffices to construct an identification
  \begin{equation*}
    \ct{\ap{f}{\refl{x}}}{H(x)}=\ct{H(x)}{\ap{g}{\refl{x}}}
  \end{equation*}
  since $\ap{f}{\refl{x}}\jdeq \refl{f(x)}$ and $\ap{g}{\refl{x}}\jdeq\refl{g(x)}$, and since $\ct{\refl{f(x)}}{H(x)}\jdeq H(x)$, we see that the path $\rightunit(H(x))^{-1}$ is of the asserted type.
\end{constr}

\begin{defn}\label{defn:retraction_swap}
Consider $f:A\to A$ and $H: f\htpy \idfunc[A]$. We construct an identification $H(f(x))=\ap{f}{H(x)}$, for any $x:A$.
\end{defn}

\begin{constr}
By the naturality of homotopies with respect to identifications the square
\begin{equation*}
\begin{tikzcd}[column sep=large]
ff(x) \arrow[d,swap,equals,"\ap{f}{H(x)}"] \arrow[r,equals,"H(f(x))"] & f(x) \arrow[d,equals,"H(x)"] \\
f(x) \arrow[r,swap,equals,"H(x)"] & x
\end{tikzcd}
\end{equation*}
commutes. This gives the desired identification $H(f(x))=\ap{f}{H(x)}$.
\end{constr}

\begin{lem}\label{lem:coherently-invertible}
  Let $f:A\to B$ be a map equipped with an inverse, i.e., consider
  \begin{align*}
    g & : B \to A \\
    G & : f \circ g \htpy \idfunc \\
    H & : g \circ f \htpy \idfunc.
  \end{align*}
  Then there is a homotopy $G':f\circ g\htpy \idfunc$ equipped with a further homotopy
  \begin{equation*}
    K : f\cdot H \htpy G'\cdot f.
  \end{equation*}
  Thus we obtain a map $\hasinverse(f)\to\iscohinvertible(f)$.\index{has-inverse(f)@{$\hasinverse(f)$}!has-inverse(f) to is-coh-invertible(f)@{$\hasinverse(f)\to\iscohinvertible(f)$}}
\end{lem}

\begin{proof}
  For each $y:B$, we construct the identification $G'(y)$ as the concatenation
  \begin{equation*}
    \begin{tikzcd}
      fg(y) \arrow[r,equals,"{G(fg(y))}^{-1}"] &[2.5em] fgfg(y) \arrow[r,equals,"\ap{f}{H(g(y))}"] &[2.5em] fg(y) \arrow[r,equals,"G(y)"] & y.
\end{tikzcd}
  \end{equation*}
  In order to construct a homotopy $f\cdot H \htpy G'\cdot f$, it suffices to show that the square
  \begin{equation*}
    \begin{tikzcd}[column sep=8em]
      fgfgf(x) \arrow[r,equals,"{G(fgf(x))}"] \arrow[d,equals,swap,"\ap{f}{H(gf(x))}"] & fgf(x) \arrow[d,equals,"\ap{f}{H(x)}"] \\
      fgf(x) \arrow[r,equals,swap,"G(f(x))"] & f(x)
    \end{tikzcd}
  \end{equation*}
  commutes for every $x:A$.
  Recall from \cref{defn:retraction_swap} that we have $H(gf(x))=\ap{gf}{H(x)}$. Using this identification, we see that it suffices to show that the square
  \begin{equation*}
    \begin{tikzcd}[column sep=8em]
      fgfgf(x) \arrow[r,equals,"(G\cdot f)(gf(x))"] \arrow[d,equals,swap,"\ap{fgf}{H(x)}"] & fgf(x) \arrow[d,equals,"\ap{f}{H(x)}"] \\
      fgf(x) \arrow[r,equals,swap,"(G\cdot f)(x)"] & f(x)
    \end{tikzcd}
  \end{equation*}
  commutes. Now we observe that this is just a naturality square the homotopy $G\cdot f:fgf\htpy f$, which commutes by \cref{defn:htpy_nat}.
\end{proof}

Now we put the pieces together to conclude that any equivalence has contractible fibers.

\begin{thm}\label{thm:contr_equiv}
Any equivalence is a contractible map.\index{equivalence!is a contractible map}\index{is a contractible map!equivalence}\index{is contractible!fiber of an equivalence}
\end{thm}

\begin{proof}
  We have seen in \cref{lem:contr-inv} that any coherently invertible map is a contractible map. Moreover, any equivalence has the structure of an invertible map by \cref{lem:inv_equiv}, and any invertible map is coherently invertible by \cref{lem:coherently-invertible}.
\end{proof}

The following corollary is very similar to \cref{thm:total_path}, which asserts that the type $\sm{x:A}a=x$ is contractible. However, we haven't yet established that the equivalence $(a=x)\simeq (x=a)$ induces an equivalence on total spaces. However, using the fact that equivalences are contractible maps we can give a direct proof.

\begin{cor}\label{cor:contr_path}
Let $A$ be a type, and let $a:A$. Then the type\index{is contractible!total space of opposite identity type}
\begin{equation*}
\sm{x:A}x=a
\end{equation*}
is contractible.
\end{cor}

\begin{proof}
By \cref{thm:id_equiv}, the identity function is an equivalence. Therefore, the fibers of the identity function are contractible by \cref{thm:contr_equiv}. Note that $\sm{x:A}x=a$ is exactly the fiber of $\idfunc[A]$ at $a:A$.
\end{proof}
\index{contractible map|)}

\begin{exercises}
  \exitem \label{ex:prop_contr}Show that if $A$ is contractible, then for any $x,y:A$ the identity type $x=y$ is also contractible.\index{contractible type!characterization of identity type}\index{is contractible!identity type of contractible type}\index{identity type!of a contractible type}\index{characterization of identity type!of a contractible type}
  \exitem \label{ex:contr_retr}Suppose that $A$ is a retract of $B$. Show that\index{contractible type!closed under retracts}
  \begin{equation*}
    \iscontr(B)\to\iscontr(A).
  \end{equation*}
  \exitem \label{ex:contr_equiv}
  \begin{subexenum}
  \item Show that for any type $A$, the map $\const_\ttt : A\to \unit$ is an equivalence if and only if $A$ is contractible.\index{contractible type!is equivalent to 1@{is equivalent to $\unit$}}
  \item Apply \cref{ex:3_for_2} to show that for any map $f:A\to B$, if any two of the three assertions\index{contractible type!3-for-2 property}\index{3-for-2 property!of contractible types}
    \begin{enumerate}
    \item $A$ is contractible
    \item $B$ is contractible
    \item $f$ is an equivalence
    \end{enumerate}
    hold, then so does the third.
  \end{subexenum}
  \exitem \label{ex:is-not-contractible-Fin}Show that $\Fin{k}$ is not contractible for all $k\neq 1$. 
  \exitem \label{ex:is-contr-prod}Show that for any two types $A$ and $B$, the following are equivalent:
  \index{contractible type!closed under cartesian products}
  \index{is contractible!factor of contractible cartesian product}
  \begin{enumerate}
  \item Both $A$ and $B$ are contractible.
  \item The type $A\times B$ is contractible.
  \end{enumerate}
  \exitem \label{ex:contr_in_sigma} Let $A$ be a contractible type with center of contraction $a:A$. Furthermore, let $B$ be a type family over $A$. Show that the map
  \begin{equation*}
    y\mapsto\pairr{a,y}:B(a)\to\sm{x:A}B(x)
  \end{equation*}
  is an equivalence.\index{left unit law!of Sigma-types@{of $\Sigma$-types}}\index{dependent pair type!left unit law}
  \exitem \label{ex:proj_fiber}Let $B$ be a family of types over $A$, and consider the projection map 
    \begin{equation*}
      \proj 1 : \big(\sm{x:A}B(x)\big)\to A.
    \end{equation*}
  \begin{subexenum}
  \item Show that for any $a:A$, the map
    \begin{equation*}
      \lam{((x,y),p)} \tr_B(p,y) : \fib{\proj 1}{a} \to B(a),
    \end{equation*}
    is an equivalence.\index{type family!fibers of projection map}
  \item Show that the following are equivalent:%
    \index{pr 1@{$\proj 1$}!of contractible family is an equivalence}%
    \index{is an equivalence!pr 1 of contractible family@{$\proj 1$ of contractible family}}
    \begin{enumerate}
    \item The projection map $\proj 1$ is an equivalence.
    \item The type $B(x)$ is contractible for each $x:A$.
    \end{enumerate}
  \item Consider a dependent function $b:\prd{x:A}B(x)$. Show that the following are equivalent:
    \begin{enumerate}
    \item The map
    \begin{equation*}
      \lam{x}(x,b(x)) : A \to \sm{x:A}B(x)
    \end{equation*}
    is an equivalence.
    \item The type $B(x)$ is contractible for each $x:A$.
    \end{enumerate}
  \end{subexenum}
  \exitem \label{ex:fib_replacement}Construct for any map $f:A\to B$ an equivalence $e:\eqv{A}{\sm{y:B}\fib{f}{y}}$ and a homotopy $H:f\htpy \proj 1\circ e$ witnessing that the triangle
  \begin{equation*}
    \begin{tikzcd}[column sep=0em]
      A \arrow[rr,"e"] \arrow[dr,swap,"f"] & & \sm{y:B}\fib{f}{y} \arrow[dl,"\proj 1"] \\
      \phantom{\sm{y:B}\fib{f}{y}} & B
    \end{tikzcd}
  \end{equation*}
  commutes. The projection $\proj 1 : (\sm{y:B}\fib{f}{y})\to B$ is sometimes also called the \define{fibrant replacement}\index{fibrant replacement|textbf} of $f$, because first projection maps are fibrations in the homotopy interpretation of type theory.
\end{exercises}
\index{contractible type|)}



\section{The fundamental theorem of identity types}\label{chap:fundamental}
\sectionmark{The fundamental theorem}

\index{fundamental theorem of identity types|(}
\index{characterization of identity type!fundamental theorem of identity types|(}
For many types it is useful to have a characterization of their identity types. For example, we have used a characterization of the identity types of the fibers of a map in order to conclude that any equivalence is a contractible map. The fundamental theorem of identity types is our main tool to carry out such characterizations, and with the fundamental theorem it becomes a routine task to characterize an identity type whenever that is of interest. We note that the fundamental theorem also appears as Theorem 5.8.4 in \cite{hottbook}.

In our first application of the fundamental theorem of identity types we show that any equivalence is an embedding. Embeddings are maps that induce equivalences on identity types, i.e., they are the homotopical analogue of injective maps. In our second application we characterize the identity types of coproducts.

Throughout this book we will encounter many more occasions to characterize identity types. For example, we will show in \cref{thm:eq_nat} that the identity type of the natural numbers is equivalent to its observational equality, and we will show in \cref{thm:eq-circle} that the loop space of the circle is equivalent to $\Z$.

In order to prove the fundamental theorem of identity types, we first prove the basic fact that a family of maps is a family of equivalences if and only if it induces an equivalence on total spaces. 

\subsection{Families of equivalences}

\index{family of equivalences|(}
\begin{defn}
Consider a family of maps
\begin{equation*}
f : \prd{x:A}B(x)\to C(x).
\end{equation*}
We define the map\index{tot(f)@{$\tot{f}$}}
\begin{equation*}
\tot{f}:\sm{x:A}B(x)\to\sm{x:A}C(x)
\end{equation*}
by $\lam{(x,y)}(x,f(x,y))$.
\end{defn}

\begin{lem}\label{lem:fib_total}
  For any family of maps $f:\prd{x:A}B(x)\to C(x)$ and any $t:\sm{x:A}C(x)$,
  there is an equivalence\index{fiber!of tot(f)@{of $\tot{f}$}}\index{tot(f)@{$\tot{f}$}!fiber}
  \begin{equation*}
    \eqv{\fib{\tot{f}}{t}}{\fib{f(\proj 1(t))}{\proj 2(t)}}.
  \end{equation*}
\end{lem}

\begin{proof}
  We first define
  \begin{equation*}
    \varphi : \prd{t:\sm{x:A}C(x)} \fib{\tot{f}}{t}\to\fib{f(\proj 1(t))}{\proj 2(t)}
  \end{equation*}
  by pattern matching by
  \begin{equation*}
    \varphi((x,f(x,y)),((x,y),\refl{}))\defeq(y,\refl{}).
  \end{equation*}

  For the proof that $\varphi(t)$ is an equivalence, for each $t:\sm{x:A}C(x)$, we construct a map
  \begin{equation*}
    \psi(t) : \fib{f(\proj 1(t))}{\proj 2(t)}\to\fib{\tot{f}}{t}
  \end{equation*}
  equipped with homotopies $G(t):\varphi(t)\circ\psi(t)\htpy\idfunc$ and $H(t):\psi(t)\circ\varphi(t)\htpy\idfunc$. Each of these definitions is given by pattern matching, as follows:
  \begin{align*}
    \psi((x,f(x,y)),(y,\refl{})) & \defeq ((x,y),\refl{}) \\
    G((x,f(x,y)),(y,\refl{})) & \defeq \refl{} \\
    H((x,f(x,y)),((x,y),\refl{})) & \defeq \refl{}.\qedhere
  \end{align*}
\end{proof}

\begin{thm}\label{thm:fib_equiv}
  Let $f:\prd{x:A}B(x)\to C(x)$ be a family of maps. The following are equivalent:
  \index{is an equivalence!total(f) of family of equivalences@{$\tot{f}$ of family of equivalences}}
  \index{tot(f)@{$\tot{f}$}!of family of equivalences is an equivalence}\index{is family of equivalences!if total(f) is an equivalence@{iff $\tot{f}$ is an equivalence}}
\begin{enumerate}
\item For each $x:A$, the map $f(x)$ is an equivalence. In this case we say that $f$ is a \define{family of equivalences}\index{family of equivalences|textbf}\index{equivalence!family of equivalences|textbf}.
\item The map $\tot{f}:\sm{x:A}B(x)\to\sm{x:A}C(x)$ is an equivalence.
\end{enumerate}
\end{thm}

\begin{proof}
By \cref{thm:equiv_contr,thm:contr_equiv} it suffices to show that $f(x)$ is a contractible map for each $x:A$, if and only if $\tot{f}$ is a contractible map. Thus, we will show that $\fib{f(x)}{c}$ is contractible if and only if $\fib{\tot{f}}{x,c}$ is contractible, for each $x:A$ and $c:C(x)$. However, by \cref{lem:fib_total} these types are equivalent, so the result follows by \cref{ex:contr_equiv}.
\end{proof}

Now consider the situation where we have a map $f:A\to B$, and a family $C$ over $B$. Then we have the map
\begin{equation*}
  \lam{(x,z)}(f(x),z):\sm{x:A}C(f(x))\to\sm{y:B}C(y).
\end{equation*}
We claim that this map is an equivalence when $f$ is an equivalence. The technique to prove this claim is the same as the technique we used in \cref{thm:fib_equiv}: first we note that the fibers are equivalent to the fibers of $f$, and then we use the fact that a map is an equivalence if and only if its fibers are contractible to finish the proof.

The converse of the following lemma does not hold. Why not?

\begin{lem}\label{lem:total-equiv-base-equiv}
  Consider a map $f:A\to B$, and let $C$ be a type family over $B$. If $f$ is an equivalence, then the map
  \begin{equation*}
    \sigma_f(C) \defeq\lam{(x,z)}(f(x),z):\sm{x:A}C(f(x))\to\sm{y:B}C(y)
  \end{equation*}
  is an equivalence.
\end{lem}

\begin{proof}
  We claim that for each $t:\sm{y:B}C(y)$ there is an equivalence
  \begin{equation*}
    \fib{\sigma_f(C)}{t}\simeq \fib{f}{\proj 1(t)}.
  \end{equation*}
  We obtain such an equivalence by constructing the following functions and homotopies:
  \begin{align*}
    \varphi(t) & : \fib{\sigma_f(C)}{t}\to\fib{f}{\proj 1 (t)} & \varphi((f(x),z),((x,z),\refl{})) & \defeq (x,\refl{}) \\
    \psi(t) & : \fib{f}{\proj 1(t)} \to\fib{\sigma_f(C)}{t} & \psi((f(x),z),(x,\refl{})) & \defeq ((x,z),\refl{}) \\
    G(t) & : \varphi(t)\circ\psi(t)\htpy\idfunc & G((f(x),z),(x,\refl{})) & \defeq \refl{} \\
    H(t) & : \psi(t)\circ\varphi(t)\htpy\idfunc & H((f(x),z),((x,z),\refl{})) & \defeq \refl{}.
  \end{align*}
  Now the claim follows, since we see that $\varphi$ is a contractible map if and only if $f$ is a contractible map.
\end{proof}

Now we use \cref{lem:total-equiv-base-equiv} to obtain a generalization of \cref{thm:fib_equiv}.

\begin{defn}\label{defn:toto}
  Consider a map $f:A\to B$ and a family of maps
  \begin{equation*}
    g:\prd{x:A}C(x)\to D(f(x)),
  \end{equation*}
  where $C$ is a type family over $A$, and $D$ is a type family over $B$. In this situation we also say that $g$ is a \define{family of maps over $f$}. Then we define\index{totf(g)@{$\tot[f]{g}$}}
  \begin{equation*}
    \tot[f]{g}:\sm{x:A}C(x)\to\sm{y:B}D(y)
  \end{equation*}
  by $\tot[f]{g}(x,z)\defeq (f(x),g(x,z))$.
\end{defn}

\begin{thm}\label{thm:equiv-toto}
  Suppose that $g$ is a family of maps over $f$ as in \cref{defn:toto}, and suppose that $f$ is an equivalence. Then the following are equivalent:
  \begin{enumerate}
  \item The family of maps $g$ over $f$ is a family of equivalences.
  \item The map $\tot[f]{g}$ is an equivalence.
  \end{enumerate}
\end{thm}

\begin{proof}
  Note that we have a commuting triangle
  \begin{equation*}
    \begin{tikzcd}[column sep=0]
      \sm{x:A}C(x) \arrow[rr,"{\tot[f]{g}}"] \arrow[dr,swap,"\tot{g}"]& & \sm{y:B}D(y) \\
      & \sm{x:A}D(f(x)) \arrow[ur,swap,"{\lam{(x,z)}(f(x),z)}"]
    \end{tikzcd}
  \end{equation*}
  By the assumption that $f$ is an equivalence, it follows that the map
  \begin{equation*}
    \sm{x:A}D(f(x))\to \sm{y:B}D(y)
  \end{equation*}
  is an equivalence. Therefore it follows that $\tot[f]{g}$ is an equivalence if and only if $\tot{g}$ is an equivalence. Now the claim follows, since $\tot{g}$ is an equivalence if and only if $g$ if a family of equivalences.
\end{proof}
\index{family of equivalences|)}

\subsection{The fundamental theorem}

\index{identity system|(}

The fundamental theorem of identity types (\cref{thm:id_fundamental}) is a general theorem that can be used to characterize the identity type of a given type. It describes necessary and sufficient conditions on a type family $B$ over a type $A$ equipped with a point $a:A$ to obtain an equivalence $(a=x)\simeq B(x)$ for each $x:A$.

One of those conditions is that the family $B$ satisfies an induction principle that is similar to the identification elimination principle. Such families are called \emph{identity systems}, which we will introduce now.

\begin{defn}\label{defn:identity-system}
  Let $A$ be a type equipped with a term $a:A$. A \define{(unary) identity system}\index{identity system|textbf}\index{unary identity system|see {identity system}} on $A$ at $a$ consists of a type family $B$ over $A$ equipped with $b:B(a)$, such that for any family of types $P(x,y)$ indexed by $x:A$ and $y:B(x)$,
  the function
  \begin{equation*}
    h\mapsto h(a,b):\Big(\prd{x:A}\prd{y:B(x)}P(x,y)\Big)\to P(a,b)
  \end{equation*}
  has a section.  
\end{defn}

In other words, if $B$ is an identity system on $A$ at $a$ and $P$ is a family of types indexed by $x:A$ and $y:B(x)$, then there is for each $p:P(a,b)$ a dependent function
\begin{equation*}
  f:\prd{x:A}\prd{y:B(x)}P(x,y)
\end{equation*}
such that $f(a,b)=p$. This is of course a variant of identification elimination, where the computation rule is given by an identification rather than as a judgmental equality.

We will state the fundamental theorem of identity types in a way that makes it maximally applicable. The fundamental theorem starts off with assuming a type $A$ equipped with a base point $a:A$, and a type family $B$ over $A$ equipped with a point $b:B(a)$. Furthermore it assumes an arbitrary family of maps
\begin{equation*}
  f:\prd{x:A}(a=x)\to B(x)
\end{equation*}
equipped with an identification $f(a,\refl{a})=b$. The theorem asserts conditions that are equivalent to $f$ being a family of equivalences.

In the setup of the fundamental theorem of identity types we can always construct the family of maps
\begin{equation*}
  f\defeq\pathind_a(b):\prd{x:A}(a=x)\to B(x)
\end{equation*}
for which the judgmental equality $f(a,\refl{a})\jdeq b$ holds. So you may wonder why we choose to formulate the fundamental theorem of identity types using a general family of maps $f$. The reason is that it is somewhat common to apply the fundamental theorem of identity types in order to conclude that $f$ is a family of equivalences, even when $f$ is not by definition the canonical family of maps, and we want to be free to do so.

The most important implication in the fundamental theorem is that (ii) implies (i). Occasionally we will also use the third equivalent statement.

\begin{thm}[The fundamental theorem of identity types]\label{thm:id_fundamental}
Let $A$ be a type with $a:A$, and let $B$ be a type family over $A$ equipped with a point $b:B(a)$. Furthermore, consider a family of maps\index{fundamental theorem of identity types}
\begin{equation*}
  f:\prd{x:A}(a=x)\to B(x)
\end{equation*}
equipped with an identification $f(a,\refl{a})=b$. Then the following are equivalent:
\begin{enumerate}
\item The family of maps $f$ is a family of equivalences.
\item The total space\index{is contractible!total space of an identity system}
\begin{equation*}
\sm{x:A}B(x)
\end{equation*}
is contractible.
\item The family $B$ equipped with $b:B(a)$ is an identity system.
\end{enumerate}
In particular, we see that for any $b:B(a)$, the canonical family of maps
\begin{equation*}
\pathind_a(b):\prd{x:A} (a=x)\to B(x)
\end{equation*}
is a family of equivalences if and only if $\sm{x:A}B(x)$ is contractible.
\end{thm}

\begin{proof}
  First we show that (i) and (ii) are equivalent.
  By \cref{thm:fib_equiv} it follows that the family of maps $f$ is a family of equivalences if and only if it induces an equivalence
  \begin{equation*}
    \eqv{\Big(\sm{x:A}a=x\Big)}{\Big(\sm{x:A}B(x)\Big)}
  \end{equation*}
  on total spaces. We have that $\sm{x:A}a=x$ is contractible, so it follows by \cref{ex:contr_equiv} that $\tot{f}$ is an equivalence if and only if $\sm{x:A}B(x)$ is contractible.

  Now we show that (ii) and (iii) are equivalent. Note that we have the following commuting triangle
  \begin{equation*}
    \begin{tikzcd}[column sep=0]
      \prd{t:\sm{x:A}B(x)}P(t) \arrow[rr,"\evpair"] \arrow[dr,swap,"{\evpt(a,b)}"] & & \prd{x:A}\prd{y:B(x)}P(x,y) \arrow[dl,"{\lam{h}h(a,b)}"] \\
      \phantom{\prd{x:A}\prd{y:B(x)}P(x,y)} & P(a,b)
    \end{tikzcd}
  \end{equation*}
  In this diagram the top map has a section. Therefore it follows by \cref{ex:3_for_2} that the left map has a section if and only if the right map has a section. Recall from \cref{defn:singleton-induction} that the type $\sm{x:A}B(x)$ satisfies singleton induction if and only if the left map in the triangle has a section for each $P$. Therefore we conclude our proof with \cref{thm:contractible}, which shows that the type $\sm{x:A}B(x)$ satisfies singleton induction if and only if it is contractible.
\end{proof}
\index{identity system|)}

\subsection{Equality on the natural numbers}
\index{natural numbers!observational equality|(}
\index{Eq N@{$\EqN$}|(}

As a first application of the fundamental theorem of identity types, we characterize the identity type of the natural numbers. We will use the observational equality $\EqN$ on $\N$. Recall from \cref{defn:obs_nat} that $\EqN$ is defined by
\begin{align*}
  \EqN(\zeroN,\zeroN) & \defeq \unit & \EqN(\zeroN,n+1) & \defeq \emptyt \\
  \EqN(m+1,\zeroN) & \defeq \emptyt & \EqN(m+1,n+1) & \defeq \EqN(m,n).
\end{align*}
This relation is an equivalence relation. In particular, the reflexivity term $\reflEqN(m):\EqN(m,m)$ is defined inductively by
\begin{align*}
  \reflEqN(\zeroN) & \defeq \ttt \\
  \reflEqN(m+1) & \defeq \reflEqN(m).
\end{align*}
Using the reflexivity term, we obtain a canonical map
\begin{equation*}
  (m=n)\to \EqN(m,n)
\end{equation*}
for every $m,n:\N$.

\begin{thm}\label{thm:eq_nat}
  For each $m,n:\N$, the canonical map\index{natural numbers!identity type}\index{natural numbers!characterization of identity type}\index{characterization of identity type!of N@{of $\N$}}\index{identity type!of the natural numbers}
  \begin{equation*}
    (m=n)\to \EqN(m,n)
  \end{equation*}
  is an equivalence.
\end{thm}

\begin{proof}
  By \cref{thm:id_fundamental} it suffices to show that the type
  \begin{equation*}
    \sm{n:\N}\EqN(m,n)
  \end{equation*}
  is contractible, for each $m:\N$. The center of contraction is defined to be $(m,\reflEqN(m))$.

  The contraction
  \begin{equation*}
    \gamma(m):\prd{n:\N}\prd{e:\EqN(m,n)}(m,\reflEqN(m))=(n,e)
  \end{equation*}
  is defined for each $m$ by induction on $m,n:\N$. In the base case we define
  \begin{equation*}
    \gamma(\zeroN,\zeroN,\ttt)\defeq \refl{}.
  \end{equation*}
  If one of $m$ and $n$ is zero and the other is a successor, then the type $\EqN(m,n)$ is empty, so the desired path can be obtained via the induction principle of the empty type.

  The inductive step remains, in which we have to define the identification
  \begin{equation*}
    \gamma(m+1,n+1,e):(m+1,\reflEqN(m+1))=(n+1,e)
  \end{equation*}
  for each $m,n:\N$ equipped with $e:\EqN(m,n)$. We first observe that there is a map
  \begin{equation*}
    \begin{tikzcd}
      \Big(\sm{n:\N}\EqN(m,n)\Big) \arrow[r,"f"] & \Big(\sm{n:\N}\EqN(m+1,n)\Big)
    \end{tikzcd}
  \end{equation*}
  given by $(n,e)\mapsto (n+1,e)$. With this definition of $f$ we have
  \begin{equation*}
    f(m,\reflEqN(m))\jdeq (m+1,\reflEqN(m+1)).
  \end{equation*}
  Therefore we can define
  \begin{equation*}
    \gamma(m+1,n+1,e)\defeq \ap{f}{\gamma(m,n,e)}.\qedhere
  \end{equation*}
\end{proof}
\index{natural numbers!observational equality|)}
\index{Eq N@{$\EqN$}|)}

\subsection{Embeddings}
\index{embedding|(}
In our second application of the fundamental theorem we show that equivalences are embeddings. The notion of embedding is the homotopical analogue of the set theoretic notion of injective map.

\begin{defn}
An \define{embedding}\index{embedding|textbf} is a map $f:A\to B$\index{is an embedding|textbf} that satisfies the property that\index{is an equivalence!action on paths of an embedding}
\begin{equation*}
\apfunc{f}:(\id{x}{y})\to(\id{f(x)}{f(y)})
\end{equation*}
is an equivalence, for every $x,y:A$. We write $\isemb(f)$\index{is-emb(f)@{$\isemb(f)$}} for the type of witnesses that $f$ is an embedding, and we define\index{A hookrightarrow B@{$A\hookrightarrow B$}|see {embedding}}
\begin{equation*}
  A\hookrightarrow B\defeq \sm{f:A\to B}\isemb(f).
\end{equation*}
\end{defn}

Another way of phrasing the following statement is that equivalent types have equivalent identity types.

\begin{thm}
\label{cor:emb_equiv} 
Any equivalence is an embedding.\index{is an embedding!equivalence}\index{equivalence!is an embedding}
\end{thm}

\begin{proof}
Let $e:\eqv{A}{B}$ be an equivalence, and let $x:A$. Our goal is to show that
\begin{equation*}
\apfunc{e} : (\id{x}{y})\to (\id{e(x)}{e(y)})
\end{equation*}
is an equivalence for every $y:A$. By \cref{thm:id_fundamental} it suffices to show that 
\begin{equation*}
\sm{y:A}e(x)=e(y)
\end{equation*}
is contractible. Now observe that there is an equivalence
\begin{samepage}
\begin{align*}
\sm{y:A}e(x)=e(y) & \eqvsym \sm{y:A}e(y)=e(x) \\
& \jdeq \fib{e}{e(x)}
\end{align*}
\end{samepage}
by \cref{thm:fib_equiv}, since for each $y:A$ the map
\begin{equation*}
\invfunc : (e(x)=e(y))\to (e(y)= e(x))
\end{equation*}
is an equivalence by \cref{ex:equiv_grpd_ops}.
The fiber $\fib{e}{e(x)}$ is contractible by \cref{thm:contr_equiv}, so it follows by \cref{ex:contr_equiv} that the type $\sm{y:A}e(x)=e(y)$ is indeed contractible.
\end{proof}
\index{embedding|)}

\subsection{Disjointness of coproducts}

\index{disjointness of coproducts|(}
\index{characterization of identity type!of coproducts|(}
\index{identity type!of a coproduct|(}
\index{coproduct!characterization of identity type|(}
\index{coproduct!disjointness|(}
In our third application of the fundamental theorem of identity types, we characterize the identity types of coproducts. Our goal in this section is to prove the following theorem.

\begin{thm}\label{thm:id-coprod-compute}
Let $A$ and $B$ be types. Then there are equivalences\index{identity type!of a coproduct}\index{characterization of identity type!of coproducts}\index{coproduct!characterization of identity type}
\begin{align*}
(\inl(x)=\inl(x')) & \eqvsym (x = x')\\
(\inl(x)=\inr(y')) & \eqvsym \emptyt \\
(\inr(y)=\inl(x')) & \eqvsym \emptyt \\
(\inr(y)=\inr(y')) & \eqvsym (y=y')
\end{align*}
for any $x,x':A$ and $y,y':B$.
\end{thm}

In order to prove \cref{thm:id-coprod-compute}, we first define
a binary relation $\Eqcoprod_{A,B}$ on the coproduct $A+B$.

\begin{defn}
Let $A$ and $B$ be types. We define\index{Eq-coprod@{$\Eqcoprod_{A,B}$}|textbf}\index{coproduct!Eq-coprod@{$\Eqcoprod_{A,B}$}|textbf}
\begin{equation*}
\Eqcoprod_{A,B} : (A+B)\to (A+B)\to\UU
\end{equation*}
by double induction on the coproduct, postulating
\begin{align*}
\Eqcoprod_{A,B}(\inl(x),\inl(x')) & \defeq (x=x') \\
\Eqcoprod_{A,B}(\inl(x),\inr(y')) & \defeq \emptyt \\
\Eqcoprod_{A,B}(\inr(y),\inl(x')) & \defeq \emptyt \\
\Eqcoprod_{A,B}(\inr(y),\inr(y')) & \defeq (y=y').
\end{align*}
The relation $\Eqcoprod_{A,B}$ is also called the \define{observational equality of coproducts}\index{observational equality!on coproduct types}\index{coproduct!observational equality}.
\end{defn}

\begin{lem}
The observational equality relation $\Eqcoprod_{A,B}$ on $A+B$ is reflexive, and therefore there is a map
\begin{equation*}
\Eqcoprodeq:\prd{s,t:A+B} (s=t)\to \Eqcoprod_{A,B}(s,t).
\end{equation*}
\end{lem}

\begin{constr}
The reflexivity term $\rho$ is constructed by induction on $t:A+B$, using
\begin{align*}
\rho(\inl(x))\defeq \refl{x}  & : \Eqcoprod_{A,B}(\inl(x),\inl(x)) \\
\rho(\inr(y))\defeq \refl{y} & : \Eqcoprod_{A,B}(\inr(y),\inr(y)).\qedhere
\end{align*}
\end{constr}

To show that $\Eqcoprodeq$ is a family of equivalences, we will use the fundamental theorem of identity types, \cref{thm:id_fundamental}. Therefore, we need to prove the following proposition.

\begin{prp}\label{lem:is-contr-total-eq-coprod}
For any $s:A+B$ the total space
\begin{equation*}
\sm{t:A+B}\Eqcoprod_{A,B}(s,t)
\end{equation*}
is contractible.
\end{prp}

\begin{proof}
  For convenience, let us write $E\defeq \Eqcoprod_{A,B}$. By induction on $s$, it suffices to show that the total spaces
  \begin{equation*}
    \sm{t:A+B}E(\inl(x),t) \qquad\text{and}\qquad \sm{t:A+B}E(\inr(y),t)
  \end{equation*}
  are contractible. The two proofs are similar, so we only prove that the type on the left is contractible. By the laws of coproducts and $\Sigma$-types given in \cref{eg:laws-products-coproducts,eg:laws-Sigma-types}, we simply compute
  \begin{samepage}
    \begin{align*}
      & \sm{t:A+B}E(\inl(x),t) \\
      & \eqvsym \Big(\sm{x':A}E(\inl(x),\inl(x'))\Big)+\Big(\sm{y':B}E(\inl(x),\inr(y'))\Big) \\
      & \eqvsym \Big(\sm{x':A}x=x'\Big)+\Big(\sm{y':B}\emptyt\Big) \\
      & \eqvsym \sm{x':A}x=x'.
    \end{align*}%
  \end{samepage}%
  The last type in this computation is contractible by \cref{thm:total_path}, so we conclude that the total space of $E(\inl(x))$ is contractible.
\end{proof}

\begin{proof}[Proof of \cref{thm:id-coprod-compute}]
  The proof is now concluded with an application of \cref{thm:id_fundamental}, using \cref{lem:is-contr-total-eq-coprod}.
\end{proof}
\index{disjointness of coproducts|)}
\index{characterization of identity type!of coproducts|)}
\index{identity type!of a coproduct|)}
\index{coproduct!characterization of identity type|)}
\index{coproduct!disjointness|)}

\subsection{The structure identity principle}\label{sec:structure-identity-principle}
\index{structure identity principle|(}

We often encounter a type consisting of certain objects equipped with further structure. For example, the fiber of a map $f:A\to B$ at $b:B$ is the type of elements $a:A$ equipped with an identification $p:f(a)=b$. Such \emph{structure} types occur all over mathematics, and it is important to have an efficient characterization of their identity types. A general structure type is just a $\Sigma$-type, and we're asking for a characterization of its identity type.

Recall from \cref{thm:eq_sigma} that the identity type of the type $\sm{x:A}B(x)$ at a pair $(a,b)$ can be characterized as
\begin{equation*}
  ((a,b)=(x,y))\simeq \sm{p:a=x}\tr_B(p,b)=y.
\end{equation*}
However, this characterization of the identity type of $\sm{x:A}B(x)$ is not as clear and useful as we like it to be, because it uses the transport function, which is completely generic. Our plan is to use identity systems on $A$ and on $B(a)$ to arrive at a more useful characterization of the identity type of $\sm{x:A}B(x)$.

In order to abstract away this characterization of the identity type of $\sm{x:A}B(x)$, let $C:A\to\UU$ be the family of types given by $C(x)\defeq (a=x)$, and let
\begin{equation*}
  D:\prd{x:A}B(x)\to(C(x)\to\UU)
\end{equation*}
be the family of types given by $D(x,y,p)\defeq \tr_B(p,b)=y$. Then $C$ is an identity system on $A$ at $a$, and the type family $y\mapsto D(a,y,\refl{})$ is an identity system on $B(a)$ at $b$. This suggests the following definition of dependent identity systems.

\begin{defn}
  Consider a type $A$ equipped with an identity system $C$ based at $a:A$, and let $c:C(a)$. Furthermore, consider a type family $B$ over $A$. A \define{dependent identity system}\index{dependent identity system|textbf}\index{identity system!dependent identity system|textbf} over $C$ at $b:B(a)$ consists of a type family
  \begin{equation*}
    D : \prd{x:A} B(x) \to (C(x)\to \UU)
  \end{equation*}
  equipped with an element $d:D(a,b,c)$ such that $y\mapsto D(a,y,c)$ is an identity system at $b$.
\end{defn}

\begin{thm}[Structure identity principle]\label{thm:structure-identity-principle}
  Consider a type family $B$ over $A$, elements $a:A$ and $b:B(a)$, and an identity system $C$ of $A$ with $c:C(a)$. Furthermore, consider a type family
  \begin{equation*}
    D : \prd{x:A} B(x) \to (C(x)\to \UU)
  \end{equation*}
  equipped with an element $d:D(a,b,c)$. Then the following are equivalent:
  \begin{enumerate}
  \item Any family of maps
    \begin{equation*}
      (b=y)\to D(a,y,c)
    \end{equation*}
    indexed by $y:B(a)$ is a family of equivalences.
  \item The total space
    \begin{equation*}
      \sm{y:B(a)}D(a,y,c)
    \end{equation*}
    is contractible.
  \item $D$ is a dependent identity system over $C$ at $b:B(a)$.
  \item Any family of maps
    \begin{equation*}
      ((a,b)=(x,y))\to \sm{z:C(x)}D(x,y,z))
    \end{equation*}
    indexed by $(x,y):\sm{x:A}B(x)$ is a family of equivalences.
  \item The total space
    \begin{equation*}
      \sm{(x,y):\sm{x:A}B(x)}\sm{z:C(x)}D(x,y,z)
    \end{equation*}
    is contractible.
  \item The type family
    \begin{equation*}
      (x,y)\mapsto \sm{z:C(x)}D(x,y,z)
    \end{equation*}
    is an identity system at $(a,b):\sm{x:A}B(x)$.
  \end{enumerate}
\end{thm}

\begin{proof}
  The first three statements as well as the last three statements are equivalent by \cref{thm:id_fundamental}. Therefore it suffices to show that (ii) and (v) are equivalent. Note that there is an equivalence
  \begin{multline*}
    \sm{(x,y):\sm{x:A}B(x)}\sm{z:C(x)}D(x,y,z) \\
    \simeq
    \sm{(x,z):\sm{x:A}C(x)}\sm{y:B(x)}D(x,y,z).
  \end{multline*}
  This equivalence, its inverse, and the homotopies witnessing that the inverse is indeed an inverse are all straightforward to construct using pattern matching. Furthermore, notice that the type $\sm{x:A}C(x)$ is contractible with center of contraction $(a,c)$ since $C$ is assumed to be an identity system at $a:A$. Therefore it follows that
  \begin{equation*}
    \sm{(x,y):\sm{x:A}B(x)}\sm{z:C(x)}D(x,y,z)\simeq\sm{y:B(a)}D(a,y,c).\qedhere
  \end{equation*}
\end{proof}

\begin{eg}
  By the structure identity principle of \cref{thm:structure-identity-principle} in combination with the fundamental theorem of identity types (\cref{thm:id_fundamental}), it becomes completely routine to characterize identity types of structures: We only have to show that the types
  \begin{equation*}
    \sm{x:A}C(x)\qquad\text{and}\qquad\sm{y:B(a)}D(a,y,c)
  \end{equation*}
  are contractible. To illustrate this use of the structure identity principle, we give an alternative characterization of the fiber of a map $f:A \to B$ at $b:B$. We claim that\index{identity type!of a fiber}\index{fiber!characterization of identity type}\index{characterization of identity type!of the fiber of a map}
  \begin{align*}
    ((x,p)=(y,q)) & \simeq \fib{\apfunc{f}}{\ct{p}{q^{-1}}} \\
                  & \jdeq \sm{\alpha:x=y}\ap{f}{\alpha}=\ct{p}{q}^{-1}.
  \end{align*}
  To see this, we apply \cref{thm:structure-identity-principle}. Note that $\sm{y:A}x=y$ is contractible by \cref{thm:total_path} with center of contraction $(x,\refl{f(x)})$. Therefore it suffices to show that the type
  \begin{equation*}
    \sm{q:f(x)=b}\refl{f(x)}=\ct{p}{q}^{-1}
  \end{equation*}
  is contractible. Of course, this type is equivalent to $\sm{q:f(x)=b}p=q$, which is again contractible by \cref{thm:total_path}.
\end{eg}
\index{structure identity principle|)}

\begin{exercises}
  \exitem
  \begin{subexenum}
  \item \label{ex:is-emb-empty}Show that the map $\emptyt\to A$ is an embedding for every type $A$.\index{is an embedding!0 to A@{$\emptyt\to A$}}
  \item \label{ex:is-emb-inl-inr}Show that $\inl:A\to A+B$ and $\inr:B\to A+B$ are embeddings for any two types $A$ and $B$.
    \index{is an embedding!inl (for coproducts)@{$\inl$ (for coproducts)}}
    \index{is an embedding!inr (for coproducts)@{$\inr$ (for coproducts)}}
    \index{inl@{$\inl$}!is an embedding}
    \index{inr@{$\inr$}!is an embedding}
  \item Show that $\inl:A\to A+B$ is an equivalence if and only if $B$ is empty, and that $\inr : B \to A+B$ is an equivalence if and only if $A$ is empty.
  \end{subexenum}
  \exitem Consider an equivalence $e:A\simeq B$. Construct an equivalence
  \begin{equation*}
    p\mapsto \tilde{p}:(e(x)=y)\simeq(x=e^{-1}(y))
  \end{equation*}
  for every $x:A$ and $y:B$, such that the triangle
  \begin{equation*}
    \begin{tikzcd}[column sep=large]
      e(x) \arrow[r,equals,"\ap{e}{\tilde{p}}"] \arrow[dr,equals,swap,"p"] & e(e^{-1}(y)) \arrow[d,equals,"G(y)"] \\
      & y
    \end{tikzcd}
  \end{equation*}
  commutes for every $p:e(x)=y$. In this diagram, the homotopy $G:e\circ e^{-1}\htpy \idfunc$ is the homotopy witnessing that $e^{-1}$ is a section of $e$.
  \exitem Show that\index{embedding!closed under homotopies}
  \begin{equation*}
    (f\htpy g)\to (\isemb(f)\leftrightarrow\isemb(g))
  \end{equation*}
  for any $f,g:A\to B$.
  \exitem \label{ex:emb_triangle}Consider a commuting triangle
  \begin{equation*}
    \begin{tikzcd}[column sep=tiny]
      A \arrow[rr,"h"] \arrow[dr,swap,"f"] & & B \arrow[dl,"g"] \\
      & X
    \end{tikzcd}
  \end{equation*}
  with $H:f\htpy g\circ h$. 
  \begin{subexenum}
  \item Suppose that $g$ is an embedding. Show that $f$ is an embedding if and only if $h$ is an embedding.\index{is an embedding!composite of embeddings}\index{is an embedding!right factor of embedding if left factor is an embedding}
  \item Suppose that $h$ is an equivalence. Show that $f$ is an embedding if and only if $g$ is an embedding.\index{is an embedding!left factor of embedding if right factor is an equivalence}
  \end{subexenum}
  \exitem Consider two embeddings $f:A\hookrightarrow B$ and $g:B\hookrightarrow C$. Show that the following are equivalent:
  \begin{enumerate}
  \item The composite $g\circ f$ is an equivalence.
  \item Both $f$ and $g$ are equivalences.
  \end{enumerate}
  \exitem Consider two maps $f:A\to C$ and $g:B\to C$. Use \cref{ex:is-emb-inl-inr} to show that the following are equivalent:
  \begin{enumerate}
  \item The map $[f,g]:A+B\to C$ is an embedding.
  \item Both $f$ and $g$ are embeddings, and
    \begin{equation*}
      f(a)\neq g(b)
    \end{equation*}
    for all $a:A$ and $b:B$.
  \end{enumerate}
  \exitem \label{ex:is-equiv-is-equiv-functor-coprod}Consider two maps $f:A\to A'$ and $g:B \to B'$.
  \begin{subexenum}
  \item Show that if the map
    \begin{equation*}
      f+g:(A+B)\to (A'+B')
    \end{equation*}
    is an equivalence, then so are both $f$ and $g$ (this is the converse of \cref{ex:coproduct_functor_equivalence}).
  \item \label{ex:is-emb-coprod}Show that $f+g$ is an embedding if and only if both $f$ and $g$ are embeddings.
  \end{subexenum}
  \exitem \label{ex:id_fundamental_retr}
  \begin{subexenum}
  \item Let $f,g:\prd{x:A}B(x)\to C(x)$ be two families of maps. Show that
    \begin{equation*}
      \Big(\prd{x:A}f(x)\htpy g(x)\Big)\to \Big(\tot{f}\htpy \tot{g}\Big). 
    \end{equation*}
  \item Let $f:\prd{x:A}B(x)\to C(x)$ and let $g:\prd{x:A}C(x)\to D(x)$. Show that
    \begin{equation*}
      \tot{\lam{x}g(x)\circ f(x)}\htpy \tot{g}\circ\tot{f}.
    \end{equation*}
  \item For any family $B$ over $A$, show that
    \begin{equation*}
      \tot{\lam{x}\idfunc[B(x)]}\htpy\idfunc.
    \end{equation*}
  \item Let $a:A$, and let $B$ be a type family over $A$. Use \cref{ex:contr_retr} to show that if each $B(x)$ is a retract of $\id{a}{x}$, then $B(x)$ is equivalent to $\id{a}{x}$ for every $x:A$.
    \index{fundamental theorem of identity types!formulation with retractions}
  \item Conclude that for any family of maps
    \index{fundamental theorem of identity types!formulation with sections}
    \begin{equation*}
      f : \prd{x:A} (a=x) \to B(x),
    \end{equation*}
    if each $f(x)$ has a section, then $f$ is a family of equivalences.
  \end{subexenum} 
  \exitem Use \cref{ex:id_fundamental_retr} to show that for any map $f:A\to B$, if
  \begin{equation*}
    \apfunc{f} : (x=y) \to (f(x)=f(y))
  \end{equation*}
  has a section for each $x,y:A$, then $f$ is an embedding.\index{is an embedding!if the action on paths have sections}
  \exitem \label{ex:path-split}(Shulman) We say that a map $f:A\to B$ is \define{path-split}\index{path-split|textbf} if $f$ has a section, and for each $x,y:A$ the map
  \begin{equation*}
    \apfunc{f}(x,y):(x=y)\to (f(x)=f(y))
  \end{equation*}
  also has a section. We write $\pathsplit(f)$\index{path-split(f)@{$\pathsplit(f)$}|textbf} for the type
  \begin{equation*}
    \sections(f)\times\prd{x,y:A}\sections(\apfunc{f}(x,y)).
  \end{equation*}
  Show that for any map $f:A\to B$ the following are equivalent:
  \begin{enumerate}
  \item The map $f$ is an equivalence.
  \item The map $f$ is path-split.\index{is an equivalence!path-split map}
  \end{enumerate}
  \exitem \label{ex:fiber_trans}Consider a triangle
  \begin{equation*}
    \begin{tikzcd}[column sep=small]
      A \arrow[rr,"h"] \arrow[dr,swap,"f"] & & B \arrow[dl,"g"] \\
      & X
    \end{tikzcd}
  \end{equation*}
  with a homotopy $H:f\htpy g\circ h$ witnessing that the triangle commutes. 
  \begin{subexenum}
  \item Construct a family of maps
    \begin{equation*}
      \fibtriangle(h,H):\prd{x:X}\fib{f}{x}\to\fib{g}{x},
    \end{equation*}
    for which the square
    \begin{equation*}
      \begin{tikzcd}[column sep=8em]
        \sm{x:X}\fib{f}{x} \arrow[r,"\tot{\fibtriangle(h,H)}"] \arrow[d] & \sm{x:X}\fib{g}{x} \arrow[d] \\
        A \arrow[r,swap,"h"] & B
      \end{tikzcd}
    \end{equation*}
    commutes, where the vertical maps are as constructed in \cref{ex:fib_replacement}.
  \item Show that $h$ is an equivalence if and only if $\fibtriangle(h,H)$ is a family of equivalences.
  \end{subexenum}
\end{exercises}
\index{fundamental theorem of identity types|)}
\index{characterization of identity type!fundamental theorem of identity types|)}



\section{Propositions, sets, and the higher truncation levels}
\sectionmark{Truncation levels}\label{chap:hierarchy}

The set theoretic foundations of mathematics arise in two stages. The first stage is to specify the formal system of first order logic; the second stage is to give an axiomatization of set theory in this formal system. Unlike set theory, type theory is its own formal system. The logic of dependent types, as given by the inference rules, is all we need.

However, even though type theory is not built upon a separate system of logic such as first order logic, we can find logic in type theory by recognizing certain types as propositions. Note that the propositions of first order logic have a virtue that could be rather useful sometimes: First order logic does not offer any way to distinguish between any two proofs of the same proposition. Therefore we say that propositions in type theory are those types that have at most one element.

This condition can be expressed with the identity type: any two elements must be equal. Examples of such types include the empty type $\emptyt$ and the unit type $\unit$. We call such types propositions. Propositions are useful, because if we know that a certain type is a proposition, then we know that any of its inhabitants are equal. Many important conditions, such as the condition that a map is an equivalence, will turn out to be propositions. This fact implies that two equivalences $A\simeq B$ are equal if and only if their underlying maps $A\to B$ are equal. However, the claim that being an equivalence is a proposition requires function extensionality, the topic of the next section.

In this section we use the idea of propositions in a different way. After we establish some basic properties of propositions, we will introduce the \emph{sets} as the types of which the identity types are propositions. This is again reminiscent of the situation in set theory, where equality is a predicate in first order logic. We will see in \cref{eg:is-set-nat} that the type of natural numbers is a set.

Next, one might ask about the types of which the identity types are \emph{sets}. Such types are called \emph{$1$-types}. There is an entire hierarchy of special types that arises this way, where a type is said to be a $(k+1)$-type if its identity types are $k$-types. Since the identity types of the $1$-types are sets, we see that sets are in fact $0$-types. Most of mathematics takes place at this level, the level of sets. The types in higher levels, as well as types that do not belong to any finite level in this hierarchy, are studied extensively in synthetic homotopy theory.

However, we can also go a step further down: Since the identity types of sets are propositions, we see that the propositions are $(-1)$-types. Moreover, the identity types of propositions are contractible. Hence we find at the bottom of this hierarchy the contractible types as the $(-2)$-types. There is no point in going down further, since we have seen in \cref{ex:prop_contr} that the identity types of contractible types are again contractible.

\index{truncated type|(}
\index{truncation level|(}

\subsection{Propositions}

\index{proposition|(}
\begin{defn}
A type $A$ is said to be a \define{proposition}\index{proposition|textbf} if its identity types are contractible, i.e., if it comes equipped with a term of type\index{is-prop(A)@{$\isprop(A)$}|textbf}
\begin{equation*}
\isprop(A)\defeq\prd{x,y:A}\iscontr(x=y).
\end{equation*}
Given a universe $\UU$, we define $\prop_\UU$\index{Prop@{$\prop_\UU$}|textbf} to be the type of all small propositions, i.e.,
\begin{equation*}
  \prop_\UU\defeq\sm{X:\UU}\isprop(X).
\end{equation*}
\end{defn}

\begin{eg}\label{eg:prop_contr}
  Any contractible type is a proposition by \cref{ex:prop_contr}\index{contractible type!is a proposition}\index{is a proposition!contractible type}. In particular, the unit type is a proposition. The empty type is also a proposition, since we have\index{empty type!is a proposition}\index{is a proposition!empty type}
\begin{equation*}
\prd{x,y:\emptyt}\iscontr(x=y)
\end{equation*}
by the induction principle of the empty type.
\end{eg}

There are many conditions on a type $A$ that are equivalent to the condition that $A$ is a proposition. In the following proposition we state four such conditions.

\begin{prp}\label{lem:isprop_eq}
  Let $A$ be a type. Then the following are equivalent:
  \begin{enumerate}
  \item The type $A$ is a proposition.
  \item Any two terms of type $A$ can be identified, i.e., there is a dependent function of type\index{is-prop'(A)@{$\isprop'(A)$}}\index{is-prop(A)@{$\isprop(A)$}!is-prop(A) iff is-prop'(A)@{$\isprop(A)\leftrightarrow\isprop'(A)$}}
    \begin{equation*}
      \isprop'(A)\defeq\prd{x,y:A}\id{x}{y}.
    \end{equation*}
  \item The type $A$ is contractible as soon as it is inhabited, i.e., there is a function of type\index{is-prop(A)@{$\isprop(A)$}!is-prop(A) iff A to is-contr(A)@{$\isprop(A)\leftrightarrow(A\to\iscontr(A))$}}
    \begin{equation*}
      A \to \iscontr(A).
    \end{equation*}
  \item The map $\const_\ttt : A\to\unit$ is an embedding.\index{is-prop(A)@{$\isprop(A)$}!is-prop(A) iff is-emb(const star)@{$\isprop(A)\leftrightarrow\isemb(\const_\ttt)$}}
  \end{enumerate}
\end{prp}

\begin{proof}
  If $A$ is a proposition, then we can use the center of contraction of the identity types of $A$ to identify any two terms in $A$. This shows that (i) implies (ii).

  To show that (ii) implies (iii), suppose that $A$ comes equipped with $p:\prd{x,y:A}\id{x}{y}$. Then for any $x:A$ the dependent function $p(x):\prd{y:A}\id{x}{y}$ is a contraction of $A$. Thus we obtain the function
  \begin{equation*}
    \lam{x}(x,p(x)):A\to\iscontr(A).
  \end{equation*}

  To show that (iii) implies (iv), suppose that $A\to\iscontr(A)$. We first make the simple observation that
  \begin{equation*}
    (X\to \isemb(f))\to \isemb(f)
  \end{equation*}
  for any map $f:X\to Y$, so it suffices to show that $A\to\isemb(\const_\ttt)$. However, assuming we have $x:A$, it follows by assumption that $A$ is contractible. Therefore, it follows by \cref{ex:contr_equiv} that the map $\const_\ttt:A\to\unit$ is an equivalence, and any equivalence is an embedding by \cref{cor:emb_equiv}.

  To show that (iv) implies (i), note that if $A\to\unit$ is an embedding, then the identity types of $A$ are equivalent to contractible types and therefore they must be contractible.
\end{proof}

One useful feature of propositions, is that in order to construct an equivalence $e:P\simeq Q$ between propositions, it suffices to construct maps back and forth between them.

\begin{prp}\label{prp:equiv-prop}
  A map $f:P\to Q$ between two propositions $P$ and $Q$ is an equivalence if and only if there is a map $g:Q\to P$. Consequently, we have for any two propositions $P$ and $Q$ that
  \begin{equation*}
    (P\simeq Q) \leftrightarrow (P\leftrightarrow Q).
  \end{equation*}
\end{prp}

\begin{proof}
  Of course, if we have an equivalence $e:P\simeq Q$, then we get maps back and forth between $P$ and $Q$. Therefore it remains to show that
  \begin{equation*}
    (P\leftrightarrow Q) \to (P\simeq Q).
  \end{equation*}
  Suppose we have $f:P\to Q$ and $g:Q\to P$. Then we obtain the homotopies $f\circ g\htpy \idfunc$ and $g\circ f\htpy \idfunc$ by the fact that any two elements in $P$ and $Q$ can be identified. Therefore $f$ is an equivalence with inverse $g$. 
\end{proof}
\index{proposition|)}

\subsection{Subtypes}
\index{subtype|(}

  In set theory, a set $y$ is said to be a subset of a set $x$, if any element of $y$ is an element of $x$, i.e., if the condition
  \begin{equation*}
    \forall_z\, (z\in y)\to (z\in x)
  \end{equation*}
  holds. We have already noted that type theory is different from set theory in that terms in type theory come equipped with a \emph{unique} type. Moreover, in set theory the proposition $x\in y$ is well-formed for any two sets $x$ and $y$, whereas in type theory we can only judge that $a:A$ by applying the rules of inference of type theory in such a manner that we arrive at the conclusion that $a:A$. Because of these differences we must find a different way to talk about subtypes.

  Note that in set theory there is a correspondence between the subsets of a set $x$, and the \emph{predicates} on $x$. A predicate on $x$ is just a proposition $P(z)$ that varies over the elements $z\in x$. Indeed, if $y$ is a subset of $x$, then the corresponding predicate is the proposition $z\in y$. Conversely, if $P$ is a predicate on $x$, then we obtain the subset
  \begin{equation*}
    \{z\in x\mid P(z)\}
  \end{equation*}
  of $x$. This observation suggests that in type theory we should define a subtype of a type $A$ to be a family of propositions over $A$.

\begin{defn}
A type family $B$ over $A$ is said to be a \define{subtype}\index{subtype|textbf} of $A$ if for each $x:A$ the type $B(x)$ is a proposition. When $B$ is a subtype of $A$, we also say that $B(x)$ is a \define{property}\index{property|textbf} of $x:A$.
\end{defn}

One reason why subtypes are important and useful, is that for any
\begin{equation*}
  (x,p),(y,q):\sm{x:A}P(x)
\end{equation*}
in a subtype of $A$, we have $(x,p)=(y,q)$ if and only if $x=y$. In other words, two terms of a subtype of $A$ are equal if and only if they are equal as terms of $A$. This fact is properly expressed using embeddings: we claim that the projection map
\begin{equation*}
  \proj 1 : \Big(\sm{x:A}P(x)\Big)\to A
\end{equation*}
is an embedding, for any subtype $P$ of $A$. This claim can be strengthened slightly. We will prove the following two closely related facts:
\begin{enumerate}
\item A map $f:A\to B$ is an embedding if and only if its fibers are propositions.
\item A family of types $B$ over $A$ is a subtype of $A$ if and only if the projection map
  \begin{equation*}
    \Big(\sm{x:A}B(x)\Big)\to A
  \end{equation*}
  is an embedding.
\end{enumerate}
The first fact is analogous to the fact that a map is an equivalence if and only if its fibers are contractible, which we saw in \cref{thm:contr_equiv,thm:equiv_contr}. To prove the above claims, we will need that propositions are closed under equivalences.

\begin{lem}\label{lem:prop_equiv}
Let $A$ and $B$ be types, and let $e:\eqv{A}{B}$. Then we have\index{proposition!closed under equivalences}
\begin{equation*}
\isprop(A)\leftrightarrow\isprop(B).
\end{equation*}
\end{lem}

\begin{proof}
We will show that $\isprop(B)$ implies $\isprop(A)$. This suffices, because the converse follows from the fact that $e^{-1}:B\to A$ is also an equivalence. 

Since $e$ is assumed to be an equivalence, it follows by \cref{cor:emb_equiv} that
\begin{equation*}
\apfunc{e} : (x=y)\to (e(x)=e(y))
\end{equation*}
is an equivalence for any $x,y:A$. If $B$ is a proposition, then in particular the type $e(x)=e(y)$ is contractible for any $x,y:A$, so the claim follows from \cref{thm:contr_equiv}.
\end{proof}

\begin{thm}\label{thm:embedding}
  Consider a map $f:A\to B$. The following are equivalent:
  \begin{enumerate}
  \item The map $f$ is an embedding.\index{embedding}
  \item The fiber $\fib{f}{b}$ is a proposition for each $b:B$.
  \end{enumerate}
\end{thm}

\begin{proof}
  By the fundamental theorem of identity types, it follows that $f$ is an embedding if and only if
  \begin{equation*}
    \sm{x:A}f(x)=f(y)
  \end{equation*}
  is contractible for each $y:A$. In other words, $f$ is an embedding if and only if $\fib{f}{f(y)}$ is contractible for each $y:A$. Note that we obtain equivalences
  \begin{equation*}
    \fib{f}{f(y)}\simeq \fib{f}{b}
  \end{equation*}
  for any $b:B$ and $p:f(y)=b$, by transporting along $p$. Therefore it follows by \cref{lem:prop_equiv} that $\fib{f}{f(y)}$ is contractible for each $y:A$ if and only if $\fib{f}{b}$ is contractible for each $y:A$, and each $b:B$ such that $p:f(y)=b$. The latter condition holds if and only if we have
  \begin{equation*}
    \fib{f}{b}\to\iscontr(\fib{f}{b})
  \end{equation*}
  for any $b:B$, which is by \cref{lem:isprop_eq} equivalent to the condition that each $\fib{f}{b}$ is a proposition.
\end{proof}

\begin{cor}\label{cor:pr1-embedding}
  Consider a family $B$ of types over $A$. The following are equivalent:
  \begin{enumerate}
  \item The map $\proj 1 : (\sm{x:A}B(x))\to A$ is an embedding.
  \item The type $B(x)$ is a proposition for each $x:A$.
  \end{enumerate}
\end{cor}

\begin{proof}
  This corollary follows at once from \cref{ex:proj_fiber}, where we showed that
  \begin{equation*}
    \fib{\proj 1}{x}\simeq B(x).\qedhere
  \end{equation*}
\end{proof}
\index{subtype|)}

\subsection{Sets}

\index{set|(}
\begin{defn}
  A type $A$ is said to be a \define{set}\index{set|textbf} if its identity types are propositions, i.e., if it comes equipped with a term of type
  \index{is-set(A)@{$\isset(A)$}}\index{is a set}
\begin{equation*}
\isset(A)\defeq \prd{x,y:A}\isprop(\id{x}{y}).
\end{equation*}
\end{defn}

\begin{eg}\label{eg:is-set-nat}
  The type of natural numbers is a set.\index{is a set!natural numbers}\index{natural numbers!is a set}\index{N@{$\N$}!is a set} To see this, recall from \cref{thm:eq_nat} that we have an equivalence
  \begin{equation*}
    (m=n)\simeq \EqN(m,n)
  \end{equation*}
  for every $m,n:\N$. Therefore it suffices to show that each $\EqN(m,n)$ is a proposition. This follows easily by induction on both $m$ and $n$.
\end{eg}

\begin{prp}
  Consider a type $A$. The following are equivalent:
  \begin{enumerate}
  \item The type $A$ is a set.
  \item The type $A$ satisfies \define{axiom K}\index{axiom K|textbf}\index{axiom!axiom K|textbf}, i.e., if and only if it comes equipped with a term of type\index{is-set(A)@{$\isset(A)$}!is-set(A) iff axiom-K(A)@{$\isset(A)\leftrightarrow\axiomK(A)$}}\index{axiom KK(A)@{$\axiomK(A)$}|textbf}
    \begin{equation*}
      \axiomK(A)\defeq\prd{x:A}\prd{p:\id{x}{x}}\id{\refl{x}}{p}.
    \end{equation*}
  \end{enumerate}
\end{prp}

\begin{proof}
If $A$ is a set, then $\id{x}{x}$ is a proposition, so any two of its elements are equal. 
This implies axiom $K$. 

For the converse, if $A$ satisfies axiom $K$, then for any $p,q:\id{x}{y}$ we have $\id{\ct{p}{q^{-1}}}{\refl{x}}$, and hence $\id{p}{q}$. This shows that $\id{x}{y}$ is a proposition, and hence that $A$ is a set.
\end{proof}

\begin{thm}\label{lem:prop_to_id}
Let $A$ be a type, and let $R:A\to A\to\UU$ be a binary relation on $A$ satisfying
\begin{enumerate}
\item Each $R(x,y)$ is a proposition,
\item $R$ is reflexive, as witnessed by $\rho:\prd{x:A}R(x,x)$,
\item There is a map
  \begin{equation*}
    R(x,y)\to (x=y)
  \end{equation*}
  for each $x,y:A$.
\end{enumerate}
Then any family of maps
\begin{equation*}
\prd{x,y:A}(\id{x}{y})\to R(x,y)
\end{equation*}
is a family of equivalences. Consequently, the type $A$ is a set.
\end{thm}

\begin{proof}
Let $f:\prd{x,y:A}R(x,y)\to(\id{x}{y})$. 
Since $R$ is assumed to be reflexive, we also have a family of maps
\begin{equation*}
\pathind_x(\rho(x)):\prd{y:A}(\id{x}{y})\to R(x,y).
\end{equation*}
Since each $R(x,y)$ is assumed to be a proposition, it therefore follows that each $R(x,y)$ is a retract of $\id{x}{y}$. Therefore it follows that $\sm{y:A}R(x,y)$ is a retract of $\sm{y:A}x=y$, which is contractible. We conclude that $\sm{y:A}R(x,y)$ is contractible, and therefore that any family of maps
\begin{equation*}
  \prd{y:A}(x=y)\to R(x,y)
\end{equation*}
is a family of equivalences.

Now it also follows that $A$ is a set, since its identity types are equivalent to propositions, and therefore they are propositions by \cref{lem:prop_equiv}. 
\end{proof}

\begin{thm}[Hedberg]\label{thm:hedberg}
Any type with decidable equality is a set.\index{Hedberg's theorem}\index{set!Hedberg's theorem}
\end{thm}

\begin{proof}
  Let $A$ be a type, and let $d:\prd{x,y:A}(\id{x}{y})+ (x\neq y)$ be the witness that $A$ has decidable equality. Furthermore, let $\UU$ be a universe containing the type $A$. We will prove that $A$ is a set by applying \cref{lem:prop_to_id}.

  For every $x,y:A$, we first define a type family $R'(x,y):((\id{x}{y})+{(x\neq y)})\to\UU$ by
\begin{align*}
R'(x,y,\inl(p)) & \defeq \unit \\
R'(x,y,\inr(p)) & \defeq \emptyt.
\end{align*}
Note that $R'(x,y,q)$ is a proposition for each $x,y:A$ and $q:(\id{x}{y})+(x\neq y)$. 
Now we define $R(x,y)\defeq R'(x,y,d(x,y))$. Then $R$ is a reflexive binary relation on $A$, and furthermore each $R(x,y)$ is a proposition. In order to apply \cref{lem:prop_to_id}, it therefore it remains to show that $R$ implies identity. 

Since $R$ is defined as an instance of $R'$, it suffices to construct a function
\begin{equation*}
  f(q) : R'(q)\to (\id{x}{y}). 
\end{equation*}
for each $q:(\id{x}{y})+(x\neq y)$. Such a function is defined by
\begin{align*}
  f(\inl(p),r) & := p \\
  f(\inr(p),r) & := \exfalso(r).\qedhere
\end{align*}
\end{proof}
\index{set|)}

\subsection{General truncation levels}
\index{truncated type|(}
\index{truncation level|(}

Consider a type $A$ in a universe $\UU$. The conditions
\begin{align*}
  \iscontr(A) & \defeq \sm{a:A}\prd{x:A}a=x \\
  \isprop(A) & \defeq \prd{x,y:A}\iscontr(x=y) \\
  \isset(A) & \defeq \prd{x,y:A}\isprop(x=y)
\end{align*}
define the first few layers of the hierarchy of truncation levels. This hierarchy starts at the level of the contractible types, which we call level $-2$. The next level is the level of propositions, and at level $0$ we have the sets.

The indexing type of the truncation levels, which will be equivalent to the type $\Z_{\geq -2}$ of integers greater than $-2$, is an inductive type $\T$\index{T@{$\T$}|see {truncation level}} equipped with the constructors\index{-2 T@{$\negtwoT$}}\index{succT@{$\succT$}}
\begin{align*}
  \negtwoT & : \T \\
  \succT & : \T\to\T.
\end{align*}
The natural inclusion $i:\N\to \T$ is defined recursively by
\begin{align*}
  i(\zeroN) & \defeq \succT(\succT(\negtwoT)) \\
  i(\succN(n)) & \defeq \succT(i(n)).
\end{align*}
Of course, we will simply write $-2$ for $\negtwoT$ and $k+1$ for $\succT(k)$.

\begin{defn}
We define $\istrunc{} : \T\to\UU\to\UU$ recursively by\index{is-trunc k(A)@{$\istrunc{k}(A)$}}
\begin{align*}
\istrunc{-2}(A) & \defeq \iscontr(A) \\
\istrunc{k+1}(A) & \defeq \prd{x,y:A}\istrunc{k}(\id{x}{y}).\qedhere
\end{align*}
For any type $A$, we say that $A$ is \define{$k$-truncated}\index{k-truncated type@{$k$-truncated type}|see {truncated type}}\index{truncated type|textbf}, or a \define{$k$-type}\index{k-type@{$k$-type}|see {truncated type}}, if there is a term of type $\istrunc{k}(A)$. We also say that a type $A$ is a \define{proper $(k+1)$-type}\index{proper (k+1)-type@{proper $(k+1)$-type}|textbf} if $A$ is a $(k+1)$-type and not a $k$-type.

Given a universe $\UU$, we define the universe $\UU^{\leq k}$ of $k$-truncated types by\index{U leq k@{$\UU^{\leq k}$}|textbf}
\begin{equation*}
  \UU^{\leq k}\defeq\sm{X:\UU}\istrunc{k}(X).
\end{equation*}

Furthermore, we say that a map $f:A\to B$ is $k$-truncated if its fibers are $k$-truncated.\index{k-truncated map@{$k$-truncated map}|see {truncated map}}\index{truncated map|textbf}
\end{defn}

\begin{rmk}
  There is a subtlety in the definition of $\istrunc{}$ regarding universes. Note that the truncation levels are defined with respect to a universe $\UU$. To be completely precise, we should therefore write $\istrunc{k}^{\UU}(A)$ for the type $\istrunc{k}(A)$ defined with respect to the universe $\UU$. If $A$ is also contained in a second universe $\VV$, then it is legitimate to ask whether
  \begin{equation*}
    \istrunc{k}^{\UU}(A)\leftrightarrow\istrunc{k}^{\VV}(A).
  \end{equation*}
  A simple inductive argument shows that this is indeed the case, where the base case follows from the judgmental equalities
  \begin{align*}
    \istrunc{-2}^{\UU}(A) & \jdeq \sm{x:A}\prd{y:A}x=y \\
    \istrunc{-2}^{\VV}(A) & \jdeq \sm{x:A}\prd{y:A}x=y.
  \end{align*}
  We may therefore safely omit explicit reference to the universes when considering truncatedness of a type.
\end{rmk}

We show in the following theorem that the truncation levels are successively contained in one another.

\begin{prp}\label{thm:istrunc_next}
If $A$ is a $k$-type, then $A$ is also a $(k+1)$-type.\index{is-trunc k(A)@{$\istrunc{k}(A)$}!is-trunc k(A) to is-trunc k+1(A)@{$\istrunc{k}(A)\to\istrunc{k+1}(A)$}}
\end{prp}

\begin{proof}
We have seen in \cref{eg:prop_contr} that contractible types are propositions. This proves the base case.
For the inductive step, note that if any $k$-type is also a $(k+1)$-type, then any $(k+1)$-type is a $(k+2)$-type, since its identity types are $k$-types and therefore $(k+1)$-types.
\end{proof}

It is immediate from the proof of \cref{thm:istrunc_next} that the identity types of $k$-types are also $k$-types.

\begin{cor}
  If $A$ is a $k$-type, then its identity types are also $k$-types.\hfill $\square$
\end{cor}

\begin{prp}\label{thm:ktype_eqv}
If $e:\eqv{A}{B}$ is an equivalence, and $B$ is a $k$-type, then so is $A$.\index{truncated type!closed under equivalences}
\end{prp}

\begin{proof}
We have seen in \cref{ex:contr_equiv} that if $B$ is contractible and $e:\eqv{A}{B}$ is an equivalence, then $A$ is also contractible. This proves the base case.

For the inductive step, assume that the $k$-types are stable under equivalences, and consider $e:\eqv{A}{B}$ where $B$ is a $(k+1)$-type. In \cref{cor:emb_equiv} we have seen that
\begin{equation*}
\apfunc{e}:(\id{x}{y})\to(\id{e(x)}{e(y)})
\end{equation*}
is an equivalence for any $x,y$. Note that $\id{e(x)}{e(y)}$ is a $k$-type, so by the induction hypothesis it follows that $\id{x}{y}$ is a $k$-type. This proves that $A$ is a $(k+1)$-type.
\end{proof}

\begin{cor}\label{cor:emb_into_ktype}
If $f:A\to B$ is an embedding, and $B$ is a $(k+1)$-type, then so is $A$.\index{truncated type!closed under embeddings}
\end{cor}

\begin{proof}
By the assumption that $f$ is an embedding, the action on paths
\begin{equation*}
\apfunc{f}:(\id{x}{y})\to (\id{f(x)}{f(y)})
\end{equation*}
is an equivalence for every $x,y:A$. Since $B$ is assumed to be a $(k+1)$-type, it follows that $f(x)=f(y)$ is a $k$-type for every $x,y:A$. Therefore we conclude by \cref{thm:ktype_eqv} that $\id{x}{y}$ is a $k$-type for every $x,y:A$. In other words, $A$ is a $(k+1)$-type.
\end{proof}

We end this section with a theorem that characterizes $(k+1)$-truncated maps. Note that it generalizes \cref{thm:embedding}, which asserts that a map is an embedding if and only if its fibers are propositions.

\begin{thm}\label{thm:trunc_ap}
Let $f:A\to B$ be a map. The following are equivalent:
\begin{enumerate}
\item The map $f$ is $(k+1)$-truncated.\index{truncated map}
\item For each $x,y:A$, the map
\begin{equation*}
\apfunc{f} : (x=y)\to (f(x)=f(y))
\end{equation*}
is $k$-truncated. 
\end{enumerate}
\end{thm}

\begin{proof}
First we show that for any $s,t:\fib{f}{b}$ there is an equivalence
\begin{equation*}
\eqv{(s=t)}{\fib{\apfunc{f}}{\ct{\proj 2(s)}{\proj 2(t)^{-1}}}}
\end{equation*}
We do this by $\Sigma$-induction on $s$ and $t$, and then we calculate
\begin{align*}
(\pairr{x,p}=\pairr{y,q}) & \eqvsym \Eqfib_f((x,p),(y,q)) \\
  & \jdeq \sm{\alpha:x=y} p=\ct{\ap{f}{\alpha}}{q} \\
  & \eqvsym \sm{\alpha:x=y} \ct{\ap{f}{\alpha}}{q}=p \\
& \eqvsym \sm{\alpha:x=y} \ap{f}{\alpha}=\ct{p}{q^{-1}} \\
& \jdeq \fib{\apfunc{f}}{\ct{p}{q^{-1}}}.
\end{align*}
By these equivalences, it follows that if $\apfunc{f}$ is $k$-truncated, then for each $s,t:\fib{f}{b}$ the identity type $s=t$ is equivalent to a $k$-truncated type, and therefore we obtain by \cref{thm:ktype_eqv} that $f$ is $(k+1)$-truncated.

For the converse, note that we have equivalences
\begin{align*}
\fib{\apfunc{f}}{p} & \eqvsym ((x,p)=(y,\refl{f(y)})).
\end{align*}
It follows that if $f$ is $(k+1)$-truncated, then the identity type $(x,p)=(y,\refl{f(y)})$ in $\fib{f}{f(y)}$ is $k$-truncated for any $p:f(x)=f(y)$. We conclude by \cref{thm:ktype_eqv} that the fiber $\fib{\apfunc{f}}{p}$ is $k$-truncated. 
\end{proof}
\index{truncated type|)}
\index{truncation level|)}

\begin{exercises}
  \exitem \label{ex:eq_bool}Show that $\bool$ is a set\index{bool@{$\bool$}!is a set} by applying \cref{lem:prop_to_id} with the observational equality on $\bool$ defined in \cref{ex:obs_bool}.
  \exitem Recall that a \define{partially ordered set (poset)}\index{poset!is a set} is defined to be a type $A$ equipped with a relation
  \begin{equation*}
    \blank\leq\blank : A \to (A \to \prop_\UU)
  \end{equation*}
  that is reflexive, antisymmetric, and transitive. Show that the underlying type of any poset is a set.
  \exitem
  \begin{subexenum}
  \item \label{cor:is-emb-is-injective}
    Show that any injective map $f:A\to B$ into a set $B$ is an embedding, and conclude that $A$ is automatically a set in this case.\index{is an embedding!injective map into a set}\index{injective map!injective maps into sets are embeddings}
  \item Show that $n\mapsto m+n$ is an embedding, for each $m:\N$.\index{add N@{$\addN$}!add N(m) is an embedding@{$\addN(m)$ is an embedding}}\index{is an embedding!add N(m)@{$\addN(m)$}} Moreover, conclude that there is an equivalence\index{N@{$\N$}!leq@{$\leq$}}
    \begin{equation*}
      (m\leq n)\simeq \sm{k:\N}m+k=n.
    \end{equation*}
  \item Show that $n\mapsto mn$ is an embedding\index{mul N@{$\mulN$}!mul N(m) is an embedding if m>0@{$\mulN(m)$ is an embedding if $m>0$}}\index{is an embedding!mul N(m) for m > 0@{$\mulN(m)$ for $m>0$}}, for each nonzero number $m:\N$. Conclude that the divisibility relation\index{d {"|" n}@{$d\mid n$}!is a proposition if d>0@{is a proposition if $d>0$}}\index{is a proposition!d {"|" n} for d>0@{$d\mid n$ for $d>0$}}\index{divisibility on N@{divisibility on $\N$}!is a proposition}
    \begin{equation*}
      d\mid n
    \end{equation*}
    is a proposition for each $d,n:\N$ such that $d>0$. 
  \end{subexenum}
  \exitem \label{ex:set_coprod}
  \begin{subexenum}
  \item Show that for any two contractible types $A$ and $B$, the coproduct $A+B$ is not contractible.
  \item Show that for any two propositions $P$ and $Q$, we have a logical equivalence
    \begin{equation*}
      \iscontr(P+Q)\leftrightarrow P\oplus Q,
    \end{equation*}
    where the \define{exclusive disjunction}\index{exclusive disjunction|textbf} $P\oplus Q$\index{P oplus Q@{$P \oplus Q$}|see {exclusive disjunction}} is defined by
    \begin{equation*}
      P\oplus Q:= (P\times\neg Q)+(Q\times\neg P).
    \end{equation*}
  \item \label{ex:is-prop-coproduct}Show that for any two propositions $P$ and $Q$, the coproduct $P+Q$ is a proposition if and only if $P\to \neg Q$.
  \item Show that for any two $(k+2)$-types $A$ and $B$, the coproduct $A+B$ is again a $(k+2)$-type\index{coproduct!is a truncated type}\index{is a truncated type!coproduct}\index{is a set!coproduct}\index{coproduct!is a set}. Conclude that $\Z$ is a set.\index{Z@{$\Z$}!is a set}
  \end{subexenum}
  \exitem \label{ex:diagonal}Let $A$ be a type, and let the \define{diagonal}\index{diagonal of a type|textbf}\index{d  A@{$\delta_A$}|textbf}\index{d  A@{$\delta_A$}|see {diagonal, of a type}} of $A$ be the map $\delta_A:A\to A\times A$ given by $\lam{x}(x,x)$. 
  \begin{subexenum}
  \item Show that\index{is a proposition!d A is an equivalence@{$\delta_A$ is an equivalence}}
    \begin{equation*}
      {\isequiv(\delta_A)}\leftrightarrow{\isprop(A)}.
    \end{equation*}
  \item Construct an equivalence $\eqv{\fib{\delta_A}{x,y}}{(x=y)}$ for any $x,y:A$.
  \item Show that $A$ is $(k+1)$-truncated if and only if $\delta_A:A\to A\times A$ is $k$-truncated.
  \end{subexenum}
  \exitem \label{ex:istrunc_sigma}
  \begin{subexenum}
  \item Consider a type family $B$ over a $k$-truncated type $A$. Show that the following are equivalent:
    \begin{enumerate}
    \item The type $B(x)$ is $k$-truncated for each $x:A$.
    \item The type $\sm{x:A}B(x)$ is $k$-truncated.\index{is a truncated type!S-type@{$\Sigma$-type}}\index{dependent pair type!is truncated}
    \end{enumerate}
    Hint: for the base case, use \cref{ex:contr_in_sigma,ex:contr_equiv}.
  \item Consider a map $f:A\to B$ into a $k$-type $B$. Show that the following are equivalent:
    \begin{enumerate}
    \item The type $A$ is $k$-truncated.
    \item The map $f$ is $k$-truncated.
    \end{enumerate}
  \end{subexenum}
  \exitem Consider two types $A$ and $B$. Show that the following are equivalent:
    \begin{enumerate}
    \item There are functions
      \begin{align*}
        f & : B \to \istrunc{k+1}(A) \\
        g & : A \to \istrunc{k+1}(B).
      \end{align*}
    \item The type $A\times B$ is $(k+1)$-truncated.
    \end{enumerate}
    Conclude with \cref{ex:is-contr-prod} that, if both $A$ and $B$ come equipped with an element, then both $A$ and $B$ are $k$-truncated if and only if the product $A\times B$ is $k$-truncated.
  \exitem
  \begin{subexenum}
  \item \label{ex:retr_id} Consider a section-retraction pair
    \begin{equation*}
      \begin{tikzcd}
        A \arrow[r,"i"] & B \arrow[r,"r"] & A,
      \end{tikzcd}
    \end{equation*}
    with $H:r\circ i\htpy \idfunc$. Show that $\id{x}{y}$ is a retract of $\id{i(x)}{i(y)}$.\index{retract!identity type}\index{identity type!of retract is retract}
  \item Use \cref{ex:contr_retr} to show that if $A$ is a retract of a $k$-type $B$, then $A$ is also a $k$-type.\index{truncated type!closed under retracts}
  \end{subexenum}
  \exitem Consider an arbitrary type $A$. Recall that concatenation of lists was defined in \cref{ex:lists}. Show that the map\index{list A@{$\lst(A)$}}
  \begin{equation*}
    f:\lst(A)\times\lst(A)\to\lst(A).
  \end{equation*}
  given by $f(x,y)\defeq\concatlist(x,y)$ is $0$-truncated.\index{truncated map!concatenation of lists}\index{concat-list@{$\concatlist$}!is a $0$-truncated map}
  \exitem \label{ex:is-trunc-const}Show that a type $A$ is a $(k+1)$-type if and only if the map $\const_x:\unit\to A$ is $k$-truncated for every $x:A$.
  \exitem \label{ex:is-trunc-comp}Consider a commuting triangle
  \begin{equation*}
    \begin{tikzcd}[column sep=tiny]
      A \arrow[rr,"h"] \arrow[dr,swap,"f"] & & B \arrow[dl,"g"] \\
      & X
    \end{tikzcd}
  \end{equation*}
  with $H: f \htpy g \circ h$, and suppose that $g$ is $k$-truncated. Show that $f$ is $k$-truncated if and only if $h$ is $k$-truncated.
  \exitem Let $f:\prd{x:A}B(x)\to C(x)$ be a family of maps. Show that the following are equivalent:
  \begin{enumerate}
  \item For each $x:A$ the map $f(x)$ is $k$-truncated.
  \item The induced map\index{tot(f)@{$\tot{f}$}!is a truncated map}\index{is a truncated map!tot(f)@{$\tot{f}$}}
    \begin{equation*}
      \tot{f}:\Big(\sm{x:A}B(x)\Big)\to\Big(\sm{x:A}C(x)\Big)
    \end{equation*}
    is $k$-truncated.
  \end{enumerate}
  \exitem \label{ex:is-trunc-fiber-inclusion}Consider a type $A$. Show that the following are equivalent:
  \begin{enumerate}
  \item The type $A$ is $(k+1)$-truncated.
  \item For any type family $B$ over $A$ and any $a:A$, the \define{fiber inclusion}\index{fiber inclusion|textbf}
    \begin{equation*}
      i_a: B(a)\to\sm{x:A}B(x)
    \end{equation*}
    given by $y\mapsto(a,y)$ is a $k$-truncated map.\index{fiber inclusion!is a truncated map}\index{is a truncated map!fiber inclusion}
  \end{enumerate}
  In particular, if $A$ is a set then any fiber inclusion $i_a:B(a)\to\sm{x:A}B(x)$ is an embedding.\index{fiber inclusion!is an embedding}\index{is an embedding!fiber inclusion}
  \exitem \label{ex:isolated-point}Consider a type $A$ equipped with an element $a:A$. We say that $a$ is an \define{isolated element}\index{isolated element|textbf} of $A$ if it comes equipped with an element of type\index{is-isolated(a)@{$\isisolated(a)$}|textbf}
  \begin{equation*}
    \isisolated(a)\defeq\prd{x:A}(a=x)+(a\neq x).
  \end{equation*}
  \begin{subexenum}
  \item Show that $a$ is isolated if and only if the map $\const_a:\unit\to A$ has decidable fibers.
  \item Show that if $a$ is isolated, then $a=x$ is a proposition, for every $x:A$. Conclude that if $a$ is isolated, then the map $\const_a:\unit\to A$ is an embedding.
  \end{subexenum}
\end{exercises}
\index{truncated type|)}
\index{truncation level|)}


\section{Function extensionality}
\label{chap:funext}
\index{function extensionality|(}
\index{axiom!function extensionality|(}

The function extensionality axiom asserts that for any two dependent functions $f,g:\prd{x:A}B(x)$, the type of identifications $f=g$ is equivalent to the type of homotopies $f\htpy g$ from $f$ to $g$. In other words, two (dependent) functions can only be distinguished by their values. The function extensionality axiom therefore provides a characterization of the identity type of (dependent) function types. By the fundamental theorem of identity types it follows immediately that the function extensionality axiom has at least three equivalent forms. There is, however, a fourth useful equivalent form of the function extensionality axiom: the \emph{weak} function extensionality axiom. This axiom asserts that any dependent product of contractible types is again contractible. A simple consequence of the weak function extensionality axiom is that any dependent product of a family of $k$-types is again a $k$-type.

The function extensionality axiom is used to derive many important properties in type theory. One class of such properties are (dependent) universal properties. Universal properties give a characterization of the type of functions into, or out of a type. For example, the universal property of the coproduct $A+B$ characterizes the type of maps $(A+B)\to X$ as the type of pairs of maps $(f,g)$ consisting of $f:A\to X$ and $g:B\to X$, i.e., the universal property of the coproduct $A+B$ is an equivalence
\begin{equation*}
  ((A+B)\to X)\simeq (A\to X)\times (B\to X).
\end{equation*}
Note that there are function types on both sides of this equivalence. Therefore we will need function extensionality in order to construct the homotopies witnessing that the inverse map is both a left and a right inverse. In fact, we leave this particular universal property as \cref{ex:up-coproduct}. The universal properties that we do show in the main text, are the universal properties of $\Sigma$-types and of the identity type. 

We end this section with two further applications of the function extensionality axiom. In the first, \cref{ex:equiv_precomp}, we show that precomposition by an equivalence is again an equivalence. More precisely we show that $f:A\to B$ is an equivalence if and only if for every type family $P$ over $B$, the precomposition map
\begin{equation*}
  \blank\circ f :\Big(\prd{y:B}P(y)\Big)\to \Big(\prd{x:A}P(f(x))\Big)
\end{equation*}
is an equivalence. To prove this fact we will make use of coherently invertible maps, which were introduced in \cref{sec:is-contr-map-is-equiv}. In the second application, \cref{thm:strong-ind-N}, we prove the strong induction principle of the natural numbers. Function extensionality is needed in order to derive the computation rule for the strong induction principle.

Many important consequences of the function extensionality axiom are left as exercises. For example, in \cref{ex:isprop_istrunc} you are asked to show that both $\iscontr(A)$ and $\istrunc{k}(A)$ are propositions, and in \cref{ex:isprop_isequiv} you are asked to show that $
\isequiv(f)$ is a proposition. The universal properties of $\emptyt$, $\unit$, and $A+B$ are left as \cref{ex:up-emptyt,ex:up-unit,ex:up-coproduct}. A few more advanced properties, such as the fact that post-composition
\begin{equation*}
  g\circ\blank : (A\to X)\to (A\to Y)
\end{equation*}
by a $k$-truncated map $g:X\to Y$ is itself a $k$-truncated map, appear in the later exercises. We encourage you to read through all of them, and get at least a basic idea of why they are true.

\subsection{Equivalent forms of function extensionality}

The function extensionality principle characterizes the identity type of an arbitrary dependent function type. It asserts that the type $f=g$ of identifications between two dependent functions is equivalent to the type of homotopies $f\htpy g$. By \cref{thm:id_fundamental}\index{fundamental theorem of identity types} there are three equivalent ways of doing this.

\begin{prp}\label{prp:funext}
  Consider a dependent function $f:\prd{x:A}B(x)$. The following are equivalent:\index{function extensionality|textbf}\index{characterization of identity type!of P-types@{of $\Pi$-types}}\index{dependent function type!characterization of identity type}
  \begin{enumerate}
  \item The \define{function extensionality principle} holds at $f$: for each $g:\prd{x:A}B(x)$, the family of maps
    \begin{equation*}
      \htpyeq:(f=g)\to (f\htpy g)
    \end{equation*}
    defined by $\htpyeq(\refl{f}):=\reflhtpy_{f}$ is a family of equivalences.
  \item The total space
    \begin{equation*}
      \sm{g:\prd{x:A}B(x)}f\htpy g
    \end{equation*}
    is contractible.
  \item
    The principle of \define{homotopy induction}\index{homotopy induction|textbf}\index{induction principle!for homotopies|textbf}:
    for any family of types $P(g,H)$ indexed by $g:\prd{x:A}B(x)$ and $H:f\htpy g$, the evaluation function
    \begin{equation*}
      \Big(\prd{g:\prd{x:A}B(x)}\prd{H:f\htpy g}P(g,H)\Big)\to P(f,\reflhtpy_f),
    \end{equation*}
    given by $s\mapsto s(f,\reflhtpy_f)$, has a section.
  \end{enumerate}
\end{prp}

\begin{proof}
  This theorem follows directly from \cref{thm:id_fundamental}.
\end{proof}

There is, however, yet a fourth condition equivalent to the function extensionality principle: the \emph{weak} function extensionality principle. The weak function extensionality principle asserts that any dependent product of contractible types is again contractible.

The following theorem is stated with respect to an arbitrary universe $\UU$, because we will use it in \cref{thm:funext-univalence} to show that the univalence axiom implies function extensionality.

\begin{thm}\label{thm:funext_wkfunext}
  Consider a universe $\UU$. The following are equivalent:\index{function extensionality}
  \begin{enumerate}
  \item The function extensionality principle holds in $\UU$: For every type family $B$ over $A$ in $\UU$ and any $f,g:\prd{x:A}B(x)$, the map
    \begin{equation*}
      \htpyeq : (f=g)\to (f\htpy g)
    \end{equation*}
    is an equivalence.
  \item The \define{weak function extensionality principle}\index{weak function extensionality|textbf}\index{function extensionality!weak function extensionality|textbf} holds in $\UU$: For every type family $B$ over $A$ in $\UU$ one has\index{contractible type!weak function extensionality|textbf}\index{is contractible!dependent function type}
    \begin{equation*}
      \Big(\prd{x:A}\iscontr(B(x))\Big)\to\iscontr\Big(\prd{x:A}B(x)\Big).
    \end{equation*}
  \end{enumerate}
\end{thm}

\begin{proof}
  First, we show that function extensionality implies weak function extensionality, suppose that each $B(a)$ is contractible with center of contraction $c(a)$ and contraction $C_a:\prd{y:B(a)}c(a)=y$. Then we take $c\defeq \lam{a}c(a)$ to be the center of contraction of $\prd{x:A}B(x)$. To construct the contraction we have to define a term of type
  \begin{equation*}
    \prd{f:\prd{x:A}B(x)}c=f.
  \end{equation*}
  Let $f:\prd{x:A}B(x)$. By function extensionality we have a map ${(c\htpy f)}\to {(c=f)}$, so it suffices to construct a term of type $c\htpy f$. Here we take $\lam{a}C_a(f(a))$. This completes the proof that function extensionality implies weak function extensionality.
  
  It remains to show that weak function extensionality implies function extensionality. By \cref{prp:funext} it suffices to show that the type
  \begin{equation*}
    \sm{g:\prd{x:A}B(x)}f\htpy g
  \end{equation*}
  is contractible for any $f:\prd{x:A}B(x)$. In order to do this, we first note that we have a section-retraction pair
  \begin{align*}
    \Big(\sm{g:\prd{x:A}B(x)}f\htpy g\Big)
    & \stackrel{i}{\longrightarrow} \Big(\prd{x:A}\sm{b:B(x)}f(x)=b\Big) \\
    & \stackrel{r}{\longrightarrow} \Big(\sm{g:\prd{x:A}B(x)}f\htpy g\Big)
  \end{align*}
  Here we have the functions
  \begin{align*}
    i & \defeq \lam{(g,H)}\lam{x}(g(x),H(x)) \\
    r & \defeq \lam{p}\pairr{\lam{x}\proj 1(p(x)),\lam{x}\proj 2(p(x))}.
  \end{align*}
  Their composite is homotopic to the identity function by the computation rule for $\Sigma$-types and the $\eta$-rule for $\Pi$-types:
  \begin{align*}
    r(i(g,H)) & \jdeq r(\lam{x}\pairr{g(x),H(x)}) \\
              & \jdeq \pairr{\lam{x}g(x),\lam{x}H(x)} \\
              & \jdeq \pairr{g,H}.
  \end{align*}
  Now we observe that the type $\prd{x:A}\sm{b:B(x)}f(x)=b$ is a product of contractible types, so it is contractible by our assumption of the weak function extensionality principle. The claim now follows, because retracts of contractible types are contractible by \cref{ex:contr_retr}.
\end{proof}

We will henceforth assume the function extensionality principle as an axiom.

\begin{axiom}[Function Extensionality]\label{axiom:funext}
  \index{function extensionality|textbf}\index{axiom!function extensionality|textbf}\index{identity type!of a Pi-type@{of a $\Pi$-type}}\index{extensionality principle!for functions|textbf}%
  For any type family $B$ over $A$, and any two dependent functions $f,g:\prd{x:A}B(x)$, the map\index{htpy-eq@{$\htpyeq$}|textbf}\index{htpy-eq@{$\htpyeq$}!is an equivalence}\index{is an equivalence!htpy-eq@{$\htpyeq$}}
  \begin{equation*}
    \htpyeq:(f=g)\to (f\htpy g)
  \end{equation*}
  is an equivalence. We will write $\eqhtpy$\index{eq-htpy@{$\eqhtpy$}|textbf} for its inverse.
\end{axiom}

\begin{rmk}
  The function extensionality axiom is added to type theory by adding the rule
  \begin{prooftree}
    \AxiomC{$\Gamma,x:A\vdash B(x)~\type$}
    \AxiomC{$\Gamma\vdash f : \prd{x:A}B(x)$}
    \AxiomC{$\Gamma\vdash g : \prd{x:A}B(x)$}
    \TrinaryInfC{$\Gamma\vdash\funext:\isequiv(\htpyeq_{f,g})$}
  \end{prooftree}
\end{rmk}

In the following theorem we extend the weak function extensionality principle to general truncation levels.

\begin{thm}\label{thm:trunc_pi}\index{k-type@{$k$-type}}
For any type family $B$ over $A$ one has\index{truncated type!dependent function type}
\begin{equation*}
\Big(\prd{x:A}\istrunc{k}(B(x))\Big)\to \istrunc{k}\Big(\prd{x:A}B(x)\Big).
\end{equation*}
\end{thm}

\begin{proof}
The theorem is proven by induction on $k\geq -2$. The base case is just the weak function extensionality principle\index{weak function extensionality}, which was shown to follow from function extensionality in \cref{thm:funext_wkfunext}.

For the inductive step, assume that the $k$-truncated types are closed under $\Pi$-types, and consider a family $B$ of $(k+1)$-truncated types. To show that the type $\prd{x:A}B(x)$ is $(k+1)$-truncated, we have to show that the type $f=g$ is $k$-truncated for every $f,g:\prd{x:A}$. By function extensionality, the type $f=g$ is equivalent to $f\htpy g$ for any two dependent functions $f,g:\prd{x:A}B(x)$. Now observe that $f\htpy g$ is a dependent product of $k$-truncated types, and therefore it is $k$-truncated by the inductive hypothesis. Since the $k$-truncated types are closed under equivalences by \cref{thm:ktype_eqv}, it follows that the type $f=g$ is $k$-truncated.
\end{proof}

\begin{cor}\label{cor:funtype_trunc}\index{truncated type!function type}
Suppose $B$ is a $k$-type. Then $A\to B$ is also a $k$-type, for any type $A$.
\end{cor}

\begin{rmk}
  It follows that $\neg A$ is a proposition for each type $A$. Note that it requires function extensionality even just to prove that $\neg P$ is a proposition for any proposition $P$.
\end{rmk}

\subsection{Identity systems on \texorpdfstring{$\Pi$}{Π}-types}

Recall from \cref{sec:structure-identity-principle} that the \emph{structure identity principle} is a way to obtain an identity system on a $\Sigma$-type. Identity systems were defined in \cref{defn:identity-system}. In this section we will describe how to obtain identity systems on a $\Pi$-type. We will first show that $\Pi$-types distribute over $\Sigma$-types\index{distributivity!of P over S@{of $\Pi$ over $\Sigma$}}. This theorem is sometimes called the \emph{type theoretic principle of choice} because it can be seen as the Curry-Howard interpretation of the axiom of choice\index{type theoretic choice|see {distributivity, of $\Pi$ over $\Sigma$}}.

\begin{thm}\label{thm:choice}
Consider a family of types $C(x,y)$ indexed by $x:A$ and $y:B(x)$. Then the map\index{choice@{$\choice$}|textbf}\index{is an equivalence!choice@{$\choice$}}
\begin{equation*}
  \choice:\Big(\prd{x:A}\sm{y:B(x)}C(x,y)\Big)\to \Big(\sm{f:\prd{x:A}B(x)}\prd{x:A}C(x,f(x))\Big)
\end{equation*}
given by
\begin{equation*}\label{eq:choice}
  \choice(h):=(\lam{x}\proj 1(h(x)),\lam{x}\proj 2(h(x))).
\end{equation*}
is an equivalence.
\end{thm}

\begin{proof}
  We define the map\index{choice -1@{$\choice^{-1}$}|textbf}
  \begin{equation*}
    \choice^{-1}:\Big(\sm{f:\prd{x:A}B(x)}\prd{x:A}C(x,f(x))\Big)\to \prd{x:A}\sm{y:B(x)}C(x,y)
  \end{equation*}
  by $\choice^{-1}(f,g):=\lam{x}(f(x),g(x))$. Then we have to construct homotopies
  \begin{equation*}
    \choice\circ\choice^{-1}\htpy\idfunc,\qquad\text{and}\qquad
    \choice^{-1}\circ\choice\htpy\idfunc.
  \end{equation*}
  For the first homotopy it suffices to construct an identification
  \begin{equation*}
    \choice(\choice^{-1}(f,g))=(f,g)
  \end{equation*}
  for any $f:\prd{x:A}B(x)$ and any $g:\prd{x:A}C(x,f(x))$. We compute the left-hand side as follows:
  \begin{align*}
    \choice(\choice^{-1}(f,g))
    & \jdeq \choice(\lam{x}(f(x),g(x))) \\
    & \jdeq (\lam{x}f(x),\lam{x}g(x)).
  \end{align*}
  By the $\eta$-rule for $\Pi$-types we have the judgmental equalities $f\jdeq \lam{x}f(x)$ and $g\jdeq\lam{x}g(x)$. Therefore we have the identification
  \begin{equation*}
    \refl{(f,g)}:\choice(\choice^{-1}(f,g))=(f,g).
  \end{equation*}
  This completes the construction of the first homotopy.

  For the second homotopy we have to construct an identification
  \begin{equation*}
    \choice^{-1}(\choice(h))=h
  \end{equation*}
  for any $h:\prd{x:A}\sm{y:B(x)}C(x,y)$. We compute the left-hand side as follows:
  \begin{align*}
    \choice^{-1}(\choice(h))
    & \jdeq \choice^{-1}(\lam{x}\proj 1(h(x)),(\lam{x}\proj 2(h(x)))) \\
    & \jdeq \lam{x}(\proj 1(h(x)),\proj 2(h(x)))
  \end{align*}
  However, it is \emph{not} the case that $(\proj 1(h(x)),\proj 2(h(x)))\jdeq h(x)$ for any $h:\prd{x:A}\sm{y:B(x)}C(x,y)$. Nevertheless, we have the identification
  \begin{equation*}
    \eqpair(\refl{},\refl{}):(\proj 1(h(x)),\proj 2(h(x)))= h(x).
  \end{equation*}
  Therefore we obtain the required homotopy by function extensionality:
  \begin{equation*}
    \lam{h}\eqhtpy(\lam{x}\eqpair(\refl{\proj 1(h(x))},\refl{\proj 2(h(x))})):\choice^{-1}\circ\choice\htpy\idfunc.\qedhere
  \end{equation*}
\end{proof}

The fact that $\Pi$-types distribute over $\Sigma$-types has many useful consequences. The most straightforward consequence is the following.

\begin{cor}
  For any two types $A$ and $B$, and any type family $C$ over $B$, we have an equivalence
\begin{equation*}
  \Big(A\to\sm{y:B}C(y)\Big)\simeq\Big(\sm{f:A\to B}\prd{x:A}C(f(x))\Big).
\end{equation*}
\end{cor}

Another direct consequence of the distributivity of $\Pi$-types over $\Sigma$-types is the fact that
\begin{equation*}
  \prd{b:B}\fib{f}{b}\simeq\sm{g:B\to A}f\circ g\htpy \idfunc.
\end{equation*}
In the following corollary we use the distributivity of $\Pi$-types over $\Sigma$-types to show that dependent functions are sections of projection maps.

\begin{cor}\label{ex:pi_sec}
  Consider a type family $B$ over $A$, and consider the projection map
  \begin{equation*}
    \proj 1:\big(\sm{x:A}B(x)\big) \to A.  
  \end{equation*}
  Then we have an equivalence
  \begin{equation*}
    \sections(\proj 1)\simeq\prd{x:A}B(x).
  \end{equation*}
\end{cor}

\begin{proof}
  \cref{thm:choice} gives the first equivalence in the following calculation:
  \begin{align*}
    \sm{h:A\to\sm{x:A}B(x)} \proj 1\circ h\htpy \idfunc
    & \simeq \sm{(f,g):\sm{f:A\to A}\prd{x:A}B(f(x))} f\htpy \idfunc \\
    & \simeq \sm{(f,H):\sm{f:A\to A}f\htpy \idfunc}\prd{x:A}B(f(x)) \\
    & \simeq \prd{x:A}B(x)
  \end{align*}
  In the second equivalence we used \cref{ex:sigma_swap} to swap the family $f\mapsto \prd{x:A}B(f(x))$ with the family $f\mapsto f\htpy\idfunc$, and in the third equivalence we used the fact that
  \begin{equation*}
    \sm{f:A\to A}f\htpy\idfunc
  \end{equation*}
  is contractible, with center of contraction $(\idfunc,\reflhtpy)$. One way to see that it is contractible is by \cref{ex:is-equiv-inv-htpy}. A direct way to see this, is by another application of \cref{thm:choice}. This gives an equivalence
  \begin{equation*}
    \left(\sm{f:A\to A}f\htpy\idfunc\right)\simeq \left(\prd{x:A}\sm{y:A}y=x\right),
  \end{equation*}
  and the right-hand side is a product of contractible types.
\end{proof}

In the final application of distributivity of $\Pi$-types over $\Sigma$-types we obtain a general way of constructing identity systems of $\Pi$-types.

\begin{thm}\label{cor:Eq-Pi}
  Consider a family $B$ of types over $A$, and for each $b:B(a)$ consider an identity system $E(b)$ at $b$. Furthermore, consider a dependent function $f:\prd{x:A}B(x)$. Then the family of types
  \begin{equation*}
    \prd{x:A}E(f(x),g(x))
  \end{equation*}
  indexed by $g:\prd{x:A}B(x)$ is an identity system at $f$.\index{identity system!of a dependent function type}\index{dependent function type!identity system}
\end{thm}

\begin{proof}
  By \cref{thm:id_fundamental} it suffices to show that the type
  \begin{equation*}
    \sm{g:\prd{x:A}B(x)}\prd{x:A}E(f(x),g(x))
  \end{equation*}
  is contractible. By \cref{thm:choice} it follows that this type is equivalent to the type
  \begin{equation*}
    \prd{x:A}\sm{y:B(x)}E(f(x),y).
  \end{equation*}
  This is a product of contractible types because each $E(f(x))$ is an identity system at $f(x):B(x)$. This product is therefore contractible by the weak function extensionality principle.
\end{proof}

\subsection{Universal properties}
The function extensionality principle allows us to prove \emph{universal properties}. Universal properties are characterizations of all maps out of or into a given type, so they are very important. Among other applications, universal properties characterize a type up to equivalence. We prove here the universal properties of dependent pair types and of identity types. In the exercises, you are asked to prove the universal properties of $\unit$, $\emptyt$, and coproducts.

\subsubsection*{The universal property of $\Sigma$-types}
\index{universal property!of S-types@{of $\Sigma$-types}|(}
\index{dependent universal property!of S-types@{of $\Sigma$-types}|(}
\index{dependent pair type!universal property|(}
\index{dependent pair type!dependent universal property|(}

The \define{universal property of $\Sigma$-types} characterizes maps \emph{out of} a dependent pair type $\sm{x:A}B(x)$. It asserts that the map\index{ev-pair@{$\evpair$}|textbf}
\begin{equation*}
\evpair:\Big(\Big(\sm{x:A}B(x)\Big)\to X\Big)\to \Big(\prd{x:A}(B(x)\to X)\Big),
\end{equation*}
given by $f\mapsto\lam{x}\lam{y}f(x,y)$, is an equivalence for any type $X$. In fact, we will prove a slight generalization of this universal property. We will prove the \define{dependent universal property} of $\Sigma$-types, which characterizes \emph{dependent} functions out of $\sm{x:A}B(x)$.

\begin{thm}\label{thm:up-sigma}
  \index{dependent universal property!of S-types@{of $\Sigma$-types}|textbf}
  \index{dependent pair type!dependent universal property|textbf}
Let $B$ be a type family over $A$, and let $C$ be a type family over $\sm{x:A}B(x)$. Then the map\index{ev-pair@{$\evpair$}!is an equivalence}
\begin{equation*}
\evpair:\Big(\prd{z:\sm{x:A}B(x)}C(z)\Big)\to \Big(\prd{x:A}\prd{y:B(x)}C(x,y)\Big),
\end{equation*}
given by $f\mapsto\lam{x}\lam{y}f(x,y)$, is an equivalence.\index{is an equivalence!ev-pair@{$\evpair$}}
\end{thm}

\begin{proof}
The map in the converse direction is obtained by the induction principle of $\Sigma$-types. It is simply the map
\begin{equation*}
\indSigma : \Big(\prd{x:A}\prd{y:B(x)}C(x,y)\Big)\to \Big(\prd{z:\sm{x:A}B(x)}C(z)\Big).
\end{equation*}
By the computation rule for $\Sigma$-types we have the homotopy
\begin{equation*}
\reflhtpy:\evpair\circ\indSigma\htpy\idfunc.
\end{equation*}
This shows that $\indSigma$ is a section of $\evpair$.

To show that $\indSigma\circ\evpair\htpy\idfunc$ we will apply the function extensionality principle. Therefore it suffices to show that $\indSigma(\lam{x}\lam{y}f(x,y))=f$. We apply function extensionality again, so it suffices to show that
\begin{equation*}
\prd{t:\sm{x:A}B(x)}\indSigma\big(\lam{x}\lam{y}f(x,y)\big)(t)=f(t).
\end{equation*}
We obtain this homotopy by another application of $\Sigma$-induction. 
\end{proof}

\begin{cor}\label{cor:times_up_out}
  \index{universal property!of cartesian products|textbf}
  \index{cartesian product type!universal property|textbf}
Let $A$, $B$, and $X$ be types. Then the map\index{ev-pair@{$\evpair$}!is an equivalence}\index{is an equivalence!ev-pair@{$\evpair$}}
\begin{equation*}
\evpair: (A\times B \to X)\to (A\to (B\to X))
\end{equation*}
given by $f\mapsto\lam{a}\lam{b}f(a,b)$ is an equivalence.
\end{cor}
\index{universal property!of S-types@{of $\Sigma$-types}|)}
\index{dependent universal property!of S-types@{of $\Sigma$-types}|)}
\index{dependent pair type!universal property|)}
\index{dependent pair type!dependent universal property|)}

\subsubsection*{The universal property of identity types}
\index{identity type!universal property|(}
\index{identity type!dependent universal property|(}
\index{universal property!of identity types|(}
\index{dependent universal property!of identity types|(}
The universal property of identity types is the fact that families of maps out of the identity type are uniquely determined by their action on the reflexivity identification. More precisely, the map
\begin{equation*}
  \evrefl:\Big(\prd{x:A}(a=x)\to B(x)\Big)\to B(a)
\end{equation*}
given by $\lam{f} f(a,\refl{a})$ is an equivalence, for every type family $B$ over $A$. Since this result is similar to the Yoneda lemma of category theory, the universal property of identity types is sometimes referred to as the \emph{type theoretic Yoneda lemma}. We will prove the \emph{dependent} universal property of identity types, a slight generalization of the universal property.

\begin{thm}\label{thm:yoneda}
  \index{dependent universal property!of identity types|textbf}
  \index{identity type!dependent universal property|textbf}
Consider a type $A$ equipped with $a:A$, and consider a family of types $B(x,p)$ indexed by $x:A$ and $p:a=x$. Then the map\index{ev-refl@{$\evrefl$}|textbf}
\begin{equation*}
\evrefl:\Big(\prd{x:A}\prd{p:a=x}B(x,p)\Big)\to B(a,\refl{a}),
\end{equation*}
given by $\lam{f} f(a,\refl{a})$, is an equivalence.\index{is an equivalence!ev-refl@{$\evrefl$}}\index{ev-refl@{$\evrefl$}!is an equivalence}
\end{thm}

\begin{proof}
  The inverse is the function
  \begin{equation*}
    \pathind_a : B(a,\refl{a})\to \prd{x:A}\prd{p:a=x}B(x,p).
  \end{equation*}
  It is immediate from the computation rule of the path induction principle that $\evrefl\circ\pathind_a\htpy \idfunc$.

To see that $\pathind_a\circ \evrefl\htpy\idfunc$, let $f:\prd{x:A}(a=x)\to B(x,p)$. To show that $\pathind_a(f(a,\refl{a}))=f$ we apply function extensionality twice. Therefore it suffices to show that
\begin{equation*}
\prd{x:A}\prd{p:a=x} \pathind_a(f(a,\refl{a}),x,p)=f(x,p).
\end{equation*}
This follows by path induction on $p$, since $\pathind_a(f(a,\refl{a}),a,\refl{a})\jdeq f(a,\refl{a})$ by the computation rule of path induction.
\end{proof}
\index{identity type!universal property|)}
\index{identity type!dependent universal property|)}
\index{universal property!of identity types|)}
\index{dependent universal property!of identity types|)}

\subsection{Composing with equivalences}

We show in the following theorem that a map $f:A\to B$ is an equivalence if and only if precomposing by $f$ is an equivalence.

\begin{thm}\label{ex:equiv_precomp}
  \index{equivalence!precomposition}
  For any map $f:A\to B$, the following are equivalent:
  \begin{enumerate}
  \item $f$ is an equivalence.
  \item For any type family $P$ over $B$ the map
    \begin{equation*}
      \Big(\prd{y:B}P(y)\Big)\to\Big(\prd{x:A}P(f(x))\Big)
    \end{equation*}
    given by $h\mapsto h\circ f$ is an equivalence.
  \item For any type $X$ the map
    \begin{equation*}
      (B\to X)\to (A\to X)
    \end{equation*}
    given by $g\mapsto g\circ f$ is an equivalence. 
  \end{enumerate}
\end{thm}

\begin{proof}
To show that (i) implies (ii), we first recall from \cref{lem:coherently-invertible} that any equivalence is also coherently invertible. Therefore $f$ comes equipped with
\begin{align*}
g & : B \to A\\
G & : f\circ g \htpy \idfunc[B] \\
H & : g\circ f \htpy \idfunc[A] \\
K & : G\cdot f \htpy f\cdot H.
\end{align*}
Then we define the inverse of $\blank\circ f$ to be the map
\begin{equation*}
\varphi:\Big(\prd{x:A}P(f(x))\Big)\to\Big(\prd{y:B}P(y)\Big)
\end{equation*}
given by $h\mapsto \lam{y}\tr_P(G(y),h(g(y)))$. 

To see that $\varphi$ is a section of $\blank\circ f$, let $h:\prd{x:A}P(f(x))$. By function extensionality it suffices to construct a homotopy $\varphi(h)\circ f\htpy h$. In other words, we have to show that
\begin{equation*}
\tr_P(G(f(x)),h(g(f(x)))=h(x)
\end{equation*}
for any $x:A$. Now we use the additional homotopy $K$ from our assumption that $f$ is coherently invertible. Since we have $K(x):G(f(x))=\ap{f}{H(x)}$ it suffices to show that
\begin{equation*}
\tr_P(\ap{f}{H(x)},h(g(f(x))))=h(x).
\end{equation*}
A simple path-induction argument yields that
\begin{equation*}
\tr_P(\ap{f}{p})\htpy \tr_{P\circ f}(p)
\end{equation*}
for any path $p:x=y$ in $A$, so it suffices to construct an identification
\begin{equation*}
\tr_{P\circ f}(H(x),h(g(f(x))))=h(x).
\end{equation*}
We have such an identification by $\apd{h}{H(x)}$.

To see that $\varphi$ is a retraction of $\blank\circ f$, let $h:\prd{y:B}P(y)$. By function extensionality it suffices to construct a homotopy $\varphi(h\circ f)\htpy h$. In other words, we have to show that
\begin{equation*}
\tr_P(G(y),h(f(g(y))))=h(y)
\end{equation*}
for any $y:B$. We have such an identification by $\apd{h}{G(y)}$. This completes the proof that (i) implies (ii).

Note that (iii) is an immediate consequence of (ii), since we can just choose $P$ to be the constant family $X$.

It remains to show that (iii) implies (i). Suppose that
\begin{equation*}
\blank\circ f:(B\to X)\to (A\to X)
\end{equation*}
is an equivalence for every type $X$. Then its fibers are contractible by \cref{thm:contr_equiv}. In particular, choosing $X\jdeq A$ we see that the fiber
\begin{equation*}
\fib{\blank\circ f}{\idfunc[A]}\jdeq \sm{h:B\to A}h\circ f=\idfunc[A]
\end{equation*}
is contractible. Thus we obtain a function $h:B\to A$ and a homotopy $H:h\circ f\htpy\idfunc[A]$ showing that $h$ is a retraction of $f$. We will show that $h$ is also a section of $f$. To see this, we use that the fiber
\begin{equation*}
\fib{\blank\circ f}{f}\jdeq \sm{i:B\to B} i\circ f=f
\end{equation*}
is contractible (choosing $X\defeq B$). 
Of course we have $(\idfunc[B],\refl{f})$ in this fiber. However we claim that there also is an identification $p:(f\circ h)\circ f=f$, showing that $(f\circ h,p)$ is in this fiber, because
\begin{align*}
(f\circ h)\circ f & \jdeq f\circ (h\circ f) \\
& = f\circ \idfunc[A] \\
& \jdeq f
\end{align*}
From the contractibility of the fiber we obtain an identification $(\idfunc[B],\refl{f})=(f\circ h,p)$. In particular we obtain that $\idfunc[B]=f\circ h$, showing that $h$ is a section of $f$.
\end{proof}

\subsection{The strong induction principle of \texorpdfstring{$\N$}{ℕ}}
\index{strong induction principle!of N@{of $\N$}|(}
\index{natural numbers!strong induction principle|(}

In the final application of the function extensionality principle we prove the strong induction principle for the type of natural numbers. Function extensionality is used to derive the computation rules of the strong induction principle.

\begin{thm}[Strong induction for the natural numbers]\label{thm:strong-ind-N}
  \index{strong induction principle!of N@{of $\N$}|textbf}
  \index{natural numbers!strong induction principle|textbf}
  Consider a type family $P$ over $\N$ equipped with
  \begin{align*}
    p_0 & : P(0) \\
    p_S & : \prd{n:\N}\Big(\prd{m:\N}(m\leq n)\to P(m)\Big)\to P(n+1).
  \end{align*}
  Then there is a dependent function\index{strong-ind@{$\strongindN$}|textbf}
  \begin{equation*}
    \strongindN(p_0,p_S) : \prd{n:\N}P(n)
  \end{equation*}
  that satisfies the following computation rules
  \begin{align*}
    \strongindN(p_0,p_S,0) & = p_0 \\
    \strongindN(p_0,p_S,n+1) & = p_S(n,(\lam{m}\lam{p}\strongindN(p_0,p_S,m))).
  \end{align*}
\end{thm}

In order to construct $\strongindN(p_0,p_S)$, we first define the type family $\tilde{P}$ over $\N$ by
\begin{equation*}
  \tilde{P}(n)\defeq \prd{m:\N} (m\leq n)\to P(m).
\end{equation*}
The idea is then to first use $p_0$ and $p_S$ to construct
\begin{align*}
  \tilde{p}_0 & : \tilde{P}(0) \\
  \tilde{p}_S & :\prd{n:\N}\tilde{P}(n)\to\tilde{P}(n+1).
\end{align*}
The ordinary induction principle of $\N$ then gives a function
\begin{equation*}
  \indN(\tilde{p}_0,\tilde{p}_S):\prd{n:\N}\tilde{P}(n),
\end{equation*}
which can be used to define a function $\prd{n:\N}P(n)$.

Before we start by the proof of \cref{thm:strong-ind-N} we state two lemmas in which we construct $\tilde{p}_0$ and $\tilde{p}_S$ with computation rules of their own. We will assume a type family $P$ over $\N$ equipped with
  \begin{align*}
    p_0 & : P(0) \\
    p_S & : \prd{n:\N}\tilde{P}(n) \to P(n+1),
  \end{align*}
as in the hypotheses of \cref{thm:strong-ind-N}.

\begin{lem}
  There is an element $\tilde{p}_0:\tilde{P}(0)$ that satisfies the judgmental equality
  \begin{equation*}
    \tilde{p}_0(0,p)\jdeq p_0
  \end{equation*}
  for any $p:0\leq 0$.
\end{lem}

\begin{proof}
  The fact that we have such a dependent function $\tilde{p}_0$ follows immediately by induction on $m$ and $p:m\leq 0$.
\end{proof}

\begin{lem}\label{lem:succ-strong-ind-N}
  There is a function
  \begin{equation*}
    \tilde{p}_S : \prd{n:\N}\tilde{P}(n)\to\tilde{P}(n+1)
  \end{equation*}
  equipped with
  \begin{enumerate}
  \item an identification
    \begin{equation*}
      \tilde{p}_S(n,H,m,p) = H(m,q)
    \end{equation*}
    for every $H:\tilde{P}(n)$ and every $p:m\leq n+1$ and $q:m\leq n$, and
  \item an identification
    \begin{equation*}
      \tilde{p}_S(n,H,n+1,p) = p_S(n,H)
    \end{equation*}
    for every $p:n+1\leq n+1$.
  \end{enumerate}
\end{lem}

\begin{proof}
  To define the function $\tilde{p}_S(n,H)$, note that there is a function
  \begin{equation*}
    f : (m\leq n+1)\to (m\leq n)+(m=n+1)\tag{\textasteriskcentered}
  \end{equation*}
  which can be defined by induction on $n$ and $m$. Using the fact that the domain and codomain of this map are both propositions, this function is easily seen to be an equivalence. Therefore we define first a function
  \begin{equation*}
    h(n,H) :\prd{m:\N} ((m\leq n)+(m=n+1))\to P(m)
  \end{equation*}
  by case analysis on $x:(m\leq n)+(m=n+1)$. There are two cases to consider: one where we have $q:m\leq n$, and one where we have $q:m=n+1$. Note that in the second case it suffices to make a definition for $q\jdeq \refl{}$. Therefore we define
  \begin{equation*}
    h(n,H,m,x) =
    \begin{cases}
      H(m,q) & \text{if }x\jdeq\inl(q)\\
      p_S(n,H) & \text{if }x\jdeq\inr(\refl{}).
    \end{cases}
  \end{equation*}
  Now we define $\tilde{p}_S$ by
  \begin{equation*}
    \tilde{p}_S(n,H,m,p)\defeq h(n,H,m,f(p)),
  \end{equation*}
  where $f:(m\leq n+1)\to (m\leq n)+(m=n+1)$ is the map we mentioned in (\textasteriskcentered).
  
  To construct the identifications claimed in (i) and (ii), note that there is an equivalence
  \begin{equation*}
    \big(\tilde{p}_S(n,H,m,p)=y\big)\simeq \big(h(n,H,m,x)=y\big),
  \end{equation*}
  for any $y:P(m)$. This equivalence is obtained from the fact that $f(p)=x$ for any $x:(m\leq n)+(m=n+1)$, i.e., the fact that $(m\leq n)+(m=n+1)$ is a proposition. Now the identifications in (i) and (ii) are obtained as a simple consequence of the computation rule for coproducts.
\end{proof}

We are now ready to finish the proof of \cref{thm:strong-ind-N}.

\begin{proof}[Proof of \cref{thm:strong-ind-N}]
  Using $\tilde{p}_0$ and $\tilde{p}_S$, we obtain by induction on $n$ a function
  \begin{equation*}
    \tilde{s}:\prd{n:\N}\tilde{P}(n)
  \end{equation*}
  satisfying the computation rules
  \begin{align*}
    \tilde{s}(0) & \jdeq \tilde{p}_0 \\
    \tilde{s}(n+1) & \jdeq \tilde{p}_S(n,\tilde{s}(n)).
  \end{align*}
  Now we define
  \begin{equation*}
    \strongindN(p_0,p_S,n) \defeq \tilde{s}(n,n,\reflleqN(n)),
  \end{equation*}
  where $\reflleqN(n):n\leq n$ is the proof of reflexivity of $\leq$.

  It remains to show that $\strongindN$ satisfies the computation rules of the strong induction principle. The identification that computes $\strongindN$ at $0$ is easy to obtain, because we have the judgmental equalities
  \begin{align*}
    \strongindN(p_0,p_S,0) & \jdeq \tilde{s}(0,0,\reflleqN(0)) \\
                                & \jdeq \tilde{p}_{0}(0,\reflleqN(0)) \\
                                & \jdeq p_0.
  \end{align*}
  To construct the identification that computes $\strongindN$ at a successor, we start by a similar computation:
  \begin{align*}
    \strongindN(p_0,p_S,n+1) & \jdeq \tilde{s}(n+1,n+1,\reflleqN(n+1)) \\
                                   & \jdeq \tilde{p}_S(n,\tilde{s}(n),n+1,\reflleqN(n+1)) \\
    & = p_S(n,\tilde{s}(n)).
  \end{align*}
  The last identification is obtained from \cref{lem:succ-strong-ind-N} (ii).
  Therefore we see that, in order to show that
  \begin{equation*}
    p_S(n,\tilde{s}(n))=p_S(n,(\lam{m}\lam{p}\tilde{s}(m,m,\reflleqN(m)))),
  \end{equation*}
  we need to prove that
  \begin{equation*}
    \tilde{s}(n)=\lam{m}\lam{p}\tilde{s}(m,m,\reflleqN(m)).
  \end{equation*}
  Here we apply function extensionality, so it suffices to show that
  \begin{equation*}
    \tilde{s}(n,m,p)=\tilde{s}(m,m,\reflleqN(m))
  \end{equation*}
  for every $m:\N$ and $p:m\leq n$. We proceed by induction on $n:\N$. The base case is trivial. For the inductive step, we note that
  \begin{align*}
    \tilde{s}(n+1,m,p)=\tilde{p}_S(n,\tilde{s}(n),m,p)=\begin{cases}\tilde{s}(n,m,p) & \text{if }m\leq n \\
    p_S(n,\tilde{s}(n)) & \text{if }m=n+1.\end{cases}
  \end{align*}
  Therefore it follows by the inductive hypothesis that
  \begin{equation*}
    \tilde{s}(n+1,m,p)=\tilde{s}(m,m,\reflleqN(m))
  \end{equation*}
  if $m\leq n$ holds. In the remaining case, where $m=n+1$, note that we have
  \begin{align*}
    \tilde{s}(n+1,n+1,\reflleqN(n+1)) & = \tilde{p}_S(n,\tilde{s}(n),n+1,\reflleqN(n+1)) \\
    & = p_S(n,\tilde{s}(n)).
  \end{align*}
  Therefore we see that we also have an identification
  \begin{equation*}
    \tilde{s}(n+1,m,p)=\tilde{s}(m,m,\reflleqN(m))
  \end{equation*}
  when $m=n+1$. This completes the proof of the computation rules for the strong induction principle of $\N$.
\end{proof}
\index{strong induction principle!of N@{of $\N$}|)}
\index{natural numbers!strong induction principle|)}

\begin{exercises}
  \exitem \label{ex:is-equiv-inv-htpy}Show that the functions\index{inv-htpy@{$\invhtpy$}!is an equivalence}\index{concat-htpy@{$\concathtpy$}!is a family of equivalences}\index{concat-htpy'@{$\concathtpy'$}!is a family of equivalences}\index{is an equivalence!inv-htpy@{$\invhtpy$}}\index{is an equivalence!concat-htpy(H)@{$\concathtpy(H)$}}\index{is an equivalence!concat-htpy'(K)@{$\concathtpy'(K)$}}
  \begin{align*}
    \invhtpy & : (f \htpy g) \to (g \htpy f) \\
    \concathtpy(H) & : (g \htpy h) \to (f \htpy h) \\
    \concathtpy'(K) & : (f \htpy g) \to (f \htpy h)
  \end{align*}
  are equivalences for every $f,g,h : \prd{x:A}B(x)$. Here, $\concathtpy'(K)$ is the function defined by $H\mapsto \ct{H}{K}$.
  \exitem Characterize the identity types of the following types:
  \begin{subexenum}
  \item The type $\sm{h:A\to B}h(a)=b$ of \define{pointed maps}\index{pointed map|textbf}, where $a:A$ and $b:B$ are given.
  \item The type $\sm{h:A\to B}f\htpy g\circ h$ of commuting triangles
    \begin{equation*}
      \begin{tikzcd}[column sep=tiny]
        A \arrow[rr,"h"] \arrow[dr,swap,"f"] & & B \arrow[dl,"g"] \\
        & X,
      \end{tikzcd}
    \end{equation*}
    where $f:A\to X$ and $g:B\to X$ are given.
  \item The type $\sm{h:X\to Y}h\circ f\htpy g$ of commuting triangles
    \begin{equation*}
      \begin{tikzcd}[column sep=tiny]
        & A \arrow[dl,swap,"f"] \arrow[dr,"g"] \\
        X \arrow[rr,swap,"h"] & & Y,
      \end{tikzcd}
    \end{equation*}
    where $f:A\to X$ and $g:A\to Y$ are given.
  \item The type $\sm{i:A\to X}\sm{j:B\to Y}j\circ f\htpy g\circ i$ of commuting squares
    \begin{equation*}
      \begin{tikzcd}
        A \arrow[d,swap,"f"] \arrow[r,"i"] & X \arrow[d,"g"] \\
        B \arrow[r,swap,"j"] & Y,
      \end{tikzcd}
    \end{equation*}
    where $f:A\to B$ and $g:X\to Y$ are given.
  \end{subexenum}
  \exitem \label{ex:isprop_istrunc}
  \begin{subexenum}
  \item Show that for any type $A$ the type $\iscontr(A)$ is a proposition\index{is-contr(A)@{$\iscontr(A)$}!is a proposition}\index{is contractible!is a property}\index{is a proposition!is-contr(A)@{$\iscontr(A)$}}.
  \item Show that for any type $A$ and any $k\geq-2$, the type $\istrunc{k}(A)$ is a proposition.\index{istrunc@{$\istrunc{k}$}!is a proposition}\index{is a proposition!istrunc(A)@{$\istrunc{k}(A)$}}
  \end{subexenum}
  \exitem \label{ex:isprop_isequiv}Let $f:A\to B$ be a function.
  \begin{subexenum}
  \item Show that if $f$ is an equivalence, then the type $\sm{g:B\to A}f\circ g\htpy \idfunc$ of sections of $f$ is contractible.
  \item Show that if $f$ is an equivalence, then the type $\sm{h:B\to A}h\circ f\htpy \idfunc$ of retractions of $f$ is contractible.
  \item Show that $\isequiv(f)$ is a proposition.\index{is-equiv(f)@{$\isequiv(f)$}!is a proposition}\index{is a proposition!is-equiv(f)@{$\isequiv(f)$}}
  \item Show that for any two equivalences $e,e':A\simeq B$, the canonical map
    \begin{equation*}
      (e=e')\to (e\htpy e')
    \end{equation*}
    is an equivalence.
  \item Show that the type $A\simeq B$ is a $k$-type if both $A$ and $B$ are $k$-types.\index{is a truncated type!A simeq B@{$A\simeq B$}}
  \end{subexenum}
  \exitem
  \begin{subexenum}
  \item Show that $\pathsplit(f)$\index{path-split!is a proposition}\index{is a proposition!is-path-split(f)@{$\pathsplit(f)$}} and $\iscohinvertible(f)$\index{is-coh-invertible(f)@{$\iscohinvertible(f)$}!is a proposition}\index{is a proposition!is-coh-invertible(f)@{$\iscohinvertible(f)$}} are propositions for any map $f:A\to B$. Conclude that we have equivalences\index{is-equiv(f)@{$\isequiv(f)$}!is-equiv(f) path-split(f)@{$\isequiv(f)\eqvsym\pathsplit(f)$}}\index{is-equiv(f)@{$\isequiv(f)$}!is-equiv(f) is-coh-invertible(f)@{$\isequiv(f)\eqvsym\iscohinvertible(f)$}}
    \begin{equation*}
      \isequiv(f) \eqvsym \pathsplit(f) \eqvsym \iscohinvertible(f).
    \end{equation*}
  \item \label{ex:idfunc_autohtpy}Construct for any type $A$ an equivalence\index{has-inverse(f)@{$\hasinverse(f)$}!has-inverse(id) id htpy id@{$\hasinverse(\idfunc)\simeq (\idfunc\htpy\idfunc)$}}
    \begin{equation*}
      \eqv{\hasinverse(\idfunc[A])}{\Big(\idfunc[A]\htpy\idfunc[A]\Big)}.
    \end{equation*}
    Note: We will use this fact in \cref{ex:is_invertible_id_S1} to show that there
    are types for which $\hasinverse(\idfunc[A])\not\eqvsym\isequiv(\idfunc[A])$.
  \end{subexenum}
  \exitem \label{ex:up-emptyt}Consider a type $A$. Show that the following are equivalent:
  \begin{enumerate}
  \item The type $A$ is empty.
  \item \label{item:dup-empty}The type $\prd{x:A}P(x)$ is contractible for any family $P$ of types over $A$. This property is the \define{dependent universal property of an empty type}\index{dependent universal property!of empty types|textbf}\index{empty type!dependent universal property|textbf}.
  \item \label{item:up-empty}The type $A\to X$ is contractible for any type $X$. This property is the \define{universal property of an empty type}\index{universal property!of empty types|textbf}\index{empty type!universal property|textbf}.
  \end{enumerate}
  \exitem \label{ex:up-unit}Consider a type $A$. Show that the following are equivalent:
  \begin{enumerate}
  \item \label{item:is-contr}The type $A$ is contractible.
  \item \label{item:dup-unit}The type $A$ comes equipped with a point $a:A$, and the map
    \begin{equation*}
      \Big(\prd{x:A}P(x)\Big)\to P(a)
    \end{equation*}
    given by $f\mapsto f(a)$ is an equivalence for any type family $P$ over $A$. This property is the \define{dependent universal property of a contractible type}\index{dependent universal property!of contractible types|textbf}\index{contractible type!dependent universal property|textbf}.
  \item \label{item:up-unit}The type $A$ comes equipped with a point $a:A$, and the map
    \begin{equation*}
      (A\to X)\to X
    \end{equation*}
    given by $f\mapsto f(a)$ is an equivalence for any type $X$. This property is the \define{universal property of a contractible type}\index{universal property!of contractible types|textbf}\index{contractible type!universal property|textbf}.
  \item The type $A$ comes equipped with a point $a:A$, and the map
    \begin{equation*}
      (A\to A)\to A
    \end{equation*}
    given by $f\mapsto f(a)$ is an equivalence.
  \item \label{item:is-equiv-diag-universal}The map
    \begin{equation*}
      X\to (A\to X)
    \end{equation*}
    given by $x\mapsto\lam{y}x$ is an equivalence for any type $X$.
  \item \label{item:is-equiv-diag}The map
    \begin{equation*}
      A\to (A\to A)
    \end{equation*}
    given by $x\mapsto\lam{y}x$ is an equivalence.
  \end{enumerate}
  \exitem \label{ex:up-coproduct}
  Consider two types $A$ and $B$. Show that the map
  \begin{equation*}
    \Big(\prd{z:A+B}P(z)\Big)\to\Big(\prd{x:A}P(\inl(x))\Big)\times\Big(\prd{y:B}P(\inr(b))\Big)
  \end{equation*}
  given by $f\mapsto (f\circ \inl,f\circ \inr)$ is an equivalence for any type family $P$ over $A+B$. This property is the \define{dependent universal property of the coproduct of $A$ and $B$}\index{dependent universal property!of coproducts|textbf}\index{coproduct!dependent universal property|textbf}. Conclude that the map
  \begin{equation*}
    (A+B\to X)\to (A\to X)\times (B\to X)
  \end{equation*}
  given by $f\mapsto (f\circ \inl,f\circ \inr)$ is an equivalence for any type $X$. This latter property is the \define{universal property of the coproduct of $A$ and $B$}\index{universal property!of coproducts|textbf}\index{coproduct!universal property|textbf}.
  \exitem\label{ex:uniqueness-identity-type}Consider a type $A$ equipped with an element $a:A$ and consider a type family $B$ over $A$ equipped with an element $b:B(a)$. Show that the following are equivalent:\index{identity type!uniqueness}
  \begin{enumerate}
  \item The map
    \begin{equation*}
      \ev_b:\Big(\prd{x:A}B(x)\to C(x)\Big)\to C(a)
    \end{equation*}
    given by $\ev_b(h)\defeq h(a,b)$ is an equivalence for any type family $C$ over $A$.
  \item The map
    \begin{equation*}
      h:\prd{x:A}(a=x)\to B(x)
    \end{equation*}
    given by $h(a,\refl{})\defeq b$ is an equivalence.
  \end{enumerate}
  \exitem Prove the \define{universal property of $\N$}\index{universal property!of N@{of $\N$}|textbf}\index{natural numbers!universal property|textbf}: For any type $X$ equipped with $x:X$ and $f:X\to X$, the type
  \begin{equation*}
    \sm{h:\N\to X} (h(\zeroN)= x)\times (h\circ \succN\htpy f\circ h)
  \end{equation*}
  is contractible.
  \exitem Show that $\N$ satisfies \define{ordinal induction}\index{ordinal induction!of N@{of $\N$}}\index{natural numbers!ordinal induction|textbf}, i.e., construct for any type family $P$ over $\N$ a function $\ordindN$ of type
  \begin{equation*}
    \Big(\prd{k:\N} \Big(\prd{m:\N} (m< k) \to P(m)\Big)\to P(k)\Big) \to \prd{n:\N}P(n).
  \end{equation*}
  Moreover, prove that
  \begin{equation*}
    \ordindN(h,n)=h(n,\lam{m}\lam{p}\ordindN(h,m))
  \end{equation*}
  for any $n:\N$ and any $h:\prd{k:\N}\Big(\prd{m:\N}(m<k)\to P(m)\Big)\to P(k)$.
  \exitem \label{ex:equiv-pi} 
  \begin{subexenum}
  \item Consider a family of $k$-truncated maps $f_i:A_i\to B_i$ indexed by $i:I$. Show that the map
    \begin{equation*}
      \lam{h}\lam{i}f_i(h(i)): \Big(\prd{i:I}A_i\Big)\to\Big(\prd{i:I}B_i\Big)
    \end{equation*}
    is also $k$-truncated.
  \item Consider an equivalence $e:I\simeq J$, and a family of equivalences $f_i:A_i\simeq B_{e(i)}$ indexed by $i:I$, where $A$ is a family of types indexed by $I$ and $B$ family of types indexed by $J$. Show that the map
    \begin{equation*}
      \lam{h}\lam{j} f_{e^{-1}(j)}(h(e^{-1}(j))) : \Big(\prd{i:I}A_i\Big)\to\Big(\prd{j:J}B_j\Big)
    \end{equation*}
    is an equivalence.
  \item Consider a family of maps $f_i:A_i\to B_i$ indexed by $i:I$. Show that the following are equivalent:
    \begin{enumerate}
    \item Each $f_i$ is $k$-truncated.
    \item For every map $\alpha:X\to I$, the map
      \begin{equation*}
        \lam{h}\lam{x}f_{\alpha(x)}(h(x)):\Big(\prd{x:X}A_{\alpha(x)}\Big)\to\Big(\prd{x:X}B_{\alpha(x)}\Big)
      \end{equation*}
      is $k$-truncated.
    \end{enumerate}
  \item \label{ex:equiv-postcomp}Show that for any map $f:A\to B$ the following are equivalent:
    \begin{enumerate}
    \item The map $f$ is $k$-truncated.
    \item For every type $X$, the postcomposition function
      \begin{equation*}
        f\circ\blank : (X\to A)\to (X\to B)
      \end{equation*}
      is $k$-truncated.
    \end{enumerate}
    In particular, $f$ is an equivalence if and only if $f\circ\blank$ is an equivalence, and $f$ is an embedding if and only if $f\circ\blank$ is an embedding.
  \end{subexenum}
  \exitem Show that \emph{$\Pi$-types distribute over coproducts}\index{distributivity!of P over coproducts@{of $\Pi$ over coproducts}}\index{dependent function type!distributivity of P over coproducts@{distributivity of $\Pi$ over coproducts}}\index{coproduct!distributivity of P over coproducts@{distributivity of $\Pi$ over coproducts}}, i.e., construct for any type $X$ and any two families $A$ and $B$ over $X$ an equivalence from the type $\prd{x:X}A(x)+B(x)$ to the type
  \begin{equation*}
    \sm{f:X\to\Fin{2}}\Big(\prd{x:X}A(x)^{f(x)=0}\Big)\times\Big(\prd{x:X}B(x)^{f(x)=1}\Big).
  \end{equation*}
  \exitem \label{ex:sec_retr}Consider a commuting triangle 
  \begin{equation*}
    \begin{tikzcd}[column sep=tiny]
      A \arrow[rr,"h"] \arrow[dr,swap,"f"] & & B \arrow[dl,"g"] \\
      & X
    \end{tikzcd}
  \end{equation*}
  with $H:f\htpy g\circ h$.
  \begin{subexenum}
  \item Show that if $h$ has a section, then $\sections(g)$ is a retract of $\sections(f)$.
  \item Show that if $g$ has a retraction, then $\retractions(h)$ is a retract of $\sections(f)$.
  \end{subexenum}
  \exitem \label{ex:triangle_fib}For any two maps $f:A\to X$ and $g:B\to X$, define the type of \define{morphisms from $f$ to $g$ over $X$}\index{morphism from f to g over X@{morphism from $f$ to $g$ over $X$}|textbf} by\index{hom X (f,g)@{$\homslice_X(f,g)$}|textbf}
  \begin{equation*}
    \homslice_X(f,g)\defeq \sm{h:A\to B} f\htpy g\circ h.
  \end{equation*}
  In other words, the type $\homslice_X(f,g)$ is the type of maps $h:A\to B$ equipped with a homotopy witnessing that the triangle
  \begin{equation*}
    \begin{tikzcd}[column sep=tiny]
      A \arrow[dr,swap,"f"] \arrow[rr,dashed,"h"] & & B \arrow[dl,"g"] \\
      & X
    \end{tikzcd}
  \end{equation*}
  commutes.
  \begin{subexenum}
  \item \label{ex:pi-fib}Consider a family $P$ of types over $X$. Show that the map
    \begin{equation*}
      \Big(\prd{x:X}\fib{f}{x}\to P(x)\Big)\to\Big(\prd{a:A}P(f(a))\Big)
    \end{equation*}
    given by $h\mapsto h_{f(a)}(a,\refl{f(a)})$ is an equivalence. 
  \item Construct three equivalences $\alpha$, $\beta$, and $\gamma$ as shown in the following diagram, and show that this triangle commutes:
    \begin{equation*}
      \begin{tikzcd}[column sep=-3em]
        & \homslice_X(f,g) \arrow[dl,swap,"\alpha"] \arrow[dr,"\beta"] & \phantom{\Big(\prd{x:X}\fib{f}{x}\to\fib{g}{x}\Big)} \\
        \Big(\prd{x:X}\fib{f}{x}\to\fib{g}{x}\Big) \arrow[rr,swap,"\gamma"] & & \prd{a:A}\fib{g}{f(a)}.
      \end{tikzcd}
    \end{equation*}
    Given a morphism $(h,H):\homslice_X(f,g)$ over $X$, we also say that $\alpha(h,H)$ is its \define{action on fibers}\index{action on fibers|textbf}\index{morphism from f to g over X@{morphism from $f$ to $g$ over $X$}!action on fibers|textbf}.
  \item \label{ex:fam-equiv}Given $(h,H):\homslice_X(f,g)$, show that the following are equivalent:
    \begin{enumerate}
    \item The map $h:A\to B$ is an equivalence.
    \item The action on fibers
      \begin{equation*}
        \alpha(h,H):\prd{x:X}\fib{f}{x}\to\fib{g}{x}
      \end{equation*}
      is a family of equivalences.
    \item The precomposition function
      \begin{equation*}
        \blank\circ (h,H) : \homslice_X(g,i)\to\homslice_X(f,i)
      \end{equation*}
      given by $(k,K)\circ (h,H) \defeq (k\circ h,\ct{H}{(K\cdot h)})$ is an equivalence for each map $i:C\to X$.
    \end{enumerate}
    Conclude that the type $\sm{h:\eqv{A}{B}} f\htpy g\circ h$ is equivalent to the type of families of equivalences
    \begin{equation*}
      \prd{x:X}\fib{f}{x}\eqvsym\fib{g}{x}.
    \end{equation*} 
  \end{subexenum}
\exitem \label{ex:iso_equiv}Let $A$ and $B$ be sets. Show that type $\eqv{A}{B}$ of equivalences from $A$ to $B$ is equivalent to the type $A\cong B$ of \define{isomorphisms}\index{isomorphism!of sets|textbf}\index{set!isomorphism|textbf} from $A$ to $B$, i.e., the type of quadruples $(f,g,H,K)$ consisting of
  \begin{align*}
    f & : A\to B \\
    g & : B\to A \\
    H & : f\circ g = \idfunc[B] \\
    K & : g\circ f = \idfunc[A].
  \end{align*}
\exitem Suppose that $A:I\to \UU$ is a type family over a set $I$ with decidable equality. Show that
  \begin{equation*}
    \Big(\prd{i:I}\iscontr(A_i)\Big)\leftrightarrow \iscontr\Big(\prd{i:I}A_i\Big).
  \end{equation*}
  \exitem \label{ex:retracts-as-limits}(Shulman) Consider a section-retraction pair
  \begin{equation*}
    \begin{tikzcd}
      A \arrow[r,"i"] & X \arrow[r,"r"] & A
    \end{tikzcd}
  \end{equation*}
  with $H:r\circ i\htpy \idfunc$ and define $f\defeq i\circ r$. Construct an equivalence
  \begin{equation*}
    A\simeq\sm{x:\N\to X}\prd{n:\N}f(x_{n+1})=x_n.
  \end{equation*}
\end{exercises}
\index{function extensionality|)}
\index{axiom!function extensionality|)}


\section{Propositional truncations}\label{sec:propositional-truncation}\label{chap:propositional-truncation}
\index{propositional truncation|(}

It is common in mathematics to express the property that a certain type of objects is inhabited, without imposing extra structure on those objects. For example, when we assert the property that a set is finite, then we only claim that there exists a bijection to a standard finite set $\{0,\ldots,n-1\}$ for some $n$, not that the set is equipped with such a bijection. There is indeed a conceptual difference between a finite set and a set equipped with a bijection to a standard finite set. The latter concept is that of a finite \emph{totally ordered} set. This difference is due to the fact that finiteness is a property, whereas there may be many different bijections to a standard finite set.

A similar observation can be made in the case of the image of a map. Note that being in the image of a given map $f:A\to B$ is a property. When we claim that $b:B$ is in the image of $f$, then we only claim that the type of $a:A$ such that $f(a)=b$ is inhabited. On the other hand, we saw in \cref{ex:fib_replacement} that the type of $b:B$ equipped with an $a:A$ such that $f(a)=b$ is equivalent to the type $A$, i.e., we have an equivalence
\begin{equation*}
  A\simeq \sm{b:B}\sm{a:A}f(a)=b.
\end{equation*}
Something is clearly off here, because the type $A$ is often not a subtype of the type $B$, while we would expect the image of $f$ to be a subtype of $B$. Therefore we see that the type $\sm{a:A}f(a)=b$ does not quite capture the concept of $b$ being in the image of $f$. The difference is again due to the fact that $\fib{f}{b}$ is often not a proposition, whereas we are looking to express the proposition that the preimage of $f$ at $b$ is inhabited.

To correctly capture the concepts of finiteness and the image of a map in type theory, and many further mathematical concepts, we need a way to assert the \emph{proposition} that a type is inhabited. The proposition that a type $A$ is inhabited is called the propositional truncation of $A$.

\subsection{The universal property of propositional truncations}\label{sec:propositional-truncation-up}

The propositional truncation of a type $A$ is a proposition $\brck{A}$ equipped with a map
\begin{equation*}
  \eta:A\to \brck{A}.
\end{equation*}
This map ensures that if we have an element $a:A$, then the proposition $\brck{A}$ that $A$ is inhabited holds. The complete specification of the propositional truncation includes the universal property of the map $\eta$. In this section we will specify in full generality when a map $f:A\to P$ into a proposition $P$ is a propositional truncation.

\begin{defn}
Let $A$ be a type, and let $f:A\to P$ be a map into a proposition $P$. We say that $f$ \define{is a propositional truncation}\index{is a propositional truncation|textbf}\index{propositional truncation!to be a propositional truncation|textbf} of $A$ if for every proposition $Q$, the precomposition map
\begin{equation*}
\blank\circ f:(P\to Q)\to (A\to Q)
\end{equation*}
is an equivalence. This property of $f$ is called the \define{universal property of the propositional truncation of $A$}\index{universal property!of propositional truncations|textbf}\index{propositional truncation!universal property|textbf}.
\end{defn}

\begin{rmk}\label{ex:prop_equiv}
  Using the fact that equivalences are maps that have contractible fibers, we can reformulate the universal property of the propositional truncation. Note that the fiber of the precomposition map $\blank\circ f:(P\to Q) \to (A \to Q)$ at a map $g:A\to Q$ is the type.
  \begin{equation*}
    \sm{h:P\to Q}h\circ f=g
  \end{equation*}
  Therefore we see that if $f$ satisfies the universal property of the propositional truncation, then these fibers are contractible. In other words, for each map $g:A\to Q$ into a proposition $Q$ there is a unique map $h:P\to Q$ for which $h\circ f=g$. We also say that every map $g:A\to Q$ into a proposition $Q$ \emph{extends} uniquely along $f$, as indicated in the diagram
  \begin{equation*}
    \begin{tikzcd}
      A \arrow[d,swap,"f"] \arrow[dr,"g"] \\
      P \arrow[r,dashed] & Q.
    \end{tikzcd}
  \end{equation*}
\end{rmk}

\begin{rmk}\label{rmk:simplified-up-trunc-Prop}
  For any two propositions $P$ and $P'$, a map $f:P\to P'$ is an equivalence if and only if there is a function $g:P'\to P$. To see this, simply note that any such function $g$ is an inverse of $f$, because any two elements in $P$ and any two elements in $P'$ are equal. 
  
  Note that the type $X\to Q$ is a proposition, for any type $X$ and any proposition $Q$. Using the previous observation, it therefore follows that the map $(P\to Q)\to (A\to Q)$ is an equivalence as soon as there is a map in the converse direction. In other words, to prove that a map $f:A\to P$ into a proposition $P$ satisfies the universal property of the propositional truncation of $A$, it suffices to construct a function
  \begin{equation*}
    (A\to Q)\to (P\to Q)
  \end{equation*}
  for every proposition $Q$.
\end{rmk}

In the following proposition we show that the propositional truncation of a type $A$ is uniquely determined up to equivalence, if it exists. In other words, any two propositional truncations of a type $A$ must be equivalent.

\begin{prp}
  Let $A$ be a type, and consider two maps
  \begin{equation*}
    f:A\to P \qquad\text{and}\qquad f':A\to P'
  \end{equation*}
  into two propositions $P$ and $P'$. If any two of the following three assertions hold, so does the third:\index{3-for-2 property!of propositional truncations}\index{propositional truncation!3-for-2 property}
  \begin{enumerate}
  \item\label{item:f-up-trunc-Prop} The map $f$ is a propositional truncation of $A$.
  \item\label{item:f-up-trunc-Prop'} The map $f'$ is a propositional truncation of $A$.
  \item\label{item:equiv-Prop} There is a (unique) equivalence $P\simeq P'$.
  \end{enumerate}
\end{prp}

\begin{proof}
  We first show that \ref{item:f-up-trunc-Prop} and \ref{item:f-up-trunc-Prop'} together imply \ref{item:equiv-Prop}. If $f$ and $f'$ are both propositional truncations of $A$, then we have maps $P\to P'$ and $P'\to P$ by the universal properties of $f$ and $f'$. Since $P$ and $P'$ are both propositions, it follows that $P\simeq P'$. For the uniqueness claim, note that the type $P\simeq P'$ is itself a proposition.

  Finally we show that \ref{item:equiv-Prop} implies that \ref{item:f-up-trunc-Prop} holds if and only if \ref{item:f-up-trunc-Prop'} holds. Suppose we have an equivalence $P\simeq P'$, let $Q$ be an arbitrary proposition, and consider the triangle
  \begin{equation*}
    \begin{tikzcd}[column sep=-1em]
      \phantom{(P'\to Q)} & (A\to Q) \arrow[dl,dashed] \arrow[dr,dashed] \\
      (P\to Q) \arrow[rr,<->] & & (P'\to Q),
    \end{tikzcd}
  \end{equation*}
  where the fact that $(P\to Q)\leftrightarrow (P'\to Q)$ holds follows from the assumption that $P$ is equivalent to $P'$. We see from this triangle that
  \begin{equation*}
    \Big((A\to Q)\to (P\to Q)\Big)\leftrightarrow\Big((A \to Q) \to (P'\to Q)\Big),
  \end{equation*}
  and this implies that \ref{item:f-up-trunc-Prop} holds if and only if \ref{item:f-up-trunc-Prop'} holds.
\end{proof}

\begin{rmk}
  One might be tempted to think that a type is inhabited if and only if it is nonempty. Recall that a type $A$ is nonempty if it satisfies the property $\neg\neg A$. Indeed, the type $\neg\neg A$ is a proposition, and it comes equipped with a map $A\to\neg\neg A$. It is therefore natural to wonder whether the map $A\to\neg\neg A$ satisfies the universal property of the propositional truncation.

  Recall that we have shown in \cref{ex:dn-monad} that any map $A\to\neg\neg Q$ extends to a map $\neg\neg A\to\neg\neg Q$, as indicated in the diagram
\begin{equation*}
  \begin{tikzcd}
    A \arrow[d] \arrow[dr] \\
    \neg\neg A \arrow[r,dashed] & \neg\neg Q.
  \end{tikzcd}
\end{equation*}
It follows that the natural map
\begin{equation*}
  (\neg\neg A\to\neg\neg Q)\to (A\to \neg\neg Q)
\end{equation*}
given by precomposition by $A\to\neg\neg A$ is an equivalence. However, this only gives us a universal property with respect to doubly negated propositions and there is no way to prove the more general universal property of the propositional truncation for the map $A\to\neg\neg A$. In fact, propositional truncations are not guaranteed to exist in Martin L\"of's dependent type theory, the way it is set up in \cref{chap:type-theory}. We will therefore add new rules to the type theory to ensure their existence.
\end{rmk}

\subsection{Propositional truncations as higher inductive types}\label{sec:propositional-truncation-hit}
\index{higher inductive type!propositional truncation|(}
\index{propositional truncation!as higher inductive type|(}

We have given a specification of the propositional truncation of a type $A$, and we have seen that this specification by a universal property determines the propositional truncation up to equivalence if it exists. However, the propositional truncation is not guaranteed to exist, so we will add new rules to the type theory that ensure that any type has a propositional truncation. We do this by presenting the propositional truncation of a type $A$ as a higher inductive type. The propositional truncation $\brck{A}$ of a type $A$ was one of the first examples of a higher inductive type, along with the circle, which we will discuss in \cref{sec:circle,sec:circle-universal-cover}.

The idea of higher inductive types is similar to the idea of ordinary inductive types, with the added feature that constructors of higher inductive types can also be used to generate \emph{identifications}. In other words, higher inductive types may be specified by two kinds of constructors:
\begin{enumerate}
\item The \emph{point constructors} are used to generate elements of the higher inductive types.
\item The \emph{path constructors} are used to generate identifications between elements of the higher inductive type.
\end{enumerate}
The induction principle of the higher inductive type then tells us how to construct sections of families over it. The rules for higher inductive types therefore come in four sets, just as the rules for ordinary inductive types in \cref{sec:inductive}: the formation rule, the constructors, the induction principle, and the computation rules.

\subsubsection{The formation rules and the constructors}
The formation rule of the propositional truncation postulates that for every type $A$ we can form the propositional truncation of $A$. The formation rule is therefore as follows:\index{[[A]]@{$\brck{A}$}|see {propositional truncation}}
  \begin{prooftree}
    \AxiomC{$\Gamma\vdash A~\type$}
    \UnaryInfC{$\Gamma\vdash \brck{A}~\type$.}
  \end{prooftree}
  Furthermore, we will assume that all universes are closed under propositional truncations. In other words, for any universe $\UU$ we will assume the rules

  \medskip
  \begin{minipage}{.4\textwidth}
    \begin{prooftree}
      \AxiomC{}
      \UnaryInfC{$X:\UU\vdash \brckcheck{X}:\UU$}
    \end{prooftree}
  \end{minipage}
  \begin{minipage}{.5\textwidth}
    \begin{prooftree}
      \AxiomC{}
      \UnaryInfC{$X:\UU\vdash \mathcal{T}(\brckcheck{X})\jdeq\brck{\mathcal{T}(X)}~\type$}
    \end{prooftree}
  \end{minipage}

  \medskip
The constructors of a (higher) inductive type tell what structure the type comes equipped with. In the case of a higher inductive type there may be point constructors and path constructors. The point constructors generate elements of the higher inductive type, and the path constructors generate identifications between those elements. In the case of the propositional truncation, there is one point constructor and one path constructor:\index{propositional truncation!eta : A -> [[A]]@{$\eta:A\to\brck{A}$}}\index{propositional truncation!alpha : P (x y : A) x = y@{$\alpha:\prd{x,y:A}x=y$}}
\begin{align*}
  \eta & : A \to \brck{A}\\*
  \alpha & : \prd{x,y:\brck{A}}x=y.
\end{align*}
The point constructor $\eta$ is sometimes called the \define{unit}\index{propositional truncation!unit of propositional truncation}\index{unit of propositional truncation} of the propositional truncation. It gives us that any element of $A$ also generates an element of $\brck{A}$. The path constructor $\alpha$ simply identifies any two elements of $\brck{A}$. Therefore it follows immediately that $\brck{A}$ is a proposition.

\begin{lem}
  For any type $A$, the type $\brck{A}$ is a proposition.\index{propositional truncation!is a proposition}\index{is a proposition!propositional truncation}\hfill $\square$ 
\end{lem}

\subsubsection{The induction principle and computation rules}
\index{induction principle!of propositional truncation|(}
\index{propositional truncation!induction principle|(}
The induction principle for the propositional truncation tells us how to construct dependent functions
\begin{equation*}
  h:\prd{t:\brck{A}}Q(t).
\end{equation*}
The induction principle will imply that such a dependent function $h$ is entirely determined by its behavior on the constructors of $\brck{A}$. The type $\brck{A}$ has two constructors: a point constructor $\eta$ and a path constructor $\alpha$, so we have two cases to consider:
\begin{enumerate}
\item Applying $h$ to points of the form $\eta(a)$ gives us a dependent function
  \begin{equation*}
    h\circ \eta : \prd{a:A}Q(\eta(a)).
  \end{equation*}
  The induction principle of $\brck{A}$ has therefore the requirement that we can construct
  \begin{equation*}
    f:\prd{a:A}Q(\eta(a))
  \end{equation*}
\item To apply $h$ to the paths $\alpha(x,y)$, we need to use the dependent action on paths from \cref{defn:apd}. For each $x,y:\brck{A}$ we obtain an identification
  \begin{equation*}
    \apd{h}{\alpha(x,y)}:\tr_Q(\alpha(x,y),h(x))=h(y)
  \end{equation*}
  in the type $Q(y)$. Note, however, that $h(x)$ and $h(y)$ are not determined by our choice of $f:\prd{a:A}Q(\eta(a))$. The second requirement of the induction principle of $\brck{A}$ is therefore that, no matter what values $h$ takes, they must always be related via the dependent action on paths of $h$. This second requirement is therefore that
  \begin{equation*}
    \tr_P(\alpha(x,y),u)=v
  \end{equation*}
  for any $u:Q(x)$ and $v:Q(y)$. 
\end{enumerate}

\begin{defn}
  The \define{induction principle}\index{induction principle!of propositional truncation|textbf}\index{propositional truncation!induction principle|textbf} of the propositional truncation $\brck{A}$ of $A$ asserts that for any family $Q$ of types over $\brck{A}$, if we have
  \begin{equation*}
    f:\prd{a:A}Q(\eta(a))
  \end{equation*}
  and if we can construct identifications
  \begin{equation*}
    \tr_Q(\alpha(x,y),u)=v
  \end{equation*}
  for any $u:Q(x)$, $v:Q(y)$ and any $x,y:\brck{A}$, then we obtain a dependent function
  \begin{equation*}
    h:\prd{t:\brck{A}}Q(t)
  \end{equation*}
  equipped with a homotopy $h\circ\eta\htpy f$.
\end{defn}

\begin{rmk}
  In fact, a family $Q$ over $\brck{A}$ satisfies the second requirement in the induction principle of the propositional truncation if and only if $Q$ is a family of propositions. To see this, simply note that transporting along $\alpha(x,y)$ is an embedding. Therefore we have
  \begin{equation*}
    (\tr_Q(\alpha(x,y),u)=\tr_Q(\alpha(x,y),v))\simeq (u=v)
  \end{equation*}
  for any $u,v:Q(x)$. By assumption, there is an identification on the left hand side, so any two elements $u$ and $v$ in $Q(x)$ are equal.

  Since the induction principle of the propositional truncation is only applicable to families of propositions over $\brck{A}$, it also follows that there are no interesting computation rules to state: any identification in a proposition just holds.
\end{rmk}
\index{induction principle!of propositional truncation|)}
\index{propositional truncation!induction principle|)}
\index{higher inductive type!propositional truncation|)}
\index{propositional truncation!as higher inductive type|)}

\subsubsection{The universal property}
\index{universal property!of propositional truncations|(}
\index{propositional truncation!universal property|(}
We have now completed the description of the propositional truncation as a higher inductive type, so it is time to show that it meets the specification we gave for the propositional truncations. In other words, we have to show that the map $\eta:A\to\brck{A}$ satisfies the universal property of the propositional truncation.

\begin{thm}
  The map $\eta:A\to\brck{A}$ satisfies the universal property of the propositional truncation.
\end{thm}

\begin{proof}
  In order to prove that $\eta:A\to\brck{A}$ satisfies the universal property of the propositional truncation of $A$, it suffices to construct a map
  \begin{equation*}
    (A\to Q)\to (\brck{A}\to Q)
  \end{equation*}
  for any proposition $Q$. Consider a map $f:A\to Q$. Then we will construct a function $\brck{A}\to Q$ by the induction principle of the propositional truncation. We have to provide a function $A\to Q$, which we have assumed already, and we have to show that
  \begin{equation*}
    \tr_{\lam{x}Q}(\alpha(x,y),u)=v.
  \end{equation*}
  for any $u,v:Q$ and any $x,y:\brck{A}$. However, we have such identifications by the assumption that $Q$ is a proposition, so the proof is complete.
\end{proof}

 One simple application of the universal property of the propositional truncation is that $\brck{\blank}$ acts on functions in a functorial way.

\begin{prp}
  There is a map\index{functorial action!of propositional truncations|textbf}\index{propositional truncation!functorial action|textbf}
  \begin{equation*}
    \brck{\blank}:(A\to B)\to (\brck{A}\to\brck{B})
  \end{equation*}
  for any two types $A$ and $B$, such that
  \begin{align*}
    \brck{\idfunc} & \htpy \idfunc \\
    \brck{g\circ f} & \htpy \brck{g}\circ\brck{f}.
  \end{align*}
\end{prp}

\begin{proof}
  For any $f:A\to B$, the map $\brck{f}:\brck{A}\to\brck{B}$ is defined to be the unique extension
  \begin{equation*}
    \begin{tikzcd}
      A \arrow[d,swap,"\eta"] \arrow[r,"f"] & B \arrow[d,"\eta"] \\
      \brck{A} \arrow[r,dashed,swap,"\brck{f}"] & \brck{B}.
    \end{tikzcd}
  \end{equation*}
  To see that $\brck{\blank}$ preserves identity maps and compositions, simply note that $\idfunc[\brck{A}]$ is an extension of $\idfunc[A]$, and that $\brck{g}\circ\brck{f}$ is an extension of $g\circ f$. Hence the homotopies are obtained by uniqueness.
\end{proof}
\index{universal property!of propositional truncations|)}
\index{propositional truncation!universal property|)}

\subsection{Logic in type theory}\label{sec:logic}
\index{logic|(}

In \cref{sec:modular-arithmetic} we interpreted logic in type theory via the Curry-Howard correspondence, which stipulates that disjunction ($\lor$) is interpreted by coproducts and the existential quantifier ($\exists$) is interpreted by $\Sigma$-types. However, when the existential quantifier is interpreted by $\Sigma$-types, then it is not possible to express certain concepts correctly, such as finiteness of a type or being in the image a map, and therefore we will add a second interpretation of logic in type theory, where logical propositions are interpreted by type theoretic propositions, i.e., the types of truncation level $-1$.

We have seen that the propositions are closed under cartesian products, implication, and dependent products indexed by arbitrary types. However, they are not closed under coproducts, and if $P$ is a family of propositions over a type $A$, then it is not necessarily the case that $\sm{x:A}P(x)$ is a proposition. We will therefore use propositional truncations to interpret disjunctions and existential quantifiers in type theory.

\begin{defn}
  Given two propositions $P$ and $Q$, we define their \define{disjunction}\index{disjunction|textbf}\index{P v Q@{$P\vee Q$}|see {disjunction}}
  \begin{equation*}
    P\vee Q \defeq \brck{P+Q}.
  \end{equation*}
\end{defn}

\begin{prp}
  Consider two propositions $P$ and $Q$. Then the disjunction $P\vee Q$ comes equipped with maps $i:P\to P\vee Q$ and $j:Q\to P\vee Q$. Moreover, the proposition $P\vee Q$ satisfies the universal property of the disjunction\index{universal property!of disjunction|textbf}\index{disjunction!universal property|textbf}: For any proposition $R$, we have
  \begin{equation*}
    (P\vee Q\to R)\leftrightarrow ((P\to R)\times (Q\to R)).
  \end{equation*}
\end{prp}

\begin{proof}
  The maps $i$ and $j$ are defined by
  \begin{align*}
    i & \defeq \eta\circ\inl  \\
    j & \defeq \eta\circ\inr.
  \end{align*}
  Now consider the following composition of maps, for an arbitrary proposition $R$:
  \begin{equation*}
    \begin{tikzcd}
      (P\vee Q\to R) \arrow[r,"\blank\circ\eta"] & (P+Q\to R) \arrow[r,"{h\,\mapsto\,(h\circ \inl,h\circ \inr)}"] &[3.6em] (P\to R)\times (Q\to R).
    \end{tikzcd}
  \end{equation*}
  The first map is an equivalence by the universal property of the propositional truncation, and the second map is an equivalence by the universal property of coproducts (\cref{ex:up-coproduct}).
\end{proof}

\begin{defn}
  Given a family $P$ of propositions over a type $A$, we define the \define{existential quantification}\index{existential quantification|textbf}\index{E (x:A) P(x)@{$\exists_{(x:A)}P(x)$}|see {existential quantification}}\index{E (x:A) P(x)@{$\exists_{(x:A)}P(x)$}|textbf}
  \begin{equation*}
    \exists_{(x:A)}P(x)\defeq \Brck{\sm{x:A}P(x)}.
  \end{equation*}
\end{defn}

\begin{prp}
  Consider a family $P$ of propositions over a type $A$. Then the existential quantification $\exists_{(x:A)}P(x)$ comes equipped with a dependent function
  \begin{equation*}
    \prd{a:A} \big(P(a)\to \exists_{(x:A)}P(x)\big).
  \end{equation*}
  Furthermore, the proposition $\exists_{(x:A)}P(x)$ satisfies the universal property of the existential quantification\index{universal property!of existential quantification|textbf}\index{existential quantification!universal property|textbf}: For any proposition $Q$, we have
  \begin{equation*}
    \Big(\Big(\exists_{(x:A)}P(x)\Big)\to Q\Big)\leftrightarrow\Big(\prd{x:A}P(x)\to Q\Big).
  \end{equation*}
\end{prp}

\begin{proof}
  The dependent function $\varepsilon : \prd{a:A} \big(P(a)\to \exists_{(x:A)}P(x)\big)$ is given by $\varepsilon(a,p):=\eta(a,p)$. Now consider the following composition of maps
  \begin{equation*}
    \begin{tikzcd}[column sep=small]
      \big(\big(\exists_{(x:A)}P(x)\big)\to Q\big) \arrow[r] &
      \big(\big(\sm{x:A}P(x)\big)\to Q\big) \arrow[r] &
      \big(\prd{x:A}P(x)\to Q\big).
    \end{tikzcd}
  \end{equation*}
  The first map in this composite is an equivalence by the universal property of the propositional truncation, and the second map is an equivalence by the universal property of $\Sigma$-types (\cref{thm:up-sigma}).
\end{proof}

In the following table we give an overview of the interpretation of the logical connectives using the propositions in type theory.

\begin{center}
  \begin{tabular}{ll}
    \toprule
    logical connective & interpretation in type theory\index{interpretation of logic in type theory}\index{logic!interpretation of logic in type theory} \\
    \midrule
    $\top$\index{T@{$\top$}|textbf} & $\unit$ \\
    $\bot$\index{T@{$\bot$}|textbf} & $\emptyt$ \\
    $P\Rightarrow Q$\index{implication|textbf} & $P\to Q$ \\
    $P\land Q$\index{conjunction|textbf} & $P\times Q$ \\
    $P\lor Q$\index{disjunction} & $\brck{P+Q}$ \\
    $P\Leftrightarrow Q$\index{bi-implication|textbf} & $P\leftrightarrow Q$ \\
    $\exists_{(x:A)}P(x)$\index{existential quantification} & $\Brck{\sm{x:A}P(x)}$ \\
    $\forall_{(x:A)}P(x)$\index{universal quantification|textbf}\index{A (x:A) P(x)@{$\forall_{(x:A)}P(x)$}|see {universal quantification}} & $\prd{x:A}P(x)$ \\
    \bottomrule
  \end{tabular}
\end{center}
\index{logic|)}

\subsection{Mapping propositional truncations into sets}

The universal property of the propositional truncation only applies when we want to define a map into a proposition. However, in some situations we might want to map the propositional truncation into a type that is not a proposition. Here we will see what we might do in such a case.

One strategy, if we want to define a map $\brck{A}\to X$, is to find a type family $P$ over $X$ such that the type $\sm{x:X}P(x)$ is a proposition. In that case, we may use the universal property of the propositional truncation to obtain a map $\brck{A}\to \sm{x:X}P(x)$ from a map $A\to \sm{x:X}P(x)$, and then we simply compose with the projection map.

\begin{eg}\label{eg:global-choice-decidable-subtype-N}
  Consider a \define{decidable subtype}\index{decidable subtype|textbf} $P$ of the natural numbers, i.e., a subtype $P:\N\to\prop_\UU$ such that each $P(n)$ is decidable. We claim that there is a function
  \begin{equation*}
    \Brck{\sm{x:\N}P(x)}\to\sm{x:\N}P(x).
  \end{equation*}
  Of course, we cannot directly use the universal property of the propositional truncation here. However, there is at most one \emph{minimal} natural number $x$ in $P$. In other words, we claim that the type
  \begin{equation}\label{eq:is-prop-minimal-element}
    \sm{x:\N}P(x)\times\islowerbound_P(x)\tag{\textasteriskcentered}
  \end{equation}
  is a proposition. To see this, note that the type $\islowerbound_P(x)$ is a proposition. By the assumption that each $P(x)$ is a proposition, it now follows that any two natural numbers $x,y:\N$ that are in $P$ and that are both lower bounds of $P$ are equal as elements in the type of \cref{eq:is-prop-minimal-element} if and only if they are equal as natural numbers. Furthermore, since both $x$ and $y$ are lower bounds of $P$, it follows that $x\leq y$ and $y\leq x$, so indeed $x=y$ holds.

  By the observation that the type in \cref{eq:is-prop-minimal-element} is a proposition, we may define a map
  \begin{equation*}
    \Brck{\sm{x:\N}P(x)}\to \sm{x:\N}P(x)\times\islowerbound_P(x)
  \end{equation*}
  by the universal property of the propositional truncation. A map
  \begin{equation*}
    \sm{x:\N}P(x)\to\sm{x:\N}P(x)\times\islowerbound_P(x)
  \end{equation*}
  was constructed in \cref{thm:well-ordering-principle-N} using the decidability of $P$.

  As a corollary of this observation, we observe that there is also a map
  \begin{equation*}
    \Brck{\sm{x:\Fin{k}}P(x)}\to\sm{x:\Fin{k}}P(x)
  \end{equation*}
  for any decidable subtype $P$ over $\Fin{k}$.
\end{eg}

\begin{rmk}\label{rmk:global-choice}
  The function of type
  \begin{equation*}
    \Brck{\sm{x:\N}P(x)}\to\sm{x:\N}P(x)
  \end{equation*}
  we constructed in \cref{eg:global-choice-decidable-subtype-N} for decidable subtypes of $\N$ is a rare case in which it is possible to obtain a function
  \begin{equation*}
    \brck{A}\to A.
  \end{equation*}
  We say that the type $A$ satisfies the \define{principle of global choice} if there is such a function $\brck{A}\to A$. Using the univalence axiom, we will see in \cref{cor:no-global-choice} that not every type satisfies the principle of global choice. 
\end{rmk}

More generally, we may wish to define a map $\brck{A}\to B$ where the type $B$ is a set. In this situation it is helpful to think of the propositional truncation of $A$ as the quotient of the type $A$ by the equivalence relation that relates every two elements of $A$ with each other. Propositional truncations can therefore also be characterized by the universal property of this quotient, which can be used to extend maps $f:A\to B$ to maps $\brck{A}\to B$ when the type $B$ is a set. The idea is that a map $f:A\to B$ into a set $B$ extends to a map $\brck{A}\to B$ if it satisfies $f(x)=f(y)$ for all $x,y:A$.

\begin{defn}\label{defn:weakly-constant}
  A map $f:A\to B$ is said to be \define{weakly constant}\index{weakly constant map|textbf}\index{constant map!weakly constant map|textbf} if it comes equipped with an element of type\index{is-weakly-constant(f)@{$\isweaklyconstant(f)$}|textbf}
  \begin{equation*}
    \isweaklyconstant(f) \defeq \prd{x,y:A}f(x)=f(y).
  \end{equation*}
\end{defn}

\begin{rmk}
  A constant map $A\to B$ is a map of the form $\const_b$. A map $f:A\to B$ is therefore constant if comes equipped with an element $b:B$ and a homotopy $f\htpy \const_b$. This is a stronger notion than the notion of weakly constant maps, which doesn't require there to be an element in $B$.

  One of the differences between constant maps and weakly constant maps manifests itself as follows: A type $A$ is contractible if and only if the identity map on $A$ is constant, while a type $A$ is a proposition if and only if the identity map on $A$ is weakly constant.
\end{rmk}

\begin{lem}
  Consider a commuting triangle
  \begin{equation*}
    \begin{tikzcd}
      A \arrow[d,swap,"\eta"] \arrow[dr,"f"] \\
      \brck{A} \arrow[r,swap,"g"] & B      
    \end{tikzcd}
  \end{equation*}
  where $B$ is an arbitrary type. Then the map $f$ is weakly constant.
\end{lem}

\begin{proof}
  Since $f$ is assumed to be homotopic to $g\circ \eta$, it suffices to show that $g\circ\eta$ is weakly constant. For any $x,y:A$, we have the identification $\alpha(x,y):\eta(x)=\eta(y)$ in $\brck{A}$. Using the action on paths of $g$, we obtain the identification
  \begin{equation*}
    \ap{g}{\alpha(x,y)}:g(\eta(x))=g(\eta(y))
  \end{equation*}
  in $B$.
\end{proof}

We now show, in a theorem due to Kraus \cite{Kraus}, that any weakly constant map $f:A\to B$ into a set $B$ extends uniquely to a map $\brck{A}\to B$. We therefore conclude that, in order to define a map $\brck{A}\to B$ into a set $B$ it suffices to define a map $f:A\to B$ and show that it is weakly constant.

\begin{thm}[Kraus]\label{ex:weakly-constant-map}
  Let $A$ be a type and let $B$ be a set. Then the map\index{universal property!of propositional truncations into sets|textbf}\index{propositional truncation!universal property into sets|textbf}
  \begin{equation*}
    (\brck{A}\to B)\to \sm{f:A\to B}\prd{x,y:A}f(x)=f(y)
  \end{equation*}
  given by $g\mapsto (g\circ\eta,\lam{x}\lam{y}\ap{g}{\alpha(x,y)})$ is an equivalence.
\end{thm}

\begin{proof}
  Consider a map $f:A\to B$ equipped with $H:\prd{x,y:A}f(x)=f(y)$. We first show that $f$ extends in at most one way to a map $\brck{A}\to B$. Let $g,h:\brck{A}\to B$ be two maps equipped with homotopies $f\htpy g\circ\eta$ and $f\htpy h\circ\eta$. In order to construct a homotopy $g\htpy h$, note that each identity type $g(x)=h(x)$ is a proposition by the assumption that $B$ is a set. We can therefore construct a homotopy $g\htpy h$ by the induction principle of propositional truncations. By the induction principle, it suffices to construct a homotopy $g\circ \eta\htpy h\circ\eta$, which we obtain from the homotopies $f\htpy g\circ\eta$ and $f\htpy h\circ\eta$.

  Since we've already proven uniqueness, it remains to construct an extension of the map $f$. We first claim that the type
  \begin{equation*}
    \sm{b:B}\Brck{\sm{x:A}f(x)=b}
  \end{equation*}
  is a proposition. To see this, consider two elements $b$ and $b'$ in this subtype of $B$. It suffices to show that $b=b'$. Since $B$ is assumed to be a set, the identity type $b=b'$ is a proposition. Therefore we may assume an element $x:A$ equipped with $p:f(x)=b$ and an element $x':A$ equipped with $p':f(x')=b'$. Using the assumption that $f$ is weakly constant, we obtain the identification
  \begin{equation*}
    \begin{tikzcd}
      b \arrow[r,equals,"p^{-1}"] & f(x) \arrow[r,equals,"{H(x,x')}"] & f(x') \arrow[r,equals,"{p'}"] & b'.
    \end{tikzcd}
  \end{equation*}
  Now we observe that the map $f:A\to B$ factors uniquely as follows
  \begin{equation*}
    \begin{tikzcd}[column sep=-1.5em]
      A \arrow[rr,dashed,"g"] \arrow[dr,swap,"f"] & & \sm{b:B}\Brck{\sm{x:A}f(x)=b} \arrow[dl,"\proj 1"] \\
      \phantom{\sm{b:B}\Brck{\sm{x:A}f(x)=b}} & B.
    \end{tikzcd}
  \end{equation*}
  Indeed, the map $g$ is given by $x\mapsto(f(x),\eta(x,\refl{}))$. Since the codomain of $g$ is a proposition, we obtain via the universal property of the propositional truncation of $A$ a unique map $h:\brck{A}\to\sm{b:B}\left\|\sm{x:A}f(x)=b\right\|$ equipped with a homotopy $g\htpy h\circ\eta$. Now we obtain the map $\proj 1\circ h:\brck{A}\to B$ equipped with the concatenated homotopy
  \begin{equation*}
    (\proj 1\circ h)\circ\eta \jdeq \proj 1\circ (h\circ\eta) \htpy \proj 1\circ g \htpy f.\qedhere
  \end{equation*}
\end{proof}

\begin{exercises}
  \exitem \label{ex:propositional-truncations-drill}Let $A$ be a type. Show that
  \begin{subexenum}
  \item $\brck{\brck{A}}\leftrightarrow\brck{A}$.
  \item $\brck{\isdecidable(A)}\leftrightarrow\isdecidable\brck{A}$.
  \item $\isdecidable(A)\to (\brck{A}\to A)$.
  \item $\neg\neg\brck{A}\leftrightarrow\neg\neg A$.
  \item $\brck{A}\vee\brck{B}\leftrightarrow \brck{A+B}$.
  \item $\exists_{(x:A)}\brck{B(x)}\leftrightarrow \brck{\sm{x:A}B(x)}$.
  \item $\neg\neg(\brck{A}\to A)$.
  \end{subexenum}
  \exitem Show that the \define{mere equality}\index{mere equality|textbf} relation given by $x,y\mapsto\brck{x=y}$ is an equivalence relation on any type.
  \exitem \label{ex:product-propositional-truncation}Consider two maps $f:A\to P$ and $g:B\to Q$ into propositions $P$ and $Q$. Show that if both $f$ and $g$ are propositional truncations then the map $f\times g : A\times B\to P\times Q$ is also a propositional truncation. Conclude that\index{propositional truncation!distributes over cartesian products}\index{distributivity!of propositional truncations over cartesian products}\index{cartesian product type!distributivity of propositional truncations over products}
  \begin{equation*}
    \brck{A\times B}\simeq \brck{A}\times\brck{B}. 
  \end{equation*}
  \exitem \label{ex:dup-trunc-prop}Consider a map $f:A\to P$ into a proposition $P$. We say that $f$ satisfies the \define{dependent universal property of the propositional truncation}\index{dependent universal property!of propositional truncations|textbf}\index{propositional truncation!dependent universal property|textbf} of $A$, if for any family $Q$ of propositions over $P$, the precomposition function
  \begin{equation*}
    \blank\circ f : \Big(\prd{p:P}Q(p)\Big)\to\Big(\prd{x:A}Q(f(x))\Big)
  \end{equation*}
  is an equivalence. Show that the following are equivalent:
  \begin{enumerate}
  \item\label{item:up-dup-trunc-Prop} The map $f$ is a propositional truncation.
  \item\label{item:dup-dup-trunc-Prop} The map $f$ satisfies the dependent universal property of the propositional truncation.
  \end{enumerate}
  \exitem Consider a map $f:A\to P$ into a proposition $P$.
  \begin{subexenum}
  \item Show that if there is a map $g:P\to A$, then $f$ is a propositional truncation. Conclude that for any type $A$ equipped with a point $a:A$, the constant map
    \begin{equation*}
      \const_\ttt: A\to\unit
    \end{equation*}
    is a propositional truncation of $A$.
  \item Show that if $A$ is a proposition, then $f$ is a propositional truncation if and only if $f$ is an equivalence. Conclude that if $A$ is a proposition, then the identity function $\idfunc:A\to A$ is a propositional truncation.
  \end{subexenum}
  \exitem Consider a type $A$ equipped with an element $d:\isdecidable(A)$.\index{decidable type}
  \begin{subexenum}
  \item Define a function $f:A\to \sm{x:A}\inl(x)=d$ and show that $f$ is a propositional truncation of $A$.
  \item Consider the function $\pi:\isdecidable(A)\to\Fin{2}$ defined by
    \begin{align*}
      \pi(\inl(x)) & := 1 \\
      \pi(\inr(x)) & := 0.
    \end{align*}
    Define a function $g:A\to (\pi(d)=1)$ and show that $g$ is a propositional truncation of $A$.
  \end{subexenum}
  \exitem Consider a family $B$ of $(k+1)$-truncated types over the propositional truncation $\brck{A}$ of a type $A$. Show that the map
  \begin{equation*}
    \Big(\prd{x:\brck{A}}B(x)\Big)\to\Big(\prd{x:A}B(x)\Big)
  \end{equation*}
  given by $f\mapsto f\circ\eta$ is a $k$-truncated map.
  \exitem Consider a universe $\UU$, let $P_1$ and $P_2$ be propositions in $\UU$, and furthermore, let $P$ be a family of propositions in $\UU$ over a type $A$ in $\UU$. Construct the following equivalences:
    \begin{align*}
      \top & \simeq \prd{Q:\prop_\UU}Q\to Q,\\
      \bot & \simeq \prd{Q:\prop_\UU}Q,\\
      \brck{A} & \simeq \prd{Q:\prop_\UU}(A\to Q)\to Q, \\
      P_1\lor P_2 & \simeq \prd{Q:\prop_\UU}(P_1\to Q) \to ((P_2\to Q)\to Q), \\
      P_1\land P_2 & \simeq \prd{Q:\prop_\UU}(P_1\to (P_2\to Q))\to Q, \\
      P_1\Rightarrow P_2 & \simeq \prd{Q:\prop_\UU}P_1\to ((P_2\to Q)\to Q), \\
      \neg A & \simeq \prd{Q:\prop_\UU} A\to Q, \\
      \exists_{(x:A)}P(x) & \simeq \prd{Q:\prop_\UU} \Big(\prd{x:A}P(x)\to Q\Big)\to Q, \\
      \forall_{(x:A)}P(x) & \simeq \prd{Q:\prop_\UU}\prd{x:A}(P(x)\to Q)\to Q\\
      \brck{a = x} & \simeq \prd{Q:A\to\prop_\UU} Q(a)\to Q(x).
    \end{align*}
    These are the \define{impredicative encodings}\index{impredicative encodings}\index{logic!impredicative encodings}\index{T@{$\top$}!impredicative encoding}\index{T@{$\bot$}!impredicative encoding}\index{propositional truncation!impredicative encoding}\index{disjunction!impredicative encoding}\index{conjunction!impredicative encoding}\index{implication!impredicative encoding}\index{negation!impredicative encoding}\index{existential quantification!impredicative encoding}\index{universal quantification!impredicative encoding}\index{mere equality!impredicative encoding} of the logical operators. \emph{Note:} It has the appearance that we could have defined $\brck{A}$ by its impredicative encoding. There is, however, a subtle issue if we take this as a definition: The map
    \begin{equation*}
      A\to\prd{Q:\prop_\UU}(A\to Q)\to Q
    \end{equation*}
    only satisfies the universal property of the propositional truncation with respect to propositions that are equivalent to propositions in $\UU$. 
  \exitem In this exercise we introduce the \define{interval}\index{interval|textbf} as a higher inductive type $\I$\index{I@{$\I$}|see{interval}}\index{I@{$\I$}|textbf}\index{higher inductive type!interval|textbf}, equipped with two point constructors and one path constructor
  \begin{align*}
    \source,\target & : \I \\
    \pathI & : \source=\target.
  \end{align*}
  The induction principle of $\I$ asserts that for any type family $P$ over $\I$, if we have
  \begin{align*}
    u & : P(\source) \\
    v & : P(\target) \\
    p & : \tr_P(\pathI,u)=v,
  \end{align*}
  then there is a section $f:\prd{x:\I}P(x)$ equipped with identifications
  \begin{align*}
    \alpha & : f(\source) = u \\
    \beta & : f(\target) = v 
  \end{align*}
  and an identification $\gamma$ witnessing that the square
  \begin{equation*}
    \begin{tikzcd}[column sep=6em]
      \tr_P(\pathI,f(\source)) \arrow[r,equals,"\ap{\tr_P(\pathI)}{\alpha}"] \arrow[d,equals,swap,"\apd{f}{\pathI}"] & \tr_P(\pathI,u) \arrow[d,equals,"p"] \\
      f(\target) \arrow[r,equals,swap,"\beta"] & v
    \end{tikzcd}
  \end{equation*}
  commutes. Note that the constructors of $\I$ induce a map
    \begin{equation*}
      \varepsilon: \Big(\prd{x:\I}P(x)\Big)\to \Big(\sm{u:P(\source)}\sm{v:P(\target)}\tr_P(\pathI,u)=v\Big).
    \end{equation*}
    given by $f\mapsto (f(\source),f(\target),\apd{f}{\pathI})$.
  \begin{subexenum}
  \item Characterize the identity types of the codomain of the map $\varepsilon$ in the following way: Construct an equivalence from the type $(u,v,q)=(u',v',q')$ to the type
    \begin{equation*}
      \sm{\alpha:u=u'}\sm{\beta:v=v'} \ct{q}{\beta}=\ct{\ap{\tr_P(\pathI)}{\alpha}}{q'},
    \end{equation*}
    for any $(u,v,q)$ and $(u',v',q')$ in the codomain of $\varepsilon$.
  \item Prove the dependent universal property of $\I$, i.e., show that the map $\varepsilon$ is an equivalence. 
  \item Show that $\I$ is contractible.
  \end{subexenum}
\end{exercises}
\index{propositional truncation|)}


\section{Image factorizations}\label{chap:image}

The image of a map $f:A\to X$ can be thought of as the least subtype of $X$ that contains all the values of $f$. More precisely, the image of $f$ is an embedding $i:\im(f)\hookrightarrow X$ that fits in a commuting triangle
\begin{equation*}
  \begin{tikzcd}[column sep=tiny]
    A \arrow[rr,"q"] \arrow[dr,swap,"f"] & & \im(f) \arrow[dl,hook,"i"] \\
    \phantom{\im(f)} & X
  \end{tikzcd}
\end{equation*}
and satisfies the \emph{universal property} of the image of $f$, which states that if a subtype $B\hookrightarrow X$ contains all the values of $f$, then it contains the image of $f$.

\subsection{The image of a map}\label{sec:image-construction}
 
\subsubsection*{The universal property of the image}
\index{universal property!of the image of a map|(}
\index{image of a map!universal property|(}

Recall from \cref{ex:triangle_fib} that we made the following definition:

\begin{defn}
  Let $f:A\to X$ and $g:B\to X$ be maps. A \define{morphism from $f$ to $g$ over $X$}\index{morphism from f to g over X@{morphism from $f$ to $g$ over $X$}|textbf} consists of a map $h:A\to B$ equipped with a homotopy $H:f\htpy g\circ h$ witnessing that the triangle
\begin{equation*}
\begin{tikzcd}[column sep=tiny]
A \arrow[rr,"h"] \arrow[dr,swap,"f"] & & B \arrow[dl,"g"] \\
& X
\end{tikzcd}
\end{equation*}
commutes. Thus, we define the type\index{hom X (f,g)@{$\homslice_X(f,g)$}|see {morphism from $f$ to $g$ over $X$}}
\begin{equation*}
\homslice_X(f,g)\defeq\sm{h:A\to B}f\htpy g\circ h.
\end{equation*}
Composition of morphisms over $X$ is defined by
\begin{equation*}
  (k,K)\circ (h,H) \defeq (k\circ h,\ct{H}{(K\cdot h)}).
\end{equation*}
\end{defn}

\begin{defn}
Consider a commuting triangle
\begin{equation*}
\begin{tikzcd}[column sep=tiny]
A \arrow[rr,"q"] \arrow[dr,swap,"f"] & & I \arrow[dl,"i"] \\
& X
\end{tikzcd}
\end{equation*}
with $H:f\htpy i\circ q$, where $i$ is an embedding\index{embedding}.
We say that $i$ satisfies the \define{universal property of the image of $f$}\index{universal property!of the image of a map|textbf} if the precomposition function
\begin{equation*}
\blank\circ(q,H) : \homslice_X(i,m)\to\homslice_X(f,m)
\end{equation*}
is an equivalence for every embedding $m:B\hookrightarrow X$. 
\end{defn}

\begin{lem}
For any $f:A\to X$ and any embedding\index{embedding} $m:B\to X$, the type $\homslice_X(f,m)$ is a proposition.\index{hom X (f,g)@{$\homslice_X(f,g)$}!is a proposition}
\end{lem}

\begin{proof}
  Recall from \cref{ex:triangle_fib} that the type $\homslice_X(f,m)$ is equivalent to the type
  \begin{equation*}
    \prd{a:A}\fib{m}{f(a)}.
  \end{equation*}
  Furthermore, recall from \cref{thm:embedding} that a map is an embedding if and only if its fibers are propositions.
  Thus we see that the type $\prd{a:A}\fib{m}{f(a)}$ is a product of propositions, hence it is a proposition by \cref{thm:trunc_pi}.
\end{proof}

\begin{prp}\label{prp:simplifly-universal-property-image}
  Consider a commuting triangle
  \begin{equation*}
    \begin{tikzcd}[column sep=tiny]
      A \arrow[rr,"q"] \arrow[dr,swap,"f"] & & I \arrow[dl,"i"] \\
      & X
\end{tikzcd}
  \end{equation*}
  with $H:f\htpy i\circ q$, where $i$ is an embedding. Then the following are equivalent:
  \begin{enumerate}
  \item The embedding $i$ satisfies the universal property of the image inclusion of $f$.
  \item For every embedding $m:B\to X$ there is a map
    \begin{equation*}
      \homslice_X(f,m)\to\homslice_X(i,m).
    \end{equation*}
  \end{enumerate}
\end{prp}

\begin{proof}
Since $\homslice_X(f,m)$ is a proposition for every embedding $m:B\to X$, the claim follows immediately by the observation made in \cref{ex:prop_equiv}.
\end{proof}
\index{universal property!of the image of a map|)}
\index{image of a map!universal property|)}

\subsubsection*{The existence of the image}
\index{image of a map!existence|(}

The image of a map $f:A\to X$ can be defined using the propositional truncation.

\begin{defn}\label{defn:im}
For any map $f:A\to X$ we define the \define{image}\index{image of a map|textbf} of $f$ to be the type\index{im f@{$\im(f)$}|see {image of a map}}\index{im f@{$\im(f)$}|textbf}
\begin{equation*}
\im(f) \defeq \sm{x:X}\brck{\fib{f}{x}}.
\end{equation*}
Furthermore, we define
\begin{enumerate}
\item the \define{image inclusion}\index{image inclusion|textbf}
  \begin{equation*}
    i_f:\im(f)\to X
  \end{equation*}
  to be the projection $\proj 1$,
\item the map
  \begin{equation*}
    q_f:A\to\im(f)
  \end{equation*}
  to be the map given by $q_f(x)\defeq(f(x),\eta(x,\refl{f(x)}))$, and
\item the homotopy $I_f:f\htpy i_f\circ q_f$ witnessing that the triangle
  \begin{equation*}
    \begin{tikzcd}[column sep=tiny]
      A \arrow[rr,"q_f"] \arrow[dr,swap,"f"] & & \im(f) \arrow[dl,"i_f"] \\
      \phantom{\im(f)} & X
    \end{tikzcd}
  \end{equation*}
  commutes, to be given by $I_f(x)\defeq\refl{f(x)}$.
\end{enumerate}
\end{defn}

\begin{prp}
  The image inclusion $i_f:\im(f)\to X$ of any map $f:A\to X$ is an embedding.\index{image inclusion!is an embedding}\index{is an embedding!image inclusion}
\end{prp}

\begin{proof}
  The claim follows directly by \cref{cor:pr1-embedding}, because the type $\brck{\fib{f}{x}}$ is a proposition for each $x:X$.
\end{proof}

\begin{thm}\label{thm:im}
  The image inclusion $i_f:\im(f)\to X$ of any map $f:A\to X$ satisfies the universal property of the image inclusion of $f$.
\end{thm}

\begin{proof}
  Consider an embedding $m:B\hookrightarrow X$. Note that we have a commuting square
  \begin{equation*}
    \begin{tikzcd}[column sep=6em]
      \homslice_X(i_f,m) \arrow[d] \arrow[r] & \homslice_X(f,m) \arrow[d] \\
      \Big(\prd{x:X}\fib{i_f}{x}\to\fib{m}{x}\Big) \arrow[r,swap,"h\mapsto{\lam{x}h_x\circ\varphi_x}"] & \Big(\prd{x:X}\fib{f}{x}\to\fib{m}{x}\Big)
    \end{tikzcd}
  \end{equation*}
  in which all four types are propositions, and the vertical maps are equivalences. Therefore it suffices to construct a map
  \begin{equation*}
    \Big(\prd{x:X}\fib{f}{x}\to\fib{m}{x}\Big)\to\Big(\prd{x:X}\fib{i_f}{x}\to\fib{m}{x}\Big)
  \end{equation*}
  The fiber $\fib{i_f}{x}$ is equivalent to the propositional truncation $\brck{\fib{f}{x}}$ and the type $\fib{m}{x}$ is a proposition by the assumption that $m$ is an embedding. Therefore we obtain the desired map by the universal property of the propositional truncation.
\end{proof}
\index{image of a map!existence|)}

\subsubsection*{The uniqueness of the image}
\index{image of a map!uniqueness|(}

We will now show that the universal property of the image implies that the image is determined uniquely up to equivalence.

\begin{thm}\label{thm:uniqueness-image}
  Let $f$ be a map, and consider two commuting triangles
  \begin{equation*}
    \begin{tikzcd}[column sep=tiny]
      A \arrow[dr,swap,"f"] \arrow[rr,"q"] & & B \arrow[dl,"i"] &[2em] A \arrow[dr,swap,"f"] \arrow[rr,"{q'}"] & & B' \arrow[dl,"{i'}"] \\
      \phantom{B'} & X & \phantom{B'} & \phantom{B'} & X
    \end{tikzcd}
  \end{equation*}
  with $I:f\htpy i\circ q$ and $I':f\htpy i'\circ q'$, in which $i$ and $i'$ are assumed to be embeddings. Then, if any two of the following three properties hold, so does the third:
  \begin{enumerate}
  \item The embedding $i$ satisfies the universal property of the image inclusion of $f$.
  \item The embedding $i'$ satisfies the universal property of the image inclusion of $f$.
  \item The type of equivalences $e:B\simeq B'$ equipped with a homotopy witnessing that the triangle
    \begin{equation*}
      \begin{tikzcd}
        B \arrow[dr,swap,"i"] \arrow[rr,"e"] & & B' \arrow[dl,"{i'}"] \\
        \phantom{B'} & X
      \end{tikzcd}
    \end{equation*}
    commutes is contractible.
  \end{enumerate}
\end{thm}

\begin{proof}
  First, we show that if (i) and (ii) hold, then (iii) holds. Note that the type $\homslice_X(i,i')$ is a proposition, since $i'$ is assumed to be an embedding. Therefore it suffices to show that the unique map $h:B\to B'$ such that the triangle
  \begin{equation*}
    \begin{tikzcd}[column sep=tiny]
      B \arrow[dr,swap,"i"] \arrow[rr,"h"] & & B' \arrow[dl,"{i'}"] \\
      & X
    \end{tikzcd}
  \end{equation*}
  commutes, is an equivalence. To see this, note that by \cref{ex:triangle_fib} it suffices to show that the action on fibers
  \begin{equation*}
    \fib{i}{x}\to\fib{i'}{x}
  \end{equation*}
  is an equivalence for each $x:X$. This follows from the universal property of $i'$, since we similarly obtain a family of maps
  \begin{equation*}
    \fib{i'}{x}\to\fib{i}{x}
  \end{equation*}
  indexed by $x:X$, and the types $\fib{i}{x}$ and $\fib{i'}{x}$ are propositions by the assumptions that $i$ and $i'$ are embeddings.
  
  Now we will show that (iii) implies that (i) holds if and only if (ii) holds. We will assume a morphism $(e,H):\homslice_X(i,i')$ such that the map $e$ is an equivalence. Furthermore, consider an embedding $m:C\to X$. Then the fact that (i) holds if and only if (ii) holds follows from the equivalence
  \begin{equation*}
    \big(\homslice_X(f,m)\to\homslice_X(i,m)\big)\simeq\big(\homslice_X(f,m)\to\homslice_X(i',m)\big).\qedhere
  \end{equation*}
\end{proof}
\index{image of a map!uniqueness|)}

\subsection{Surjective maps}\label{subsec:surjective}

A map $f:A\to B$ is surjective if for every $b:B$ there is an \emph{unspecified} element $a:A$ that maps to $b$. We define this property using the propositional truncation.

\begin{defn}
A map $f:A\to B$ is said to be \define{surjective}\index{surjective map|textbf} if there is an element of type\index{is-surj f@{$\issurj(f)$}|see {surjective map}}
\begin{equation*}
\issurj(f)\defeq \prd{b:B}\brck{\fib{f}{b}}.
\end{equation*}
\end{defn}

\begin{eg}
  Any equivalence is a surjective map, since its fibers are contractible. More generally, any map that has a section is surjective. Those are sometimes called \define{split epimorphisms}. Note that having a section is stronger than surjectivity, since in general we don't have a function $\brck{\fib{f}{b}}\to\fib{f}{b}$.
\end{eg}

In \cref{ex:dup-trunc-prop} we showed the dependent universal property of the propositional truncation: a map $f:A\to B$ into a proposition $B$ satisfies the universal property of the propositional truncation if and only if for every family of propositions $P$ over $B$, the precomposition map
\begin{equation*}
  \blank\circ f : \Big(\prd{b:B}P(b)\Big)\to\Big(\prd{a:A}P(f(a))\Big)
\end{equation*}
is an equivalence. In the following proposition we show that, if we omit the condition that $B$ is a proposition, then $f$ satisfies this dependent universal property if and only if $f$ is surjective.

\begin{prp}\label{prp:surjective}
  Consider a map $f:A\to B$. Then the following are equivalent:
  \begin{enumerate}
  \item \label{prp-item:surjective}The map $f:A\to B$ is surjective.
  \item \label{prp-item:is-equiv-precomp-surjective}The map $f:A\to B$ satisfies the \define{dependent universal property of a surjective map}\index{dependent universal property!of surjective maps|textbf}\index{surjective map!dependent universal property|textbf}: For any family $P$ of propositions over $B$, the precomposition map
    \begin{equation*}
      \blank\circ f : \Big(\prd{y:B}P(y)\Big)\to\Big(\prd{x:A}P(f(x))\Big)
    \end{equation*}
    is an equivalence. In other words, any subtype of $B$ that contains all the elements of the form $f(x)$ contains all the elements of $B$.
  \item \label{prp-item:is-trunc-map-precomp-surjective}For any $k\geq-2$, and for any family $P$ of $(k+1)$-truncated types over $B$, the precomposition map
    \begin{equation*}
      \blank\circ f : \Big(\prd{y:B}P(y)\Big)\to\Big(\prd{x:A}P(f(x))\Big)
    \end{equation*}
    is a $k$-truncated map.
  \end{enumerate}
\end{prp}

\begin{proof}
  To prove that \ref{prp-item:surjective} implies \ref{prp-item:is-equiv-precomp-surjective}, suppose first that $f$ is surjective, and consider the commuting square
  \begin{equation*}
    \begin{tikzcd}[column sep=4.4em]
      \Big(\prd{y:B}P(y)\Big) \arrow[r,"\blank\circ f"] \arrow[d,swap,"h\mapsto\lam{y}\const_{h(y)}"] & \Big(\prd{x:A}P(f(x))\Big)  \\
      \Big(\prd{y:B}\brck{\fib{f}{y}}\to P(y)\Big) \arrow[r,swap,"h\mapsto h(\blank)\circ\eta"] & \Big(\prd{y:B}\fib{f}{y}\to P(y)\Big) \arrow[u,swap,"{h\mapsto\lam{x}h(f(x),(x,\refl{}))}"]
    \end{tikzcd}
  \end{equation*}
  In this square, the bottom map is an equivalence by \cref{ex:equiv-pi} and by the universal property of the propositional truncation of $\fib{f}{y}$. The map on the right is an equivalence by \cref{ex:pi-fib}. Furthermore, the map on the left is an equivalence by \cref{ex:equiv-pi,ex:up-unit}, because the type $\brck{\fib{f}{y}}$ is contractible by the assumption that $f$ is surjective. Therefore it follows that the top map is an equivalence, which completes the proof that \ref{prp-item:surjective} implies \ref{prp-item:is-equiv-precomp-surjective}.

  The proof that \ref{prp-item:is-equiv-precomp-surjective} implies \ref{prp-item:is-trunc-map-precomp-surjective} is by induction on $k$. The base case holds by assumption. For the inductive step, it suffices by \cref{thm:trunc_ap} to show that $\apfunc{\blank\circ f}$ is $k$-truncated for any $g,h:\prd{y:B}P(y)$. Notice that we have a commuting square
  \begin{equation*}
    \begin{tikzcd}[column sep=large]
      (g=h) \arrow[r,"\apfunc{\blank\circ f}"] \arrow[d,swap,"\htpyeq"] & (g\circ f = h\circ f) \arrow[d,"\htpyeq"] \\
      \prd{y:B}g(y)=h(y) \arrow[r,swap,"\blank\circ f"] & \prd{x:A}g(f(x))=h(f(x))
    \end{tikzcd}
  \end{equation*}
  The vertical maps on the left and right are equivalences by function extensionality, and the bottom map is $k$-truncated by the inductive hypothesis. This implies that $\apfunc{\blank\circ f}$ is $k$-truncated.

  To prove that \ref{prp-item:is-trunc-map-precomp-surjective} implies \ref{prp-item:surjective}, note that the assumption in \ref{prp-item:is-trunc-map-precomp-surjective} implies that the precomposition function
  \begin{equation*}
    \blank\circ f : \Big(\prd{y:B}\brck{\fib{f}{y}}\Big)\to\Big(\prd{x:A}\brck{\fib{f}{f(x)}}\Big)
  \end{equation*}
  is an equivalence. Hence it suffices to construct an element of type $\brck{\fib{f}{f(x)}}$ for each $x:A$. This is easy, because we have
  \begin{equation*}
    \eta(x,\refl{f(x)}):\brck{\fib{f}{f(x)}}.\qedhere
  \end{equation*}
\end{proof}

As a corollary we obtain that any surjective map into a proposition satisfies the universal property of the propositional truncation.

\begin{cor}
  For any map $f:A\to P$ into a proposition $P$, the following are equivalent:
  \begin{enumerate}
  \item The map $f$ satisfies the universal property of the propositional truncation of $A$.
  \item The map $f$ is surjective.
  \end{enumerate}
\end{cor}

Using the characterization of surjective maps of \cref{prp:surjective}, we can also give a new characterization of the image of a map. 

\begin{thm}\label{thm:surjective}
Consider a commuting triangle
\begin{equation*}
\begin{tikzcd}[column sep=tiny]
A \arrow[rr,"q"] \arrow[dr,swap,"f"] & & B \arrow[dl,"m"] \\
& X
\end{tikzcd}
\end{equation*}
in which $m$ is an embedding. Then the following are equivalent:
\begin{enumerate}
\item The embedding $m$ satisfies the universal property of the image inclusion of $f$.\index{image of a map!universal property}\index{surjective map!universal property of the image of a map}
\item The map $q$ is surjective.
\end{enumerate}
\end{thm}

\begin{proof}
  First assume that $m$ satisfies the universal property of the image inclusion of $f$, and consider the composite function
  \begin{equation*}
    \begin{tikzcd}
      \Big(\sm{y:B}\brck{\fib{q}{y}}\Big) \arrow[r,"\proj 1"] & B \arrow[r,"m"] & X.
    \end{tikzcd}
  \end{equation*}
  Note that $m\circ\proj 1$ is a composition of embeddings, so it is an embedding. By the universal property of $m$ there is a unique map $h$ for which the triangle
  \begin{equation*}
    \begin{tikzcd}[column sep=0]
      B \arrow[dr,swap,"m"] \arrow[rr,dashed,"h"] & & \sm{y:B}\brck{\fib{q}{y}} \arrow[dl,"m\circ\proj 1"] \\
      \phantom{\sm{y:B}\brck{\fib{q}{y}}} & X
    \end{tikzcd}
  \end{equation*}
  commutes. Now note that $\proj 1\circ h$ is a map such that $m\circ (\proj 1\circ h)\htpy m$. The identity function is another map for which we have $m\circ\idfunc\htpy m$, so it follows by uniqueness that $\proj 1\circ h\htpy \idfunc$. In other words, the map $h$ is a section of the projection map. Therefore we obtain by \cref{ex:pi_sec} a dependent function
  \begin{equation*}
    \prd{b:B}\brck{\fib{q}{b}},
  \end{equation*}
  showing that $q$ is surjective.

  For the converse, suppose that $q$ is surjective. To prove that $m$ satisfies the universal property of the image factorization of $f$, it suffices to construct a map
  \begin{equation*}
    \homslice_X(f,m')\to\homslice_X(m,m'),
  \end{equation*}
  for any embedding $m':B'\to X$. To see that there is such an equivalence, we make the following calculation
  \begin{align*}
    \homslice_X(m,m') &  \simeq \prd{b:B}\fib{m'}{m(b)} \tag{By \cref{ex:triangle_fib}}\\
                         & \simeq \prd{a:A}\fib{m'}{m(q(a))} \tag{By \cref{prp:surjective}}\\
                         & \simeq \prd{a:A}\fib{m'}{f(a)} \tag{By $f\htpy m\circ q$}\\
                         & \simeq \homslice_X(f,m').\tag{By \cref{ex:triangle_fib}}
  \end{align*}
\end{proof}

\begin{cor}
  Every map factors uniquely as a surjective map followed by an embedding.\index{surjective map!factorization}\index{embedding!factorization}
\end{cor}

\begin{proof}
  Consider a map $f:A\to X$, and two factorizations
  \begin{equation*}
    \begin{tikzcd}[column sep=tiny]
      A \arrow[rr,"q"] \arrow[dr,swap,"f"] & & B \arrow[dl,"i"] &[3em] A \arrow[rr,"{q'}"] \arrow[dr,swap,"f"] & & B' \arrow[dl,"{i'}"] \\
      & X & & & X
    \end{tikzcd}
  \end{equation*}
  of $f$ where $m$ and $m'$ are embeddings, and $q$ and $q'$ are surjective. Then both $m$ and $m'$ satisfy the universal property of the image factorization of $f$ by \cref{thm:surjective}. Now it follows by \cref{thm:uniqueness-image} that the type of $(e,H):\homslice_X(i,i')$ in which $e$ is an equivalence, equipped with an identification
  \begin{equation*}
    (e,H)\circ(q,I)=(q',I')
  \end{equation*}
  in $\homslice_X(f,i')$, is contractible.
\end{proof}

\subsection{Cantor's diagonal argument}
\index{Cantor's diagonal argument|(}

Now that we have introduced surjective maps, we are in position to give Cantor's famous diagonal argument, which he used to show that there are infinite sets of different cardinality. The diagonal argument gives a proof that there is no surjective map from $X$ to its power set $\mathcal{P}(X)$. The power set of a type $X$ is of course defined with respect to a universe $\UU$, as the type of families of propositions in $\UU$ indexed by $X$.

\begin{defn}
  Consider a type $X$, and a universe $\UU$. We define the \define{$\UU$-power set}\index{power set|textbf} of $X$ to be\index{P U X@{$\mathcal{P}_\UU(X)$}|see {power set}}
  \begin{equation*}
    \mathcal{P}_{\mathcal{U}}(X)\defeq X\to\prop_\UU.
  \end{equation*}
\end{defn}

\begin{thm}
  For any type $X$ and any universe $\UU$, there is no surjective function
  \begin{equation*}
    f : X \to \mathcal{P}_{\mathcal{U}}(X)
  \end{equation*}
\end{thm}

\begin{proof}
  Consider a function $f:X\to (X\to \prop_\UU)$, and suppose that $f$ is surjective. Following Cantor's diagonalization argument, we define the subset $P:X\to\prop_\UU$ by
  \begin{equation*}
    P(x)\defeq \neg(f(x,x)).
  \end{equation*}
  Our goal is to reach a contradiction and $f$ is assumed to be surjective. Therefore, it suffices to show that
  \begin{equation*}
    \Brck{\sm{x:X}f(x)=P}\to\emptyt.
  \end{equation*}
  The empty type is a proposition, so by the universal property of the propositional truncation it is equivalent to show that
  \begin{equation*}
    \Big(\sm{x:X}f(x)=P\Big)\to\emptyt.
  \end{equation*}
  Consider an element $x:X$ equipped with an identification $f(x)=P$. Our goal is to construct an element of the empty type, i.e, to reach a contradiction. By the identification $f(x)=P$ it follows that 
  \begin{equation*}
    f(x,y)\leftrightarrow P(y)
  \end{equation*}
  for all $y:X$. In particular, it follows that $f(x,x)\leftrightarrow P(x)$. However, since $P(x)$ is defined as $\neg(f(x,x))$, we obtain that $f(x,x)\leftrightarrow\neg(f(x,x))$. By \cref{ex:no-fixed-points-neg} this gives us the desired contradiction.
\end{proof}
\index{Cantor's diagonal argument|)}

\begin{exercises}
  \exitem Consider a commuting triangle
  \begin{equation*}
    \begin{tikzcd}[column sep=tiny]
      A \arrow[dr,swap,"f"] \arrow[rr,"h"] & & B \arrow[dl,"g"] \\
      & X
    \end{tikzcd}
  \end{equation*}
  where $g$ is an embedding.
  \begin{subexenum}
  \item Show that if there is a morphism
    \begin{equation*}
      \begin{tikzcd}[column sep=tiny]
        B \arrow[dr,swap,"g"] \arrow[rr,"k"] & & A \arrow[dl,"f"] \\
        & X
      \end{tikzcd}
    \end{equation*}
    over $X$, then $g$ satisfies the universal property of the image of $f$.
  \item Show that if $f$ is an embedding, then $g$ satisfies the universal property of $f$ if and only if $h$ is an equivalence.
  \end{subexenum}
  \exitem
  \begin{subexenum}
  \item Show that for any proposition $P$, the constant map\index{constant map!is an embedding}
    \begin{equation*}
      \const_\ttt : P \to \unit
    \end{equation*}
    is an embedding. Use this fact to construct an equivalence
    \begin{equation*}
      \Big(\sm{A:\UU}A\hookrightarrow\unit\Big)\simeq\prop_\UU.
    \end{equation*}
  \item Consider a map $f:A\to P$ into a proposition $P$. Show that the following are equivalent:
    \begin{enumerate}
    \item The map $f$ is a propositional truncation of $A$.\index{propositional truncation!universal property of the image of A arrow 1@{universal property of the image of $A\to \unit$}}
    \item The constant map $P\to\unit$ satisfies the universal property of the image of the constant map $A\to\unit$.
    \end{enumerate}
  \end{subexenum}
  \exitem \label{ex:is-equiv-is-emb-is-surjective}Consider a map $f:A\to B$. Show that the following are equivalent:
  \begin{enumerate}
  \item $f$ is an equivalence.\index{is an equivalence!is surjective and an embedding}
  \item $f$ is both surjective and an embedding.
  \end{enumerate}
  \exitem Consider a commuting triangle
  \begin{equation*}
    \begin{tikzcd}[column sep=tiny]
      A \arrow[rr,"h"] \arrow[dr,swap,"f"] & & B \arrow[dl,"g"] \\
      & X
    \end{tikzcd}
  \end{equation*}
  with $H:f\htpy g\circ h$.
  \begin{subexenum}
  \item Show that if $f$ is surjective, then $g$ is surjective.
  \item Show that if both $g$ and $h$ are surjective, then $f$ is surjective.
  \item As a converse to \cref{ex:is-trunc-comp}, show that if $f$ and $h$ are $k$-truncated, then $g$ is also $k$-truncated.
  \end{subexenum}
  \exitem Prove \define{Lawvere's fixed point theorem}\index{Lawvere's fixed point theorem}: For any two types $A$ and $B$, if there is a surjective map $f:A\to B^A$, then for any $h:B\to B$ there exists an $x:B$ such that $h(x)=x$, i.e., show that
  \begin{equation*}
    \Big(\exists_{(f:A\to(A\to B))}\issurj(f)\Big)\to\Big(\forall_{(h:B\to B)}\exists_{(b:B)}h(b)=b\Big).
  \end{equation*}
\end{exercises}

\section{Finite types}\label{chap:finite}

\subsection{Counting in type theory}
\index{counting|(}
When someone counts the elements of a finite set $A$, they go through the elements of $A$ one by one, at each stage keeping track of how many elements have been counted so far. This process results in the number $|A|$ of elements of the set $A$, and moreover it gives a bijection from the standard finite set with $|A|$ elements. In other words, to count the elements of $A$ is to give an equivalence from one of the standard finite sets to the set $A$. We turn this into a definition.

\begin{defn}
  For each type $A$, we define the type\index{count(A)@{$\cnt(A)$}|textbf}
  \begin{equation*}
    \cnt(A)\defeq\sm{k:\N}(\Fin{k}\simeq A).
  \end{equation*}
  The elements of $\cnt(A)$ are called \define{countings}\index{countings of a type|textbf} of $A$. When we have $(k,e):\cnt(A)$, we also say that $A$ \define{has $k$ elements}\index{has k elements@{has $k$ elements}|textbf}.
\end{defn}

Note that the type $\cnt(A)$ is often not a proposition. For instance, different equivalences of type $\Fin{k}\simeq\Fin{k}$ induce different elements of type $\cnt(\Fin{k})$.

\begin{eg}
  It follows immediately from the definition of countings that every standard finite type can be counted in a canonical way: For any $k:\N$ we have $(k,\idfunc) : \cnt(\Fin{k})$. It also follows immediately from the definition of countings that types equipped with countings are closed under equivalences.
\end{eg}

\begin{eg}
  Suppose $A$ comes equipped with a counting $(k,e):\cnt(A)$. Then $k=0$ if and only if $A$ is empty. Indeed, the inverse of $e$ is a map $e^{-1}:A\to\emptyt$. Conversely, if we have $f:\isempty(A)$, then the map $f:A\to\emptyt$ is automatically an equivalence. This shows that $\Fin{k}\simeq\emptyt$, and a short argument by induction on $k$ yields that $k=0$. 
\end{eg}

\begin{eg}
  A type $A$ has one element if and only if it is contractible. Indeed, the type $\Fin{1}$ is contractible, so it follows from the 3-for-2 property of contractible types (\cref{ex:contr_retr}) that there is an equivalence $\Fin{1}\simeq A$ if and only if $A$ is contractible.   
\end{eg}

\begin{eg}\label{rmk:count-decidable-proposition}
  A proposition $P$ comes equipped with a counting if and only if it is decidable. To see this, note that for any type $X$, if we have $(k,e):\cnt(X)$, then it follows that $X$ is decidable. This is shown by induction on $k$. In the case where $k=0$, it follows that $X$ is empty, and hence that $X$ is decidable. In the case where $k$ is a successor, the bijection $e:\Fin{k}\simeq X$ gives us the element $e(\ttt):X$. Again we conclude that $X$ is decidable.

  Conversely, if $P$ is decidable, then we can construct a counting of $P$ by case analysis on $d:P+\neg P$. If $P$ holds, then it is contractible and hence equivalent to $\Fin{1}$. If $\neg P$ holds, then $P$ is equivalent to $\Fin{0}$.  
\end{eg}

\begin{rmk}\label{rmk:has-decidable-equality-count}
  We also note that any type $A$ equipped with a counting $e:\Fin{k}\simeq A$ has decidable equality.\index{has decidable equality!type equipped with a counting} This follows from \cref{prp:has-decidable-equality-Fin}, where we showed that $\Fin{k}$ has decidable equality, for any $k:\N$.
\end{rmk}

\begin{thm}\label{thm:count}
  We make the following claims about countings:
  \begin{enumerate}
  \item\label{item:count-coprod} Consider two types $A$ and $B$. The following are equivalent:
    \begin{enumerate}
    \item Both $A$ and $B$ come equipped with a counting.
    \item The coproduct $A+B$ comes equipped with a counting.
    \end{enumerate}
  \item\label{item:count-Sigma} Consider a type family $B$ indexed by a type $A$. Consider the following three conditions:
    \begin{enumerate}
    \item \label{item:count-Sigma-count-base}The type $A$ comes equipped with a counting.
    \item \label{item:count-Sigma-count-fibers}The type $B(x)$ comes equipped with a counting, for each $x:A$.
    \item \label{item:count-Sigma-count-total}The type $\sm{x:A}B(x)$ comes equipped with a counting.
    \end{enumerate}
    If (a) holds, then (b) holds if and only if (c) holds. Furthermore, if both (b) and (c) hold and if $B$ comes equipped with a section $f:\prd{x:A}B(x)$, then (a) holds.

    Consequently, if $P$ is a subtype of a type $A$ equipped with a counting, then we have
    \begin{equation*}
      \cnt\Big(\sm{x:A}P(x)\Big)\leftrightarrow \prd{x:A}\isdecidable(P(x)).
    \end{equation*}
  \end{enumerate}
\end{thm}

\begin{proof}
  We will first prove the forward direction of \ref{item:count-coprod}. Then we will prove both claims in \ref{item:count-Sigma}, and we will prove the reverse direction of claim \ref{item:count-coprod} last.

  For the forward direction of claim \ref{item:count-coprod}, suppose we have equivalences $e:\Fin{k}\simeq A$ and $f:\Fin{l}\simeq B$. The equivalences $e$ and $f$ induce via \cref{ex:coproduct_functor,ex:laws-Fin} a composite equivalence
  \begin{equation*}
    \begin{tikzcd}
      A+B \arrow[r,"\simeq"] & \Fin{k}+\Fin{l} \arrow[r,"\simeq"] & \Fin{k+l},
    \end{tikzcd}
  \end{equation*}
  from which we obtain an element of type $\cnt(A+B)$.

  Next, we will prove the forward direction in the first claim of \ref{item:count-Sigma}, i.e., we will prove that if $A$ comes equipped with an equivalence $e:\Fin{k}\simeq A$, and if $B$ is a family of types over $A$ equipped with
  \begin{equation*}
    f:\prd{x:A}\cnt(B(x)),
  \end{equation*}
  then the total space $\sm{x:A}B(x)$ also has a counting. The proof is by induction on $k$. Note that in the base case, where $k=0$, the type $\sm{x:A}B(x)$ is empty, so it has a counting. For the inductive step, note $\Sigma$ distributes from the right over coproducts. This gives an equivalence
  \begin{align*}
    \sm{x:A}B(x) & \simeq \sm{x:\Fin{k+1}}B(e(x)) \\
    & \simeq \Big(\sm{x:\Fin{k}}B(e(\inl(x)))\Big)+ B(e(\inr(\ttt))).
  \end{align*}
  The type $\sm{x:\Fin{k}}B(e(\inl(x)))$ has a counting by the inductive hypothesis, and the type $B(e(\inr(\ttt)))$ has a counting by assumption. Therefore, it follows that the total space $\sm{x:A}B(x)$ has a counting.
  
  Now we will prove the converse direction of the first claim in \ref{item:count-Sigma}. Suppose that $A$ comes equipped with $e:\Fin{k}\simeq A$, and that $\sm{x:A}B(x)$ comes equipped with $f:\Fin{l}\simeq\sm{x:A}B(x)$. By \cref{rmk:count-decidable-proposition} it suffices to show that, for $a:A$, the type $B(a)$ is a decidable subtype of $\sm{x:A}B(x)$. Consider the map
  \begin{equation*}
    i: B(a)\to \sm{x:A}B(x)
  \end{equation*}
  given by $b\mapsto (a,b)$. For $(x,y):\sm{x:A}B(x)$, we have the equivalences
  \begin{align*}
    \fibf{i}(x,y) & \simeq \sm{b:B(a)}(a,b)=(x,y) \\
                  & \simeq \sm{b:B(a)}\sm{p:a=x}\tr_B(p,b)=y \\
                  & \simeq \sm{p:a=x}\fib{\tr_B(p)}{y} \\
                  & \simeq a=x.
  \end{align*}
  Here we used that $\tr_B(p)$ is an equivalence, and therefore has contractible fibers. Now note that the type $a=x$ is a decidable proposition by \cref{rmk:has-decidable-equality-count}.

  Next, we will prove the second claim in \ref{item:count-Sigma}. Suppose that $B$ is a family over $A$ that comes equipped with a section $b:\prd{x:A}B(x)$, and suppose that each $B(x)$ has a counting, and that the total space $\sm{x:A}B(x)$ has a counting. Then we have a map
  \begin{equation*}
    g : A\to\sm{x:A}B(x)
  \end{equation*}
  given by $a\mapsto (a,b(a))$. The fibers of $g$ can be computed by the following equivalences:
  \begin{align*}
    \fibf{g}(x,y) & \simeq \sm{a:A}(a,b(a))=(x,y) \\
                  & \simeq \sm{a:A}\sm{p:a=x}\tr_B(p,b(a))=y \\
    & \simeq \tr_B(p,b(x))=y.
  \end{align*}
  Note that the type $\tr_B(p,b(x))=y$ is a decidable proposition by \cref{rmk:has-decidable-equality-count}. Now it follows by the forward direction of the first claim in \ref{item:count-Sigma} that $A$ has a counting.

  It remains to prove the converse direction of \ref{item:count-coprod}. Note that the forward direction of the first claim in \ref{item:count-Sigma} implies that countings on a type $X$ induce countings on any decidable subtype of $X$. Note that both $A$ and $B$ are decidable subtypes of the coproduct $A+B$. Any counting of $A+B$ therefore induces countings of $A$ and of $B$.
\end{proof}

\begin{cor}\label{cor:count-prod}
  Consider two types $A$ and $B$. We make two claims:
  \begin{enumerate}
  \item If both $A$ and $B$ come equipped with a counting, then the product $A\times B$ has a counting.
  \item If the product $A\times B$ comes equipped with a counting, then we have two functions
    \begin{align*}
      B & \to \cnt(A) \\
      A & \to \cnt(B).
    \end{align*}
  \end{enumerate}
\end{cor}

\begin{proof}
  The first claim follows from \ref{item:count-Sigma-count-base} in \cref{thm:count}, and the second claim follows from \ref{item:count-Sigma-count-fibers} in \cref{thm:count}. 
\end{proof}
\index{counting|)}

\subsection{Double counting in type theory}
\index{double counting|(}

In combinatorics, counting arguments often proceed by showing that two finite sets are isomorphic---or, in the language of type theory, by showing that two finite types are equivalent. The idea here is, of course, that when we count the elements of a type twice correctly, then both countings must result in the same number. However, this is something that we must prove before we can use it. In other words, we must show that
\begin{equation*}
  (\Fin{k}\simeq\Fin{l})\to (k=l)
\end{equation*}
for any two natural numbers $k$ and $l$. We will prove this claim as a consequence of the following general fact.

\begin{prp}\label{prp:is-injective-maybe}
  For any two types $X$ and $Y$, there is a map
  \begin{equation*}
    (X+\unit\simeq Y+\unit)\to (X\simeq Y).
  \end{equation*}
\end{prp}

\begin{proof}
  We prove the claim in four steps. We will write $i$ for $\inl:X\to X+\unit$ and also for $\inl:Y\to Y+\unit$, and we will write $\star$ for $\inr(\star):X+\unit$ and also for $\inr(\star):Y+\unit$.
  \begin{enumerate}
  \item We first show that for any equivalence $e:X+\unit\simeq Y+\unit$ and any $x:X$ equipped with an identification $p:e(i(x))=\star$, that there is an element
    \begin{equation*}
      \starvalue(e,x,p):Y
    \end{equation*}
    equipped with an identification
    \begin{equation*}
      \alpha:i(\starvalue(e,x,p))=e(\star).  
    \end{equation*}
    To see this, note that the map $e$ is injective. The elements $i(x)$ and $\star$ are distinct, so it follows that the elements $e(i(x))$ and $e(\star)$ are distinct. In particular, we have $e(\star)\neq\star$. Therefore it follows that there is an element $y:Y$ equipped with an identification $i(y)=e(\star)$. 
  \item Next, we construct for every equivalence $e:X+\unit\simeq Y+\unit$ a map $f:X\to Y$ equipped with identifications
    \begin{align*}
      \beta & : \prd{y:Y} (e(i(x))=i(y))\to (f(x)=y) \\
      \gamma & : \prd{p:e(i(x))=\star} f(x)=\starvalue(e,x,p).
    \end{align*}
    In order to construct the map $f:X\to Y$, we first construct a dependent function
    \begin{equation*}
      f':\prd{x:X}\prd{u:Y+\unit}((e(i(x))=u)\to Y).
    \end{equation*}
    This function is defined by pattern matching on $u$, by
    \begin{align*}
      f'(x,i(y),p) & \defeq y \\
      f'(x,\star,p) & \defeq \starvalue(e,x,p)
    \end{align*}
    Then we define $f(x):=f'(x,e(i(x)),\refl{})$. By the definition of $f'$ it then follows that we have an identification
    \begin{align*}
      f(x) & \jdeq f'(x,e(i(x)),\refl{}) \\
           & = f'(x,i(y),p) \\
           & \jdeq y
    \end{align*}
    for any $y:Y$ and $p:e(i(x))=i(y)$, and that we have an identification
    \begin{align*}
      f(x) & \jdeq f'(x,e(i(x)),\refl{}) \\
           & = f'(x,\ttt,p) \\
           & \jdeq \starvalue(e,x,p)
    \end{align*}
    for any $p:e(i(x))=\ttt$. 
  \item The inverse function $g:Y\to X$ is constructed in the same way as the function $f:X\to Y$, using the equivalence $e^{-1}:Y+\unit\simeq X+\unit$. This function comes equipped with
    \begin{align*}
      \delta & : \prd{x:X}(e^{-1}(i(y))=i(x))\to (g(y)=x) \\
      \varepsilon & : \prd{p:e^{-1}(i(y))=\star}g(y)=\starvalue(e^{-1},y,p). 
    \end{align*}
  \item It remains to show that $f$ and $g$ are inverse to each other. The proof that $g$ is a retraction of $f$ is similar to the proof that $g$ is a section of $f$, so we will only prove the latter. In other words, we will construct an identification
    \begin{equation*}
      f(g(y))=y
    \end{equation*}
    for any $y:Y$. The proof is by case analysis on $(e^{-1}(i(y))=\ttt)+(e^{-1}(i(y))\neq \ttt)$. In the case where $p:e^{-1}(i(y))=\ttt$, we have the identification
    \begin{equation*}
      \varepsilon(p):g(y)=\starvalue(e^{-1},y,p).
    \end{equation*}
    Furthermore, we have the identification
    \begin{equation*}
      \gamma(q) : f(g(y)) = \starvalue(e,g(y),q),
    \end{equation*}
    where $q:e(i(g(y)))=\ttt$ is the composite of the identifications
    \begin{align*}
      e(i(g(y))) & = e(i(\starvalue(e^{-1},y,p))) \\
                 & = e(e^{-1}(\ttt)) \\
      & =\ttt.
    \end{align*}
    Using the identification $\gamma(q)$, we obtain
    \begin{align*}
      i(f(g(y))) & = i(\starvalue(e,g(y),q)) \\
                 & = e(\ttt) \\
                 & = e(e^{-1}(i(y))) \\
                 & = i(y).
    \end{align*}
    Since $i:Y\to Y+\unit$ is injective, it follows that $f(g(y))=y$. 
    \qedhere
  \end{enumerate}
\end{proof}

\begin{thm}\label{thm:is-injective-Fin}
  For any two natural numbers $k$ and $l$, there is a map
  \begin{equation*}
    (\Fin{k}\simeq\Fin{l})\to (k=l).
  \end{equation*}
\end{thm}

\begin{proof}
  The proof is by induction on $k$ and $l$. In the base case, where both $k$ and $l$ are zero, the claim is obvious. If $k$ is zero and $l$ is a successor, then we have $0:\Fin{l}$. Any equivalence $e:\Fin{k}\simeq \Fin{l}$ now gives us the element
  \begin{equation*}
    e^{-1}(0):\emptyt,
  \end{equation*}
  which is of course absurd. Similarly, if $k$ is a successor and $l$ is zero, we obtain $e(0):\emptyt$, which is again absurd. If both $k$ and $l$ are a successor, then we have by \cref{prp:is-injective-maybe} the composite
  \begin{equation*}
    \begin{tikzcd}[column sep=1.5em]
      (\Fin{k+1}\simeq\Fin{l+1}) \arrow[r] & (\Fin{k}\simeq\Fin{l}) \arrow[r] & (k=l) \arrow[r] & (k+1=l+1).
    \end{tikzcd}\qedhere
  \end{equation*}
\end{proof}
\index{double counting|)}

\subsection{Finite types}

The type of all finite types is the subtype of the base universe $\UU_0$ consisting of all types $X$ for which there exists an unspecified equivalence $\Fin{k}\simeq X$ for some $k:\N$.

\begin{defn}\label{defn:finite}
  A type $X$ is said to be \define{finite}\index{finite type|textbf} if it comes equipped with an element of type\index{is-finite@{$\isfinite(X)$}|textbf}
  \begin{equation*}
    \isfinite(X) \defeq \Brck{\sm{k:\N}\Fin{k}\simeq X}
  \end{equation*}
  The type $\F$ of all finite types is defined to be\index{F@{$\F$}|see {finite type}}\index{F@{$\F$}|textbf}
  \begin{equation*}
    \F:=\sm{X:\UU_0}\isfinite(X).
  \end{equation*}
  In other words, the type $\F$ of finite types is the image of the map $\Fin{} : \N \to \UU_0$.
  We also define the type $\BS_k$\index{BS n@{$\BS_n$}|textbf} of \define{$k$-element types} by
  \begin{equation*}
    \BS_k\defeq \sm{X:\UU_0}\brck{\Fin{k}\simeq X}.
  \end{equation*}
\end{defn}

\begin{rmk}
  It follows directly from the definition of finiteness that any type $X$ equipped with a counting is finite. In particular, any $\Fin{k}$ is finite. Furthermore, it follows that if $X$ is equivalent to a finite type $Y$, then $X$ is also finite. Indeed, we can use the functoriality of the propositional truncation to obtain a function
  \begin{equation*}
    \Brck{\sm{k:\N}\Fin{k}\simeq Y}\to\Brck{\sm{k:\N}\Fin{k}\simeq X}
  \end{equation*}
  from a map $\big(\sm{k:\N}\Fin{k}\simeq Y\big)\to\big(\sm{k:\N}\Fin{k}\simeq X\big)$. Given an equivalence $e:X\simeq Y$, such a map is given as the map induced on total spaces from the family of maps $f\mapsto e^{-1}\circ f$.

  Similarly, it follows that any finite type has decidable equality, and that every finite type is a set.
\end{rmk}

In the following proposition we will show that each finite type can be assigned a unique cardinality.

\begin{thm}
  For any type $X$, consider the type $\isfinite'(X)$ defined by\index{is-finite'@{$\isfinite'(X)$}|textbf}
  \begin{equation*}
    \isfinite'(X) \defeq \sm{k:\N}\brck{\Fin{k}\simeq X}.
  \end{equation*}
  Then the type $\isfinite'(X)$\index{is-finite'@{$\isfinite'(X)$}!is a proposition} is a proposition, and there is an equivalence
  \begin{equation*}
    \isfinite(X)\leftrightarrow\isfinite'(X).
  \end{equation*}
  If $X$ is a finite type, then the unique number $k$ such that $\brck{\Fin{k}\simeq X}$ is the \define{cardinality}\index{cardinality!of a finite type|textbf} of $X$. We write $|X|$\index{{"|"}X{"|"}@{$|X|$}|see {cardinality, of a finite type}} for the cardinality of $X$.
\end{thm}

\begin{proof}
  We first prove the claim that the type $\isfinite'(X)$ is a proposition. In other words, we need to show that any two natural numbers $k$ and $k'$ for which there are respective elements of the types $\brck{\Fin{k}\simeq X}$ and $\brck{\Fin{k'}\simeq X}$, can be identified.

  Since the type of natural numbers is a set, the type $k=k'$ is a proposition. Therefore, we may assume that we have equivalences $\Fin{k}\simeq X$ and $\Fin{k'}\simeq X$. Consequently, we have an equivalence $\Fin{k}\simeq\Fin{k'}$. Now it follows from \cref{thm:is-injective-Fin} that $k=k'$.

  The second claim is that the propositions $\isfinite(X)$ and $\isfinite'(X)$ are equivalent, which we will show by constructing functions back and forth.
  Since we have shown that the type $\isfinite'(X)$ is a proposition, we obtain a map $\isfinite(X)\to\isfinite'(X)$ via the universal property of the propositional truncation, from the map
  \begin{equation*}
    \Big(\sm{k:\N}\Fin{k}\simeq X\Big)\to \sm{k:\N}\brck{\Fin{k}\simeq X}
  \end{equation*}
  given by $(k,e)\mapsto (k,\eta(e))$. 
  
  To construct a map $\isfinite'(X)\to\isfinite(X)$, it suffices to construct a map
  \begin{equation*}
    \brck{\Fin{k'}\simeq X}\to \Brck{\sm{k:\N}\Fin{k}\simeq X}
  \end{equation*}
  for each $k':\N$. Again by the universal property of the propositional truncation, we obtain this map from the function
  \begin{equation*}
    (\Fin{k'}\simeq X) \to \Brck{\sm{k:\N}\Fin{k}\simeq X}
  \end{equation*}
  given by $e\mapsto \eta(k',e)$. 
\end{proof}

\begin{cor}
  There is an equivalence
  \begin{equation*}
    \F \simeq \sm{k:\N}\BS_k.
  \end{equation*}
\end{cor}

\begin{proof}
  This equivalence can be obtained by composing the equivalences
  \begin{align*}
    \sm{X:\UU_0}\isfinite(X) & \simeq \sm{X:\UU_0}\sm{k:\N}\brck{\Fin{k}\simeq X} \\
    & \simeq \sm{k:\N}\sm{X:\UU_0}\brck{\Fin{k}\simeq X}. \qedhere
  \end{align*}
\end{proof}

We now aim to extend \cref{thm:count} to obtain some closure properties of finite types. Before we do so, we prove the \textbf{principle of finite choice}.

\begin{prp}\label{prp:finite-choice}
  Consider a type family $B$ over a finite type $A$. Then there is a \define{finite choice}\index{finite choice|textbf}\index{finite type!finite choice|textbf} map
  \begin{equation*}
    \Big(\prd{x:A}\brck{B(x)}\Big)\to\Brck{\prd{x:A}B(x)}
  \end{equation*}
\end{prp}

\begin{proof}
  Note that the type $\big\|\prd{x:A}B(x)\big\|$ is a proposition. Therefore we may assume that the type $A$ comes equipped with a counting $e:\Fin{k}\simeq A$. By this equivalence, it suffices to show that for every type family $B$ over $\Fin{k}$, there is a map
  \begin{equation*}
    \Big(\prd{x:\Fin{k}}\brck{B(x)}\Big)\to\Brck{\prd{x:\Fin{k}}B(x)}.
  \end{equation*}
  We proceed by induction on $k$. In the base case, $\Fin{k}$ is empty and therefore the type $\prd{x:\Fin{k}}B(x)$ is contractible. The asserted function therefore exists.

  For the inductive step, note that by the dependent universal property of coproducts (\cref{ex:up-coproduct}) we have the equivalences
  \begin{align*}
    \Big(\prd{x:\Fin{k+1}}\brck{B(x)}\Big) & \simeq \Big(\prd{x:\Fin{k}}\brck{B(i(x))}\Big)\times \brck{B(\ttt)} \\
    \Brck{\prd{x:\Fin{k}}B(x)} & \simeq \Brck{\Big(\prd{x:\Fin{k}}B(i(x))\Big)\times B(\ttt)}.
  \end{align*}
  Recall from \cref{ex:product-propositional-truncation} that $\brck{X\times Y}\simeq \brck{X}\times\brck{Y}$ for any two types $X$ and $Y$. This fact together with the inductive hypothesis finishes the proof.
\end{proof}

\begin{thm} ~
  \begin{enumerate}
    \item \label{item:coproduct-finite-types}For any two types $X$ and $Y$, the following are equivalent:
    \begin{enumerate}
    \item Both $X$ and $Y$ are finite.
    \item The coproduct $X+Y$ is finite.
    \end{enumerate}
  \item \label{item:product-finite-types}For any two types $X$ and $Y$, we make two claims:
    \begin{enumerate}
    \item If both $X$ and $Y$ are finite, then the cartesian product $X\times Y$ is finite.
    \item If the type $X\times Y$ is finite, then we have two functions
      \begin{align*}
        Y & \to \isfinite(X) \\
        X & \to \isfinite(Y).
      \end{align*}
    \end{enumerate}
  \item \label{item:Sigma-finite-types}Consider a type family $B$ over $A$, and consider the following three conditions:
    \begin{enumerate}
    \item The type $A$ is finite.
    \item The type $B(x)$ is finite for each $x:A$.
    \item The type $\sm{x:A}B(x)$ is finite.
    \end{enumerate}
    If (a) holds, then (b) is equivalent to (c). Moreover, if (b) and (c) hold, then (a) holds if and only if $A$ is a set and the type $\sm{x:A}\neg B(x)$ is finite. Furthermore, if (b) and (c) hold and $B$ has a section, then (a) holds.
  \end{enumerate}
\end{thm}

\begin{proof}
  To prove claim \ref{item:coproduct-finite-types}, first suppose that both $X$ and $Y$ are finite. Since the type $\isfinite(X+Y)$ is a proposition, we may assume that $X$ and $Y$ come equipped with countings. It follows from \cref{thm:count} that $X+Y$ has a counting, so it is finite. Conversely, suppose that the type $X+Y$ is finite. Since the types $\isfinite(X)$ and $\isfinite(Y)$ are both propositions, we may assume that the coproduct $X+Y$ comes equipped with a counting. Again it follows from \cref{thm:count} that the types $X$ and $Y$ have countings, so they are finite.

  The proof of claim \ref{item:product-finite-types} is similar to the proof of claim \ref{item:coproduct-finite-types}, hence we omit it.

  It remains to prove claim \ref{item:Sigma-finite-types}. First, suppose that the type $A$ is finite, and that each $B(x)$ is finite. By \cref{prp:finite-choice} we have a map
  \begin{equation*}
    \Big(\prd{x:A}\isfinite(B(x))\Big)\to \Brck{\prd{x:A}\cnt(B(x))}.
  \end{equation*}
  Since our goal is to construct an element of a proposition, we may therefore assume that each $B(x)$ comes equipped with a counting. We may also assume that $A$ comes equipped with a counting. It follows from \cref{thm:count} that the type $\sm{x:A}B(x)$ has a counting, so it is finite.

  Next, assume that $A$ is finite and that the type $\sm{x:A}B(x)$ is finite, and let $a:A$. The type $\isfinite(B(a))$ is a proposition, so we may assume that the types $A$ and $\sm{x:A}B(x)$ come equipped with countings. It follows from \cref{thm:count} that $B(a)$ has a counting, so it is finite.

  The final claim has two parts. First, assume that each $B(x)$ is finite, that the type $\sm{x:A}B(x)$ is finite, and that the type family $B$ has a section $f:\prd{x:A}B(x)$. It follows that the map
  \begin{equation*}
    A\to\sm{x:A}B(x)
  \end{equation*}
  given by $x\mapsto (x,f(x))$ is a decidable embedding, because the fiber at $(x,y)$ of this map is equivalent to the identity type $f(x)=y$ in $B(x)$, which is a decidable proposition. It follows from the fact that (a) and (b) together imply (c) that $A$ is finite.

  For the remaining part of the final claim, assume that $A$ is a set. Note that the assumption that each $B(x)$ is finite implies that each $B(x)$ is either inhabited or empty. It follows that we have an equivalence
  \begin{equation*}
    A\simeq \Big(\sm{x:A}\brck{B(x)}\Big)+\Big(\sm{x:A}\neg B(x)\Big).
  \end{equation*}
  We assume that the type $\sm{x:A}\neg B(x)$ is finite. In order to show that $A$ is finite, it therefore suffices to show that the type $\sm{x:A}\brck{B(x)}$ is finite. Without loss of generality, we assume that each $B(x)$ is inhabited. To finish the proof, it suffices to show that there is an element of type
  \begin{equation*}
    \Brck{\prd{x:A}B(x)}
  \end{equation*}
  using the assumption that $\prd{x:A}\brck{B(x)}$. To construct such an element, we may assume a counting $e:\Fin{k}\simeq\sm{x:A}B(x)$. We claim that there is a function
  \begin{equation*}
    \brck{B(a)}\to B(a),
  \end{equation*}
  i.e., that the type $B(a)$ satisfies the principle of global choice of \cref{rmk:global-choice} for each $a:A$. Recall from \cref{eg:global-choice-decidable-subtype-N} that the decidable subtypes of $\Fin{k}$ satisfy global choice. Therefore it also follows that the decidable subtypes of $\sm{x:A}B(x)$ satisfy global choice. Thus, it suffices to show that $B(x)$ is a decidable subtype of $\sm{x:A}B(x)$.
  
   The assumption that $A$ is a set implies by \cref{ex:is-trunc-fiber-inclusion} that the fiber inclusion $i_a:B(a)\to\sm{x:A}B(x)$ is an embedding for each $a:A$. Furthermore, we note that we have the following equivalence computing the fibers of $i_a$ at $(x,y)$:
  \begin{equation*}
    \Big(\sm{z:B(a)}(a,z)=(x,y)\Big)\simeq (a=x).
  \end{equation*}
  The type on the left hand side is decidable, so it follows that the type $A$ has decidable equality. We conclude that each $B(a)$ is a decidable subtype of $\sm{x:A}B(x)$.
\end{proof}

\begin{exercises}
  \exitem
  \begin{subexenum}
  \item Construct an equivalence $\Fin{n^m}\simeq(\Fin{m}\to\Fin{n})$. Conclude that if $A$ and $B$ are finite types, then $A\to B$ is finite.
  \item Construct an equivalence $\Fin{n!}\simeq(\Fin{n}\simeq\Fin{n})$. Conclude that if $A$ is finite, then $A\simeq A$ is finite.
  \end{subexenum}
  \exitem Suppose that $A$ is a retract of $B$. Show that $\cnt(B)\to\cnt(A)$.
  Conclude that $\isfinite(B)\to\isfinite(A)$.\index{finite type!closed under retracts}
  \exitem 
  \begin{subexenum}
  \item Consider a family of decidable types $A_i$ indexed by a finite type $I$. Show that the dependent product
    \begin{equation*}
      \prd{i:I}A_i
    \end{equation*}
    is decidable.
  \item Show that $\isemb(f)$ is decidable, for any map $f:I\to J$ between finite types.
  \item Show that $\issurj(f)$ is decidable, for any map $f:I\to J$ between finite types.
  \item Show that $\isequiv(f)$ is decidable, for any map $f:I\to J$ between finite types.
  \end{subexenum}
  \exitem \label{item:quotient-finite-types}Consider a surjective map $f:A\to B$, and suppose that $A$ is finite. Show that the following are equivalent:
    \begin{enumerate}
    \item The type $B$ has decidable equality.
    \item The type $B$ is finite.
    \end{enumerate}
  \exitem Consider a family $B$ of types over $A$.
  \begin{subexenum}
  \item Show that if $A$ is finite and if each $B(x)$ is finite, then the type
    \begin{equation*}
      \prd{x:A}B(x)
    \end{equation*}
    is finite.
  \item Show that if $A$ is finite and if $\prd{x:A}B(x)$ is finite, then we have
    \begin{equation*}
      \Big(\prd{x:A}\brck{B(x)}\Big)\to\Big(\prd{x:A}\isfinite(B(x))\Big).
    \end{equation*}
  \item Show that if $\prd{x:A}B(x)$ is finite and if each $B(x)$ is finite, then $A$ is finite if and only if the following three conditions hold:
    \begin{enumerate}
    \item $A$ has decidable equality.
    \item The type
      \begin{equation*}
        \sm{x:A}|B(x)|\leq 1
      \end{equation*}
      is finite.
    \item The type
      \begin{equation*}
        \prd{x:A}(2\leq|B(x)|)\to B(x)
      \end{equation*}
      is finite.
    \end{enumerate}
  \end{subexenum}
  \exitem Consider two finite types $X$ and $Y$ with $m$ and $n$ elements, respectively, and let $f:X\to Y$ be a map.
  \begin{subexenum}
  \item Show that
    \begin{equation*}
      \isinj(f)\to (m\leq n).
    \end{equation*}
  \item Prove the \define{pigeonhole principle}\index{pigeonhole principle|textbf}\index{finite type!pigeonhole principle|textbf}, i.e., show that
    \begin{equation*}
      (n>m)\to \exists_{(x,x':X)}(x\neq x')\times(f(x)=f(x')).
    \end{equation*}
  \item Show that there is no embedding $\N\hookrightarrow \Fin{k}$, for any $k:\N$.
  \end{subexenum}
  \exitem Consider a finite type $X$.
  \begin{subexenum}
  \item Show that any embedding $f:X\to X$ is an equivalence. Sets $X$ such that every embedding $X\hookrightarrow X$ is an equivalence are also called \define{Dedekind finite}.\index{Dedekind finite type|textbf}\index{finite type!Dedekind finite type|textbf}
  \item Show that any surjective map $f:X\to X$ is an equivalence.
  \end{subexenum}
  \exitem Consider two arbitrary types $A$ and $B$. For any $2$-element type $X$, construct an equivalence
  \begin{equation*}
    (A+B)^X\simeq A^X+X\times (A\times B)+B^X.
  \end{equation*}
  \exitem
  \begin{subexenum}
  \item Consider a set $A$ and an arbitrary type $B$. Show that any embedding $A\hookrightarrow B$ factors uniquely through the embedding $(\unit\hookrightarrow B)\hookrightarrow B$ given by $e\mapsto e(\ttt)$. 
  \item A map $f:A\to B$ is said to be \define{decidable}\index{decidable map} if the type $\fib{f}{b}$ is decidable for all $b:B$. Write $A\demb B$\index{A hookrightarrow d B@{$A \demb B$}|textbf} for the type of decidable embeddings\index{decidable embedding} from $A$ to $B$. Show that for any type $A$ with decidable equality and an arbitrary type $B$, any decidable embedding $A\demb B$ factors uniquely through the embedding $(\unit\demb B)\emb B$.
  \item (Escard\'o) For any two types $A$ and $B$, construct an equivalence
  \begin{equation*}
    ((A+\unit)\simeq(B+\unit))\simeq (\unit\demb (B+\unit))\times(A\simeq B).
  \end{equation*}
  \end{subexenum}
  \exitem
  \begin{subexenum}
  \item For any two types $A$ and $B$, construct an equivalence
    \begin{equation*}
      ((A+\unit)\demb(B+\unit))\simeq (\unit \demb (B+\unit))\times (A\demb B).
    \end{equation*}
  \item Construct an equivalence $\Fin{\fallingfactorial{n}{m}}\simeq(\Fin{m}\hookrightarrow\Fin{n})$, where $\fallingfactorial{n}{m}$ is the \define{$m$-th falling factorial}\index{falling factorial|textbf}\index{(n) m@{$\fallingfactorial{n}{m}$}|see {falling factorial}} of $n$, which is defined recursively by
    \begin{align*}
      \fallingfactorial{0}{0} & \defeq 1 & \fallingfactorial{0}{m+1} & \defeq 0 \\*
      \fallingfactorial{n+1}{0} & \defeq 1 & \fallingfactorial{n+1}{m+1} & \defeq (n+1)\fallingfactorial{n}{m}.
    \end{align*}
    Conclude that if $A$ and $B$ are finite with cardinality $m$ and $n$, then the type $A\hookrightarrow B$ is finite with cardinality $\fallingfactorial{n}{m}$.
  \end{subexenum}
  \exitem
  \begin{subexenum}
  \item Consider an arbitrary type $A$ and a type $B$ with decidable equality. Construct an equivalence
    \begin{equation*}
      ((A+\unit)\twoheadrightarrow(B+\unit))\simeq (B+\unit)\times(A\twoheadrightarrow B)+(A\twoheadrightarrow B+\unit).
    \end{equation*}
  \item Construct an equivalence $\Fin{\numberofsurjectivemaps{m}{n}}\simeq(\Fin{m}\twoheadrightarrow\Fin{n})$, where $\numberofsurjectivemaps{m}{n}$\index{S(m,n)@{$\numberofsurjectivemaps(m,n)$}|textbf} is defined recursively by
    \begin{align*}
      \numberofsurjectivemaps{0}{0} & \defeq 1 \\*
      \numberofsurjectivemaps{0}{n+1} & \defeq 0 \\*
      \numberofsurjectivemaps{m+1}{0} & \defeq 0 \\*
      \numberofsurjectivemaps{m+1}{n+1} & \defeq (n+1)\numberofsurjectivemaps{m}{n}+\numberofsurjectivemaps{m}{n+1}.
    \end{align*}
    Conclude that if $A$ and $B$ are finite with cardinality $m$ and $n$, then the type $A\twoheadrightarrow B$ is finite with cardinality $\numberofsurjectivemaps{m}{n}$. Note: the number $\numberofsurjectivemaps{m}{n}$ is $n!\stirling{m}{n}$, where $\stirling{m}{n}$\index{{{n m}}@{$\stirling{n}{m}$}|see {Stirling number of the second kind}} is the \define{Stirling number of the second kind}\index{Stirling number of the second kind} at $(m,n)$.
  \end{subexenum}
\end{exercises}

\section{The univalence axiom}
\index{univalence axiom|(}
\index{axiom!univalence|(}

The univalence axiom characterizes the identity type of a universe. Roughly speaking, it asserts that equivalent types are equal. The univalence axiom therefore postulates the common mathematical habit of identifying equivalent objects, such as equivalent types, isomorphic groups, isomorphic rings, logically equivalent propositions, subsets with the same elements, and so on. The univalence axiom is due to Voevodsky, who also showed that it is modeled in the simplicial sets. He also showed, in one of his first applications, that the univalence axiom implies function extensionality, which we will also prove here.

One way to think about the univalence axiom is that it \emph{expands} the notion of equality to encapsulate the notion of equivalence. It asserts that for each equivalence $e$ between two types $X$ and $Y$ in a universe $\mathcal{U}$ there is a unique identification $p_e:X=Y$ in the universe $\mathcal{U}$ such that transporting along $p_e$ in the universal type family over $\mathcal{U}$ is homotopic to the original equivalence $e:X\simeq Y$. 

Since there might be many distinct equivalences between two types $X$ and $Y$, there will be equally many identifications those types. The univalence axiom is therefore inconsistent with the commonly assumed axiom that all identity types are propositions, i.e., that all types are sets. Indeed, there are two equivalences $\bool\simeq\bool$, so a univalent universe cannot be a set.

\subsection{Equivalent forms of the univalence axiom}
By the fundamental theorem of identity types, \cref{thm:id_fundamental}, it is immediate that the univalence axiom comes in three equivalent forms.

\begin{thm}\label{thm:univalence}\index{identity type!of a universe}\index{characterization of identity type!of a universe}\index{universe!characterization of identity type}
Consider a universe $\UU$. The following are equivalent:
\begin{enumerate}
\item The universe $\UU$ is \define{univalent}\index{univalent universe|textbf}: For any two types $A,B:\UU$, the map\index{equiv-eq@{$\equiveq$}|textbf}
  \begin{equation*}
    \equiveq:(A=B)\to (A\simeq B)
  \end{equation*}
  given by $\equiveq(\refl{}):=\idfunc$, is an equivalence.
\item The type
\begin{equation*}
\sm{B:\UU}\eqv{A}{B}
\end{equation*}
is contractible for each $A:\UU$.
\item For any type $A:\UU$, the family of types $A\simeq X$ indexed by $X:\UU$ is an identity system on $\UU$. In other words, the universe $\UU$ satisfies the principle of \define{equivalence induction}\index{equivalence induction|textbf}\index{induction principle!for equivalences|textbf}: For every $A:\UU$ and for every type family of types $P(X,e)$ indexed by $X:\UU$ and $e:A\simeq X$, the map
\begin{equation*}
\Big(\prd{X:\UU}\prd{e:\eqv{A}{X}}P(X,e)\Big)\to P(A,\idfunc)
\end{equation*}
given by $f\mapsto f(A,\idfunc)$ has a section.
\end{enumerate}
\end{thm}

\begin{proof}
  The claim is a special case of \cref{thm:id_fundamental}, the fundamental theorem of identity types\index{fundamental theorem of identity types}.
\end{proof}

One way to see that the univalence axiom is plausible, is by observing that all type constructors preserve equivalences. For example, in \cref{thm:equiv-toto} we showed that for any type family $B$ over $A$ and any type family $B'$ over $A'$, if we have an equivalence $e:A\simeq A'$ and family of equivalences $f:\prd{x:A}B(x)\simeq B'(e(x))$, then we obtain an equivalence
\begin{equation*}
  \Big(\sm{x:A}B(x)\Big)\simeq\Big(\sm{x':A'}B'(x')\Big).
\end{equation*}
Under the same assumptions, we showed in \cref{ex:equiv-pi} that we obtain an equivalence
\begin{equation*}
  \Big(\prd{x:A}B(x)\Big)\simeq\Big(\prd{x':A'}B'(x')\Big).
\end{equation*}
Furthermore, for any two elements $x,y:A$ any equivalence $e:A\simeq A'$ induces an equivalence $(x=y)\simeq (e(x)=e(y))$ by \cref{cor:emb_equiv}. In other words, all the standard type formers within a universe $\UU$ are \emph{equivalence invariant}. Since identity types are not assumed to be propositions, we have the possibility to postulate the univalence axiom.

\begin{axiom}[The univalence axiom]\label{axiom:univalence}\index{univalence axiom|textbf}\index{axiom!univalence|textbf}
  We will assume that all the universes generated by \cref{enough-universes}\index{enough universes}\index{universe!enough universes} are univalent. Given a univalent universe $\UU$, we will write $\eqequiv$\index{eq-equiv@{$\eqequiv$}|textbf} for the inverse of $\equiveq$.
\end{axiom}

As a first application of the univalence axiom, let us show that for any type $A$ the type of types in a univalent universe $\UU$ that are equivalent to $A$ is a proposition.

\begin{defn}\label{defn:small-types}
  Consider a univalent universe $\UU$. A type $X$ is said to be \define{$\UU$-small}\index{small type|textbf}\index{U-small type@{$\UU$-small type}|textbf} if it comes equipped with an element of type\index{is-small@{$\issmall_\UU(A)$}|textbf}
  \begin{equation*}
    \issmall_\UU(A)\defeq\sm{X:\UU}A\simeq X.
  \end{equation*}
  Similarly, a map $f:A\to B$ is said to be \define{$\UU$-small}\index{small map|textbf}\index{U-small map@{$\UU$-small map}|textbf} if all of its fibers are $\UU$-small.
\end{defn}

\begin{eg}
  ~
  \begin{enumerate}
  \item Any type in $\UU$ is $\UU$-small.
  \item Any contractible type is $\UU$-small with respect to any universe $\UU$. \index{contractible type!is U-small@{is $\UU$-small}}\index{small type!contractible type}
  \item For any family $P$ of $\UU$-small types over a $\UU$-small type $A$, the dependent product $\prd{x:A}B(x)$ is $\UU$-small.\index{small type!dependent function type}\index{dependent function type!is U-small@{is $\UU$-small}}
  \item The type of $\UU$-small types in $\VV$ is equivalent to the type of $\VV$-small types in $\UU$. This follows from the equivalence
    \begin{equation*}
      \Big(\sm{Y:\VV}\sm{X:\UU}Y\simeq X\Big) \simeq \Big(\sm{X:\UU}\sm{Y:\VV}X\simeq Y\Big).
    \end{equation*}
  \item Any finite type is $\UU$-small for any universe $\UU$.\index{finite type!is U-small@{is $\UU$-small}}\index{small type!finite type} Consequently, we get equivalences
  \begin{equation*}
    \Big(\sm{X:\UU}\isfinite(X)\Big)\simeq\Big(\sm{Y:\VV}\isfinite(Y)\Big)
  \end{equation*}
  for any two univalent universes $\UU$ and $\VV$. This observation is the reason why we usually write $\F$ for the type of finite types (in $\UU$), without referring to its universe.
  \item In \cref{thm:russell} we will show that $\UU$ cannot be $\UU$-small, i.e., that there cannot be a type $U:\UU$ equipped with an equivalence $U\simeq \UU$.
  \end{enumerate}
\end{eg}

\begin{prp}\label{prp:small}
  For any univalent universe $\UU$ and any type $A$, the type $\issmall_\UU(A)$ is a proposition.\index{is-small@{$\issmall_\UU(A)$}!is a proposition}
\end{prp}

\begin{proof}
  By \cref{lem:isprop_eq} it suffices to show that
  \begin{equation*}
    \issmall_\UU(A)\to\iscontr(\issmall_\UU(A)).
  \end{equation*}
  Let $X:\UU$ be a type equipped with $e:A\simeq X$. Then we have an equivalence
  \begin{equation*}
    \Big(\sm{Y:\UU}A\simeq Y\Big)\simeq\Big(\sm{Y:\UU}X\simeq Y\Big).
  \end{equation*}
  The latter type is contractible by \cref{thm:univalence}.
\end{proof}

\begin{cor}
  Consider a univalent universe $\UU$ and a univalent universe $\VV$ containing all types in $\UU$. Then the universe inclusion $i:\UU\to\VV$ is an embedding.
\end{cor}

\begin{proof}
  Since $\VV$ is assumed to be univalent, it follows that
  \begin{equation*}
    \fib{i}{A}\simeq\issmall_\UU(A)
  \end{equation*}
  for any type $A:\VV$. The type $\issmall_\UU(A)$ is a proposition since $\UU$ is univalent. Hence the claim follows by \cref{thm:embedding}.
\end{proof}

\subsection{Propositional extensionality}
\index{propositional extensionality|(}

An important direct consequence of the univalence axiom is the principle of propositional extensionality. This principle asserts that any two logically equivalent propositions $P$ and $Q$ can be identified. Propositional extensionality is an important principle on its own, which is sometimes assumed in formal systems without the univalence axiom.

In order to prove propositional extensionality, we first observe that the univalence axiom also characterizes the identity type of any subuniverse.

\begin{prp}\label{prp:univalence-subuniverse}
  Consider a universe $\UU$, and let $P$ be a family of propositions over $\UU$. Then the family of maps
  \begin{equation*}
    \equiveq:(A=B)\to (\proj 1(A) \simeq \proj 1(B))
  \end{equation*}
  indexed by $A,B:\sm{X:\UU}P(X)$, given by $\equiveq(\refl{}):=\idfunc$ is an equivalence.
\end{prp}

\begin{proof}
  Since $P$ is a subuniverse, it follows from \cref{cor:pr1-embedding} that the projection map is an embedding. Therefore we see that the asserted map is the composite of the equivalences
  \begin{equation*}
    \begin{tikzcd}
      (A=B) \arrow[r,"\apfunc{\proj 1}"] & (\proj 1(A)=\proj 1(B)) \arrow[r,"\equiveq"] &[2em] (\proj 1(A)\simeq \proj 1(B)).
    \end{tikzcd}\qedhere
  \end{equation*}
\end{proof}

\begin{rmk}
  Often, when $P$ is a subuniverse, i.e., a subtype of the a universe $\UU$, we will also write $A$ for the type $\proj 1(A)$ if $A:\sm{X:\UU}P(X)$. Using this shorthand notation, the equivalence in \cref{prp:univalence-subuniverse} is displayed as
  \begin{equation*}
    (A=B)\simeq (A\simeq B).
  \end{equation*}
\end{rmk}

Important examples of subuniverses include the subuniverse $\prop_\UU$ of propositions in $\UU$, the subuniverse $\Set_\UU$ of sets in $\UU$, and the subuniverse $\UU^{\leq k}$ of $k$-truncated types in $\UU$. The subuniverse $\F$ of finite types in $\UU_0$, and the subuniverses $\BS_k$ of $k$-element types are further important subuniverses to which \cref{prp:univalence-subuniverse} applies. Note that by the univalence axiom, any subuniverse is automatically closed under equivalences\index{subuniverse!closed under equivalences}. Indeed, if we have $X\simeq Y$, then we have $P(X)\to P(Y)$ by transporting along the equality $X=Y$ induced by univalence.

\begin{thm}\label{prp:propositional-extensionality}
  Propositions satisfy \define{propositional extensionality}\index{propositional extensionality|textbf}\index{extensionality principle!for propositions|textbf}:
  For any two propositions $P$ and $Q$, the canonical map\index{bi-implication}\index{iff-eq@{$\iffeq$}|textbf}
  \begin{equation*}
    \iffeq:(P=Q)\to (P\iffprop Q)
  \end{equation*}
  defined by $\iffeq(\refl{}):=(\idfunc,\idfunc)$ is an equivalence. It follows that the type $\prop_\UU$ of propositions in $\UU$ is a set.\index{Prop@{$\prop_\UU$}!is a set}\index{univalence axiom!implies propositional extensionality}
\end{thm}

\begin{proof}
  Recall from \cref{ex:isprop_istrunc} that $\isprop(X)$ is a proposition for any type $X$. \cref{prp:univalence-subuniverse} therefore applies, which gives
  \begin{equation*}
    (P=Q)\simeq (P\simeq Q)\simeq (P\leftrightarrow Q).
  \end{equation*}
  The last equivalence follows from \cref{prp:equiv-prop}, using the fact that $(P\simeq Q)$ is a proposition by \cref{ex:isprop_isequiv}.
\end{proof}

\begin{cor}\label{cor:decidable-Prop}
  The type\index{DProp@{$\decidableProp_\UU$}|see {decidable proposition}}\index{DProp@{$\decidableProp_\UU$}|textbf}\index{decidable proposition|textbf}\index{DProp@{$\decidableProp_\UU$}!is equivalent to bool@{is equivalent to $\bool$}}
  \begin{equation*}
    \decidableProp_\UU \defeq \sm{P:\prop_\UU}\isdecidable(P)
  \end{equation*}
  of decidable propositions in any universe $\UU$ is equivalent to $\bool$.
\end{cor}

\begin{proof}
  Note that $\Sigma$ distributes from the left over coproducts, so we have an equivalence
  \begin{equation*}
    \Big(\sm{P:\prop_\UU}P+\neg P\Big)\simeq \Big(\sm{P:\prop_\UU}P\Big)+\Big(\sm{Q:\prop_\UU}\neg Q\Big). 
  \end{equation*}
  Therefore it suffices to show that both $\sm{P:\prop_\UU}P$ and $\sm{Q:\prop_\UU}\neg Q$ are contractible. At the centers of contraction we have $(\unit,\ttt)$ and $(\emptyt,\idfunc)$, respectively. For the contractions, note that both types are subtypes of the types of propositions. Therefore it suffices to show that $\unit=P$ for any proposition $P$ equipped with $p:P$, and that $\emptyt=Q$ for any proposition $Q$ equipped with $q:\neg Q$. Both identifications are obtained immediately from propositional extensionality.
\end{proof}
\index{propositional extensionality|)}

\subsection{Univalence implies function extensionality}
One of the first applications of the univalence axiom was Voevodsky's theorem that the univalence axiom on a universe $\UU$ implies function extensionality for types in $\UU$. The proof uses the fact that weak function extensionality implies function extensionality. We will also make use of the following lemma. 

\begin{lem}\label{lem:postcomp-equiv}
  For any equivalence $e:\eqv{X}{Y}$ in a univalent universe $\UU$, and any type $A$, the post-composition map
  \begin{equation*}
    e\circ\blank : (A \to X) \to (A\to Y)
  \end{equation*}
  is an equivalence.
\end{lem}

Note that this statement was also part of \cref{ex:equiv-postcomp}. That exercise is solved using function extensionality. However, since our present goal is to derive function extensionality from the univalence axiom, we cannot make use of that exercise. Therefore we give a new proof, using the univalence axiom.

\begin{proof}
  Since $\UU$ is assumed to be a univalent universe, it satisfies by \cref{thm:univalence} the principle of equivalence induction. Therefore, it suffices to show that the post-composition map
  \begin{equation*}
    \idfunc\circ\blank : (A\to X)\to (A\to X)
  \end{equation*}
  is an equivalence. This post-composition map is of course just the identity map on $A\to X$, so it is indeed an equivalence.
\end{proof}

\begin{thm}\label{thm:funext-univalence}\index{univalence axiom!implies function extensionality}\index{function extensionality!univalence implies function extensionality}
  For any universe $\UU$, the univalence axiom on $\UU$ implies function extensionality on $\UU$.
\end{thm}

\begin{proof}
  Note that by \cref{thm:funext_wkfunext}\index{weak function extensionality} it suffices to show that univalence implies weak function extensionality. We note that the proof of \cref{thm:funext_wkfunext} also goes through when it is restricted to types in $\UU$.
  
Suppose that $B:A\to \UU$ is a family of contractible types. Our goal is to show that the product $\prd{x:A}B(x)$ is contractible.
Since each $B(x)$ is contractible, the projection map $\proj 1:\big(\sm{x:A}B(x)\big)\to A$ is an equivalence by \cref{ex:proj_fiber}.

Now it follows by \cref{lem:postcomp-equiv} that $\proj1\circ\blank$ is an equivalence. Consequently, it follows from \cref{thm:contr_equiv} that the fibers of
\begin{equation*}
\proj 1\circ\blank : \Big(A\to \sm{x:A}B(x)\Big)\to (A\to A)
\end{equation*}
are contractible. In particular, the fiber at $\idfunc[A]$ is contractible. Therefore it suffices to show that $\prd{x:A}B(x)$ is a retract of $\sm{f:A\to\sm{x:A}B(x)}\proj 1\circ f=\idfunc[A]$. In other words, we will construct a section-retraction pair
\begin{equation*}
\begin{tikzcd}[column sep=1em]
\Big(\prd{x:A}B(x)\Big) \arrow[r,"i"] & \Big(\sm{f:A\to\sm{x:A}B(x)}\proj 1\circ f=\idfunc[A]\Big) \arrow[r,"r"] & \Big(\prd{x:A}B(x)\Big),
\end{tikzcd}
\end{equation*}
with $H:r\circ i\htpy \idfunc$.

We define the function $i$ by
\begin{equation*}
  i(f) \defeq (\lam{x}(x,f(x)),\refl{\idfunc}).
\end{equation*}
To see that this definition is correct, we need to know that
\begin{equation*}
  \lam{x}\proj 1(x,f(x))\jdeq \idfunc.
\end{equation*}
This is indeed the case, by the rule $\lambda$-eq for $\Pi$-types, on \cpageref{page:lambda-eq}.

Next, we define the function $r$. Consider a function $h:A\to \sm{x:A}B(x)$ equipped with an identification $p:\proj 1 \circ h = \idfunc$. Then we have the homotopy $\htpyeq(p):\proj 1 \circ h \htpy \idfunc$. Furthermore, we obtain $\proj 2(h(x)):B(\proj 1(h(x)))$. Using these ingredients, we define $r$ by
\begin{equation*}
  r((h,p),x)\defeq \tr_B(\htpyeq(p,x),\proj 2(h(x))).
\end{equation*}

It remains to construct a homotopy $H:r\circ i\htpy \idfunc$. We simply compute
\begin{align*}
  r(i(f)) & \jdeq r(\lam{x}(x,f(x)),\refl{}) \\
          & \jdeq \tr_B(\htpyeq(\refl{},x),\proj 2(x,f(x))) \\
          & \jdeq \tr_B(\refl{},f(x)) \\
          & \jdeq f(x).
\end{align*}
Thus we see that $r\circ i\jdeq \idfunc$ by an application of the $\eta$-rule for $\Pi$-types. Therefore we simply define $H(f)\defeq\refl{}$.
\end{proof}

\subsection{Maps and families of types}

Using the univalence axiom, we can establish a fundamental relation between maps into a type $A$, and families of types indexed by $A$. A special case of this relation asserts that the type of all pairs $(X,e)$ consisting of a type $X$ and an embedding $e:X\emb A$ is equivalent to the type of all subtypes of $A$, i.e., the type of all families $P$ of propositions indexed by $A$.

\begin{thm}\label{thm:object-classifier}\index{type family}
  For any type $A$ and any univalent universe $\UU$ containing $A$, the map
  \begin{equation*}
    \Big(\sm{X:\UU}X\to A\Big)\to (A\to \UU)
  \end{equation*}
  given by $(X,f)\mapsto\fibf{f}$ is an equivalence.\index{fib f b@{$\fib{f}{b}$}}
\end{thm}

\begin{proof}
  The map in the converse direction is given by
  \begin{equation*}
    B\mapsto \Big(\sm{x:A}B(x),\proj 1\Big). 
  \end{equation*}
  To verify that this map is a section of the asserted map, we have to prove that
  \begin{equation*}
    \fibf{\proj 1}=B
  \end{equation*}
  for any $B:A\to\UU$. By function extensionality and the univalence axiom, this is equivalent to
  \begin{equation*}
    \prd{x:A}\fib{\proj 1}{x}\simeq B(x).
  \end{equation*}
  Such a family of equivalences was constructed in \cref{ex:proj_fiber}.

  It remains to verify that
  \begin{equation*}
    (X,f)=\Big(\sm{x:A}\fib{f}{x},\proj 1\Big). 
  \end{equation*}
  Before we do this, we claim that the identity type
  \begin{equation*}
    (X,f)=(Y,g)
  \end{equation*}
  in the type $\sm{X:\UU}X\to A$ is equivalent to the type of pairs $(e,f)$ consisting of an equivalence $e:X\simeq Y$ equipped with a homotopy $f\htpy g\circ e$. This fact follows from \cref{thm:id_fundamental}, because the type
  \begin{equation*}
    \sm{Y:\UU}\sm{g:Y\to A}\sm{e:X\simeq Y}f\htpy g\circ e
  \end{equation*}
  is contractible by the structure identity principle, \cref{thm:structure-identity-principle}.

  To finish the proof, it therefore suffices to construct an equivalence
  \begin{equation*}
    e:X\simeq \sm{a:A}\fib{f}{a}
  \end{equation*}
  equipped with a homotopy $f\htpy \proj 1\circ e$. Such an equivalence $e$ equipped with a homotopy was constructed in \cref{ex:fib_replacement}.
\end{proof}

The following corollary is so important, that we call it again a theorem.

\begin{thm}\label{thm:classifier-subuniverse}
  Consider a type $A$ and a univalent universe $\UU$ containing $A$. Furthermore, let $P$ be a family of types indexed by $\UU$, and write
  \begin{equation*}
    \UU_P\defeq\sm{X:\UU}P(X).
  \end{equation*}
  Then the map
  \begin{equation*}
    \Big(\sm{X:\UU}\sm{f:X\to A}\prd{a:A}P(\fib{f}{a})\Big)\to (A\to\UU_P)
  \end{equation*}
  given by $(X,f,p)\mapsto \lam{a}(\fib{f}{a},p(a))$ is an equivalence.
\end{thm}

\begin{proof}
  The asserted map is homotopic to the composition of the equivalences
  \begin{align*}
    & \sm{X:\UU}\sm{f:X\to A}\prd{a:A}P(\fib{f}{a}) \\
    & \simeq \sm{(X,f):\sm{X:\UU}X\to A}\prd{a:A}P(\fib{f}{a}) \\
    & \simeq \sm{B:A\to \UU}\prd{a:A}P(B(a)) \\
    & \simeq A\to\sm{X:\UU}P(X).\qedhere
  \end{align*}
\end{proof}

\cref{thm:classifier-subuniverse} applies to any subuniverse. Examples include the subuniverse of $k$-types, for any truncation level $k$, the subuniverse of decidable propositions, the subuniverse of finite types, the subuniverse of inhabited types, and so on. It also applies to type families over $\UU$ that aren't families of propositions. The families $P\defeq\isdecidable$ and $P\defeq\cnt$ are examples.

\begin{cor}\label{cor:subtype}
  Consider a type $A$ and a univalent universe $\UU$ containing $A$. Then the map\index{embedding}\index{subtype}\index{fib f b@{$\fib{f}{b}$}}
  \begin{equation*}
    \Big(\sm{X:\UU}X\emb A\Big)\to (A\to\prop_\UU)
  \end{equation*}
  given by $(X,f)\mapsto \fibf{f}$ is an equivalence.
\end{cor}

In other words, a subtype of a type $A$ is equivalently described as a type $X$ equipped with an embedding $e:X\hookrightarrow A$. This brings us to an important point about equality of subtypes.

\begin{rmk}
  By function extensionality and propositional extensionality, it follows that two subtypes $P,Q:A\to\prop_\UU$ are the same if and only if\index{characterization of identity type!of subtypes of A@{of subtypes of $A$}}\index{subtype!characterization of identity type}
\begin{equation*}
  P(a)\leftrightarrow Q(a)
\end{equation*}
holds for all $a:A$. In other words, two subtypes of $A$ are the same if and only if they contain the same elements of $A$.

On the other hand, by \cref{cor:subtype} we can also consider two types $X$ and $Y$ equipped with embeddings $f:X\emb A$ and $g:Y\emb A$ as subtypes of $A$. Using the structure identity principle, \cref{thm:structure-identity-principle}, we see that the identity type $(X,f)=(Y,g)$ in the type $\sm{X:\UU}X\emb A$ is equivalent to the type
\begin{equation*}
  \sm{e:X\simeq Y}f\htpy g\circ e.
\end{equation*}
In other words, two subtypes $(X,f)$ and $(Y,g)$ of $A$ are equal if and only if there is an equivalence $X\simeq Y$ that is compatible with the embeddings $f:X\emb A$ and $g:Y\emb X$. Indeed, this condition is equivalent to the previous condition that two subtypes are the same if and only if they have the same elements.

We see that the combination of the structure identity principle\index{structure identity principle} and the univalence axiom automatically characterizes equality of subtypes in the most natural way, and we will see similar natural characterizations of identity types throughout the remainder of this book.
\end{rmk}

\subsection{Classical mathematics with the univalence axiom}

In classical mathematics, the axiom of choice asserts that for any collection $X$ of nonempty sets, there is a choice function $f$ such that $f(x)\in x$ for each $x\in X$. The univalence axiom is consistent with the axiom of choice, but we have to be careful in our formulation of the axiom of choice to make it about sets. A naive interpretation that would be applicable to all types, such as the assertion that every family $B$ of inhabited types has a section, is not consistent with univalence. We will use the type $\BS_2$ of $2$-element types for a counterexample.

\begin{prp}\label{prp:Eq-F2}
  The type\index{BS 2@{$\BS_2$}!characterization of identity type}\index{characterization of identity type!of BS 2@{of $\BS_2$}}
  \begin{equation*}
    \sm{X:\BS_2}X
  \end{equation*}
  of pointed $2$-element types is contractible. Consequently, the canonical family of maps
  \begin{equation*}
    (\Fin{2}= X) \to X
  \end{equation*}
  indexed by $X:\BS_2$, is a family of equivalences.
\end{prp}

\begin{proof}
  By the univalence axiom it follows that the type $\sm{X:\BS_2}\Fin{2}\simeq X$ is contractible. In order to show that $\sm{X:\BS_{2}}X$ is contractible, it therefore suffices to show that the map
  \begin{equation*}
    f: (\Fin{2}\simeq X)\to X
  \end{equation*}
  given by $f(e)\defeq e(\star)$, is an equivalence. Since being an equivalence is a proposition by \cref{ex:isprop_isequiv}, we may assume an equivalence $\alpha:\Fin{2}\simeq X$, and we proceed by equivalence induction on $\alpha$. Therefore, it suffices to show that the map
  \begin{equation*}
    f: (\Fin{2}\simeq \Fin{2})\to\Fin{2}
  \end{equation*}
  give by $f(e)\defeq e(\star)$ is an equivalence. Using the notation from \cref{sec:Fin}, we define the inverse map $g$ by
  \begin{align*}
    g(\star) & \defeq \idfunc \\
    g(i(\star)) & \defeq \succFin_2,
  \end{align*}
  and it is a straightforward verification that $f$ and $g$ are inverse to each other.
\end{proof}

\begin{cor}\label{cor:no-section-F2}
  There is no dependent function\index{BS 2@{$\BS_2$}!is not contractible}
  \begin{equation*}
    \prd{X:\BS_2}X.
  \end{equation*}
\end{cor}

\begin{proof}
  By \cref{prp:Eq-F2,ex:equiv-pi}, we have an equivalence
  \begin{equation*}
    \Big(\prd{X:\BS_2}\Fin{2}=X\Big)\simeq \Big(\prd{X:\BS_2}X\Big).
  \end{equation*}
  Note that $\prd{X:\BS_2}\Fin{2}=X$ is the type of contractions of $\BS_2$, using the center of contraction $\Fin{2}$. Therefore it suffices to show that $\BS_2$ is not contractible. Recall from \cref{ex:prop_contr} that the identity types of contractible types are contractible. On the other hand, it follows from \cref{prp:Eq-F2} that the identity type $\Fin{2}=\Fin{2}$ in $\BS_2$ is equivalent to $\Fin{2}$. This type isn't contractible by \cref{ex:is-not-contractible-Fin}. We conclude that $\BS_2$ is not contractible.
\end{proof}

The family $X\mapsto X$ over $\BS_2$ is therefore an example of a family of nonempty types for which there are provably no sections. In the following corollary we conclude more generally that there is no way to construct an element of an arbitrary inhabited type.

\begin{cor}\label{cor:no-global-choice}
  If $\UU$ is a univalent universe, then there is no \define{global choice}\index{global choice|textbf} function
  \begin{equation*}
    \prd{A:\UU}\brck{A}\to A.
  \end{equation*}
\end{cor}

\begin{proof}
  Suppose $f:\prd{A:\UU}\brck{A}\to A$. By restricting $f$ to the type of $2$-element types in $\UU$, we obtain a function
  \begin{equation*}
    \prd{A:\BS_2}\brck{A}\to A.
  \end{equation*}
  Note that every $2$-element type is inhabited, i.e., there is an element of type $\brck{A}$ for every $2$-element type $A$. To see this, consider a type $A:\UU$ such that $\brck{\Fin{2}\simeq A}$. To obtain an element of type $\brck{A}$, we may assume an equivalence $e:\Fin{2}\simeq A$. Then we have $\eta(e(0)):\brck{A}$.

  Since every $2$-element type is inhabited, we obtain a function $\prd{A:\BS_2}A$, which is impossible by \cref{cor:no-section-F2}.
\end{proof}

\cref{cor:no-global-choice} is of philosophical importance. It shows that the \define{principle of global choice} is incompatible with the univalence axiom, i.e., that there is no way to obtain construct a function $\brck{A}\to A$ for all types $A$. In other words, we cannot obtain an element of $A$ merely from the assumption that the type $A$ is inhabited. What is the obstruction? It is the fact that no such choice of an element of $A$ can be invariant under the automorphisms on $A$, i.e., under the self-equivalences on $A$. Indeed, in the example where $A$ is the $2$-element type $\Fin{2}$ there are no fixed point of the equivalence $\succFin_2:\Fin{2}\simeq\Fin{2}$. By the univalence axiom, there is an identification $p:\Fin{2}=\Fin{2}$ in $\UU$, such that $\tr(p,x)=\succFin_2(x)$. If we had a function
\begin{equation*}
  f:\prd{X:\UU}\brck{X}\to X,
\end{equation*}
the dependent action on paths of $f$ would give an identification
\begin{equation*}
  \apd{f}{p}:\succFin_2(f(\Fin{2},p,H))=f(\Fin{2},p,\eta(0)).
\end{equation*}
In other words, it would give us a fixed point for the successor function on $\Fin{2}$.

This is perhaps a good moment to stress that the axiom of choice is really an axiom about \emph{sets}, not about more general types. And indeed, when we restrict the axiom of choice to sets, it turns out to be consistent with the univalence axiom and therefore safe to assume. In this book, however, we will not have many applications for the axiom of choice and therefore we will not assume it, unless we explicitly say otherwise.

\begin{defn}
  The \define{axiom of choice}\index{axiom of choice|textbf}\index{axiom!axiom of choice|textbf} asserts that for any family $B$ of inhabited sets indexed by a set $A$, the type of sections of $B$ is also inhabited, i.e., it asserts that there is an element of type
  \begin{equation*}
    \AC_{\UU}(A,B)\defeq \Big(\prd{x:A}\brck{B(x)}\Big)\to\Brck{\prd{x:A}B(x)},
  \end{equation*}
  for every $A:\Set_\UU$ and $B:A\to\Set_\UU$.
\end{defn}

Similar care has to be taken with the type theoretic formulation of the law of excluded middle. It is again inconsistent to assume that every type is decidable.

\begin{thm}
  There is no \define{global decidability function}\index{global decidability|textbf}
  \begin{equation*}
    \prd{X:\UU}\isdecidable(X). 
  \end{equation*}
\end{thm}

\begin{proof}
  Suppose there is such a dependent function $d$. By restricting $d$ to the subuniverse of $2$-element types, we obtain a dependent function
  \begin{equation*}
    d:\prd{X:\BS_2}\isdecidable(X).
  \end{equation*}
  However, each $2$-element type $X$ is inhabited. By \cref{ex:propositional-truncations-drill} we obtain a function
  \begin{equation*}
    \isdecidable(X)\to X
  \end{equation*}
  for each $2$-element type $X$. Therefore, we obtain from $d$ a dependent function $\prd{X:\BS_2}X$, which does not exist by \cref{cor:no-section-F2}.
\end{proof}

The law of excluded middle is really an axiom of propositional logic, and it is indeed consistent with the univalence axiom that every \emph{proposition} is decidable.

\begin{defn}
  The \define{law of excluded middle}\index{law of excluded middle|textbf}\index{axiom!law of excluded middle|textbf} asserts that every proposition is decidable, i.e.,
  \begin{equation*}
    \LEM_\UU\defeq \prd{P:\prop_\UU}\isdecidable(P).
  \end{equation*}
\end{defn}

We will again not assume the law of excluded middle, unless we say otherwise. Nevertheless, we have seen in \cref{sec:decidability} that some propositions are already decidable without assuming the law of excluded middle, and decidability remains an important concept in type theory and mathematics.

\subsection{The binomial types}
\index{binomial type|(}

To wrap up this section on univalence, we will use the univalence axiom to construct for any two types $A$ and $B$ a type $\dbinomtype{A}{B}$ that has properties similar to the binomial coefficients $\dbinomtype{n}{k}$. Indeed, we will show that if $A$ is an $n$-element type and $B$ is a $k$-element type, then $\dbinomtype{A}{B}$ is an $\binom{n}{k}$-element type. The binomial types are defined using decidable embeddings.

\begin{defn}
  A map $f:A\to B$ is said to be \define{decidable}\index{decidable map|textbf} if it comes equipped with an element of type
  \begin{equation*}
    \isdecidable(f) \defeq \prd{b:B}\isdecidable(\fib{f}{b}).
  \end{equation*}
  We will write $A\demb B$ for the type of \define{decidable embeddings}\index{decidable embedding|textbf}\index{A hookrightarrow d B@{$A \demb B$}|see {decidable embedding}} from $A$ to $B$, i.e., for the type of embeddings that are also decidable maps. 
\end{defn}

\begin{defn}
  Consider a type $A$ and a universe $\UU$. We define the \define{connected component}\index{universe!connected component|textbf}\index{connected component!of a universe|textbf} of $\UU$ at $A$ by\index{U A@{$\UU_A$}|textbf}
  \begin{equation*}
    \UU_A\defeq \sm{X:\UU}\brck{A\simeq X}.
  \end{equation*}
\end{defn}

\begin{eg}
  Note that type $\UU_{\Fin{n}}$ is the type $\BS_n$ of all $n$-element types. Note also that if $A\simeq B$, then $\UU_A\simeq\UU_B$.\index{BS n@{$\BS_n$}}
\end{eg}

\begin{defn}\label{defn:binomial-type}
  Consider two types $A$ and $B$ and a universe $\UU$ containing both $A$ and $B$. We define the \define{binomial type}\index{binomial type|textbf} $\dbinomtype[\UU]{A}{B}$\index{(A B)@{$\dbinomtype{A}{B}$}|see {binomial type}} by
  \begin{equation*}
    \dbinomtype[\UU]{A}{B} \defeq \sm{X:\UU_B}X\demb A.
  \end{equation*}
\end{defn}

\begin{rmk}
  We define the binomial types using decidable embeddings because the usual properties of binomial coefficients generalize most naturally under the extra assumption of decidability. In particular the binomial theorem for types, which is stated as \cref{ex:binomial-theorem} and generalized in \cref{ex:distributive-pi-coprod}, rely on the use of decidable embeddings.
\end{rmk}

\begin{prp}\label{prp:equiv-binom-type}
  Consider two types $A$ and $B$, and a universe $\UU$ containing both $A$ and $B$. Then we have an equivalence
  \begin{equation*}
    \dbinomtype[\UU]{A}{B}\simeq \sum_{(P:A\to\decidableProp_\UU)}\Brck{B\simeq\sm{a:A}P(a)}.
  \end{equation*}
  from the binomial type $\dbinomtype[\UU]{A}{B}$ to the type of decidable subtypes\index{decidable subtype} of $A$ that are merely equivalent to $B$.
\end{prp}

\begin{proof}
  This equivalence follows from \cref{thm:classifier-subuniverse}, by which we have
  \begin{equation*}
    \Big(\sm{X:\UU}X\demb A\Big)\simeq (A\to\decidableProp_\UU).\qedhere
  \end{equation*}
\end{proof}

\begin{rmk}
  Combining \cref{cor:decidable-Prop,prp:equiv-binom-type}, we obtain an equivalence
  \begin{equation*}
    \dbinomtype[\UU]{A}{B}\simeq \sum_{(f:A\to\bool)}\Brck{B\simeq\sm{a:A}f(a)=\btrue}.
  \end{equation*}
  for any universe $\UU$ that contains both $A$ and $B$. This equivalence is important, because the right hand side doesn't depend on the universe $\UU$. Therefore we will simply write $\dbinomtype{A}{B}$ for $\dbinomtype[\UU]{A}{B}$, if the universe $\UU$ contains both $A$ and $B$.
\end{rmk}
  
\begin{lem}\label{prp:binomtype-recursion}
  For any two types $A$ and $B$, we have equivalences\index{binomial type!recursive relations}
  \begin{align*}
    \dbinomtype{\emptyt}{\emptyt} & \simeq \unit & \dbinomtype{A+\unit}{\emptyt} & \simeq \unit \\
    \dbinomtype{\emptyt}{B+\unit} & \simeq \emptyt & \dbinomtype{A+\unit}{B+\unit} & \simeq \dbinomtype{A}{B}+\dbinomtype{A}{B+\unit}.
  \end{align*}
\end{lem}

\begin{proof}
  For the first two equivalences, we prove that $\dbinomtype{X}{\emptyt}$ is contractible for any type $X$. To see this, we first note that the type $\UU_\emptyt$ is contractible. Indeed, at the center of contraction we have the empty type, and any two types that are merely equivalent to the empty type are empty and hence equivalent. Therefore it follows that
  \begin{equation*}
    \dbinomtype{X}{\emptyt}\simeq \emptyt\demb X.
  \end{equation*}
  The type of decidable embeddings $\emptyt\demb X$ is contractible, because the type $\emptyt\to X$ is contractible with the map $\exfalso:\emptyt\to X$ at the center of contraction, which is of course a decidable embedding.

  Next, the fact that the binomial type $\dbinomtype{\emptyt}{B+\unit}$ is empty follows from the fact that the type of maps $X\to\emptyt$ is empty for any type $X$ merely equivalent to $B+\unit$. 

  For the last equivalence we will use \cref{prp:equiv-binom-type}. Using the universal property of $A+\unit$, we see that
  \begin{equation*}
    \dbinomtype{A+\unit}{B+\unit}\simeq \sm{P:A\to\decidableProp_\UU}\sm{Q:\decidableProp_\UU}\brck{(B+\unit)\simeq \sm{a:A}P(a)+ Q}.
  \end{equation*}
  Using the fact that $\decidableProp_\UU\simeq\Fin{2}$, observe that we have an equivalence
  \begin{align*}
    & \sm{Q:\decidableProp_\UU}\brck{(B+\unit)\simeq(\sm{a:A}P(a)+Q)} \\
    & \qquad\simeq \brck{(B+\unit)\simeq(\sm{a:A}P(a)+\unit)}+\brck{(B+\unit)\simeq\sm{a:A}P(a)}.
  \end{align*}
  Furthermore, note that it follows from \cref{prp:is-injective-maybe} that
  \begin{equation*}
    \brck{(B+\unit)\simeq(\sm{a:A}P(a)+\unit)}\simeq\brck{B\simeq\sm{a:A}P(a)}.
  \end{equation*}
  Thus we see that
  \begin{align*}
    \dbinomtype{A+\unit}{B+\unit} & \simeq \Big(\sm{P:A\to\decidableProp_\UU}\brck{B\simeq\sm{a:A}P(a)}\Big) \\
    & \qquad\qquad +\Big(\sm{P:A\to\decidableProp_\UU}\brck{(B+\unit)\simeq\sm{a:A}P(a)}\Big).\qedhere
  \end{align*}
\end{proof}

\begin{thm}
  If $A$ and $B$ are finite types of cardinality $n$ and $k$, respectively, then the type $\dbinomtype{A}{B}$ is finite of cardinality $\binom{n}{k}$.\index{finite type!binomial type}\index{binomial type!is finite}
\end{thm}

\begin{proof}
  The claim that the type $\binomtype{A}{B}$ is finite of cardinality $\binom{n}{k}$ is a proposition, so we may assume $e:\Fin{n}\simeq A$ and $f:\Fin{k}\simeq B$. The claim now follows by induction on $n$ and $k$, using \cref{prp:binomtype-recursion}.
\end{proof}

\begin{rmk}
  It is perhaps remarkable that the type $\sm{X:\UU_B}X\demb A$ is a good generalisation of the binomial coefficients to types. Note that when $A$ and $B$ are finite types of cardinality $n$ and $k$, respectively, then the type $B\demb A$ has a factor $k!$ too many elements. When we seemingly enlarge it by the type $\UU_B$ of all types merely equivalent to $B$, it turns out that we obtain the correct generalisation of the binomial coefficients.

  One reason why it works is that the identity type of $\sm{X:\UU_B}X\demb A$ is characterized, via the univalence axiom, by
  \begin{equation*}
    ((X,f)=(Y,g))\simeq\sm{e:X\simeq Y}f\htpy g\circ e. 
  \end{equation*}
  Therefore it follows that for any two decidable embeddings $f,g:B\demb A$, if $f$ and $g$ are the same up to a permutation on $B$, then we get an identification $(B,f)=(B,g)$ in the type $\sm{X:\UU_B}X\demb A$.

  From a group theoretic perspective we may observe that the automorphism group $B\simeq B$ acts freely on the set of decidable embeddings $\B\demb A$, and the type $\sm{X:\UU_B}X\hookrightarrow A$ can be viewed as the type of orbits of that action. Since this action of $\Aut(B)$ on $B\demb A$ is free, we see that the number of orbits is $\frac{1}{k!}$ times the number of elements in $B\demb A$. 
\end{rmk}
\index{binomial type|)}

\begin{exercises}
  \exitem \label{ex:istrunc_UUtrunc}
  \begin{subexenum}
  \item Use the univalence axiom to show that the type $\sm{A:\UU}\iscontr(A)$ of all contractible types in $\UU$ is contractible.\index{universe!of contractible types}\index{contractible type}
  \item Use the univalence axiom and \cref{ex:isprop_istrunc,ex:isprop_isequiv} to show that the universe of $k$-types\index{universe!of k-types@{of $k$-types}}\index{U leq k@{$\UU^{\leq k}$}}\index{k-type@{$k$-type}!universe of k-types@{universe of $k$-types}}\index{truncated type!universe of k-types@{universe of $k$-types}}\index{universe!of k-types@{of $k$-types}!is a k-type@{is a $k$-type}}\index{U leq k@{$\UU^{\leq k}$}!is a k-type@{is a $k$-type}}
    \begin{equation*}
      \UU^{\leq k}\defeq \sm{X:\UU}\istrunc{k}(X)
    \end{equation*}
    is a $(k+1)$-type, for any $k\geq -2$.
  \item Show that $\prop_\UU$ is not a proposition.\index{universe!of propositions}\index{Prop@{$\prop_\UU$}}
  \item Show that the universe $\Set_\UU$ of sets\index{universe!of sets} in $\UU$ is not a set. \index{Set@{$\Set_\UU$}!is not a set}
  \end{subexenum}
  \exitem Give an example of a type family $B$ over a type $A$ for which the implication
  \begin{equation*}
    \neg\Big(\prd{x:A}B(x)\Big) \to \Big(\sm{x:A}\neg B(x)\Big)
  \end{equation*}
  is false.
  \exitem Show that the law of excluded middle holds if and only if every set has decidable equality.\index{law of excluded middle}\index{axiom!law of excluded middle}
  \exitem Consider a type $A$ and a univalent universe $\UU$ containing $A$. Construct an equivalence
  \begin{equation*}
    A\simeq\sm{B:A\to\UU}\iscontr\left(\sm{a:A}B(a)\right).
  \end{equation*}
  \exitem \label{ex:surjective-precomp}Consider a map $f:A\to B$. Show that the following are equivalent:
  \begin{enumerate}
  \item The map $f$ is surjective.\index{surjective map}
  \item For every set $C$, the precomposition function
    \begin{equation*}
      \blank\circ f:(B\to C)\to (A\to C)
    \end{equation*}
    is an embedding.
  \end{enumerate}
  Hint: To show that (ii) implies (i), use the assumption with the set $C\defeq\prop_\UU$, where $\UU$ is a univalent universe containing both $A$ and $B$.
  \exitem (Escard\'o)\label{ex:idtype-is-emb} Consider a type $A$ in $\UU$. Show that the identity type, seen as a function\index{Id A@{$\idtypevar{A}$}!is an embedding}\index{is an embedding!Id A@{$\idtypevar{A}$}}
  \begin{equation*}
    \idtypevar{} : A \to (A\to\UU),
  \end{equation*}
  is an embedding.
  \exitem
  \begin{subexenum}
  \item For any type $A$ in $\UU$, consider the function\index{S A@{$\Sigma_A$}|see {dependent pair type}}\index{S A@{$\Sigma_A$}!is a k-truncated map@{is a $k$-truncated map}}\index{truncated map!S A@{$\Sigma_A$}}
    \begin{equation*}
      \Sigma_A : (A \to \UU) \to \UU,
    \end{equation*}
    which takes a family $B$ of $\UU$-small types to its $\Sigma$-type. Show that the following are equivalent:
    \begin{enumerate}
    \item The type $A$ is $k$-truncated.
    \item The map $\Sigma_A$ is $k$-truncated.
    \end{enumerate}
    Hint: Construct an equivalence $\fib{\Sigma_A}{X}\simeq (X\to A)$.\index{fiber!of S A@{of $\Sigma_A$}}
  \item Show that the map ${+}:\UU\times\UU\to\UU$, which takes $(A,B)$ to the coproduct $A+B$, is $0$-truncated.\index{A + B@{$A+B$}!is a 0-truncated map@{${\blank}+{\blank}$ is a $0$-truncated map}}\index{truncated map!-+-@{${\blank}+{\blank}$}}
  \end{subexenum}
  \exitem (Escard\'o) Consider a proposition $P$ and a universe $\UU$ containing $P$. Show that the map
  \begin{equation*}
    \Pi_P : (P\to \UU)\to\UU,
  \end{equation*}
  given by $A\mapsto\prd{p:P}A(p)$, is an embedding.\index{P P@{$\Pi_P$}!is an embedding}
  \exitem Consider two types $A$ and $B$ and a universe $\UU$ containing both $A$ and $B$. A \define{binary correspondence}\index{binary correspondence|textbf} $R:A\to(B\to\UU)$ is said to be a \define{function}\index{binary correspondence!function|textbf}\index{function!binary correspondence|textbf} if it satisfies the condition\index{is-function@{$\isfunction(R)$}|textbf}
  \begin{equation*}
    \isfunction(R)\defeq\prd{a:A}\iscontr\Big(\sm{b:B}R(a,b)\Big),
  \end{equation*}
  and $R$ is said to be an \define{opposite function}\index{function!opposite function|textbf}\index{opposite function|textbf} if the \define{opposite correspondence}\index{correspondence!opposite correspondence}\index{opposite correspondence}\index{R op@{$\op{R}$}|see {opposite correspondence}} $\op{R}:B\to(A\to\UU)$ given by $\op{R}(b,a)\defeq R(a,b)$ is functional.
  \begin{subexenum}
  \item Construct an equivalence
    \begin{equation*}
      (A\to B)\simeq \sm{R:A\to(B\to\UU)}\isfunction(R).
    \end{equation*}
  \item Construct an equivalence
    \begin{equation*}
      (A\simeq B)\simeq\sm{R:A\to (B\to\UU)}\isfunction(R)\times\isfunction(\op{R}).
    \end{equation*}
  \end{subexenum}
  \exitem
  \begin{subexenum}
  \item For any $k:\N$, show that the type\index{BS n@{$\BS_n$}}
    \begin{equation*}
      \sm{X:\BS_{k+1}}\Fin{k}\hookrightarrow X
    \end{equation*}
    is contractible.
  \item More generally, construct for any $k,l:\N$ and any $k$-element type $A$ an equivalence
    \begin{equation*}
      \Big(\sm{X:\BS_{k+l}}A\hookrightarrow X\Big)\simeq \BS_l
    \end{equation*}
  \end{subexenum}
  \exitem \label{ex:complement-Fk}
  \begin{subexenum}
  \item For any type $A$, construct an equivalence
  \begin{equation*}
    \UU_A \simeq \sum_{(X:\UU_{A+\unit})}\dbinomtype{X}{\unit}.
  \end{equation*}
  \item For any $k:\N$, construct an equivalence\index{BS n@{$\BS_n$}}
    \begin{equation*}
      \Big(\sm{X:\BS_{k+1}}X\Big)\simeq \BS_k.
    \end{equation*}
    In other words, show that the type of $(k+1)$-element types equipped with a point is equivalent to the type of $k$-element types. Conclude that the type of pointed finite types is equivalent to the type of finite types, i.e., conclude that we have an equivalence\index{F@{$\F$}}
  \begin{equation*}
    \Big(\sm{X:\F}X\Big)\simeq \F.
  \end{equation*}
  \end{subexenum}
  \exitem 
  \begin{subexenum}
  \item Show that for $k\neq 2$, the type\index{BS n@{$\BS_n$}}
    \begin{equation*}
      \prd{X:\BS_k}X\to X
    \end{equation*}
    is contractible. Conclude that the type $\prd{X:\BS_k}X\simeq X$ is also contractible. Hint: Use \cref{ex:complement-Fk}.
  \item Show that the type\index{BS 2@{$\BS_2$}}
    \begin{equation*}
      \prd{X:\BS_2}X\to X
    \end{equation*}
    is equivalent to $\Fin{2}$. Conclude that the type $\prd{X:\BS_2}X\simeq X$ is also equivalent to $\Fin{2}$.
  \end{subexenum}
  \exitem Consider a type $A$.
  \begin{subexenum}
  \item Recall from \cref{ex:isolated-point} that an element $a:A$ is isolated\index{isolated point} if and only if the map $\const_a:\unit\to A$ is a decidable embedding. Construct an equivalence
  \begin{equation*}
    \dbinomtype{A}{\unit}\simeq \sm{a:A}\isisolated(a).
  \end{equation*}
  \item Construct an equivalence
    \begin{equation*}
      \dbinomtype{A}{\unit}\simeq\Big(\sm{X:\UU}(X+\unit)\simeq A\Big).
    \end{equation*}
    Conclude that the map $X\mapsto X+\unit$ on a univalent universe $\UU$ is $0$-truncated.\index{truncated map!-+1@{${\blank}+\unit$}}
  \item More generally, construct an equivalence\index{binomial type}
    \begin{equation*}
      \dbinomtype{A}{B} \simeq \sm{X:\UU_B}\sm{Y:\UU}(X+Y\simeq A).
    \end{equation*}
  \end{subexenum}
  \exitem \label{ex:binomial-theorem}For any $(X,i):\dbinomtype{A}{B}$, we define $A\setminus(X,i)\defeq (A\setminus X,A\setminus i):\dbinomtype{A}{B}$, where
  \begin{align*}
    A\setminus X & \defeq \sm{a:A}\neg(\fib{i}{a}) \\
    A\setminus i & \defeq \proj 1.
  \end{align*}
  Now consider a finite type $X$ and two arbitrary types $A$ and $B$. Construct an equivalence\index{binomial theorem}
  \begin{equation*}
    (A+B)^X\simeq\sm{k:\N}\sm{(Y,i):\dbinomtype{X}{\Fin{k}}}A^Y\times B^{X\setminus Y}.
  \end{equation*}
  \exitem Let $\UU$ be a univalent universe.
  \begin{subexenum}
  \item Consider a section-retraction pair
    \begin{equation*}
      \begin{tikzcd}
        A \arrow[r,"i"] & X \arrow[r,"r"] & A
      \end{tikzcd}
    \end{equation*}
    with $H:r\circ i\htpy \idfunc$. Show that if $X$ is $\UU$-small, then so is $A$. Hint: Use \cref{ex:retracts-as-limits}.\index{small type!retracts of small types are small}\index{retract!retracts of small types are small}
  \item Consider two inhabited types $A$ and $B$. Show that if the product $A\times B$ is $\UU$-small, then so are the types $A$ and $B$.\index{small type!products of small types are small}\index{cartesian product type!products of small types are small}
  \end{subexenum}
  \exitem Consider a finite type $X$ and a univalent universe $\UU$ containing $X$. Show that the type\index{Retr X@{$\Retr_\UU(X)$}!is finite}\index{finite type!Retr X of a finite type X@{$\Retr_\UU(X)$ of a finite type $X$}}
  \begin{equation*}
    \Retr_\UU(X)\defeq \sm{A:\UU}\sm{i:A\to X}\sm{r:X\to A}r\circ i\htpy\idfunc
  \end{equation*}
  of all retracts of $X$ is finite.
  \exitem Consider a $k$-truncated type $X$ and a univalent universe $\UU$ containing $X$. Show that the type
  \begin{equation*}
    \mathsf{Retr}_\UU(X)
  \end{equation*}
  of all retracts of $X$ is $k$-truncated.\index{Retr X@{$\Retr_\UU(X)$}!is a k-truncated type@{is a $k$-truncated type}}
  \exitem \label{ex:surjection-into-k-type}For any type $A$ and any $k\geq-1$, show that the type
  \begin{equation*}
    \sm{X:\typele{k}}A\twoheadrightarrow X
  \end{equation*}
  of $k$-truncated types $X$ equipped with a surjective map $A\twoheadrightarrow X$ is $k$-truncated, even though the type $\typele{k}$ itself is $(k+1)$-truncated.
  \exitem
  \begin{subexenum}
    \item Show that for $k\geq 3$, the type\index{BS n@{$\BS_n$}}
    \begin{equation*}
      \prd{X:\BS_k}(X+X)\emb (X\times X)+\unit
    \end{equation*}
    is empty, even though the inequality $2k\leq k^2+1$ holds for all $k:\N$.
  \item Show that the type\index{BS 2@{$\BS_2$}}
    \begin{equation*}
      \prd{X:\BS_2}(X+X)\emb (X\times X)+\unit
    \end{equation*}
    is equivalent to $\Fin{8}$.
  \end{subexenum}
  \exitem \label{ex:prime}For any natural number $n$ consider the type
  \begin{equation*}
    \tilde{D}_n\defeq\sm{X:\BS_2}\sm{Y:X\to\F}\Big(\Fin{n}\simeq\prd{x:X}Y(x)\Big).
  \end{equation*}
  \begin{subexenum}
  \item Show that $\tilde{D}_{1}\simeq \BS_2$.\index{BS 2@{$\BS_2$}}
  \item Show that $\tilde{D}_{n}$ is contractible if and only if $n$ is prime\index{prime number}.
  \item Show that $\tilde{D}_n$ is a set if and only if $n$ is not a square.
  \end{subexenum}
\end{exercises}
\index{univalence axiom|)}
\index{axiom!univalence|)}


\section{Set quotients}\label{sec:set-quotients}
\index{set quotient|(}

In this section we construct the quotient of a type by an equivalence relation. By an equivalence relation we understand a binary relation $R$ which is reflexive, symmetric, transitive, and moreover, we require that the type $R(x,y)$ relating $x$ and $y$ is a proposition. Therefore, if $\UU$ is a universe that contains $R(x,y)$ for each $x,y:A$, then we can view $R$ as a map
\begin{equation*}
  R:A\to(A\to\prop_\UU).
\end{equation*}
The quotient $A/R$ is constructed as the type of equivalence classes, which is just the image of the map $R:A\to (A\to\prop_\UU)$. This construction of the quotient by an equivalence relation is very much like the construction of a quotient set in classical set theory. Examples of set quotients are abundant in mathematics. We cover two of them in this section: the type of rational numbers and the set truncation of a type.

There is, however, a subtle issue with our construction of the set quotient as the image of the map $R:A\to(A\to\prop_\UU)$. What universe is the quotient $A/R$ in? Note that $\prop_\UU$ is a type in the successor universe $\UU^+$, constructed in \cref{defn:successor-universe}. Therefore the function type $A\to \prop_\UU$ as well as the quotient $A/R$ are also types in $\UU^+$. That seems unfortunate, because in Zermelo-Fraenkel set theory the quotient of a set by an equivalence relation is an ordinary set, and not a more general class.

To address the size issues of set quotients, we will introduce the type theoretic replacement axiom. This axiom is analogous to the replacement axiom in Zermelo-Fraenkel set theory, which asserts that the image of a set under any function is again a set. The type theoretic replacement property asserts that for any map $f:A\to B$ from a type $A$ in $\UU$ to a type $B$ of which the \emph{identity types} are equivalent to types in $\UU$, the image of $f$ is also equivalent to a type in $\UU$. The replacement axiom can either be assumed, or it can be proven from the assumption that universes are closed under certain \emph{higher inductive types}, and it is therefore considered to be a very mild assumption.

\subsection{Equivalence relations and the replacement axiom}

\begin{defn}\label{defn:eq_rel}
Consider a type $A$ and a universe $\UU$. Let $R:A\to (A\to\prop_\UU)$ be a binary relation on $A$ valued in the propositions in $\UU$\index{relation!valued in propositions}. We say that $R$ is an \define{equivalence relation}\index{equivalence relation|textbf}\index{relation!equivalence relation|textbf} if $R$ comes equipped with
\begin{align*}
\rho & : \prd{x:A}R(x,x) \\*
\sigma & : \prd{x,y:A} R(x,y)\to R(y,x) \\*
\tau & : \prd{x,y,z:A} R(x,y)\to (R(y,z)\to R(x,z)),
\end{align*}
witnessing that $R$ is reflexive\index{relation!reflexive}, symmetric\index{relation!symmetric}, and transitive\index{relation!transitive}. We write $\eqrel_\UU(A)$\index{Eq-Rel@{$\eqrel_\UU(A)$}|see {equivalence relation}}\index{Eq-Rel@{$\eqrel_\UU(A)$}|textbf} for the type of all equivalence relations on $A$ valued in the propositions in $\UU$.
\end{defn}

\begin{defn}
  Let $R:A\to (A\to\prop_\UU)$ be an equivalence relation. A subtype $P:A\to \prop_\UU$ is said to be an \define{equivalence class}\index{equivalence class|textbf}\index{equivalence relation!equivalence class|textbf} if it satisfies the condition
  \begin{equation*}
    \isequivalenceclass(P)\defeq\exists_{(x:A)}\forall_{(y:A)}P(y)\leftrightarrow R(x,y).
  \end{equation*}
  We define $A/R$\index{A/R@{$A/R$}|see {set quotient}}\index{set quotient|textbf} to be the type of equivalence classes, i.e., we define
  \begin{equation*}
    A/R\defeq \sm{P:A\to\prop_\UU}\isequivalenceclass(P).
  \end{equation*}
  Furthermore, we define \define{equivalence class of $x:A$}\index{[x] R@{$[x]_R$}|textbf}\index{equivalence relation![x] R@{$[x]_R$}|textbf}\index{set quotient![x] R@{$[x]_R$}|textbf} to be
  \begin{equation*}
    [x]_R\defeq R(x),
  \end{equation*}
  which is indeed an equivalence class. Sometimes we will write $q_R:A\to A/R$ for the map $x\mapsto [x]_R$.\index{q R@{$q_R$}|textbf}
\end{defn}

In other words, $A/R$ is the image of the map $R:A\to (A\to\prop_\UU)$. In the following proposition we characterize the identity type of $A/R$. As a corollary, we obtain equivalences
\begin{equation*}
  ([x]_R=[y]_R)\simeq R(x,y),
\end{equation*}
justifying that the quotient $A/R$ is defined to be the type of equivalence classes. Note that in our characterization of the identity type of $A/R$ we make use of propositional extensionality.

\begin{prp}\label{prp:eq-quotient}
  Let $R:A\to (A\to\prop_\UU)$ be an equivalence relation. Furthermore, consider $x:A$ and an equivalence class $P$. Then the canonical map\index{characterization of identity type!of set quotients}\index{set quotient!characterization of identity type}\index{identity type!of A/R@{of $A/R$}}
  \begin{equation*}
    ([x]_R=P)\to P(x)
  \end{equation*}
  is an equivalence.
\end{prp}

\begin{proof}
  By \cref{thm:id_fundamental} it suffices to show that the total space
  \begin{equation*}
    \sm{P:A/R}P(x)
  \end{equation*}
  is contractible. The center of contraction is of course $[x]_R$, which satisfies $[x]_R(x)$ by reflexivity of $R$. It remains to construct a contraction. Since $\sm{P:A/R}P(x)$ is a subtype of $A/R$, we construct a contraction by showing that
  \begin{equation*}
    [x]_R=P
  \end{equation*}
  whenever $P(x)$ holds. Since $P$ is an equivalence class there exists an element $y:A$ such that $P=[y]_R$. Note that our goal is a proposition, so we may assume that we have such a $y$. From the assumption that $P(x)$ holds, it follows that $R(x,y)$ holds. To complete the proof, it therefore is suffices to show that
  \begin{equation*}
    [x]_R=[y]_R,
  \end{equation*}
  assuming that $R(x,y)$ holds. By function extensionality and propositional extensionality, it is equivalent to show that
  \begin{equation*}
    \prd{z:A}R(x,z)\leftrightarrow R(y,z),
  \end{equation*}
  which follows directly from the assumption that $R$ is an equivalence relation.
\end{proof}

\begin{cor}\label{cor:eq-quotient}
  Consider an equivalence relation $R$ on a type $A$, and let $x,y:A$. Then there is an equivalence
  \begin{equation*}
    ([x]_R=[y]_R)\simeq R(x,y).
  \end{equation*}
\end{cor}

\begin{rmk}
  Notice that type of equivalence classes of an equivalence relation in $\UU$ is a type in the universe $\UU^+$ that contains $\UU$ and every type in $\UU$, or indeed in any universe $\VV$ containing $\UU$ and every type in $\UU$. Indeed, the type
  \begin{equation*}
    \prop_\UU\jdeq\sm{X:\UU}\isprop(X)
  \end{equation*}
  of propositions in $\UU$ is a type in $\VV$. It follows that the type $A\to\prop_\UU$ is a type in $\VV$. The type of equivalence classes of an equivalence relation $R$ on $A$ in $\UU$ is a subtype of $A\to\prop_\UU$ in $\UU$, so we conclude that $A/R$ is a type in $\VV$.
\end{rmk}

In classical mathematics, on the other hand, we consider the class of equivalence classes of an equivalence relation to be a (small) set. We will introduce the replacement axiom in order to ensure that set quotients in type theory are small.

Recall that in set theory, the replacement axiom asserts that for any family of sets $\{X_i\}_{i\in I}$ indexed by a set $I$, there is a set $X[I]$ consisting of precisely those sets $x$ for which there exists an $i\in I$ such that $x\in X_i$. In other words: the image of a set-indexed family of sets is again a set. Without the replacement axiom, $X[I]$ would be a class.

In type theory, we may similarly ask whether the image of a map $X:I\to\UU$ is $\UU$-small, assuming that $I$ is $\UU$-small. The replacement axiom settles a more general variant of this question. The key observation is that the identity types of $\UU$ are $\UU$-small by the univalence axiom. In other words, univalent universes are \emph{locally small} in the following sense.

\begin{defn}\label{defn:locally-small-type}
  Consider a universe $\UU$. A type $A$ is said to be \define{locally $\UU$-small}\index{locally small type|textbf} if the identity type $x=y$ is $\UU$-small for every $x,y:A$.
    We write\index{is-locally-small(A)@{$\islocallysmall_\UU(A)$}}
    \begin{equation*}
      \islocallysmall_\UU(A)\defeq \prd{x,y:A}\issmall_\UU(x=y).
    \end{equation*}
    Similarly, a map $f:A\to B$ is said to be \define{locally $\UU$-small}\index{locally small map|textbf} if all of its fibers are locally $\UU$-small.
\end{defn}

\begin{eg}
  ~
  \begin{enumerate}
  \item Any $\UU$-small type is also locally $\UU$-small.\index{small type!small types are locally small}\index{locally small type!small types are locally small}
  \item Any proposition is locally small with respect to any universe $\UU$.\index{proposition!propositions are locally U-small@{propositions are locally $\UU$-small}}\index{locally small type!propositions are locally small}
  \item Any univalent universe $\UU$ is locally $\UU$-small, because by the univalence axiom we have an equivalence\index{locally small type!U is locally U-small@{$\UU$ is locally $\UU$-small}}\index{universe!U is locally U-small@{$\UU$ is locally $\UU$-small}}
    \begin{equation*}
      (A=B)\simeq (A\simeq B)
    \end{equation*}
    for each $A,B:\UU$, and the type $A\simeq B$ is in $\UU$.
  \item For any family $B$ of locally $\UU$-small types over a $\UU$-small type $A$, the dependent product $\prd{x:A}B(x)$ is locally $\UU$-small.\index{dependent function type!is locally U-small@{is locally $\UU$-small}}\index{locally small type!dependent function type is locally small}
  \end{enumerate}
\end{eg}

We are now ready to assume the replacement axiom.

\begin{axiom}[The replacement axiom]\label{axiom:replacement}
  For any universe $\UU$, we assume that for any map $f:A\to B$ from a $\UU$-small type $A$ into a locally $\UU$-small type $B$, the image of $f$ is $\UU$-small.\index{replacement axiom|textbf}\index{axiom!replacement axiom|textbf}
\end{axiom}

\begin{eg}
  For any type $A:\UU$, the type $\UU_A$ of all types in $\UU$ merely equivalent to $A$ is equivalent to the image of the constant map $\const_A:\unit\to \UU$ is small. Since $\unit$ is small and $\UU$ is locally $\UU$-small, it follows from the replacement axiom that $\UU_A$ is $\UU$-small.\index{U A@{$\UU_A$}!is U-small@{is $\UU$-small}}
\end{eg}

\begin{eg}
  The type $\F$ of all finite types in $\UU$ is equivalent to be the image of the map
  \begin{equation*}
    \Fin{} : \N\to\UU.
  \end{equation*}
  Since $\N$ is $\UU$-small and $\UU$ is locally $\UU$-small, it follows from the replacement axiom that $\F$ is $\UU$-small.\index{F@{$\F$}!is U-small@{is $\UU$-small}}\index{small type!F is U-small@{$\F$ is $\UU$-small}}
\end{eg}

\begin{eg}
  Consider a type $A$ in $\UU$ and an equivalence relation $R$ on $A$ in $\UU$. Then the type $A/R$ is $\UU$-small, since it is equivalent to the image of
  \begin{equation*}
    R:A\to (A\to\prop_\UU),
  \end{equation*}
  which maps the $\UU$-small type $A$ into the locally $\UU$-small type $A\to\prop_\UU$.\index{set quotient!is U-small@{is $\UU$-small}}\index{small type!set quotient is U-small@{set quotient is $\UU$-small}}
\end{eg}

\subsection{The universal property of set quotients}
\index{universal property!of set quotients|(}
\index{set quotient!universal property|(}

The quotient $A/R$ is constructed as the image of $R$, so we obtain a commuting triangle
\begin{equation*}
  \begin{tikzcd}[column sep=-1em]
    A \arrow[rr,"q_R"] \arrow[dr,swap,"R"] & & A/R \arrow[dl,hook,"i_R"] \\
    \phantom{A/R} & \prop_\UU^A,
  \end{tikzcd}
\end{equation*}
and the embedding $i_R:A/R\to\prop_\UU^A$ satisfies the universal property of the image of $R$. This universal property is, however, not the usual universal property of the quotient.

\begin{defn}
  Consider a map $q:A\to B$ into a set $B$ satisfying the property that
  \begin{equation*}
    R(x,y)\to (q(x)=q(y))
  \end{equation*}
  for all $x,y:A$. We say that $q:A\to B$ \define{is a set quotient}\index{is set quotient}\index{set quotient|textbf} of $R$, or that $q$ satisfies the \define{universal property of the set quotient by $R$}\index{universal property!of set quotients|textbf}\index{set quotient!universal property|textbf}, if for every map $f:A\to X$ into a set $X$ such that $f(x)=f(y)$ whenever $R(x,y)$ holds, there is a unique extension
  \begin{equation*}
    \begin{tikzcd}
      A \arrow[d,swap,"q"] \arrow[dr,"f"] \\
      B \arrow[r,dashed] & X.
    \end{tikzcd}
  \end{equation*}
\end{defn}

\begin{rmk}
  Formally, we express the universal property of the quotient by $R$ as follows. Consider a map $q:A\to B$ that satisfies the property that
  \begin{equation*}
    H:\prd{x,y:A}R(x,y)\to (f(x)=f(y)).
  \end{equation*}
  Then there is for any set $X$ a map
  \begin{equation*}
    q^\ast:(B\to X) \to \Big(\sm{f:A\to X}\prd{x,y:A}R(x,y)\to (f(x)=f(y))\Big).
  \end{equation*}
  This map takes a function $h:B\to X$ to the pair
  \begin{equation*}
    q^\ast(h)\defeq(h\circ q,\lam{x}\lam{y}\lam{r}\ap{h}{H_{x,y}(r)}).
  \end{equation*}
  The universal property of the set quotient of $R$ asserts that the map $q^\ast$ is an equivalence for every set $X$. It is important to note that the universal property of set quotients is formulated with respect to sets.
\end{rmk}

\begin{thm}\label{thm:quotient_up}
  Consider a type $A$ and a universe $\UU$ containing $A$. Furthermore, let $R:A\to (A\to \prop_\UU)$ be an equivalence relation\index{equivalence relation}, and consider a map $q:A\to B$ into a set $B$, not necessarily in $\UU$. Then the following are equivalent.
  \begin{enumerate}
  \item \label{item:thm-quotient-up}The map $q$ satisfies the property that
    \begin{equation*}
      q(x)=q(y)
    \end{equation*}
    for every $x,y:A$ for which $R(x,y)$ holds, and moreover $q$ satisfies the universal property of the set quotient of $R$.
  \item \label{item:thm-quotient-effective}The map $q$ is surjective and \define{effective}\index{effective map|textbf}\index{equivalence relation!effective map|textbf}, which means that for each $x,y:A$ we have an equivalence
    \begin{equation*}
      (q(x)=q(y))\simeq R(x,y).
    \end{equation*}
  \item \label{item:thm-quotient-up-image}The map $R:A\to (A\to \prop_\UU)$ extends along $q$ to an embedding
    \begin{equation*}
      \begin{tikzcd}[column sep=tiny]
        A \arrow[rr,"q"] \arrow[dr,swap,"R"] & & B \arrow[dl,dashed,"i"] \\
        & \prop_\UU^A
      \end{tikzcd}
    \end{equation*}
    and the embedding $i$ satisfies the universal property of the image inclusion of $R$.
  \end{enumerate}
\end{thm}

In \cref{thm:quotient_up} we don't assume that $B$ is in the same universe as $A$ and $R$, because we want to apply it to $B\defeq\im(R)$. As we will see below, this extra generality only affects the proof that \ref{item:thm-quotient-effective} implies \ref{item:thm-quotient-up-image}.

\begin{proof}
  We first show that \ref{item:thm-quotient-effective} is equivalent to \ref{item:thm-quotient-up-image}, since this is the easiest part. After that, we will show that \ref{item:thm-quotient-up} is equivalent to \ref{item:thm-quotient-effective}.

  Assume that \ref{item:thm-quotient-up-image} holds. Then $q$ is surjective by \cref{thm:surjective}. Moreover, we have
  \begin{align*}
    R(x,y) & \simeq R(x)=R(y) \\
           & \simeq i(q(x))=i(q(y)) \\
           & \simeq q(x)=q(y)
  \end{align*}
  In this calculation, the first equivalence holds by \cref{cor:eq-quotient}; the second equivalence holds since we have a homotopy $R\htpy i\circ q$; and the third equivalence holds since $i$ is an embedding. This completes the proof that \ref{item:thm-quotient-up-image} implies \ref{item:thm-quotient-effective}.

  Next, we show that \ref{item:thm-quotient-effective} implies \ref{item:thm-quotient-up-image}. First, we want to define a map 
  \begin{equation*}
    i:B\to\prop_\UU^A.
  \end{equation*}
  We would like to define $i(b,a):=(b=q(a))$. This direct definition does not go through, however, because the type $B$ is not assumed to be in $\UU$. Nevertheless, observe that by the assumption that $q$ is surjective and effective, the type $B$ is locally $\UU$-small. To see this, first note that $\issmall_U(X)$ is a proposition for any type $X$ by \cref{prp:small}. Using the assumption that $q$ is surjective, it follows from \cref{prp:surjective} that it suffices to show that $q(a)=q(a')$ is $\UU$-small for each $a,a':A$. This follows by the assumption that $q$ is effective. In particular, the identity type $b=q(a)$ is a $\UU$-small proposition, for every $b:B$ and $a:A$. Let us write $s(b,a)$ for the element of type $\issmall_\UU(b=q(a))$.

  Now consider a universe $\VV$ containing $B$. Then we can define a map
  \begin{equation*}
    j:B\to\Big(A\to\sm{P:\prop_\VV}\issmall_\UU(P)\Big)
  \end{equation*}
  by $j(b,a):=(b=q(a),s(b,a))$, and now we obtain $i$ from $j$ by defining
  \begin{equation*}
    i(b,a):=\proj 1(s(b,a)).
  \end{equation*}
  Note that we have an equivalence $i(b,a)\simeq (b=q(a))$ for every $b:B$ and $a:A$.
  Then the triangle
  \begin{equation*}
    \begin{tikzcd}[column sep=tiny]
      A \arrow[rr,"q"] \arrow[dr,swap,"R"] & & B \arrow[dl,"i"] \\
      & \prop_\UU^A
    \end{tikzcd}
  \end{equation*}
  commutes, since we have an equivalence
  \begin{equation*}
    i(q(a),a') \simeq (q(a)=q(a')) \simeq R(a,a')
  \end{equation*}
  for each $a,a':A$. To show that $i$ is an embedding, recall from \cref{cor:is-emb-is-injective} that it suffices to show that $i$ is injective, i.e., that
  \begin{equation*}
    \prd{b,b':B}(i(b)=i(b'))\to (b=b'),
  \end{equation*}
  since the codomain of $i$ is a set by \cref{prp:propositional-extensionality}. Note that injectivity of a map into a set is a property, and that $q$ is assumed to be surjective. Hence by \cref{prp:surjective} it is sufficient to show that
  \begin{equation*}
    \prd{a,a':A}(i(q(a))=i(q(a')))\to (q(a)=q(a')).
  \end{equation*}
  Since $R\htpy i\circ q$, and $q(a)=q(a')$ is assumed to be equivalent to $R(a,a')$, it suffices to show that
  \begin{equation*}
    \prd{a,a':A}(R(a)=R(a'))\to R(a,a'),
  \end{equation*}
  which follows directly from \cref{cor:eq-quotient}. Thus we have shown that the factorization $R\htpy i\circ q$ factors $R$ as a surjective map followed by an embedding. We conclude by \cref{thm:surjective} that the embedding $i$ satisfies the universal property of the image factorization of $R$, which finishes the proof that \ref{item:thm-quotient-effective} implies \ref{item:thm-quotient-up-image}.
  
  Now we show that \ref{item:thm-quotient-up} implies \ref{item:thm-quotient-effective}. To see that $q$ is surjective if it satisfies the assumptions in \ref{item:thm-quotient-up}, consider the image factorization
  \begin{equation*}
    \begin{tikzcd}[column sep=tiny]
      A \arrow[dr,swap,"q"] \arrow[rr,"q_q"] & & \im(q) \arrow[dl,"i_q"] \\
      \phantom{\im(q)} & B.
    \end{tikzcd}
  \end{equation*}
  We claim that the map $i_q$ has a section. To see this, we first note that we have
  \begin{equation*}
    q_q(x)=q_q(y)
  \end{equation*}
  for any $x,y:A$ satisfying $R(x,y)$, because if $R(x,y)$ holds, then $q(x)=q(y)$ and hence $i_q(q_q(x))=i_q(q_q(y))$ holds and $i_q$ is an embedding. Since $\im(q)$ is a set, we may apply the universal property of $q$ and we obtain a unique extension of $q_q$ along $q$
  \begin{equation*}
    \begin{tikzcd}
      A \arrow[d,swap,"q"] \arrow[dr,"q_q"] \\
      B \arrow[r,dashed,swap,"h"] & \im(q).
    \end{tikzcd}
  \end{equation*}
  Now we observe that the composite $i_q\circ h$ is an extension of $q$ along $q$, so it must be the identity function by uniqueness. Thus we have established that $h$ is a section of $i_q$. Since $i_q$ is an embedding with a section, it follows that $i_q$ is an equivalence. We conclude that $q$ is surjective, because $q$ is the composite $i_q\circ q_q$ of a surjective map followed by an equivalence.

  Now we have to show that the map $q$ is effective, i.e., that $q(x)=q(y)$ is equivalent to $R(x,y)$ for every $x,y:A$. We first apply the universal property of $q$ to obtain for each $x:A$ an extension of $R(x)$ along $q$
  \begin{equation*}
    \begin{tikzcd}
      A \arrow[d,swap,"q"] \arrow[dr,"R(x)"] \\
      B \arrow[r,dashed,swap,"\tilde{R}(x)"] & \prop_\UU.
    \end{tikzcd}
  \end{equation*}
  Since the triangle commutes, we have an equivalence $\tilde{R}(x,q(x'))\simeq R(x,x')$ for each $x':A$. Now we apply \cref{thm:id_fundamental} to see that the canonical family of maps
  \begin{equation*}
    \prd{y:B}(q(x)=y)\to \tilde{R}(x,y)
  \end{equation*}
  is a family of equivalences. Thus, we need to show that the type $\sm{y:B}\tilde{R}(x,y)$ is contractible. For the center of contraction, note that we have $q(x):B$, and the type $\tilde{R}(x,q(x))$ is equivalent to the type $R(x,x)$, which is inhabited by reflexivity of $R$. To construct the contraction, it suffices to show that
  \begin{equation*}
    \prd{y:B}\tilde{R}(x,y)\to (q(x)=y).
  \end{equation*}
  Since this is a property, and since we have already shown that $q$ is a surjective map, we may apply \cref{prp:surjective}, by which it suffices to show that
  \begin{equation*}
    \prd{x':A}\tilde{R}(x,q(x'))\to (q(x)=q(x')).
  \end{equation*}
  Since $\tilde{R}(x,q(x'))\simeq R(x,x')$, this is immediate from our assumption on $q$. Thus we obtain the contraction, and we conclude that we have an equivalence $\tilde{R}(x,y)\simeq (q(x)=y)$ for each $y:B$. It follows that we have an equivalence
  \begin{equation*}
    R(x,y)\simeq (q(x)=q(y))
  \end{equation*}
  for each $x,y:A$, which completes the proof that \ref{item:thm-quotient-up} implies \ref{item:thm-quotient-effective}.
  
  It remains to show that \ref{item:thm-quotient-effective} implies \ref{item:thm-quotient-up}. Assume \ref{item:thm-quotient-effective}, and let $f:A\to X$ be a map into a set $X$, satisfying the property that
  \begin{equation*}
    \prd{a,a':A}R(a,a')\to (f(a)=f(a')).
  \end{equation*}
  Our goal is to show that the type of extensions of $f$ along $q$ is contractible. By \cref{ex:surjective-precomp} it follows that there is at most one such an extension, so it suffices to construct one.

  In order to construct an extension, we will construct for every $b:B$ a term $x:X$ satisfying the property
  \begin{equation*}
    P(x)\defeq \exists_{(a:A)}(f(a)=x)\land (q(a)=b).
  \end{equation*}
  Before we make this construction, we first observe that there is at most one such $x$, i.e., that the type of $x:X$ satisfying $P(x)$ is in fact a proposition. To see this, we need to show that $x=x'$ for any $x,x':X$ satisfying $P(x)$ and $P(x')$. Since $X$ is assumed to be a set, our goal of showing that $x=x'$ is a property. Therefore we may assume that we have $a,a':A$ satisfying
  \begin{align*}
    f(a) & = x & q(a) & = b \\
    f(a') & = x' & q(a') & = b.
  \end{align*}
  It follows from these assumptions that $q(a)=q(a')$, and hence that $R(a,a')$ holds. This in turn implies that $f(a)=f(a')$, and hence that $x=x'$.

  Now let $b:B$. Our goal is to construct an $x:X$ that satisfies the property
  \begin{equation*}
    \exists_{(a:A)}(f(a)=x)\land (q(a)=b).
  \end{equation*}
  Since $q$ is assumed to be surjective, we have $\brck{\fib{q}{b}}$. Moreover, since we have shown that at most one $x:X$ exists with the asserted property, we get to assume that we have $a:A$ satisfying $q(a)=b$. Now we see that $x\defeq f(a)$ satisfies the desired property.

  Thus, we obtain a function $h:B\to X$ satisfying the property that for all $b:B$ there exists an $a:A$ such that
  \begin{equation*}
    f(a)=h(b)\qquad\text{and}\qquad q(a)=b.
  \end{equation*}
  In particular, it follows that $h(q(a))=f(a)$ for all $a:A$, which completes the proof that \ref{item:thm-quotient-effective} implies \ref{item:thm-quotient-up}.  
\end{proof}

\begin{cor}
  Consider an equivalence relation $R$ over a type $A$. Then the quotient map
  \begin{equation*}
    q:A\to A/R
  \end{equation*}
  is surjective and effective, and it satisfies the universal property of the set quotient.
\end{cor}

\cref{thm:quotient_up} can be used to show that the type of equivalence relations is equivalent to the type of sets $X$ equipped with a surjective map $A\twoheadrightarrow X$. This may seem remarkable if you haven't tried \cref{ex:surjection-into-k-type} yet, because at first glance one might think that the type of sets $X$ equipped with a surjective map $A\twoheadrightarrow X$ is a $1$-type, while the type of equivalence relations on $A$ is a set.

\begin{thm}\label{thm:eqrel-surj}
  For any type $A$ and any universe $\UU$ containing $A$, we have an equivalence
  \index{surjective map!surjective maps into sets are set quotients}\index{set quotient!surjective maps into sets are set quotients}
  \begin{equation*}
    \eqrel_\UU(A)\simeq\sm{X:\Set_\UU}A\twoheadrightarrow X.
  \end{equation*}
\end{thm}

\begin{proof}
  Given an equivalence relation $R:A\to(A\to\prop_\UU)$ on $A$ we first use the replacement axiom, by which the set quotient $A/R$ is $\UU$-small, to obtain a set $Q(R):\Set_\UU$, an equivalence $e:Q(R)\simeq A/R$, and a surjective map $f:A\to Q(R)$ such that the triangle
  \begin{equation*}
    \begin{tikzcd}[column sep=1em]
      \phantom{Q(R)} & A \arrow[dl,swap,"f"] \arrow[dr,"q"] & \phantom{Q(R)}\\
      Q(R) \arrow[rr,swap,"e"] & & A/R
    \end{tikzcd}
  \end{equation*}
  commutes. This defines a map
  \begin{equation*}
    \mathcal{Q}_A:\eqrel_\UU(A)\to\sm{X:\Set_\UU}A\twoheadrightarrow X.
  \end{equation*}
  The map $\mathcal{K}_A:\big(\sm{X:\Set_\UU}A\twoheadrightarrow X\big)\to\eqrel_\UU(A)$ is given by
  \begin{equation*}
    \mathcal{K}_A(X,f,x,y)\defeq K_f(x,y) \defeq (f(x)=f(y)).
  \end{equation*}
  Note that $K_f$ is valued in propositions because $X$ is assumed to be a set, and obviously it is an equivalence relation.

  To see that $\mathcal{K}_A(\mathcal{Q}_A(R))=R$ note that by function extensionality and propositional extensionality it follows that two equivalence relations $R$ and $S$ on $A$ are equal if and only if $R(x,y)\leftrightarrow S(x,y)$ for all $x,y:A$. Note that $\mathcal{K}_A(\mathcal{Q}_A(R))(x,y)\leftrightarrow R(x,y)$ holds for all $x,y:A$ if and only if $(q_R(x)=q_R(y))\leftrightarrow R(x,y)$ holds for all $x,y:A$. This follows from \cref{cor:eq-quotient}.

  It remains to show that $\mathcal{Q}_A(\mathcal{K}_A(X,f))=(X,f)$. Note that the type of identifications $(Y,g)=(X,f)$ is by the univalence axiom equivalent to the type
  \begin{equation*}
    \sm{e:Y\simeq X}e\circ g\htpy f.
  \end{equation*}
  Therefore it suffices to construct a commuting triangle
  \begin{equation*}
    \begin{tikzcd}[column sep=1em]
      \phantom{A/K_f} & A \arrow[dl,swap,"q_{K_f}"] \arrow[dr,"f"] & \phantom{A/K_f} \\
      A/K_f \arrow[rr] & & X
    \end{tikzcd}
  \end{equation*}
  We obtain such an equivalence by combining \cref{thm:quotient_up} and \cref{thm:uniqueness-image}.
\end{proof}
\index{universal property!of set quotients|)}
\index{set quotient!universal property|)}

\subsection{Partitions}
\index{partition|(}
\index{set quotient!partition|(}

There are many equivalent ways of stating what an equivalence relation is. We saw in \cref{thm:eqrel-surj} that the type of equivalence relations on $A$ is equivalent to the type of surjective maps out of $A$ into a set. Here we will show that the type of equivalence relations on $A$ is equivalent to the type of partitions of $A$. Another type that is equivalent to the type of equivalence relations of $A$ is the type of set-indexed $\Sigma$-decompositions of $A$, i.e., the type of triples $(X,Y,e)$ consisting of a set $X$, a family $Y$ of inhabited types indexed by $X$, and an equivalence $e:A\simeq \sm{x:X}Y(x)$. The fact that the type of equivalence relations on $A$ is equivalent to the type of set-indexed $\Sigma$-decompositions of $A$ is stated as \cref{ex:sigmadecompositions}

In this section we show that equivalence relations on $A$ are partitions of $A$. Recall that the type of inhabited subtypes of $A$ is defined to be
\begin{equation*}
  \mathcal{P}_{\mathcal{U}}^+(A)\defeq\sm{Q:A\to\prop_\UU}\Brck{\sm{a:A}Q(a)}.
\end{equation*}
The equivalence of equivalence relations and partitions requires some finesse regarding universes. This is why we set up the definition of partitions in the following way.

\begin{defn}
  Let $A$ be a type and let $\UU$ and $\VV$ be two universes. A \define{$(\UU,\VV)$-partition}\index{partition|textbf}\index{set quotient!partition|textbf} of a type $A$ is a subset
  \begin{equation*}
    P:\mathcal{P}_{\UU}^+(A)\to\prop_{\VV}
  \end{equation*}
  of the type of inhabited subsets of $A$ such that for each $x:A$ there is a unique inhabited subset $Q$ of $A$ in $P$ that contains $x$, i.e., if it comes equipped with an element of type
  \begin{equation*}
    \ispartition(P):=\prd{x:A}\iscontr\Big(\sm{Q:\mathcal{P}_{\UU}^+(A)}P(Q)\times Q(x)\Big)
  \end{equation*}
  The type of all $(\UU,\VV)$-partitions of $A$ is defined by
  \begin{equation*}
    \partition_{\UU,\VV}(A)\defeq\sm{P:\mathcal{P}_{\UU}^+(A)\to\prop_{\VV}}\ispartition(P)
  \end{equation*}
\end{defn}

\begin{thm}
  Consider a type $A$, a universe $\UU$, and consider a universe $\VV$ containing both $A$ and every type in $\UU$. Then we have an equivalence
  \begin{equation*}
    \eqrel_{\UU}(A)\simeq\partition_{\UU,\VV}(A).
  \end{equation*}
\end{thm}

\begin{proof}
  Consider an equivalence relation $R$ on $A$. Then we define
  \begin{equation*}
    P:\mathcal{P}_{\UU}^+(A)\to\prop_\VV
  \end{equation*}
  by $P(Q)\defeq\exists_{(x:A)}\forall_{(y:A)}Q(y)\leftrightarrow R(x,y)$. In other words, $P$ is the subtype of equivalence classes of $R$, which are all inhabited. To show that $P$ is a partition of $A$, let $x:A$. The type
  \begin{equation*}
    \sm{Q:\mathcal{P}_{\UU}^+(A)}P(Q)\times Q(x)
  \end{equation*}
  is equivalent to the type
  \begin{equation*}
    \sm{Q:\mathcal{P}_{\UU}^+(A)}\prd{y:A}Q(y)\leftrightarrow R(x,y)
  \end{equation*}
  since the proposition $\exists_{(z:A)}\forall_{(y:A)}Q(y)\leftrightarrow R(z,y)$ is equivalent to the type $\prd{y:A}Q(y)\leftrightarrow R(x,y)$, given an element $q:Q(x)$. By univalence it follows that the latter type is equivalent to the identity type $Q=R(x)$ in $\mathcal{P}_{\UU}^+(A)$, so the total space is contractible. Thus we obtain a map
  \begin{equation*}
    \psi:\eqrel_{\mathcal{U}}(A)\to\partition_{\UU,\VV}(A).
  \end{equation*}
  
  For the converse map, we first define for any $(\UU,\VV)$-partition $P$ of $A$ a binary relation $R_P$ such that $R_P(x)$ is at the center of contraction in the type
  \begin{equation*}
    \sm{Q:\mathcal{P}^+_{\mathcal{U}}(A)}P(Q)\times Q(x).
  \end{equation*}
  In other words, $R_P(x)$ is defined to be the unique block in the partition $P$ such that $R_P(x,x)$ holds. It is immediate from its definition that $R_P(x,y)$ is a proposition in $\UU$. To see that $R_P$ is symmetric, note that if $R_P(x,y)$ holds, then $R_P(x)$ is an element of type
  \begin{equation*}
    \sm{Q:\mathcal{P}^+_{\mathcal{U}}(A)}P(Q)\times Q(y).
  \end{equation*}
  By contractibility, this implies that $R(x)=R(y)$, from which we obtain that $R(y,x)$ holds. To see that $R_P$ is transitive we observe similarly that if $R(x,y)$ and $R(y,z)$ hold, then we have an identification $R(x)=R(y)$ and it follows that $R(x,z)$ holds.  Thus we obtain a map
  \begin{equation*}
    \varphi:\partition_{\UU,\VV}(A)\to\eqrel_{\UU}(A).
  \end{equation*}
  It remains to prove that the maps $\psi$ and $\varphi$ are inverse to each other, first let $R$ be an equivalence relation. In order to show that $\varphi(\psi(R))=R$ it suffices by univalence to show that the equivalence relation obtained from the partition induced by $R$ is given by
  \begin{equation*}
    R'(x,y):=\sm{Q:\mathcal{P}_{\UU}^+(A)}\Big(\exists_{(u:A)}\forall_{(v:A)}Q(v)\leftrightarrow R(u,v)\Big)\times Q(x)\times Q(y).
  \end{equation*}
  is equivalent to $R$. Observe that the proposition $R'(x,y)$ is equivalent to $R(x,x)\times R(x,y)$, which is equivalent to $R(x,y)$. This shows that the composite
  \begin{equation*}
    \begin{tikzcd}
      \eqrel_{\UU}(A) \arrow[r,"\psi"] & \partition_{\UU,\VV}(A) \arrow[r,"\varphi"] & \eqrel_{\UU}(A)
    \end{tikzcd}
  \end{equation*}
  is homotopic to the identity function.

  Finally, we have to show that for any partition $P$ of $A$ and any inhabited subtype $Q$ of $A$ we have $\psi(\varphi(P))(Q)\leftrightarrow P(Q)$. Note that this is a proposition, so we may assume an element $x:A$ such that $Q(x)$ holds. By univalence it follows that $\psi(\varphi(P))(Q)$ holds if and only if $Q=R_P(x)$, where $R_P$ is the equivalence relation constructed in the definition of the map $\varphi$. Now we see that $P(Q)$ holds if and only if $Q$ is in the contractible type
  \begin{equation*}
    \sm{Q':\mathcal{P}^+_{\mathcal{U}}(A)}P(Q')\times Q'(x),
  \end{equation*}
  which is the case if and only if $Q=R_P(x)$. This shows that the composite
  \begin{equation*}
    \begin{tikzcd}
      \partition_{\UU,\VV}(A) \arrow[r,"\varphi"] & \eqrel_{\UU}(A) \arrow[r,"\psi"] & \partition_{\UU,\VV}(A)
    \end{tikzcd}
  \end{equation*}
  is homotopic to the identity function.
\end{proof}
\index{partition|)}
\index{set quotient!partition|)}

\subsection{Unique representatives of equivalence classes}
\index{equivalence class!choice of unique representatives|(}
\index{equivalence relation!choice of unique representatives|(}

A common way to construct set quotients is by showing that the equivalence classes of an equivalence relation have a choice of unique representatives. In this section we show that if there is a choice of unique representatives, then the set quotient can be constructed as the type of those representatives. An important reason to define set quotients as the type of canonical representatives, if that is possible, is that the universe level of the set quotient can be kept as low as possible without needing to appeal to the replacement axiom.

\begin{defn}
  Consider an equivalence relation $R$ on a type $A$, and consider a family of types $C(x)$ indexed by $x:A$. We say that $C$ is a \define{choice of (unique) representatives}\index{choice of unique representatives|textbf}\index{equivalence class!choice of unique representatives|textbf}\index{equivalence relation!choice of unique representatives|textbf} of the equivalence classes of $R$ if $C$ comes equipped with an element of type
  \begin{equation*}
    \ischoiceofrepresentatives(C) \defeq \prd{x:A}\iscontr\Big(\sm{y:A}C(y)\times R(x,y)\Big).
  \end{equation*}
\end{defn}

\begin{thm}\label{thm:choice-of-representatives}
  Consider an equivalence relation $R$ on a type $A$, and let $C$ be a choice of representatives of the equivalence classes of $R$, with $(h(x),c(x),r(x))$ at the center of contraction of $\sm{y:A}C(y)\times R(x,y)$. Then the map
  \begin{equation*}
    q:A\to\sm{x:A}C(x)
  \end{equation*}
  given by $q(x)\defeq(h(x),c(x))$ is a map into a set such that $q(x)=q(y)$ for every $x,y:A$ such that $R(x,y)$ holds, and moreover $q$ satisfies the universal property of the set quotient of $A$ by $R$. 
\end{thm}

\begin{proof}
  First, we will use \cref{lem:prop_to_id} to show that the type $\sm{y:A}C(y)$ is a set, such that
  \begin{equation*}
    ((x,c)=(y,d))\simeq R(x,y)
  \end{equation*}
  for any $(x,c)$ and $(y,d)$ in $\sm{y:A}C(y)$. Note that we have a function
  \begin{equation*}
    R(x,y)\to ((x,c)=(y,d)),
  \end{equation*}
  since for any $r:R(x,y)$ both $(x,c,r)$ and $(y,d,r)$ are elements of the contractible type ${\sm{y:A}C(y)\times R(x,y)}$. Since $R$ is a reflexive relation valued in propositions, the claim follows. In particular, it follows that
  \begin{equation*}
    (q(x)=q(y))\simeq R(x,y)
  \end{equation*}
  for any $x,y:A$, i.e., $q$ is effective.
  
  To prove the universal property of set quotients, note that by characterization \ref{item:thm-quotient-effective} in \cref{thm:quotient_up} it suffices to show that $q$ is surjective and effective. We have already shown above that $q$ is effective, so it remains to show that $q$ is surjective. In fact, we will prove the stronger claim that the projection map
  \begin{equation*}
    \proj 1:\sm{x:A}C(x)\to A  
  \end{equation*}
  is a section of $q$. Let $x:A$ and $c:C(x)$. Then $(x,c,\rho(x))$ is an element of the type
  \begin{equation*}
    \sm{y:A}C(y)\times R(x,y),
  \end{equation*}
  which is contractible with center of contraction $(h(x),c(x),r(x))$. Therefore it follows that $q(x)\jdeq (h(x),c(x))=(x,c)$. In particular, we see that $q(\proj 1(x,c))=(x,c)$, i.e., that $\proj 1$ is a section of $q$.
\end{proof}

\begin{eg}
  In \cref{prp:congruence-eqrel} we constructed the congruence relations $x\equiv y \mod k$ on the natural numbers for every natural number $k$, and in \cref{thm:effective-mod-k,thm:issec-nat-Fin} we showed that the map
  \begin{equation*}
    x\mapsto [x]_{k+1}:\N\to\Fin{k+1}
  \end{equation*}
  is effective and split surjective. By \cref{thm:quotient_up} it follows that the map
  \begin{equation*}
    x\mapsto [x]_{k+1}:\N\to\Fin{k+1}
  \end{equation*}
  satisfies the universal property of the set quotient of the equivalence relation $x,y\mapsto x\equiv y\mod k+1$.

  We also claim that there is a choice of representatives of the congruence relations.\index{congruence relations on N@{congruence relations on $\N$}!choice of unique representatives}\index{choice of unique representatives!for the congruence relations on N@{for the congruence relations on $\N$}} We define our choice of representatives by
  \begin{equation*}
    C(y)\defeq \fib{\natFin}{y},
  \end{equation*}
  where $\natFin:\Fin{k+1}\to\N$ is the inclusion of $\Fin{k+1}$ into $\N$ constructed in \cref{defn:natFin}. To see that $C$ is a choice of representatives, we have to prove that
  \begin{equation*}
    \sm{y:\N}C(y)\times (x\equiv y\mod k+1)
  \end{equation*}
  is contractible for each $x:\N$. At the center of contraction we have the triple $(\natFin([x]_{k+1}),([x]_{k+1},\refl{}),p)$ where $p:x\equiv\natFin([x]_{k+1})\mod k+1$ is the proof obtained via \cref{thm:effective-mod-k,thm:issec-nat-Fin}. In order to construct the contraction, note that both $C(y)$ and $x\equiv y\mod k+1$ are propositions for each $y:\N$. Therefore it suffices to prove that for any $y:\N$ such that $C(y)$ and $x\equiv y\mod k+1$ hold, we have
  \begin{equation*}
    \natFin([x]_{k+1})=y.
  \end{equation*}
  Since $C(y)$ holds, we see that $y=\natFin([y]_{k+1})$. Therefore it suffices to prove that $[x]_{k+1}=[y]_{k+1}$. This follows from \cref{thm:effective-mod-k}, since we assumed $x\equiv y\mod k+1$.
\end{eg}

\begin{eg}
  Consider the type of \define{(integer) fractions}\index{fraction|textbf}\index{integers!integer fractions|textbf}
  \begin{equation*}
    Q\defeq \Z\times\sm{y:\Z}y\neq 0.
  \end{equation*}
  We define an equivalence relation on $Q$ by
  \begin{equation*}
    ((x,y)\sim (x',y'))\defeq (xy'=x'y).
  \end{equation*}
  This equivalence relation has a choice of representatives defined by
  \index{choice of unique representatives!for integer fractions}
  \index{fraction!choice of unique representatives}
  \begin{equation*}
    C(x,y)\defeq (y>0)\land (\gcd(x,y)=1). 
  \end{equation*}
  In other words, we say that $(x,y)$ is a \define{reduced fraction} if $y>0$ and $x$ and $y$ are coprime. 

  To see that $C$ defines a choice of unique representatives, we first need to construct the center of contraction of
  \begin{equation*}
    \sm{q:Q}C(q)\times ((x,y)\sim q).
  \end{equation*}
  Note that if $y<0$ then $(x,y)\sim (-x,-y)$, and we have $-y>0$. It is therefore safe to assume that $y>0$. We claim that
  \begin{equation*}
    (x/\gcd(x,y),y/\gcd(x,y)):Q
  \end{equation*}
  satisfies $C$ and is equivalent to $(x,y)$. It is immediate that $y/\gcd(x,y)>0$ and that $(x,y)\sim(x/\gcd(x,y),y/\gcd(x,y))$. The fact that $x/\gcd(x,y)$ and $y/\gcd(x,y)$ are coprime follows from the fact that
  \begin{equation*}
    \gcd(x/d,y/d)=\gcd(x,y)/d
  \end{equation*}
  for any common divisor $d$ of $x$ and $y$. 
  
  To construct a contraction, let $(x',y'):Q$ such that $C(x',y')$ and $(x,y)\sim (x',y')$. Since $C(q)$ and $(x,y)\sim q$ are propositions for every $q:Q$ it suffices to show that
  \begin{equation*}
    x'=x/\gcd(x,y)\qquad\text{and}\qquad y'=y/\gcd(x,y).
  \end{equation*}
  Since $x'$ and $y'$ are assumed to be coprime, it follows from the equation
  \begin{equation*}
    x'y/\gcd(x,y)=xy'/\gcd(x,y)
  \end{equation*}
  that $x'$ divides $x/\gcd(x,y)$. Similarly $x/\gcd(x,y)$ and $y/\gcd(x,y)$ are coprime, it follows from the same equation that $x/\gcd(x,y)$ divides $x'$, so we conclude that $ux'=x/\gcd(x,y)$ for some $u=\pm 1$. The fact that $vy'=y/\gcd(x,y)$ for some $v=\pm 1$ is proven similarly. However, since both $y$ and $y'$ are positive, and the $\gcd(x,y)$ of any two integers is positive, it follows that $v=1$. Using the assumption that $x'y/\gcd(x,y)=xy'/\gcd(x,y)$, this allows us to deduce that also $u=1$.

  We define the type of \define{rational numbers}\index{Q@{$\Q$}|see {rational numbers}}\index{rational numbers|textbf} by
  \begin{equation*}
    \Q\defeq \sm{(x,y):Q}(y>0)\land \gcd(x,y)=1,
  \end{equation*}
  and we define the quotient map $(x,y)\mapsto x/y:Q\to \Q$ to be the quotient map $q$ in \cref{thm:choice-of-representatives}. By \cref{thm:choice-of-representatives} it also follows that $(x,y)\mapsto x/y$ satisfies the universal property of the set quotient of the equivalence relation $\sim$ on $Q$.
\end{eg}
\index{equivalence class!choice of unique representatives|)}
\index{equivalence relation!choice of unique representatives|)}

\subsection{Set truncations}
\index{set truncation|(}
\index{set quotient!set truncation|(}

An important instance of set quotients in the univalent foundations of mathematics is the notion of set truncation. Analogous to the propositional truncation, the set truncation of a type $A$ is a map $\eta:A\to \trunc{0}{A}$ into a set $\trunc{0}{A}$ such that any map $f:A\to X$ into a set $X$ extends uniquely along $\eta$:
\begin{equation*}
  \begin{tikzcd}
    A \arrow[dr,"f"] \arrow[d,swap,"\eta"] \\
    \trunc{0}{A} \arrow[r,dashed] & X.
  \end{tikzcd}
\end{equation*}
In other words, the set truncation $\eta:A\to\trunc{0}{A}$ is the universal way of mapping $A$ into a set. We first specify what it means for a map $f:A\to B$ into a set $B$ to be a set truncation of $A$.

\begin{defn}
  We say that a map $f:A\to B$ into a set $B$ is a \define{set truncation}\index{set truncation|textbf}\index{set quotient!set truncation|textbf}\index{universal property!of set truncations|textbf}\index{set truncation!universal property|textbf} if the precomposition function
  \begin{equation*}
    \blank\circ f : (B\to X)\to (A\to X)
  \end{equation*}
  is an equivalence for every set $X$. 
\end{defn}

In the following theorem we prove several conditions that are equivalent to being a set truncation.

\begin{thm}\label{thm:set-truncation}
  Consider a map $f:A\to B$ into a set $B$. Then the following are equivalent:
  \begin{enumerate}
  \item\label{item:is-set-truncation} The map $f$ is a set truncation.
  \item\label{item:dup-set-truncation} The map $f$ satisfies the \textbf{dependent universal property}\index{set truncation!dependent universal property|textbf}\index{dependent universal property!of set truncations|textbf} of the set truncation: For every family $X$ of sets over $B$, the precomposition function
    \begin{equation*}
      \blank\circ f : \Big(\prd{b:B}X(b)\Big)\to\Big(\prd{a:A}X(f(a))\Big)
    \end{equation*}
    is an equivalence.
  \item\label{item:is-quotient-set-truncation} The map $f$ is surjective and effective with respect to the equivalence relation $x,y\mapsto\brck{x=y}$, i.e., we have equivalences
    \begin{equation*}
      (f(x)=f(y))\simeq \brck{x=y}
    \end{equation*}
    for every $x,y:A$.
  \end{enumerate}
\end{thm}

\begin{proof}
  The fact that \ref{item:dup-set-truncation} implies \ref{item:is-set-truncation} is immediate. Moreover, the fact that \ref{item:is-set-truncation} is equivalent to \ref{item:is-quotient-set-truncation} follows from the fact that any map $h:A\to X$ into a set $X$ comes equipped with a function
  \begin{equation*}
    \brck{x=y}\to (h(x)=h(y))
  \end{equation*}
  for every $x,y:A$. 
  
  It remains to prove that \ref{item:is-set-truncation} implies \ref{item:dup-set-truncation}. Consider a family $X$ of sets over $B$, and consider the commuting square
  \begin{equation*}
    \begin{tikzcd}[column sep=7em]
      \sm{g:B\to B}\prd{b:B}X(g(b)) \arrow[d,swap,"\simeq"] \arrow[r,"{(g,s)\mapsto (g\circ f,s\circ f)}"] & \sm{h:A\to B}\prd{a:A}X(h(a)) \arrow[d,"\simeq"] \\
      (B\to\sm{b:B}X(b)) \arrow[r,swap,"\blank\circ f"] & (A\to\sm{b:B}X(b))
    \end{tikzcd}
  \end{equation*}
  The side maps are equivalences by the distributivity of $\Pi$ over $\Sigma$, and the bottom map is an equivalence by the assumption that $f$ is a set truncation. Therefore it follows that the top map is an equivalence. Furthermore, note that the map
  \begin{equation*}
    \blank\circ f : (B\to B)\to (A\to B)
  \end{equation*}
  is an equivalence by the assumption that $f$ is a set truncation. Therefore it follows from \cref{thm:equiv-toto} that the map
  \begin{equation*}
    \blank\circ f : \Big(\prd{b:B}X(g(b))\Big)\to \Big(\prd{a:A}X(g(f(a)))\Big)
  \end{equation*}
  is an equivalence for every $g:B\to B$. Now we take $g\defeq \idfunc$ to complete the proof that \ref{item:is-set-truncation} implies \ref{item:dup-set-truncation}.
\end{proof}

\begin{cor}
  On any universe $\UU$, there is an operation $\trunc{0}{\blank}:\UU\to\Set_\UU$\index{[[A]] 0@{$\trunc{0}{A}$}|see {set truncation}} such that every type $A$ in $\UU$ comes equipped with a map
  \begin{equation*}
    \eta:A\to\trunc{0}{A}
  \end{equation*}
  that satisfies the universal property of the set truncation. The set $\trunc{0}{A}$ is called the \define{set truncation}\index{set truncation|textbf} of $A$.
\end{cor}

\begin{proof}
  By \cref{thm:set-truncation} it follows that a map $f:A\to B$ into a set $B$ is a set truncation if and only if it is a quotient map with respect to the equivalence relation $x,y\mapsto\brck{x=y}$. Given a type $A$ in $\UU$, the quotient of $A$ by $x,y\mapsto\brck{x=y}$ is equivalent to a type in $\UU$ by the replacement axiom.
\end{proof}

\begin{cor}
  The set truncation $\eta:A\to\trunc{0}{A}$ is surjective and effective with respect to the equivalence relation $x,y\mapsto\brck{x=y}$, i.e., we have an equivalence
  \begin{equation*}
    (\eta(x)=\eta(y))\simeq \brck{x=y}
  \end{equation*}
  for each $x,y:A$. 
\end{cor}

By this corollary, we may think of the set truncation $\trunc{0}{A}$ of $A$ as the set of connected components of $A$. Indeed, if we have an unspecified identification $\brck{x=y}$ in $A$, then we think of $x$ and $y$ as being in the same connected component. For example, any $k$-element set is a type that is in the same connected component of $\UU$ as the type $\Fin{k}$.

\begin{defn}
  A type $A$ is said to be \define{connected}\index{connected type|textbf} if its set truncation $\trunc{0}{A}$ is contractible. We define\index{is-conn@{$\isconn(A)$}|textbf}
  \begin{equation*}
    \isconn(A)\defeq\iscontr\trunc{0}{A}.
  \end{equation*}
  Furthermore, we say that a map $f:A\to B$ is \define{connected}\index{connected map|textbf} if all its fibers are connected.
\end{defn}

\begin{rmk}
  In particular, every connected type is inhabited, because if $\trunc{0}{A}$ is contractible, then we have equivalences\index{inhabited type!connected types are inhabited}\index{connected type!connected types are inhabited}
  \begin{equation*}
    \brck{A}\simeq (\trunc{0}{A}\to\brck{A}) \simeq (A\to \brck{A}),
  \end{equation*}
  and the latter type contains the unit of the propositional truncation.
\end{rmk}

Using the notion of connectivity, we can add one more property to the list of equivalent characterizations of set truncations given in \cref{thm:set-truncation}.

\begin{thm}\label{thm:unit-set-truncation-connected}
  Consider a map $f:A\to B$ into a set $B$. Then the following are equivalent:
  \begin{enumerate}
  \item \label{item:unit-set-truncation-connected-i}The map $f$ is a set truncation.
  \item \label{item:unit-set-truncation-connected-ii}The map $f$ is connected.
  \end{enumerate}
\end{thm}

\begin{proof}
  First, suppose that $f$ is a set truncation, and consider $b:B$. Our goal is to show that the type
  \begin{equation*}
    \trunc{0}{\fib{f}{b}}
  \end{equation*}
  is contractible. Since $f$ is surjective by \cref{thm:set-truncation}, there exists an element $a:A$ equipped with an identification $f(a)=b$. We are proving a proposition, so it suffices to show that $\trunc{0}{\fib{f}{f(a)}}$ is contractible. At the center of contraction we have
  \begin{equation*}
    \eta(a,\refl{}):\trunc{0}{\fib{f}{f(a)}}.
  \end{equation*}
  In order to construct the contraction, we use the dependent universal property of the set truncation, by which it suffices to construct a function
  \begin{equation*}
    \prd{x:A}\prd{p:f(x)=f(a)}\eta(a,\refl{})=\eta(x,p)
  \end{equation*}
  Recall from \cref{thm:set-truncation} that the map $f$ is effective, so we have an equivalence $e:\brck{x=a}\simeq (f(x)=f(a))$ for every $x:A$. Furthermore, equality in set truncations are propositions, so we may even eliminate the propositional truncation from $\brck{x=a}$. Therefore it suffices to prove
  \begin{equation*}
    \prd{x:A}\prd{p:x=a}\eta(a,\refl{})=\eta(x,e(\eta(p)))
  \end{equation*}
  This is immediate, since $e(\eta(\refl{}))=\refl{}$. This completes the proof of \ref{item:unit-set-truncation-connected-i} implies \ref{item:unit-set-truncation-connected-ii}.

  For the converse, suppose that $f$ is connected, and consider a set $X$. Note that we have a commuting square
  \begin{equation*}
    \begin{tikzcd}[column sep=8em]
      \Big(\prd{b:B}\trunc{0}{\fib{f}{b}}\to X\Big) \arrow[r,"{h\mapsto\lam{b}{t}h(b,\eta(t))}"] & \Big(\prd{b:B}\fib{f}{b}\to X\Big) \arrow[d,swap,"{h\mapsto\lam{a}h(f(a),(a,\refl{}))}"] \\
      (B\to X) \arrow[r,swap,"\blank\circ f"] \arrow[u,"h\mapsto\lam{b}\lam{u}h(b)"] & (A\to X)
    \end{tikzcd}
  \end{equation*}
  In this commuting square, the map on the left is an equivalence since $\trunc{0}{\fib{f}{b}}$ is contractible for each $b:B$. The top map is an equivalence because $X$ is a set, and the right map is an equivalence by \cref{ex:pi-fib}. Therefore it follows that the bottom map is an equivalence, which completes the proof that \ref{item:unit-set-truncation-connected-ii} implies \ref{item:unit-set-truncation-connected-i}.
\end{proof}

\begin{rmk}
  There are truncation operations for every truncation level. That is, we can define for every type $A$ a map $\eta:A\to\trunc{k}{A}$ such that the map
\begin{equation*}
  \blank\circ\eta : (\trunc{k}{A}\to X)\to (A\to X)
\end{equation*}
is an equivalence for every $k$-truncated type $X$. To learn more about general $k$-truncations, we refer to Chapter 7 of \cite{hottbook}.
\end{rmk}
\index{set truncation|)}
\index{set quotient!set truncation|)}

\begin{exercises}
  \exitem Consider a proposition $P$, and define the relation $\sim_P$\index{~ P@{$\sim_P$}|textbf} on $\bool$ by
  \begin{align*}
    (\btrue\sim_P\btrue) & \defeq \unit & (\btrue\sim_P\bfalse) & \defeq P \\
    (\bfalse\sim_P\btrue) & \defeq P & (\bfalse\sim_P\bfalse) & \defeq \unit
  \end{align*}
  \begin{subexenum}
  \item  Show that $\sim_P$ is an equivalence relation.
  \item Consider a universe $\UU$ containing the proposition $P$. Construct an embedding ${\bool/{\sim}_P}\hookrightarrow\prop_\UU$.  
  \item Use the quotient $\bool/\sim_P$ to show that the axiom of choice implies the law of excluded middle.\index{axiom of choice!implies law of excluded middle}
  \end{subexenum}
  \exitem Consider an equivalence relation $R:A\to(A\to\prop_\UU)$ on $A$, where $\UU$ is a universe containing $A$. Show that type type
  \begin{equation*}
    \sm{X:\UU}\sm{f:A\twoheadrightarrow X}\prd{x,y:A}(f(x)=f(y))\simeq R(x,y)
  \end{equation*}
  is contractible.
  \exitem For any type $A$, show that the type of equivalence relations equipped with a choice of representatives of its equivalence classes is equivalent to the type of \define{set-based retracts} of $A$, i.e., the type\index{retract!set-based retract|textbf}\index{set-based retract|textbf}\index{Retr Set U X@{$\Retr_{\Set_\UU}(X)$}|textbf}\index{Retr Set U X@{$\Retr_{\Set_\UU}(X)$}!is equivalent to equivalence relations with canonical representatives}
  \begin{equation*}
    \Retr_{\Set_\UU}(A) \defeq \sm{X:\Set_\UU}\sm{i:X\to A}\sm{q:A\to X} q\circ i\htpy \idfunc.
  \end{equation*}
  \exitem  \label{ex:sigmadecompositions}A \define{$\Sigma$-decomposition}\index{S-decomposition@{$\Sigma$-decomposition}|textbf} of a type $A$ consists of a type $X$ (the \define{indexing type}\index{S-decomposition@{$\Sigma$-decomposition}!indexing type|textbf} of the $\Sigma$-decomposition) equipped with a family $Y$ of inhabited types indexed by $X$ and an equivalence
  \begin{equation*}
    e:A\simeq \sm{x:X}Y(x).
  \end{equation*}
  In other words, the type of all $\Sigma$-decompositions of $A$ is defined by
  \begin{equation*}
    \Sigmadecomposition_\UU(A) \defeq \sm{X:\UU}\sm{Y:X\to\sm{Z:\UU}\brck{Z}}A\simeq\sm{x:X}Y(x).
  \end{equation*}
  \begin{subexenum}
  \item Construct an equivalence
    \begin{equation*}
      \Sigmadecomposition_\UU(A)\simeq \sm{X:\UU}A\twoheadrightarrow X.
    \end{equation*}
  \item A $\Sigma$-decomposition is said to be \define{set-indexed}\index{S-decomposition@{$\Sigma$-decomposition}!set indexed|textbf}\index{set-indexed S-decomposition@{set-indexed $\Sigma$-decomposition}|textbf} if its indexing type is a set. We will write $\Sigmadecomposition_{\Set_\UU}(A)$ for the type of all set-indexed $\Sigma$-decompositions of $A$ in $\UU$. Construct an equivalence
    \begin{equation*}
      \eqrel_\UU(A)\simeq \Sigmadecomposition_{\Set_\UU}(A).
    \end{equation*}
  \end{subexenum}
  \exitem \label{ex:is-surjective-fiber-inclusion}Consider a type $A$ equipped with an element $a:A$. Show that the following are equivalent:
  \begin{enumerate}
  \item The type $A$ is connected.
  \item There is an element of type $\brck{a=x}$ for any $x:A$.
  \item For any family $B$ over $A$, the fiber inclusion\index{fiber inclusion}
    \begin{equation*}
      i_a:B(a)\to\sm{x:A}B(x)
    \end{equation*}
    defined in \cref{ex:is-trunc-fiber-inclusion} is surjective.
  \end{enumerate}
  \exitem \label{ex:poset-reflection}Consider a preorder $(A,\leq)$, and define for any $a:A$ the order preserving map
  \begin{equation*}
    y_a : \preord(\op{(A,\leq)},{(\prop_\UU,{\to})})
  \end{equation*}
  by $y_a(x)\defeq(x\leq a)$. Furthermore, define the \define{poset reflection}\index{poset reflection|textbf}\index{poset!poset reflection}\index{preorder!poset reflection} $\posetreflection{A}$\index{[[A]] Pos@{$\posetreflection{A}$}|see {poset reflection}} to be the image of the map
  \begin{equation*}
    a\mapsto y_a : A\to \preord(\op{(A,\leq)},{(\prop_\UU,{\to})}).
  \end{equation*}
  \begin{subexenum}
  \item Show that the image of the map $a\mapsto y_a$ satisfies the universal property of the set quotient of the equivalence relation\index{set quotient!poset reflection}\index{poset reflection!is set quotient}
    \begin{equation*}
      x,y\mapsto (x\leq y)\land (y\leq x).
    \end{equation*}
  \item Equip the type $\posetreflection{A}$ with the structure of a poset and construct an order preserving map $\eta : A \to \posetreflection{A}$ that satisfies the following universal property: For any poset $P$, any order preserving map $f:A\to P$ extends uniquely along $\eta$ to an order preserving map $g:\posetreflection{A}\to P$, as indicated in the following diagram:\index{universal property!of poset reflections|textbf}\index{poset reflection!universal property|textbf}
  \begin{equation*}
    \begin{tikzcd}
      A \arrow[r,"f"] \arrow[d,swap,"\eta"] & P. \\
      \posetreflection{A} \arrow[ur,dashed]
    \end{tikzcd}
  \end{equation*}
  \end{subexenum}
  \exitem Consider a map $f:A\to B$.
  \begin{subexenum}
  \item Show that the type of maps $\trunc{0}{f}:\trunc{0}{A}\to\trunc{0}{B}$ equipped with a homotopy witnessing that the square
    \begin{equation*}
      \begin{tikzcd}
        A \arrow[r,"f"] \arrow[d,swap,"\eta"] & B \arrow[d,"\eta"] \\
        \trunc{0}{A} \arrow[r,swap,"\trunc{0}{f}"] & \trunc{0}{B}
      \end{tikzcd}
    \end{equation*}
    commutes is contractible.\index{functorial action!set truncation|textbf}\index{set truncation!functorial action|textbf}
  \item Show that if $f$ is injective, then $\trunc{0}{f}:\trunc{0}{A}\to\trunc{0}{B}$ is injective.
  \item Show that the following are equivalent
    \begin{enumerate}
    \item The map $f$ is surjective.
    \item the map $\trunc{0}{f}:\trunc{0}{A}\to\trunc{0}{B}$ is surjective.
    \end{enumerate}
  \item Construct a map $h:\im(f)\to\im\trunc{0}{f}$ such that the squares
    \begin{equation*}
      \begin{tikzcd}
        A \arrow[r,"q_f"] \arrow[d,swap,"\eta"] & \im(f) \arrow[d,swap,"h"] \arrow[r,"i_f"] & B \arrow[d,"\eta"] \\
        \trunc{0}{A} \arrow[r,swap,"q_{\trunc{0}{f}}"] & \im\trunc{0}{f} \arrow[r,swap,"i_{\trunc{0}{f}}"] & \trunc{0}{B}
      \end{tikzcd}
    \end{equation*}
    commute, and show that $h$ is a set truncation of $\im(f)$.
  \end{subexenum}
  \exitem Consider a type $A$, and suppose that $\trunc{0}{A}$ is a finite type with $k$ elements. Show that there exists a map $f:\Fin{k}\to A$ such that $\eta\circ f$ is an equivalence, i.e., prove the proposition
  \begin{equation*}
    \exists_{(f:\Fin{k}\to A)}\isequiv(\eta\circ f).
  \end{equation*}
  \exitem Consider a type $A$ and a universe $\UU$ containing $A$. Let
  \begin{equation*}
    \tau:A\to ((A\to\prop_\UU)\to\prop_\UU)
  \end{equation*}
  be the map defined by $\tau(a)\defeq\lam{f}f(a)$. Show that the map
  \begin{equation*}
    q_\tau : A\to\im(\tau)
  \end{equation*}
  obtained from the image factorization of $A$ is a set truncation of $A$.
  \exitem \label{ex:weakly-path-constant}A map $f:A \to B$ is called \define{weakly path-constant}\index{weakly path constant map|textbf} if it comes equipped with an element of type
  \begin{equation*}
    \isweaklypathconstant(f) : \prd{x,y:A}\prd{p,q:x=y}\ap{f}{p}=\ap{f}{q}.
  \end{equation*}
  In other words, $f$ is weakly path-constant if for each $x,y:A$ the map $\apfunc{f}:(x=y)\to (f(x)=f(y))$ is weakly constant in the sense of \cref{defn:weakly-constant}.
  \begin{subexenum}
  \item Show that every map $\trunc{0}{A}\to B$ is weakly path-constant. Use this to obtain a map
    \begin{equation*}
      \alpha : \Big(\trunc{0}{A}\to B\Big)\to\Big(\sm{f:A\to B}\isweaklypathconstant(f)\Big).
    \end{equation*}
  \item Show that if $B$ is a $1$-type, then the map $\alpha$ is an equivalence. In other words, show that every weakly path-constant map $f:A\to B$ into a $1$-type $B$ has a unique extension
    \begin{equation*}
      \begin{tikzcd}
        A \arrow[r,"f"] \arrow[d,swap,"\eta"] & B. \\
        \trunc{0}{A} \arrow[ur,dashed]
      \end{tikzcd}
    \end{equation*}
  \end{subexenum}
  \exitem Consider two universes $\UU$ and $\VV$. Use the type theoretic replacement axiom to show that the type of locally $\VV$-small types in $\UU$ is equivalent to the type\index{locally small type}
  \begin{equation*}
    \sum_{(Y:\mathsf{Locally\usc{}}\VV\mathsf{\usc{}Small\usc{}Set}_\UU)}\Big(Y\to \sm{Z:\VV}\isconn(Z)\Big)
  \end{equation*}
  of locally $\VV$-small sets $Y$ in $\UU$ equipped with a family of connected types in $\VV$.
  \exitem Show that every finite type is uniquely a product of finitely many finite types of prime cardinality, in the sense that the type
  \begin{equation*}
    \sm{X:\F}\sm{Y:X\to\sm{p:\primeN}\BS_p}\Brck{A\simeq\prd{x:X}Y(x)}
  \end{equation*}
  is connected for every finite type $A$.
  \exitem \label{ex:stirling-type-of-the-second-kind}Consider two types $A$ and $B$. The \define{Stirling type of the second kind}\index{Stirling type of the second kind|textbf}\index{{{A B}}@{$\stirling{A}{B}$}|see {Stirling type of the second kind}} is the type
  \begin{equation*}
    \stirling{A}{B}:=\sm{X:\UU_B}A\twoheadrightarrow X. 
  \end{equation*}
  \begin{subexenum}
  \item Show that if $B$ is a $k$-type, then the type $\stirling{A}{B}$ is also a $k$-type.\index{Stirling type of the second kind!is truncated}\index{truncated type!Stirling type of the second kind}
  \item Suppose that $B$ has decidable equality. Construct an equivalence
  \begin{equation*}
    \stirling{A+\unit}{B+\unit}\simeq (B+\unit)\times\stirling{A}{B+\unit}+\stirling{A}{B}
  \end{equation*}
  \item Suppose that $A$ and $B$ are finite types of cardinality $n$ and $k$. Show that the Stirling type $\stirling{A}{B}$ of the second kind is a finite type of cardinality $\stirling{n}{k}$, where $\stirling{n}{k}$ is the \define{Stirling number of the second kind}.\index{Stirling type of the second kind!is finite}\index{finite type!Stirling type of the second kind}
  \end{subexenum}
  \exitem \label{ex:distributive-pi-coprod}In this exercise we extend the definition of the binomial types to $\trunc{0}{\UU}$ as follows: For a type $X:\UU$ and $k:\trunc{0}{\UU}$, we define\index{binomial type!extended definition|textbf}
  \begin{equation*}
    \binomtype{X}{k}\defeq \sm{Y:\fib{\eta}{k}}Y\demb X.
  \end{equation*}
  Furthermore, for $(Y,i):\binomtype{X}{k}$, define
  \begin{align*}
    X\setminus Y & \defeq \sm{x:X}\neg(\fib{i}{x}). \\
    \complement{i} & \defeq \proj 1.
  \end{align*}
  Now consider a type $X$ and two type families $A$ and $B$ over $X$, and let $\UU$ be a universe containing $X$, $A$, and $B$. Show that the type $\prd{x:X}A(x)+B(x)$ is equivalent to the type\index{distributivity!of P over coproducts@{of $\Pi$ over coproducts}}\index{dependent function type!distributivity of P over coproducts@{distributivity of $\Pi$ over coproducts}}
  \begin{equation*}
    \sm{k:\trunc{0}{\UU}}\sm{(Y,i):\dbinomtype{X}{k}}\Big(\prd{y:Y}A(i(y))\Big)\times\Big(\prd{y:X\setminus Y}B(\complement{i}(y))\Big).
  \end{equation*}
\end{exercises}
\index{set quotient|)}


\section{Groups in univalent mathematics}\label{sec:groups}
\index{group|(}

In this section we demonstrate a very common way to use the univalence axiom\index{univalence axiom}, showing that isomorphic groups can be identified. When you introduce a certain kind of structure in type theory, such as groups or rings, you automatically obtain the type of all such structures. In other words, we define what a group is by defining the type of all groups, we define what a ring is by defining the type of all rings, and so on. The elements of the type of all groups are of course groups, such as the group of integers, integers modulo $k$, automorphism groups, and so on. The next important question is how two elements in the type of groups can be identified. This question is answered with the help of the univalence axiom: isomorphic groups can be identified. This is an instance of the \emph{structure identity principle}\index{structure identity principle}, which we covered in \cref{sec:structure-identity-principle}.

Identifiying isomorphic groups is a common \emph{informal} practice in classical mathematics. For example, by the third isomorphism theorem we have an isomorphism
\begin{equation*}
  (G/N)/(K/N)\cong (G/K)
\end{equation*}
for any sequence $N \trianglelefteq K \trianglelefteq G$ of normal subgroups of $G$, and it is common to simply write $(G/N)/(K/N)=G/K$. Of course, classical mathematicians know that this convention is incompatible with the axioms of Zermelo-Fraenkel set theory, but that does not stop them from applying this useful abuse of notation. In univalent mathematics we make this informal practice precise and formal.

\subsection{The type of all groups}

In order to efficiently characterize the identity type of the type of all groups in a universe $\UU$, we introduce the type of groups in two stages: first we introduce the type of \emph{semigroups}, and then we introduce groups as semigroups that possess a unit element and inverses. Since semigroups can have at most one unit element and since elements of semigroups can have at most one inverse, it follows that the type of groups is a subtype of the type of semigroups, and this will help us with the characterization of the identity type of the type of all groups.

\begin{rmk}
  In order to show that isomorphic (semi)groups can be identified, it has to be part of the definition of a (semi)group that its underlying type is a set. This is an important observation: in many branches of algebra the objects of study are \emph{set-level} structures\index{set-level structure}.

  A notable exception is formed by categories, which are objects at truncation level $1$, i.e., at the level of \emph{groupoids}. We will not cover categories in this book. For more about categories we recommend Chapter 9 of \cite{hottbook}.
\end{rmk}

\begin{defn}
  A \define{semigroup}\index{semigroup|textbf} in a universe $\UU$ is a triple $(G,\mu,\alpha)$ consisting of a set $G$ in $\UU$ equipped with a binary operation $\mu:G\to (G\to G)$ and a homotopy
  \begin{equation*}
    \alpha : \prd{x,y,z:G}\mu(\mu(x,y),z)=\mu(x,\mu(y,z))
  \end{equation*}
  witnessing that $\mu$ is \define{associative}\index{associative|textbf}.
  We write $\semigroup_\UU$\index{Semigroup@{$\semigroup_\UU$}|textbf} for the type of all semigroups in $\UU$, i.e., for the type
  \begin{equation*}
    \sm{G:\Set_\UU}\sm{\mu:G\to(G\to G)}\prd{x,y,z:G}\mu(\mu(x,y),z)=\mu(x,\mu(y,z)).
  \end{equation*}
\end{defn}

\begin{defn}
  A semigroup $G$ is said to be \define{unital}\index{semigroup!unital}\index{unital semigroup} if it comes equipped with a \define{unit}\index{unit!of a unital semigroup} $e:G$ that satisfies the left and right unit laws\index{unit laws!for a unital semigroup}
  \begin{align*}
    \leftunit : \prd{y:G}\mu(e,y)=y \\
    \rightunit : \prd{x:G}\mu(x,e)=x.
  \end{align*}
  We write $\isunital(G)$\index{is-unital@{$\isunital$}} for the type of such triples $(e,\leftunit,\rightunit)$. Unital semigroups are also called \define{monoids}\index{monoid|textbf}, so we define\index{Monoid@{$\monoid_\UU$}}
  \begin{equation*}
    \monoid_\UU\defeq\sm{G:\semigroup_\UU}\isunital(G).
  \end{equation*}
\end{defn}

The unit of a semigroup is of course unique once it exists. In univalent mathematics we express this fact by asserting that the type $\isunital(G)$ is a proposition for each semigroup $G$. In other words, being unital is a \emph{property} of semigroups rather than structure on it. This is typical for univalent mathematics: we express that a structure is a property by proving that this structure is a proposition.

\begin{lem}
  For a semigroup $G$ the type $\isunital(G)$ is a proposition.\index{is-unital@{$\isunital$}!is a proposition}
\end{lem}

\begin{proof}
  Let $G$ be a semigroup. Note that since $G$ is a set, it follows that the types of the left and right unit laws are propositions. Therefore it suffices to show that any two elements $e,e':G$ satisfying the left and right unit laws can be identified. This is easy:
  \begin{equation*}
    e = \mu(e,e') = e'.\qedhere
  \end{equation*}
\end{proof}

\begin{defn}
  Let $G$ be a unital semigroup. We say that $G$ \define{has inverses}\index{unital semigroup!has inverses}\index{semigroup!has inverses} if it comes equipped with an operation $x\mapsto x^{-1}$ of type $G\to G$, satisfying the left and right inverse laws\index{inverse laws!for a group}
  \begin{align*}
    \leftinv & : \prd{x:G}\mu(x^{-1},x)=e \\
    \rightinv & : \prd{x:G}\mu(x,x^{-1}) = e.
  \end{align*}
  We write $\isgroup'(G,e)$\index{is-group'@{$\isgroup'$}|textbf} for the type of such triples $((\blank)^{-1},\leftinv,\rightinv)$, and we write\index{is-group@{$\isgroup$}|textbf}
  \begin{equation*}
    \isgroup(G)\defeq\sm{e:\isunital(G)}\isgroup'(G,e)
  \end{equation*}
  A \define{group}\index{group|textbf} is a unital semigroup with inverses. We write $\group$\index{Group@{$\group_\UU$}|textbf} for the type of all groups in $\UU$.
\end{defn}

\begin{lem}
  For any semigroup $G$ the type $\isgroup(G)$ is a proposition.\index{is-group@{$\isgroup$}!is a proposition}
\end{lem}

\begin{proof}
  We have already seen that the type $\isunital(G)$ is a proposition. Therefore it suffices to show that the type $\isgroup'(G,e)$ is a proposition\index{is-group'@{$\isgroup'$}!is a proposition} for any $e:\isunital(G)$.

  Since a semigroup $G$ is assumed to be a set, we note that the types of the inverse laws are propositions. Therefore it suffices to show that any two inverse operations satisfying the inverse laws are homotopic.

  Let $x\mapsto x^{-1}$ and $x\mapsto x^{-1'}$ be two inverse operations on a unital semigroup $G$, both satisfying the inverse laws. Then we have the following identifications
  \begin{align*}
    x^{-1} & = \mu(e,x^{-1}) \\
    & = \mu(\mu(x^{-1'},x),x^{-1}) \\
    & = \mu(x^{-1'},\mu(x,x^{-1})) \\
    & = \mu(x^{-1'},e) \\
    & = x^{-1'}
  \end{align*}
  for any $x:G$. Thus the two inverses of $x$ are the same, and the claim follows.
\end{proof}

\begin{eg}
  The type $\Z$ of integers\index{Z@{$\Z$}!is a group}\index{group!Z@{$\Z$}} has the structure of a group, with the group operation being addition. The fact that $\Z$ is a set was shown in \cref{ex:set_coprod}, and the group laws were shown in \cref{ex:int_group_laws}. 
\end{eg}

\begin{eg}
  Given a set $X$, we define  the \define{automorphism group}\index{automorphism group|textbf}\index{group!automorphism group of a set|textbf}\index{set!automorphism group|textbf} of $X$ by\index{Aut(X)@{$\Aut(X)$}|see {automorphism group}}
  \begin{equation*}
    \Aut(X)\defeq (X\simeq X).
  \end{equation*}
  The group operation of $\Aut(X)$ is given by composition of equivalences, and the unit of the group is the identity function. An important special case of the automorphism groups is the \define{symmetric group}\index{symmetric group|textbf}\index{S n@{$S_n$}|see {symmetric group}}\index{group!S n@{$S_n$}|textbf}
  \begin{equation*}
    S_n\defeq \Aut(\Fin{n}).
  \end{equation*}
\end{eg}

\subsection{Group homomorphisms}

\begin{defn}
  Let $G$ and $H$ be (semi)groups. A \define{homomorphism}\index{homomorphism!of semigroups|textbf}\index{semigroup!homomorphism|textbf}\index{homomorphism!of groups|textbf}\index{group!homomorphism|textbf}\index{group homomorphism|textbf}\index{semigroup homomorphism|textbf} of (semi)groups from $G$ to $H$ is a pair $(f,\mu_f)$ consisting of a function $f:G\to H$ between their underlying types, and a homotopy
  \begin{equation*}
    \mu_f:\prd{x,y:G} f(\mu_G(x,y))=\mu_H(f(x),f(y))
  \end{equation*}
  witnessing that $f$ preserves the binary operation of $G$. We will write\index{hom(G,H) for semigroups@{$\hom(G,H)$ for semigroups}|textbf}\index{hom(G,H) for groups@{$\hom(G,H)$ for groups}|textbf}
  \begin{equation*}
    \hom(G,H)
  \end{equation*}
  for the type of all (semi)group homomorphisms from $G$ to $H$.
\end{defn}

\begin{rmk}\label{rmk:is-set-hom-semigroup}
  Since it is a property for a function to preserve the multiplication of a semigroup, it follows easily that equality of semigroup homomorphisms is equivalent to the type of homotopies between their underlying functions. In particular, it follows that the type of homomorphisms of semigroups is a set.
\end{rmk}

\begin{rmk}\label{rmk:category-semigroup}
  The \define{identity homomorphism}\index{identity homomorphism!of semigroups|textbf}\index{identity homomorphism!of groups|textbf} on a (semi)group $G$ is defined to be the pair consisting of
  \begin{align*}
    \idfunc & : G \to G \\
    \lam{x}\lam{y}\refl{} & : \prd{x,y:G} \mu_G(x,y) = \mu_G(x,y).
  \end{align*}
  Let $f:G\to H$ and $g:H\to K$ be (semi)group homomorphisms. Then the composite function $g\circ f:G\to K$ is also a (semi)group homomorphism\index{composition!of semigroup homomorphisms|textbf}\index{composition!of group homomorphisms|textbf}, since we have the identifications
  \begin{equation*}
    \begin{tikzcd}
      {g(f(\mu_G(x,y)))} \arrow[r,equals] & {g(\mu_H(f(x),f(y)))} \arrow[r,equals] & {\mu_K(g(f(x)),g(f(y)))}.
    \end{tikzcd}
  \end{equation*}
  Since the identity type of (semi)group homomorphisms is equivalent to the type of homotopies between (semi)group homomorphisms it is easy to see that (semi)group homomorphisms satisfy the laws of a category, i.e., that we have the identifications
  \begin{align*}
    \idfunc\circ f & = f \\
    g\circ \idfunc & = g \\
    (h\circ g) \circ f & = h \circ (g \circ f)
  \end{align*}
  for any composable (semi)group homomorphisms $f$, $g$, and $h$.
\end{rmk}

\begin{defn}
Let $h:\hom(G,H)$ be a homomorphism of (semi)groups. Then $h$ is said to be an \define{isomorphism}\index{group!isomorphism|textbf}\index{isomorphism!of groups|textbf}\index{semigroup!isomorphism|textbf}\index{isomorphism!of semigroups|textbf} if it comes equipped with an element of type $\isiso(h)$\index{is-iso@{$\isiso(h)$}!for semigroup homomorphisms|textbf}\index{is-iso@{$\isiso(h)$}!for group homomorphisms|textbf}, consisting of triples $(h^{-1},p,q)$ consisting of a homomorphism $h^{-1}:\hom(H,G)$ of semigroups and identifications
\begin{equation*}
p:h^{-1}\circ h=\idfunc[G]\qquad\text{and}\qquad q:h\circ h^{-1}=\idfunc[H]
\end{equation*}
witnessing that $h^{-1}$ satisfies the inverse laws\index{inverse laws!for semigroup isomorphisms}\index{inverse laws!for group isomorphisms}We write $G\cong H$ for the type of all isomorphisms of semigroups from $G$ to $H$, i.e.,
\begin{equation*}
G\cong H \defeq \sm{h:\hom(G,H)}\sm{k:\hom(H,G)} (k\circ h = \idfunc[G])\times (h\circ k=\idfunc[H]).
\end{equation*}
\end{defn}

If $f$ is an isomorphism, then its inverse is unique. In other words, being an isomorphism is a property.

\begin{lem}
  For any semigroup homomorphism $h:\hom(G,H)$, the type
  \begin{equation*}
    \isiso(h)
  \end{equation*}
  is a proposition.\index{is-iso@{$\isiso(h)$}!is a proposition} It follows that the type $G\cong H$ is a set for any two semigroups $G$ and $H$.
\end{lem}

\begin{proof}
  Let $k$ and $k'$ be two inverses of $h$. In \cref{rmk:is-set-hom-semigroup} we have observed that the type of semigroup homomorphisms between any two semigroups is a set. Therefore it follows that the types $h\circ k=\idfunc$ and $k\circ h=\idfunc$ are propositions, so it suffices to check that $k=k'$. In \cref{rmk:is-set-hom-semigroup} we also observed that the equality type $k=k'$ is equivalent to the type of homotopies $k\htpy k'$ between their underlying functions. We construct a homotopy $k\htpy k'$ by the usual argument:
  \begin{equation*}
    \begin{tikzcd}
      k(y) \arrow[r,equals] & k(h(k'(y)) \arrow[r,equals] & k'(y).
    \end{tikzcd}\qedhere
  \end{equation*}
\end{proof}

\subsection{Isomorphic groups are equal}

\begin{lem}\label{lem:grp_iso}
  A (semi)group homomorphism $h:\hom(G,H)$ is an isomorphism if and only if its underlying map is an equivalence. Consequently, there is an equivalence
  \begin{equation*}
    (G\cong H)\simeq \sm{e:G\simeq H}\prd{x,y:G}e(\mu_G(x,y))=\mu_H(e(x),e(y))
  \end{equation*}
\end{lem}

\begin{proof}
  If $h:\hom(G,H)$ is an isomorphism, then the inverse semigroup homomorphism also provides an inverse of the underlying map of $h$. Thus we obtain that $h$ is an equivalence. For the converse, suppose that the underlying map of $f:G\to H$ is an equivalence. Then its inverse is also a semigroup homomorphism, since we have
  \begin{align*}
    f^{-1}(\mu_H(x,y)) & = f^{-1}(\mu_H(f(f^{-1}(x)),f(f^{-1}(y)))) \\
               & = f^{-1}(f(\mu_G(f^{-1}(x),f^{-1}(y)))) \\
               & = \mu_G(f^{-1}(x),f^{-1}(y)). \qedhere
  \end{align*}
\end{proof}

\begin{defn}
Let $G$ and $H$ be a semigroups in a univalent universe $\UU$. We define the family of maps\index{iso-eq for semigroups@{$\isoeq$ for semigroups}}
\begin{equation*}
\isoeq : (G=H)\to (G\cong H)
\end{equation*}
indexed by $H:\semigroup_\UU$ by $\isoeq(\refl{})\defeq\idfunc[G]$.
\end{defn}

\begin{thm}\label{thm:iso-eq-semigroup}
Consider a semigroup $G$ in a univalent universe $\UU$. Then the family of maps\index{identity type!of Semigroup@{of $\semigroup_\UU$}}\index{Semigroup@{$\semigroup_\UU$}!identity type}\index{characterization of identity type!of Semigroup@{of $\semigroup_\UU$}}
\begin{equation*}
\isoeq : (G=H)\to (G\cong H)
\end{equation*}
indexed by $H:\semigroup_\UU$ is a family of equivalences.
\end{thm}

\begin{proof}
By the fundamental theorem of identity types \cref{thm:id_fundamental}\index{fundamental theorem of identity types} it suffices to show that the total space
\begin{equation*}
\sm{H:\semigroup_\UU}G\cong H
\end{equation*}
is contractible. Since the type of isomorphisms from $G$ to $H$ is equivalent to the type of equivalences from $G$ to $H$ it suffices to show that the type
\begin{equation*}
  \sm{H:\semigroup_\UU}\sm{e:\eqv{G}{H}}\prd{x,y:G}e(\mu_G(x,y))=\mu_{H}(e(x),e(y)))
\end{equation*}
is contractible. Since $\semigroup_\UU\jdeq\sm{H:\Set_\UU}\hasassociativemul(H)$ we are in position to apply the structure identity principle stated in \cref{thm:structure-identity-principle}. Note that $H\mapsto G\simeq H$ is an identity system on $\Set_\UU$ at the set $G$. By condition (v) of \cref{thm:structure-identity-principle} it therefore suffices to show that the type
\begin{equation*}
  \sm{\mu':\hasassociativemul(G)}\prd{x,y:G}\mu_G(x,y)=\mu'(x,y)
\end{equation*}
is contractible. This follows by function extensionality, since associativity of a binary operation on a set is a proposition.
\end{proof}

\begin{cor}
The type $\semigroup_\UU$ is a $1$-type.\index{Semigroup@{$\semigroup_\UU$}!is a 1-type@{is a $1$-type}}
\end{cor}

\begin{proof}
  The identity types of $\semigroup_\UU$ are sets because they are equivalent to the sets of isomorphisms between semigroups.
\end{proof}

We now turn to the proof that isomorphic groups are equal. Analogously to the map $\isoeq$ of semigroups, we have a map $\isoeq$ of groups. Note, however, that the domain of this map is now the identity type $G=H$ of the \emph{groups} $G$ and $H$, so the maps $\isoeq$ of semigroups and groups are not exactly the same maps.

\begin{defn}
  Let $G$ and $H$ be groups in a univalent universe $\UU$. We define the family of maps\index{iso-eq for groups@{$\isoeq$ for groups}}
  \begin{equation*}
    \isoeq : (G=H)\to (G\cong H)
  \end{equation*}
  indexed by $H:\Group_\UU$ by $\isoeq(\refl{})\defeq\idfunc[G]$.
\end{defn}

\begin{thm}
  For any two groups $G$ and $H$ in a univalent universe $\UU$, the map\index{identity type!of Group@{of $\group_\UU$}}\index{Group@{$\group_\UU$}!characterization of identity type}\index{characterization of identity type!of Group@{of $\Group_\UU$}}
  \begin{equation*}
    \isoeq:(G=H)\to (G\cong H)
  \end{equation*}
  is an equivalence.
\end{thm}

\begin{proof}
  Let $G$ and $H$ be groups in $\UU$, and write $UG$ and $UH$ for their underlying semigroups, respectively. Then we have a commuting triangle
  \begin{equation*}
    \begin{tikzcd}[column sep=0]
      (G=H) \arrow[rr,"\apfunc{\proj 1}"] \arrow[dr,swap,"\isoeq"] & & (UG=UH) \arrow[dl,"\isoeq"] \\
      \phantom{(UG=UH)} & (G\cong H)
    \end{tikzcd}
  \end{equation*}
  Since being a group is a property of semigroups it follows that the projection map $\group_\UU\to\semigroup_\UU$ forgetting the unit and inverses, is an embedding. Thus the top map in this triangle is an equivalence. The map on the right is an equivalence by \cref{thm:iso-eq-semigroup}, so the claim follows by the 3-for-2 property.
\end{proof}

\begin{cor}
  The type of groups is a $1$-type.\index{Group@{$\group_\UU$}!is a 1-type@{is a $1$-type}}
\end{cor}

\subsection{Homotopy groups of types}
\index{homotopy group|(}

Since the identity type gives every type groupoidal structure, we can construct for every type $A$ equipped with a base point $a:A$ a sequence of groups $\pi_n(A,a)$ indexed by $n\geq 1$. In order to construct this sequence of groups, we first define the \emph{loop space} operation, which takes pointed types to pointed types.

\begin{defn}
  The type of \define{pointed types}\index{pointed type|textbf} in a universe $\UU$ is defined as\index{U *@{$\UU_\ast$}|textbf}
  \begin{equation*}
    \UU_\ast\defeq\sm{X:\UU}X.
  \end{equation*}
  Given two pointed types $A$ and $B$ with base points $a$ and $b$ respectively, we define the type of \define{pointed maps}\index{pointed map|textbf}\index{A arrow* B@{$A\to_\ast B$}|see {pointed map}}
  \begin{equation*}
    (A\to_\ast B)\defeq\sm{f:A\to B}f(a)=b.
  \end{equation*}
\end{defn}

\begin{defn}\label{defn:loop-spaces}
  Consider a universe $\UU$. We define the \define{loop space}\index{loop space|textbf}\index{O (A)@{$\loopspace{A}$}|see {loop space}} operation
  \begin{equation*}
    \loopspacesym : \UU_\ast\to\UU_\ast
  \end{equation*}
  by $\loopspace{A,a}\defeq(a=a,\refl{})$. Furthermore, we define for every $A:\UU_\ast$ the \define{iterated loop space}\index{iterated loop space|textbf}\index{O n(A)@{$\loopspace[n]{A}$}|see {iterated loop space}} $\loopspace[n]{A}$ recursively by
  \begin{align*}
    \loopspace[0]{A}\defeq A \\
    \loopspace[n+1]{A}\defeq\loopspace{\loopspace[n]{A}}.
  \end{align*}
\end{defn}

\begin{eg}\label{eg:loop-spaces}
  If $A$ is a pointed $1$-type\index{pointed 1-type@{pointed $1$-type}}, then the loop space $\loopspace{A}$ is a set. Furthermore, it has the structure of a group. Its unit is $\refl{}$, and the group operation is given by concatenation of identifications. This satisfies the group laws, since the group laws are just a special case of the groupoid laws for identity types, constructed in \cref{sec:groupoid}. Thus we see that the loop space of a pointed $1$-type is a group\index{loop space!of a 1-type is a group@{of a $1$-type is a group}}\index{group!loop space of a 1-type@{loop space of a $1$-type}}.
\end{eg}

If $A$ is a pointed type, but not assumed to be $1$-truncated, then we can still get 

\begin{defn}
  Consider a pointed type $A$ with base point $a:A$, and let $n\geq 1$. Then we define the \define{$n$-th homotopy group}\index{homotopy group|textbf}\index{group!homotopy group|textbf} $\pi_n(A)$\index{p  n(A)@{$\pi_n(A)$}|see {homotopy group}} of $A$ at $a$ to be the group with underlying set
  \begin{equation*}
    \pi_n(A)\defeq\trunc{0}{\loopspace[n]{A}}
  \end{equation*}
  The unit of the group is $\eta(\refl{})$ and the group operation is the unique binary operation such that
  \begin{equation*}
    \eta(r)\eta(s)=\eta(\ct{r}{s})
  \end{equation*}
  for every $r,s:\loopspace[n]{A}$. The group $\pi_1(A)$\index{p  1(A)@{$\pi_1(A)$}|see {fundamental group}} of a pointed type is called the \define{fundamental group}\index{fundamental group|textbf}\index{group!fundamental group|textbf} of $A$ at its base point $a:A$.
\end{defn}

\begin{rmk}
  Note that for $n=0$, we can still define the set
  \begin{equation*}
    \pi_0(A)\defeq\trunc{0}{A}.
  \end{equation*}
  However, this set does not necessarily come equipped with the structure of a group.
\end{rmk}

\begin{prp}\label{prp:homotopy-group-loop-space}
  For any pointed type $A$ and any $n\geq 1$ we have an isomorphism
  \begin{equation*}
    \pi_{n+1}(A)\cong \pi_n(\loopspace{A}).
  \end{equation*}
\end{prp}

\begin{proof}
  First, observe that we have a pointed equivalence
  \begin{equation*}
    \loopspace{\loopspace[n]{A}}\equiv_\ast\loopspace[n]{\loopspace{A}}.
  \end{equation*}
  This equivalence is constructed by induction on $n$, and also preserves the concatenation operation. Using this equivalence, we obtain a group isomorphism
  \begin{equation*}
    \pi_{n+1}(A)\jdeq \trunc{0}{\loopspace{\loopspace[n]{A}}}\cong\trunc{0}{\loopspace[n]{\loopspace{A}}} \jdeq \pi_n(\loopspace{A}).\qedhere
  \end{equation*}
\end{proof}

Homotopy groups are important algebraic invariants of a type. For example, they can be used to show that two pointed types $A$ and $B$ are not equivalent by showing that two types $A$ and $B$ have non-isomorphic homotopy groups. The study of homotopy groups of types is an intricate and complicated subject, analogous to algebraic topology. Since the homotopy groups of types are obtained in such a canonical manner from the identity types, which are inductively generated by just the reflexivity identification, the subject of studying homotopy groups of types is also called \emph{synthetic homotopy theory}. In the final section of this book we will show that the fundamental group of the circle, which is introduced as a \emph{higher inductive type}, is $\Z$. In this section we will show that equivalent types have isomorphic homotopy groups, and that the homotopy groups $\pi_n(A)$ are abelian if $n\geq 2$.

\begin{defn}
  Consider a pointed map $f:A\to_\ast B$ between two pointed types $A$ and $B$, where $p:f(a)=b$. Then we define the pointed map\index{functorial action!of O@{of $\loopspacesym$}|textbf}
  \begin{equation*}
    \loopspace{f}:\loopspace{A}\to_\ast\loopspace{B}
  \end{equation*}
  by $\loopspace{f}(r)\defeq \ct{(\ct{p^{-1}}{\ap{f}{r}})}{p}$. The identification witnessing that this is indeed a pointed map is obtained from the fact that $\ap{f}{\refl{}}\jdeq\refl{}$ and $\ct{p^{-1}}{p}=\refl{}$.

  Similarly, we define $\loopspace[n]{f}:\loopspace[n]{A}\to_\ast\loopspace[n]{B}$ recursively by\index{functorial action!of O n@{of $\loopspacesym^n$}|textbf}
  \begin{align*}
    \loopspace[0]{f} & \defeq f \\
    \loopspace[n+1]{f} & \defeq \loopspace{\loopspace[n]{f}}.
  \end{align*}
  The functorial action of $\Omega^n$ together with the functorial action of set truncation yield a functorial action\index{functorial action!of p n@{of $\pi_n$}|textbf}
  \begin{equation*}
    \pi_n(f):\pi_n(A)\to\pi_n(B)
  \end{equation*}
  for every pointed map $f:A\to_\ast B$. 
\end{defn}

\begin{rmk}
  Since action of paths preserves path concatenation, it follows that $\Omega^n(f)$ preserves path concatenation, for each $n\geq 1$. Consequently, the maps
  \begin{equation*}
    \pi_n(f):\pi_n(A)\to\pi_n(B)
  \end{equation*}
  are group homomorphisms.
\end{rmk}

\begin{prp}
  Consider a pointed equivalence $e:A\simeq_\ast B$ between two pointed types $A$ and $B$. Then we obtain group isomorphisms
  \begin{equation*}
    \pi_n(e):\pi_n(A)\cong\pi_n(B)
  \end{equation*}
  for all $n\geq 1$.
\end{prp}

\begin{proof}
  For any pointed equivalence $e:A\simeq_\ast B$ it follows that $\pi_n(e)$ is also an equivalence. Using \cref{lem:grp_iso}, the claim now follows.
\end{proof}

\subsection{The Eckmann-Hilton argument}

\index{Eckmann-Hilton argument|(}
The Eckmann-Hilton argument is used to show that $\pi_n(A)$ is an abelian group for all $n\geq 2$. This is achieved by constructing an identification
\begin{equation*}
  \ct{p}{q}=\ct{q}{p}
\end{equation*}
for all $p,q:\loopspace[2]{A}$. Note that identification elimination is not immediately applicable here, since both $p$ and $q$ are identifications of type $\refl{a}=\refl{a}$ with neither endpoint free. Therefore, we must come up with something else.

\begin{defn}
  Consider a binary operation $f:A\to(B\to C)$. The \define{binary action on paths}\index{binary action on paths|textbf}\index{action on paths!binary action on paths|textbf}\index{action on paths!ap-binary@{$\apbinary_f$}|textbf} of $f$ is the family of functions\index{ap-binary@{$\apbinary_f$}|textbf}
  \begin{equation*}
    \apbinary_f:(x=x')\to ((y=y') \to (f(x,y)=f(x',y'))
  \end{equation*}
  indexed by $x,x':A$ and $y,y':B$ given by $\apbinary_f(\refl{},\refl{})\defeq\refl{}$.
\end{defn}

\begin{lem}\label{lem:laws-ap-binary}
  The binary action on paths of $f:A\to(B\to C)$ satisfies the following laws:
  \begin{align*}
    \apbinary_f(\refl{},q) & = \ap{f(x)}{q} \\
    \apbinary_f(p,\refl{}) & = \ap{f(\blank,y)}{p}
  \end{align*}
  and moreover both triangles in the following diagram commute:
  \begin{equation*}
    \begin{tikzcd}[column sep=10em,row sep=4em]
      f(x,y) \arrow[r,equals,"{\ap{f(\blank,y)}{p}}"] \arrow[d,equals,swap,"{\ap{f(x,\blank)}{q}}"] \arrow[dr,equals,"{\apbinary_f(p,q)}"] & f(x',y) \arrow[d,equals,"{\ap{f(x',\blank)}{q}}"] \\
      f(x,y') \arrow[r,equals,swap,"{\ap{f(\blank,y')}{p}}"] & f(x',y')
    \end{tikzcd}
  \end{equation*}
\end{lem}

\begin{proof}
  The proof is immediate by identification elimination on $p$ and $q$, where applicable.
\end{proof}

\begin{eg}
  One particular binary operation to which we can apply the binary action on paths is concatenation of identifications
  \begin{equation*}
    \ct{\blank}{\blank}:(x=y)\to((y=z)\to (x=z))
  \end{equation*}
  This results in the \define{horizontal concatenation}\index{identity type!horizontal concatenation|textbf}\index{horizontal concatenation|textbf} operation\index{r .h s@{$\ct[h]{r}{s}$}|see {horizontal concatenation}}
  \begin{equation*}
    \ct[h]{\blank}{\blank} : (p=p')\to ((q=q') \to (\ct{p}{q}=\ct{p'}{q'})).
  \end{equation*}
  In other words, for any two identifications $r:p=p'$ and $s:q=q'$ as in the diagram
  \begin{equation*}
    \begin{tikzcd}[column sep=huge]
      x \arrow[r,equals,bend left=30,"p",""{name=A,below}] \arrow[r,equals,bend right=30,""{name=B,above},"{p'}"{below}] \arrow[from=A,to=B,phantom,"r\Downarrow"] & y \arrow[r,equals,bend left=30,"q",""{name=C,below}] \arrow[r,equals,bend right=30,""{name=D,above},"{q'}"{below}] \arrow[from=C,to=D,phantom,"s\Downarrow"] & z.
    \end{tikzcd}
  \end{equation*}
  we obtain $\ct[h]{r}{s}\defeq\apbinary_{\ct{\blank}{\blank}}(r,s):\ct{p}{q}=\ct{p'}{q'}$. The \define{vertical concatenation}\index{identity type!vertical concatenation|textbf}\index{vertical concatenation|textbf} operation, which concatenates $r:p=p'$ and $r':p'=p''$ as in the diagram
  \begin{equation*}
    \begin{tikzcd}[column sep=7em]
      x \arrow[r,equals,bend left=60,"p",""{name=A,below}] \arrow[r,equals,""{name=B},""{name=E,below},"{p'}"{near end}] \arrow[r,equals,bend right=60,"{p''}"{below},""{name=F,above}] \arrow[from=A,to=B,phantom,"r\Downarrow"] \arrow[from=E,to=F,phantom,"{r'\Downarrow}"] 
      & y
    \end{tikzcd}
  \end{equation*}
  is given by ordinary concatenation of identifications.
\end{eg}

\begin{lem}\label{lem:unit-laws-horizontal-concat}
  Horizontal concatenation satisfies the following left and right unit laws:\index{unit laws!for horizontal concatenation}\index{horizontal concatenation!unit laws}
  \begin{align*}
    \ct[h]{\refl{\refl{}}}{s} & = s \\
    \ct[h]{r}{\refl{\refl{}}} & = r.
  \end{align*}
\end{lem}

\begin{proof}
  This follows by identification elimination on $r$ and $s$, or alternatively via \cref{lem:laws-ap-binary}.
\end{proof}

In the following lemma we establish the \define{interchange law} for horizontal and vertical concatenation.

\begin{lem}\label{lem:interchange-law}
Consider a diagram of the form\index{interchange law!of horizontal and vertical concatenation}\index{horizontal concatenation!interchange law}\index{vertical concatenation!interchange law}
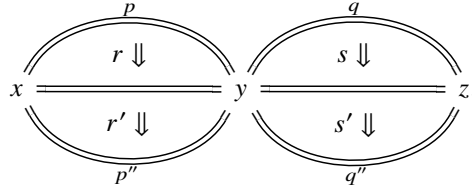
\begin{equation*}
\begin{tikzcd}[column sep=7em]
x \arrow[r,equals,bend left=60,"p",""{name=A,below}] \arrow[r,equals,""{name=B},""{name=E,below}] \arrow[r,equals,bend right=60,"{p''}"{below},""{name=F,above}] \arrow[from=A,to=B,phantom,"r\Downarrow"] \arrow[from=E,to=F,phantom,"{r'\Downarrow}"] 
& y \arrow[r,equals,bend left=60,"q",""{name=C,below}] \arrow[r,equals,""{name=G,above},""{name=H,below}] \arrow[r,equals,bend right=60,""{name=D,above},"{q''}"{below}] \arrow[from=C,to=G,phantom,"s\Downarrow"] \arrow[from=H,to=D,phantom,"{s'\Downarrow}"] & z.
\end{tikzcd}
\end{equation*}
Then there is an identification
\begin{equation*}
  \ct[h]{(\ct{r}{r'})}{(\ct{s}{s'})}=\ct{(\ct[h]{r}{s})}{(\ct[h]{r'}{s'})}.
\end{equation*}
\end{lem}

\begin{proof}
  We use path induction on both $r$ and $s$. Then it suffices to show that
  \begin{equation*}
    \ct[h]{(\ct{\refl{}}{r'})}{(\ct{\refl{}}{s'})}=\ct{(\ct[h]{\refl{}}{\refl{}})}{(\ct[h]{r'}{s'})}
  \end{equation*}
  Using the unit laws for ordinary concatenation, we see that both sides reduce to $\ct[h]{r'}{s'}$.
\end{proof}

\begin{thm}
  Consider a pointed type $A$, and let $r,s:\loopspace[2]{A}$. Then there is an identification
  \begin{equation*}
    \ct{r}{s}=\ct{s}{r}
  \end{equation*}
\end{thm}

\begin{proof}
  First we observe that $\ct{r}{s}=\ct[h]{r}{s}$ by the following calculation using the unit laws from \cref{lem:unit-laws-horizontal-concat} and the interchange law from \cref{lem:interchange-law}:
  \begin{align*}
    \ct{r}{s} & = \ct{(\ct[h]{r}{\refl{\refl{}}})}{(\ct[h]{\refl{\refl{}}}{s})} \\
              & = \ct[h]{(\ct{r}{\refl{\refl{}}})}{(\ct{\refl{\refl{}}}{s})} \\
              & = \ct[h]{r}{s}
  \end{align*}
  Similarly, we observe that $\ct[h]{r}{s}=\ct{s}{r}$ by the following calculation:
  \begin{align*}
    \ct[h]{r}{s} & = \ct[h]{(\ct{\refl{\refl{}}}{r})}{(\ct{s}{\refl{\refl{}}})} \\
                 & = \ct{(\ct[h]{\refl{\refl{}}}{s})}{(\ct[h]{r}{\refl{\refl{}}})} \\
                 & = \ct{s}{r}.
  \end{align*}
  These two calculations combined prove the claim.
\end{proof}

\begin{cor}
For $n\geq 2$, the $n$-th homotopy group of any pointed type is abelian.\index{homotopy group!is abelian for n geq 2@{is abelian for $n\geq 2$}}
\end{cor}

\begin{proof}
  By \cref{prp:homotopy-group-loop-space} it follows that $\pi_n(A)$ is isomorphic to the second homotopy group of some pointed type, for every $n\geq 2$. Therefore it suffices to prove the claim for $\pi_2(A)$ for every pointed type $A$.
  
  Our goal is to show that 
  \begin{equation*}
    \prd{r,s:\pi_2(A)} rs=sr.
  \end{equation*}
  Since we are constructing an identification in a set, we can use the dependent universal property of $0$-truncation on both $r$ and $s$, stated in \cref{thm:set-truncation}. Therefore it suffices to show that
  \begin{equation*}
    \prd{r,s:\loopspace[2]{A}} \eta(r)\eta(s)=\eta(s)\eta(r).
  \end{equation*}
  The claim now follows, because
  \begin{equation*}
    \eta(r)\eta(s)=\eta(\ct{r}{s})=\eta(\ct{s}{r})=\eta(s)\eta(r).\qedhere
  \end{equation*}
\end{proof}
\index{homotopy group|)}
\index{Eckmann-Hilton argument|)}

\subsection{Concrete versus abstract groups in univalent mathematics}

In univalent mathematics there is another exciting perspective on group theory. We won't be able to go in full details here, but we can sketch some of key ideas. To learn more about this beautiful univalent perspective on group theory, I recommend the forthcoming \emph{Symmetry} book \cite{symmetrybook}.

We saw in \cref{eg:loop-spaces} that for every pointed connected $1$-type $X$ we obtain a group with underlying type $\loopspace{X}$. All groups can be constructed in this way. In fact, for every group $G$ in $\UU$ the type
\begin{equation*}
  \sm{B:\mathsf{Pointed\usc{}Connected\usc{}}1\mathsf{\usc{}Type}_\UU}G\cong\loopspace{B}
\end{equation*}
of pointed connected $1$-types $B$ equipped with a group isomorphism from $G$ to $\loopspace{B}$ is contractible. We write $BG$ for the unique pointed connected $1$-type whose loop space is isomorphic to $G$. The pointed type $BG$ is also called the \define{delooping}\index{delooping|textbf}\index{group!delooping|textbf} of $G$, or the \define{classifying type}\index{classifying type|textbf}\index{group!classifying type|textbf} of $G$. The fact that the above type is contractible is of course heavily reliant on the univalence axiom.

\begin{eg}
  We have already seen that
  \begin{equation*}
    S_n\cong\loopspace{\BS_n},
  \end{equation*}
  i.e., that the loop space of the type of all finite types of cardinality $n$ is equivalent to the symmetric group $S_n$. The type $\BS_n$ is of course a pointed connected $1$-type, so $BS_n$ is indeed the classifying type of the symmetric group $S_n$.\index{BS n@{$\BS_n$}!is classifying type of symmetric group}
\end{eg}

Since the map
\begin{equation*}
  \loopspacesym:\mathsf{Pointed\usc{}Connected\usc{}}1\mathsf{\usc{}Type}_\UU\to\Grp_\UU
\end{equation*}
is an equivalence, we obtain two perspectives on the type of all groups. The elements of the type $\Grp_\UU$ are groups according to the traditional definition of groups. We call such groups \define{abstract groups}\index{group!abstract group|textbf}\index{abstract group|textbf}. On the other hand, pointed connected $1$-types $B$ are \define{concrete groups}\index{concrete group|textbf}\index{group!concrete group|textbf} in the sense that the contain an object $\ast:B$, and the group $B$ represents is the group of self-identifications (i.e., symmetries) of the base point $\ast:B$. Thus we see that when we present a group as a pointed connected $1$-type, then we \emph{concretely} manifest that group as the group of symmetries of some object.

We can also bring group homomorphisms into the mix: for every group homomorphism $f:G\to H$ the type of pointed maps $b:BG\to_\ast BH$ equipped with a homotopy witnessing 
\begin{equation*}
  \begin{tikzcd}
    G \arrow[d,swap,"\cong"] \arrow[r,"f"] & H \arrow[d,"\cong"] \\
    \loopspace{BG} \arrow[r,swap,"\loopspace{b}"] & \loopspace{BH}
  \end{tikzcd}
\end{equation*}
commutes is contractible. In other words, every group homomorphism $f:G\to H$ has a unique \define{delooping} $Bf:BG\to BH$.

We can do all of group theory in this way. For example, traditionally a $G$-set is defined to be a set $X$ equipped with a group homomorphism $G\to\Aut(X)$. That is, the type of \define{abstract $G$-sets}\index{abstract G-set@{abstract $G$-set}|textbf}\index{group!abstract G-set@{abstract $G$-set}|textbf} is defined to be\index{G-Set@{$G\mathsf{\usc{}}\Set_\UU$}|see {abstract $G$-set}}\index{G-Set@{$G\mathsf{\usc{}}\Set_\UU$}|textbf}
\begin{equation*}
  G\mathsf{\usc{}}\Set_\UU\defeq\sm{X:\Set_\UU}\hom(G,\Aut(X)).
\end{equation*}
However, this definition is equivalent to family $X:BG\to\Set_\UU$ of sets indexed by the classifying type $BG$. Therefore we define \define{concrete $G$-sets}\index{concrete G-set@{concrete $G$-set}|textbf}\index{group!concrete G-set@{concrete $G$-set}|textbf} to be type families $X:BG\to\Set_\UU$. Given a concrete $G$-set $X:BG\to\Set_\UU$, the set being acted upon is the set $X(\ast)$, and the action of $G$ on $X(\ast)$ is given by transport, since the elements of $G$ are equivalent to loops in $BG$.

The type of \define{orbits}\index{orbit}\index{concrete G-set@{concrete $G$-set}!orbit|textbf} of a concrete $G$-set $X:BG\to\Set_\UU$ can then be defined as
\begin{equation*}
  X/G\defeq \sm{u:BG}X(u)
\end{equation*}
and the type of \define{fixed points}\index{fixed point!of a concrete G-set@{of a concrete $G$-set}|textbf}\index{concrete G-set@{concrete $G$-set}!fixed point|textbf} of $X$ can be defined as
\begin{equation*}
  X_G\defeq \prd{u:BG}X(u).
\end{equation*}
To see that these definitions make sense, note that the fiber inclusion $X(\ast)\to X/G$ maps each element in the $G$-set $X$ to its orbit. The fiber inclusion is surjective by \cref{ex:is-surjective-fiber-inclusion}, and it maps two elements $x,y:X(\ast)$ to the same orbit precisely when there is a group element $g$ such that $gx=y$. Similarly, for the type of fixed points notice that each $x:X_G$ determines an element $x_\ast:X(\ast)$, which comes equipped with an identification
\begin{equation*}
  \apd{x}{g}:gx_\ast=x_\ast
\end{equation*}
since the group action of $G$ on $X$ is given by transport.

Also, notice that a subgroup $H$ of $G$ determines an inclusion homomorphism $i:H\to G$, and this inclusion function corresponds uniquely to a pointed map $Bi:BG\to BH$. Since $\loopspace{Bi}$ is an embedding, we note that $Bi$ must be a $0$-truncated map. Therefore, a concrete subgroup of a concrete group $BG$ is defined to be a concrete $G$-set $X$ such that the type of orbits is connected. Such concrete $G$-sets are called \define{transitive}\index{transitive concrete G-set@{transitive concrete $G$-set}|textbf}\index{concrete G-set@{concrete $G$-set}!transitive concrete G-set@{transitive concrete $G$-set}|textbf}.

Dually, we say that a concrete $G$-set $X$ is \define{free}\index{concrete G-set@{concrete $G$-set}!free concrete G-set@{free concrete $G$-set}|textbf}\index{free concrete G-set@{free concrete $G$-set}|textbf} if the type of orbits $X/G$ is a set. To see that this definition makes sense, we use the following generalization of the fundamental theorem of identity types:

\begin{thm}\label{thm:truncated-fundamental}
  Consider a connected type $A$ equipped with an element $a:A$, and consider a family of types $B(x)$ indexed by $x:A$. Then the following are equivalent:\index{fundamental theorem of identity types!generalization to truncated maps}\index{generalized fundamental theorem of identity types!truncated maps}
  \begin{enumerate}
  \item Every family of maps
    \begin{equation*}
      f:\prd{x:A}(a=x)\to B(x)
    \end{equation*}
    is a family of $k$-truncated maps.
  \item The total space
    \begin{equation*}
      \sm{x:A}B(x)
    \end{equation*}
    is $(k+1)$-truncated.
  \end{enumerate}
\end{thm}

\begin{proof}
  Recall from \cref{ex:is-trunc-const} that the total space $\sm{x:A}B(x)$ is $(k+1)$-truncated if and only if the base point inclusion
  \begin{equation*}
    (x,y):\unit\to\sm{x:A}B(x)
  \end{equation*}
  is $k$-truncated for every $(x,y):\sm{x:A}B(x)$. Since the type $A$ is assumed to be connected, this is equivalent to the condition that every base point inclusion of the form
  \begin{equation*}
    (a,y):\unit\to\sm{x:A}B(x)
  \end{equation*}
  is $k$-truncated. Base point inclusions of this form are homotopic to $\tot{f}$, where
  \begin{equation*}
    f:\prd{x:A}(a=x)\to B(x)
  \end{equation*}
  is given by $f(a,\refl{})\defeq y$. The condition that $\tot{f}$ is $k$-truncated is by \cref{lem:fib_total} equivalent to the condition that $f$ is a family of $k$-truncated maps. Furthermore, every family of maps $f:\prd{x:A}(a=x)\to B(x)$ is of the above form by the type theoretic Yoneda lemma (\cref{thm:yoneda}), completing the proof. 
\end{proof}

By the previous theorem it follows that if the type of orbits of a concrete $G$-set $X$ is a set, then the map $g\mapsto gx$ must be an embedding for every $x:X(\ast)$. In other words, the action of $G$ on $X$ is free.

\begin{rmk}
  \cref{thm:truncated-fundamental} can be generalized further. We include this generalization in \cref{ex:connected-fundamental}.
\end{rmk}

\begin{eg}
  Consider two sets $A$ and $B$, and a universe $\UU$ containing both of them. Then the automorphism group $\Aut(B)$ acts on the decidable embeddings $B\demb A$ by precomposition. Its type of orbits is the binomial type\index{binomial type}
  \begin{equation*}
    \dbinomtype{A}{B}\defeq\sm{X:\UU_B}X\demb A,
  \end{equation*}
  which we introduced in \cref{defn:binomial-type}. By \cref{prp:equiv-binom-type} it follows that $\dbinomtype{A}{B}$ is a set, so the action of $\Aut(B)$ on $B\demb A$ is free. Note that we didn't need to assume that $A$ and $B$ are sets: the action of $\Aut(B)$ on $B\demb A$ is always free.

  Similarly, we have an action of the automorphism group $\Aut(B)$ on the surjective maps $A\twoheadrightarrow B$ by postcomposition. Its type of orbits is the stirling type of the second kind\index{Stirling type of the second kind}
  \begin{equation*}
    \stirling{A}{B}\defeq\sm{X:\UU_B}A\twoheadrightarrow X,
  \end{equation*}
  which we introduced in \cref{ex:stirling-type-of-the-second-kind}. Assuming that $B$ is a set, it was shown in \cref{ex:stirling-type-of-the-second-kind} that $\stirling{A}{B}$ is a set. In other words, the action of $\Aut(B)$ on $A\twoheadrightarrow B$ is free.
\end{eg}

\begin{eg}
  In \cref{ex:prime} we introduced the type\index{prime number}
  \begin{equation*}
    \tilde{D}_n\defeq \sm{X:\BS_2}\sm{Y:X\to\F}\Big(\Fin{n}\simeq\prd{x:X}Y(x)\Big).
  \end{equation*}
  Notice that this type is the type of orbits of the $\Z/2$-set $D_n$ given by
  \begin{equation*}
    D_n(X)\defeq \sm{Y:X\to\F}\Big(\Fin{n}\simeq\prd{x:X}Y(x)\Big).
  \end{equation*}
  The fact that this is a family of sets is a nice exercise. Note that there is a surjective morphism of $\Z/2$-sets from $D_n(\Fin{2})$ to the $\Z/2$ set of divisors of $n$, where the action is given by $d\mapsto n/d$. The concrete $\Z/2$-action $D_n$ is transitive precisely when $n$ is either $1$ or a prime, and it is is free precisely when $n$ is not a square. Combining these two observations, we see that $n$ is prime if and only if this action is both transitive and free. In other words, $n$ is prime if and only if the type $\tilde{D}_n$ of orbits is contractible.
\end{eg}

$G$-sets which are both transitive and free are very special. Such $G$-sets are called \define{$G$-torsors}\index{torsor|textbf}\index{group!torsor|textbf}. Note that a $G$-set $X$ is a $G$-torsor if and only if the type of orbits $X/G$ is contractible. By the fundamental theorem of identity types, this implies that the family of maps
\begin{equation*}
  \prd{v:BG}(u=v)\to X(v)
\end{equation*}
is a family of equivalences, where $(u,x)$ is the center of contraction of $X/G$. It follows that a concrete $G$-set $X:BG\to\Set_\UU$ is a $G$-torsor if and only if it is in the image of
\begin{equation*}
  \idtypevar{} : BG\to (BG\to \Set_\UU).
\end{equation*}
However, we know from \cref{ex:idtype-is-emb} that this map is an embedding, so it follows that the type of concrete $G$-torsors is equivalent to $BG$. On the other hand, the type of concrete $G$-torsors is equivalent to the type of abstract $G$-torsors. This suggests that the classifying type $BG$ of any group $G$ can be constructed as the type of abstract $G$-torsors, and this is indeed one way to construct the classifying type of a group $G$.

\begin{exercises}
  \exitem Consider a set $X$ equipped with an associative binary operation $\mu:X\to (X\to X)$, and suppose that
  \begin{enumerate}
  \item The type $X$ is inhabited, i.e., $\brck{X}$ holds.
  \item The maps $\mu(x,\blank)$ and $\mu(\blank,y)$ are equivalences, for each $x,y:X$.
  \end{enumerate}
  Show that $X$ is a group.\index{group}
  \exitem Let $f:\hom(G,H)$ be a group homomorphism\index{homomorphism!of groups}. Show that $f$ preserves units and inverses, i.e., show that\index{group homomorphism!preserves units and inverses}
  \begin{align*}
    f(e_G) & = e_H \\
    f(x^{-1}) & = f(x)^{-1}.
  \end{align*}
  \exitem \label{ex:groupop-embedding}
  Consider a group $G$. Show that the function
  \begin{equation*}
    \mu_G:G\to (G\simeq G)
  \end{equation*}
  is an injective group homomorphism.
  \exitem Let $X$ be a set. Show that the map\index{equiv-eq@{$\equiveq$}!is a group isomorphism}
  \begin{equation*}
    \equiveq : (X=X)\to (\eqv{X}{X})
  \end{equation*}
  is a group isomorphism.
  \exitem Consider a group $G$. Show that the map\index{Z@{$\Z$}!is the free group with one generator|textbf}\index{group!free group with one generator|textbf}
  \begin{equation*}
    \Grp(\Z,G)\to G
  \end{equation*}
  given by $h\mapsto h(\oneZ)$, is an equivalence. In other words, the group $\Z$ satisfies the universal property of the \define{free group on one generator}\index{free group with one generator|textbf}.
  \exitem Give a direct proof and a proof using the univalence axiom of the fact that all semigroup isomorphisms between unital semigroups preserve the unit. Conclude that isomorphic monoids are equal.\index{isomorphism!of semigroups!preserves unit}\index{characterization of identity type!of Monoid@{of $\monoid_\UU$}}\index{identity type!of Monoid@{of $\monoid_\UU$}}
  \exitem \label{ex:dihedral-group}Consider an abelian group $A$, and let $D_A\defeq A+A$\index{D A@{$D_A$}|see {generalized dihedral group}}\index{D A@{$D_A$}|textbf} be the set equipped with $1\defeq\inl(0)$, the binary operation ${\blank}\cdot{\blank}:D_A\to (D_A\to D_A)$ defined by
    \begin{align*}
    \inl(x)\cdot\inl(y) & \defeq \inl(x+y) \\
    \inl(x)\cdot\inr(y) & \defeq \inr(-x+y) \\
    \inr(x)\cdot\inl(y) & \defeq \inr(x+y) \\
    \inr(x)\cdot\inr(y) & \defeq \inl(-x+y),
  \end{align*}
  and the unary operation $(\blank)^{-1}:D_A\to D_A$ defined by
  \begin{align*}
    \inl(x)^{-1} & \defeq \inl(-x) \\
    \inr(x)^{-1} & \defeq \inr(x).
  \end{align*}
  Show that $D_A$ equipped with these operations is a group. The group $D_A$ is called the \define{generalized dihedral group}\index{generalized dihedral group|textbf}\index{group!generalized dihedral group|textbf} on $A$. The (ordinary) \define{dihedral group}\index{dihedral group|textbf}\index{group!dihedral group|textbf} $D_k$\index{D k@{$D_k$}|see {dihedral group}}\index{D k@{$D_k$}|textbf} is defined to be $D_k\defeq D_{\Z/k}$.
  \exitem Recall that a \define{subgroup}\index{subgroup|textbf}\index{group!subgroup|textbf} of a group $G$ in $\UU$ consists of a subtype
  \begin{equation*}
    P:G\to\prop_\UU
  \end{equation*}
  such that $P$ contains the unit and is closed under the group operation and under inverses.
  \begin{subexenum}
  \item Consider a proposition $P$, and let $N_P$ be the subtype of $\Z/2$ given by
    \begin{equation*}
      N_P(x)\defeq (x=0)\vee P.
    \end{equation*}
    Show that $N_P$ is a subgroup of $\Z/2$.
  \item Show that the map $P\mapsto N_P$ is an embedding
    \begin{equation*}
      \prop_\UU\hookrightarrow\subgroup_\UU(\Z/2).
    \end{equation*}
  \end{subexenum}
  \exitem Recall that a \define{normal subgroup}\index{normal subgroup|textbf}\index{group!normal subgroup|textbf} $H$ of a group $G$ is a subgroup of $G$ such that $xyx^{-1}$ is in $H$ for every $y:H$ and $x:G$. Show that the type of normal subgroups of $G$ in $\UU$ is equivalent to the type
  \begin{equation*}
    \sm{H:\Grp_\UU}\sm{f:\hom(G,H)}\issurj(f).
  \end{equation*}
  \exitem \label{ex:commutative-binary-operations}For any type $A$, we define the type of \define{commutative binary operations}\index{commutative binary operation|textbf} on $A$ to be
  \begin{equation*}
    \Big(\sm{X:\BS_2}A^X\Big)\to A.
  \end{equation*}
  If $A$ is a set, show that the map
  \begin{equation*}
    \Big(\Big(\sm{X:\BS_2}A^X\Big)\to A\Big)\to\Big(\sm{f:A\to (A\to A)}\prd{x,y:A}f(x,y)=f(y,x)\Big)
  \end{equation*}
  given by $h\mapsto \lam{x}\lam{y}h(\Fin{2},(x,y))$ is an equivalence. In other words, show that every commutative operation $f:A\to(A\to A)$ extends uniquely along the map $f\mapsto(\Fin{2},f)$ as in the diagram
  \begin{equation*}
    \begin{tikzcd}
      A^{\Fin{2}} \arrow[dr,"\mu"] \arrow[d,swap,"f\mapsto{(\Fin{2},f)}"] \\
      \sm{X:\BS_2}A^X \arrow[r,dashed] & A.
    \end{tikzcd}
  \end{equation*}
  Give an informal explanation of this fact in terms fixed points of the concrete $\Z/2$-action on the set of binary operations $A\to (A\to A)$.
  \exitem Consider a commutative monoid $M$. Define an operation
  \begin{equation*}
    \prd{X:\mathbb{F}} M^X\to M
  \end{equation*}
  that extends the (binary) monoid operation to the finite unordered $n$-tuples of elements in $M$.
  \exitem Show that the type of $3$-element groups is equivalent to the type of $2$-element types.
  \exitem Show that the number of connected components in the type of all groups of order $n$ is as follows, for $n\leq 8$:
  \begin{center}
    \begin{tabular}{rllllllll}
      \emph{order:} & 1 & 2 & 3 & 4 & 5 & 6 & 7 & 8 \\
      \midrule
      \emph{number of groups:} & 1 & 1 & 1 & 2 & 1 & 2 & 1 & 5
    \end{tabular}
  \end{center}
  \exitem \label{ex:connected-fundamental}Consider a subtype\index{fundamental theorem of identity types!full generalization}\index{generalized fundamental theorem of identity types}
  \begin{equation*}
    P:\UU\to\prop_\VV
  \end{equation*}
  of a universe $\UU$. We say that a type $A:\UU$ is a \define{$P$-type}\index{P-type@{$P$-type}|textbf} if $P(A)$ holds, we say that a map $f:A\to B$ is a \define{$P$-map}\index{P-map@{$P$-map}|textbf} if its fibers are $P$-types, and we say that $A$ is \define{$P$-separated}\index{P-separated type@{$P$-separated type}|textbf} if its identity types are $P$-types.
  
  Now consider a connected type $A:\UU$ equipped with an element $a:A$, and consider a family of types $B(x):\UU$ indexed by $x:A$. Show that the following are equivalent:
  \begin{enumerate}
  \item Every family of maps
    \begin{equation*}
      f:\prd{x:A}(a=x)\to B(x)
    \end{equation*}
    is a family of $P$-maps.
  \item The total space
    \begin{equation*}
      \sm{x:A}B(x)
    \end{equation*}
    is $P$-separated.
  \end{enumerate}
  For readers familiar with the notion of $k$-connectedness: Conclude that every $f:\prd{x:A}(a=x)\to B(x)$ is a family of $k$-connected maps if and only if $\sm{x:A}B(x)$ is a $(k+1)$-connected type.
  \exitem Consider a group $G$ in a universe $\UU$ and a pointed connected $1$-type $B$. In analogy with \cref{thm:quotient_up}, show that the following are equivalent:
  \begin{enumerate}
  \item The pointed connected $1$-type $B$ comes equipped with a group homomorphism\index{universal property!of the classifying type of a group|textbf}\index{classifying type!universal property|textbf}
    \begin{equation*}
      \varphi:G \to \loopspace{B}
    \end{equation*}
    and for every pointed connected $1$-type $C$ that comes equipped with a group homomorphism $\psi:G\to \loopspace{C}$ there is a unique pointed map $f:B\to_\ast C$ equipped with a homotopy witnessing that the triangle
    \begin{equation*}
      \begin{tikzcd}[column sep=tiny]
        & G \arrow[dl,swap,"\varphi"] \arrow[dr,"\psi"] \\
        \loopspace{B} \arrow[rr,swap,"\loopspace{f}"] & & \loopspace{C}
      \end{tikzcd}
    \end{equation*}
    commutes.
  \item The pointed connected $1$-type $B$ comes equipped with a group isomorphism
    \begin{equation*}
      \varphi:G\cong \loopspace{B}.
    \end{equation*}
  \item There is an embedding $i:B\hookrightarrow G\mathsf{\usc{}}\Set_\UU$ such that the triangle
    \begin{equation*}
      \begin{tikzcd}[column sep=tiny]
        \unit \arrow[rr] \arrow[dr,swap,"\mathsf{Pr}_G"] & & B \arrow[dl,"i"] \\
        & G\mathsf{\usc{}}\Set_\UU
      \end{tikzcd}
    \end{equation*}
    commutes, where $\mathsf{Pr}_G$ is the \define{principal $G$-set}\index{principal G-set@{principal $G$-set}|textbf}\index{group!principal G-set@{principal $G$-set}|textbf}, i.e., $G$ acting on itself from the left.
  \end{enumerate}
  \exitem Consider a group $G$ and a pointed connected $1$-type $B$ equipped with a group isomorphism
  \begin{equation*}
    \varphi:G\cong \loopspace{B}.
  \end{equation*}
  \begin{subexenum}
  \item Show that the map
    \begin{equation*}
      \ev_\ast:(B\to\Set_\UU)\to \sm{X:\Set_\UU}\hom(G,\Aut(X))
    \end{equation*}
    sending concrete $G$-sets\index{concrete G-set@{concrete $G$-set}}\index{group!concrete G-set@{concrete $G$-set}} to abstract $G$-sets\index{abstract G-set@{abstract $G$-set}}\index{group!abstract G-set@{abstract $G$-set}} defined by
    \begin{equation*}
      \ev_\ast(X)\defeq (X(\ast),g\mapsto \tr_X(\varphi(g)))
    \end{equation*}
    is an equivalence. In the remainder of this exercise we will write $gx$ for $\tr_X(\varphi(g),x)$.
  \item Show that the type $X_G\defeq\prd{u:BG}X(u)$ of concrete fixed points of $X$\index{fixed point!of a concrete G-set@{of a concrete $G$-set}}\index{concrete G-set@{concrete $G$-set}!fixed point} is equivalent to the type
    \begin{equation*}
      \sm{x:X(\ast)}gx=x
    \end{equation*}
    of \define{fixed points} of the abstract $G$-set $\ev_\ast(X)$.\index{fixed point!of an abstract G-set@{of an abstract $G$-set}|textbf}\index{abstract G-set@{abstract $G$-set}!fixed point|textbf}
  \item Show that the type $X/G$ of orbits\index{concrete G-set@{concrete $G$-set}!orbit}\index{orbit} of $X$ is connected if and only if the abstract $G$-set $\ev_\ast(X)$ is transitive in the sense that\index{abstract G-set@{abstract $G$-set}!transitive abstract G-set@{transitive abstract $G$-set}|textbf}\index{transitive abstract G-set@{transitive abstract $G$-set}|textbf}
    \begin{equation*}
      \forall_{(x:X(\ast))}\issurj(g\mapsto gx)
    \end{equation*}
  \item Show that the type $X/G$ of orbits\index{concrete G-set@{concrete $G$-set}!orbit}\index{orbit} of $X$ is a set if and only if the abstract $G$-set $\ev_\ast(X)$ is free\index{abstract G-set@{abstract $G$-set}!free abstract G-set@{free abstract $G$-set}|textbf}\index{free abstract G-set@{free abstract $G$-set}|textbf} in the sense that
    \begin{equation*}
      \forall_{(x:X(\ast))}\isinj(g\mapsto gx).
    \end{equation*}
  \item Show that the type of abstract $G$-torsors\index{torsor}\index{group!torsor} is equivalent to the type of families $X:B\to\Set_\UU$ with contractible total space.\index{torsor}\index{group!torsor}
  \end{subexenum}
  \exitem (Buchholtz) Consider a group $G$ with classifying type $BG$ equipped with a group isomorphism
  \begin{equation*}
    \varphi:G\cong \loopspace{BG}.
  \end{equation*}
  Define the $G$-type $\concretesubgroup_\UU(G) : BG\to\UU$ of \define{concrete subgroups} of $G$ by\index{concrete subgroup|textbf}\index{group!concrete subgroup|textbf}\index{Concrete-Subgroup@{$\concretesubgroup_\UU(G,u)$}|textbf}
  \begin{equation*}
    \concretesubgroup_\UU(G,u)\defeq \sum_{(X:BG\to\Set_\UU)}\sum_{(x:X(u))}\isconn(X/G).
  \end{equation*}
  \begin{subexenum}
  \item Construct an equivalence\index{subgroup}\index{group!subgroup}
    \begin{equation*}
      \concretesubgroup_\UU(G,\ast)\simeq\subgroup_\UU(G).
    \end{equation*}
  \item Show that $G$ acts on $\concretesubgroup_\UU(G,\ast)$ by conjugation\index{conjugation|textbf}\index{group!conjugation|textbf}, i.e., show that for any $g:G$ we have a commuting square
    \begin{equation*}
      \begin{tikzcd}[column sep=1.6em]
        \concretesubgroup_\UU(G,\ast) \arrow[r,"g"] \arrow[d,swap,"\simeq"] & \concretesubgroup_\UU(G,\ast) \arrow[d,"\simeq"] \\
        \subgroup_\UU(G) \arrow[r,swap,"H\mapsto\{ghg^{-1}\mid h\in H\}"] & \subgroup_\UU(G)
      \end{tikzcd}
    \end{equation*}
  \item Conclude that the type of normal subgroups of a group $G$\index{normal subgroup}\index{group!normal subgroup} is equivalent to the type of \define{concrete normal subgroups}\index{concrete normal subgroup|textbf}\index{group!concrete normal subgroup|textbf}
    \begin{equation*}
      \prd{u:BG}\concretesubgroup_\UU(G,u).
    \end{equation*}
  \item Show that the type of normal subgroups of a group $G$ is also equivalent to the type
    \begin{equation*}
      \sm{BH:\concretegroup_\UU}\sm{f:BG\to_\ast BH}\isconn(f)
    \end{equation*}
  \end{subexenum}
\end{exercises}
\index{group|)}


\section{General inductive types}\label{sec:w-types}

\index{inductive type|(}
\index{W-type|(}

Most inductive types we have seen in this book have a finite number of constructors with finite arities. For example, the type $\N$ has two constructors: one constant $\zeroN$ and one unary constructor $\succN$. However, there is no objection to having an nonfinite amount of constructors, possibly with nonfinite arities. W-types are general inductive types that have a \emph{type} of constructors, whose arities are \emph{types}. W-types are therefore specified by a type $A$ of \emph{symbols} for the constructors, and a type family $B$ over $A$ specifying the arities of the constructors that the symbols represent.

An example of a W-type is the type of finitely branching rooted trees. This inductive type has a constructor with arity $X$ for each finite type $X$. In other words, a finitely branching rooted tree is obtained by attaching a finitely many finitely branching rooted trees to a root. The root itself is therefore a finitely branching tree, obtained from the $0$-ary constructor corresponding to the empty type, and if we have any finite family finitely branching rooted trees, we can combine them all into one finitely branching rooted tree by attaching them to a new root.

\subsection{The type of well-founded trees}

\begin{defn}
  Consider a type family $B$ over $A$. The \define{W-type}\index{W-type|textbf} $\W(A,B)$\index{W(A,B)@{$\W(A,B)$}|see {W-type}}\index{W(A,B)@{$W(A,B)$}|textbf} is defined as the inductive type with constructor\index{tree@{$\collect$}|textbf}\index{W-type!tree@{$\collect$}|textbf}\index{inductive type!W-type}
  \begin{equation*}
    \collect : \prd{x:A} (B(x)\to \W(A,B))\to \W(A,B).
  \end{equation*}
  The induction principle\index{W-type!induction principle|textbf}\index{induction principle!of W-types|textbf} of the W-type $\W(A,B)$ asserts that, for any type family $P$ over $\W(A,B)$, any dependent function
  \begin{equation*}
    h : \prd{x:A}\prd{\alpha:B(x)\to \W(A,B)} \Big(\prd{y:B(x)}P(\alpha(x))\Big)\to P(\collect(x,\alpha))
  \end{equation*}
  determines a dependent function\index{ind_W@{$\indW$}|textbf}\index{W-type!ind_W@{$\indW$}|textbf}
  \begin{equation*}
    \indW(h):\prd{x:\W(A,B)}P(x)
  \end{equation*}
  that satisfies the judgmental equality\index{computation rules!for W-types|textbf}\index{W-type!computation rule|textbf}
  \begin{equation*}
    \indW(h,\collect(x,\alpha))\jdeq h(x,\alpha,\lam{y}\indW(h,\alpha(y))).
  \end{equation*}
  The elements of W-types are called \define{(well-founded) trees}\index{well-founded trees|textbf}.
\end{defn}

\begin{rmk}
  Some authors write $\mathsf{sup}$ for the constructor of a W-type. The intuition that $\collect(a,\alpha)$ is a supremum of the family of elements $\alpha(y)$ indexed by $y:B(a)$ is, however, somewhat misleading, because $\collect(a,\alpha)$ does not satisfy the defining properties of a supremum.
\end{rmk}

\begin{rmk}
  When we define a dependent function
  \begin{equation*}
    f:\prd{x:\W(A,B)}P(x)
  \end{equation*}
  via the induction principle of W-types, we will often display that definition by pattern matching\index{pattern matching!for W-types}\index{W-type!pattern matching}. Such definitions are then displayed as
  \begin{equation*}
    f(\collect(x,\alpha))\defeq h(x,\alpha,\lam{y}f(\alpha(y))),
  \end{equation*}
  which contains all the information to carry out the construction via the induction principle of W-types. The advantage of definitions by pattern matching is that they directly display the defining judgmental equality the function being defined.
\end{rmk}

\begin{rmk}\label{rmk:constant-W}
  For any $x:A$, the function
  \begin{equation*}
    \collect(x):(B(x)\to\W(A,B))\to\W(A,B)
  \end{equation*}
  takes a family of elements $\alpha(y):\W(A,B)$ indexed by $y:B(x)$ and collects them into an element $\collect(x,\alpha):\W(A,B)$. Since the element $\collect(x,\alpha)$ has been constructed out of a family $\alpha(y)$ of elements of $\W(A,B)$ indexed by $y:B(x)$, we say that the type $B(x)$ is the \define{arity}\index{arity of constructor W-type|textbf}\index{W-type!arity of constructor|textbf} of $\collect(x,\alpha)$. In other words, there is a function\index{arity@{$\arity$}|textbf}\index{W-type!arity@{$\arity$}|textbf}
  \begin{equation*}
    \arity : \W(A,B)\to\UU
  \end{equation*}
  given by $\arity(\collect(x,\alpha))\defeq B(x)$. The element $x:A$ is the \define{symbol}\index{symbol of a constructor of a W-type|textbf}\index{W-type!symbol of a constructor|textbf} of the operation $\collect(x):(B(x)\to\W(A,B))\to\W(A,B)$. Note that there might be many different symbols $x,y:A$ for which the operations $\collect(x)$ and $\collect(y)$ have equivalent arities, i.e., for which $B(x)\simeq B(y)$.

  Furthermore, the \define{components}\index{component of an element pf a W-type|textbf}\index{W-type!component of an element|textbf} of $\collect(x,\alpha)$ are the elements $\alpha(y):\W(A,B)$ indexed by $y:B(x)$. In other words, we have
  \begin{align*}
    \component & : \prd{w:\W(A,B)} \arity(w) \to \W(A,B),
  \end{align*}
  given by $\component(\collect(x,\alpha))\defeq\alpha$.
  
  In the special case where $B(x)$ is empty, there is exactly one family of elements $\alpha(y):\W(A,B)$ indexed by $y:B(x)$. Therefore, it follows that any $x:A$ such that $B(x)$ is empty induces a constant in the W-type $\W(A,B)$. More precisely, if we are given a map $h:B(x)\to \emptyt$, then we can define the \define{constant}\index{constant element in a W-type|textbf}\index{W-type!constant element|textbf}
  \begin{equation*}
    c_x(h)\defeq \collect(x,\exfalso\circ h).
  \end{equation*}
  The elements of $w:\W(A,B)$ for which the type $B(\arity(w))$ is empty are called the \define{constants} of $\W(A,B)$. In other words, the predicate\index{is-constant@{$\isconstantW$}|textbf}\index{W-type!is-constant@{$\isconstantW$}|textbf}
  \begin{equation*}
    \isconstantW : \W(A,B)\to\prop_\UU
  \end{equation*}
  is defined by $\isconstantW(w)\defeq\isempty(B(\arity(w)))$.
  
  On the other hand, if each type $B(x)$ is inhabited, then there are no such constants and we will see in the following proposition that the W-type $\W(A,B)$ is empty in this case.  
\end{rmk}

\begin{prp}\label{prp:is-empty-W}
  Consider a family $B$ of types over $A$. Then the following are equivalent:
  \begin{enumerate}
  \item For each $x:A$, the type $B(x)$ is nonempty.
  \item The $W$-type $\W(A,B)$ is empty.\index{is empty!W-type}\index{W-type!is empty}
  \end{enumerate}
  In particular, if each $B(x)$ is inhabited, then $\W(A,B)$ is empty.
\end{prp}

\begin{proof}
  To prove that (i) implies (ii), assume that $\neg\neg(B(x))$ holds for each $x:A$. Our goal is to construct a function $f:\W(A,B)\to \emptyt$. By the induction principle of W-types it suffices to construct a function of type
  \begin{equation*}
    \prd{x:A}\prd{\alpha:B(x)\to\W(A,B)}\Big(\prd{y:B(x)}\emptyt\Big)\to\emptyt.
  \end{equation*}
  This type is judgmentally equal to the type
  \begin{equation*}
    \prd{x:A}\prd{\alpha:B(x)\to\W(A,B)}\neg\neg(B(x)),
  \end{equation*}
  so we obtain the desired function from the assumption that $\neg\neg(B(x))$ holds for every $x:A$.

  To prove that (ii) implies (i), suppose that $\W(A,B)$ is empty and let $x:A$. To show that $\neg\neg(B(x))$ holds, assume that $\neg(B(x))$ holds. In other words, assume a function $h:B(x)\to\emptyt$. Then we have the constant element $c_x(h):\W(A,B)$. This is impossible, since $\W(A,B)$ was assumed to be empty.
\end{proof}

\begin{eg}\label{eg:Nat-W}
  Consider the type family $P$ over $\bool$ given by
  \begin{equation*}
    P(\bfalse) \defeq \emptyt \qquad\text{and}\qquad P(\btrue) \defeq \unit.
  \end{equation*}
  We claim that the W-type $N\defeq \W(\bool,P)$\index{natural numbers!as W-type}\index{W-type!natural numbers} is equivalent to $\N$. The idea is that the constructor $\collect$ of $\W(\bool,P)$ splits into one nullary constructor with symbol $\bfalse$ and arity $P(\bfalse)\jdeq\emptyt$, and one unary constructor with symbol $\btrue$ and arity $P(\btrue)\jdeq\unit$.

  More formally, we define the zero element $z:N$ and the successor function $s:N\to N$ by
  \begin{equation*}
    z\defeq \collect(\bfalse,\exfalso) \qquad\text{and}\qquad s(x)\defeq \collect(\btrue,\const_x).
  \end{equation*}
  Thus, we obtain a function $f:\N\to N$ that satisfies $f(\zeroN)\jdeq z$ and $f(\succN(n))\jdeq s(f(n))$. It's inverse $g:N\to \N$ is defined via the induction principle of W-types by
  \begin{align*}
    g(\collect(\bfalse,\alpha)) & \defeq \zeroN \\
    g(\collect(\btrue,\alpha)) & \defeq \succN(g(\alpha(\ttt))).
  \end{align*}
  It is immediate from these definitions that $g(f(n))=n$ for all $n:\N$. It remains to construct an identification $p(x):f(g(x))=x$ for all $x:N$. Such an identification is constructed inductively. First, there is an identification
  \begin{equation*}
    p(\collect(\bfalse,\alpha)) : \collect(\bfalse,\exfalso)=\collect(\bfalse,\alpha)
  \end{equation*}
  by the fact that $\exfalso=\alpha$ for any $\alpha:\emptyt\to N$. Second, there is an identification
  \begin{equation*}
    p(\collect(\btrue,\alpha)) : \collect(\btrue,\const_{\alpha(\ttt)})=\collect(\btrue,\alpha)
  \end{equation*}
  by the fact that $\const_{\alpha(\ttt)}=\alpha$ for any map $\alpha:\unit\to N$. This completes the construction of the equivalence $\N\simeq N$.
\end{eg}

\begin{eg}\label{eg:planar-binary-tree-W}
  Consider the type family $B$ over $\bool$ given by
  \begin{equation*}
    B(\bfalse) \defeq \emptyt \qquad\text{and}\qquad B(\btrue) \defeq \bool.
  \end{equation*}
  Then the W-type $\W(\bool,B)$ is equivalent to the type of \define{oriented binary rooted trees}\index{oriented binary rooted tree|textbf}\index{tree!oriented binary rooted tree|textbf}\index{W-type!oriented binary rooted trees|textbf}, which is the inductive type with constructors\index{oriented binary rooted tree!node@{$\node$}|textbf}\index{oriented binary rooted tree![-,-]@{$[\blank,\blank]$}|textbf}\index{inductive type!oriented binary rooted trees|textbf}
  \begin{align*}
    \node & : \planarBinTree \\
    {[\blank,\blank]} & : \planarBinTree\to (\planarBinTree \to \planarBinTree).
  \end{align*}
  We leave the construction of the equivalence $\planarBinTree\simeq\W(\bool,B)$ as \cref{ex:oriented-bin-tree}. The reason we call the elements of $\planarBinTree$ oriented binary rooted trees is that in a tree of the form $[T_1,T_2]$ we can see by inspection which branch is on the left and which branch is on the right.
\end{eg}

\begin{eg}\label{eg:binary-tree-W}
  Consider the type $A\defeq \unit+\BS_2$, where $\BS_2$ is the type of $2$-element types. We define the family $B$ over $A$ by pattern matching:
  \begin{align*}
    B(\inl(x)) & \defeq \emptyt \\
    B(\inr(X)) & \defeq X.
  \end{align*}
  The type of \define{binary rooted trees}\index{binary rooted tree|textbf}\index{tree!binary rooted tree|textbf}\index{W-type!binary rooted trees|textbf} is the W-type $\W(A,B)$ for this choice of $A$ and $B$. We can also present the type of binary rooted trees as an inductive type with the following constructors:\index{binary rooted tree!node@{$\node$}|textbf}\index{binary rooted tree!bin-tree@{$\collectBinTree$}|textbf}\index{bin-tree@{$\collectBinTree$}|textbf}\index{Bin-Tree@{$\BinTree$}|textbf}\index{inductive type!binary rooted trees|textbf}
  \begin{align*}
    \node & : \BinTree \\
    \collectBinTree & : \prd{X:\BS_2} \BinTree^X\to \BinTree.
  \end{align*}
  There is an important qualitative difference between the type of oriented binary rooted trees and the type of binary rooted trees. Given two distinct oriented binary rooted trees $T_1$ and $T_2$, the two oriented binary rooted trees $[T_1,T_2]$ and $[T_2,T_1]$ will also be distinct. On the other hand, given two binary rooted trees $T_1$ and $T_2$, the binary rooted trees
  \begin{align*}
    & \collectBinTree (\bool,\indbool(T_1,T_2)) \\
    & \collectBinTree(\bool,\indbool(T_2,T_1))
  \end{align*}
  can always be identified. In the terminology of \cref{ex:commutative-binary-operations}, the constructor $\collectBinTree$ of $\BinTree$ is equivalently described as a commutative binary operation on $\BinTree$.
\end{eg}

\begin{eg}\label{eg:finitely-branching-tree-W}
  The W-type $\W(\N,\Fin{})$ is the type of \define{oriented finitely branching rooted trees}\index{tree!oriented finitely branching rooted tree|textbf}\index{oriented finitely branching rooted tree|textbf}\index{W-type!oriented finitely branching rooted trees|textbf}. On the other hand, we define the type of \define{(unoriented) finitely branching rooted trees}\index{tree!finitely branching rooted tree|textbf}\index{finitely branching rooted tree|textbf}\index{W-type!finitely branching rooted trees|textbf} to be the W-type $\W(\F,\mathcal{T})$. The qualitive difference between the types of oriented and unoriented finitely branching rooted trees is similar to the qualitative difference between types of oriented and unoriented binary rooted trees. In the type of oriented finitely branching rooted trees, we record the ordering of the branches while in the type of unoriented finitely branching rooted trees there are identifications between trees that have the same branches up to permutation.
\end{eg}

\subsection{Observational equality of W-types}

\index{observational equality!on W-types|(}
\index{W-type!observational equality|(}

Each element $x:\W(A,B)$ has symbol $\prearity(x):A$ and a family of components $\component(x):B(\prearity(x))\to\W(A,B)$. Therefore, we have a map
\begin{equation*}
  \eta : \W(A,B)\to \sm{x:A}(B(x)\to\W(A,B))
\end{equation*}
given by $\eta(x)\defeq(\prearity(x),\component(x))$.

\begin{prp}\label{prp:algebra-W}
  The map $\eta:\W(A,B)\to\sm{x:A}(B(x)\to\W(A,B))$ is an equivalence.
\end{prp}

\begin{proof}
  We define
  \begin{equation*}
    \varepsilon : \Big(\sm{x:A}(B(x)\to\W(A,B))\Big)\to\W(A,B)
  \end{equation*}
  by $\varepsilon(x,\alpha)\defeq\collect(x,\alpha)$. The fact that $\varepsilon$ is an inverse of $\eta$ follows easily.
\end{proof}

The fact that we have an equivalence
\begin{equation*}
  \W(A,B)\simeq\sm{x:A}(B(x)\to\W(A,B)),
\end{equation*}
suggests a way to characterize the identity type of $\W(A,B)$. Indeed, any equivalence is an embedding, and therefore we also have
\begin{equation*}
  (x=y)\simeq (\eta(x)=\eta(y)).
\end{equation*}
The latter is an identity type in a $\Sigma$-type, which can be characterized as a $\Sigma$-type of identity types. We therefore define the following observational equality relation on $\W(A,B)$.

\begin{defn}
  Suppose $A$ and each $B(x)$ are in $\UU$. We define a binary relation\index{Eq W@{$\EqW$}|textbf}\index{W-type!Eq W@{$\EqW$}|textbf}
  \begin{equation*}
    \EqW : \W(A,B)\to \W(A,B)\to \UU
  \end{equation*}
  recursively by
  \begin{equation*}
    \EqW(\collect(x,\alpha),\collect(y,\beta)) \defeq \sm{p:x=y}\prd{z:B(x)}\,\alpha(z)=\beta(\tr_B(p,z))
  \end{equation*}
\end{defn}

\begin{thm}\label{thm:EqW}
  The observational equality relation $\EqW$ on $\W(A,B)$ is reflexive, and the canonical map\index{characterization of identity type!of W-types}\index{W-type!characterization of identity type}\index{identity type!of W(A,B)@{of $\W(A,B)$}}
  \begin{equation*}
    (x=y)\to \EqW(x,y)
  \end{equation*}
  is an equivalence for each $x,y:\W(A,B)$. 
\end{thm}

\begin{proof}
  The element $\reflEqW(x):\EqW(x,x)$ is defined recursively as
  \begin{equation*}
    \reflEqW(\collect(x,\alpha))\defeq (\refl{x},\reflhtpy_\alpha).
  \end{equation*}
  This proof of reflexivity induces the canonical map $(x=y)\to\EqW(x,y)$. To show that it is an equivalence for each $x,y:\W(A,B)$, we apply the fundamental theorem of identity types, by which it suffices to show that the type
  \begin{equation*}
    \sm{y:\W(A,B)}\EqW(x,y)
  \end{equation*}
  is contractible for each $x:\W(A,B)$. The center of contraction is the pair $(x,\reflEqW(x))$. For the contraction, we have to construct a function
  \begin{equation*}
    h:\prd{y:\W(A,B)}\prd{p:\EqW(x,y)}\,(x,\reflEqW(x))=(y,p).
  \end{equation*}
  By the induction principle of W-types, it suffices to define
  \begin{equation*}
    h(\collect(y,\beta),(p,H))\defeq (x,(\refl{},\reflhtpy))=(y,(p,H)).
  \end{equation*}
  Here we proceed by identification elimination on $p:x=y$, followed by homotopy induction on the homotopy $H:\alpha\htpy \beta$. Thus, it suffices to construct an identification
  \begin{equation*}
    (x,(\refl{},\reflhtpy))=(x,(\refl{},\reflhtpy)),
  \end{equation*}
  which we have by reflexivity.
\end{proof}

\begin{thm}
  Consider a type family $B$ over a type $A$, and let $k:\T$ be a truncation level. If $A$ is a $(k+1)$-type, then so is $\W(A,B)$.\index{W-type!is truncated}\index{is truncated!W-type}
\end{thm}

\begin{proof}
  Suppose that $A$ is a $(k+1)$-type. In order to show that $\W(A,B)$ is a $(k+1)$-type, we have to show that its identity types are $k$-types. The proof is by induction on $x,y:\W(A,B)$. For $x\jdeq\collect(a,\alpha)$ and $y\jdeq\collect(b,\beta)$, we have the equivalence
  \begin{equation*}
    (\collect(a,\alpha)=\collect(b,\beta))\simeq\sm{p:a=b}\prd{z:B(a)}\,\alpha(z)=\beta(\tr_B(p,z))
  \end{equation*}
  Note that the type $a=b$ is a $k$-type by the assumption that $A$ is a $(k+1)$-type. Furthermore, the type $\alpha(z)=\beta(\tr_B(p,z))$ is a $k$-type by the induction hypothesis. Therefore it follows that the type on the right-hand side of the displayed equivalence is a $k$-type, and this completes the proof.
\end{proof}
\index{observational equality!on W-types|)}
\index{W-type!observational equality|)}

\subsection{Functoriality of W-types}
\index{functorial action!of W-types|(}
\index{W-type!functorial action|(}

\begin{defn}
  Consider a type family $B$ over $A$, and a type family $B'$ over $A'$. Furthermore, consider a map $f:A'\to A$ and a family of equivalences
  \begin{equation*}
    e_x:B'(x)\simeq B(f(x))
  \end{equation*}
  indexed by $x:A'$. Then we define the map $\W(f,e):\W(A',B')\to\W(A,B)$\index{W(f,e)@{$\W(f,e)$}|textbf}\index{W(f,e)@{$\W(f,e)$}|see {W-type, functorial action}}\index{W-type!functorial action|textbf}\index{functorial action!of W-types|textbf} of W-types inductively by
  \begin{equation*}
    \W(f,e)(\collect(x,\alpha))\defeq\collect(f(x),\W(f,g)\circ \alpha\circ e_x^{-1}).
  \end{equation*}
\end{defn}

\begin{lem}\label{lem:fib-W}
  For any morphism $\W(f,e):\W(A',B')\to\W(A,B)$ of W-types and any $\collect(x,\alpha):\W(A,B)$, there is an equivalence\index{fiber!of W(f,e)@{of $\W(f,e)$}}\index{W(f,e)@{$\W(f,e)$}!fiber}
  \begin{equation*}
    \fib{\W(f,e)}{\collect(x,\alpha)} \simeq \fib{f}{x}\times\prd{b:B(x)}\fib{\W(f,e)}{\alpha(b)}.
  \end{equation*}
\end{lem}

\begin{proof}
  First, note that by the characterization in \cref{thm:EqW} of the identity type of $\W(A,B)$, there is an equivalence between the fiber $\fib{\W(f,e)}{\collect(x,\alpha)}$ and the type
  \begin{align*}
    & \sm{x':A'}\sm{\alpha':B'(x')\to\W(A',B')}\sm{p:f(x')=x} \\*
    & \phantom{\sm{x':A'}}\prd{b:B(f(x'))}\W(f,e)(\alpha'(e_{x'}^{-1}(b)))=\alpha(\tr_{B}(p,b)). \\
    \intertext{By rearranging the $\Sigma$-type, we see that this type is equivalent to the type}
    & \sm{(x',p):\fib{f}{x}}\sm{\alpha':B'(x')\to\W(A',B')} \\*
    & \phantom{\sm{x':A'}}\prd{b:B(f(x'))}\W(f,e)(\alpha'(e_{x'}^{-1}(b)))=\alpha(\tr_{B}(p,b)).
  \end{align*}
  Therefore, it suffices to show for each $(x',p):\fib{f}{x}$, that the type
  \begin{equation*}
    \sm{\alpha':B'(x')\to\W(A',B')}\prd{b:B(f(x'))}\W(f,e)(\alpha'(e_{x'}^{-1}(b)))=\alpha(\tr_{B}(p,b))
  \end{equation*}
  is equivalent to the type $\prd{b:B(x)}\fib{\W(f,e)}{\alpha(b)}$. Since we have an identification $p:f(x')=x$ and an equivalence $e_{x'}:B'(x')\simeq B(f(x'))$, it follows that the type above is equivalent to the type
  \begin{equation*}
    \sm{\alpha':B(x)\to\W(A',B')}\prd{b:B(x)}\W(f,e)(\alpha'(b))=\alpha(b).
  \end{equation*}
  By distributivity of $\Pi$ over $\Sigma$, i.e., by \cref{thm:choice}, this type is equivalent to the type
  \begin{equation*}
    \prd{b:B(x)}\sm{w:\W(A',B')}\W(f,e)(w)=\alpha(b),
  \end{equation*}
  completing the proof.
\end{proof}

\begin{thm}
  Consider a morphism $\W(f,e):\W(A,B)\to\W(A',B')$ of W-types. If the map $f:A\to A'$ is $k$-truncated, then so is the map $\W(f,e)$. In particular, if $f$ is an equivalence or an embedding, then so is $\W(f,e)$.\index{W(f,e)@{$\W(f,e)$}!is a truncated map}\index{W(f,e)@{$\W(f,e)$}!is an embedding}\index{W(f,e)@{$\W(f,e)$}!is an equivalence}\index{is an equivalence!W(f,e)@{$\W(f,e)$}}\index{is an embedding!W(f,e)@{$\W(f,e)$}}
\end{thm}

\begin{proof}
  Suppose that the map $f$ is $k$-truncated. We will prove recursively that the fibers of the morphism $\W(f,e)$ on W-types is $k$-truncated. We saw in \cref{lem:fib-W} that there is an equivalence
  \begin{equation*}
    \fib{\W(f,e)}{\collect(x,\alpha)}\simeq \fib{f}{x}\times\prd{b:B(x)}\fib{\W(f,e)}{\alpha(b)}.
  \end{equation*}
  The type $\fib{f}{x}$ is $k$-truncated by assumption, and each of the types
  \begin{equation*}
    \fib{\W(f,e)}{\alpha(b)}
  \end{equation*}
  is $k$-truncated by the inductive hypothesis, so the claim follows.
\end{proof}
\index{functorial action!of W-types|)}
\index{W-type!functorial action|)}

\subsection{The elementhood relation on W-types}
\index{elementhood relation on W-types|(}
\index{W-type!elementhood relation|(}

The elements of a W-type $\W(A,B)$ are constructed out of families of elements of $\W(A,B)$ indexed by a type $B(x)$ for some $x:A$. More precisely, for each $\collect(x,\alpha):\W(A,B)$ we have a family of elements
\begin{equation*}
  \alpha(y):\W(A,B)
\end{equation*}
indexed by $y:B(x)$. Thus, we could say that $\alpha(y)$ is in $\collect(x,\alpha)$, for each $y:B(x)$. More abstractly, we can define an elementhood relation on $\W(A,B)$.

\begin{defn}
  Given a W-type $\W(A,B)$ and a universe $\UU$ containing both $A$ and each type in the family $B$, we define a type-valued relation\index{e@{$\in$}|see {elementhood relation on W-types}|textbf}\index{e@{$\in$}|textbf}\index{W-type!e@{$\in$}|textbf}
  \begin{equation*}
    {\in}:\W(A,B)\to\W(A,B)\to \UU
  \end{equation*}
  by $(x\in \collect(a,\alpha))\defeq \sm{y:B(a)}\alpha(y)=x$. 
\end{defn}

Using the elementhood relation on $\W(A,B)$, we can reformulate the induction principle to, perhaps, a more recognizable form:

\begin{thm}
  For any family $P$ of types over $\W(A,B)$, there is a function\index{induction principle!of W-types}\index{W-type!induction principle}
  \begin{equation*}
    i : \Big(\prd{x:\W(A,B)}\Big(\prd{y:\W(A,B)}(y\in x)\to P(y)\Big)\to P(x)\Big)\to \Big(\prd{x:X}P(x)\Big)
  \end{equation*}
  that comes equipped with an identification
  \begin{equation*}
    i(h,x)=h(x,\lam{y}\lam{e}i(h,y))
  \end{equation*}
  for every $h:\prd{x:\W(A,B)}\Big(\prd{y:\W(A,B)}(y\in x)\to P(y)\Big)\to P(x)$, and every $x:\W(A,B)$.
\end{thm}

\begin{proof}
  For any type family $P$ over $\W(A,B)$, we first define a new type family $\square P$ over $\W(A,B)$ given by
  \begin{equation*}
    \square P(x):=\prd{y:\W(A,B)}(y\in x)\to P(y).
  \end{equation*}
  The family $\square P(x)$ comes equipped with a map
  \begin{equation*}
    \eta : \Big(\prd{x:\W(A,B)}P(x)\Big)\to \Big(\prd{x:\W(A,B)}\square P(x)\Big)
  \end{equation*}
  given by $\eta(f,x,y,e)\defeq f(y)$. Conversely, there is a map
  \begin{equation*}
    \varepsilon(h) : \Big(\prd{y:\W(A,B)}\square P(y)\Big) \to \Big(\prd{x:\W(A,B)}P(x)\Big)
  \end{equation*}
  for every $h:\prd{y:\W(A,B)}\square P (y)\to P(y)$, given by $\varepsilon(h,g,x)\defeq h(x,g(x))$. Note that the induction principle can now be stated as
  \begin{equation*}
    i : \Big(\prd{y:\W(A,B)}\square P (y)\to P(y)\Big)\to\Big(\prd{x:\W(A,B)}P(x)\Big),
  \end{equation*}
  and the computation rule states that
  \begin{equation*}
    i(h,x)=h(x,\eta(i(h),x)).
  \end{equation*}
  Before we prove the induction principle, we prove the intermediate claim that there is a function
  \begin{equation*}
    i' : \Big(\prd{y:\W(A,B)}\square P (y)\to P(y)\Big) \to \Big(\prd{x:\W(A,B)}\square P(x)\Big)
  \end{equation*}
  equipped with an identification
  \begin{equation*}
    j'(h,x,y,e) : i'(h,x,y,e) = h(y,i'(h,y))
  \end{equation*}
  for every $h:\prd{y:\W(A,B)}\square P(y)\to P(y)$ and every $x,y:\W(A,B)$ equipped with $e:y\in x$. Both $i'$ and $j'$ are defined by pattern matching:
  \begin{align*}
    i'(h,\collect(a,f),f(b),(b,\refl{})) & := h(f(b),i'(h,f(b))) \\
    j'(h,\collect(a,f),f(b),(b,\refl{})) & := \refl{}.
  \end{align*}
  Now we define $i(h):=\varepsilon(h,i'(h))$. Note that we have the judgmental equalities
  \begin{align*}
    i(h,x) & \jdeq \varepsilon(h,i'(h),x) \\
           & \jdeq h(x,i'(h,x)), \\
    \intertext{and}
    h(x,\lam{y}\lam{e}i(h,y))
           & \jdeq h(x,\lam{y}\lam{e}\varepsilon(h,i'(h),y)) \\
           & \jdeq h(x,\lam{y}\lam{e}h(y,i'(h,y))).
  \end{align*}
  The computation rule is therefore satisfied by the identification
  \begin{equation*}
    \begin{tikzcd}[column sep=12.5em]
      h(x,i'(h,x)) \arrow[r,equals,"\ap{h(x)}{\eqhtpy(\lam{y}\eqhtpy(j'(h,x,y)))}"] & h(x,\lam{y}\lam{e} h(y,i'(h,y))).
    \end{tikzcd}
    \qedhere
  \end{equation*}
\end{proof}

\subsection{Extensional W-types}\label{sec:extensional-W-types}

It is tempting to think that an element $w:\W(A,B)$ is completely determined by the elements $z:\W(A,B)$ equipped with a proof $z\in w$. However, this may not be the case. For instance, a W-type $\W(A,B)$ might have \emph{two} unary constructors, e.g., when $A\defeq \unit+\bool$ and the family $B$ over $A$ is given by
\begin{align*}
  B(\inl(x)) & \defeq \emptyt \\
  B(\inr(y)) & \defeq \unit.
\end{align*}
If we write $f$ and $g$ for the two unary constructors of $\W(A,B)$, then we see that for any element $w:\W(A,B)$, the elements
\begin{equation*}
  u\defeq\collect(\inr(\bfalse),\const_w)\qquad\text{and}\qquad v\defeq\collect(\inr(\btrue),\const_w)
\end{equation*}
both only contain the element $w$. However, the elements $u$ and $v$ are distinct in $\W(A,B)$.

Something similar happens in the type of oriented binary rooted trees. Given two binary rooted trees $S$ and $T$, there are two ways to combine $S$ and $T$ into a new binary tree: we have $[S,T]$ and $[T,S]$. Both contain precisely the elements $S$ and $T$, but they are distinct. Nevertheless, there are many important W-types in which the elements $w$ are uniquely determined by the elements $z\in w$. Such W-types are called extensional.

\begin{defn}
  We say that a W-type $\W(A,B)$ is \define{extensional}\index{extensional W-type|textbf}\index{W-type!extensionality|textbf}\index{extensionality principle!for W-types|textbf} if the canonical map
  \begin{equation*}
    (x=y)\to\prd{z:\W(A,B)}(z\in x)\simeq (z\in y)
  \end{equation*}
  is an equivalence.
\end{defn}

In the following theorem we give a precise characterization of the inhabited extensional W-types. 

\begin{thm}\label{thm:extensional-W}
  Consider an inhabited W-type $\W(A,B)$. Then the following are equivalent:
  \begin{enumerate}
  \item The W-type $\W(A,B)$ is extensional.
  \item The family $B$ is \define{univalent}\index{type family!univalent type family|textbf}\index{univalent type family|textbf} in the sense that the map
  \begin{equation*}
    \tr_B:(x=y)\to (B(x)\simeq B(y))
  \end{equation*}
  is an equivalence, for every $x,y:A$.
  \end{enumerate}
\end{thm}

\begin{rmk}
  Note that if the W-type $\W(A,B)$ is empty, then it is vacuously extensional. However, we saw in \cref{prp:is-empty-W} that any family $B$ of inhabited types over $A$ gives rise to an empty W-type $\W(A,B)$, so there is no hope of showing that $B$ is a univalent family if $\W(A,B)$ is empty.

  We also note that a type family $B$ over $A$ is univalent if and only if the map $B:A\to \UU$ is an embedding. In other words, the claim in \cref{thm:extensional-W} is that an inhabited W-type $\W(A,B)$ is extensional if and only if $B$ is the canonical type family over a subuniverse $A$ of $\UU$.
\end{rmk}

\begin{proof}
  We will first show that (ii) is equivalent to the following property:
  \begin{enumerate}
  \item[(ii')] The map
    \begin{equation*}
      \tr_B : (\prearity(x)=y)\to (B(\prearity(x))\simeq B(y))
    \end{equation*}
    is an equivalence for every $x:\W(A,B)$ and every $y:A$.
  \end{enumerate}
  Clearly, (ii) implies (ii'). For the converse we use the assumption that $\W(A,B)$ is inhabited. Since the property in (ii) is a proposition, we may assume an element $w:\W(A,B)$. Using $w$, we obtain for every $x:A$ the element
  \begin{equation*}
    \collect(x,\const_w):\W(A,B)
  \end{equation*}
  The symbol of $\collect(x,\const_w)$ is $x$, and therefore the hypothesis that (ii') holds implies that the map $(x=y)\to (B(x)\simeq B(y))$ is an equivalence. This concludes the proof that (ii) is equivalent to (ii'). It remains to show that (i) is equivalent to (ii').

  Let $x:\W(A,B)$. By the fundamental theorem of identity types, the W-type $\W(A,B)$ is extensional if and only if the total space
  \begin{equation*}
    \sm{y:\W(A,B)}\prd{z:\W(A,B)}(z\in x)\simeq (z\in y)
  \end{equation*}
  is contractible, for any $x:\W(A,B)$. When $x$ is of the form $\collect(a,\alpha)$, the type $z\in x$ is just the fiber $\fib{\alpha}{z}$. Using this observation, we see that the above type is equivalent to the type
  \begin{equation*}
    \sm{b:A}\sm{\beta:B(b)\to\W(A,B)}\prd{z:\W(A,B)}\fib{\alpha}{z}\simeq\fib{\beta}{z}.\tag{\textasteriskcentered}
  \end{equation*}
  By \cref{ex:fam-equiv} it follows that this type is equivalent to the type
  \begin{equation*}
    \sm{y:A}\sm{\beta:B(y)\to\W(A,B)}\sm{e:B(x)\simeq B(y)}\alpha\htpy e\circ \beta.
  \end{equation*}
  Note that the type $\sm{\beta:B(y)\to\W(A,B)}\alpha\htpy e\circ\beta$ is contractible for any equivalence $e:B(x)\simeq B(y)$. Therefore, it follows that the above type is contractible if and  only if the type
  \begin{equation*}
    \sm{y:A}B(x)\simeq B(y)
  \end{equation*}
  is contractible, which is the case if and only if the map $(x=y)\to(B(x)\simeq B(y))$ is an equivalence for all $y:A$.
\end{proof}

\begin{eg}
  The type $N$ of \cref{eg:Nat-W}\index{natural numbers!is an extensional W-type}, the type of binary rooted trees \cref{eg:binary-tree-W}\index{binary rooted tree!is an extensional W-type}, and the type of finitely branching rooted trees \cref{eg:finitely-branching-tree-W}\index{finitely branching rooted tree!is an extensional W-type} are all examples extensional W-types. On the other hand, the type of oriented binary rooted trees of \cref{eg:planar-binary-tree-W} and the type of oriented finitely branching rooted trees of \cref{eg:finitely-branching-tree-W} are not extensional.
\end{eg}
\index{elementhood relation on W-types|)}
\index{W-type!elementhood relation|)}

\subsection{Russell's paradox in type theory}\label{subsec:russell}

\index{Russell's paradox|(}
\index{multiset|(}
Russell's paradox tells us that there cannot be a set of all sets. If there were such a set $S$, then we could form the set
\begin{equation*}
  R\defeq \{x\in S\mid x\notin x\},
\end{equation*}
for which we have $R\in R\leftrightarrow R\notin R$, a contradiction. To reproduce Russell's paradox in type theory, we first recall a crucial difference between the type theoritic judgment $a:A$ and the set theoretic proposition $x\in y$. Although the judgment $a:A$ plays a similar role in type theory as the elementhood relation, types and their elements are fundamentally different entities, whereas in Zermelo-Fraenkel set theory there are only sets, and the proposition $x\in y$ can be formed for any two sets $x$ and $y$. In type theory, there is no relation on the universe that is similar to the elementhood relation.

However, we have seen in \cref{sec:extensional-W-types} that it is possible to define an elementhood relation on arbitrary W-types. We will use this elementhood relation on the W-type $\W(\UU,\Ty)$ to derive a paradox analogous to Russell's paradox, and we will see that $\UU$ cannot be equivalent to a type in $\UU$.

The type $\W(\UU,\Ty)$ possesses a lot of further structure. In fact, it can be used to encode constructive set theory in type theory. There is, however, one significant difference with ordinary set theory: the elementhood relation is type-valued. In other words, there may be many ways in which $x\in y$ holds. The type $\W(\UU,\Ty)$ is therefore also called the type of \define{multisets}. It was first studied by Aczel in \cite{AczelCZF}, with refinements in \cite{AczelGambinoCZF}, and in the setting of univalent mathematics it has been studied extensively by Gylterud in \cite{GylterudMultisets}.

\begin{defn}
  Consider a $\UU$ with universal type family $\Ty$. We define the type\index{M_U@{$\multiset{\UU}$}|textbf}\index{M_U@{$\multiset{\UU}$}|see {multiset}}
  \begin{equation*}
    \multiset{\UU} \defeq \W(\UU,\Ty),
  \end{equation*}
  and the elements of $\multiset{\UU}$ are called \define{multisets in $\UU$}\index{multiset|textbf}. We will write\index{{f(x)"|"x:A}@{$\{f(x)\mid x:A\}$}|textbf}\index{multiset!{f(x)"|"x:A}@{$\{f(x)\mid x:A\}$}|textbf}
  \begin{equation*}
    \{f(x)\mid x:A\}
  \end{equation*}
  for the multiset in $\UU$ of the form $\collect(A,f)$. More generally, given an element $t(x_0,\ldots,x_n):\multiset{\UU}$ in context $x_0:A_0,\ldots,x_n:A_n(x_0,\ldots,x_{n-1})$, where each $A_i$ is in $\UU$, we will write\index{{t(x0,...,xn)"|"x0:A0,...,xn:An}@{$\{t(x_0,\ldots,x_n)\mid x_0:A_0,\ldots,x_n:A_n\}$}|textbf}\index{multiset!{t(x0,...,xn)"|"x0:A0,...,xn:An}@{$\{t(x_0,\ldots,x_n)\mid x_0:A_0,\ldots,x_n:A_n\}$}|textbf}
  \begin{equation*}
    \{t(x_0,\ldots,x_n)\mid x_0:A_0,\ldots,x_n:A_n(x_0,\ldots,x_{n-1})\}
  \end{equation*}
  for the multiset in $\UU$ of the form
  \begin{equation*}
    \collect\left(\sm{x_0:A_0}\cdots A_n(x_0,\ldots,x_{n-1}),\lam{(x_0,\ldots,x_n)}t(x_0,\ldots,x_n)\right).
  \end{equation*}

  Given a multiset $X\jdeq \{f(x)\mid x:A\}$ in $\UU$, the \define{cardinality}\index{cardinality!of a multiset|textbf}\index{multiset!cardinality|textbf} of $X$ is the type $A$, and the \define{elements}\index{element of a multiset|textbf}\index{multiset!element|textbf} of $X$ are the multisets $f(x)$ in $\UU$, for each $x:A$.
\end{defn}

In the notation of multisets, the elementhood relation ${\in}:\multiset{\UU}\to\multiset{\UU} \to\UU^+$ is defined by
\begin{equation*}
  (X\in \{g(y)\mid y:B\}) \jdeq \sm{y:B} g(y)=X.
\end{equation*}
In other words, a multiset $X$ is in a multiset of the form $\{g(y)\mid y:B\}$ if and only if $X$ comes equipped with an element $y:B$ and an identification $g(y)=X$. The W-type of multisets is extensional by \cref{thm:extensional-W} and the univalence axiom.

Recall from \cref{defn:small-types}\index{small type}\index{U-small type@{$\UU$-small type}} that for a universe $\UU$, we say that a type $A$ is (essentially) $\UU$-small if $A$ comes equipped with an element of type\index{is-small@{$\issmall_\UU(A)$}}
\begin{equation*}
  \issmall_{\UU}(A)\defeq \sm{X:\UU}A\simeq X. 
\end{equation*}
Our goal in this section is to show, via Russell's paradox, that the universe $\UU$ is not $\UU$-small, i.e., that there cannot be a type $U:\UU$ equipped with an equivalence $\UU\simeq U$. We will use a similar condition of smallness for multisets.

\begin{defn}
  Let $\UU$ and $\VV$ be universes. We say that a multiset $\{f(x)\mid x:A\}$ in $\VV$ \define{is $\UU$-small}\index{small multiset|textbf}\index{U-small multiset@{$\UU$-small multiset}|textbf}\index{multiset!small multiset|textbf} if the type $A$ is $\UU$-small and if each mulitset $f(x)$ in $\VV$ is $\UU$-small. In other words, the type family\index{is-small-M@{$\issmallmultiset{\UU}$}|textbf}\index{multiset!is-small-M@{$\issmallmultiset{\UU}$}|textbf}
  \begin{equation*}
    \issmallmultiset{\UU} : \multiset{\VV}\to \VV\sqcup\UU^+ 
  \end{equation*}
  is defined recursively by
  \begin{equation*}
    \issmallmultiset{\UU}(\{f(x)\mid x:A\}) \defeq \issmall_{\UU}(A)\times \prd{x:A}\issmallmultiset{\UU}(f(x)).
  \end{equation*}
\end{defn}

We will need quite a few properties of smallness before we can reproduce Russell's paradox. We begin with a simple lemma.

\begin{lem}\label{lem:is-small-comprehension-multiset}
  Consider a $\UU$-small multiset $\{f(x)\mid x:A\}$ in $\VV$, and let $B$ be a family of $\UU$-small types over $A$. Then the multiset
  \begin{equation*}
    \{f(x)\mid x:A, y:B(x)\}
  \end{equation*}
  is again $\UU$-small.
\end{lem}

\begin{proof}
  If the multiset $\{f(x)\mid x:A\}$ is $\UU$-small, then the type $A$ is $\UU$-small. By the assumption that $B$ is a family of $\UU$-small types together with the fact that $\UU$-small types are closed under formation of $\Sigma$-types, it follows that the type
  \begin{equation*}
    \sm{x:A}B(x)
  \end{equation*}
  is $\UU$-small. Furthermore, since each $f(x)$ is $\UU$-small, we conclude that the multiset $\{f(x)\mid x:A,y:B(x)\}$ is $\UU$-small. 
\end{proof}

The main purpose of the following lemma is to know that the elementhood relation takes values in the $\UU$-small types, when it is applied to $\UU$-small multisets. We will use the univalence axiom to prove this fact. 

\begin{prp}\label{prp:is-small-elementhood-multiset}
  Consider two univalent universes $\UU$ and $\VV$, and let $X$ and $Y$ be $\UU$-small multisets in $\VV$. We make two claims:
  \begin{enumerate}
  \item The type $X=Y$ is $\UU$-small.
  \item The type $X\in Y$ is $\UU$-small.
  \end{enumerate}
\end{prp}

\begin{proof}
  For the first claim, let $X\jdeq\{f(x)\mid x : A\}$ and let $Y\jdeq\{g(y)\mid y:B\}$. The proof is by induction. Via \cref{thm:EqW} it follows that the type $X=Y$ is equivalent to the type
  \begin{equation*}
    \sm{p:A=B}\prd{x:A}f(x)=g(\equiveq(p)).
  \end{equation*}
  The type $A=B$ is $\UU$-small because it is equivalent to the type $A\simeq B$, which is $\UU$-small. Therefore it suffices to show that the type
  \begin{equation*}
    \prd{x:A}f(x)=g(\equiveq(p))
  \end{equation*}
  is $\UU$-small, for every $p:A=B$. Here we proceed by identification elimination, and the type $\prd{x:A}f(x)=g(x)$ is a product of $\UU$-small types by the induction hypothesis. This concludes the proof of the first claim.

  For the second claim, let $Y\jdeq\{g(y)\mid y:B\}$. Then the type
  \begin{equation*}
    \sm{y:B}g(y)=X
  \end{equation*}
  is a dependent sum of $\UU$-small types, indexed by an $\UU$-small type, which is again $\UU$-small.
\end{proof}

The condition that a multiset $\{f(x)\mid x:A\}$ in $\VV$ is $\UU$-small suggests that there is an `equivalent' multiset in $\UU$. 

\begin{defn}\label{defn:inclusion-small-multisets}
  Given two universes $\UU$ and $\VV$, we define an inclusion function
  \begin{equation*}
    i : \Big(\sm{X:\multiset{\VV}}\issmallmultiset{\UU}(X)\Big)\to\multiset{\UU},
  \end{equation*}
  of the $\UU$-small multisets in $\VV$ into the multisets in $\UU$, inductively by
  \begin{equation*}
    i(\{f(x)\mid x:A\})\defeq \{i(f(e^{-1}(y))) \mid y:B\}.
  \end{equation*}
  for any multiset $\{f(x)\mid x:A\}$ of which the type $A$ is equipped with an equivalence $e:A\simeq B$ for some $B$ in $\UU$, and such that the multiset $f(x)$ in $\VV$ is $\UU$-small for each $x:A$.
\end{defn}

\begin{prp}\label{prp:is-embedding-inclusion-small-multisets}
  The inclusion function $i$ of $\UU$-small multisets in $\VV$ into the multisets in $\UU$ satisfies the following properties
  \begin{enumerate}
  \item For each $\UU$-small multiset $X$ in $\VV$, the multiset $i(X)$ in $\UU$ is $\VV$-small.
  \item The induced map
    \begin{equation*}
      \Big(\sm{X:\multiset{\VV}}\issmallmultiset{\UU}(X)\Big)\to\Big(\sm{Y:\multiset{\UU}}\issmallmultiset{\VV}(Y)\Big)
    \end{equation*}
    is an equivalence.
  \end{enumerate}
  Consequently, the inclusion function $i$ is an embedding.
\end{prp}

\begin{proof}
  To see that $i(\{f(x)\mid x:A\})$ is $\VV$-small for each $\UU$-small multiset $\{f(x)\mid x:A\}$ in $\VV$, note that the assumption that $\{f(x)\mid x:A\}$ is $\UU$-small gives us an equivalence $e:A\simeq B$ and an element $H(x):\issmallmultiset{\UU}(f(x))$ for each $x:A$. The type $B$ is the indexing type of $i(\{f(x)\mid x:A\})$, and $B$ is $\VV$-small because it is equivalent to the type $A$ in $\VV$. Furthermore, each multiset $i(f(e^{-1}(y)))$ is $\VV$-small by the inductive hypothesis. This completes the proof of the first claim.

  We therefore have inclusion functions
  \begin{equation*}
    \begin{tikzcd}
      \Big(\sm{X:\multiset{\VV}}\issmallmultiset{\UU}(X)\Big) \arrow[r,yshift=-.7ex,swap,"i"] &
      \Big(\sm{Y:\multiset{\UU}}\issmallmultiset{\VV}(Y)\Big) \arrow[l,yshift=.7ex,swap,"i"]
    \end{tikzcd}
  \end{equation*}
  To see that the maps $i$ and $i$ are mutual inverses, it suffices to show that $i(i(X))=X$. This follows by induction from the following calculation, where we assume an equivalence $e:A\simeq B$ into a $B$ in $\UU$.
  \begin{align*}
    i(i(\{f(x)\mid x :A\})) & \jdeq i(\{i(f(e^{-1}(y)))\mid y:B\}) \\
                            & \jdeq \{i(i(f(e^{-1}(e(x))))) \mid x:A\} \\
                            & = \{i(i(f(x)))\mid x:A\} \\
                            & = \{f(x)\mid x:A\}.
  \end{align*}
  
  For the last claim, note that we have factored $i$ as an equivalence followed by an embedding
  \begin{equation*}
    \begin{tikzcd}[column sep=2em]
      \Big(\sm{X:\multiset{\VV}}\issmallmultiset{\UU}(X)\Big) \arrow[r] &
      \Big(\sm{Y:\multiset{\UU}}\issmallmultiset{\VV}(Y)\Big) \arrow[r] &
      \multiset{\VV},
    \end{tikzcd}
  \end{equation*}
  and therefore $i$ is an embedding.
\end{proof}

Furthermore, the embedding $i$ induces equivalences on the elementhood relation on multisets.

\begin{prp}\label{prp:elementhood-small-multisets}
  Consider a multiset $X$ in $\UU$ and a multiset $Y$ in $\VV$. Furthermore, suppose that $X$ is $\VV$-small and that $Y$ is $\UU$-small. Then we have
  \begin{equation*}
    (i(X)\in Y)\simeq (X\in i(Y)).
  \end{equation*}
\end{prp}

\begin{proof}
  Let $X\jdeq\{f(x)\mid x:A\}$ and $Y\jdeq\{g(y)\mid y:B\}$. By the assumption that $Y$ is $\UU$-small we have an equivalence $e:B\simeq B'$ to a type $B'$ in $\UU$. Then we have the equivalences
  \begin{align*}
    i(X) \in \{g(y)\mid y:B\} & \jdeq \sm{y:B}g(y)=i(X) \\
                              & \simeq \sm{y:B}i(g(y))=X \\
                              & \simeq \sm{y':B'}i(g(e^{-1}(y')))=X \\
                              & \jdeq X\in i(Y).\qedhere
  \end{align*}
\end{proof}

We are now almost in position to reproduce Russell's paradox. We will need one more ingredient: the universal tree, i.e., the multiset of all multisets in $\UU$. 

\begin{defn}
  Let $\UU$ be a universe. Then we define the \define{universal tree}\index{universal tree|textbf}\index{tree!universal tree|textbf}\index{multiset!universal tree|textbf} $\yggdrasil$ to be the multiset\index{Y_U@{$\yggdrasil$}|see {universal tree}}
  \begin{equation*}
    \yggdrasil:=\{i(X) \mid X:\multiset{\UU}\}
  \end{equation*}
  in $\UU^{+}$, where $i:\multiset{\UU}\to\multiset{\UU^+}$ is the inclusion of the multisets in $\UU$ to the multisets in $\UU^+$ given by the fact that each multiset in $\UU$ is $\UU^+$-small.
\end{defn}

\begin{prp}\label{prp:is-small-universal-tree}
  Consider two universes $\UU$ and $\VV$, and suppose that $\UU$ as well as each $X:\UU$ are $\VV$-small. Then the universal tree $\yggdrasil$ is also $\VV$-small.
\end{prp}

\begin{proof}
  To show that the universal tree $\{i(X)\mid X:\multiset{\UU}\}$ is $\VV$-small, we first have to show that the type $\multiset{\UU}$ is $\VV$-small. This follows from the more general fact that the subuniverse of $\VV$-small types is closed under the formation of W-types. Indeed, if a type $A$ is $\VV$-small, and if $B(x)$ is $\VV$-small for each $x:A$, then we have an equivalence $\alpha:A\simeq A'$ to a type $A'$ in $\VV$, and for each $x':A'$ we have an equivalence $B(\alpha^{-1}(x'))\simeq B'(x')$ in $\VV$. These equivalences induce an equivalence
  \begin{equation*}
    \W(A,B)\simeq \W(A',B')
  \end{equation*}
  into the type $W(A',B')$, which is in $\VV$. This concludes the proof that $\multiset{\UU}$ is $\VV$-small.

  It remains to show that the multiset $i(X)$ in $\UU^+$ is $\VV$-small, for each $X:\multiset{\UU}$. Equivalently, we have to show that each multiset $X$ in $\UU$ is $\VV$-small. This follows by recursion: given a multiset $\{f(x)\mid x:A\}$, the type $A$ is $\VV$-small by assumption, and the multiset $f(x)$ is $\VV$-small by the induction hypothesis.
\end{proof}

We are finally ready to employ \define{Russell's paradox} to prove that a univalent universe cannot be equivalent to any type it contains.\index{Russell's paradox|textbf}

\begin{thm}\label{thm:russell}
  Consider a univalent universe $\UU$. Then $\UU$ cannot be $\UU$-small.\index{universe!U is not U-small@{$\UU$ is not $\UU$-small}}\index{U-small type@{$\UU$-small type}!U is not U-small@{$\UU$ is not $\UU$-small}}\index{small type!U is not U-small@{$\UU$ is not $\UU$-small}}
\end{thm}

\begin{proof}
  Suppose that $\UU$ is $\UU$-small, and consider the multiset
  \begin{equation*}
    R\defeq \{i(X) \mid X:\multiset{\UU}, H : X\notin X\}
  \end{equation*}
  in $\UU^+$, where $i:\multiset{\UU}\to\multiset{\UU^+}$ is the inclusion of the multisets in $\UU$ to the multisets in $\UU^+$ given by the fact that each multiset in $\UU$ is $\UU^+$-small.

  First, we note that $R$ is $\UU$-small. This follows from \cref{lem:is-small-comprehension-multiset}, using the fact that the universal tree $\{i(X)\mid X:\multiset{\UU}\}$ is $\UU$-small by \cref{prp:is-small-universal-tree}, and the fact that $X\in X$ is $\UU$-small by \cref{prp:is-small-elementhood-multiset}.

  Since $R$ is $\UU$-small, there is a multiset $R':\multiset{\UU}$ such that $i(R')=R$. Now it follows that
  \begin{align*}
    R\in R & \simeq \sm{X:\multiset{\UU}}\sm{H:X\notin X} i(X)=R \\
           & \simeq \sm{X:\multiset{\UU}}\sm{H:X\notin X} X=R' \\
           & \simeq R'\notin R' \\
           & \simeq R\notin R.
  \end{align*}
  In the second step we used \cref{prp:is-embedding-inclusion-small-multisets}, where we showed that $i$ is an embedding, and in the last step we used \cref{prp:elementhood-small-multisets}. Now we obtain a contradiction, because it follows from \cref{ex:no-fixed-points-neg} that no type is (logically) equivalent to its own negation.
\end{proof}
\index{Russell's paradox|)}
\index{multiset|)}

\begin{exercises}
  \exitem
  \begin{subexenum}
  \item \label{ex:oriented-bin-tree}Let $B:\bool\to\UU$ be the type family defined in \cref{eg:planar-binary-tree-W}. Construct an equivalence\index{oriented binary rooted tree}\index{tree!oriented binary rooted tree}
    \begin{equation*}
      \planarBinTree\simeq\W(\bool,B).
    \end{equation*}
  \item Prove that $\W(\bool,B)$ is not extensional.
  \end{subexenum}
  \exitem Show that for any univalent universe $\UU$ there is no type $U:\UU$ equipped with a surjection $\UU\twoheadrightarrow U$.
  \exitem For a type family $B$ over $A$, suppose that each $B(x)$ is empty. Show that the type $\W(A,B)$ is equivalent to the type $A$.
  \exitem
  \begin{subexenum}
  \item Show that the elementhood relation ${\in}$ on $\W(A,B)$ is irreflexive, for any type family $B$ over any type $A$.\index{elementhood relation on W-types!is irreflexive}\index{W-type!elementhood relation!is irreflexive}
  \item Use the previous fact along with \cref{prp:elementhood-small-multisets} to give a second proof of the fact that there can be no type $U:\UU$ equipped with an equivalence $\UU\simeq U$.
  \end{subexenum}
  \exitem \label{ex:le-W} For each $x:\W(A,B)$, let ${x <(\blank)}:\W(A,B)\to\UU$ be the type family generated inductively by the following constructors:\index{W-type!strict ordering|textbf}\index{strict ordering!on W-types|textbf}
  \begin{align*}
    i & : \prd{y:\W(A,B)}(x\in y) \to (x < y) \\
    j & : \prd{y,z:\W(A,B)} (y\in z) \to ((x<y) \to (x<z)).
  \end{align*}
  \begin{subexenum}
  \item Show that the type-valued relation $<$ is transitive and irreflexive.\index{W-type!strict ordering!is transitive}\index{strict ordering!on W-types!is transitive}\index{W-type!strict ordering!is irreflexive}\index{strict ordering!on W-types!is irreflexive}
  \item Suppose that the type $\W(A,B)$ is inhabited and suppose that there exists an element $a:A$ for which $B(a)$ is inhabited. Show that the following are equivalent:
    \begin{enumerate}
    \item The type $x<y$ is a proposition for all $x,y:\W(A,B)$.
    \item The type $x\in y$ is a proposition for all $x,y:\W(A,B)$.
    \item The type $A$ is a set and the type $B(a)$ is a proposition for all $a:A$.
    \end{enumerate}
    Thus, in general it is not the case that $<$ is a relation valued in propositions.
  \item Show that $\W(A,B)$ satisfies the following \define{strong induction principle}\index{strong induction principle!of W-types|textbf}\index{W-type!strong induction principle|textbf}: For any type family $P$ over $\W(A,B)$, if there is a function
    \begin{equation*}
      h:\prd{x:\W(A,B)}\Big(\prd{y:\W(A,B)} (y<x)\to P(y)\Big)\to P(x),
    \end{equation*}
    then there is a function $f:\prd{x:\W(A,B)}P(x)$ equipped with an identification
    \begin{equation*}
      f(x)=h(x,\lam{y}\lam{p}f(y))
    \end{equation*}
    for all $x:\W(A,B)$.
  \item Show that there can be no sequence of elements $x:\N\to\W(A,B)$ such that $x_{n+1}< x_n$ for all $n:\N$.\index{W-type!no infinitely descending sequences}
  \end{subexenum}
  \exitem (Awodey, Gambino, Sojakova \cite{AwodeyGambinoSojakova}) For any type family $B$ over $A$, the \define{polynomial endofunctor}\index{polynomial endofunctor|textbf} $P_{A,B}$ acts on types by
  \begin{equation*}
    P_{A,B}(X) \defeq \sm{x:A}X^{B(x)},
  \end{equation*}
  and it takes a map $h:X\to Y$ to the map\index{functorial action!polynomial endofunctor|textbf}
  \begin{equation*}
    P_{A,B}(h) : P_{A,B}(X)\to P_{A,B}(Y)
  \end{equation*}
  defined by $P_{A,B}(h,(x,\alpha)) \defeq (x,h\circ \alpha)$. Furthermore, there is a canonical map
  \begin{equation*}
    (h\htpy h') \to (P_{A,B}(h)\htpy P_{A,B}(h'))
  \end{equation*}
  taking a homotopy $H:h\htpy h'$ to a homotopy $P_{A,B}(H):P_{A,B}(h)\htpy P_{A,B}(h')$. 

  A type $X$ is said to be equipped with the \define{structure of an algebra}\index{polynomial endofunctor!algebra|textbf}\index{algebra of a polynomial endofunctor|textbf} for the polynomial endofunctor $P_{A,B}$ if $X$ comes equipped with a map
  \begin{equation*}
    \mu: P_{A,B}(X)\to X.
  \end{equation*}
  Thus, \define{algebras} for the polynomial endofunctor $P_{A,B}$ are pairs $(X,\mu)$ where $X$ is a type and $\mu:P_{A,B}(X)\to X$. Note that $\W(A,B)$ comes equipped with the structure of an algebra for $P_{A,B}$ by \cref{prp:algebra-W}.
  
  Given two algebras $X$ and $Y$ for the polynomial endofunctor $P_{A,B}$, we say that a map $h:X\to Y$ is equipped with the \define{structure of a homomorphism}\index{polynomial endofunctor!morphism of algebras|textbf}\index{morphism of algebras!polynomial endofunctor|textbf} of algebras if it comes equipped with a homotopy witnessing that the square
  \begin{equation*}
    \begin{tikzcd}[column sep=large]
      P_{A,B}(X) \arrow[d,swap,"\mu_X"] \arrow[r,"P_{A,B}(h)"] & P_{A,B}(Y) \arrow[d,"\mu_Y"] \\
      X \arrow[r,swap,"h"] & Y
    \end{tikzcd}
  \end{equation*}
  commutes. The type $\hom((X,\mu_X),(Y,\mu_Y))$ of homomorphisms of algebras for $P_{A,B}$ is therefore defined as
  \begin{equation*}
    \hom((X,\mu_X),(Y,\mu_Y))\defeq \sm{h:X\to Y}h\circ\mu_X\htpy \mu_Y\circ P_{A,B}(h).
  \end{equation*}
  \begin{subexenum}
  \item For any $(x,\alpha),(y,\beta):P_{A,B}(X)$, construct an equivalence
    \begin{equation*}
      ((x,\alpha)=(y,\beta)) \simeq \sm{p:x=y}\alpha\htpy \beta\circ\tr_B(p).
    \end{equation*}
  \item For any two morphisms $(f,K),(g,L):\hom((X,\mu_X),(Y,\mu_Y))$ of algebras for $P_{A,B}$, construct an equivalence
    \begin{equation*}
      ((f,K)=(g,L))\simeq \sm{H:f\htpy g} \ct{K}{(\mu_Y\cdot P_{A,B}(H))}\htpy \ct{(H\cdot \mu_X)}{L}.
    \end{equation*}
  \item Show that the W-type $\W(A,B)$ equipped with the canonical structure $\varepsilon$ of a $P_{A,B}$-algebra, constructed in \cref{prp:algebra-W}, is a \define{(homotopy) initial $P_{A,B}$-algebra} in the sense that the type
    \begin{equation*}
      \hom((\W(A,B),\varepsilon),(X,\mu))
    \end{equation*}
    is contractible, for each $P_{A,B}$-algebra $(X,\mu)$.\index{W-type!is initial algebra of polynomial endofunctor}
  \end{subexenum}
  \exitem Consider the \define{rank comparison relation}\index{W-type!rank comparison relation|textbf}\index{rank comparison relation!W-type|textbf} ${\preceq} : \W(A,B)\to (\W(A,B)\to\prop_\UU)$ defined recursively by\index{preceq@{$\preceq$}|see {W-type, rank comparison relation}}
  \begin{equation*}
    (\collect(a,\alpha)\preceq\collect(b,\beta)) \defeq \forall_{(x:B(a))}\exists_{(y:B(b))}\,\alpha(x)\preceq \beta(y).
  \end{equation*}
  If $x\preceq y$ holds, we say that $x$ has \define{lower rank} than $y$. Furthermore, we define the \define{strict rank comparison relation}\index{W-type!strict rank comparison relation|textbf}\index{strict rank comparison relation!W-type|textbf} ${\prec}$\index{prec@{$\prec$}|see {W-type, strict rank comparison relation}} on $\W(A,B)$ by
  \begin{equation*}
    (x\prec y)\defeq \exists_{(z\in y)}x\preceq z.
  \end{equation*}
  If $x\prec y$ holds, we say that $x$ has \define{strictly lower rank} than $y$.
  \begin{subexenum}
  \item Show that the rank comparison relation defines a preordering on $\W(A,B)$, i.e., show that $\preceq$ is reflexive and transitve\index{W-type!rank comparison relation!is a preordering}\index{rank comparison relation!W-type!is a preordering}. Furthermore, prove the following properties, in which $<$ is the strict ordering on $\W(A,B)$ defined in \cref{ex:le-W}:
    \begin{enumerate}
    \item $(x \preceq y) \leftrightarrow \forall_{(x'<x)}\exists_{(y'<y)}\,x'\preceq y'$
    \item $(x < y)\to (x\preceq y)$
    \item $(x < y) \to (y \npreceq x)$
    \item $\isconstantW(x)\leftrightarrow \forall_{(y:\W(A,B))}\,x \preceq y$.
    \end{enumerate}
  \item Show that the relation $\prec$ on $\W(A,B)$ is a strict ordering on $\W(A,B)$, i.e., show that it is irreflexive and transitive\index{W-type!strict rank comparison relation!is a strict ordering}\index{strict rank comparison relation!W-type!is a strict ordering}. Furthermore, prove the following properties:
    \begin{enumerate}
    \item $(x < y)\to (x\prec y)$
    \item $(x \prec y)\to (x\preceq y)$
    \item $\forall_{(y\preceq y')}\forall_{(x'\preceq x)}(x\prec y)\to (x'\prec y')$.
    \end{enumerate}
  \end{subexenum}
  Since $\preceq$ defines a preordering on $\W(A,B)$, it follows that the preorder $(\W(A,B),\preceq)$ has a poset reflection, in the sense of \cref{ex:poset-reflection}. We will write
  \begin{equation*}
    \eta : (\W(A,B),\preceq)\to(\rank(A,B),\preceq)
  \end{equation*}
  for the poset reflection of $(\W(A,B),\preceq)$\index{R(A,B)@{$\rank(A,B)$}|see {W-type, rank poset}}\index{W-type!rank poset|textbf}\index{rank poset!W-type|textbf} and its quotient map. We will call the poset $(\mathcal{R}(A,B),\preceq)$ the \define{rank poset} of the W-type $\W(A,B)$.
  \begin{subexenum}[resume]
  \item Show that if each $B(x)$ is finite, then the rank poset $(\rank(A,B),\preceq)$ is either the empty poset, the poset with one element, or it is isomorphic to the poset $(\N,\leq)$. 
  \item Show that the strict ordering $\prec$\index{strict rank comparison relation!W-type!extends to rank poset} extends to a relation $\prec$ on $\rank(A,B)$ with the following properties:
    \begin{enumerate}
    \item We have $(x\prec y)\leftrightarrow (\eta(x)\prec\eta(y))$ for every $x,y:\W(A,B)$.
    \item We have $(x\prec y)\to (x\preceq y)$ for every $x,y:\rank(A,B)$.
    \item The relation $\prec$ is transitive and irreflexive on $\rank(A,B)$.
    \end{enumerate}
    We will call the strictly ordered set $(\mathcal{R}(A,B),\prec)$ the \define{(strict) rank}\index{rank!of W-type}\index{W-type!rank} of the W-type $\W(A,B)$. 
  \item A \define{strictly ordered set}\index{strictly ordered set|textbf} $(X,<)$, i.e., a set $X$ equipped with a transitive, irreflexive relation $<$ valued in the propositions, is said to be \define{well-founded}\index{well-founded relation|textbf} if for any family $P$ of propositions over $X$, the implication
  \begin{equation*}
    \Big(\forall_{(x:X)}\Big(\forall_{(y<x)}P(y)\Big)\to P(x)\Big)\to \forall_{(x:X)}P(x).
  \end{equation*}
  holds. Show that the rank $(\rank(A,B),\prec)$ of $\W(A,B)$ is well-founded.\index{rank!of W-type!is well-founded}\index{W-type!rank!is well-founded}
  \item A strictly ordered set $(X,<)$ is said to be \define{extensional}\index{strict ordering!extensional|textbf}\index{extensional strict ordering|textbf} if the logical equivalence
  \begin{equation*}
    (x=y)\leftrightarrow\forall_{(z:X)}\,(z<x)\leftrightarrow(z<y)
  \end{equation*}
  holds for any $x,y:X$. Show that the rank $(\rank(A,B),\prec)$ of $\W(A,B)$ is extensional.\index{rank!of W-type!is extensional}\index{W-type!rank!is extensional}
  \end{subexenum}
\end{exercises}
\index{inductive type|)}
\index{W-type|)}


\cleardoublepage


\chapter{The circle}\label{chapter-synthetic-homotopy-theory}

Many spaces familiar from classical topology have a counterpart in homotopy type theory. For example, there are types representing spheres, projective spaces, and many other CW-complexes. Often such types are constructed as \emph{higher inductive types}. The study of such types leads to the subject of \emph{synthetic homotopy theory}. We conclude this book with a short window into this exciting new subject by introducing the circle as a higher inductive type.

We have seen many examples of inductive types. Those are specified by their (point) constructors, which tell us how we can construct their elements, and an induction principle that allows us to construct sections of type families over them. Inductive types are freely generated by their constructors. Higher inductive types are specified not only by their point constructors, but also by \emph{path constructors}, which are identifications between the point constructors.

The presence of path constructors in a higher inductive type make it slightly more complicated to specify its induction principle. With ordinary inductive types the induction principle only took the point constructors into account. However, when we specify a higher inductive type, we also have to take the path constructors into account. We will see how this works in the example of the circle.

Univalence has a prominent role in the study of higher inductive types. In particular, when we want to characterize the identity type of a higher inductive type we will have to construct a type family over that higher inductive type using the induction principle. Any type family over a higher inductive type must be compatible with the path constructors, and this is where univalence comes in.

The prime example of a type family constructed over a higher inductive type using the univalence axiom is the \emph{universal cover} of the circle. We will see how the univalence axiom can be used to construct the universal cover, and how the universal cover can be used to construct the famous identification
\begin{equation*}
  \pi_1(\sphere{1})=\Z
\end{equation*}
in the type of groups, which asserts that the fundamental group of the circle is the group of integers.


\section{The circle}\label{sec:circle}
\index{circle|(}
\index{inductive type!circle|(}

\subsection{The induction principle of the circle}

The \emph{circle} is specified as a higher inductive type $\sphere{1}$\index{S 1@{$\sphere{1}$}|see {circle}} that comes equipped with\index{base@{$\base$}}\index{loop@{$\lloop$}}\index{circle!base@{$\base$}}\index{circle!loop@{$\lloop$}}
\begin{align*}
\base & : \sphere{1} \\
\lloop & : \id{\base}{\base}.
\end{align*}
Just like for ordinary inductive types, the induction principle for higher inductive types provides us with a way of constructing sections of dependent types. However, we need to take the \emph{path constructor}\index{path constructor} $\lloop$ into account in the induction principle. 

The induction principle of the circle tells us how to define a section
\begin{equation*}
  f:\prd{x:\sphere{1}}P(x)
\end{equation*}
of an arbitrary type family $P$ over $\sphere{1}$. To see what the induction principle of the circle should be, we start with an arbitrary section $f:\prd{x:\sphere{1}}P(x)$ and see how it acts on the constructors of $\sphere{1}$. By applying $f$ to the base point of the circle, we obtain an element $f(\base):P(\base)$. Moreover, using the dependent action on paths\index{dependent action on paths} of $f$ of \cref{defn:apd} we also obtain an identification
\begin{align*}
\apd{f}{\lloop} & : \id{\tr_P(\lloop,f(\base))}{f(\base)}
\end{align*}
in the type $P(\base)$. In other words, we obtain a \emph{dependent action on generators} for every section of a family of types.

\begin{defn}
Let $P$ be a type family over the circle. The \define{dependent action on generators}\index{dependent action on generators for the circle|textbf}\index{circle!dependent action on generators|textbf} is the map\index{dgen_S1@{$\dgen_{\sphere{1}}$}|see {circle, dependent action on generators}}\index{dgen_S1@{$\dgen_{\sphere{1}}$}|textbf}
\begin{equation}\label{eq:dgen_circle}
\dgen_{\sphere{1}}:\Big(\prd{x:\sphere{1}}P(x)\Big)\to\Big(\sm{u:P(\base)}\id{\tr_P(\lloop,u)}{u}\Big)
\end{equation}
given by $\dgen_{\sphere{1}}(f)\defeq\pairr{f(\base),\apd{f}{\lloop}}$.
\end{defn}

The induction principle of the circle states that in order to construct a section $f:\prd{x:\sphere{1}}P(x)$, it suffices to provide an element $u:P(\base)$ and an identification
\begin{equation*}
  \tr_P(\lloop,u)=u.
\end{equation*}
More precisely, the induction principle of the circle is formulated as follows:

\begin{defn}
The \define{circle}\index{circle|textbf}\index{inductive type!circle|textbf} is a type $\sphere{1}$\index{S 1@{$\sphere{1}$}|see {circle}}\index{S 1@{$\sphere{1}$}|textbf} that comes equipped with\index{base@{$\base$}|textbf}\index{loop@{$\lloop$}|textbf}\index{circle!base@{$\base$}|textbf}\index{circle!loop@{$\lloop$}|textbf}\index{higher inductive type!circle|textbf}
\begin{align*}
\base & : \sphere{1} \\
\lloop & : \id{\base}{\base},
\end{align*}
and satisfies the \define{induction principle of the circle}\index{induction principle!of the circle|textbf}\index{circle!induction principle|textbf}, which provides for each type family $P$ over $\sphere{1}$ a map\index{ind S 1@{$\ind{\sphere{1}}$}|textbf}
\begin{equation*}
\ind{\sphere{1}}:\Big(\sm{u:P(\base)}\id{\tr_P(\lloop,u)}{u}\Big)\to \Big(\prd{x:\sphere{1}}P(x)\Big),
\end{equation*}
and a homotopy witnessing that $\ind{\sphere{1}}$ is a section of $\dgen_{\sphere{1}}$
\begin{equation*}
\comphtpy{\sphere{1}}:\dgen_{\sphere{1}}\circ \ind{\sphere{1}}\htpy \idfunc
\end{equation*}
for the computation rules\index{computation rules!for the circle|textbf}\index{circle!computation rules|textbf}.
\end{defn}

\begin{rmk}\label{rmk:circle-induction}
  The type of identifications $(u,p)=(u',p')$ in the type
  \begin{equation*}
    \sm{u:P(\base)}\tr_P(\lloop,u)=u
  \end{equation*}
  is equivalent to the type of pairs $(\alpha,\beta)$ consisting of an identification $\alpha:u=u'$, and an identification $\beta$ witnessing that the square
  \begin{equation*}
    \begin{tikzcd}[column sep=6em]
      \tr_P(\lloop,u) \arrow[d,equals,swap,"p"] \arrow[r,equals,"\ap{\tr_P(\lloop)}{\alpha}"] & \tr_P(\lloop,u') \arrow[d,equals,"{p'}"] \\
      u \arrow[r,equals,swap,"\alpha"] & u'
    \end{tikzcd}
  \end{equation*}
  commutes. Therefore it follows from the induction principle of the circle that for any $(u,p):\sm{u:P(\base)}\tr_P(\lloop,u)=u$, there is a dependent function $f:\prd{x:\sphere{1}}P(x)$ equipped with an identification
  \begin{equation*}
    \alpha : f(\base)=u,
  \end{equation*}
  and an identification $\beta$ witnessing that the square
  \begin{equation*}
    \begin{tikzcd}[column sep=6em]
      \tr_P(\lloop,f(\base)) \arrow[d,equals,swap,"{\apd{f}{\lloop}}"] \arrow[r,equals,"\ap{\tr_P(\lloop)}{\alpha}"] & \tr_P(\lloop,u) \arrow[d,equals,"{p}"] \\
      f(\base) \arrow[r,equals,swap,"\alpha"] & u
    \end{tikzcd}
  \end{equation*}
  commutes.  
\end{rmk}

\subsection{The (dependent) universal property of the circle}
\subsectionmark{The universal property of the circle}
\index{circle!dependent universal property|(}
\index{circle!universal property|(}
\index{universal property!of the circle|(}
\index{dependent universal property!of the circle|(}

We will now use the induction principle of the circle to derive the \emph{dependent universal property} and the \emph{universal property} of the circle. The universal property of the circle states that, for any type $X$ the canonical map
\begin{equation*}
  \Big(\sphere{1}\to X\Big)\to\Big(\sm{x:X}x=x\Big)
\end{equation*}
given by $f\mapsto(f(\base),\ap{f}{\lloop})$ is an equivalence. The type $\sm{x:X}x=x$ is also called the type of \define{free loops}\index{free loop|textbf} in $X$. In other words, the universal property of the circle states that a map $\sphere{1}\to X$ is the same thing as a free loop in $X$.

The \emph{dependent universal property} of the circle similarly states that for any type family $P$ over the circle, the canonical map
\begin{equation*}
  \dgen_{\sphere{1}}:\Big(\prd{x:\sphere{1}}P(x)\Big)\to\Big(\sm{y:P(\base)}\tr_P(\lloop,y)=y\Big)
\end{equation*}
given by $f\mapsto(f(\base),\apd{f}{\lloop})$ is an equivalence. Note that the induction principle already states that this map has a section. The dependent universal property therefore improves on this by stating that this map also has a retraction.

\begin{thm}[The dependent universal property of the circle]\label{thm:circle-dependent-universal-property}
  \index{dependent universal property!of the circle|textbf}
  \index{circle!dependent universal property|textbf}
  For any type family $P$ over the circle, the map
  \begin{equation*}
    \dgen_{\sphere{1}}:
    \Big(\prd{x:\sphere{1}}P(x)\Big)
    \to
    \Big(\sm{y:P(\base)}\tr_P(\lloop,y)=y\Big)
  \end{equation*}
  given by $f\mapsto(f(\base),\apd{f}{\lloop})$ is an equivalence.
\end{thm}

\begin{proof}
  By the induction principle of the circle we know that the map has a section, i.e., we have
  \begin{align*}
    \ind{\sphere{1}} & : \Big(\sm{y:P(\base)}\tr_P(\lloop,y)=y\Big) \to \Big(\prd{x:\sphere{1}}P(x)\Big) \\
    \comphtpy{\sphere{1}} & : \dgen_{\sphere{1}}\circ\ind{\sphere{1}}\htpy\idfunc
  \end{align*}
  Therefore it remains to construct a homotopy
  \begin{equation*}
    \ind{\sphere{1}}\circ\dgen_{\sphere{1}}\htpy\idfunc.
  \end{equation*}
  Thus, for any $f:\prd{x:\sphere{1}}P(x)$ our task is to construct an identification
  \begin{equation*}
    \ind{\sphere{1}}(\dgen_{\sphere{1}}(f))=f.
  \end{equation*}
  By function extensionality it suffices to construct a homotopy
  \begin{equation*}
    \prd{x:\sphere{1}} \ind{\sphere{1}}(\dgen_{\sphere{1}}(f))(x)= f(x).
  \end{equation*}
  We proceed by the induction principle of the circle using the family of types $E_{g,f}(x)\defeq g(x)=f(x)$ indexed by $x:\sphere{1}$, where $g$ is the function
  \begin{equation*}
    g\defeq\ind{\sphere{1}}(\dgen_{\sphere{1}}(f)).
  \end{equation*}
  Thus, it suffices to construct
  \begin{align*}
    \alpha & : g(\base)=f(\base)\\
    \beta  & : \tr_{E_{g,f}}(\lloop,\alpha)=\alpha. 
  \end{align*}
  An argument by path induction on $p$ yields that
  \begin{equation*}
    \Big(\ct{\apd{g}{p}}{r}=\ct{\ap{\tr_P(p)}{q}}{\apd{f}{p}}\Big)\to\Big(\tr_{E_{g,f}}(p,q)=r\Big),
  \end{equation*}
  for any $f,g:\prd{x:X}P(x)$ and any $p:x=x'$, $q:g(x)=f(x)$ and $r:g(x')=f(x')$.
  Therefore it suffices to construct an identification $\alpha:g(\base)=f(\base)$ equipped with an identification $\beta$ witnessing that the square
  \begin{equation*}
    \begin{tikzcd}[column sep=6em]
      \tr_P(\lloop,g(\base)) \arrow[d,equals,swap,"\apd{g}{\lloop}"] \arrow[r,equals,"\ap{\tr_P(\lloop)}{\alpha}"] & \tr_P(\lloop,f(\base)) \arrow[d,equals,"\apd{f}{\lloop}"] \\
      g(\base) \arrow[r,equals,swap,"\alpha"] & f(\base)"
    \end{tikzcd}
  \end{equation*}
  commutes. Notice that we get exactly such a pair $(\alpha,\beta)$ from the computation rule of the circle, by \cref{rmk:circle-induction}.
\end{proof}

As a corollary we obtain the following uniqueness principle for dependent functions defined by the induction principle of the circle.

\begin{cor}
  Consider a type family $P$ over the circle, and let
  \begin{align*}
    y & : P(\base) \\
    p & : \tr_{P}(\lloop,y)=y.
  \end{align*}
  Then the type of functions $f:\prd{x:\sphere{1}}P(x)$ equipped with an identification
  \begin{equation*}
    \alpha: f(\base)=y
  \end{equation*}
  and an identification $\beta$ witnessing that the square
  \begin{equation*}
    \begin{tikzcd}[column sep=6em]
      \tr_P(\lloop,f(\base)) \arrow[d,equals,swap,"{\apd{f}{\lloop}}"] \arrow[r,equals,"\ap{\tr_P(\lloop)}{\alpha}"] & \tr_P(\lloop,y) \arrow[d,equals,"{p}"] \\
      f(\base) \arrow[r,equals,swap,"\alpha"] & y
    \end{tikzcd}
  \end{equation*}
  commutes, is contractible.
\end{cor}

Now we use the dependent universal property to derive the ordinary universal property of the circle. It would be tempting to say that it is a direct corollary, but we need to address the transport that occurs in the dependent universal property.

\begin{thm}[The universal property of the circle]\label{thm:circle_up}
  \index{universal property!of the circle|textbf}
  \index{circle!universal property|textbf}
For each type $X$, the \define{action on generators}\index{action on generators for the circle|textbf}\index{gen_S1@{$\mathsf{gen}_{\sphere{1}}$}|textbf}
\begin{equation*}
\mathsf{gen}_{\sphere{1}}:(\sphere{1}\to X)\to \sm{x:X}x=x
\end{equation*}
given by $f\mapsto (f(\base),\ap{f}{\lloop})$ is an equivalence.
\end{thm}

\begin{proof}
  We prove the claim by constructing a commuting triangle
  \begin{equation*}
    \begin{tikzcd}[column sep=-2em]
      \phantom{\Big(\sm{x:X}\tr_{\const_X}(\lloop,x)=x\Big)} & (\sphere{1}\to X) \arrow[dl,swap,"\gen_{\sphere{1}}"] \arrow[dr,"\dgen_{\sphere{1}}"] \\
      \Big(\sm{x:X}x=x\Big) \arrow[rr,swap,"\simeq"] & & \Big(\sm{x:X}\tr_{\const_X}(\lloop,x)=x\Big)
    \end{tikzcd}
  \end{equation*}
  in which the bottom map is an equivalence. Indeed, once we have such a triangle, we use the fact from \cref{thm:circle-dependent-universal-property} that $\dgen_{\sphere{1}}$ is an equivalence to conclude that $\gen_{\sphere{1}}$ is an equivalence.

  To construct the bottom map, we first observe that for any constant type family $\const_B$ over a type $A$, any $p:a=a'$ in $A$, and any $b:B$, there is an identification
  \begin{equation*}
    \mathsf{tr\usc{}const}_B(p,b):\mathsf{tr}_{\const_B}(p,b)=b.
  \end{equation*}
  This identification is easily constructed by path induction on $p$. Now we construct the bottom map as the induced map on total spaces of the family of maps
  \begin{equation*}
    l\mapsto \ct{\mathsf{tr\usc{}const}_X(\lloop,x)}{l},
  \end{equation*}
  indexed by $x:X$. Since concatenating by a path is an equivalence, it follows by \cref{thm:fib_equiv} that the induced map on total spaces is indeed an equivalence.

  To show that the triangle commutes, it suffices to construct for any $f:\sphere{1}\to X$ an identification witnessing that the triangle
  \begin{equation*}
    \begin{tikzcd}[column sep=1em]
      \tr_{\const_X}(\lloop,f(\base)) \arrow[dr,equals,swap,"\apd{f}{\lloop}"] \arrow[rr,equals,"{\mathsf{tr\usc{}const}_X(\lloop,f(\base))}"] & & f(\base) \arrow[dl,equals,"\ap{f}{\lloop}"] \\
      & f(\base) & \phantom{\tr_{\const_X}(\lloop,f(\base))}
    \end{tikzcd}
  \end{equation*}
  commutes. This again follows from general considerations: for any $f:A\to B$ and any $p:a=a'$ in $A$, the triangle
  \begin{equation*}
    \begin{tikzcd}[column sep=1em]
      \tr_{\const_B}(p,f(a)) \arrow[dr,equals,swap,"\apd{f}{p}"] \arrow[rr,equals,"{\mathsf{tr\usc{}const}_B(p,f(a))}"] & & f(a) \arrow[dl,equals,"\ap{f}{p}"] \\
      & f(a') & \phantom{\tr_{\const_B}(p,f(a))}
    \end{tikzcd}
  \end{equation*}
  commutes by path induction on $p$.
\end{proof}

\begin{cor}
  For any loop $l:x=x$ in a type $X$, the type of maps $f:\sphere{1}\to X$ equipped with an identification
  \begin{equation*}
    \alpha : f(\base)=x 
  \end{equation*}
  and an identification $\beta$ witnessing that the square
  \begin{equation*}
    \begin{tikzcd}
      f(\base) \arrow[r,equals,"\alpha"] \arrow[d,equals,swap,"\ap{f}{\lloop}"] & x \arrow[d,equals,"l"] \\
      f(\base) \arrow[r,equals,swap,"\alpha"] & x
    \end{tikzcd}
  \end{equation*}
  commutes, is contractible.
\end{cor}
\index{circle!dependent universal property|)}
\index{circle!universal property|)}
\index{universal property!of the circle|)}
\index{dependent universal property!of the circle|)}

\subsection{Multiplication on the circle}
\label{sec:mulcircle}
\index{circle!H-space structure|(}
\index{circle!multiplication|(}
\index{multiplication!on the circle|(}

One way the circle arises classically, is as the set of complex numbers at distance $1$ from the origin. It is an elementary fact that $|xy|=|x||y|$ for any two complex numbers $x,y\in\mathbb{C}$, so it follows that when we multiply two complex numbers that both lie on the unit circle, then the result lies again on the unit circle. This operation puts a group structure on the classical circle.

This suggests that it should also be possible to construct a multiplication on the higher inductive type $\sphere{1}$. More precisely, we will equip $\sphere{1}$ with an \emph{H-space structure}, and in the exercises you will be asked to show that this multiplicative structure is associative, commutative, and has inverses.

\begin{defn}
  Consider a pointed type $A$ with a base point $\pt$. An \textbf{H-space structure}\index{H-space!structure|textbf} on $(A,\pt)$ consists of a binary operation $\mu:A\to (A\to A)$ satisfying the following \define{coherent unit laws}\index{coherent unit laws|textbf}\index{unit laws!coherent unit laws|textbf}\index{H-space!coherent unit laws|textbf}:
  \begin{align*}
    \leftunit_\mu(y) & : \mu(\pt,y)= y \\
    \rightunit_\mu(x) & : \mu(x,\pt)= x \\
    \cohunit_\mu    & : \leftunit_\mu(\pt)=\rightunit_\mu(\pt).
  \end{align*}
  An \define{H-space}\index{H-space|textbf} is a pointed type equipped with an H-space structure.
\end{defn}

\begin{rmk}\label{rmk:hspace}
  The data of an H-space structure is equivalently described by a family of base point preserving maps
  \begin{equation*}
    \mu : \prd{x:A}\sm{f:A\to A}f(\pt)=x
  \end{equation*}
  equipped with an identification $\mu_\pt=(\idfunc,\refl{})$. The data $\mu(a,\pt)=a$ corresponds to the right unit law for $\mu$, whereas the data $\mu_\pt=(\idfunc,\refl{})$ combines the left unit law and the coherence in one single identification.

  Note that for any identification $\alpha:x=y$ in $A$ and two base-point preserving functions $(f,p):\sm{f:A\to A}f(\pt)=x$ and $(g,q):\sm{f:A\to A}f(\pt)=y$, we have
  \begin{equation*}
    \tau:\Big(\sm{H:f\htpy g}\ct{p}{\alpha}=\ct{H(\pt)}{q}\Big) \to \tr(\alpha,(f,p))=(g,q)
  \end{equation*}
  This function is easily constructed by identification elimination on $\alpha$. We will be using this in our construction of the H-space structure on the circle.
\end{rmk}

\begin{thm}\label{defn:hspace-circle}
  There is an H-space structure\index{H-space!circle|textbf}\index{circle!multiplication|textbf}\index{circle!H-space structure|textbf}\index{unit laws!for multiplication on S1@{for multiplication on $\sphere{1}$}}\index{mul S 1@{$\mulcircle$}}
  \begin{align*}
    \mulcircle & : \sphere{1}\to(\sphere{1}\to\sphere{1}) \\
    \leftunit_{\sphere{1}} & : \prd{y:\sphere{1}}\mulcircle(\base,y)=y \\
    \rightunit_{\sphere{1}} & : \prd{x:\sphere{1}}\mulcircle(x,\base)=x \\
    \cohunit_{\sphere{1}} & : \leftunit_{\sphere{1}}(\base)=\rightunit_{\sphere{1}}(\base).
  \end{align*}
  on the circle.
\end{thm}

\begin{proof}[Construction]
  By \cref{rmk:hspace} it suffices to construct a dependent function
  \begin{equation*}
    \mu:\prd{x:\sphere{1}}\sm{f:\sphere{1}\to\sphere{1}}f(\base)=x
  \end{equation*}
  such that $\mu(\base)=(\idfunc,\refl{})$. This provides us with a useful shortcut, because the identification will follow from the computation rule of the induction principle of the circle.

  Let $P$ be the family of types given by $P(x):=\sm{f:\sphere{1}\to\sphere{1}}f(\base)=x$. By the dependent universal property of the circle there is a unique
  \begin{equation*}
    \mu :\prd{x:\sphere{1}}\sm{f:\sphere{1}\to\sphere{1}}f(\base)=x
  \end{equation*}
  equipped with an identification $\alpha:\mu(\base)=(\idfunc,\refl{})$ and an identification witnessing that the square
  \begin{equation*}
    \begin{tikzcd}[column sep=6em]
      \tr_P(\lloop,\mu(\base)) \arrow[d,equals,swap,"\apd{\mu}{\lloop}"] \arrow[r,equals,"\ap{\tr_P(\lloop)}{\alpha}"] & \tr_P(\lloop,(\idfunc,\refl{})) \arrow[d,equals,"{\tau(\htpyidcircle,r)}"] \\
      \mu(\base) \arrow[r,equals,swap,"\alpha"] & (\idfunc,\refl{})
    \end{tikzcd}
  \end{equation*}
  commutes. In this square, $\tau$ is the function from \cref{rmk:hspace}, and the homotopy $\htpyidcircle:\idfunc\htpy\idfunc$ equipped with an identification $r:\lloop = \ct{\htpyidcircle(\base)}{\refl{}}$ remain to be defined.

  We use the dependent universal property of the circle with respect to the family $E_{\idfunc,\idfunc}$ given by
  \begin{equation*}
    E_{\idfunc,\idfunc}(x) \defeq (x=x),
  \end{equation*}
  to define $\htpyidcircle$ as the unique homotopy equipped with an identification
  \begin{equation*}
    \basehtpyidcircle : \htpyidcircle(\base)=\lloop
  \end{equation*}
  and an identification $\loophtpyidcircle$ witnessing that the square
  \begin{equation*}
    \begin{tikzcd}[column sep=8em]
      \tr_{E_{\idfunc,\idfunc}}(\lloop,\htpyidcircle(\base)) \arrow[r,equals,"\ap{\tr_{E_{\idfunc,\idfunc}}(\lloop)}{\basehtpyidcircle}"] \arrow[d,equals,swap,"\apd{\htpyidcircle}{\lloop}"] & \tr_{E_{\idfunc,\idfunc}}(\lloop,\lloop) \arrow[d,equals,"\gamma"] \\
      \htpyidcircle(\base) \arrow[r,equals,swap,"\basehtpyidcircle"] & \lloop
    \end{tikzcd}
  \end{equation*}
  commutes. Now it remains to define the path $\gamma:\tr_{E_{\idfunc,\idfunc}}(\lloop,\lloop)=\lloop$ in the above square. To proceed, we first observe that a simple path induction argument yields a function
  \begin{equation*}
    \Big(\ct{p}{r}=\ct{q}{p}\Big)\to\Big(\tr_{E_{\idfunc,\idfunc}}(p,q)=r\Big),
  \end{equation*}
  for any $p:\base=x$, $q:\base=\base$ and $r:x=x$. In particular, we have a function
  \begin{equation*}
    \Big(\ct{\lloop}{\lloop}=\ct{\lloop}{\lloop}\Big)\to\Big(\tr_{E_{\idfunc,\idfunc}}(\lloop,\lloop)=\lloop\Big).
  \end{equation*}
  Now we apply this function to $\refl{\ct{\lloop}{\lloop}}$ to obtain the desired identification
  \begin{equation*}
    \gamma:\tr_{E_{\idfunc,\idfunc}}(\lloop,\lloop)=\lloop.\qedhere
  \end{equation*}
\end{proof}

\begin{rmk}
  For some of the exercises below it may be useful to know that the binary operation $\mulcircle$ is the unique map $\sphere{1}\to(\sphere{1}\to\sphere{1})$ equipped with an identification
  \begin{equation*}
    \basemulcircle :\mulcircle(\base)=\idfunc
  \end{equation*}
  and an identification $\loopmulcircle$ witnessing that the square
  \begin{equation*}
    \begin{tikzcd}[column sep=huge]
      \mulcircle(\base) \arrow[r,equals,"\basemulcircle"] \arrow[d,equals,swap,"\ap{\mulcircle}{\lloop}"] & \idfunc \arrow[d,equals,"\eqhtpy(\htpyidcircle)"] \\
      \mulcircle(\base) \arrow[r,equals,swap,"\basemulcircle"] & \idfunc
  \end{tikzcd}
  \end{equation*}
  commutes, where the homotopy $\htpyidcircle:\idfunc\htpy\idfunc$ is the one constructed in \cref{defn:hspace-circle}.
\end{rmk}

\index{circle!H-space structure|)}
\index{circle!multiplication|)}
\index{multiplication!on the circle|)}

\begin{exercises}
  \exitem \label{ex:circle-constant}
  \begin{subexenum}
  \item Show that for any type $X$ and any $x:X$, the map
    \begin{equation*}
      \ind{\sphere{1}}(x,\refl{x}):\sphere{1}\to X
    \end{equation*}
    is homotopic to the constant map $\mathsf{const}_x$.\index{circle!constant maps}\index{constant map!on the circle}
  \item Show that
    \begin{equation*}
      \ind{\sphere{1}}(\base,\lloop) : \sphere{1}\to\sphere{1}
    \end{equation*}
    is homotopic to the identity function.
  \item Consider a map $f:X\to Y$ and a free loop\index{free loop} $(x,l)$ in $X$. Construct a homotopy
    \begin{equation*}
      \ind{\sphere{1}}(f(x),\ap{f}{l})\htpy f\circ \ind{\sphere{1}}(x,l).
    \end{equation*}
  \end{subexenum}
  \exitem \label{ex:circle-connected}
  \begin{subexenum}
  \item Show that the circle is connected.\index{circle!is connected}\index{is connected!circle}
  \item Let $P:\sphere{1}\to\prop$ be a family of propositions over the circle. Show that
    \begin{equation*}
      P(\base)\to\prd{x:\sphere{1}}P(x).
    \end{equation*}
  \item Show that any embedding $m:\sphere{1}\to\sphere{1}$ is an equivalence.
  \item Show that for any embedding $m:X\to\sphere{1}$, there is a proposition $P$ and an equivalence $e:\eqv{X}{\sphere{1}\times P}$ for which the triangle
    \begin{equation*}
      \begin{tikzcd}[column sep=0]
        X \arrow[dr,swap,"m"] \arrow[rr,"e"] & & \sphere{1}\times P \arrow[dl,"\proj 1"] \\
        \phantom{\sphere{1}\times P} & \sphere{1}
      \end{tikzcd}
    \end{equation*}
    commutes. In other words, all the embeddings into the circle are of the form $\sphere{1}\times P\to \sphere{1}$.
  \end{subexenum}
  \exitem \label{ex:mulcircle-is-equiv}
  \begin{subexenum}
  \item Show that for any $x:\sphere{1}$, both functions
    \begin{equation*}
      \mulcircle(x,\blank)\qquad\text{and}\qquad\mulcircle(\blank,x)
    \end{equation*}
    are equivalences.\index{is an equivalence!mul S 1 (x,-)@{$\mulcircle(x,\blank)$}}\index{is an equivalence!mul S 1 (-,x)@{$\mulcircle(\blank,x)$}}
  \item Show that the function\index{is an embedding!mul S 1@{$\mulcircle$}}
    \begin{equation*}
      \mulcircle : \sphere{1}\to(\sphere{1}\to\sphere{1})
    \end{equation*}
    is an embedding.
  \end{subexenum}
  \exitem \label{ex:circle_connected}
  \begin{subexenum}
  \item Show that a type $X$ is a set\index{set} if and only if the map
    \begin{equation*}
      \lam{x}\lam{t} x : X \to (\sphere{1}\to X)
    \end{equation*}
    is an equivalence.
  \item Show that a type $X$ is a set\index{set} if and only if the map
    \begin{equation*}
      \lam{f}f(\base) : (\sphere{1}\to X)\to X
    \end{equation*}
    is an equivalence.
  \end{subexenum}
  \exitem Show that the multiplicative operation on the circle is associative and commutative, i.e.~construct an identifications\index{circle!associativity of multiplication}\index{associativity!of multiplication on S 1@{of multiplication on $\sphere{1}$}}\index{circle!commutativity of multiplication}\index{commutativity!of multiplication on S 1@{of multiplication on $\sphere{1}$}}
  \begin{align*}
    \mulcircle(\mulcircle(x,y),z) & = \mulcircle(x,\mulcircle(y,z)) \\
    \mulcircle(x,y) & =\mulcircle(y,x).
  \end{align*}
  for every $x,y,z:\sphere{1}$.
  \exitem Show that the circle, equipped with the multiplicative operation $\mulcircle$ is an abelian group, i.e.~construct an inverse operation
  \begin{equation*}
    \invcircle : \sphere{1}\to\sphere{1}
  \end{equation*}
  and construct identifications
  \begin{align*}
    \leftinv_{\sphere{1}} & : \mulcircle(\invcircle(x),x) = \base \\
    \rightinv_{\sphere{1}} & : \mulcircle(x,\invcircle(x)) = \base.
  \end{align*}
  Moreover, show that the square
  \begin{equation*}
    \begin{tikzcd}
      \invcircle(\base) \arrow[d,equals] \arrow[r,equals] & \mulcircle(\base,\invcircle(\base)) \arrow[d,equals] \\
      \mulcircle(\invcircle(\base),\base) \arrow[r,equals] & \base
    \end{tikzcd}
  \end{equation*}
  commutes.
  \exitem Show that for any multiplicative operation
  \begin{equation*}
    \mu:\sphere{1}\to(\sphere{1}\to\sphere{1})
  \end{equation*}
  that satisfies the condition that $\mu(x,\blank)$ and $\mu(\blank,x)$ are equivalences for any $x:\sphere{1}$, there is an element $e:\sphere{1}$ such that
  \begin{equation*}
    \mu(x,y)=\mulcircle(x,\mulcircle(\bar{e},y))
  \end{equation*}
  for every $x,y:\sphere{1}$, where $\bar{e}\defeq\invcircle(e)$ is the complex conjugation of $e$ on $\sphere{1}$.
  \exitem Consider a pointed type $(A,\pt)$ equipped with a \define{noncoherent H-space structure}\index{H-space!noncoherent H-space|textbf}\index{noncoherent H-space|textbf}\index{H-space!coherence theorem} $(\mu,H,K)$ consisting of a binary operation $\mu:A\to (A\to A)$ and homotopies
  \begin{align*}
    H & : \prd{y:A}\mu(\pt,y)=y \\
    K & : \prd{x:A}\mu(x,\pt)=x.
  \end{align*}
  Show that the homotopy $K$ can be adjusted to a new homotopy $K':\prd{x:A}\mu(x,\pt)=x$ in such a way that
  \begin{equation*}
    H(\pt)=K'(\pt)
  \end{equation*}
  holds. In other words, any noncoherent H-space structure can be improved to an H-space structure with the same underlying binary operation. Hint: Take some inspiration from \cref{lem:coherently-invertible}, where one of the homotopies of the invertibility of a map was adjusted to obtain coherent invertibility.
\end{exercises}


\section{The universal cover of the circle}\label{sec:circle-universal-cover}
\index{circle!universal cover|(}
\index{universal cover of S 1@{universal cover of $\sphere{1}$}|(}

In this section we use the univalence axiom to construct the \emph{universal cover} of the circle and show that the loop space of the circle is equivalent to $\mathbb{Z}$. The universal cover of the circle is a family of sets over the circle with contractible total space.
Classically, the universal cover is described as a map $\mathbb{R}\to\sphere{1}$ that winds the real line around the circle. In homotopy type theory the universal cover is constructed as a map $\sphere{1}\to\Set$ into the univalent type of all sets, and we will use the dependent universal property of the circle to show that its total space is contractible.

\subsection{The universal cover of the circle}

The type of small families over $\sphere{1}$ is just the function type $\sphere{1}\to\UU$. Therefore, we may use the universal property of the circle to construct type families over the circle.

By the universal property, $\UU$-small type families over $\sphere{1}$ are equivalently described as pairs $(X,p)$ consisting of a type $X:\UU$ and an identification $p:X=X$. The univalence axiom\index{univalence axiom!families over $\sphere{1}$} implies that the map
\begin{equation*}
\mathsf{eq\usc{}equiv}_{X,X}:(\eqv{X}{X})\to (X=X)
\end{equation*}
is an equivalence. Therefore, type families over the circle are equivalently described as pairs $(X,e)$, consisting of a type $X$ and an equivalence $e:\eqv{X}{X}$. The type $\sm{X:\UU}\eqv{X}{X}$ is also called the type of \define{descent data}\index{descent data for the circle|textbf}\index{circle!descent data|textbf} for the circle.

\begin{defn}\label{defn:circle_descent}
Consider a type $X$ and an equivalence $e:\eqv{X}{X}$.
We will construct a dependent type $\mathcal{D}(X,e):\sphere{1}\to\UU$\index{D (X,e)@{$\mathcal{D}(X,e)$}|textbf}\index{D (X,e)@{$\mathcal{D}(X,e)$}|see {type family, over $\sphere{1}$}}\index{type family!over S 1@{over $\sphere{1}$}|textbf}\index{circle!type family over S 1@{type family over $\sphere{1}$}|textbf} equipped with an equivalence $x\mapsto x_{\mathcal{D}}:\eqv{X}{\mathcal{D}(X,e,\base)}$ for which the square
\begin{equation*}
\begin{tikzcd}
X \arrow[r,"\eqvsym"] \arrow[d,swap,"e"] & \mathcal{D}(X,e,\base) \arrow[d,"\mathsf{tr}_{\mathcal{D}(X,e)}(\lloop)"] \\
X \arrow[r,swap,"\eqvsym"] & \mathcal{D}(X,e,\base)
\end{tikzcd}
\end{equation*}
commutes. We will write $d\mapsto d_{X}$ for the inverse of this equivalence, so that the relations
\begin{samepage}%
\begin{align*}
(x_{\mathcal{D}})_X & =x & (e(x)_{\mathcal{D}}) & = \mathsf{tr}_{\mathcal{D}(X,e)}(\lloop,x_{\mathcal{D}}) \\
(d_X)_{\mathcal{D}} & =d & (\mathsf{tr}_{\mathcal{D}(X,e)}(d))_X & = e(d_X)
\end{align*}
\end{samepage}%
hold.
\end{defn}

\begin{constr}
  An easy path induction argument reveals that
\begin{equation*}
\mathsf{equiv\usc{}eq}(\ap{P}{\lloop})=\mathsf{tr}_P(\lloop)
\end{equation*}
for each dependent type $P:\sphere{1}\to\UU$. Therefore we see that the triangle\index{desc_S1@{$\mathsf{desc}_{\sphere{1}}$}|textbf}
\begin{equation*}
\begin{tikzcd}
& (\sphere{1}\to \UU) \arrow[dl,swap,"\mathsf{gen}_{\sphere{1}}"] \arrow[dr,"\mathsf{desc}_{\sphere{1}}"] \\
\sm{X:\UU}X=X \arrow[rr,swap,"\tot{\lam{X}\mathsf{equiv\usc{}eq}_{X,X}}"] & & \sm{X:\UU}\eqv{X}{X}
\end{tikzcd}
\end{equation*}
commutes, where the map $\mathsf{desc}_{\sphere{1}}$ is given by $P\mapsto\pairr{P(\base),\mathsf{tr}_P(\lloop)}$ and the bottom map is an equivalence by the univalence axiom and \cref{thm:fib_equiv}.
Now it follows by the 3-for-2 property that $\mathsf{desc}_{\sphere{1}}$ is an equivalence, since $\mathsf{gen}_{\sphere{1}}$ is an equivalence by \cref{thm:circle_up}.
This means that for every type $X$ and every $e:\eqv{X}{X}$ there is a type family $\mathcal{D}(X,e):\sphere{1}\to\UU$ equipped with an identification
\begin{equation*}
\pairr{\mathcal{D}(X,e,\base),\mathsf{tr}_{\mathcal{D}(X,e)}(\lloop)}=\pairr{X,e}.
\end{equation*}
For convenience, we invert this identification. Now we observe that the type of identifications in $\sm{X:\UU}\eqv{X}{X}$ can be characterized by
\begin{equation*}
  ((X,e)=(X',e'))\simeq \sm{\alpha:X\simeq X'} e'\circ \alpha\htpy \alpha\circ e'.
\end{equation*}
This implies that we obtain an equivalence $x\mapsto x_{\mathcal{D}}:X\simeq D(X,e,\base)$ such that the square
\begin{equation*}
\begin{tikzcd}[column sep=huge]
X \arrow[d,swap,"e"] \arrow[r,"x\mapsto x_{\mathcal{D}}"] & \mathcal{D}(X,e,\base) \arrow[d,"\mathsf{tr}_{\mathcal{D}(X,e)}(\lloop)"] \\
X \arrow[r,swap,"x\mapsto x_{\mathcal{D}}"] & \mathcal{D}(X,e,\base)
\end{tikzcd}
\end{equation*}
commutes.
\end{constr}

Recall from \cref{eg:is-equiv-succ-Z} that the successor function $\succZ :\Z\to \Z$ is an equivalence. Its inverse is the predecessor function defined in \cref{ex:int_pred}. 

\begin{defn}\label{defn:universal-cover-circle}
  The \define{universal cover}\index{universal cover of S 1@{universal cover of $\sphere{1}$}|textbf}\index{circle!universal cover|textbf} of the circle is defined via \cref{defn:circle_descent} to be the unique dependent type $\universalcovercircle\defeq\mathcal{D}(\Z,\succZ ):\sphere{1}\to\UU$.\index{Z@{$\Z$}!universal cover of S 1@{universal cover of $\sphere{1}$}|textbf}\index{E S 1@{$\universalcovercircle$}|see {universal cover of $\sphere{1}$}}\index{E S 1@{$\universalcovercircle$}|textbf} equipped with an equivalence $x\mapsto x_\mathcal{E}:\Z\to\universalcovercircle(\base) $ and a homotopy witnessing that the square
  \begin{equation*}
    \begin{tikzcd}[column sep=large]
      \mathbb{Z} \arrow[r,"x\mapsto x_\mathcal{E}"] \arrow[d,swap,"\succZ "] & \universalcovercircle(\mathsf{base}) \arrow[d,"\mathsf{tr}_{\universalcovercircle}(\mathsf{loop})"] \\
      \mathbb{Z} \arrow[r,swap,"x\mapsto x_\mathcal{E}"] & \universalcovercircle(\mathsf{base})
    \end{tikzcd}
  \end{equation*}
  commutes. We will occasionally write $y\mapsto y_\Z$ for the inverse of $x\mapsto x_{\mathcal{E}}$.
\end{defn}

The picture of the universal cover is that of a helix\index{helix} over the circle. This picture emerges from the path liftings of $\mathsf{loop}$ in the total space. The segments of the helix connecting $k$ to $k+1$ in the total space of the helix, are constructed in the following lemma.

\begin{lem}
For any $k:\Z$, there is an identification\index{segment-helix@{$\segmenthelix_k$}|textbf}\index{universal cover of S 1@{universal cover of $\sphere{1}$}!segment-helis@{$\segmenthelix_k$}}
\begin{equation*}
\segmenthelix_k:(\base,k_{\mathcal{E}})=(\base,\succZ (k)_{\mathcal{E}})
\end{equation*}
in the total space $\sm{t:\sphere{1}}\mathcal{E}(t)$.\index{universal cover of S 1@{universal cover of $\sphere{1}$}!total space}\index{total space!universal cover of S 1@{universal cover of $\sphere{1}$}}
\end{lem}

\begin{proof}
By \cref{thm:eq_sigma} it suffices to show that
\begin{equation*}
\prd{k:\Z} \sm{\alpha:\base=\base} \mathsf{tr}_{\mathcal{E}}(\alpha,k_{\mathcal{E}})= \succZ (k)_{\mathcal{E}}.
\end{equation*}
We just take $\alpha\defeq\lloop$. Then we have $\mathsf{tr}_{\mathcal{E}}(\alpha,k_{\mathcal{E}})= \succZ (k)_{\mathcal{E}}$ by the commuting square provided in the definition of $\mathcal{E}$.
\end{proof}

\subsection{Working with descent data}
\index{circle!descent data|(}
\index{descent data for the circle|(}

The equivalence
\begin{equation*}
  (\sphere{1}\to\UU)\simeq \sm{X:\UU}X\simeq X
\end{equation*}
yields that for any type family $A$ over the circle the type of descent data $(X,e)$ equipped with an equivalence $\alpha:X\simeq A(\base)$ and a homotopy $H$ witnessing that the square
\begin{equation*}
  \begin{tikzcd}
    X \arrow[r,"\alpha"] \arrow[d,swap,"e"] & A(\base) \arrow[d,"\tr_A(\lloop)"] \\
    X \arrow[r,swap,"\alpha"] & A(\base)
  \end{tikzcd}
\end{equation*}
commutes is contractible. In the remainder of this section we study arbitrary type families over the circle equipped with such descent data, which will put us in a good position to prove things about the universal cover of the circle.

\begin{prp}
  Consider a type family $A$ over the circle and consider descent data $(X,e)$ equipped with an equivalence $\alpha:X\simeq A(\base)$ and a homotopy witnessing that the square
  \begin{equation*}
    \begin{tikzcd}
      X \arrow[r,"\alpha"] \arrow[d,swap,"e"] & A(\base) \arrow[d,"\tr_A(\lloop)"] \\
      X \arrow[r,swap,"\alpha"] & A(\base)
    \end{tikzcd}
  \end{equation*}
  commutes. Furthermore, consider two elements $x,y:X$. Then we have an equivalence
  \begin{equation*}
    \bar{\alpha}:(e(x)=y)\simeq (\tr_{A}(\lloop,\alpha(x))=\alpha(y)).
  \end{equation*}
\end{prp}

\begin{proof}
  Note that the commutativity of the square implies that
  \begin{equation*}
    \tr_A(\lloop,\alpha(x))=\alpha(e(x)).
  \end{equation*}
  By \cref{thm:id_fundamental} it therefore suffices to prove that the total space
  \begin{equation*}
    \sm{y:X}\tr_A(\lloop,\alpha(x))=\alpha(y)
  \end{equation*}
  is contractible. This type is equivalent to $\fib{\alpha}{\tr_A(\lloop,\alpha(x))}$, which is contractible because $\alpha$ is an equivalence.
\end{proof}

In the following proposition we show that sections of a type family $A$ equipped with descent data $(X,e)$ are equivalently described as fixed points for $e:X\simeq X$. 

\begin{prp}
  Consider a type family $A$ over the circle and descent data $(X,e)$ equipped with an equivalence $\alpha:X\simeq A(\base)$ and a homotopy witnessing that the square
  \begin{equation*}
    \begin{tikzcd}
      X \arrow[r,"\alpha"] \arrow[d,swap,"e"] & A(\base) \arrow[d,"\tr_A(\lloop)"] \\
      X \arrow[r,swap,"\alpha"] & A(\base)
    \end{tikzcd}
  \end{equation*}
  commutes. Then there is a commuting square
  \begin{equation*}
    \begin{tikzcd}
      \prd{t:\sphere{1}}A(t) \arrow[d,swap,"\ev_\base"] \arrow[r] & \sm{x:X}e(x)=x \arrow[d,"\proj 1"] \\
      A(\base) \arrow[r,swap,"\alpha^{-1}"] & X
    \end{tikzcd}
  \end{equation*}
  in which the top map is an equivalence.
\end{prp}

\begin{proof}
  By the dependent universal property of the circle we have an equivalence
  \begin{equation*}
    \Big(\prd{t:\sphere{1}}A(t)\Big)\simeq \sm{x:A(\base)}\tr_A(\lloop,x)=x.
  \end{equation*}
  This equivalence fits in a commuting triangle
  \begin{equation*}
    \begin{tikzcd}[column sep=-2em]
      \phantom{\sm{x:A(\base)}\tr_A(\lloop,x)=x} & \prd{t:\sphere{1}}A(t) \arrow[dl] \arrow[dr,"\dgen_{\sphere{1}}"] \\
      \sm{x:X}e(x)=x \arrow[rr,swap,"{\tot[\alpha]{\bar{\alpha}}}"] & & \sm{x:A(\base)}\tr_A(\lloop,x)=x 
    \end{tikzcd}
  \end{equation*}
  where the map on the left is given by $s\mapsto(\alpha^{-1}(s(\base)),\bar{\alpha}^{-1}(\apd{s}{\lloop}))$. The bottom map and the map on the right are equivalences, so it follows by the 3-for-2 property of equivalences that the map on the left is an equivalence.
\end{proof}

The following corollary can be used to compare type families over the circle. In particular, we will use it to compare the identity type of the circle with the universal cover.

\begin{cor}
  Consider two type families $A$ and $B$ over the circle equipped with descent data $(X,e)$ and $(Y,f)$, equivalences $\alpha:X\simeq A(\base)$ and $\beta:Y\simeq B(\base)$, and homotopies $H$ and $K$ witnessing that the squares
  \begin{equation*}
    \begin{tikzcd}
      X \arrow[r,"\alpha"] \arrow[d,swap,"e"] & A(\base) \arrow[d,"\tr_A(\lloop)"] & Y \arrow[r,"\beta"] \arrow[d,swap,"f"] & B(\base) \arrow[d,"\tr_B(\lloop)"] \\
      X \arrow[r,swap,"\alpha"] & A(\base) & Y \arrow[r,swap,"\beta"] & B(\base)
    \end{tikzcd}
  \end{equation*}
  commute, respectively. Then there is a commuting square
  \begin{equation*}
    \begin{tikzcd}
      \Big(\prd{t:\sphere{1}}A(t)\to B(t)\Big) \arrow[r] \arrow[d,swap,"\ev_\base"] & \sm{h:X\to Y}h\circ e\htpy f\circ h \arrow[d,"\proj 1"] \\
      (A(\base)\to B(\base)) \arrow[r,swap,"h\mapsto \beta^{-1}\circ h\circ \alpha"] & (X\to Y)
    \end{tikzcd}
  \end{equation*}
  in which the top map is an equivalence.
\end{cor}

\begin{proof}
  The claim follows once we observe that $(Y^X,\lam{h}f\circ h\circ e^{-1})$ is descent data for the family of types $(A(t)\to B(t))$ indexed by $t:\sphere{1}$. Indeed, we have the equivalence $h\mapsto \beta\circ h\circ\alpha^{-1} : Y^X\simeq B(\base)^{A(\base)}$ for which the square
  \begin{equation*}
    \begin{tikzcd}[column sep=6em]
      Y^X \arrow[r,"h\mapsto\beta\circ h\circ\alpha^{-1}"] \arrow[d,swap,"h\mapsto f\circ h\circ e^{-1}"] & B(\base)^{A(\base)} \arrow[d,"\tr_{t\mapsto A(t)\to B(t)}(\lloop)"] \\
      Y^X \arrow[r,swap,"h\mapsto\beta\circ h\circ\alpha^{-1}"] & B(\base)^{A(\base)}
    \end{tikzcd}
  \end{equation*}
  commutes. 
\end{proof}

\begin{cor}\label{cor:compute-families-of-maps-universal-cover}
  Consider a type family $A$ over the circle and descent data $(X,e)$ equipped with an equivalence $\alpha:X\simeq A(\base)$ and a homotopy witnessing that the square
  \begin{equation*}
    \begin{tikzcd}
      X \arrow[r,"\alpha"] \arrow[d,swap,"e"] & A(\base) \arrow[d,"\tr_A(\lloop)"] \\
      X \arrow[r,swap,"\alpha"] & A(\base)
    \end{tikzcd}
  \end{equation*}
  commutes. Then there is a commuting square
  \begin{equation*}
    \begin{tikzcd}[column sep=3em]
      \Big(\prd{t:\sphere{1}}\universalcovercircle(t)\to A(t)\Big) \arrow[r] \arrow[d,swap,"\ev_\base"] & \sm{h:\Z \to X} h\circ \succZ \htpy e\circ h \arrow[d,"\proj 1"] \\
      (\universalcovercircle(\base)\to A(\base)) \arrow[r,swap,"h\mapsto \alpha^{-1}\circ h\circ (k\mapsto k_{\mathcal{E}})"] & (\Z \to X)
    \end{tikzcd}
  \end{equation*}
  in which the top map is an equivalence.
\end{cor}

In other words, a family of maps $\universalcovercircle(t)\to A(t)$ indexed by $t:\sphere{1}$ is equivalently described as a map $h:\Z\to X$ for which the square
\begin{equation*}
  \begin{tikzcd}
    \Z \arrow[r,"h"] \arrow[d,swap,"\succZ"] & X \arrow[d,"e"] \\
    \Z \arrow[r,swap,"h"] & X
  \end{tikzcd}
\end{equation*}
commutes. It is now time to prove the universal property of the integers.
\index{circle!descent data|)}
\index{descent data for the circle|)}  

\subsection{The (dependent) universal property of the integers}
\index{Z@{$\Z$}!dependent universal property|(}
\index{Z@{$\Z$}!universal property|(}
\index{dependent universal property!of Z@{of $\Z$}|(}
\index{universal property!of Z@{of $\Z$}|(}

The dependent universal property precisely characterizes sections of families over the integers, for those families $A(k)$ indexed by $k:\Z$ that come equipped with families of equivalences $A(k)\simeq A(k+1)$ for all $k:\Z$.

\begin{lem}\label{lem:elim-Z}
Let $B$ be a family over $\Z$, equipped with an element $b_0:B(0)$, and an equivalence
\begin{equation*}
e_k : B(k)\eqvsym B(\succZ(k))
\end{equation*}
for each $k:\Z$. Then there is a dependent function $f:\prd{k:\Z}B(k)$ equipped with identifications $f(0)=b_0$ and
\begin{equation*}
f(\succZ(k))=e_k(f(k))
\end{equation*}
for any $k:\Z$.
\end{lem}

\begin{proof}
The map is defined using the induction principle for the integers, stated in \cref{lem:Z_ind}. First we take
\begin{align*}
f(-1) & \defeq e_{-1}^{-1}(b_0) \\
f(0) & \defeq b_0 \\
f(1) & \defeq e_0(b_0).
\end{align*}
For the induction step on the negative integers we use
\begin{equation*}
\lam{n}e_{\inneg(\succN(n))}^{-1} : \prd{n:\N} B(\inneg(n))\to B(\inneg(\succN(n)))
\end{equation*}
For the induction step on the positive integers we use
\begin{equation*}
\lam{n}e_{\inpos(n)} : \prd{n:\N} B(\inpos(n))\to B(\inpos(\succN(n))).
\end{equation*}
The computation rules follow in a straightforward way from the computation rules of $\Z$-induction and the fact that $e^{-1}$ is an inverse of $e$. 
\end{proof}

\begin{eg}
For any type $A$, we obtain a map $f:\Z\to A$ from any $x:A$ and any equivalence $e:\eqv{A}{A}$, such that $f(0)=x$ and the square
\begin{equation*}
\begin{tikzcd}
\Z \arrow[d,swap,"\succZ"] \arrow[r,"f"] & A \arrow[d,"e"] \\
\Z \arrow[r,swap,"f"] & A
\end{tikzcd}
\end{equation*}
commutes. In particular, if we take $A\defeq (x=x)$ for some $x:X$, then for any $p:x=x$ we have the equivalence $\lam{q}\ct{p}{q}:(x=x)\to (x=x)$. This equivalence induces a map
\begin{equation*}
k\mapsto p^k : \Z \to (x=x),
\end{equation*}
for any $p:x=x$. This induces the \define{degree $k$ map}\index{degree k map@{degree $k$ map}|textbf}\index{circle!degree k map@{degree $k$ map}|textbf} on the circle\index{deg(k)@{$\mathsf{deg}(k)$}|textbf}\index{deg(k)@{$\mathsf{deg}(k)$}|see {degree $k$ map}}\index{circle!deg(k)@{$\mathsf{deg}(k)$}|textbf}\index{circle!deg(k)@{$\mathsf{deg}(k)$}|see {degree $k$ map}}
\begin{equation*}
\mathsf{deg}(k) : \sphere{1}\to\sphere{1},
\end{equation*}
for any $k:\mathbb{Z}$, see \cref{ex:circle_degk}.
\end{eg}

In the following proposition we show that the dependent function constructed in \cref{lem:elim-Z} is unique. This is the \define{dependent universal property of the integers}\index{dependent universal property!of Z@{of $\Z$}|textbf}\index{Z@{$\Z$}!dependent universal property|textbf}.

\begin{prp}\label{prp:unique-elim-Z}
  Consider a type family $B:\mathbb{Z}\to\UU$ equipped with $b:B(0)$ and a family of equivalences
  \begin{equation*}
    e:\prd{k:\Z} \eqv{B(k)}{B(\succZ (k))}.
  \end{equation*}
  Then the type
  \begin{equation*}
    \sm{f:\prd{k:\Z}B(k)}(f(0)=b)\times\prd{k:\Z}f(\succZ (k))=e_k(f(k))
  \end{equation*}
  is contractible.
\end{prp}

\begin{proof}
  In \cref{lem:elim-Z} we have already constructed an element of the asserted type.
  Therefore it suffices to show that any two elements of this type can be identified.
  Note that the type $(f,p,H)=(f',p',H')$ is equivalent to the type of triples $(K,\alpha,\beta)$ consisting of
  \begin{align*}
    K & : f\htpy f' \\
    \alpha & : K(0)=\ct{p}{(p')^{-1}} \\
    \beta & : \prd{k:\Z}K(\succZ (k))=\ct{(\ct{H(k)}{\ap{e_k}{K(k)}})}{H'(k)^{-1}}.
  \end{align*}
  We obtain such a triple by applying \cref{lem:elim-Z} to the family $C$ over $\Z$ given by $C(k)\defeq f(k)=f'(k)$, which comes equipped with the base point
  \begin{equation*}
    \ct{p}{(p')^{-1}} : C(0),
  \end{equation*}
  and the family of equivalences
  \begin{equation*}
    \prd{k:\Z}\eqv{C(k)}{C(\succZ (k))}
  \end{equation*}
  given by $r\mapsto \ct{(\ct{H(k)}{\ap{e_k}{r}})}{H'(k)^{-1}}$.
\end{proof}

The \define{universal property of the integers}\index{universal property!of Z@{of $\Z$}|textbf}\index{Z@{$\Z$}!universal property|textbf} is a simple corollary of the dependent universal property. One way of phrasing it is that $\Z$ is the \emph{initial type equipped with a point and an automorphism}\index{integers!initial type with a point and an automorphism}.

\begin{cor}
  For any type $X$ equipped with a base point $x_0:X$ and an automorphism $e:\eqv{X}{X}$, the type
  \begin{equation*}
    \sm{f:\Z\to X}(f(0)=x_0)\times ((f \circ \succZ )\htpy(e\circ f))
  \end{equation*}
  is contractible.
\end{cor}

Using the fact that equivalences are contractible maps, we can reformulate the dependent universal property of the integers as follows.

\begin{thm}
  For any type family $A$ over $\Z$ equipped with a family of equivalences
  \begin{equation*}
    e:\prd{k:\Z}A(k)\simeq A(\succZ(k)),
  \end{equation*}
  the map
  \begin{equation*}
    \ev_0:\Big(\sm{f:\prd{k:\Z}A(k)}\prd{k:\Z}f(\succZ(k))=e_k(f(k))\Big)\to A(0)
  \end{equation*}
  given by $(f,H)\mapsto f(0)$ is an equivalence.
\end{thm}

\begin{proof}
  Note that the fibers of $\ev_0$ are equivalent to the types that are shown to be contractible in \cref{prp:unique-elim-Z}.
\end{proof}

The following corollary will be used to prove that the fundamental cover of the circle is equivalent to the identity type based at $\base:\sphere{1}$.

\begin{cor}\label{cor:universal-property-Z}
  For any type $X$ equipped with an equivalence $e:X\simeq X$, the map
  \begin{equation*}
    \Big(\sm{f:\Z\to X} f\circ \succZ \htpy e\circ f\Big)\to X
  \end{equation*}
  given by $(f,H)\mapsto f(0)$ is an equivalence.
\end{cor}
\index{Z@{$\Z$}!dependent universal property|)}
\index{Z@{$\Z$}!universal property|)}
\index{dependent universal property!of Z@{of $\Z$}|)}
\index{universal property!of Z@{of $\Z$}|)}

\subsection{The fundamental group of the circle}
\index{characterization of identity type!of the circle|(}
\index{circle!characterization of identity type|(}

We have two goals remaining in this book. The first goal is to prove that the universal cover of the circle is an identity system at $\base:\sphere{1}$, in the sense of \cref{defn:identity-system}. Since the universal cover is a family of sets over the circle, this implies that the circle is a $1$-type.

\begin{thm}
  \label{thm:eq-circle}%
  The universal cover of the circle is an identity system at $\base:\sphere{1}$.\index{universal cover of S 1@{universal cover of $\sphere{1}$}!is an identity system}\index{identity system!universal cover of S 1@{universal cover of $\sphere{1}$}}
\end{thm}

\begin{proof}
  By \cref{ex:uniqueness-identity-type} it suffices to show that the map
  \begin{equation*}
    f\mapsto f(0_{\mathcal{E}}) : \Big(\prd{t:\sphere{1}}\universalcovercircle(t)\to A(t)\Big)\to A(\base)
  \end{equation*}
  is an equivalence for every type family $A$ over the circle. Note that we have a commuting triangle
  \begin{equation*}
    \begin{tikzcd}[column sep=5em]
      \Big(\prd{t:\sphere{1}}\universalcovercircle(t)\to A(t)\Big) \arrow[d] \arrow[dr,"f\mapsto f(0_{\mathcal{E}})"] \\
      \sm{h:\Z\to A(\base)}h\circ\succZ\htpy \tr_A(\lloop)\circ h \arrow[r,swap,"{(h,H)\mapsto h(0)}"] & A(\base)
    \end{tikzcd}
  \end{equation*}
  in which the left map is the equivalence obtained in \cref{cor:compute-families-of-maps-universal-cover} and the bottom map is an equivalence by \cref{cor:universal-property-Z}.
\end{proof}

\begin{cor}
  The circle is a $1$-type and not a $0$-type.\index{circle!is a 1-type@{is a $1$-type}}\index{circle!is not a set}
\end{cor}

\begin{proof}
  To see that the circle is a $1$-type we have to show that $s=t$ is a $0$-type for every $s,t:\sphere{1}$. By \cref{ex:circle-connected} it suffices to show that the loop space of the circle is a $0$-type. This is indeed the case, because $\Z$ is a $0$-type, and we have an equivalence $(\base=\base)\simeq \Z$.

  Furthermore, since $\Z$ is a $0$-type and not a $(-1)$-type, it follows that the circle is a $1$-type and not a $0$-type.
\end{proof}

Our second goal is to construct a group isomorphism
\begin{equation*}
  \pi_1(\sphere{1})\cong \Z.
\end{equation*}
However, \cref{thm:eq-circle} doesn't immediately show that the fundamental group of the circle is $\Z$. It only gives us an equivalence
\begin{equation*}
  \loopspace{\sphere{1}}\simeq \Z.
\end{equation*}
In order to compute the fundamental group of the circle we augment the fundamental theorem of identity types with the following proposition.\index{fundamental theorem of identity types}

\begin{prp}\label{prp:fundamental-theorem-id-with-operation}
  Consider a type $A$ equipped with a point $a:A$, and consider an identity system\index{identity system} $B$ on $A$ at $a$ equipped with $b:B(a)$. Furthermore, suppose that there is a binary operation
  \begin{equation*}
    \mu:B(a)\to (B(x)\to B(x))
  \end{equation*}
  for every $x:A$, equipped with a homotopy $\mu(\blank,b)\htpy \idfunc$. Then we have
  \begin{equation*}
    f(\ct{p}{q})=\mu(f(p),f(q))
  \end{equation*}
  for the unique family of maps
  \begin{equation*}
    f:\prd{x:A}(a=x)\to B(x)
  \end{equation*}
  such that $f(\refl{})=b$, and for every $p:a=a$ and $q:a=x$.
\end{prp}

\begin{proof}
  Consider a family of maps $f:(a=x)\to B(x)$ indexed by $x:A$ such that $f(\refl{})=b$, and let $p:a=a$ and $q:a=x$. By induction on $q$ it suffices to show that
  \begin{equation*}
    f(p)=\mu(f(p),f(\refl{}))
  \end{equation*}
  This follows, since $f(\refl{})=b$ and $\mu(f(p),b)=f(p)$.
\end{proof}

We are now ready to prove that the fundamental group of the circle is $\Z$. Recall from \cref{defn:universal-cover-circle} that we write $y\mapsto y_\Z$ for the inverse of the equivalence
\begin{equation*}
  x\mapsto x_{\mathcal{E}}:\Z\simeq\universalcovercircle(\base).
\end{equation*}

\begin{thm}\label{thm:fundamental-group-circle}
  There is a group isomorphism\index{circle!p 1 S 1 cong Z@{$\pi_1(\sphere{1})\cong\Z$}}\index{p  1 S 1 cong Z@{$\pi_1(\sphere{1})\cong \Z$}}\index{circle!fundamental group}\index{fundamental group!of the circle}
  \begin{equation*}
    \pi_1(\sphere{1})\cong \Z.
  \end{equation*}
\end{thm}

\begin{proof}
  First we observe that, since the circle is a $1$-type, we have an isomorphism of groups $\pi_1(\sphere{1})\cong\loopspace{\sphere{1}}$. In order to show that the group $\loopspace{\sphere{1}}$ is isomorphic to $\Z$, we prove that the family of equivalences
  \begin{equation*}
    \alpha:\prd{t:\sphere{1}} (\base=t)\to \universalcovercircle(t)
  \end{equation*}
  given by $\alpha(\refl{})\defeq 0_{\mathcal{E}}$ satisfies
  \begin{equation*}
    \alpha(\ct{p}{q})_\Z=\alpha(p)_\Z+\alpha(q)_\Z
  \end{equation*}
  for every $p,q:\loopspace{\sphere{1}}$.
  
  To see that the claim holds, note that by \cref{prp:fundamental-theorem-id-with-operation} it suffices to construct a binary operation
  \begin{equation*}
    \mu : \universalcovercircle(\base)\to(\universalcovercircle(x)\to\universalcovercircle(x))
  \end{equation*}
  equipped with a homotopy $\mu(\blank,0_{\mathcal{E}})\htpy\idfunc$, such that
  \begin{equation*}
    \mu(k_{\mathcal{E}},l_{\mathcal{E}})=(k+l)_{\mathcal{E}}
  \end{equation*}
  holds for every $k,l:\Z$. Equivalently, it suffices to construct for each $k:\Z$ a function
  \begin{equation*}
    \mu(k_{\mathcal{E}}):\universalcovercircle(x)\to\universalcovercircle(x)
  \end{equation*}
  indexed by $x:\sphere{1}$ equipped with an identification $\mu(k_{\mathcal{E}},l_{\mathcal{E}})=(k+l)_{\mathcal{E}}$ for each $k,l:\Z$. Since we have
  \begin{equation*}
    k+(l+1)=(k+l)+1
  \end{equation*}
  for all $k,l:\Z$, such a function is obtained at once from \cref{cor:compute-families-of-maps-universal-cover}.
\end{proof}

In order to prove that the fundamental group of the circle is $\Z$, we first had to use the univalence axiom to construct the universal cover of the circle. This proof was originally discovered by Mike Shulman in 2011, and later published in \cite{LicataShulman}. Its importance of this proof to the field of homotopy type theory is hard to overestimate. The proof led to the discovery of the \emph{encode-decode method}, which we presented in this book as the fundamental theorem of identity types, and it was the start of the field that is now sometimes called \emph{synthetic homotopy theory}, where the induction principle for identity types and the univalence axiom are used along with methods from algebraic topology in order to compute algebraic invariants of types.

\index{characterization of identity type!of the circle|)}
\index{circle!characterization of identity type|)}

\begin{exercises}
  \exitem
  \begin{subexenum}
  \item Show that
    \begin{equation*}
      \prd{x:\sphere{1}}\brck{\base=x}.
    \end{equation*}
  \item On the other hand, use the universal cover of the circle to show that
    \begin{equation*}
      \neg\Big(\prd{x:\sphere{1}}\base=x\Big).
    \end{equation*}
  \item Use the circle to conclude that
    \begin{equation*}
      \neg\Big(\prd{X:\UU} \brck{X}\to X\Big).
    \end{equation*}
  \end{subexenum}
  \exitem \label{ex:circle_degk}
\begin{subexenum}
\item Show that for every $x:X$, we have an equivalence
\begin{equation*}
\eqv{\Big(\sm{f:\sphere{1}\to X}f(\base)= x \Big)}{(x=x)}
\end{equation*}
\item Show that for every $t:\sphere{1}$, we have an equivalence
\begin{equation*}
\eqv{\Big(\sm{f:\sphere{1}\to \sphere{1}}f(\base)= t \Big)}{\Z}
\end{equation*}
The base point preserving map $f:\sphere{1}\to\sphere{1}$ corresponding to $k:\Z$ is the degree $k$ map\index{circle!degree k map@{degree $k$ map}}\index{degree k map@{degree $k$ map}} on the circle.
\item Show that for every $t:\sphere{1}$, we have an equivalence
\begin{equation*}
\eqv{\Big(\sm{e:\eqv{\sphere{1}}{\sphere{1}}}e(\base)= t \Big)}{\bool}
\end{equation*}
\end{subexenum}
\exitem \label{ex:circle_double_cover} The \define{(twisted) double cover}\index{twisted double cover of S 1@{twisted double cover of $\sphere{1}$}|textbf}\index{circle!twisted double cover|textbf}\index{double cover of S 1@{double cover of $\sphere{1}$}|textbf}\index{circle!double cover|textbf} of the circle is defined as the type family $\mathcal{T}\defeq\mathcal{D}(\bool,\negbool):\sphere{1}\to\UU$, where $\negbool:\eqv{\bool}{\bool}$ is the negation equivalence of \cref{ex:neg_equiv}.
\begin{subexenum}
\item Show that $\neg(\prd{t:\sphere{1}}\mathcal{T}(t))$.
\item Construct an equivalence $e:\eqv{\sphere{1}}{\sm{t:\sphere{1}}\mathcal{T}(t)}$ for which the triangle
\begin{equation*}
\begin{tikzcd}[column sep=tiny]
\sphere{1} \arrow[rr,"e"] \arrow[dr,swap,"\mathsf{deg}(2)"] & & \sm{t:\sphere{1}}\mathcal{T}(t) \arrow[dl,"\proj 1"] \\
\phantom{\sm{t:\sphere{1}}\mathcal{T}(t)} & \sphere{1}
\end{tikzcd}
\end{equation*}
commutes.
\end{subexenum}
\exitem Construct an equivalence $\eqv{(\eqv{\sphere{1}}{\sphere{1}})}{\sphere{1}+\sphere{1}}$ for which the triangle
\begin{equation*}
  \begin{tikzcd}
    (\eqv{\sphere{1}}{\sphere{1}}) \arrow[rr,"\simeq"] \arrow[dr,swap,"\evbase"] & & (\sphere{1}+\sphere{1}) \arrow[dl,"\fold"] \\
    & \sphere{1}
  \end{tikzcd}
\end{equation*}
commutes. Conclude that a univalent universe containing a circle is not a $1$-type.
\exitem \label{ex:is_invertible_id_S1}
\begin{subexenum}
\item Construct a family of equivalences
\begin{equation*}
\prd{t:\sphere{1}} \big(\eqv{(t=t)}{\Z}\big).
\end{equation*}
\item Use \cref{ex:circle_connected} to show that $\eqv{(\idfunc[\sphere{1}]\htpy\idfunc[\sphere{1}])}{\Z}$.
\item Use \cref{ex:idfunc_autohtpy} to show that
\begin{equation*}
\eqv{\mathsf{has\usc{}inverse}(\idfunc[\sphere{1}])}{\Z},
\end{equation*}
and conclude that ${\mathsf{has\usc{}inverse}}(\idfunc[\sphere{1}])\not\simeq{\isequiv(\idfunc[\sphere{1}])}$. 
\end{subexenum}
\exitem Consider a map $i:A \to \sphere{1}$, and assume that $i$ has a retraction. Construct an element of type
  \begin{equation*}
    \iscontr(A)+\isequiv(i).
  \end{equation*}
  \exitem
  \begin{subexenum}
  \item Show that the multiplicative operation on the circle is associative\index{associativity!of multiplication on S 1@{of multiplication on $\sphere{1}$}}\index{circle!associativity of multiplication}, i.e.~construct an identification
    \begin{equation*}
      \assoc_{\sphere{1}}(x,y,z) :
      \mulcircle(\mulcircle(x,y),z)=\mulcircle(x,\mulcircle(y,z))
    \end{equation*}
    for any $x,y,z:\sphere{1}$.
  \item Show that the associator satisfies unit laws, in the sense that the following triangles commute:
    \begin{equation*}
      \begin{tikzcd}[column sep=-1em]
        \mulcircle(\mulcircle(\base,x),y) \arrow[rr,equals] \arrow[dr,equals] & & \mulcircle(\base,\mulcircle(x,y)) \arrow[dl,equals] \\
        & \mulcircle(x,y)
      \end{tikzcd}
    \end{equation*}
    \begin{equation*}
      \begin{tikzcd}[column sep=-1em]
        \mulcircle(\mulcircle(x,\base),y) \arrow[rr,equals] \arrow[dr,equals] & & \mulcircle(x,\mulcircle(\base,y)) \arrow[dl,equals] \\
        & \mulcircle(x,y)
      \end{tikzcd}
    \end{equation*}
    \begin{equation*}
      \begin{tikzcd}[column sep=-1em]
        \mulcircle(\mulcircle(x,y),\base) \arrow[rr,equals] \arrow[dr,equals] & & \mulcircle(x,\mulcircle(y,\base)) \arrow[dl,equals] \\
        & \mulcircle(x,y).
      \end{tikzcd}
    \end{equation*}
  \item State the laws that compute
    \begin{align*}
      & \assoc_{\sphere{1}}(\base,\base,x) \\
      & \assoc_{\sphere{1}}(\base,x,\base) \\
      & \assoc_{\sphere{1}}(x,\base,\base) \\
      & \assoc_{\sphere{1}}(\base,\base,\base).
    \end{align*}
    Note: the first three laws should be $3$-cells and the last law should be a $4$-cell. The laws are automatically satisfied, since the circle is a $1$-type.
  \end{subexenum}
  \exitem For convenience, we will write $x\cdot_{\sphere{1}}y\defeq\mulcircle(x,y)$ in this exercise. Construct the \define{Mac Lane pentagon}\index{Mac Lane pentagon}\index{circle!Mac Lane pentagon} for the circle, i.e.~show that the pentagon
  \begin{equation*}
    \begin{tikzcd}[column sep=-2em]
      &[-6em] ((x\cdot_{\sphere{1}} y)\cdot_{\sphere{1}} z)\cdot_{\sphere{1}} w \arrow[rr,equals] \arrow[dl,equals] & & (x\cdot_{\sphere{1}} y)\cdot_{\sphere{1}} (z\cdot_{\sphere{1}} w) \arrow[dr,equals] &[-6em] \\
      (x\cdot_{\sphere{1}} (y\cdot_{\sphere{1}} z))\cdot_{\sphere{1}} w \arrow[drr,equals] & & & & x\cdot_{\sphere{1}} (y\cdot_{\sphere{1}} (z \cdot_{\sphere{1}} w)) \\
      & & x\cdot_{\sphere{1}} ((y\cdot_{\sphere{1}} z)\cdot_{\sphere{1}} w) \arrow[urr,equals]
    \end{tikzcd}
  \end{equation*}
  commutes for every $x,y,z,w:\sphere{1}$.
  \exitem Recall from \cref{ex:surjective-precomp} that if $f:A\to B$ is a surjective map, then the precomposition map
  \begin{equation*}
    \blank\circ f : (B\to C)\to (A\to C)
  \end{equation*}
  is an embedding for every set $C$. 
  Give an example of a surjective map $f:A\to B$, such that the precomposition function
  \begin{equation*}
    \blank\circ f:(B\to \sphere{1})\to (A\to \sphere{1})
  \end{equation*}
  is \emph{not} an embedding, showing that the condition that $C$ is a set is essential.
  \exitem In this exercise we give an alternative proof that the total space of $\universalcovercircle$\index{universal cover of S 1@{universal cover of $\sphere{1}$}!total space} is contractible.
  \begin{subexenum}
  \item Construct a function
    \begin{equation*}
      h : \prd{k:\Z}(\base,0_{\mathcal{E}})=(\base,k_{\mathcal{E}})
    \end{equation*}
    equipped with a homotopy
    \begin{equation*}
      H : \prd{k:\Z}h(\succZ (k)_{\mathcal{E}})=\ct{h(k)}{\segmenthelix(k)}.
    \end{equation*}
  \item Show that the total space $\sm{t:\sphere{1}}\universalcovercircle(t)$ of the universal cover of the circle is contractible.
  \end{subexenum}
  \exitem Consider the type $\mathcal{C}$ of families $A:\sphere{1}\to\Set$ of sets over the circle equipped with a point $a_0:A(\base)$, such that the total space
  \begin{equation*}
    \sm{t:\sphere{1}}A(t)
  \end{equation*}
  is connected.
  \begin{subexenum}
  \item For any type family $A$ over the circle equipped with $a_0:A(\base)$, show that the total space $\sm{t:\sphere{1}}A(t)$ is connected if and only if $\tr_A(\lloop):A(\base)\to A(\base)$ has a single orbit in the sense that the map $k\mapsto \tr_A(\lloop)^k(a_0):\Z\to A(\base)$ is surjective. 
  \item Let $(A,a_0)$ and $(B,b_0)$ be in $\mathcal{C}$. Show that the type
    \begin{equation*}
      ((A,a_0)\leq (B,b_0))\defeq \sm{f:\prd{t:\sphere{1}}A(t)\to B(t)}f(\base,a_0)=b_0
    \end{equation*}
    is a proposition. Furthermore, show that this inequality relation gives $\mathcal{C}$ the structure of a poset.
  \item Show that the poset $\mathcal{C}$ is isomorphic to the poset of subgroups of $\Z$.
  \end{subexenum}
\end{exercises}

\index{circle!universal cover|)}
\index{universal cover of S 1@{universal cover of $\sphere{1}$}|)}
\index{circle|)}
\index{inductive type!circle|)}



\backmatter

\printindex

\printbibliography

\end{document}